\documentclass[10pt]{amsbook}

\usepackage[all]{xy}

\usepackage{latexsym}
\usepackage{amsmath}
\usepackage{amssymb}
\usepackage{amsfonts}

\makeindex

\numberwithin{subsection}{section}

\numberwithin{section}{chapter}

\numberwithin{chapter}{part}

\newtheorem{thm}{Theorem}[section]
\newtheorem{prop}[thm]{Proposition}
\newtheorem{lem}[thm]{Lemma}
\newtheorem{sublem}[thm]{Sub-lemma}
\newtheorem{df}[thm]{Definition}
\newtheorem{cor}[thm]{Corollary}

\newtheorem{ass}[thm]{Assumption}
\newtheorem{rmk}[thm]{Remark}
\newtheorem{ex}[thm]{Example}
\newtheorem{Q}[thm]{Question}
\newtheorem{conv}[thm]{Convention}

\author[B. To\"en]{Bertrand To\"en}

\address{Laboratoire Emile Picard UMR CNRS 5580 \\Universit\'{e} Paul Sabatier, Toulouse}

\email{toen@picard.ups-tlse.fr}

\author[G. Vezzosi]{Gabriele Vezzosi}

\address{Dipartimento di Matematica Applicata
``G. Sansone'' \\ Universit\`a di Firenze }

\email{gabriele.vezzosi@unifi.it}

\title[HAGII]{Homotopical Algebraic Geometry II:\\
geometric stacks and applications}

\begin{document}

\maketitle

\begin{quote}\centering{This work is dedicated to Alexandre Grothendieck}\end{quote}

\frontmatter

\setcounter{page}{6}

\tableofcontents

\bigskip
\bigskip

\chapter*{Abstract}\footnote{2000 Mathematics Subject Classification: 14A20, 18G55, 18F10, 55U40, 55P42, 55P43, 18F20, 18D10, 18E30 18G35, 18G30, 13D10, 55N34.\\
Keywords: Algebraic stacks, higher stacks, derived algebraic geometry\\ Received by the Editor June 14 2005}
This is the second part of a series of papers called ``HAG'', and
devoted to develop the foundations of \textit{homotopical algebraic
geometry}. We start by defining and studying generalizations of standard
notions of linear algebra in an abstract monoidal model category, such
as derivations, \'etale and smooth morphisms, flat and projective modules, etc.
We then use our theory of stacks over model categories, 
introduced in \cite{hagI}, in order to define
a general notion of geometric stack over a base symmetric monoidal model
category $C$, and prove that this notion satisfies the expected properties.

The rest of the paper consists in specializing $C$ in order to give
various examples of applications in several different contexts.
First of all, when $C=k-Mod$ is the category of $k$-modules with the
trivial model structure, we show that our notion gives back the
algebraic $n$-stacks of C. Simpson. Then
we set $C=sk-Mod$, the model category of simplicial $k$-modules, and obtain
this way a notion of geometric $D^{-}$\emph{-stack} which are the main geometric
objects of \emph{derived algebraic geometry}. We give several
examples of derived version of classical moduli stacks, as 
the $D^{-}$-stack of local systems on a space, the $D^{-}$-stack
of algebra structures over an operad, the $D^{-}$-stack of flat bundles
on a projective complex manifold, etc. We also present
the cases where $\mathcal{C}=C(k)$ is the model category of unbounded complexes of
$k$-modules, and $C=Sp^{\Sigma}$ the model category of symmetric spectra.
In these two contexts we give some examples of geometric stacks such as
the stack of associative dg-algebras, the stack of dg-categories, and
a geometric stack constructed using topological modular forms.\\

\bigskip
\bigskip
\bigskip

\rightline{\scriptsize{\texttt{There are more things in heaven and earth, Horatio,}}}
\rightline{\scriptsize{\texttt{than are dreamt of in our philosophy. But come...}}}
\smallskip
\rightline{\scriptsize{W. Shakespeare, \textit{Hamlet},  Act 1, Sc. 5.}}
\bigskip
\bigskip

\bigskip
\bigskip

\rightline{\scriptsize{\texttt{Mon cher Cato, il faut en convenir, les
forces de l'ether nous p\'en\`etrent,}}}
\rightline{\scriptsize{\texttt{et ce fait d\'eliquescent il nous faut
l'appr\'ehender coute que coute.}}}
\smallskip
\rightline{\scriptsize{P. Sellers, \textit{Quand la  panth\`ere rose
s'en m\^ele.}}}
\bigskip
\bigskip
\bigskip

\mainmatter

\chapter*{Introduction}

This is the second part of a series of papers called ``HAG'',
devoted to start the development of \textit{homotopical algebraic
geometry}. The first part \cite{hagI} was concerned with the homotopical generalization
of sheaf theory, and contains the notions of model topologies, model sites, stacks
over model sites and model topoi, all of these being homotopical versions of the usual notions
of Grothendieck topologies, sites, sheaves and topoi.
The purpose of the present
work is to use these new concepts in some
specific situations, in order to introduce a very general notion of
\emph{geometric stacks}, a far reaching homotopical generalization of
the notion of algebraic stacks introduced by P. Deligne, D. Mumford and M. Artin.
This paper includes the general study and the
standard properties of geometric stacks, as well as various examples
of applications in the contexts of algebraic geometry and algebraic topology. \\

\bigskip

\subsection*{Reminders on abstract algebraic geometry}

A modern point of view on algebraic geometry consists of viewing algebraic varieties and schemes
through their functors of points. In this
functorial point of view, schemes are certain sheaves of sets on
the category of commutative rings endowed with a suitable topology (e.g. the Zariski topology).
Keeping this in mind, it turns out that the whole theory of schemes can be completely reconstructed
starting from the symmetric monoidal category $\mathbb{Z}-Mod$ of $\mathbb{Z}$-modules alone. Indeed,
the category of commutative rings is reconstructed by taking the category of
commutative monoids in $\mathbb{Z}-Mod$. Flat morphisms can be recognized via the
exactness property of the base change functor on the category of modules. Finitely presented
morphisms are recognized via the usual categorical characterization in terms of
commutation of mapping sets with respect to filtered colimits. Finally, Zariski open
immersion can be defined as flat morphisms of finite presentation $A \longrightarrow B$ such that
$B\simeq B\otimes_{A}B$. Schemes are then reconstructed as certain Zariski sheaves on the
opposite category of commutative rings, which are obtained by gluing affine schemes
via Zariski open immersions (see for example the first chapters of \cite{dg}).

The fact that the notion of schemes has such a purely categorical interpretation has
naturally lead to the theory of \emph{relative algebraic geometry}, in which the
base symmetric monoidal category $\mathbb{Z}-Mod$ is replaced by an abstract
base symmetric monoidal category $\mathcal{C}$, and under reasonable assumptions on $\mathcal{C}$
the notion
of \emph{schemes over $\mathcal{C}$} can be made meaningful
as well as useful (see for example \cite{del1,hak}
for some applications).

The key observation of this work is that one can generalize further the theory of
relative algebraic geometry by requiring $\mathcal{C}$ to be endowed with
an additional \textit{model category} structure, compatible with its monoidal structure
(relative algebraic geometry is then recovered by taking the trivial
model structure), in such a way that the notions of \emph{schemes} and more
generally of \emph{algebraic spaces} or \emph{algebraic stacks} still have a natural and useful meaning,
compatible with the homotopy theory carried by $\mathcal{C}$.
In this work, we present this general theory, and show how this enlarges
the field of applicability by
investigating several examples not covered by the standard theory of relative algebraic geometry. The most important of these applications
is the existence of foundations for \emph{derived algebraic geometry}, 
a global counter part of the derived deformation theory of
V. Drinfel'd, M. Kontsevich and al. \\

\subsection*{The setting}

Our basic datum is a symmetric monoidal model category $\mathcal{C}$
(in the sense of \cite{ho}), on which certain conditions are imposed
(see assumptions \ref{ass-1}, \ref{ass1}, \ref{ass0} and \ref{ass2}).
We briefly discuss these requirements here.
The model category $\mathcal{C}$ is assumed to satisfy some reasonable
additional properties (as for example being proper,
or that cofibrant objects are flat for the monoidal structure). These
assumptions are only made for the convenience of certain constructions, and
may clearly be omitted. The model category $\mathcal{C}$ is also assumed to be
\emph{combinatorial} (see e.g. \cite{du2}), making it
reasonably behaved with respect to localization techniques.
The first really important assumption on $\mathcal{C}$ states that it is pointed (i.e. that the final
and initial object coincide) and that its homotopy category $\mathrm{Ho}(\mathcal{C})$ is additive.
This makes the model category $\mathcal{C}$ \emph{homotopically additive}, which is
a rather strong condition, but is used all along this work
and seems difficult to avoid (see however
\cite{sousz}). Finally, the last condition we make on $\mathcal{C}$ is
also rather strong, and states that the theory of commutative monoids in $\mathcal{C}$, and
the theory of modules over them, both possess reasonable model category structures.
This last condition is of course far from being satisfied in general (as for
example it is not satisfied when $\mathcal{C}$ is the model category of complexes over
some ring which is not of characteristic zero), but all the examples we have in mind
can be treated in this setting\footnote{Alternatively, one could switch to $E_{\infty}$-algebras (and modules over them) for which useful model and semi-model structures are known to exist, thanks to the work of M. Spitzweck \cite{sp}, in  much more general situations than for the case of commutative
monoids (and modules over them).}. The model categories of simplicial modules,
of complexes over a ring of characteristic zero, and of symmetric spectra are three
important examples of symmetric monoidal model category satisfying all our assumptions.
More generally, the model categories of sheaves with values in any of these
three fundamental categories provide additional examples.  \\

\subsection*{Linear and commutative
algebra in a symmetric monoidal model category}

An important consequence of our assumptions on the base symmetric monoidal model
category $\mathcal{C}$ is the existence of reasonable generalizations of
general constructions and results from standard linear and commutative algebra.
We have gathered some of these notions (we do not claim to be exhaustive) in
\S \ref{partI.1}. For example, we give definitions of derivations as well as of cotangent complexes
representing them, we define formally \'etale morphisms, flat morphisms, open Zariski
immersions, formally unramified morphisms, finitely presented morphism
of commutative monoids and modules, projective and flat modules, Hochschild cohomology, etc.
They are all generalizations of the well known
notions in the sense that when applied to the case where $\mathcal{C}=\mathbb{Z}-Mod$ with the
trivial model structure we find back the usual notions.
However, there are sometimes
several nonequivalent generalizations, as for example there exist at least two, nonequivalent
reasonable generalizations of smooth morphisms which restrict to the usual one when
$\mathcal{C}=\mathbb{Z}-Mod$. This is why we have tried to give an overview of several
possible generalizations, as we think all definitions could have their own interest
depending both on the context and on what one wants to do with them. Also we wish to mention that
all these notions depend heavily on the base model category $\mathcal{C}$, in the sense that
the same object in $\mathcal{C}$, when considered in different model categories
structures on $\mathcal{C}$, might not
behave the same way. For example, a commutative ring can also be considered
as a simplicial commutative ring, and the notion of finitely presented
morphisms is not the same in the two cases.
We think that keeping track of the base model category $\mathcal{C}$ is rather important,
since playing with the change of base categories might be very useful, and is also
an interesting feature of the theory.

The reader will immediately notice that several notions behave in a much better way when
the base model category satisfies certain stability assumptions (e.g. is
a stable model category, or when the suspension functor is fully faithful, see
for example Prop. \ref{p6}, Cor. \ref{c1}). We think this is one of the main
features of homotopical algebraic geometry: linear and commutative algebra notions
tend to be better behaved as the base model category tend to be ``more'' stable. We do not claim that
everything becomes simpler in the stable situation, but that certain difficulties
encountered can be highly simplified by enlarging the base model category to a more
stable one. \\

\subsection*{Geometric stacks}

In \S \ref{partI.3} we present the general notions of geometric stacks
relative to our base model category $\mathcal{C}$. Of course, we start by
defining $Aff_{\mathcal{C}}$, the model category of affine objects over $\mathcal{C}$, as the
opposite of the model category $Comm(\mathcal{C})$ of commutative monoids in $\mathcal{C}$.
We assume we are given a \emph{model (pre-)topology} $\tau$ on $Aff_{\mathcal{C}}$, in the
sense we have given to this expression in \cite[Def. 4.3.1]{hagI} (see also Def. \ref{modtop}). We also assume that
this model topology satisfies certain natural assumptions, as
quasi-compactness and the descent property for modules. The model category
$Aff_{\mathcal{C}}$ together with its model topology $\tau$ is a model site in the sense
of \cite[Def. 4.3.1]{hagI} or Def. \ref{modtop}, and it gives rise to a model category of stacks
$Aff_{\mathcal{C}}^{\sim,\tau}$. The homotopy category of $Aff_{\mathcal{C}}^{\sim,\tau}$ will simply be denoted by $\mathrm{St}(\mathcal{C}, \tau)$. The Yoneda embedding for model categories allows
us to embed the homotopy category $\mathrm{Ho}(Aff_{\mathcal{C}})$ into $\mathrm{St}(\mathcal{C}, \tau)$, 
and this gives a notion of representable stack, our analog of
the notion of affine scheme. Geometric stacks will result from
a certain kind of gluing representable stacks.

Our notion of geometric stack
is relative to a class of morphisms \textbf{P} in $Aff_{\mathcal{C}}$, satisfying
some compatibility conditions with respect to the topology $\tau$, essentially stating that the
notion of morphisms in \textbf{P} is local for the topology $\tau$.
With these two notions, $\tau$ and \textbf{P}, we define by induction on
$n$ a notion of $n$-geometric stack (see \ref{d11}). The precise definition
is unfortunately too long to be reproduced here, but one can
roughly say that $n$-geometric stacks are stacks $F$ whose diagonal
is $(n-1)$-representable (i.e. its fibers over representable stacks
are $(n-1)$-geometric stacks), and which admits a covering
by representable stacks $\coprod_{i}U_{i} \longrightarrow F$, such that
all morphisms $U_{i} \longrightarrow F$ are in \textbf{P}.

The notion of $n$-geometric stack satisfies all the expected basic properties.
For example, geometric stacks are stable by (homotopy) fiber products and disjoint unions,
and being an $n$-geometric stack is a local property (see Prop. \ref{p11}, \ref{p14}).
We also present a way to produce $n$-geometric stacks as certain quotients of groupoid actions, in the
same way that algebraic stacks (in groupoids) can always be presented as quotients
of a scheme by a smooth groupoid action (see Prop. \ref{p13}). When a property
\textbf{Q} of morphisms in $Aff_{\mathcal{C}}$ satisfies a certain compatibility
with both \textbf{P} and $\tau$, there exists a natural notion of \textbf{Q}-morphism between
$n$-geometric stacks, satisfying all the expected properties (see Def. \ref{d23}, Prop. \ref{p20}).
We define the stack of quasi-coherent modules, as well its sub-stacks of
vector bundles and of perfect modules (see Thm. \ref{t2}, Cor. \ref{ct2}).
These also behave as expected, and for example the stack of vector bundles
is shown to be $1$-geometric as soon as the class \textbf{P} contains the
class of smooth morphisms (see Cor. \ref{cp24}). \\

\subsection*{Infinitesimal theory}

In \S \ref{partI.4}, we investigate
the infinitesimal properties of geometric stacks. For this we define a notion of
derivation of a stack $F$ with coefficients in a module, and the notion of
cotangent complex is defined via the representability of the functor of derivations
(see Def. \ref{d16}, \ref{d17}, \ref{d18}). The object representing
the derivations, the cotangent complex, is not in general
an object in the base model category $\mathcal{C}$, but belongs to the
stabilization of $\mathcal{C}$ (this is of course related to the well known fact that
cotangent spaces of algebraic stacks are not vector spaces but rather
complexes of vector spaces). This is why these notions will be only defined
when the suspension functor of $\mathcal{C}$ is fully faithful, or equivalently when
the stabilization functor from $\mathcal{C}$ to its stabilization is fully faithful
(this is again an incarnation of the fact, mentioned above, that homotopical algebraic
geometry seems to prefer stable situations). In a way, this explains from a conceptual point of view
the fact that
the infinitesimal study of usual algebraic stacks in the sense of Artin is already
part of homotopical algebraic geometry, and does not really belong to standard
algebraic geometry. We also define stacks having an \emph{obstruction theory}
(see Def. \ref{d24}), a notion which controls obstruction to lifting morphisms along
a first order deformation in terms of the cotangent complex.
Despite its name, having an obstruction theory is a property of a stack and not
an additional structure. Again, this notion is really well behaved when the suspension
functor of $\mathcal{C}$ is fully faithful, and this once again explains the relevance
of derived algebraic geometry with respect to infinitesimal deformation theory.
Finally, in the last section we give sufficient conditions
(that we called \emph{Artin's conditions}) insuring that any
$n$-geometric stack has an obstruction theory (Thm. \ref{t1}). This last result
can be considered as a far reaching generalization of the
exixtence of cotangent complexes for algebraic stacks
as presented in \cite{lm}. \\

\subsection*{Higher Artin stacks (after C. Simpson)}

As a first example of application, we show how our general notion of
geometric stacks specializes to C. Simpson's algebraic $n$-stacks
introduced in \cite{s4}. For this, we let $\mathcal{C}=k-Mod$, be the symmetric
monoidal category of $k$-modules (for some fixed commutative ring $k$), endowed
with its trivial model structure. The topology $\tau$ is chosen to be
the \'etale (\'et) topology, and \textbf{P} is chosen to be
the class of smooth morphisms. We denote by $\mathrm{St}(k)$ the corresponding
homotopy category of \'et-stacks. Then, our definition of $n$-geometric stack
gives back the notion of algebraic $n$-stack introduced in \cite{s4}
(except that the two $n$'s might differ); these stacks will be called
\emph{Artin $n$-stacks} as they contain the usual algebraic stacks in the sense of Artin
as particular cases (see Prop. \ref{pII-2}). However, all our infinitesimal
study (cotangent complexes and obstruction theory) does not apply here as the
suspension functor on $k-Mod$ is the zero functor. This should not be viewed
as a drawback of the theory; on the contrary we rather think this explains why
deformation theory and obstruction theory in fact already belong to the realm of
derived algebraic geometry, which is our next application. \\

\subsection*{Derived algebraic geometry: $D^{-}$-stacks}

Our second application is the so-called \emph{derived algebraic
geometry}. The base model category $\mathcal{C}$ is chosen to be
$sk-Mod$, the symmetric monoidal model category of simplicial
commutative $k$-modules, $k$ being some fixed base ring. The category
of affine objects is $k-D^{-}Aff$, the opposite model category of
the category of commutative simplicial $k$-algebras.
In this setting, our general notions of flat, \'etale, smooth morphisms and
Zariski open immersions all have
explicit descriptions in terms of standard notions (see Thm. \ref{tII-1}).
More precisely, we prove that a morphism of simplicial commutative $k$-algebras
$A \longrightarrow B$ is flat (resp. smooth, resp. \'etale, resp. a Zariski open immersion)
in the general sense we have given to these notions in \S \ref{partI.2}, if and only if
it satisfies the following two conditions
\begin{itemize}
\item The induced morphism of affine schemes
$$Spec\, \pi_{0}(B) \longrightarrow Spec\, \pi_{0}(A)$$
is flat (resp. smooth, resp. \'etale, resp. a Zariski open immersion)
in the usual sense.
\item The natural morphism
$$\pi_{*}(A)\otimes_{\pi_{0}(A)}\pi_{0}(B) \longrightarrow
\pi_{*}(B)$$
is an isomorphism.
\end{itemize}

We endow $k-D^{-}Aff$ with the \'etale model
topology, a natural extension of the \'etale topology for affine schemes; the corresponding 
homotopy category of $D^{-}$-stacks is simply denoted by $D^{-}\mathrm{St}(k)$.
The class \textbf{P} is taken to be the class of smooth morphisms.
The $n$-geometric stacks in this context will be called
\emph{$n$-geometric $D^{-}$-stacks},
where the notation $D^{-}$ is meant to bring to mind the negative
bounded derived category\footnote{Recall that
the homotopy theory of simplicial $k$-modules is
equivalent to the homotopy theory of
negatively graded cochain complexes of $k$-modules. Therefore,
derived algebraic geometry can also be considered
as algebraic geometry over the category of negatively
graded complexes.}. An important consequence of the above descriptions
of \'etale and smooth morphisms is that the natural inclusion functor from
the category of $k$-modules to the category of simplicial $k$-modules, induces
a full embedding of the category of Artin $n$-stacks
into the category of $n$-geometric $D^{-}$-stack. This inclusion
functor $i$ has furthermore a right adjoint, called the truncation functor $t_{0}$ (see Def. \ref{dII-8}),
and the adjunction morphism $it_{0}(F) \longrightarrow F$ provides a closed embedding
of the classical Artin $n$-stack $it_{0}(F)$ to its derived version
$F$, which behaves like a formal thickening (see Prop. \ref{pII-6}).
This is a global counterpart of the common picture of derived deformation theory
of a formal classical moduli space sitting as a closed sub-space in
the corresponding formal derived moduli space.

We also prove that our general conditions for the
existence of an obstruction theory are satisfied, and so any $n$-geometric
$D^{-}$-stack has an obstruction theory (see Prop. \ref{cII-4}). An important particular case is when
this result is applied to the image of an Artin $n$-stack via the natural inclusion functor $i$; we
obtain in this way the existence of an obstruction theory for
\textit{any} Artin $n$-stack, and in particular the existence of a cotangent complex.
This is a very good instance of our principle that things simplifies
when the base model category $\mathcal{C}$ becomes more stable: the infinitesimal study of classical
objects of algebraic geometry (such as schemes, algebraic stacks or Artin $n$-stacks) 
becomes conceptually clearer and behaves much better when we consider these objects as geometric $D^{-}$-stacks.

Finally, we give several examples of $D^{-}$-stacks being
derived versions of some well known classical moduli problems.
First of all the $D^{-}$-stack of \textit{local systems} on a topological space, and the
$D^{-}$-stack of \textit{algebra structures} over a given operad are shown to be
$1$-geometric (see Lem. \ref{lII-11}, Prop. \ref{pII-15}).
We also present derived versions of the \textit{scheme of morphisms} between two projective schemes, and
of the moduli stack of \textit{flat bundles} on a projective complex manifold (Cor. \ref{cpII-17} and
Cor. \ref{cpII-17'}). The proofs that these
last two stacks are geometric rely on a special version of J. Lurie's representability theorem
(see \cite{lu} and Appendix C).\\

\subsection*{Complicial algebraic geometry: $D$-stacks}

What we call \emph{complicial algebraic geometry} is an unbounded version
of derived algebraic geometry in which the base model category
is $C(k)$ the category of complexes over some commutative ring $k$ (of
characteristic zero), and is presented in \S \ref{IIunb}.
It turns out that linear algebra
over $C(k)$ behaves rather differently than over the category of
simplicial $k$-modules (corresponding to complexes in non-positive degrees).
Indeed, the smooth, \'etale and Zariski open immersion can not
be described using a simple description on homotopy groups anymore. For
example, a usual ring $A$ may have Zariski open localizations
$A \longrightarrow A'$ in the context of complicial algebraic geometry such that 
$A'$ is not cohomologically concentrated in degree $0$ anymore. Also, 
a usual non affine scheme might be affine when considered 
as a scheme over $C(k)$: for example any quasi-compact open subscheme
of a usual affine scheme is representable by an affine scheme over $C(k)$ (see example
\ref{exudag}).

This makes the complicial theory rather different from
derived algebraic geometry for which the geometric intuition was instead 
quite close to the usual one, and constitutes a very interesting new feature
of complicial algebraic geometry.
The category $Aff_{\mathcal{C}}$ is here the opposite of the category of
unbounded commutative differential graded algebras over $k$. It is
endowed with a \emph{strong \'etale topology}, and the corresponding 
homotopy category of $D$-stacks is simply denoted by $D\mathrm{St}(k)$.
A new feature is here the existence of several interesting choices for the class
\textbf{P}. We will present two of them, one for which
\textbf{P} is taken to be the class of perfect morphisms, a rather
weak notion of smoothness, and a second one for which \textbf{P}
is taken to be the class of \emph{fip-smooth} morphism,
a definitely stronger notion behaving similarly
to usual smooth morphisms with respect to lifting properties.
We check that such choices satisfy the required properties in
order for $n$-geometric stacks (called
\emph{weakly $n$-geometric $D$-stacks} and
\emph{$n$-geometric $D$-stacks}, according to the choice of \textbf{P})
to make sense. Furthermore, for our second
choice, we prove that Artin's conditions are satisfied, and thus
that $n$-geometric $D$-stacks have a good infinitesimal
theory.
We give several examples of weakly geometric $D$-stacks,
the first one being the
$D$-stack of \textit{perfect modules} $\mathbf{Perf}$. We also show that the
$D$-stack of \textit{associative algebra structures} $\mathbf{Ass}$
is a weakly $1$-geometric $D$-stack. Finally,
the $D$-stack of \textit{connected} $dg$\textit{-categories} $\mathbf{Cat}_{*}$
is shown to  be weakly $2$-geometric. It is important to note that
these $D$-stacks can not be reasonably described
as geometric $D^{-}$-stacks, and provide
examples of truly ``exotic'' geometric objects.

Suitable slight modifications of the $D$-stacks
$\mathbf{Perf}$, $\mathbf{Ass}$ and $\mathbf{Cat}_{*}$ are given
and shown to be geometric. This allows us to
study their tangent complexes, and show in particular that
the infinitesimal theory of a certain class of dg-algebras and
dg-categories is controlled by derivations and
Hochschild cohomology, respectively (see Cor. \ref{ccwass} and Cor. \ref{ctdg-cat}).
We also show that
Hochschild cohomology does not control deformations of
general dg-categories in any reasonable sense (see Cor. \ref{ctdg-cat2} and Rem. \ref{rmkcat}). This has been a true surprise, as it contradicts
some of the statements one finds in the existing literature, including
some made by the authors themselves (see e.g. \cite[Thm. 5.6]{hagdag}).\\

\subsection*{Brave new algebraic geometry: $S$-stacks}

Our last context of application, briefly presented in \S \ref{IIbnag}, is
the one where the base symmetric monoidal model category is $\mathcal{C}=Sp^{\Sigma}$, the model category
of symmetric spectra (\cite{hss, shi}), and gives rise to what we call, after F. Waldhausen,
\emph{brave new algebraic geometry}. Like in the complicial case, 
the existence of negative homotopy groups makes the general theory
of flat, smooth, \'etale morphisms and of Zariski open immersions rather
different from the corresponding one in derived algebraic geometry. Moreover, typical phenomena coming
from the existence of Steenrod operations makes the notion of smooth morphism
even more exotic and rather different from algebraic geometry; to give just a striking example,
$\mathbb{Z}[T]$ is \textit{not} smooth over $\mathbb{Z}$ in the context
of brave new algebraic geometry. Once again, we do not think this
is a drawback of the theory, but rather an interesting new feature one should
contemplate and try to understand, as it might reveal interesting new insights
also on classical objects. In brave new algebraic geometry, we also check that the
strong \'etale topology and the class \textbf{P} of fip-smooth morphisms
satisfy our general assumptions, so that $n$-geometric stacks exists
in this context. We call them $n$-geometric $S$\textit{-stacks}, while the homotopy category of 
$S$-stacks for the strong \'etale topology is simply denoted by $\mathrm{St}(S)$. As an example, we give
a construction of a $1$-geometric $S$-stack starting from the ``sheaf'' of topological modular forms
(Thm. \ref{ttmf}). \\

\subsection*{Relations with other works}

It would be rather long to present all related works, and we
apologize in advance for not mentioning all of them.

The general fact that the notion of geometric stack only depends on
a topology and a choice of the class of morphisms \textbf{P} has already been
stressed by Carlos Simpson in \cite{s4}, who
attributes this idea to C. Walter.
Our general definition of geometric $n$-stacks is
a straightforward generalization to our abstract
context of the definitions found in \cite{s4}.

Originally, derived algebraic geometry have been approached
using the notion of \textit{dg-schemes}, as introduced by M. Kontsevich,
and developed by I. Ciocan-Fontanine and M. Kapranov.
We have not tried to make a full comparison with our theory. Let us only mention
that there exists a functor from dg-schemes to our category of
$1$-geometric $D^{-}$-stacks (see \cite[\S 3.3]{hagdag}). Essentially nothing is known
about this functor: we tend to believe that it is
not fully faithful, though this question does not seem very relevant. On the contrary,
the examples of $dg$-schemes constructed in \cite{ck1,ck2} do provide
examples of $D^{-}$-stacks and we think it is interesting to
look for derived moduli-theoretic interpretations of these (i.e. describe
their functors of points).

We would like to mention that an approach to  formal derived algebraic
geometry has been settled down by V. Hinich in \cite{hin2}. As far as
we know, this is the first functorial point of view on derived algebraic
geometry that appeared in the literature.

There is a big overlap between our Chapter \ref{IIder} and Jacob Lurie's thesis \cite{lu}.
The approach to derived algebraic geometry used by J. Lurie is different
from ours as it is based on a notion of $\infty$-category, whereas
we are working with model categories. The simplicial localization techniques
of Dwyer and Kan provide a way to pass from model categories to
$\infty$-categories, and the ``strictification'' theorem
of \cite[Thm. 4.2.1]{msri} can be used to see that our approach and
Lurie's approach are in fact equivalent
(up to some slight differences, for instance
concerning the notion of descent).
We think that the present
work and \cite{lu} can not be reasonably considered as totally independent, as
their authors have been frequently communicating on the subject since
the spring 2002. It seems rather clear that we all have been influenced
by these communications and that we have greatly benefited from
the reading of the first drafts of \cite{lu}. We have to mention however that
a huge part of the material of the present paper had been already announced
in earlier papers (see e.g. \cite{web,hagdag}), and have been worked
out since the summer 2000 at the time were our project has started.
The two works are also rather disjoint and complementary, as
\cite{lu} contains much more materials on derived algebraic
geometry than what we have included in \S \ref{IIder} (e.g. a wonderful
generalization of Artin's representability theorem, to state only the
most striking result). On the other hand, our ``HAG'' project has also been
motivated by rather exotic contexts of applications, as
the ones exposed for example in \S \ref{IIunb} and \S \ref{IIbnag} , and which are not covered
by the framework of \cite{lu}.

The work of K. Behrend on $dg$-schemes and $dg$-stacks \cite{be1,be2} has been
done while we were working on our project, and therefore
\S \ref{IIder} also has some overlaps with his work. However, the two approaches
are rather different and nonequivalent, as K. Behrend uses
a $2$-truncated version of our notions of stacks, with the effect
of killing some higher homotopical information. We have not investigated
a precise comparison between these two approaches in this work,
 but we would like to mention that there
exists a functor from our category of $D^{-}$-stacks to K. Behrend's category
of $dg$-sheaves. This functor is extremely badly behaved: it is not
full, nor faithful, nor essentially surjective, nor even injective
on isomorphism classes of objects. The only good property is that
it sends $1$-geometric Deligne-Mumford $D^{-}$-stacks to
Deligne-Mumford $dg$-stack in Behrend's sense.  However,
there are non geometric $D^{-}$-stacks that become geometric
objects in Behrend's category of $dg$-sheaves.

Some notions of \'etale and smooth morphisms of commutative $S$-algebras have been
introduced in \cite{min2}, and they seem to be related to the general
notions we present in \S \ref{partI.2}. However a precise comparison is not so easy.
Moreover,
\cite{min2} contains some wrong statements
like the fact that $thh$-smoothness generalizes
smoothness for discrete algebras (right after
Definition 4.2) or like Lemma 4.2 (2). The proof
of Theorem 6.1 also contains an important gap,
since the local equivalences at the end of the proof
are not checked to glue together.

Very recentely, J. Rognes has proposed a brave new version of
Galois theory, including brave new notions of \'etaleness
which are very close to our notions
(see \cite{ro}).

A construction of the moduli of dg-algebras and
dg-categories appears in \cite{kos}. These moduli are
only formal moduli by construction, and we propose
our $D$-stacks $\mathbf{Ass}$ and $\mathbf{Cat}_{*}$ as
their global geometrical counterparts.

We wish to mention the work of M. Spitzweck \cite{sp}, in which he
proves the existence of model category structures for $E_{\infty}$-algebras and
modules in a rather general context. This work can therefore be used in order
to suppress our assumptions on the existence of model category of commutative
monoids. Also, a nice symmetric monoidal model category of motivic complexes
is defined in \cite{sp}, providing a new interesting context to
investigate. It has been suggested to us to consider this
example of \emph{algebraic geometry over motives} by Yu. Manin, already during
spring 2000, but we do not have
at the moment interesting things to say on the subject.

Finally, J. Gorski  has recently constructed
a $D^{-}$-stack version of the Quot functor (see \cite{go}), providing
this way a functorial interpretation of the
derived Quot scheme of \cite{ck1}. A geometric
$D^{-}$-stack classifying objects in a dg-category
has been recently constructed by the first author
and M. Vaqui\'e in \cite{tv}. \\

\subsection*{Acknowledgments}

First of all, we are very grateful to C. Simpson for sharing with us all
his ideas on higher stacks, and for encouraging us to pursue our project.
We owe him a lot.

We are very grateful to J. Lurie for various communications on the subject, and
for sharing with us his work \cite{lu}. We have learned a lot about derived
algebraic geometry from him.

We thank M. Vaqui\'e for reading a former 
version of the present work, and for his comments and
suggestions. 

For various conversations on the subject we thank M. Anel,
J. Gorski,  A. Hirschowitz,
A. Joyal, M. Kontsevich, L. Katzarkov, H. Miller,
T. Pantev, C. Rezk, J. Rognes, S. Schwede, B. Shipley, M. Spitzweck,
N. Strickland and J. Tapia.

Finally, we thank both referees for their careful reading of the manuscript and
for their interesting and useful remarks and suggestions.

\subsection*{Notations and conventions}

We will use the word \textit{universe} in the sense of
\cite[Exp.I, Appendice]{sga4}. Universes will be denoted by $\mathbb{U} \in
\mathbb{V} \in \mathbb{W} \dots$. For any universe $\mathbb{U}$ we
will assume that $\mathbb{N} \in \mathbb{U}$. The category of sets
(resp. simplicial sets, resp. \dots) belonging to a universe
$\mathbb{U}$ will be denoted by $Set_{\mathbb{U}}$ (resp.
$SSet_{\mathbb{U}}$, resp. \dots). The objects of
$Set_{\mathbb{U}}$ (resp. $SSet_{\mathbb{U}}$, resp. \dots) will
be called $\mathbb{U}$-sets (resp. $\mathbb{U}$-simplicial sets,
resp. \dots). We will use the expression
\textit{$\mathbb{U}$-small set} (resp. \textit{$\mathbb{U}$-small
simplicial set}, resp. \dots) to mean \textit{a set isomorphic to
a set in $\mathbb{U}$} (resp. \textit{a simplicial set isomorphic
to a simplicial set in $\mathbb{U}$}, resp. \dots).
A unique exception concerns categories.
The expression \textit{$\mathbb{U}$-category} refers to the
usual notion of \cite[$I Def. 1.2$]{sga4}, and denotes a category $\mathcal{C}$ such that
for any two objects $x$ and $y$ in $\mathcal{C}$ the set $Hom_{\mathcal{C}}(x,y)$ is $\mathbb{U}$-small.
In the same way, a category $\mathcal{C}$ is \textit{$\mathbb{U}$-small} is it is isomorphic
to some element in $\mathbb{U}$.

Our references for model categories are \cite{ho} and \cite{hi}.
By definition, our model categories will always be \textit{closed}
model categories, will have all \textit{small} limits and colimits
and the functorial factorization property. The word
\textit{equivalence} will always mean \textit{weak equivalence}
and will refer to a model category structure.

The homotopy category of a model category $M$ is $W^{-1}M$ (see
\cite[Def. $1.2.1$]{ho}), where $W$ is the subcategory of
equivalences in $M$, and it will be denoted as $\mathrm{Ho}(M)$.
The sets of morphisms in $\mathrm{Ho}(M)$ will be denoted by
$[-,-]_{M}$, or simply by $[-,-]$ when the reference to the model
category $M$ is clear. We will say that two objects in a model
category $M$ are equivalent if they are isomorphic in
$\mathrm{Ho}(M)$. We say that two model categories are
\textit{Quillen equivalent} if they can be connected by a finite
string of Quillen adjunctions each one being a Quillen
equivalence. The mapping space of morphisms between two objects
$x$ and $y$ in a model category $M$ is denoted by $Map_{M}(x,y)$
(see \cite[\S 5]{ho}),
or simply $Map(x,y)$ if the reference to $M$ is clear. The simplicial
set depends on the choice of cofibrant and fibrant resolution
functors, but is well defined as an object in the
homotopy category of simplicial sets $\mathrm{Ho}(SSet)$.
If $M$ is a $\mathbb{U}$-category, then $Map_{M}(x,y)$ is
a $\mathbb{U}$-small simplicial set.

The homotopy fiber product (see \cite[$13.3$, $19.5$]{hi}, \cite[Ch.
XIV]{dkh} or \cite[10]{ds} ) of a diagram $\xymatrix{x \ar[r] & z & \ar[l] y}$ in a
model category $M$ will be denoted by $x\times_{z}^{h}y$. In the
same way, the homotopy push-out of a diagram $\xymatrix{x & \ar[l]
z  \ar[r] & y}$ will be denoted by $x\coprod_{z}^{\mathbb{L}}y$.
For a pointed model category $M$, the suspension and loop functors
functor will be denoted by
$$S : \mathrm{Ho}(M) \longrightarrow \mathrm{Ho}(M) \qquad \mathrm{Ho}(M) \longleftarrow
\mathrm{Ho}(M) : \Omega.$$
Recall that $S(x):=*\coprod_{x}^{\mathbb{L}}*$, and
$\Omega(x):=*\times_{x}^{h}*$.

When a
model category $M$ is a simplicial model category, its simplicial
sets of morphisms will be denoted by $\underline{Hom}_{M}(-,-)$, and
their derived functors by $\mathbb{R}\underline{Hom}_{M}$ (see
\cite[1.3.2]{ho}), or simply $\underline{Hom}(-,-)$ and
$\mathbb{R}\underline{Hom}(-,-)$ when the reference to $M$ is clear.
When $M$ is a symmetric monoidal
model category in the sense of \cite[\S 4]{ho}, the derived
monoidal structure will be denoted by
$\otimes^{\mathbb{L}}$.

For the notions of $\mathbb{U}$-cofibrantly generated,
$\mathbb{U}$-combinatorial and $\mathbb{U}$-cellular model
category, we refer to \cite{ho, hi, du2}, or to 
 \cite[Appendix]{hagI}, where
the basic definitions and crucial properties are recalled in a way
that is suitable for our needs.

As usual, the standard simplicial category will be denoted by
$\Delta$. The category of simplicial objects
in a category $\mathcal{C}$ will be denoted by $s\mathcal{C}:=\mathcal{C}^{\Delta^{op}}$.
In the same way, the category of co-simplicial objects
in $\mathcal{C}$ will be denoted by $cs\mathcal{C}$.
For any simplicial object $F \in s\mathcal{C}$ in a
category $\mathcal{C}$, we will use the notation $F_{n}:=F([n])$. Similarly,
for any co-simplicial object $F \in \mathcal{C}^{\Delta}$, we will
use the notation $F_{n}:=F([n])$. Moreover, when
$\mathcal{C}$ is a model category, we will use the notation
$$|X_{*}|:=Hocolim_{[n]\in \Delta^{op}}X_{n}$$
for any $X_{*}\in s\mathcal{C}$.

A sub-simplicial set $K\subset L$ will be called
\emph{full} is $K$ is a union of connected components
of $L$.  We will also say that a morphism
$f : K \longrightarrow L$ of simplicial sets is
\emph{full} if it induces an
equivalence bewteen $K$ and a full sub-simplicial set of $L$.
In the same way, we will use the expressions \emph{full sub-simplicial presheaf},
 and \emph{full morphisms of simplicial presheaves} for the
levelwise extension of the above notions to presheaves of simplicial sets.

For a Grothendieck site $(\mathcal{C},\tau)$ in a universe $\mathbb{U}$, we
will denote by $Pr(\mathcal{C})$ the category of presheaves of
$\mathbb{U}$-sets on $\mathcal{C}$, $Pr(\mathcal{C}):=\mathcal{C}^{Set_{\mathbb{U}}^{op}}$. The
subcategory of sheaves on $(\mathcal{C},\tau)$ will be denoted by
$Sh_{\tau}(\mathcal{C})$, or simply by $Sh(\mathcal{C})$ if the topology $\tau$ is
unambiguous.

All complexes will be cochain complexes (i.e. with differential increasing the degree by one) and therefore will look like
$$\xymatrix{
\dots \ar[r] & E^{n} \ar[r]^-{d_{n}} & E^{n+1} \ar[r] & \dots \ar[r] & E^{0} \ar[r] &
E^{1} \ar[r] & \dots.}$$

The following notations concerning various homotopy categories of stacks are defined in the main text, 
and recalled here for readers' convenience (see also the Index at the end of the book). \\

$$\mathrm{St}(\mathcal{C},\tau):=\mathrm{Ho}(Aff_{\mathcal{C}}^{\sim,\tau})$$

$$\mathrm{St}(k):=\mathrm{Ho}(k-Aff^{\sim,\textrm{\'et}})$$

$$D^{-}\mathrm{St}(k):=\mathrm{Ho}(k-D^{-}Aff^{\sim, \textrm{\'et}})$$

$$D\mathrm{St}(k):=\mathrm{Ho}(k-DAff^{\sim, \textrm{s-\'et}})$$

$$\mathrm{St}(S):=\mathrm{Ho}(SAff^{\sim, \textrm{s-\'et}})$$

\part{General theory of geometric stacks}

\chapter*{Introduction to Part 1}
In this first part we will study the general theory
of stacks and geometric stacks over a base symmetric
monoidal model category $\mathcal{C}$. For this, we will
start in \S \ref{partI.1} by introducing the notion of
a \emph{homotopical algebraic context} (HA \emph{context} for short),
which consists of a triple $(\mathcal{C}, \mathcal{C}_{0}, \mathcal{A})$ where
$\mathcal{C}$ is our base monoidal model category, $\mathcal{C}_{0}$ is a sub-category of $\mathcal{C}$, and $\mathcal{A}$ is a sub-category of the category $Comm(\mathcal{C})$
of commutative monoids in $\mathcal{C}$; we also require that the triple $(\mathcal{C}, \mathcal{C}_{0}, \mathcal{A})$ satisfies certain compatibility conditions. Although this might look like a rather unnatural
and complicated definition, it will be shown in 
\S \ref{partI.2} that this data precisely allows us to
define abstract versions of standard notions such as
derivations, unramified, \'etale, smooth and flat morphisms.
In other words a HA context describes an abstract context in which
the basic notions of linear and commutative algebra can
be developed. 

The first two chapters are only concerned with purely
algebraic notions and the geometry only starts in the third one, \S 
\ref{partI.3}. We start by some reminders on the notions
of model topology and of model topos (developed in \cite{hagI}), which are homotopical versions of the notions
of Grothendieck topology and of Grothendieck topos and which will be used all along this work.
Next, we introduce a notion of a \emph{homotopical
algebraic geometry context} (HAG \emph{context} for short),
consisting of a HA context together with
two additional data, $\tau$ and \textbf{P}, satsfying some compatiblity conditions.
The first datum
$\tau$ is a model topology on $Aff_{\mathcal{C}}$, the
opposite model category of commutative monoids in $\mathcal{C}$.
The second datum \textbf{P} consists of a class of morphisms in
$Aff_{\mathcal{C}}$ which behaves well with respect to $\tau$.
The model topology $\tau$ gives a category of stacks over $Aff_{\mathcal{C}}$
 (a homotopical generalization of the category of sheaves on affine schemes) 
in which everything is going to be embedded by means of a Yoneda lemma.
The class of morphisms \textbf{P} will then be used in order to define
\textit{geometric stacks} and more generally $n$\textit{-geometric stacks},
by considering successive quotient stacks of objects of $Aff_{\mathcal{C}}$
by action of groupoids whose structural morphisms are in \textbf{P}.
The compatibility axioms between $\tau$ and \textbf{P}
will insure that this notion of geometricity behaves well,
and satisfies the basic expected properties (stability by
homotopy pullbacks, gluing and certain quotients).

In \S \ref{partI.4}, the last chapter of part I, we will go
more deeply into the study of geometric stacks by
introducing infinitesimal constructions such as derivations,
cotangent complexes and obstruction theories. The main result of this last
chapter states that any geometric stack has an obstruction
theory (including a cotangent complex) as soon as
the HAG context satisfies suitable additional conditions.\\

\chapter{Homotopical algebraic context}\label{partI.1}

The purpose of this chapter is to
fix once for all our base model category as well
as several general assumptions it should satisfy. \\

All along this chapter, we refer
to \cite{ho} for the general
definition of monoidal model categories, and to
\cite{schw-shi} for general results about monoids and modules
in monoidal model categories. \\

From now on, and all along this work, we fix
three universes $\mathbb{U}\in \mathbb{V}\in \mathbb{W}$ (see, e.g. \cite[Exp.I, Appendice]{sga4}).
We also
let $(\mathcal{C},\otimes,\mathbf{1})$ be a symmetric monoidal model category in the sense
of \cite[\S 4]{ho}.
We assume that $\mathcal{C}$ is a $\mathbb{V}$-small
category, and that it is $\mathbb{U}$-combinatorial in the
sense of \cite[Appendix]{hagI}. \\

We make a first assumption on the base model category
$\mathcal{C}$, making it closer to an additive category. Recall
that we denote by $Q$ a cofibrant replacement functor and
by $R$ a fibrant replacement functor in $M$.

\begin{ass}\label{ass-1}
\begin{enumerate}
\item
The model category $\mathcal{C}$ is proper, pointed (i.e. the final object
is also an initial object) and
for any two object $X$ and $Y$ in $\mathcal{C}$ the natural morphisms
$$QX\coprod QY \longrightarrow X\coprod Y
\longrightarrow RX\times RY$$
are all equivalences.
\item The homotopy category $\mathrm{Ho}(\mathcal{C})$
is an additive category.
\end{enumerate}
\end{ass}

Assumption \ref{ass-1} implies in particular that finite homotopy coproducts
are also finite homotopy products in $\mathcal{C}$. It is always
satisfied when $\mathcal{C}$ is furthermore a stable model category in the
sense of \cite[\S 7]{ho}. Note that \ref{ass-1} implies in
particular that for any two objects $x$ and $y$ in $\mathcal{C}$, the set
$[x,y]$ has a natural abelian group structure.  \\

As $(\mathcal{C},\otimes,\mathbf{1})$ is a symmetric monoidal category, which is closed and
has $\mathbb{U}$-small limits and colimits, all the standard
notions and constructions of linear algebra makes sense
in $\mathcal{C}$ (e.g. monoids, modules over monoids, operads, algebra
over an operad \dots). The category of all associative, commutative and unital monoids
in $\mathcal{C}$ will be denoted by $Comm(\mathcal{C})$\index{$Comm(\mathcal{C})$!commutative monoids in $\mathcal{C}$}. Objects of $Comm(\mathcal{C})$ will simply be called
\emph{commutative monoids in $\mathcal{C}$}, or \emph{commutative monoids} if $\mathcal{C}$
is clear. In the same way, one defines
$Comm_{nu}(\mathcal{C})$\index{$Comm_{nu}(\mathcal{C})$!non-unital commutative monoids in $\mathcal{C}$} to be the category of
non-unital commutative monoids in $\mathcal{C}$. Therefore, our convention
will be that monoids are unital unless the contrary is specified.

The categories $Comm(\mathcal{C})$
and $Comm_{nu}(\mathcal{C})$ are again
$\mathbb{U}$-categories which are $\mathbb{V}$-small categories, and
possess all $\mathbb{U}$-small limits and colimits. They
come equipped with natural
forgetful functors
$$Comm(\mathcal{C}) \longrightarrow \mathcal{C} \qquad Comm_{nu}(\mathcal{C}) \longrightarrow \mathcal{C},$$
possessing left adjoints
$$F : \mathcal{C} \longrightarrow Comm(\mathcal{C}) \qquad F_{nu} : \mathcal{C} \longrightarrow Comm_{nu}(\mathcal{C})$$
sending an object of $\mathcal{C}$ to the free commutative monoid it generates.
We recall that for $X \in \mathcal{C}$ one has
$$F(X)=\coprod_{n\in \mathbb{N}}X^{\otimes n}/\Sigma_{n}$$
$$F_{nu}(X)=\coprod_{n\in \mathbb{N}-\{0\}}X^{\otimes n}/\Sigma_{n},$$
where $X^{\otimes n}$ is the $n$-tensor power of $X$, $\Sigma_{n}$
acts on it by permuting the factors and $X^{\otimes n}/\Sigma_{n}$ denotes
the quotient of this action in $\mathcal{C}$.  \\

Let $A \in Comm(\mathcal{C})$ be a commutative monoid. We will denote by
$A-Mod$\index{$A-Mod$!unital left $A$-modules} the category of unital left $A$-modules in $\mathcal{C}$. The category $A-Mod$ is again
a $\mathbb{U}$-category which is a $\mathbb{V}$-small category, and
has all $\mathbb{U}$-small limits and colimits. The objects in
$A-Mod$ will simply be called \emph{$A$-modules}.
It comes equiped with a natural
forgetful functor
$$A-Mod \longrightarrow \mathcal{C},$$
possessing a left adjoint
$$A\otimes - : \mathcal{C} \longrightarrow A-Mod$$
sending an object of $\mathcal{C}$ to the free $A$-module it generates.
We also recall that the category $A-Mod$ has a natural
symmetric monoidal structure $-\otimes_{A}-$. For two
$A$-modules $X$ and $Y$, the object $X\otimes_{A}Y$
is defined as the coequalizer in $\mathcal{C}$ of the two natural morphisms
$$X\otimes A \otimes Y \longrightarrow X\otimes Y \qquad
X\otimes A \otimes Y \longrightarrow X\otimes Y.$$
This symmetric monoidal structure is furthermore closed, and for
two $A$-modules $X$ and $Y$ we will denoted by $\underline{Hom}_{A}(X,Y)$\index{$\underline{Hom}_{A}(X,Y)$}
the $A$-module of morphisms. One has the usual adjunction isomorphisms
$$Hom(X\otimes_{A}Y,Z)\simeq Hom(X,\underline{Hom}_{A}(Y,Z)),$$
as well as isomorphisms of $A$-modules
$$\underline{Hom}_{A}(X\otimes_{A}Y,Z)\simeq
\underline{Hom}(X,\underline{Hom}_{A}(Y,Z)).$$

We define a morphism
in $A-Mod$ to be a \emph{fibration} or an \emph{equivalence} if it is so
on the underlying objects in $\mathcal{C}$.

\begin{ass}\label{ass1}
Let $A\in Comm(\mathcal{C})$ be any commutative monoid in $\mathcal{C}$.
Then, the above notions of equivalences and fibrations
makes $A-Mod$ into a $\mathbb{U}$-combinatorial proper model
category. The monoidal structure $-\otimes_{A}-$ makes furthermore
$A-Mod$ into a symmetric monoidal model category in the sense
of \cite[\S 4]{ho}.
\end{ass}

Using the assumption \ref{ass1} one sees that the homotopy category
$\mathrm{Ho}(A-Mod)$ has a natural symmetric monoidal structure
$\otimes^{\mathbb{L}}_{A}$, and derived internal $Hom's$
associated to it $\mathbb{R}\underline{Hom}_{A}$, satisfying the usual
adjunction rule
$$[X\otimes^{\mathbb{L}}_{A}Y,Z]\simeq
[X,\mathbb{R}\underline{Hom}_{A}(Y,Z)].$$

\begin{ass}\label{ass0}
Let $A$ be a commutative monoid in $\mathcal{C}$.
For any cofibrant object $M\in A-Mod$, the
functor
$$-\otimes_{A}M : A-Mod \longrightarrow A-Mod$$
preserves equivalences.
\end{ass}

Let us still denote by $A$ a commutative monoid in $\mathcal{C}$.
We have categories of commutative monoids in
$A-Mod$, and non-unital commutative monoids in
$A-Mod$, denoted respectively by $A-Comm(\mathcal{C})$\index{$A-Comm(\mathcal{C})$!commutative $A$-algebras} and
$A-Comm_{nu}(\mathcal{C})$\index{$A-Comm_{nu}(\mathcal{C})$!non-unital commutative $A$-algebras}, and
whose objects
will be called \emph{commutative $A$-algebras}
and \emph{non-unital commutative $A$-algebras}.
They
come equipped with natural
forgetful functors
$$A-Comm(\mathcal{C}) \longrightarrow A-Mod \qquad A-Comm_{nu}(\mathcal{C}) \longrightarrow A-Mod,$$
possessing left adjoints
$$F_{A} : A-Mod \longrightarrow A-Comm(\mathcal{C}) \qquad F_{A}^{nu} :
A-Mod \longrightarrow A-Comm_{nu}(\mathcal{C})$$
sending an object of $A-Mod$ to the free commutative monoid it generates.
We recall that for $X \in A-Mod$ one has
$$F_{A}(X)=\coprod_{n\in \mathbb{N}}X^{\otimes_{A} n}/\Sigma_{n}$$
$$F_{A}^{nu}(X)=\coprod_{n\in \mathbb{N}-\{0\}}X^{\otimes_{A} n}/\Sigma_{n},$$
where $X^{\otimes_{A} n}$ is the $n$-tensor power of $X$
in $A-Mod$, $\Sigma_{n}$
acts on it by permuting the factors and $X^{\otimes_{A} n}/\Sigma_{n}$ denotes
the quotient of this action in $A-Mod$.

Finally, we define a morphism in $A-Comm(\mathcal{C})$ or in
$A-Comm_{nu}(\mathcal{C})$ to be a \emph{fibration}
(resp. an \emph{equivalence}) if it is so as a morphism in
the category $\mathcal{C}$ (or equivalently as a morphism
in $A-Mod$).

\begin{ass}\label{ass2}
Let $A$ be any commutative monoid in $\mathcal{C}$.
\begin{enumerate}
\item The above classes 
of equivalences and fibrations make
the categories $A-Comm(\mathcal{C})$ and $A-Comm_{nu}(\mathcal{C})$ into
$\mathbb{U}$-combinatorial proper model
categories.
\item If $B$ is a cofibrant object in  $A-Comm(\mathcal{C})$,  then
the functor
$$B\otimes_{A} - : A-Mod \longrightarrow B-Mod$$
preserves equivalences.
\end{enumerate}
\end{ass}

\begin{rmk}\label{remnonunit}
\emph{One word concerning non-unital algebras. We will not really use
this notion in the sequel, except at one point in order to prove the existence
of a cotangent complex (so the reader is essentially allowed to forget about this 
unfrequently used notion). 
In fact, by our assumptions, the model category of
non-unital commutative $A$-algebras is Quillen equivalent to the model category
of augmented commutative $A$-algebra. However, the categories themselves
are not equivalent, since the category $\mathcal{C}$ is not assumed to be
strictly speaking additive, but only additive up to homotopy (e.g. 
it could be the model category of symmetric spectra of \cite{hss}). Therefore, we do not
think that the existence of the model structure on $A-Comm(\mathcal{C})$
implies the existence of the model structure on $A-Comm_{nu}(\mathcal{C})$; 
this explains why we had to add condition (1) on $A-Comm_{nu}(\mathcal{C})$ in Assumption \ref{ass2}.
Furthermore, the passage 
from augmented $A$-algebras to non-unital $A$-algebras will be in any case  
necessary to construct a certain Quillen adjunction during the proof
of Prop. \ref{p1}, because such a Quillen adjunction does not exist from
the model category of augmented $A$-algebras (as it is a composition 
of a left Quillen functor by a right Quillen equivalence).}
\end{rmk}

An important consequence of assumption \ref{ass2} $(2)$ is that
for $A$ a commutative monoid in $\mathcal{C}$, and
$B$, $B'$ two commutative $A$-algebras, the natural morphism
in $\mathrm{Ho}(A-Mod)$
$$B\coprod^{\mathbb{L}}_{A}B' \longrightarrow
B\otimes^{\mathbb{L}}_{A}B'$$
is an isomorphism (here the object on the left is
the homotopy coproduct in $A-Comm(\mathcal{C})$, and the one
on the right is
the derived tensor product in $A-Mod$).

An important remark we will use implicitly very often in this paper
is that the category $A-Comm(\mathcal{C})$ is naturally equivalent to
the comma category $A/Comm(\mathcal{C})$, of objects under $A$.
Moreover, the model structure on $A-Comm(\mathcal{C})$
coincides through this equivalence with the comma model
category $A/Comm(\mathcal{C})$. \\

We will also fix a full subcategory
$\mathcal{C}_{0}$ of $\mathcal{C}$, playing essentially the role of a kind of ``$t$-structure''
on $\mathcal{C}$ (i.e. essentially defining which are the
``non-positively graded objects'', keeping in mind that in this work
we use the cohomological grading when concerned with complexes). More precisely, we will fix a
subcategory $\mathcal{C}_{0}\subseteq \mathcal{C}$ satisfying the following conditions.

\begin{ass}\label{ass7}
\begin{enumerate}
\item $\mathbf{1} \in \mathcal{C}_{0}$.
\item The full subcategory
$\mathcal{C}_{0}$ of $\mathcal{C}$
is stable by equivalences and
by $\mathbb{U}$-small homotopy colimits.
\item
The full subcategory
$\mathrm{Ho}(\mathcal{C}_{0})$ of $\mathrm{Ho}(\mathcal{C})$
is stable by the monoidal structure $-\otimes^{\mathbb{L}}-$
(i.e. for $X$ and $Y$ in $\mathrm{Ho}(\mathcal{C}_{0})$ we have
$X\otimes^{\mathbb{L}}Y\in \mathrm{Ho}(\mathcal{C}_{0})$).
\end{enumerate}
\end{ass}

Recall that as $\mathcal{C}$ is a pointed model category one can
define its \textit{suspension functor}\index{suspension functor}
$$\begin{array}{cccc}
S : & \mathrm{Ho}(\mathcal{C}) & \longrightarrow & \mathrm{Ho}(\mathcal{C}) \\
& x & \mapsto & 
*\coprod_{x}^{\mathbb{L}}* \end{array}$$
left adjoint to the \textit{loop functor}\index{loop functor}
$$\begin{array}{cccc}
\Omega  : & \mathrm{Ho}(\mathcal{C}) & \longrightarrow & \mathrm{Ho}(\mathcal{C}) \\
& x & \mapsto & 
:=*\times^{h}_{x}*.
\end{array}$$
We  set $\mathcal{C}_{1}$ to be the full subcategory
of $\mathcal{C}$ consisting of all objects
equivalent to the suspension of some object in $\mathcal{C}_{0}$.
The full subcategory $\mathcal{C}_{1}$ of $\mathcal{C}$
is also closed by equivalences, homotopy colimits
and the derived tensor structure.
We will denote by
$Comm(\mathcal{C})_{0}$ the full subcategory of
$Comm(\mathcal{C})$ consisting of commutative monoids
whose underlying $\mathcal{C}$-object lies in $\mathcal{C}_{0}$.
In the same way, for $A\in Comm(\mathcal{C})$ we denote by
$A-Mod_{0}$\index{$A-Mod_{0}$} (resp. $A-Mod_{1}$\index{$A-Mod_{1}$}, resp. 
$A-Comm(\mathcal{C})_{0}$\index{$A-Comm(\mathcal{C})_{0}$}) the full subcategory of
$A-Mod$ consisting of $A$-modules whose underlying $\mathcal{C}$-object lies
in $\mathcal{C}_{0}$ (resp. of
$A-Mod$ consisting of $A$-modules whose underlying $\mathcal{C}$-object lies
in $\mathcal{C}_{1}$, resp.
of $A-Comm(\mathcal{C})$ consisting of commutative $A$-algebras
whose underlying $\mathcal{C}$-object lies in $\mathcal{C}_{0}$).

An important consequence of Assumption \ref{ass7}
is that for any morphism $A \longrightarrow B$
in $Comm(\mathcal{C})_{0}$, and any
$M\in A-Mod_{0}$, we have
$B\otimes_{A}^{\mathbb{L}}M \in \mathrm{Ho}(B-Mod_{0})$. Indeed, 
any such $A$-module can be written as a homotopy colimit 
of $A$-modules of the form $A^{\otimes^{\mathbb{L}}n}\otimes^{\mathbb{L}}M$, 
for which we have
$$B\otimes_{A}^{\mathbb{L}}A^{\otimes^{\mathbb{L}}n}\otimes^{\mathbb{L}}M\simeq
A^{\otimes^{\mathbb{L}}(n-1)}\otimes^{\mathbb{L}}M.$$
In particular the full subcategory
$Comm(\mathcal{C})_{0}$ is closed by homotopy push-outs
in $Comm(\mathcal{C})$. Passing to the suspension we
also see that for any $M\in A-Mod_{1}$, one also has
$B\otimes_{A}^{\mathbb{L}}M \in \mathrm{Ho}(B-Mod_{1})$.

\begin{rmk}\label{r1}
\begin{enumerate}
\item
\emph{The reason for introducing the subcategory
$\mathcal{C}_{0}$ is to be able to consider
reasonable infinitesimal lifting properties; these infinitesimal lifting properties will be used
to develop the abstract obstruction
theory of geometric stacks in \S \ref{partI.4}. It is useful to keep in mind that
$\mathcal{C}_{0}$  plays a role analogous to a kind of $t$-structure
on $\mathcal{C}$, in that it morally defines which are the
non-positively graded objects (a typical example will appear
in \S \ref{IIunb} where $\mathcal{C}$ will be the model category of unbounded complexes
and $\mathcal{C}_{0}$ the subcategory of complexes with vanishing positive cohomology). 
Different  choices of $\mathcal{C}_{0}$ will then give different
notions of formal smoothness (see \S \ref{Iinf}), and thus possibly
different notions of geometric stacks. We think that
playing with $\mathcal{C}_{0}$ as a degree of freedom
is an interesting feature of our abstract infinitesimal theory.}

\item
\emph{Assumptions \ref{ass-1}, \ref{ass1} and \ref{ass0} are not really
serious, and are only useful to avoid taking
too many fibrant and cofibrant replacements. With some
care, they can be omitted. On the other hand,
the careful reader will probably be surprised by
assumption \ref{ass2},
as it is known not to be satisfied
in several interesting examples (e.g. when $\mathcal{C}$ is the model category of
complexes over some commutative ring $k$, not of characteristic zero). Also,
it is well known that in some situations the notion of
commutative monoid is too strict and it is often preferable to use
the weaker notion of $E_{\infty}$-monoid.}
\emph{Two reasons has led us to assume  \ref{ass2}.
First of all, for all contexts of application of the general theory we will present
in this work, there is always a base model category $\mathcal{C}$
for which this condition is satisfied and gives rise to the correct
theory. Moreover, if one replaces commutative
monoids by $E_{\infty}$-monoids then assumption
\ref{ass2} is almost always satisfied, as shown in
\cite{sp}, and we think that translating our general constructions should then
be a rather academic exercise. Working with commutative monoids instead of
$E_{\infty}$-monoids simplifies a lot the notations and certain constructions,
and in our opinion this theory already captures the real essence of the subject.}

\emph{Finally, in partial defense of our choice, let us also mention that contrary to what one could think at first sight, 
working with $E_{\infty}$-monoids would not strictly speaking increase the degree
of generality of the theory. Indeed, one of our major application is to the category 
of simplicial $k$-modules, whose category of commutative monoids is the category 
of simplicial commutative $k$-algebras. However, if $k$ has non-zero characteristic,  
the homotopy theory of simplicial commutative $k$-algebras is \textit{not} 
equivalent to the homotopy theory of $E_{\infty}$-monoids in simplicial
$k$-modules (the latter is equivalent to the homotopy theory 
of $E_{\infty}$-$k$-algebras in non positive degrees). Therefore, using $E_{\infty}$-monoids throughout 
would prevent us from 
developing derived algebraic geometry as presented in \S 2.2.}
\end{enumerate}
\end{rmk}

We list below some important examples of symmetric monoidal model categories
$\mathcal{C}$ satisfying the four above assumptions, and of crucial importance for our
applications.

\begin{enumerate}
\item Let $k$ be any commutative ring, and $\mathcal{C}=k-Mod$ be the category
of $\mathbb{U}$-$k$-modules, symmetric
monoidal for the tensor product $\otimes_{k}$, and
endowed with its trivial
model structure (i.e. equivalences are isomorphisms
and all morphisms are fibrations and cofibrations). Then assumptions
\ref{ass-1}, \ref{ass1}, \ref{ass0} and \ref{ass2}
are satisfied. For $\mathcal{C}_{0}$ one can take
the whole $\mathcal{C}$ for which Assumption \ref{ass7} is
clearly satisfied.

\item Let $k$ be a commutative ring of
characteristic $0$, and $\mathcal{C}=C(k)$ be the
category of (unbounded) complexes of $\mathbb{U}$-$k$-modules,
symmetric monoidal for the tensor product of complexes
$\otimes_{k}$, and endowed
with its projective model structure for which equivalences are quasi-isomorphisms and
fibrations are epimorphisms (see \cite{ho}). Then, the category
$Comm(\mathcal{C})$ is the category of  commutative differential graded
$k$-algebras belonging to $\mathbb{U}$. For $A\in Comm(\mathcal{C})$, the category $A-Mod$ is then the category
of differential graded $A$-modules. It is well known that as $k$ is of characteristic
zero then assumptions
\ref{ass-1}, \ref{ass1}, \ref{ass0} and \ref{ass2}
are satisfied (see e.g; \cite{hin}).
For $\mathcal{C}_{0}$ one can take either
the whole $\mathcal{C}$, or the full subcategory
of $\mathcal{C}$ consisting of complexes $E$ such that
$H^{i}(E)=0$ for any $i>0$, for which
Assumption \ref{ass7} is satisfied.

A similar example
is given by non-positively graded complexes $\mathcal{C}=C^{-}(k)$.

\item Let $k$ be any commutative ring, and $\mathcal{C}=sMod_k$ be the category of
$\mathbb{U}$-simplicial $k$-modules, endowed with the levelwise tensor
product and the usual model structure for which equivalences and
fibrations are defined on the underlying simplicial sets (see
e.g. \cite{gj}). The category $Comm(\mathcal{C})$ is then the category of
simplicial commutative $k$-algebras, and for $A \in sMod_k$,
$A-Mod$ is the category of simplicial modules over the simplicial
ring $A$. Assumptions \ref{ass-1}, \ref{ass1},
\ref{ass0} and \ref{ass2} are again well
known to be satisfied. For $\mathcal{C}_{0}$ one can take
the whole $\mathcal{C}$.

Dually, one could also
let $\mathcal{C}=csMod_k$ be the category of co-simplicial $k$-modules, and 
our assumptions would again be satisfied. In this case,
$\mathcal{C}_{0}$ could be for example the full subcategory
of co-simplicial modules $E$ such that $\pi_{-i}(E)=H^{i}(E)=0$ for any $i>0$. This
subcategory is stable under homotopy colimits as the functor
$E \mapsto H^{0}(E)$ is right Quillen and right adjoint to the inclusion 
functor $k-Mod \longrightarrow csMod_{k}$. 

\item Let $Sp^{\Sigma}$ be the category of $\mathbb{U}$-symmetric spectra
and its smash product, endowed with
the positive stable model structure of \cite{shi}. Then, $Comm(\mathcal{C})$ is the category of
commutative symmetric ring spectra, and assumptions
\ref{ass-1}, \ref{ass1}, \ref{ass0} and \ref{ass2}
are known to be satisfied
(see \cite[Thm. 3.1, Thm. 3.2, Cor. 4.3]{shi}).
The two canonical choices for $\mathcal{C}_{0}$ are
the whole $\mathcal{C}$ or the full subcategory of
connective spectra.

It is important to note that
one can also take for $\mathcal{C}$ the category of $Hk$-modules in $Sp^{\Sigma}$ for some
commutative ring $k$. This will give a model for the
homotopy theory of $E_{\infty}$-$k$-algebras that were not
provided by our example $2$.

\item Finally, the above three examples can be sheafified over some
Grothendieck site, giving the corresponding relative theories over a base Grothendieck topos.

In few words, let $\mathcal{S}$ be a $\mathbb{U}$-small Grothendieck site, and $\mathcal{C}$ be one the three
symmetric monoidal model category $C(k)$, $sMod_k$, $Sp^{\Sigma}$ discussed
above. One considers the corresponding categories of presheaves on $\mathcal{S}$,
$Pr(\mathcal{S},C(k))$, $Pr(\mathcal{S},sMod_k)$,
$Pr(\mathcal{S},Sp^{\Sigma})$. They can be endowed with the
projective model structures for which fibrations and equivalences are
defined levelwise. This first model structure does
not depend on the topology on $\mathcal{S}$, and will be called
the \emph{strong model structure}: its (co)fibrations and equivalences will be called
\emph{global (co)fibrations} and \emph{global equivalences}.

The next step is to introduce notions of \emph{local equivalences} in the
model categories $Pr(\mathcal{S},C(k))$, $Pr(\mathcal{S},sMod_k)$,
$Pr(\mathcal{S},Sp^{\Sigma})$. This notion
is defined by first defining reasonable \emph{homotopy sheaves}, as done
for the notion of local equivalences in the theory of simplicial presheaves, and
then define a morphism to be a local equivalence if it induces
isomorphisms on all homotopy sheaves (for various choice of base points,
see \cite{jo,ja} for more details). The final model structures on
$Pr(\mathcal{S},C(k))$, $Pr(\mathcal{S},sMod_k)$, $Pr(\mathcal{S},Sp^{\Sigma})$ are the one for which
equivalences are the local equivalences, and
cofibrations are the global cofibrations. The proof that this indeed defines
model categories is not given here and is very similar to the proof of the
existence of the local projective model structure on
the category of simplicial presheaves (see for example
\cite{hagI}).

Finally, the symmetric monoidal structures on the categories
$C(k)$, $sMod_k$ and $Sp^{\Sigma}$ induces natural symmetric monoidal structures
on the categories $Pr(\mathcal{S},C(k))$, $Pr(\mathcal{S},sMod_k)$,
$Pr(\mathcal{S},Sp^{\Sigma})$. These symmetric monoidal structures
make them into
symmetric monoidal model categories when $\mathcal{S}$ has finite products. One can also check that
the symmetric monoidal model categories
$Pr(\mathcal{S},C(k))$, $Pr(\mathcal{S},sMod_k)$, $Pr(\mathcal{S},Sp^{\Sigma})$
constructed that way all satisfy the assumptions
\ref{ass-1}, \ref{ass1}, \ref{ass0} and \ref{ass2}.

An important example is the following. Let $\mathcal{O}$ be a sheaf of commutative rings
on the site $\mathcal{S}$, and let $H\mathcal{O} \in Pr(\mathcal{S},Sp^{\Sigma})$ be the presheaf of
symmetric spectra it defines. The object $H\mathcal{O}$
is a commutative monoid in $Pr(\mathcal{S},Sp^{\Sigma})$ and
one can therefore consider the
model category $H\mathcal{O}-Mod$, of $H\mathcal{O}$-modules.
The category $H\mathcal{O}-Mod$ is a symmetric monoidal model category and its
homotopy category is equivalent to the unbounded derived category $D(\mathcal{S},\mathcal{O})$ of
$\mathcal{O}$-modules on the site $\mathcal{S}$. This gives a way to define all the standard constructions
as derived tensor products, derived internal $Hom's$ etc., in the context of
unbounded complexes of $\mathcal{O}$-modules.

\end{enumerate}

Let $f : A \longrightarrow B$ be a morphism of commutative monoids in $\mathcal{C}$. We deduce
an adjunction between the categories of modules
$$f^{*} : A-Mod \longrightarrow B-Mod \qquad
A-Mod \longleftarrow B-Mod : f_{*},$$
where $f^{*}(M):=B\otimes_{A}M$, and $f_{*}$ is the forgetful functor that sees 
a $B$-module as an $A$-module through the morphism $f$.
Assumption \ref{ass1} tells us that this adjunction
is a Quillen adjunction, and  assumption \ref{ass0} implies it
is furthermore a Quillen equivalence when
$f$ is an equivalence (this is one of the main reasons for assumption \ref{ass0}).

The morphism $f$ induces a pair of adjoint derived functors
$$\mathbb{L}f^{*} : \mathrm{Ho}(A-Mod) \longrightarrow \mathrm{Ho}(B-Mod) \qquad
\mathrm{Ho}(A-Mod) \longleftarrow \mathrm{Ho}(B-Mod) : \mathbb{R}f_{*}\simeq f_*,$$
and, as usual, we will also use the notation
$$\mathbb{L}f^{*}(M)=:B\otimes_{A}^{\mathbb{L}}M \in \mathrm{Ho}(B-Mod).$$

Finally, let
$$\xymatrix{A \ar[r]^{f} \ar[d]_{p} & B \ar[d]^{p'} \\
A' \ar[r]_{f'} & B'}$$
be a homotopy cofiber square
in $Comm(\mathcal{C})$. Then, for any $A'$-module $M$ we have the
well known base change morphism
$$\mathbb{L}f^{*}p_{*}(M)\longrightarrow (p')_*\mathbb{L}(f')^{*}(M).$$

\begin{prop}\label{ptrans}
Let us keep the notations as above. Then, the morphism
$$\mathbb{L}f^{*}p_{*}(M)\longrightarrow (p')_*\mathbb{L}(f')^{*}(M)$$
is an isomorphism in $\mathrm{Ho}(B-Mod)$ for any $A'$-module $M$.
\end{prop}

\begin{proof} As the homotopy categories
of modules are invariant under equivalences of
commutative monoids, one can suppose that
the diagram
$$\xymatrix{A \ar[r]^{f} \ar[d]_{p} & B \ar[d]^{p'} \\
A' \ar[r]_{f'} & B'}$$
is cocartesian in $Comm(\mathcal{C})$, and
consists of cofibrations in $Comm(\mathcal{C})$. Then, using
\ref{ass0} and \ref{ass2} $(2)$ one sees that
the natural morphisms
$$M\otimes^{\mathbb{L}}_{A}B \longrightarrow M\otimes_{A}B
\qquad M\otimes^{\mathbb{L}}_{A'}B' \longrightarrow M\otimes_{A'}B',$$
are isomorphism in $\mathrm{Ho}(B-Mod)$. Therefore, the proposition
follows from the natural isomorphism of $B$-modules
$$M\otimes_{A'}B'\simeq M\otimes_{A'}(A'\otimes_{A}B)\simeq
M\otimes_{A}B.$$
\end{proof}

\begin{rmk}\label{remtrans}
\emph{The above base change formula will be extremely important 
in the sequel, and most often used implicitly. It should be noticed that 
it implies that the homotopy coproduct in the model category of commutative
monoids is given by the derived tensor product. This last property is only
satisfied because we have used commutative monoids, and is clearly wrong 
for simply associative monoids. This is one major reason why our setting cannot be
used, at least without some modifications, to develop truly
non-commutative geometries. Even partially commutative
structures, like $E_{n}$-monoids for $n>1$ would not satisfy the base
change formula, and one really needs $E_{\infty}$-monoids at least. 
This fact also prevents us to generalize our setting by replacing 
the category of commutative monoids by more general categories, like some
category of algebras over more general operads. }
\end{rmk}

We now consider $A\in Comm(\mathcal{C})$ and the natural inclusion
$A-Mod_{0} \longrightarrow A-Mod$. We consider
the restricted Yoneda embedding
$$\mathbb{R}\underline{h}_{0}^{-} : \mathrm{Ho}(A-Mod^{op}) \longrightarrow \mathrm{Ho}((A-Mod_{0}^{op})^{\wedge}),$$
sending an $A$-module $M$ to the functor
$$Map(M,-) : A-Mod_{0}^{op} \longrightarrow SSet_{\mathbb{V}}.$$
We recall here from \cite[\S 4.1]{hagI} that 
for a model category 
$M$, and a full sub-category stable by equivalences $M_{0}\subset M$, 
$M_{0}^{\wedge}$ is the left Bousfield localization 
of $SPr(M_{0})$ along equivalences in $M_{0}$. 

\begin{df}\label{dgood}
We will say that $A\in Comm(\mathcal{C})$
is \emph{good with respect to $\mathcal{C}_{0}$}
(or simply \emph{$\mathcal{C}_{0}$-good})\index{good!$\mathcal{C}_{0}$-good} if
the functor
$$\mathrm{Ho}(A-Mod^{op}) \longrightarrow \mathrm{Ho}((A-Mod_{0}^{op})^{\wedge})$$
is fully faithful.
\end{df}

In usual category theory, a full subcategory $D \subset C$ is called
\emph{dense} if the restricted Yoneda functor
$C \longrightarrow Pr(D):=\underline{Hom}(D^{op},Ens)$ is 
fully faithful (in \cite{sga4} this notion is equivalent to the fact that
\emph{$D$ generates $C$ through strict epimorphisms}). This implies for example that any object
of $C$ is the colimit of objects of $D$, but is a slightly 
stronger condition because any object $x\in C$ is in fact isomorphic to 
the colimit of the canonical diagram $D/x \longrightarrow C$ (see e.g. \cite[ExpI-Prop. 7.2]{sga4}). Our notion 
of being good (Def. \ref{dgood}) essentially means that 
$A-Mod_{0}^{op}$ is \emph{homotopically dense} in $A-Mod^{op}$. This of course
implies that any object in $A-Mod$ is equivalent to a homotopy limit of
objects in $A-Mod_{0}$, and is equivalent to the fact that 
any cofibrant object $M\in A-Mod$ is equivalent to the homotopy limit of the
natural diagram $(M/A-Mod_{0})^{c} \longrightarrow A-Mod$,  where
$(M/A-Mod_{0})^{c}$ denotes the category of cofibrations under $M$. Dually, one could 
say that $A$ being good with respect to $\mathcal{C}_{0}$ means that 
\emph{$A-Mod_{0}$ cogenerates $A-Mod$ through strict monomorphisms} in a homotopical sense. \\

\bigskip

We finish this first chapter
by the following definition, gathering
our assumptions \ref{ass-1},
\ref{ass1}, \ref{ass0}, \ref{ass2} and \ref{ass7}
all together.

\begin{df}\label{dha}
A \emph{Homotopical Algebraic context}\index{HA context!Homotopical Algebraic context}
(or simply \emph{HA context}) is a
triplet $(\mathcal{C},\mathcal{C}_{0},\mathcal{A})$,
consisting of a symmetric monoidal model category
$\mathcal{C}$, two full sub-categories stable by equivalences
$$\mathcal{C}_{0}\subset \mathcal{C}
\qquad \mathcal{A}\subset Comm(\mathcal{C}),$$
such that any $A\in \mathcal{A}$ is $\mathcal{C}_{0}$-good, and
assumptions \ref{ass-1},
\ref{ass1}, \ref{ass0}, \ref{ass2}, \ref{ass7} are satisfied.
\end{df}

\chapter{Preliminaries on linear and commutative algebra in an HA context}\label{partI.2}

All along this chapter we fix once for all
a HA context $(\mathcal{C},\mathcal{C}_{0},\mathcal{A})$, in the sense
of Def. \ref{dha}. The purpose of this chapter is to show that
the assumptions of the last chapter imply that many
general notions of linear and commutative algebra generalize
in some reasonable sense in our base category $\mathcal{C}$. \\

\section{Derivations and the cotangent complex}\label{Ider}

This section is nothing else than a rewriting
of the first pages of \cite{ba}, which stay
valid in our general context. \\

Let $A\in Comm(\mathcal{C})$ be a commutative monoid in $\mathcal{C}$, and
$M$ be an $A$-module. We define a new commutative monoid
$A\oplus M$ in the following way. The underlying object
of $A\oplus M$ is the coproduct $A\coprod M$. The multiplicative
structure is defined by the morphism
$$(A\coprod M)\otimes (A\coprod M) \simeq A\otimes A\coprod A\otimes M \coprod
A\otimes M \coprod M\otimes M \longrightarrow A\coprod M$$
given by the three morphisms
$$\mu \coprod * : A\otimes A \longrightarrow A\coprod M$$
$$*\coprod \rho : A\otimes M \longrightarrow A\coprod M$$
$$* : M\otimes M \longrightarrow M,$$
where $\mu : A \otimes A \longrightarrow A$ is the multiplicative
structure of $A$, and $\rho : A\otimes M \longrightarrow M$
is the module structure of $M$. The monoid $A\oplus M$ is commutative and unital,
and defines an object in $Comm(\mathcal{C})$. It comes furthermore with a natural morphism
of commutative monoids $id\coprod * : A\oplus M \longrightarrow A$, which
has a natural section $id\coprod * : A\longrightarrow A\oplus M$.

Now, if $A \longrightarrow B$ is a morphism in $Comm(\mathcal{C})$, and $M$ is a $B$-module,
the morphism $B\oplus M \longrightarrow B$ can be seen as a morphism
of commutative $A$-algebras. In other words,
$B\oplus M$ can be seen as an object of the double comma model category
$A-Comm(\mathcal{C})/B$.

\begin{df}\label{d1}
Let $A\longrightarrow B$ be a morphism of commutative monoids, and
$M$ be a $B$-module. The simplicial set of \emph{derived $A$-derivations from $B$ to $M$}\index{derived derivations},
is the object
$$\mathbb{D}er_{A}(B,M):=Map_{A-Comm(\mathcal{C})/B}(B,B\oplus M)\in \mathrm{Ho}(SSet_{\mathbb{U}}).$$
\end{df}

Clearly, $M\mapsto \mathbb{D}er_{A}(B,M)$ defines a functor
from the homotopy category of $B$-module $\mathrm{Ho}(B-Mod)$ to the
homotopy category of simplicial sets $\mathrm{Ho}(SSet)$. More precisely,
the functoriality of the construction of mapping spaces
implies that one can also construct a genuine functor
$$\mathbb{D}er_{A}(B,-) : B-Mod \longrightarrow SSet_{\mathbb{U}},$$
lifting the previous functor on the homotopy categories.
This last functor will be
considered as an object
in the model category of pre-stacks $(B-Mod^{op})^{\wedge}$
as defined in \cite[\S 4.1]{hagI}. Recall from
\cite[\S 4.2]{hagI} that there exists a Yoneda embedding
$$\mathrm{Ho}(B-Mod)^{op} \longrightarrow \mathrm{Ho}((B-Mod^{op})^{\wedge})$$
sending a $B$-module $M$ to the simplicial presheaf
$N \mapsto Map_{B-Mod}(M,N)$, and objects in the essential image will be called
\emph{co-representable}.

\begin{prop}\label{p1}
For any morphism $A \longrightarrow B$ in $Comm(\mathcal{C})$, there
exists a $B$-module $\mathbb{L}_{B/A}$, and
an element $d\in \pi_{0}(\mathbb{D}er_{A}(B,\mathbb{L}_{B/A}))$, such that
for any $B$-module $M$, the natural morphism obtained by composing with $d$
$$d^{*} : Map_{B-Mod}(\mathbb{L}_{B/A},M)\longrightarrow \mathbb{D}er_{A}(B,M)$$
is an isomorphism in $\mathrm{Ho}(SSet)$.
\end{prop}

\begin{proof} The proof is the same as
in \cite{ba}, and uses our assumptions \ref{ass-1},
\ref{ass1}, \ref{ass0} and \ref{ass2}. We will reproduce
it for the reader convenience. \\

We first consider
the Quillen adjunction
$$-\otimes_{A}B : A-Comm(\mathcal{C})/B \longrightarrow B-Comm(\mathcal{C})/B \qquad
A-Comm(\mathcal{C})/B \longleftarrow B-Comm(\mathcal{C})/B : F,$$
where $F$ is the forgetful functor. This induces an adjunction
on the level of homotopy categories
$$-\otimes^{\mathbb{L}}_{A}B : \mathrm{Ho}(A-Comm(\mathcal{C})/B) \longrightarrow
\mathrm{Ho}(B-Comm(\mathcal{C})/B)$$
$$\mathrm{Ho}(A-Comm(\mathcal{C})/B) \longleftarrow \mathrm{Ho}(B-Comm(\mathcal{C})/B) : F.$$
We consider a second Quillen adjunction
$$K : B-Comm_{nu}(\mathcal{C}) \longrightarrow B-Comm(\mathcal{C})/B
\qquad B-Comm_{nu}(\mathcal{C}) \longleftarrow B-Comm(\mathcal{C})/B : I,$$
where $B-Comm_{nu}(\mathcal{C})$ is the category of
non-unital commutative $B$-algebras (i.e. non-unital
commutative monoids in $B-Mod$). The functor
$I$ takes a diagram of commutative
monoids
$\xymatrix{B \ar[r]^{s} &  C \ar[r]^{p} & B}$ to
the kernel of $p$ computed in the category
of non-unital commutative $B$-algebras. In the other direction,
the functor $K$ takes a non-unital commutative
$B$-algebra $C$ to the trivial extension of $B$ by $C$
(defined as our $B\oplus M$ but taking into account
the multiplication on $C$).
Clearly, $I$ is a right Quillen functor, and the adjunction defines
an adjunction on the homotopy categories
$$\mathbb{L}K : \mathrm{Ho}(B-Comm_{nu}(\mathcal{C})) \longrightarrow \mathrm{Ho}(B-Comm(\mathcal{C})/B)$$
$$ \mathrm{Ho}(B-Comm_{nu}(\mathcal{C})) \longleftarrow \mathrm{Ho}(B-Comm(\mathcal{C})/B) : \mathbb{R}I.$$

\begin{lem}\label{l1}
The adjunction $(\mathbb{L}K,\mathbb{R}I)$ is an equivalence.
\end{lem}

\begin{proof} This follows easily from our assumption \ref{ass-1}.
Indeed, it implies that for any fibration in $\mathcal{C}$,
$f : X \longrightarrow Y$, which has a section $s : Y \longrightarrow X$,
the natural morphism
$$i\coprod s : F\coprod Y \longrightarrow X,$$
where $i : F \longrightarrow X$ is the fiber of $f$, is an equivalence.
It also implies that the homotopy fiber of the natural
morphism $id\coprod * : X\coprod Y \longrightarrow X$
is naturally equivalent to $Y$. These two facts imply the lemma. \end{proof}

Finally, we consider a third adjunction
$$Q : B-Comm_{nu}(\mathcal{C}) \longrightarrow B-Mod \qquad
B-Comm_{nu}(\mathcal{C}) \longleftarrow B-Mod : Z,$$
where $Q$ of an object $C \in B-Comm_{nu}(\mathcal{C})$
is the push-out of $B$-modules
$$\xymatrix{
C\otimes_{B} C \ar[r]^{\mu} \ar[d] & C \ar[d] \\
\bullet \ar[r] & Q(C),}$$
and $Z$ sends a $B$-module $M$ to the non-unital
$B$-algebra $M$ endowed with the zero multiplication.
Clearly, $(Q,Z)$ is a Quillen adjunction and gives rise
to an adjunction on the homotopy categories
$$\mathbb{L}Q : \mathrm{Ho}(B-Comm_{nu}(\mathcal{C})) \longrightarrow \mathrm{Ho}(B-Mod) $$
$$\mathrm{Ho}(B-Comm_{nu}(\mathcal{C})) \longleftarrow \mathrm{Ho}(B-Mod) : Z.$$
We can now conclude the proof of the proposition by
chaining up
the various adjunction to get a string of isomorphisms in
$\mathrm{Ho}(SSet)$
$$\mathbb{D}er_{A}(B,M)\simeq
Map_{A-Comm(\mathcal{C})/B}(B,F(B\oplus M))\simeq
Map_{B-Comm(\mathcal{C})/B}(B\otimes_{A}^{\mathbb{L}}B,B\oplus M)$$
$$\simeq Map_{B-Comm_{nu}(\mathcal{C})}(\mathbb{R}I(B\otimes_{A}^{\mathbb{L}}B),\mathbb{R}I(B\oplus M))\simeq
Map_{B-Comm_{nu}(\mathcal{C})}(\mathbb{R}I(B\otimes_{A}^{\mathbb{L}}B),Z(M))$$
$$\simeq Map_{B-Mod}(\mathbb{L}Q\mathbb{R}I(B\otimes_{A}^{\mathbb{L}}B),M).$$
Therefore, $\mathbb{L}_{B/A}:=\mathbb{L}Q\mathbb{R}I(B\otimes_{A}^{\mathbb{L}}B)$
and the image of $id \in Map_{B-Mod}(\mathbb{L}_{B/A},\mathbb{L}_{B/A})$
gives what we were looking for. \end{proof}

\begin{rmk}\label{r2}
\emph{Proposition \ref{p1} implies that the two functors
$$M \mapsto Map_{B-Mod}(\mathbb{L}_{B/A},M) \qquad
M \mapsto \mathbb{D}er_{A}(B,M)$$
are isomorphic as objects in
$\mathrm{Ho}((B-Mod^{op})^{\wedge})$. In other words,
Prop. \ref{p1} implies that the
functor $\mathbb{D}er_{A}(B,-)$ is} co-representable \emph{
in the sense of \cite{hagI}.}
\end{rmk}

\begin{df}\label{d2}
Let $A \longrightarrow B$ be a morphism in $Comm(\mathcal{C})$.
\begin{enumerate}
\item The
$B$-module $\mathbb{L}_{B/A} \in \mathrm{Ho}(B-Mod)$ is called the
\emph{cotangent
complex of $B$ over $A$}\index{cotangent complex!relative of commutative monoids}.
\item When $A=\mathbf{1}$, we will use the following notation
$$\mathbb{L}_{B}:=\mathbb{L}_{B/\mathbf{1}},$$
and $\mathbb{L}_{B}$ will be called the \emph{cotangent complex of $B$}.
\end{enumerate}
\end{df}

Using the definition and proposition \ref{p1}, 
it is  easy to check the following facts.

\begin{prop}\label{p2}
\begin{enumerate}
\item Let $A \longrightarrow B \longrightarrow C$ be two morphisms
in $Comm(\mathcal{C})$. Then, there is a homotopy cofiber sequence in $C-Mod$
$$\mathbb{L}_{B/A}\otimes^{\mathbb{L}}_{B}C \longrightarrow
\mathbb{L}_{C/A} \longrightarrow  \mathbb{L}_{C/B}.$$

\item
Let
$$\xymatrix{A\ar[r] \ar[d] & B \ar[d] \\
A' \ar[r] & B'}$$
be a homotopy cofiber square in $Comm(\mathcal{C})$. Then, the natural
morphism
$$\mathbb{L}_{B/A}\otimes^{\mathbb{L}}_{B}B'\longrightarrow \mathbb{L}_{B'/A'}$$
is an isomorphism in $\mathrm{Ho}(B'-Mod)$. Furthermore, the natural morphism
$$\mathbb{L}_{B/A}\otimes_{B}^{\mathbb{L}}B' \coprod
\mathbb{L}_{A'/A}\otimes_{A'}^{\mathbb{L}}B' \longrightarrow
\mathbb{L}_{B'/A}$$
is an isomorphism in $\mathrm{Ho}(B'-Mod)$.

\item Let
$$\xymatrix{A\ar[r] \ar[d] & B \ar[d] \\
A' \ar[r] & B'}$$
be a homotopy cofiber square in $Comm(\mathcal{C})$. Then the following
square is homotopy cocartesian in $B'-Mod$
$$\xymatrix{\mathbb{L}_{A}\otimes^{\mathbb{L}}_{A}B' \ar[r] \ar[d] &
\mathbb{L}_{B}\otimes^{\mathbb{L}}_{B}B' \ar[d] \\
\mathbb{L}_{A'}\otimes^{\mathbb{L}}_{A'}B' \ar[r] &
\mathbb{L}_{B'}.}$$

\item For any commutative monoid $A$ and any $A$-module $M$, one has
a natural isomorphism in $\mathrm{Ho}(A-Mod)$
$$\mathbb{L}_{A\oplus M}\otimes_{A\oplus M}^{\mathbb{L}}A\simeq
\mathbb{L}_{A}\coprod \mathbb{L}QZ(M),$$
where
$$Q : A-Comm_{nu}(\mathcal{C}) \longrightarrow A-Mod \qquad
A-Comm_{nu}(\mathcal{C}) \longleftarrow A-Mod : Z$$
is the Quillen adjunction used during the proof of
\ref{p1}.
\end{enumerate}
\end{prop}

\begin{proof} $(1)$ to $(3)$ are simple exercises, using the definitions and
that for any morphism of commutative monoids $A \longrightarrow B$,
and any $B$-module $M$, the following square is homotopy cartesian
in $Comm(\mathcal{C})$ (because of our assumption \ref{ass-1})
$$\xymatrix{
A\oplus M \ar[r] \ar[d] & B\oplus M \ar[d] \\
A\ar[r] & B.}$$

$(4)$ We note that for any commutative monoid $A$, and any
$A$-modules $M$ and $N$, one has a natural
homotopy fiber sequence
$$Map_{A-Comm(\mathcal{C})/A}(A\oplus M,A\oplus N) \longrightarrow 
Map_{Comm(\mathcal{C})/A}(A\oplus M,A\oplus N) \rightarrow $$
$$\longrightarrow Map_{Comm(\mathcal{C})/A}(A,A\oplus N),$$
or equivalently using lemma \ref{l1}
$$Map_{A-Mod}(\mathbb{L}QZ(M),N) \longrightarrow
\mathbb{D}er(A\oplus M,N) \longrightarrow \mathbb{D}er(A,N).$$
This
implies the existence of a natural homotopy cofiber sequence
of $A$-modules
$$\mathbb{L}_{A} \longrightarrow \mathbb{L}_{A\oplus M}\otimes_{A\oplus M}^{\mathbb{L}}A \longrightarrow \mathbb{L}QZ(M).$$
Clearly this sequence splits in $\mathrm{Ho}(A-Mod)$ and gives rise to a
natural isomorphism
$$\mathbb{L}_{A\oplus M}\otimes_{A\oplus M}^{\mathbb{L}}A\simeq
\mathbb{L}_{A}\coprod \mathbb{L}QZ(M).$$
\end{proof}

The importance of derivations come from the fact that
they give rise to infinitesimal extensions in the following way.
Let $A \longrightarrow B$ be a morphism of commutative monoids
in $\mathcal{C}$, and $M$ be a $B$-module. Let $d : \mathbb{L}_{B/A} \longrightarrow
M$ be a morphism in $\mathrm{Ho}(B-Mod)$, corresponding to
a derivation $d\in \pi_{0}(\mathbb{D}er_{A}(B,M))$.
This derivation can be seen as a section
$d : B\longrightarrow B\oplus M$
of the morphism of commutative $A$-algebras
$B\oplus M \longrightarrow B$. We consider
the following homotopy cartesian diagram in the
category of commutative $A$-algebras
$$\xymatrix{
C \ar[r] \ar[d] & B\ar[d]^-{d} \\
B \ar[r]_-{s} & B\oplus M}$$
where $s : B \longrightarrow B\oplus M$
is the natural section corresponding
to the zero morphism $\mathbb{L}_{B/A} \longrightarrow
M$. Then, $C \longrightarrow B$ is a morphism
of commutative $A$-algebras such that its fiber
is a non-unital commutative $A$-algebra
isomorphic in $\mathrm{Ho}(A-Comm_{nu}(\mathcal{C}))$ to the
loop $A$-module $\Omega M:=*\times^{h}_{M}*$
with the zero multiplication. In other words,
$C$ is a \emph{square zero extension of $B$ by
$\Omega M$}. It will be denoted by
$B\oplus_{d}\Omega M$. The most important case is
of course when $A=B$, and we make the following definition.

\begin{df}\label{dext}
Let $A$ be a commutative monoid, $M$ and
$A$-module and $d\in \pi_{0}\mathbb{D}er(A,M)$ be
a derivation given by a morphism
in $d : A \longrightarrow A\oplus M$
in $\mathrm{Ho}(Comm(\mathcal{C})/A)$. The \emph{square zero
extension associated to $d$}\index{square zero
extension}, denoted
by $A\oplus_{d}\Omega M$, is
defined as the homotopy pullback diagram
of commutative monoids
$$\xymatrix{
A\oplus_{d}\Omega M \ar[r] \ar[d] & A \ar[d]^-{d} \\
A \ar[r]_-{s} & A\oplus M,}$$
where $s$ is the natural morphism corresponding to
the zero derivation. The top horizontal morphism
$A\oplus_{d}\Omega M \longrightarrow A$ will be called
the natural projection.
\end{df}

\section{Hochschild homology}

For a commutative monoid $A \in Comm(\mathcal{C})$, we set
$$THH(A):=S^{1}\otimes^{\mathbb{L}}A \in \mathrm{Ho}(Comm(\mathcal{C})),$$
where $S^{1}\otimes^{\mathbb{L}}$ denotes the
derived external product of object in $Comm(\mathcal{C})$ by
the simplicial circle $S^{1}:=\Delta^{1}/\partial \Delta^{1}$.
Presenting the circle $S^{1}$ has the
homotopy push-out
$$*\coprod^{\mathbb{L}}_{*\coprod *}*$$
one gets that
$$THH(A)\simeq A\otimes^{\mathbb{L}}_{A\otimes^{\mathbb{L}}A}A.$$
The natural point $* \longrightarrow S^{1}$ induces a natural
morphism in $\mathrm{Ho}(Comm(\mathcal{C}))$
$$A \longrightarrow THH(A)$$
making $THH(A)$ as a commutative $A$-algebra, and as a natural
object in $\mathrm{Ho}(A-Comm(\mathcal{C}))$.

\begin{df}\label{d2'}
Let $A$ be a commutative monoid in $\mathcal{C}$. The
\emph{topological Hochschild homology of $A$}\index{topological Hochschild homology}
(or simply \emph{Hochschild homology})\index{Hochschild homology}
is the commutative $A$-algebra $THH(A):=S^{1}\otimes^{\mathbb{L}}A$.

More generally, if $A \longrightarrow B$ is a morphism
of commutative monoids in $\mathcal{C}$, the \emph{relative
topological Hochschild homology of $B$ over $A$}
(or simply \emph{relative Hochschild homology})
is the commutative $A$-algebra
$$THH(B/A):=THH(B)\otimes^{\mathbb{L}}_{THH(A)}A.$$
\end{df}

By definition, we have for any commutative monoid $B$,
$$Map_{Comm(\mathcal{C})}(THH(A),B)\simeq Map_{SSet}(S^{1},Map_{Comm(\mathcal{C})}(A,B)).$$
This implies that if $f : A \longrightarrow B$ is a morphism
of commutative monoids in $\mathcal{C}$, then we have
$$Map_{A-Comm(\mathcal{C})}(THH(A),B)\simeq \Omega_{f}Map_{Comm(\mathcal{C})}(A,B),$$
where $\Omega_{f}Map_{Comm(\mathcal{C})}(A,B)$ is the loop space
of $Map_{Comm(\mathcal{C})}(A,B)$ at the point $f$.
More generally, if $B$ and $C$ are commutative $A$-algebras,
then
$$Map_{A-Comm(\mathcal{C})}(THH(B/A),C)\simeq Map_{SSet}(S^{1},Map_{A-Comm(\mathcal{C})}(B,C)),$$
and for a morphism $f : B \longrightarrow C$ of commutative $A$-algebras
$$Map_{B-Comm(\mathcal{C})}(THH(B/A),C)\simeq \Omega_{f}Map_{A-Comm(\mathcal{C})}(B,C).$$

\begin{prop}\label{p2'}
\begin{enumerate}
\item Let $A \longrightarrow B \longrightarrow C$ be two morphisms
in $Comm(\mathcal{C})$. Then, the natural morphism
$$THH(C/A)\otimes_{THH(B/A)}^{\mathbb{L}}B\longrightarrow THH(C/B)$$
is an isomorphism in $\mathrm{Ho}(B-Comm(\mathcal{C}))$.

\item
Let
$$\xymatrix{A\ar[r] \ar[d] & B \ar[d] \\
A' \ar[r] & B'}$$
be a homotopy cofiber square in $Comm(\mathcal{C})$. Then, the natural
morphism
$$THH(B/A)\otimes_{A}^{\mathbb{L}}THH(A'/A)\longrightarrow THH(B'/A)$$
is an isomorphism in $\mathrm{Ho}(A-Comm(\mathcal{C}))$.

\end{enumerate}
\end{prop}

\begin{proof} Exercise.  \end{proof}

\section{Finiteness conditions}

We present here two different finiteness
conditions for objects in model categories. The first one
is valid in any model category, and is a homotopy analog
of the notion of finitely
presented object in a category. The second
one is only valid for symmetric monoidal categories,
and is a homotopy generalization of the notion
of rigid objects in monoidal categories.

\begin{df}\label{d3}
A morphism $x\longrightarrow y$ in a proper model category
$M$
is \emph{finitely presented}\index{finitely presented morphism} (we also say that
\emph{$y$ is finitely presented over $x$})
if for any filtered diagram of
objects under $x$, $\{z_{i}\}_{i\in I}
\in x/M$,
the natural morphism
$$Hocolim_{i\in I} Map_{x/M}(y,z_{i})
\longrightarrow Map_{x/M}(y,Hocolim_{i\in I}z_{i})$$
is an isomorphism in $\mathrm{Ho}(SSet)$.
\end{df}

\begin{rmk}\label{r3}
\emph{If the model category $M$ is not proper the definition
\ref{d3} has to be modified by replacing $Map_{x/M}$
with $Map_{Qx/M}$, where $Qx$ is a cofibrant model
for $x$. By our assumption \ref{ass-1} all the model categories we
will use are proper.}
\end{rmk}

\begin{prop}\label{p3}
Let $M$ be a proper model category.
\begin{enumerate}
\item
Finitely presented morphisms in $M$
are stable by equivalences. In other words,
if one has a commutative diagram in $M$
$$\xymatrix{x \ar[r]^{f} \ar[d]_{p} & y\ar[d]^{q} \\
x' \ar[r]_{f'} & y'}$$
such that $p$ and $q$ are equivalences, then
$f$ is finitely presented if and only
if $f'$ is finitely presented.

\item Finitely presented morphisms in $M$
are stable by compositions and retracts.

\item Finitely presented morphisms in $M$ are stable by
homotopy push-outs. In other words, if
one has a homotopy push-out diagram in $M$
$$\xymatrix{x \ar[r]^{f} \ar[d]_{p} & y\ar[d]^{q} \\
x' \ar[r]_{f'} & y'}$$
then $f'$ is finitely presented if $f$ is so.
\end{enumerate}
\end{prop}

\begin{proof} $(1)$ is clear as $Map_{x/M}(a,b)$
only depends on the isomorphism class
of $a$ and $b$ as objects in the homotopy category
$\mathrm{Ho}(x/M)$.

$(2)$ Let $x \longrightarrow y\longrightarrow z$ be two
finitely presented morphisms in $M$, and let  $\{z_{i}\}_{i\in I}
\in x/M$ be a filtered diagrams of objects. Then, one has
for any object $t \in x/M$ a fibration sequence of
simplicial sets
$$Map_{y/M}(z,t) \longrightarrow
Map_{x/M}(z,t) \longrightarrow Map_{x/M}(y,t).$$
As fibration sequences are stable by filtered homotopy colimits,
one gets a morphism of fibration sequences
$$\xymatrix{
\mbox{\em Hocolim}_{i\in I}\mbox{\em Map}_{y/M}(z,z_{i}) \ar[r] \ar[d] &
\mbox{\em Hocolim}_{i\in I}\mbox{\em Map}_{x/M}(z,z_{i}) \ar[r] \ar[d] &
\mbox{\em Hocolim}_{i\in I}\mbox{\em Map}_{x/M}(y,z_{i}) \ar[d] \\
\mbox{\em Map}_{y/M}(z,\mbox{\em Hocolim}_{i\in I}z_{i}) \ar[r] &
\mbox{\em Map}_{x/M}(z,\mbox{\em Hocolim}_{i\in I}z_{i}) \ar[r] &
\mbox{\em Map}_{x/M}(y,\mbox{\em Hocolim}_{i\in I}z_{i}), }$$
and the five lemma tells us that the
vertical arrow in the middle is an isomorphism in
$\mathrm{Ho}(SSet)$. This implies that $z$ is finitely presented over $x$.

The assertion concerning retracts is clear since, 
if $x\longrightarrow y$ is a retract of
$x' \longrightarrow y'$, for any $z\in x/M$ the simplicial set
$Map_{x/M}(y,z)$ is a retract of  $Map_{x'/M}(y',z)$. \\

$(3)$ This is clear since we have for any
object $t \in x'/M$, a natural equivalence
$Map_{x'/M}(y',t)\simeq Map_{x/M}(y,t)$.  
\end{proof}

Let us now fix $I$, a set of generating cofibrations
in $M$.

\begin{df}\label{dI-cell}
\begin{enumerate}
\item An object $X$ is a \emph{strict finite $I$-cell object}\index{object!strict finite cell},
if there exists a finite sequence
$$\xymatrix{
X_{0}=\emptyset \ar[r] & X_{1} \ar[r] & \dots \ar[r] & X_{n}=X,}$$
and for any $0\leq i<n$ a push-out square
$$\xymatrix{
X_{i} \ar[r] & X_{i+1} \\
A \ar[u] \ar[r]_-{u_{i}} & B, \ar[u]}$$
with $u_{i} \in I$.
\item An object $X$ is a \emph{finite $I$-cell object}\index{object!finite cell}
(or simply a \emph{finite cell object} when $I$ is clear)
if it is equivalent to a strict finite $I$-cell object.
\item
The model category $M$ is \emph{compactly generated}\index{model!category!compactly generated}
if it satisfies the following conditions.
\begin{enumerate}
\item The model category $M$ is cellular (in the sense
of \cite[\S 12]{hi}).
\item  There exists a set of generating cofibrations
$I$, and generating trivial cofibrations $J$
whose domains and codomains are cofibrant,
$\omega$-compact (in the sense of \cite[\S 10.8]{hi}) and
$\omega$-small with respect to the whole category $M$.
\item Filtered colimits commute with finite limits in $M$.
\end{enumerate}
\end{enumerate}
\end{df}

The following proposition identifies 
finitely presented objects when $M$ is 
compactly generated. 

\begin{prop}\label{pI-cell}
Let $M$ be a compactly generated model category, and
$I$ be a set of generating cofibrations
whose domains and codomains are cofibrant, $\omega$-compact
and $\omega$-small with respect to the whole category $M$.
\begin{enumerate}
\item A filtered colimit of fibrations (resp. trivial
fibrations) is a fibration (resp. a trivial fibration).
\item For any filtered diagram $X_{i}$ in $M$, the natural
morphism
$$Hocolim_{i}X_{i} \longrightarrow Colim_{i}X_{i}$$
is an isomorphism in $\mathrm{Ho}(M)$.

\item Any object $X$ in $M$ is equivalent to a filtered
colimit of strict finite $I$-cell objects.

\item An object $X$ in $M$ is finitely
presented if and only if it is equivalent to a
retract, in $\mathrm{Ho}(M)$, of a strict finite $I$-cell object.

\end{enumerate}
\end{prop}

\begin{proof} $(1)$ By assumption the domain and codomain
of morphisms of $I$ are $\omega$-small, so $M$ is finitely generated
in the sense of \cite[\S 7]{ho}. Property $(1)$ is then proved in 
\cite[\S 7]{ho}. \\

$(2)$ For a filtered category $A$, the colimit functor
$$Colim : M^{A} \longrightarrow M$$
is a left Quillen functor for the levelwise projective model structure on $M$. 
By $(1)$ we know that $Colim$ preserves trivial fibrations, and thus
that it also preserves equivalences. We therefore have
isomorphisms of functors 
$Hocolim\simeq \mathbb{L}Colim\simeq Colim$. \\

$(3)$ The small object argument (e.g. \cite[Thm. 2.1.14]{ho}) gives that
any object $X$ is equivalent to a $I$-cell complex $Q(X)$.
By $\omega$-compactness of the domains and codomains of $I$,
$Q(X)$ is the filtered colimit of its finite sub-$I$-cell complexes.
This implies that $X$ is equivalent to a filtered colimit
of strict finite $I$-cell objects. \\

$(4)$ Let $A$ be a filtered category, and $Y \in M^{A}$ be a $A$-diagram.
Let $c(Y) \longrightarrow R_{*}(Y)$ be a Reedy fibrant replacement
of the constant simplicial object $c(Y)$ with values $Y$ 
(in the model category of simplicial objects in $M^{A}$, see \cite[\S 5.2]{ho}).
By $(2)$, the induced morphism
$$c(Colim_{a\in A} Y_{a}) \longrightarrow Colim_{a\in A} R_{*}(Y_{a})$$
is an equivalence of simplicial objects in $M^{A}$. Moreover, 
$(1)$ and the exactness of filtered colimits implies that 
$Colim_{a\in A} R_{*}(Y_{a})$ is a Reedy fibrant object in the model
categroy of simplicial objects in $M$ (as filtered colimits
commute with matching objects for the Reedy category $\Delta^{op}$, see \cite[\S 5.2]{ho}). 
This implies that for any cofibrant and $\omega$-small object $K$ in $M$, 
we have
$$Hocolim_{a\in A}Map(K,Y_{a})\simeq
Colim_{a\in A}Map(K,Y_{a}) \simeq 
Colim_{a\in A}Hom(K,R_{*}(Y_{a}))\simeq $$
$$Hom(K,Colim_{a\in A}R_{*}(Y_{a}))\simeq Map(K,Colim_{a\in A}Y_{a}).$$
This implies that the domains and codomains of 
$I$ are homotopically finitely presented. 

As filtered
colimits of simplicial sets preserve homotopy pull-backs, we deduce that
any finite cell objects is also homotopically finitely presented,
as they are constructed from domains and codomains of $I$
by iterated homotopy push-outs (we use
here that domains and codomains of $I$ are cofibrant).
This implies that
any retract of a finite cell object is
homotopically finitely presented. Conversely,
let $X$ be a homotopically finitely presented object in $Ho(M)$, and
by $(3)$ let us write it as $Colim_{i}X_{i}$, where
$X_{i}$ is a filtered diagram of finite cell objects. Then,
$[X,X]\simeq Colim_{i}[X,X_{i}]$, which
implies that this identity of $X$ factors through some $X_{i}$, or in other words
that $X$ is a retract in $Ho(M)$ of some $X_{i}$. 
\end{proof}

Now, let $M$ be a symmetric monoidal model category in the sense
of \cite[\S 4]{ho}. We remind that this implies in particular that
the monoidal structure on $M$ is closed, and therefore possesses
$Hom$'s objects $\underline{\mathbf{Hom}}_{M}(x,y)\in M$ satisfying the usual adjunction rule
$$Hom(x,\underline{\mathbf{Hom}}_{M}(y,z)))\simeq Hom(x\otimes y,z).$$
The internal structure can be derived, and gives
on one side a symmetric monoidal structure $-\otimes^{\mathbb{L}}-$
on $\mathrm{Ho}(M)$, as well as $Hom$'s objects
$\mathbb{R}\underline{\mathbf{Hom}}_{M}(x,y) \in \mathrm{Ho}(M)$ satisfying the
derived version of the previous adjunction
$$[x,\mathbb{R}\underline{\mathbf{Hom}}_{M}(y,z))]\simeq [x\otimes^{\mathbb{L}} y,z].$$
In particular, if $\mathbf{1}$ 
is the unit of the monoidal structure of $M$, then
$$[\mathbf{1},\mathbb{R}\underline{\mathbf{Hom}}_{M}(x,y)]\simeq [x,y],$$
and more generally
$$Map_{M}(\mathbf{1},\mathbb{R}\underline{\mathbf{Hom}}_{M}(x,y))\simeq Map_{M}(x,y).$$
Moreover, the adjunction between $-\otimes^{\mathbb{L}}$ and
$\mathbb{R}\underline{\mathbf{Hom}}_{M}$ extends naturally to an adjunction isomorphism
$$\mathbb{R}\underline{\mathbf{Hom}}_{M}(x,\mathbb{R}\underline{\mathbf{Hom}}_{M}(y,z)))\simeq
\mathbb{R}\underline{\mathbf{Hom}}_{M}(x\otimes^{\mathbb{L}} y,z).$$
The derived dual of an object $x\in M$ will be denoted by
$$x^{\vee}:=\mathbb{R}\underline{\mathbf{Hom}}_{M}(x,\mathbf{1}).$$

\begin{df}\label{d4}
Let $M$ be a symmetric monoidal model category.
An object $x\in M$ is called \emph{perfect}\index{object!perfect}
if the natural morphism
$$x\otimes^{\mathbb{L}}x^{\vee} \longrightarrow
\mathbb{R}\underline{\mathbf{Hom}}_{M}(x,x)$$
is an isomorphism in $\mathrm{Ho}(M)$.
\end{df}

\begin{prop}\label{p4}
Let $M$ be a symmetric monoidal model category.
\begin{enumerate}
\item If $x$ and $y$ are perfect objects in
$M$, then so is $x\otimes^{\mathbb{L}}y$.
\item If $x$ is a perfect object in $M$, then
for any objects $y$ and $z$, the natural morphism
$$\mathbb{R}\underline{\mathbf{Hom}}_{M}(y,x\otimes^{\mathbb{L}}z)
\longrightarrow
\mathbb{R}\underline{\mathbf{Hom}}_{M}(y\otimes^{\mathbb{L}}x^{\vee},z)$$
is an isomorphism in $\mathrm{Ho}(M)$.
\item If $x$ in perfect in $M$, and $y \in M$, then the natural morphism
$$\mathbb{R}\underline{\mathbf{Hom}}_{M}(x,y) \longrightarrow x^{\vee}\otimes^{\mathbb{L}}y$$
is an isomorphism in $\mathrm{Ho}(M)$.

\item If $\mathbf{1}$ is finitely presented in $M$, then
so is any perfect object.

\item If $M$ is furthermore a stable model category
then perfect objects are stable by
homotopy push-outs and homotopy pullbacks. In other words,
if $\xymatrix{x\ar[r]& y \ar[r] & z}$ is a homotopy fiber
sequence in $M$, and if two of the objects $x$, $y$, and $z$
are perfect then so is the third.

\item Perfect objects are stable by retracts in $\mathrm{Ho}(M)$.

\end{enumerate}
\end{prop}

\begin{proof} $(1)$, $(2)$ and $(3)$ are standard, as perfect objects
are precisely the strongly dualizable objects of the closed
monoidal category $\mathrm{Ho}(M)$ (see for example \cite{m}). \\

$(4)$ Let $x$ be a perfect object in $M$, and $\{z_{i}\}_{i\in I}$ be a filtered diagram
of objects in $M$. Let $x^{\vee}:=\mathbb{R}\underline{Hom}(x,\mathbf{1})$ the dual
of $x$ in $\mathrm{Ho}(M)$. Then, we have
$$Map_{M}(x,Hocolim_{i}z_{i})\simeq Map_{M}(\mathbf{1},x^{\vee}\otimes^{\mathbb{L}}Hocolim_{i}z_{i})\simeq
Map_{M}(\mathbf{1},Hocolim_{i}x^{\vee}\otimes^{\mathbb{L}}z_{i})$$
$$\simeq Hocolim_{i}Map_{M}(\mathbf{1},x^{\vee}\otimes^{\mathbb{L}}z_{i})
\simeq Hocolim_{i}Map_{M}(x,z_{i}).$$

$(5)$ Let $\xymatrix{x\ar[r]& y \ar[r] & z}$ be a homotopy fiber
sequence in $M$. It is enough to prove that if $y$ and $z$
are perfect then so is $x$. For this, let $x^{\vee}$, $y^{\vee}$ and
$z^{\vee}$ the duals of $x$, $y$ and $z$.

One has a morphism of homotopy fiber sequences
$$\xymatrix{
x\otimes^{\mathbb{L}} z^{\vee} \ar[d] \ar[r] & x\otimes^{\mathbb{L}} y^{\vee} \ar[r] \ar[d] &
x\otimes^{\mathbb{L}} x^{\vee} \ar[d] \\
\mathbb{R}\underline{\mathbf{Hom}}_{M}(z,x) \ar[r] & \mathbb{R}\underline{\mathbf{Hom}}_{M}(y,x) \ar[r] &
\mathbb{R}\underline{\mathbf{Hom}}_{M}(x,x).}$$
The five lemma and point $(3)$ implies that the last vertical
morphism is isomorphism, and that $x$ is
perfect. \\

$(6)$ If $x$ is a retract of $y$, then the natural morphism
$$x\otimes^{\mathbb{L}}x^{\vee} \longrightarrow \mathbb{R}\underline{\mathbf{Hom}}_{M}(x,x)$$
is a retract of
$$y\otimes^{\mathbb{L}}y^{\vee} \longrightarrow \mathbb{R}\underline{\mathbf{Hom}}_{M}(y,y).$$
\end{proof}

The following corollary gives a condition under which 
perfect and finitely presented objects are the same.

\begin{cor}\label{cI-cell}
Suppose that $M$ is a stable and compactly generated 
symmetric monoidal model category with
$\mathbf{1}$ being $\omega$-compact and cofibrant. We assume
that the set $I$ of morphisms of the form
$$S^{n}\otimes \mathbf{1} \longrightarrow \Delta^{n+1}\otimes
\mathbf{1}$$
is a set of generating cofibration for $M$.  
Then an object $x$ in $M$ is perfect if and only 
if it is finitely presented, and if and only if
it is a retract of a finite $I$-cell object. 
\end{cor}

\begin{proof}
This essentially follows from Prop. \ref{pI-cell} and
Prop. \ref{p4}, the only statement which  remains to be proved
is that retract of finite $I$-cell objects are perfect. 
But this follows from the fact that 
$\mathbf{1}$ is always perfect, and 
from Prop. \ref{p4} $(5)$ and $(6)$. 
\end{proof}

\begin{rmk}\label{r4}
\emph{The notion of finitely presented morphisms and
perfect objects depend on the model structure
and not only on the underlying category $M$. They specialize to
the corresponding usual categorical notions when $M$ is endowed
with the trivial model structure.}
\end{rmk}

\section{Some properties of modules}

In this Section we give some general notions
of flatness and projectiveness of
modules over a commutative monoid in $\mathcal{C}$.

\begin{df}\label{d7-}
Let $A \in Comm(\mathcal{C})$ be a commutative monoid, and
$M$ be an $A$-module.
\begin{enumerate}
\item The $A$-module $M$ is \emph{flat}\index{flat!module} if
the functor
$$-\otimes_{A}^{\mathbb{L}}M : \mathrm{Ho}(A-Mod) \longrightarrow \mathrm{Ho}(A-Mod)$$
preserves homotopy pullbacks.

\item The $A$-module $M$ is \emph{projective}\index{projective!module}
if it is a retract in $\mathrm{Ho}(A-Mod)$ of
$\coprod_{E}^{\mathbb{L}}A$ for some $\mathbb{U}$-small
set $E$.

\end{enumerate}
\end{df}

\begin{prop}\label{pd7-}
Let $A \longrightarrow B$ be a morphism
of commutative monoids in $\mathcal{C}$.
\begin{enumerate}
\item The free $A$-module $A^{n}$ of rank  $n$ is flat. 
Moreover, if infinite direct sums in $\mathrm{Ho}(A-Mod)$ 
commute with homotopy pull-backs, then for any $\mathbb{U}$-small set $E$, 
the free $A$-module $\coprod_{E}^{\mathbb{L}}A$
is flat. 
\item Flat modules in $\mathrm{Ho}(A-Mod)$ are stable by derived tensor
products, finite coproducts and retracts.
\item Projective modules in $\mathrm{Ho}(A-Mod)$
are stable by derived tensor products, finite coproducts
and retracts.
\item If $M$ is a flat (resp. projective) $A$-module, then
$M\otimes_{A}^{\mathbb{L}}B$ is a flat
(resp. projective) $B$-module.
\item A perfect $A$-module is flat.
\item Let us suppose that
$\mathbf{1}$ is a finitely presented object in $\mathcal{C}$.
Then, a projective $A$-module
is finitely presented if and only if it
is a retract of $\coprod_{E}^{\mathbb{L}}A$
for some finite set $E$.
\item Let us assume that $\mathbf{1}$ is a finitely presented
object in $\mathcal{C}$. Then, a projective $A$-module
is perfect if and only it is finitely presented.
\end{enumerate}
\end{prop}

\begin{proof} $(1)$ to $(4)$ are easy and follow from
the definitions. \\

$(5)$ Let $M$ be a perfect $A$-module, and
$M^{\vee}:=\mathbb{R}\underline{Hom}_{A}(M,A)$ be its dual.
Then, for any homotopy cartesian square of $A$-modules
$$\xymatrix{
P \ar[r] \ar[d] & Q \ar[d] \\
P'\ar[r] & Q',}$$
the diagram
$$\xymatrix{
P\otimes_{A}^{\mathbb{L}}M \ar[r] \ar[d] & Q\otimes_{A}^{\mathbb{L}}M \ar[d] \\
P'\otimes_{A}^{\mathbb{L}}M\ar[r] & Q'\otimes_{A}^{\mathbb{L}}M}$$
is equivalent to
$$\xymatrix{
\mathbb{R}\underline{Hom}_{A}(M^{\vee},P) \ar[r] \ar[d] &
\mathbb{R}\underline{Hom}_{A}(M^{\vee},Q) \ar[d] \\
\mathbb{R}\underline{Hom}_{A}(M^{\vee},P') \ar[r]&
\mathbb{R}\underline{Hom}_{A}(M^{\vee},Q'),}$$
which is again homotopy cartesian by the general properties
of derived internal $Hom$'s. This shows that
$-\otimes^{\mathbb{L}}M$ preserves homotopy pullbacks, and hence
that $M$ is flat.  \\

$(6)$ Clearly, if $E$ a finite set, then
$\coprod_{E}^{\mathbb{L}}A$ is finitely presented object, as
for any $A$-module $M$ one has
$$Map_{A-Mod}(\coprod_{E}^{\mathbb{L}}A,M)\simeq
Map_{\mathcal{C}}(\mathbf{1},M)^{E}.$$
Therefore, a retract of $\coprod_{E}^{\mathbb{L}}A$
is also finitely presented.

Conversely, let $M$ be a  projective $A$-module which
is also finitely presented. Let $i : M \longrightarrow
\coprod_{E}^{\mathbb{L}}A$ be a morphism which admits
a retraction. As $\coprod_{E}^{\mathbb{L}}A$ is the colimit
of $\coprod_{E_{0}}^{\mathbb{L}}A$, for $E_{0}$ running over the
finite subsets of $E$, the morphism
$i$ factors as
$$M \longrightarrow \coprod_{E_{0}}^{\mathbb{L}}A \longrightarrow
\coprod_{E}^{\mathbb{L}}A$$
for some finite subset $E_{0} \subset E$. This shows that
$M$ is in fact a retract of $\coprod_{E}^{\mathbb{L}}A$. \\

$(7)$ Using Prop. \ref{p4} $(4)$ one sees that if
$M$ is perfect then it is finitely presented. Conversely, let
$M$ be a finitely presented  projective $A$-module.
By $(6)$ we know that $M$ is a retract of
$\coprod_{E}^{\mathbb{L}}A$ for some finite set $E$.
But, as $E$ is finite, $\coprod_{E}^{\mathbb{L}}A$
is perfect, and therefore so is $M$ as a retract
of a perfect module. \end{proof}

\section{Formal coverings}

The following notion will be highly used, and is
a categorical version of faithful morphisms
of affine schemes.

\begin{df}\label{dcov}
A family of morphisms of commutative
monoids $\{f_{i} : A \rightarrow B_{i}\}_{i\in I}$ is
a \emph{formal covering}\index{formal covering (family)} if
the family of functors
$$\{\mathbb{L}f_{i}^{*} : \mathrm{Ho}(A-Mod) \longrightarrow \mathrm{Ho}(B_{i}-Mod)\}_{i\in I}$$
is conservative (i.e. a morphism $u$ in $\mathrm{Ho}(A-Mod)$ is an isomorphism
if and only if all the $\mathbb{L}f_{i}^{*}(u)$ are isomorphisms).
\end{df}

The formal covering families are stable by
equivalences, homotopy push-outs and compositions
and therefore do form a model topology in the
sense of \cite[Def. 4.3.1]{hagI} (or Def. \ref{modtop}).

\begin{prop}\label{pcov}
Formal covering families form a
model topology (Def. \ref{modtop}) on the model category
$Comm(\mathcal{C})$.
\end{prop}

\begin{proof} Stability by equivalences and composition
is clear. The stability by homotopy push-outs
is an easy consequence of the
transfer formula Prop. \ref{ptrans}.  \end{proof}

\section{Some properties of morphisms}

In this part we review several classes of morphisms
of commutative monoids in $\mathcal{C}$, generalizing the
usual notions of Zariski open immersions, unramified, \'etale,
smooth and flat morphisms of affine schemes. It is interesting to
notice that these notions make sense in our general context, but
specialize to very different notions in various
specific cases (see the examples at the beginning
of each section \S \ref{IIsim}, \S \ref{IIder}, \S \ref{IIunb} and \S \ref{IIbnag}). \\

\begin{df}\label{d5}
Let $f : A \longrightarrow B$ be a morphism
of commutative monoids in $\mathcal{C}$.
\begin{enumerate}

\item The morphism $f$ is an \emph{epimorphism}\index{morphism of commutative monoids!epimorphism} if for any
commutative $A$-algebra $C$ the simplicial set
$Map_{A-Comm(\mathcal{C})}(B,C)$ is either empty or
contractible.

\item The morphism $f$ is \emph{flat}\index{morphism of commutative monoids!flat} if
the induced functor
$$\mathbb{L}f^{*} : \mathrm{Ho}(A-Mod) \longrightarrow \mathrm{Ho}(B-Mod)$$
commutes with finite homotopy limits.

\item The morphism $f$ is a \emph{formal Zariski
open immersion}\index{morphism of commutative monoids!formal Zariski
open immersion} if it is flat and if the functor
$$f_{*} : \mathrm{Ho}(B-Mod) \longrightarrow \mathrm{Ho}(A-Mod)$$
is fully faithful.

\item The morphism $f$ is \emph{formally unramified}\index{morphism of commutative monoids!formally unramified}
if $\mathbb{L}_{B/A}\simeq 0$ in $\mathrm{Ho}(B-Mod)$.

\item The morphism $f$ is \emph{formally \'etale}\index{morphism of commutative monoids!formally \'etale}
if the natural morphism
$$\mathbb{L}_{A}\otimes^{\mathbb{L}}_{A}B \longrightarrow
\mathbb{L}_{B}$$
is an isomorphism in $\mathrm{Ho}(B-Mod)$.

\item The morphism $f$ is \emph{formally thh-\'etale}\index{morphism of commutative monoids!formally thh-\'etale}
if the natural morphism $$B \longrightarrow THH(B/A)$$
is an isomorphism in $\mathrm{Ho}(Comm(\mathcal{C}))$.

\end{enumerate}
\end{df}

\begin{rmk}\label{explepi}\emph{One remark concerning our notion of epimorphism of commutative monoids is in order. First of all, 
in a category $\mathcal{C}$ (without any model structure) having fiber products, 
a morphism $x \longrightarrow y$ is a monomorphism if and only 
the diagonal morphism $x \longrightarrow x\times_{y}x$ is an isomorphism.
A natural generalization of this fact gives a notion of 
monomorphism in any model category $M$, as a morphism $x \longrightarrow y$ whose diagonal 
$x \longrightarrow x\times^{h}_{y}x$ is an isomorphism in $\mathrm{Ho}(M)$. 
Equivalently, the morphism $x \longrightarrow y$ is a monomorphism
if and and only if for any $z\in M$ the induced morphism of simplicial sets
$Map_{M}(z,x) \longrightarrow Map_{M}(z,y)$ is 
a monomorphism in the model category $SSet$. Furthermore, it is easy to check that 
a morphism $f : K \longrightarrow L$ is a monomorphism in the model category 
$SSet$ if and only if for any $s\in L$ the homotopy fiber of $f$ at $s$ is either
empty or contractible. Therefore we see that 
a morphism $A\longrightarrow B$ in $Comm(\mathcal{C})$ is an epimorphism in the sense
of Def. \ref{d5} $(1)$ if and only if it is a monomorphism when considered
as a morphism in the model category $Comm(\mathcal{C})^{op}$, or equivalently if
the induced morphism $B\otimes_{A}^{\mathbb{L}}B \longrightarrow B$
is an isomorphism in $\mathrm{Ho}(Comm(\mathcal{C}))$.
This justifies our
terminology, and moreover shows that when the model structure on the category $M$ is trivial, 
an epimorphism in the sense of our definition is nothing else than an epimorphism in $M$ 
in the usual categorical sense.} 
\end{rmk}

\begin{prop}\label{p4'}
Epimorphisms, flat morphisms, formal Zariski open
immersions, formally unramified morphisms, formally \'etale morphisms,
and formally $thh$-\'etale morphisms,  are
all stable by compositions, equivalences and homotopy push-outs.
\end{prop}

\begin{proof} This is a simple exercise using the definitions
and Propositions \ref{ptrans}, \ref{p2}, \ref{p2'},
and \ref{p4}. 
\end{proof}

The relations between all these notions are given by
the following proposition.

\begin{prop}\label{p5}
\begin{enumerate}
\item A morphism $f : A \longrightarrow B$ is
an epimorphism if and only if the functor
$$f_{*} : \mathrm{Ho}(B-Mod) \longrightarrow \mathrm{Ho}(A-Mod)$$
is fully faithful.

\item A formal Zariski open immersion is
an epimorphism. A flat epimorphism is a formal Zariski
open immersion.

\item A morphism $f : A \longrightarrow B$
of commutative monoids is formally $thh$-\'etale
if and only if for any commutative $A$-algebra
$C$ the simplicial set $$Map_{A-Comm(\mathcal{C})}(B,C)$$
is discrete (i.e. equivalent to a set).

\item A formally \'etale morphism is
formally unramified.

\item An epimorphism is
formally unramified and
formally $thh$-\'etale.

\item A morphism $f : A \longrightarrow B$ in $Comm(\mathcal{C})$
is formally unramified if and only the morphism
$$B\otimes_{A}^{\mathbb{L}}B \longrightarrow B$$
is formally \'etale.

\end{enumerate}
\end{prop}

\begin{proof} $(1)$
Let $f : A \longrightarrow B$ be
a morphism such that the right Quillen
functor $f_{*} : B-Mod \longrightarrow A-Mod$
induces a fully faithful functor on the homotopy categories.
Therefore, the adjunction morphism
$$M\otimes_{A}^{\mathbb{L}}B \longrightarrow M$$
is an isomorphism for any $M\in \mathrm{Ho}(B-Mod)$.
In particular, the functor
$$f_{*} : \mathrm{Ho}(B-Comm(\mathcal{C})) \longrightarrow \mathrm{Ho}(A-Comm(\mathcal{C}))$$
is also fully faithful. Let $C$ be a commutative $A$-algebra, and
let us suppose that $Map_{A-Comm(\mathcal{C})}(B,C)$ is not empty.
This implies that $C$ is isomorphic in $\mathrm{Ho}(A-Comm(\mathcal{C}))$
to some $f_{*}(C')$ for $C'\in \mathrm{Ho}(B-Comm(\mathcal{C}))$.
Therefore, we have
$$Map_{A-Comm(\mathcal{C})}(B,C)\simeq Map_{A-Comm(\mathcal{C})}(f_{*}(B),f_{*}(C'))\simeq
Map_{B-Comm(\mathcal{C})}(B,C')\simeq *,$$
showing that $f$ is a epimorphism.

Conversely, let $f : A \longrightarrow B$ be epimorphism.
For any commutative $A$-algebra $C$, we have
$$Map_{A-Comm(\mathcal{C})}(B,C)\simeq Map_{A-Comm(\mathcal{C})}(B,C)\times Map_{A-Comm(\mathcal{C})}(B,C),$$
showing that the natural morphism $B\otimes_{A}^{\mathbb{L}}B\longrightarrow B$
is an isomorphism in $\mathrm{Ho}(A-Comm(\mathcal{C}))$.
This implies that for any $B$-module $M$, we have
$$M\otimes^{\mathbb{L}}_{A}B\simeq M\otimes_{B}^{\mathbb{L}}(B\otimes_{A}^{\mathbb{L}}B)\simeq M,$$
or in other words, that the adjunction morphism
$M \longrightarrow f_{*}\mathbb{L}f^{*}(M)\simeq M\otimes^{\mathbb{L}}_{A}B$
is an isomorphism in $\mathrm{Ho}(B-Mod)$. This means that
$f_{*} : \mathrm{Ho}(B-Mod) \longrightarrow \mathrm{Ho}(A-Mod)$
is fully faithful. \\

$(2)$ This is clear by $(1)$ and the definitions. \\

$(3)$ For any morphism of commutative $A$-algebras,
$f : B\longrightarrow C$ we have
$$Map_{B-Comm(\mathcal{C})}(THH(B/A),C))\simeq \Omega_{f}Map_{A-Comm(\mathcal{C})}(B,C).$$
Therefore, $B\longrightarrow THH(B/A)$ is an equivalence if and only
if for any such $f : B \longrightarrow C$, the simplicial set
$\Omega_{f}Map_{A-Comm(\mathcal{C})}(B,C)$ is contractible. Equivalently,
$f$ is formally thh-\'etale if and only if $Map_{A-Comm(\mathcal{C})}(B,C)$
is discrete. \\

$(4)$ Let $f : A \longrightarrow B$ be
a formally \'etale morphism of commutative monoids in $\mathcal{C}$.
By Prop. \ref{p2} $(1)$, there is a homotopy cofiber sequence of $B$-modules
$$\mathbb{L}_{A}\otimes_{A}^{\mathbb{L}}B \longrightarrow \mathbb{L}_{B}
\longrightarrow \mathbb{L}_{B/A},$$
showing that if the first morphism is an isomorphism
then $\mathbb{L}_{B/A}\simeq *$ and therefore that
$f$ is formally unramified. \\

$(5)$  Let $f : A \longrightarrow B$ an epimorphism. By definition
and $(3)$ we know that $f$ is formally $thh$-\'etale. Let us prove that
$f$ is also formally unramified.
Let $M$ be a $B$-module. As we have seen before, for any commutative
$B$-algebra $C$, the adjunction morphism
$C \longrightarrow C\otimes^{\mathbb{L}}_{A}B$ is
an isomorphism. In particular,
the functor
$$f_{*} : \mathrm{Ho}(B-Comm(\mathcal{C})/B) \longrightarrow \mathrm{Ho}(A-Comm(\mathcal{C})/B)$$
is fully faithful. Therefore we have
$$\mathbb{D}er_{A}(B,M)\simeq Map_{A-Comm(\mathcal{C})/B}(B,B\oplus M)\simeq
Map_{B-Comm(\mathcal{C})/B}(B,B\oplus M)\simeq *,$$
showing that $\mathbb{D}er_{A}(B,M)\simeq *$ for any
$B$-module $M$, or equivalently that
$\mathbb{L}_{B/A}\simeq *$. \\

$(6)$ Finally, Prop. \ref{p2} $(2)$
shows that the morphism $B\otimes_{A}^{\mathbb{L}}B \longrightarrow B$
is formally \'etale if and only if the natural morphism
$$\mathbb{L}_{B/A}\coprod \mathbb{L}_{B/A} \longrightarrow
\mathbb{L}_{B/A}$$
is an isomorphism in $\mathrm{Ho}(B-Mod)$. But
this is equivalent to
$\mathbb{L}_{B/A}\simeq *$.  \end{proof}

In order to state the next results we recall
that as $\mathcal{C}$ is a pointed model category
one can define a suspension functor (see \cite[\S 7]{ho})
$$\begin{array}{cccc}
S : & \mathrm{Ho}(\mathcal{C}) & \longrightarrow & \mathrm{Ho}(\mathcal{C}) \\
 & X & \mapsto & S(X):=*\coprod^{\mathbb{L}}_{X}*.
\end{array}$$
This functor possesses a right adjoint, the loop
functor
$$\begin{array}{cccc}
\Omega : & \mathrm{Ho}(\mathcal{C}) & \longrightarrow & \mathrm{Ho}(\mathcal{C}) \\
 & X & \mapsto & \Omega(X):=*\times^{h}_{X}*.
\end{array}$$

\begin{prop}\label{p6}
Assume that the base model category $\mathcal{C}$ is such that
the suspension functor
$S : \mathrm{Ho}(\mathcal{C}) \longrightarrow \mathrm{Ho}(\mathcal{C})$ is fully faithful.
\begin{enumerate}
\item A morphism of commutative monoids in $\mathcal{C}$ is
formally \'etale if and only if it is formally
unramified.
\item A formally $thh$-\'etale morphism of commutative
monoids in $\mathcal{C}$ is a formally \'etale morphism.
\item An epimorphism of commutative monoids in $\mathcal{C}$ is formally \'etale.
\end{enumerate}
\end{prop}

\begin{proof} $(1)$ By the last proposition we only need to prove that
a formally unramified morphism is also formally \'etale.
Let $f : A \longrightarrow B$ be such a morphism. 
By Prop. \ref{p2} $(1)$ there is a homotopy cofiber sequence
of $B$-modules
$$\mathbb{L}_{A}\otimes^{\mathbb{L}}_{A}B \longrightarrow
\mathbb{L}_{B} \longrightarrow *.$$
This implies that for any $B$-module $M$, the homotopy fiber of the morphism
$$\mathbb{D}er(B,M) \longrightarrow \mathbb{D}er(A,M)$$
is contractible, and in particular that this morphism induces
isomorphisms on all higher homotopy groups.
It remains to show that this morphism induces also
an isomorphism on $\pi_{0}$. For this, we can use the
hypothesis on $\mathcal{C}$ which implies that the suspension functor
on $\mathrm{Ho}(B-Mod)$ is fully faithful. Therefore, we have
$$\pi_{0}(\mathbb{D}er(B,M))\simeq \pi_{0}(\mathbb{D}er(B,\Omega S(M)))
\simeq \pi_{1}(\mathbb{D}er(B,SM))\simeq $$
$$\simeq \pi_{1}(\mathbb{D}er(A,SM))
\simeq \pi_{0}(\mathbb{D}er(A,M)).$$

$(2)$ Let $A \longrightarrow B$
be a formally $thh$-\'etale morphism in $Comm(\mathcal{C})$.\\
As $Map_{A-Comm(\mathcal{C})}(B,C)$ is discrete for any
commutative $A$-algebra $C$, the simplicial set
$\mathbb{D}er_{A}(B,M)$ is discrete for any
$B$-module $M$. Using the hypothesis on $\mathcal{C}$
we get that for any $B$-module $M$
$$\pi_{0}(\mathbb{D}er_{A}(B,M))\simeq \pi_{0}(\mathbb{D}er_{A}(B,\Omega S(M)))
\simeq \pi_{1}(\mathbb{D}er_{A}(B,SM))\simeq 0,$$
showing that $\mathbb{D}er_{A}(B,M)\simeq *$, and therefore
that $\mathbb{L}_{B/A}\simeq *$.  This implies that
$f$ is formally unramified, and therefore is
formally \'etale by the first part of the proposition. \\

$(3)$ This follows from $(2)$ and Prop. \ref{p5} $(5)$.
\end{proof}

The hypothesis of Proposition \ref{p6} saying that the suspension is fully faithful will appear 
in many places in the sequel. It is essentially equivalent to saying that the homotopy 
theory of $\mathcal{C}$ can be embedded in a stable homotopy theory 
in such a way that homotopy colimits are preserved (it could be in fact called
\emph{left semi-stable}). It is a very natural condition as many statements will 
then simplify, as shown for example by our infinitesimal theory
in \S 3. 

\begin{cor}\label{c1}
Assume furthermore that $\mathcal{C}$ is a stable model category, and
let $f : A \longrightarrow B$ be a morphism of commutative monoids
in $\mathcal{C}$. The following are equivalent.
\begin{enumerate}
\item The morphism $f$  is a formal Zariski open immersion.
\item The morphism $f$  is an epimorphism.
\item The morphism $f$  is a formally \'etale epimorphism.
\end{enumerate}
\end{cor}

\begin{proof} Indeed, in the stable case all
base change functors commute with limits, so
Prop. \ref{p5} $(1)$ and $(2)$ shows that
formal Zariski open immersions are exactly
the epimorphisms. Furthermore, by
Prop. \ref{p5} $(5)$ and Prop. \ref{p6}
all epimorphisms are formally \'etale.
\end{proof}

\begin{df}\label{d6}
A morphism of commutative monoids in $\mathcal{C}$ is
a \emph{Zariski open immersion}\index{Zariski open immersion! between commutative monoids} (resp. \emph{unramified}\index{morphism of commutative monoids!unramified},
resp. \emph{\'etale}\index{morphism of commutative monoids!\'etale}, resp.
\emph{$thh$-\'etale})\index{morphism of commutative monoids!thh-\'etale}  if it is
finitely presented (as a morphism in
the model category $Comm(\mathcal{C})$) and is
a formal Zariski open immersion (resp. formally unramified,
resp. formally \'etale, resp.
formally $thh$-\'etale).
\end{df}

Clearly, using what we have seen before,
Zariksi open immersions, unramified morphisms, \'etale
morphisms, and $thh$-\'etale morphisms
are all stable by equivalences, compositions and push-outs. \\

\section{Smoothness}

We define two general notions of smoothness, both
different generalizations of the usual notion, and both
useful in certain contexts.  A third, and still
different, notion of smoothness will be given in the next section.

\begin{df}\label{d7}
Let $f : A\longrightarrow B$
be a morphism of commutative algebras.

\begin{enumerate}

\item The morphism $f$ is \emph{formally perfect}
(or simply \emph{fp})\index{morphism of commutative monoids!fp}\index{morphism of commutative monoids!formally perfect}
if the $B$-module $\mathbb{L}_{B/A}$ is perfect
(in the sense of Def. \ref{d4}).
\item The morphism $f$ is \emph{formally smooth}\index{morphism of commutative monoids!formally smooth}
if the $B$-module $\mathbb{L}_{B/A}$ is projective
(in the sense of Def. \ref{d7-}) and if the
morphism
$$\mathbb{L}_{A}\otimes_{A}^{\mathbb{L}}B \longrightarrow 
\mathbb{L}_{B}$$
has a retraction in $\mathrm{Ho}(B-Mod)$.
\end{enumerate}
\end{df}

\begin{df}\label{d7'}
Let $f : B \longrightarrow C$
be a morphism of commutative algebras.
The morphism $f$ is
\emph{perfect}\index{morphism of commutative monoids!perfect}, or
simply \emph{p},\index{morphism of commutative monoids!p} (resp. \emph{smooth})\index{morphism of commutative monoids!smooth}  if it is
finitely presented (as a morphism in
the model category $Comm(\mathcal{C})$) and is
fp (resp. formally smooth).
\end{df}

Of course, (formally) \'etale morphisms are (formally) smooth morphisms
as well as (formally) perfect morphisms.

\begin{prop}\label{p7'}
The fp, perfect, formally smooth and smooth morphisms are all stable by compositions,
homotopy push outs and equivalences.
\end{prop}

\begin{proof} Exercise. \end{proof}

\section{Infinitesimal lifting properties}\label{Iinf}

While the first two notions of smoothness only depend on the underlying symmetric monoidal model category $\mathcal{C}$,
the third one, to be defined below, will also depend on the HA context we are working in.\\
Recall from Definition \ref{dha} that an \textit{HA context} is a
triplet $(\mathcal{C},\mathcal{C}_{0},\mathcal{A})$,
consisting of a symmetric monoidal model category
$\mathcal{C}$, two full sub-categories stable by equivalences
$$\mathcal{C}_{0}\subset \mathcal{C}
\qquad \mathcal{A}\subset Comm(\mathcal{C}),$$
such that:
\begin{itemize}
  \item $\mathbf{1} \in \mathcal{C}_{0}$, $\mathcal{C}_{0}$ is closed under
   by $\mathbb{U}$-small homotopy colimits, and 
   $X\otimes^{\mathbb{L}}Y\in \mathrm{Ho}(\mathcal{C}_{0})$ if $X$ and $Y$ in $\mathrm{Ho}(\mathcal{C}_{0})$.
	\item  any $A\in \mathcal{A}$ is $\mathcal{C}_{0}$-good (i.e. the functor
         $$\mathrm{Ho}(A-Mod) \longrightarrow \mathrm{Ho}((A-Mod_{0}^{op})^{\wedge})$$
         is fully faithful);
  \item assumptions \ref{ass-1},
  \ref{ass1}, \ref{ass0}, \ref{ass2}, \ref{ass7} are satisfied.     
 \end{itemize} 

\smallskip

Recall also that $\mathcal{C}_{1}$ is the full subcategory
of $\mathcal{C}$ consisting of all objects
equivalent to suspensions of objects in $\mathcal{C}_{0}$, 
$Comm(\mathcal{C})_{0}$ the full subcategory of
$Comm(\mathcal{C})$ consisting of commutative monoids
whose underlying $\mathcal{C}$-object is in $\mathcal{C}_{0}$, and, for $A\in Comm(\mathcal{C})$,
$A-Mod_{0}$\index{$A-Mod_{0}$} (resp. $A-Mod_{1}$\index{$A-Mod_{1}$}, 
resp. $A-Comm(\mathcal{C})_{0}$\index{$A-Comm(\mathcal{C})_{0}$}) is the full subcategory of
$A-Mod$ consisting of $A$-modules whose underlying $\mathcal{C}$-object is
in $\mathcal{C}_{0}$ (resp. of
$A-Mod$ consisting of $A$-modules whose underlying $\mathcal{C}$-object is
in $\mathcal{C}_{1}$, resp.
of $A-Comm(\mathcal{C})$ consisting of commutative $A$-algebras
whose underlying $\mathcal{C}$-object is in $\mathcal{C}_{0}$).\\

\begin{df}\label{dismooth}
Let $f : A \longrightarrow B$ be a morphism in $Comm(\mathcal{C})$.
\begin{enumerate}
\item
The morphism $f$ is called\index{formally infinitesimally
smooth!morphism between commutative monoids} \emph{formally infinitesimally
smooth relative to the HA context $(\mathcal{C },\mathcal{C}_{0},\mathcal{A})$} (or simply \emph{formally i-smooth} when
the HA context is clear)
if for any $R\in \mathcal{A}$, any morphism $A\longrightarrow R$
of commutative monoids, any
$M\in R-Mod_{1}$, and any $d\in \pi_{0}(\mathbb{D}er_{A}(R,M))$, 
the natural morphism
$$\pi_{0}\left( Map_{A-Comm(\mathcal{C})}(B,R\oplus_{d}\Omega M)\right)
\longrightarrow
\pi_{0}\left( Map_{A-Comm(\mathcal{C})}(B,R) \right)$$
is surjective.
\item The morphism $f$ is  called \emph{i-smooth}\index{morphism of commutative monoids!i-smooth} if it is
formally i-smooth  and finitely presented.
\end{enumerate}
\end{df}

The following proposition is immediate from the
definition.

\begin{prop}
Formally i-smooth
and i-smooth morphisms are stable
by equivalences, composition and homotopy push-outs.
\end{prop}

The next result provides a criterion for formally i-smooth morphisms
in terms of their cotangent complexes.

\begin{prop}\label{pismooth}
A morphism $f : A \longrightarrow B$ is
formally i-smooth if and only if
for any morphism $B\longrightarrow R$
with $R\in \mathcal{A}$, and
any $R$-module
$M\in R-Mod_{1}$,
the natural morphism
$$[\mathbb{L}_{R/A},M]\simeq \pi_{0}(\mathbb{D}er_{A}(R,M))) \longrightarrow \pi_{0}(\mathbb{D}er_{A}(B,M))\simeq
[\mathbb{L}_{B/A},M]_{B-Mod}$$
is zero.
\end{prop}

\begin{proof} Let us first assume that the
condition of the proposition is satisfied.
Let us
consider $R\in \mathcal{A}$, a morphism $A\longrightarrow R$, an
$R$-module
$M\in R-Mod_{1}$
and $d\in \pi_{0}(\mathbb{D}er_{A}(R,M))$.
The
homotopy fiber of the natural morphism
$$Map_{A-Comm(\mathcal{C})}(B,R\oplus_{d}\Omega M)
\longrightarrow
Map_{A-Comm(\mathcal{C})}(B,R)$$
taken at some morphism
$B \longrightarrow R$ in $\mathrm{Ho}(Comm(\mathcal{C}))$, can be
identified with the path space
$Path_{0,d'}\mathbb{D}er_{A}(B,M)$ from $0$ to $d'$ in
$\mathbb{D}er_{A}(B,M)$, where $d'$ is the image
of $d$ under the natural morphism
$$\pi_{0}(\mathbb{D}er_{A}(R,M))\longrightarrow \pi_{0}(\mathbb{D}er_{A}(B,M)).$$
By assumption, $0$ and $d'$ belong to the same
connected component of
$\mathbb{D}er_{A}(B,M)$, and thus we see that
$Path_{0,d'}\mathbb{D}er_{A}(B,M)$ is non-empty.
We have thus shown that
$$Map_{A-Comm(\mathcal{C})}(B,R\oplus_{d}\Omega M)
\longrightarrow
Map_{A-Comm(\mathcal{C})}(B,R)$$
has non-empty homotopy fibers and therefore that
$$\pi_{0}\left( Map_{A-Comm(\mathcal{C})}(B,R\oplus_{d}\Omega M)\right)
\longrightarrow
\pi_{0}\left( Map_{A-Comm(\mathcal{C})}(B,R) \right)$$
is surjective. The morphism $f$ is therefore
formally i-smooth.

Conversely, let $R\in \mathcal{A}$, $A\longrightarrow R$ a morphism, 
$M\in R-Mod_{1}$ and
$d\in \pi_{0}(\mathbb{D}er_{A}(R,M))$. Let $A \rightarrow R$ be a
morphism of commutative monoids.
We consider the diagram
$$\xymatrix{
A \ar[r] \ar[d] & R[\Omega_{d}M] \ar[d] \ar[r] & R \ar[d] \\
B \ar[r] & R \ar[r]_-{s} & R\oplus M. }$$
The homotopy fiber of
$$Map_{A-Comm(\mathcal{C})}(B,R\oplus_{d}\Omega M)
\longrightarrow
Map_{A-Comm(\mathcal{C})}(B,R)$$
is non-empty because $f$ is formally i-smooth.
By definition this means that the image of $d$ by the morphism
$$\pi_{0}(\mathbb{D}er_{A}(R,M))\longrightarrow \pi_{0}(\mathbb{D}er_{A}(B,M))$$
is zero. As this is true for any $d$, this finishes the proof of the proposition.
\end{proof}

\begin{cor}\label{cpismooth}
\begin{enumerate}
\item
Any formally unramified morphism is formally i-smooth.
\item
Assume that the suspension functor $S$ is fully faithful, that $\mathcal{C}_{1}=\mathcal{C}$
(so that in particular $\mathcal{C}$ is stable), and $\mathcal{A}=Comm(\mathcal{C})$. Then
formally i-smooth morphisms are
precisely the formally \'etale morphisms.
\end{enumerate}
\end{cor}

\begin{proof} It follows immediately from
\ref{pismooth}. \end{proof}

\begin{cor}\label{cpismooth2}
We assume that for any $M\in \mathcal{C}_{1}$
one has $[\mathbf{1},M]=0$.
Then any formally smooth morphism
is formally i-smooth.
\end{cor}

\begin{proof} By Prop. \ref{pismooth} and definition of
formal smoothness, it is
enough to show that for any commutative
monoid $A\in \mathcal{C}$, any
$A$-module $M\in A-Mod_{1}$
and any
projective $A$-module $P$ we have
$[P,M]=0$. By assumption this is true for
$P=A$, and thus also true for free $A$-modules and their retracts.
\end{proof}

\section{Standard localizations and Zariski open immersions}

For any object $A\in \mathcal{C}$, we can define
its \emph{underlying space} as
$$|A|:=Map_{\mathcal{C}}(\mathbf{1},A) \in \mathrm{Ho}(SSet_{\mathbb{U}}).$$
The model category $\mathcal{C}$ being pointed, the simplicial set
$|A|$ has a natural base point $* \in |A|$, and one can therefore
define the homotopy groups of $A$
$$\pi_{i}(A):=\pi_{i}(|A|,*).$$
When $A$ is the underlying object of a commutative monoid
$A\in Comm(\mathcal{C})$ one has by adjunction
$$\pi_{0}(A)=\pi_{0}(Map_{\mathcal{C}}(\mathbf{1},A))\simeq [\mathbf{1},A]_{\mathcal{C}}\simeq [A,A]_{A-Mod}.$$
Since $A$ is the unit of the monoidal structure on the additive category $\mathrm{Ho}(A-Mod)$,
the composition of endomorphisms endows $\pi_{0}(A)$ with a multiplicative structure making
it into a commutative ring. More generally, the category
$\mathrm{Ho}(A-Mod)$ has a natural structure of a graded category (i.e. has a natural
enrichment into the symmetric monoidal category of $\mathbb{N}$-graded
abelian groups), defined by
$$[M,N]_{*}:=\oplus_{i}[S^{i}(N),M]\simeq \oplus_{i}[N,\Omega^{i}(M)],$$
where $S^{i}$ is the $i$-fold iterated suspension functor, and
$\Omega^{i}$ the $i$-fold iterated loop functor.
Therefore, the graded endomorphism ring of the unit
$A$ has a natural structure of a graded commutative ring. As this
endomorphism ring is naturally isomorphic to
$$\pi_{*}(A):=\oplus_{i}\pi_{i}(A),$$
we obtain this way a natural structure of a graded
commutative ring on $\pi_{*}(A)$. This clearly defines a functor
$$\pi_{*} : \mathrm{Ho}(Comm(\mathcal{C})) \longrightarrow GComm,$$
from $\mathrm{Ho}(Comm(\mathcal{C}))$ to the category of $\mathbb{N}$-graded
commutative rings.

In the same way, for a commutative monoid $A\in Comm(\mathcal{C})$, and
a $A$-module $M$, one can define
$$\pi_{*}(M):=\pi_{*}(|M|)=\pi_{*}(Map_{A-Mod}(A,M)),$$
which has a natural structure of a graded $\pi_{*}(A)$-module.
This defines a functor
$$\pi_{*} : \mathrm{Ho}(A-Mod) \longrightarrow \pi_{*}(A)-GMod,$$
from $\mathrm{Ho}(A-Mod)$ to the category of $\mathbb{N}$-graded
$\pi_{*}(A)$-modules.

\begin{prop}\label{ploc}
Let $A\in Comm(\mathcal{C})$ be a commutative monoid in $\mathcal{C}$, and
$a\in \pi_{0}(A)$. There exists an epimorphism
$A\longrightarrow A[a^{-1}]$, such that for any commutative $A$-algebra $C$,
the simplicial set
$Map_{A-Comm(\mathcal{C})}(A[a^{-1}],C)$
is non-empty (and thus contractible) if and only if
the image of $a$ in $\pi_{0}(C)$ by the morphism
$\pi_{0}(A) \rightarrow \pi_{0}(C)$ is an invertible element.
\end{prop}

\begin{proof} We represent the element $a$ as a morphism in $\mathrm{Ho}(A-Mod)$
$a : A \longrightarrow A$. Taking the image of $a$ by the
left derived functor of the free commutative $A$-algebra functor
(which is left Quillen)
$$F_{A} : A-Mod \longrightarrow A-Comm(\mathcal{C})$$
we find a morphism in $\mathrm{Ho}(A-Comm(\mathcal{C}))$
$$a : \mathbb{L}F_{A}(A) \longrightarrow \mathbb{L}F_{A}(A).$$
As the model category $A-Comm(\mathcal{C})$ is a $\mathbb{U}$-combinatorial model
category, one can apply the localization techniques of
in order to invert any $\mathbb{U}$-small set of
morphisms in $A-Comm(\mathcal{C})$ (see \cite{sm,du2}). We
let $L_{a}A-Comm(\mathcal{C})$ be the left Bousfield localization
of $A-Comm(\mathcal{C})$ along the set of morphisms (with one element)
$$S_{a}:=\{a : \mathbb{L}F_{A}(A) \longrightarrow \mathbb{L}F_{A}(A)\}.$$ We define
$A[a^{-1}] \in \mathrm{Ho}(A-Comm(\mathcal{C}))$ as a local model
of $A$ in $L_{a}A-Comm(\mathcal{C})$, the left Bousfield localization of
$A-Comm(\mathcal{C})$ along $S_{a}$.

First of all, the $S_{a}$-local objects are the commutative $A$-algebras $B$ such that
the induced morphism
$$a^{*} : Map_{A-Mod}(A,B) \longrightarrow Map_{A-Mod}(A,B)$$
is an equivalence. Equivalently, the multiplication by $a\in \pi_{0}(A)$
$$\times a : \pi_{*}(B) \longrightarrow \pi_{*}(B)$$
is an isomorphism. This shows that the $S_{a}$-local objects
are the commutative $A$-algebras $A
\longrightarrow B$ such that the image of $a$ by $\pi_{0}(A)
\longrightarrow \pi_{0}(B)$ is invertible.

Suppose now that $C\in A-Comm(\mathcal{C})$ is such that
$Map_{A-Comm(\mathcal{C})}(A[a^{-1}],C)$ is not empty. The morphism
$\pi_{0}(A) \longrightarrow \pi_{0}(C)$ then factors through
$\pi_{0}(A[a^{-1}])$, and thus the image of $a$ is invertible in
$\pi_{0}(C)$. Therefore, $C$ is a $S_{a}$-local object, and thus
$$Map_{A-Comm(\mathcal{C})}(A[a^{-1}],C)\simeq Map_{A-Comm(\mathcal{C})}(A,C)\simeq *.$$
This implies that $A \longrightarrow A[a^{-1}]$ is an epimorphism.
It only remain to prove that if $C\in A-Comm(\mathcal{C})$ is such that the
image of $a$ is invertible in $\pi_{0}(C)$, then
$Map_{A-Comm(\mathcal{C})}(A[a^{-1}],C)$ is non-empty. But such a $C$ is a
$S_{a}$-local object, and therefore
$$*\simeq Map_{A-Comm(\mathcal{C})}(A,C)\simeq Map_{A-Comm(\mathcal{C})}(A[a^{-1}],C).$$
\end{proof}

\begin{df}\label{dloc}
Let $A\in Comm(\mathcal{C})$ and $a\in \pi_{0}(A)$. The commutative
$A$-algebra $A[a^{-1}]$\index{$A[a^{-1}]$} is called
\emph{the standard localization of $A$ with respect to $a$}\index{standard localization}.
\end{df}

A useful property of standard localizations is given by the
following corollary of  the proof of Prop. \ref{ploc}. In order to state it,
we will use the following notations. For any
$A\in Comm(\mathcal{C})$ and $a\in \pi_{0}(A)$, we represent
$a$ as a morphism in $\mathrm{Ho}(A-Mod)$,
$A \longrightarrow A$. Tensoring this morphism with $M$ gives
a morphism in $\mathrm{Ho}(A-Mod)$, denoted by
$$\times a : M \longrightarrow M.$$

\begin{cor}\label{cloc}
Let $A\in Comm(\mathcal{C})$, $a\in \pi_{0}(A)$ and
let $f : A \longrightarrow A[a^{-1}]$ be as in
Prop. \ref{ploc}.
\begin{enumerate}
\item The functor
$$\mathrm{Ho}(A[a^{-1}]-Comm(\mathcal{C})) \longrightarrow \mathrm{Ho}(A-Comm(\mathcal{C}))$$
is fully faithful, and its essential image consists of
all commutative $B$-algebras such that the image of
$a$ is invertible in
$\pi_{0}(B)$.
\item
The functor
$$f_{*} : \mathrm{Ho}(A[a^{-1}]-Mod) \longrightarrow \mathrm{Ho}(A-Mod)$$
is fully faithful and its essential image
consists of all $A$-modules $M$ such that
the multiplication by $a$
$$\times a : M \longrightarrow M$$
is an isomorphism in $\mathrm{Ho}(A-Mod)$.
\end{enumerate}
\end{cor}

\begin{proof} $(1)$ The fact that the functor $f_{*}$ is fully 
faithful is immediate
as $f$ is an epimorphism. The fact that the functor $f_{*}$ 
takes its values in the
required subcategory is clear by functoriality of the construction $\pi_{*}$.

Let $B$ be a commutative $A$-algebra such that the image of $a$ is invertible
in $\pi_{0}(B)$. Then, we know that $B$ is a $S_{a}$-local object. Therefore, one has
$$Map_{A-Comm(\mathcal{C})}(A[a^{-1}],B)\simeq Map_{A-Comm(\mathcal{C})}(A,B)\simeq *,$$
showing that $B$ is in the image of $f_{*}$. \\

$(2)$ The fact that $f_{*}$ is fully faithful follows from
Prop. \ref{p5} $(1)$. Let $M\in A[a^{-1}]-Mod$ and let us prove that
the morphism $\times a : M \longrightarrow M$ is an isomorphism
in $\mathrm{Ho}(A-Mod)$. Using that $A \longrightarrow A[a^{-1}]$ is an epimorphism, one finds
$M\simeq M\otimes^{\mathbb{L}}_{A}A[a^{-1}]$, which reduces the problem to the case
where $M=A[a^{-1}]$. But then, the morphism
$\times a : A[a^{-1}] \longrightarrow A[a^{-1}]$, as a morphism in $\mathrm{Ho}(A[a^{-1}]-Mod)$ lives in
$[A[a^{-1}],A[a^{-1}]]\simeq \pi_{0}(A[a^{-1}])$, and correspond to the image of
$a$ by the morphism $\pi_{0}(A) \longrightarrow \pi_{0}(A[a^{-1}])$, which is then
invertible. In other words, $\times a : A[a^{-1}] \longrightarrow A[a^{-1}]$ is
an isomorphism.

Conversely, let $M$ be an $A$-module such that the morphism
$\times a : M \longrightarrow M$ is an isomorphism in $\mathrm{Ho}(A-Mod)$.
We need to show that
the adjunction morphism
$$M\longrightarrow M\otimes^{\mathbb{L}}_{A}A[a^{-1}]$$
is an isomorphism. For this, we use that the morphism
$A\longrightarrow A[a^{-1}]$ can be constructed using a small object
argument with respect to the horns over the morphism
$a : \mathbb{L}F_{A}(A) \longrightarrow \mathbb{L}F_{A}(A)$
(see \cite[4.2]{hi}). Therefore, a transfinite induction argument shows that it is
enough to prove that the  morphism induced by tensoring
$$a\otimes Id : \mathbb{L}F_{A}(A)\otimes^{\mathbb{L}}_{A}M\simeq \mathbb{L}F_{A}(M) \longrightarrow
\mathbb{L}F_{A}(A)\otimes^{\mathbb{L}}_{A}M\simeq \mathbb{L}F_{A}(M)$$
is an isomorphism in $\mathrm{Ho}(A-Mod)$. But this morphism
is the image by
the functor $\mathbb{L}F_{A}$ of the morphism $\times a : M \longrightarrow M$, and is
therefore an isomorphism. \end{proof}

\begin{prop}\label{ploc2}
Let $A \in Comm(\mathcal{C})$, $a\in \pi_{0}(A)$ and
$A\longrightarrow A[a^{-1}]$ the standard localization with respect to $a$.
Assume that the model category $\mathcal{C}$ is finitely generated (in the sense
of \cite{ho}).
\begin{enumerate}
\item The morphism $A \longrightarrow A[a^{-1}]$ is a
formal Zariski open immersion.

\item If $1$ is a finitely presented object in $\mathcal{C}$, then
$A \longrightarrow A[a^{-1}]$ is a Zariski open immersion.

\end{enumerate}
\end{prop}

\begin{proof} $(1)$ It only remains to show that the morphism
$A\longrightarrow A[a^{-1}]$ is flat. In other words, we need to prove that the
functor $M \mapsto M\otimes^{\mathbb{L}}_{A}A[a^{-1}]$, preserves homotopy fiber
sequences of $A$-modules.

The model category $A-Mod$ is $\mathbb{U}$-combinatorial and
finitely generated. In particular, there exists a
$\mathbb{U}$-small set $G$ of $\omega$-small cofibrant objects in $A-Mod$, such that
a morphism $N \longrightarrow P$ is an equivalence in $A-Mod$ if and only if
for any $X\in G$ the induced morphism $Map_{A-Mod}(X,N)\longrightarrow
Map_{A-Mod}(X,P)$ is an isomorphism in $\mathrm{Ho}(SSet)$.
Furthermore, filtered homotopy colimits preserve
homotopy fiber sequences.
For any $A$-module $M$, we let
$M_{a}$ be the transfinite homotopy colimit
$$\xymatrix{
M\ar[r]^-{\times a} & M \ar[r]^-{\times a} \ar[r] & \dots & }$$
where the morphism $\times a$ is composed with itself $\omega$-times.

The functor
$M \mapsto M_{a}$ commutes with homotopy fiber sequences. Therefore, it only remain to show that
$M_{a}$ is naturally isomorphic in $\mathrm{Ho}(A-Mod)$ to $M\otimes^{\mathbb{L}}_{A}A[a^{-1}]$.
For this, it is enough to check that the natural morphism
$M \longrightarrow M_{a}$ induces an isomorphism in $\mathrm{Ho}(A-Mod)$
$$M\otimes^{\mathbb{L}}_{A}A[a^{-1}] \simeq (M_{a})\otimes^{\mathbb{L}}_{A}A[a^{-1}],$$
and that the natural morphism
$$(M_{a})\otimes^{\mathbb{L}}_{A}A[a^{-1}] \longrightarrow M_{a}$$
is an isomorphism in $\mathrm{Ho}(A-Mod)$. The first assumption follows easily from the fact
that $-\otimes^{\mathbb{L}}_{A}A[a^{-1}]$ commutes with homotopy colimits, and the fact that
$\times a :  A[a^{-1}] \longrightarrow A[a^{-1}]$ is an isomorphism in $\mathrm{Ho}(A-Mod)$.
For the second assumption we use Cor. \ref{cloc} $(2)$, which tells us that
it is enough to check that the morphism
$$\times a : M_{a} \longrightarrow M_{a}$$
is an isomorphism in $\mathrm{Ho}(A-Mod)$. For this, we need to show that for any $X\in G$
the induced morphism
$$Map_{A-Mod}(X,M_{a})\longrightarrow
Map_{A-Mod}(X,M_{a})$$
is an isomorphism in $\mathrm{Ho}(SSet)$.
But, as the objects in $G$ are cofibrant and $\omega$-small,
$Map_{A-Mod}(X,-)$ commutes with $\omega$-filtered homotopy colimits, and the
morphism
$$Map_{A-Mod}(X,M_{a})\longrightarrow
Map_{A-Mod}(X,M_{a})$$
is then obviously an isomorphism in $\mathrm{Ho}(SSet)$ by the construction of $M_{a}$. \\

$(2)$ When $1$ is finitely presented, one has for any filtered diagram
of commutative $A$-algebras $B_{i}$ an isomorphism
$$Colim_{i}\pi_{*}(B_{i})\simeq \pi_{*}(Hocolim_{i}B_{i}).$$
Using that $Map_{A-Comm(\mathcal{C})}(A[a^{-1}],B)$ is either empty or contractible, depending
whether of not $a$ goes to a unit in $\pi_{*}(B)$, we easily deduce that
$$Hocolim_{i}Map_{A-Comm(\mathcal{C})}(A[a^{-1}],B_{i})\simeq Map_{A-Comm(\mathcal{C})}(A[a^{-1}],Hocolim_{i}B_{i}).$$
\end{proof}

We can also show that the natural morphism $A \longrightarrow A[a^{-1}]$
is a formally \'etale morphism in the sense of Def. \ref{d5}.

\begin{prop}\label{ploc3}
Let $A\in Comm(\mathcal{C})$ and $a\in \pi_{0}(A)$. Then, the natural morphism
$A \longrightarrow A[a^{-1}]$ is  formally \'etale.
\end{prop}

\begin{proof} Let $M$ be any $A[a^{-1}]$-module. We need to show that the
natural morphism
$$Map_{A-Comm(\mathcal{C})/A[a^{-1}]}(A[a^{-1}],A[a^{-1}]\oplus M) \longrightarrow
Map_{A-Comm(\mathcal{C})/A[a^{-1}]}(A,A[a^{-1}]\oplus M)$$
is an isomorphism in $\mathrm{Ho}(SSet)$. Using the universal property of
$A\longrightarrow A[a^{-1}]$ given by Prop. \ref{ploc} we see that it is enough to
prove that for any $B\in Comm(\mathcal{C})$, and any $B$-module $M$ the natural
projection $\pi_{0}(B\oplus M) \longrightarrow \pi_{0}(B)$ reflects invertible elements
(i.e. an element in $\pi_{0}(B\oplus M)$ is invertible if and only
its image in $\pi_{0}(B)$ is so). But, clearly,
$\pi_{0}(B\oplus M)$ can be identified with the trivial square zero extension
of the commutative ring $\pi_{0}(B)$ by $\pi_{0}(M)$, which implies the required result. 
\end{proof}

\begin{cor}\label{cloc2}
Assume that the model category $\mathcal{C}$ is finitely presented, and that
the unit $\mathbf{1}$ is finitely presented in $\mathcal{C}$. Then for any
$A\in Comm(\mathcal{C})$ and $a\in \pi_{0}(A)$, the
morphism $A\longrightarrow A[a^{-1}]$ is an \'etale, flat epimorphism.
\end{cor}

\begin{proof} Put \ref{ploc2} and \ref{ploc3} together. \end{proof}

\section{Zariski open immersions and perfect modules}

Let $A$ be a commutative monoid in $\mathcal{C}$ and $K$ be
a perfect $A$-module in the sense of Def. \ref{d4}. We are going to
define a Zariski open immersion $A \longrightarrow A_{K}$,
which has to be thought as the complement of
the support of the $A$-module $K$.

\begin{prop}\label{psupp}
Assume that $\mathcal{C}$ is stable model category. Then there exists
a formal Zariski open immersion $A \longrightarrow A_{K}$,
such that for any commutative $A$-algebra $C$,
the simplicial set $$Map_{A-Comm(\mathcal{C})}(A_{K},C)$$
is non-empty (and thus contractible) if and only if
$K\otimes_{A}^{\mathbb{L}}C\simeq *$ in $\mathrm{Ho}(C-Mod)$.
If the unit $\mathbf{1}$ is furthermore finitely presented, then
$A \longrightarrow A_{K}$ is finitely presented and thus
is a Zariski open immersion.
\end{prop}

\begin{proof} The commutative $A$-algebra is constructed
using a left Bousfield localization of the model
category $A-Comm(\mathcal{C})$, as done in the proof of
Prop. \ref{ploc}.

We let $I$ be a generating $\mathbb{U}$-small set of
cofibrations in $A-Comm(\mathcal{C})$, and
$K^{\vee}:=\mathbb{R}\underline{Hom}_{A-Mod}(K,A)$
be the dual of $K$ in $\mathrm{Ho}(A-Mod)$. For any morphism
$X \longrightarrow Y$ in $I$, we consider
the morphism of free commutative $A$-algebras
$$\mathbb{L}F_{A}(K^{\vee}\otimes^{\mathbb{L}}_{A}X) \longrightarrow
\mathbb{L}F_{A}(K^{\vee}\otimes^{\mathbb{L}}_{A}Y),$$
where $F_{A} : A-Mod \longrightarrow A-Comm(\mathcal{C})$ is the
left Quillen functor sending an $A$-module to the
free commutative $A$-algebra it generates.
When $X \longrightarrow Y$ varies in $I$ this gives
a $\mathbb{U}$-small set of morphisms
denoted by $S_{K}$ in $A-Comm(\mathcal{C})$.
We consider $L_{K}A-Comm(\mathcal{C})$, the left Bousfield
localization of $A-Comm(\mathcal{C})$ along the set $S_{K}$.
By definition, $A\longrightarrow A_{K}$ is
an $S_{K}$-local model of $A$ in the localized model
category $L_{K}A-Comm(\mathcal{C})$.

\begin{lem}\label{l2}
The $S_{K}$-local objects in $L_{K}A-Comm(\mathcal{C})$
are precisely the commutative $A$-algebras $B$ such that
$K\otimes_{A}^{\mathbb{L}}B\simeq *$ in $\mathrm{Ho}(A-Mod)$.
\end{lem}

\begin{proof} First of all, one has an adjunction
isomorphism in $\mathrm{Ho}(SSet)$
$$Map_{A-Comm(\mathcal{C})}(\mathbb{L}F_{A}(K^{\vee}\otimes^{\mathbb{L}}_{A}X),B)
\simeq Map_{A-Mod}(X,K\otimes^{\mathbb{L}}_{A}B).$$
This implies that an object $B \in A-Comm(\mathcal{C})$ is
$S_{K}$-local if and only if for all morphism
$X \longrightarrow Y$ in $I$ the induced morphism
$$Map_{A-Mod}(Y,K\otimes^{\mathbb{L}}_{A}B)\longrightarrow
Map_{A-Mod}(X,K\otimes^{\mathbb{L}}_{A}B)$$
is an isomorphism in $\mathrm{Ho}(SSet)$. As $I$ is a set of
generating cofibrations in $A-Mod$ this implies that
$B$ is $S_{K}$-local if and only if
$K\otimes^{\mathbb{L}}_{A}C\simeq *$ in $\mathrm{Ho}(A-Mod)$. \end{proof}

We now finish the proof of proposition \ref{psupp}. First of all,
Lem. \ref{l2} implies that
for any commutative $A$-algebra $B$, if the mapping space
$Map_{A-Comm(\mathcal{C})}(A_{K},B)$ is non-empty then
$$B\otimes_{A}^{\mathbb{L}}K\simeq B\otimes_{A_{K}}^{\mathbb{L}}(A_{K}\otimes_{A}^{\mathbb{L}}K)
\simeq *,$$ and thus
$B$ is an $S_{K}$-local object. This shows that if
$Map_{A-Comm(\mathcal{C})}(A_{K},B)$ is non-empty then one has
$$Map_{A-Comm(\mathcal{C})}(A_{K},B)\simeq Map_{A-Comm(\mathcal{C})}(A,B)\simeq *.$$
In other words $A \longrightarrow A_{K}$ is a formal
Zariski open immersion by
Cor. \ref{c1} and the stability assumption on
$\mathcal{C}$. It only remains to prove that
$A\longrightarrow A_{K}$ is also finitely presented when $\mathbf{1}$ is
a finitely presented object in $\mathcal{C}$.

For this, let $\{C_{i}\}_{i\in I}$ be a filtered
diagram of commutative $A$-algebras and $C$ be its homotopy colimit.
By the property of $A \longrightarrow A_{K}$, we need
to show that if $K\otimes^{\mathbb{L}}_{A}C\simeq *$
then there is an $i\in I$ such that
$K\otimes^{\mathbb{L}}_{A}C_{i}\simeq *$. For this, we consider
the two elements $Id$ and $*$ in
$[K\otimes^{\mathbb{L}}_{A}C,K\otimes^{\mathbb{L}}_{A}C]_{C-Mod}$.
As $K$ is perfect and $\mathbf{1}$ is a finitely presented object,
$K$ is a finitely presented $A$-module by Prop. \ref{p4} $(4)$.
Therefore, one has
$$*\simeq [K\otimes^{\mathbb{L}}_{A}C,K\otimes^{\mathbb{L}}_{A}C]\simeq
Colim_{i\in I}[K,K\otimes^{\mathbb{L}}_{A}C_{i}]_{A-Mod}.$$
As the two elements $Id$ and $*$ becomes equal
in the colimit, there is an $i$ such that
they are equal as elements in
$$[K,K\otimes^{\mathbb{L}}_{A}C_{i}]_{A-Mod}\simeq
[K\otimes^{\mathbb{L}}_{A}C_{i},K\otimes^{\mathbb{L}}_{A}C_{i}]_{C_{i}-Mod},$$
showing that $K\otimes^{\mathbb{L}}_{A}C_{i}\simeq *$ in
$\mathrm{Ho}(C_{i}-Mod)$.  \end{proof}

\begin{cor}\label{csupp}
Assume that $\mathcal{C}$ is a stable model category, and
let $f : A \longrightarrow A_{K}$ be as in
Prop. \ref{psupp}.
\begin{enumerate}
\item The functor
$$\mathrm{Ho}(A_{K}-Comm(\mathcal{C})) \longrightarrow \mathrm{Ho}(A-Comm(\mathcal{C}))$$
is fully faithful, and its essential image consists of
all commutative $B$-algebras such that
$B\otimes_{A}^{\mathbb{L}}K\simeq *$.
\item
The functor
$$f_{*} : \mathrm{Ho}(A_{K}-Mod) \longrightarrow \mathrm{Ho}(A-Mod)$$
is fully faithful and its essential image
consists of all $A$-modules $M$ such that
$M\otimes_{A}^{\mathbb{L}}K\simeq *$.
\end{enumerate}
\end{cor}

\begin{proof} As the morphism
$A\longrightarrow A_{K}$ is an epimorphism,
we know by Prop. \ref{p5} that the functors
$$\mathrm{Ho}(A_{K}-Comm(\mathcal{C})) \longrightarrow \mathrm{Ho}(A-Comm(\mathcal{C})) \qquad \mathrm{Ho}(A_{K}-Mod) \longrightarrow \mathrm{Ho}(A-Mod)$$
are both fully faithful.
Furthermore, a commutative $A$-algebra $B$ is in the essential image
of the first one if and only if $Map_{A-Comm(\mathcal{C})}(A_{K},B)\simeq *$, and
therefore if and only if $B\otimes_{A}^{\mathbb{L}}K\simeq *$. For the second functor, its clear
that if $M$ is a $A_{K}$-module, then
$$M\otimes_{A}^{\mathbb{L}}K\simeq
M\otimes_{A_{K}}^{\mathbb{L}}A_{K}\otimes_{A}^{\mathbb{L}}K\simeq *.$$
Conversely, let $M$ be an $A$-module such that $M\otimes_{A}^{\mathbb{L}}K\simeq *$.

\begin{lem}\label{l3}
For any $A$-module $M$,
$K\otimes^{\mathbb{L}}_{A}M\simeq *$ if and only if
$K^{\vee}\otimes^{\mathbb{L}}_{A}M\simeq *$.
\end{lem}

\begin{proof} Indeed, as $K$ is perfect,
$K^{\vee}$ is a retract of  $K^{\vee}\otimes_{A}^{\mathbb{L}}K\otimes_{A}^{\mathbb{L}}K^{\vee}$
in $\mathrm{Ho}(A-Mod)$. This implies that $K^{\vee}\otimes^{\mathbb{L}}_{A}M$
is a retract of
$K^{\vee}\otimes_{A}^{\mathbb{L}}K\otimes_{A}^{\mathbb{L}}K^{\vee}
\otimes_{A}^{\mathbb{L}}M$, showing that
$$\left(K\otimes^{\mathbb{L}}_{A}M\simeq * \right) \Rightarrow
\left(K^{\vee}\otimes^{\mathbb{L}}_{A}M\simeq * \right).$$
By symmetry this proves the lemma. \end{proof}

We need to prove that the adjunction morphism
$$M\longrightarrow M\otimes_{A}^{\mathbb{L}}A_{K}$$
is an isomorphism in $\mathrm{Ho}(A-Mod)$. For this, we use the fact that
the morphism $A\longrightarrow A_{K}$ can be constructed using a
small object argument on the set of horns on the set $S_{K}$ (see
\cite[4.2]{hi}). By a transfinite induction we are therefore reduced to show that
for any morphism
$X \longrightarrow Y$ in $\mathcal{C}$,
the natural morphism
$$\mathbb{L}F_{A}(K^{\vee}\otimes^{\mathbb{L}}_{A}X)\otimes_{A}^{\mathbb{L}}M \longrightarrow
\mathbb{L}F_{A}(K^{\vee}\otimes^{\mathbb{L}}_{A}Y)\otimes_{A}^{\mathbb{L}}M$$
is an isomorphism in $\mathrm{Ho}(\mathcal{C})$. But, using that
$M\otimes_{A}^{\mathbb{L}}K\simeq *$ and lemma \ref{l3}, this is clear by the explicit description
of the functor $F_{A}$. \end{proof}

\section{Stable modules}

We recall
that as $\mathcal{C}$ is a pointed model category
one can define a suspension functor (see \cite[\S 7]{ho})
$$\begin{array}{cccc}
S : & \mathrm{Ho}(\mathcal{C}) & \longrightarrow & \mathrm{Ho}(\mathcal{C}) \\
 & X & \mapsto & S(X):=*\coprod^{\mathbb{L}}_{X}*.
\end{array}$$
This functor possesses a right adjoint, the loop
functor
$$\begin{array}{cccc}
\Omega : & \mathrm{Ho}(\mathcal{C}) & \longrightarrow & \mathrm{Ho}(\mathcal{C}) \\
 & X & \mapsto & \Omega(X):=*\times^{h}_{X}*.
\end{array}$$

We fix, once for all an object $S^{1}_{\mathcal{C}}\in \mathcal{C}$, which
is a cofibrant model for $S(\mathbf{1})\in \mathrm{Ho}(\mathcal{C})$.
For any commutative monoid $A\in Comm(\mathcal{C})$, we let
$$S^{1}_{A}:=S^{1}_{\mathcal{C}}\otimes A \in A-Mod$$
be the free $A$-module on $S^{1}_{\mathcal{C}}$. It is a cofibrant
object in $A-Mod$, which is a model for the suspension $S(A)$ (note
that $S^{1}_{A}$ is cofibrant in $A-Mod$, but not in $\mathcal{C}$  
unless $A$ is itself cofibrant in $\mathcal{C}$).
The functor
$$S^{1}_{A}\otimes_{A} -  :
A-Mod \longrightarrow A-Mod$$
has a right adjoint
$$\underline{Hom}_{A}(S^{1}_{A},-) : A-Mod \longrightarrow A-Mod.$$
Furthermore, assumption \ref{ass1} implies that
$S^{1}_{A}\otimes_{A} - $ is a left Quillen functor. We can therefore
apply the general construction of \cite{ho2}
in order to produce a model category $Sp^{S_{A}^{1}}(A-Mod)$, of
spectra in $A-Mod$ with respect to the
left Quillen endofunctor $S^{1}_{A}\otimes_{A} - $.

\begin{df}\label{stabmod}
Let $A\in Comm(\mathcal{C})$ be a commutative monoid in $\mathcal{C}$.
The model category of \emph{stable $A$-modules}\index{stable modules}
is the model category $Sp^{S^{1}_{A}}(A-Mod)$\index{$Sp^{S^{1}_{A}}(A-Mod)$}, of
spectra in $A-Mod$ with respect to the left
Quillen endo-functor
$$S_{A}^{1}\otimes_{A} - : A-Mod \longrightarrow A-Mod.$$
It will simply be denoted
by $Sp(A-Mod)$, and its objects will be called
\emph{stable $A$-modules}.
\end{df}

Recall that objects in the category
$Sp(A-Mod)$ are families of objects
$M_{n} \in A-Mod$ for $n\geq 0$, together with
morphisms $\sigma_{n} : S^{1}_{A}\otimes_{A}M_{n} \longrightarrow M_{n+1}$.
Morphisms in $Sp(A-Mod)$ are simply families of morphisms
$f_{n} : M_{n} \rightarrow N_{n}$ commuting
with the morphisms $\sigma_{n}$. One starts by endowing 
$Sp(A-Mod)$ with the levelwise model structure, for which
equivalences (resp. fibrations) are the morphisms
$f : M_{*} \longrightarrow N_{*}$ such that each morphism
$f_{n} : M_{n} \longrightarrow N_{n}$
is an equivalence in $A-Mod$ (resp. a fibration). The
definitive model structure, called the stable model structure,
is the left Bousfield localization of $Sp(A-Mod)$
whose local objects are the stable modules $M_{*} \in Sp(A-Mod)$ such that
each induced morphism
$$M_{n} \longrightarrow \mathbb{R}\underline{Hom}_{A}(
S^{1}_{A},M_{n+1})$$
is an isomorphism in $\mathrm{Ho}(A-Mod)$. These local objects
will be called \emph{$\Omega$-stable $A$-modules}.
We refer to
\cite{ho2} for details concerning the existence
and the properties of this model structure.

There exists an adjunction
$$S_{A} : A-Mod \longrightarrow Sp(A-Mod) \qquad
A-Mod \longleftarrow Sp(A-Mod) : (-)_{0},$$
where the right adjoint sends a stable $A$-module $M_{*}$
to $M_{0}$. The left adjoint
is defined by
$S_{A}(M)_{n}:=(S^{1}_{A})^{\otimes_{A} n}\otimes_{A}M$,
with the natural transition morphisms. This adjunction is
a Quillen adjunction, and can be derived into an
adjunction on the level of homotopy categories
$$\mathbb{L}S_{A}\simeq S_{A} : \mathrm{Ho}(A-Mod) \longrightarrow \mathrm{Ho}(Sp(A-Mod))$$
$$\mathrm{Ho}(A-Mod) \longleftarrow \mathrm{Ho}(Sp(A-Mod)) : \mathbb{R}(-)_{0}.$$
Note that by \ref{ass1}, $S_{A}$ preserves equivalences,
so $\mathbb{L}S_{A}\simeq S_{A}$. On the contrary, the functor
$(-)_{0}$ does not preserve equivalences and must be derived
on the right. In particular, the functor
$S_{A} : \mathrm{Ho}(A-Mod) \longrightarrow \mathrm{Ho}(Sp(A-Mod))$
is not fully faithful in general.

\begin{lem}\label{lstabmod}
\begin{enumerate}
\item
Assume that the suspension functor
$$S : \mathrm{Ho}(\mathcal{C}) \longrightarrow \mathrm{Ho}(\mathcal{C})$$
is fully faithful. Then, for any commutative
monoid $A\in Comm(\mathcal{C})$, the functor
$$S_{A} : \mathrm{Ho}(A-Mod) \longrightarrow \mathrm{Ho}(Sp(A-Mod))$$
is fully faithful.
\item If furthermore, $\mathcal{C}$ is a stable
model category then $S_{A}$ is a Quillen equivalence.
\end{enumerate}
\end{lem}

\begin{proof} $(1)$ As the adjunction morphism
$M \longrightarrow S_{A}(M)_{0}$ is always an
isomorphism in $A-Mod$, it is enough to show that
for any $M\in A-Mod$ the stable $A$-module
$S_{A}(M)$ is a $\Omega$-stable $A$-module.
For this, it is enough to show that for any
$M\in A-Mod$, the adjunction morphism
$$M \longrightarrow
\mathbb{R}\underline{Hom}_{A}(
S^{1}_{A},S^{1}_{A}\otimes_{A}M)$$
is an isomorphism in $\mathrm{Ho}(\mathcal{C})$. But, one has natural
isomorphisms in $\mathrm{Ho}(\mathcal{C})$
$$S^{1}_{A}\otimes_{A}M\simeq S^{1}_{\mathcal{C}}\otimes M\simeq S(M)$$
$$\mathbb{R}\underline{Hom}_{A}(
S^{1}_{A},S^{1}_{A}\otimes_{A}M)\simeq
\mathbb{R}\underline{Hom}_{1}(
S^{1},S(M))\simeq \Omega(S(M)).$$
This shows that the above morphism is in fact isomorphic
in $\mathrm{Ho}(\mathcal{C})$ to the adjunction morphism
$$M \longrightarrow \Omega(S(M))$$
which is an isomorphism by hypothesis on $\mathcal{C}$. \\

$(2)$ As $\mathcal{C}$ is a stable model category, the functor
$S^{1}_{\mathcal{C}}\otimes - : \mathcal{C} \longrightarrow \mathcal{C}$ is a
Quillen equivalence. This also implies that for
any $A\in Comm(\mathcal{C})$, the functor
$S_{A}^{1} : A-Mod \longrightarrow A-Mod$
is a Quillen equivalence. We know by \cite{ho2}
that $S_{A} : A-Mod \longrightarrow Sp(A-Mod)$
is a Quillen equivalence. \end{proof}

The Quillen adjunction
$$S_{A} : A-Mod \longrightarrow Sp(A-Mod) \qquad
A-Mod \longleftarrow Sp(A-Mod) : (-)_{0},$$
is furthermore functorial in $A$. Indeed, for
$A \longrightarrow B$ a morphism in $Comm(\mathcal{C})$, one defines
a functor
$$-\otimes_{A}B : Sp(A-Mod) \longrightarrow Sp(B-Mod)$$
defined by
$$(M_{*}\otimes_{A}B)_{n}:=M_{n}\otimes_{A}B.$$
The transitions morphisms are given by
$$S^{1}_{B}\otimes_{B}(M_{n}\otimes_{A}B) \simeq
(S^{1}_{A}\otimes_{A}M_{n})\otimes_{A}B \longrightarrow
M_{n+1}\otimes_{A}B.$$
Clearly, the
square of left Quillen functors
$$\xymatrix{
A-Mod \ar[r]^-{S_{A}} \ar[d]_-{-\otimes_{A}B} & Sp(A-Mod) \ar[d]^-{-\otimes_{A}B} \\
 B-Mod \ar[r]_{S_{B}} & Sp(B-Mod)}$$
commutes up to a natural isomorphism. So does the square of
right Quillen functors
$$\xymatrix{
B-Mod \ar[r] \ar[d] & Sp(B-Mod) \ar[d] \\
A-Mod \ar[r] & Sp(A-Mod).}$$
Finally, using techniques of symmetric spectra, as
done in \cite{ho2}, it is possible to show that the homotopy category of
stable $A$-modules inherits from $\mathrm{Ho}(A-Mod)$ a symmetric monoidal structure,
still denoted by $-\otimes^{\mathbb{L}}_{A}-$. This makes the homotopy
category $\mathrm{Ho}(Sp(A-Mod))$ into a closed symmetric monoidal category.
In particular, for two stable $A$-modules $M_{*}$ and $N_{*}$ one can define
a stable $A$-modules of morphisms
$$\mathbb{R}\underline{Hom}_{A}^{Sp}(M_{*},N_{*}) \in
\mathrm{Ho}(Sp(A-Mod)).$$

We now consider the
category $(A-Mod_{0}^{op})^{\wedge}$ of pre-stacks over
$A-Mod_{0}^{op}$, as defined in
\cite{hagI}. Recall it is
the category of $\mathbb{V}$-simplicial presheaves
on $A-Mod_{0}^{op}$,
and that its model structure is obtained
from the projective levelwise model structure
by a left Bousfield localization inverting the equivalences
in $A-Mod$. The homotopy category
$\mathrm{Ho}((A-Mod_{0}^{op})^{\wedge})$ can be naturally identified with the
full subcategory of $\mathrm{Ho}(SPr(A-Mod_{0}^{op}))$ consisting of
functors
$$F : A-Mod_{0} \longrightarrow SSet_{\mathbb{V}},$$
sending equivalences of $A$-modules to equivalences of simplicial sets.

We  define a functor
$$\begin{array}{cccc}
\underline{h}_{s}^{-} : & Sp(A-Mod)^{op} & \longrightarrow &
(A-Mod_{0}^{op})^{\wedge} \\
 & M_{*} & \mapsto & \underline{h}_{s}^{M_{*}}
\end{array}$$
by
$$\begin{array}{cccc}
\underline{h}_{s}^{M_{*}} : & A-Mod_{0} & \longrightarrow & SSet_{\mathbb{V}} \\
 & N & \mapsto & Hom(M_{*},\Gamma_{*}(S_{A}(N))),
\end{array}$$
where $\Gamma_{*}$ is a simplicial resolution functor
on the model category $Sp(A-Mod)$.

Finally we will need some terminology. A stable $A$-module $M_{*}\in \mathrm{Ho}(Sp(A-Mod))$ is
called \emph{$0$-connective}, if it is isomorphic to some
$S_{A}(M)$ for an $A$-module $M\in \mathrm{Ho}(A-Mod_{0})$. By induction,
for an integer $n>0$, a stable
$A$-module $M_{*}\in \mathrm{Ho}(Sp(A-Mod))$ is
called \emph{$(-n)$-connective}, if it is isomorphic to
$\Omega(M_{*}')$ for some $-(n-1)$-connective
stable $A$-module $M_{*}'$ (here $\Omega$ is the loop
functor on $\mathrm{Ho}(Sp(A-Mod))$).
Note that if the suspension functor is fully faithful, connective stable modules are exactly 
connective objects with respect to the natural $t$-structure on $A$-modules.

\begin{prop}\label{pstabmod}
For any $A\in Comm(\mathcal{C})$,
the functor $\underline{h}_{s}^{-}$ has a
total right derived functor
$$\mathbb{R}\underline{h}_{s}^{-} : \mathrm{Ho}(Sp(A-Mod))^{op} \longrightarrow
\mathrm{Ho}((A-Mod_{0}^{op})^{\wedge}),$$
which commutes with homotopy limits
\footnote{This makes
sense as the functor $\mathbb{R}\underline{h}_{s}^{-}$ is naturally
defined on the level of the Dwyer-Kan simplicial localizations with respect to equivalences
$$LSp(A-Mod)^{op} \longrightarrow
L((A-Mod_{0}^{op})^{\wedge}).$$}.
If the suspension functor
$$S : \mathrm{Ho}(\mathcal{C}) \longrightarrow \mathrm{Ho}(\mathcal{C})$$
is fully faithful, and if $A\in \mathcal{A}$, then for any integer $n\geq 0$, the functor
$\mathbb{R}\underline{h}_{s}^{-}$ is fully faithful when restricted to the
full subcategory of $(-n)$-connective objects.
\end{prop}

\begin{proof} As
$S_{A} : A-Mod_{0} \longrightarrow Sp(A-Mod)$ preserves
equivalences, one checks easily that
$\underline{h}_{s}^{M_{*}}$ is a fibrant object in
$(A-Mod_{0}^{op})^{\wedge}$ when $M_{*}$ is cofibrant
in $Sp(A-Mod)$. This easily implies that
$M_{*} \mapsto \underline{h}_{s}^{QM_{*}}$ is a right derived
functor for $\underline{h}_{s}^{-}$, and the standard properties
of mapping spaces imply that it commutes with homotopy limits.

We now assume that the suspension functor
$$S : \mathrm{Ho}(\mathcal{C}) \longrightarrow \mathrm{Ho}(\mathcal{C})$$
is fully faithful, and that
$A\in \mathcal{A}$. Let $n\geq 0$ be an integer.

Let
$S^{n} : A-Mod \longrightarrow A-Mod$ be a left Quillen functor which is
a model for the
suspension functor iterated $n$ times (e.g.
$S^{n}(N):=S^{n}_{A}\otimes_{A} N$, where
$S^{n}_{A}:=(S^{1}_{A})^{\otimes_{A} n}$). There is a pullback functor
$$(S^{n})^{*} : \mathrm{Ho}((A-Mod_{0}^{op})^{\wedge}) \longrightarrow
\mathrm{Ho}((A-Mod_{0}^{op})^{\wedge})$$
defined by $(S^{n})^{*}(F)(N):=F(S^{n}(N))$ for any $N\in A-Mod_{0}$
(note that $S^{n}$ stablizes the subcategory
$A-Mod_{0}$ because of our assumption \ref{ass7}).
For an $(-n)$-connective object $M_{*}$,
we claim there exists a natural isomorphism in
$\mathrm{Ho}((A-Mod_{0}^{op})^{\wedge})$
between $\mathbb{R}\underline{h}_{s}^{M_{*}}$ and
$(S^{n})^{*}(\mathbb{R}\underline{h}^{M})$, where
$\mathbb{R}\underline{h}^{M}$ is the value at $M$ of the restricted
Yoneda embedding
$$\mathbb{R}\underline{h}^{M}_{0} :
\mathrm{Ho}(A-Mod) \longrightarrow \mathrm{Ho}((A-Mod_{0}^{op})^{\wedge}).$$

Indeed, let us write $M$ as
$\Omega^{n}(S_{A}(M))$ for some object
$M\in \mathrm{Ho}(A-Mod)$, where $\Omega^{n}$ is the loop
functor of $Sp(A-Mod)$, iterated $n$ times.
Then, using our lemma Lem. \ref{lstabmod}, for any
$N\in \mathrm{Ho}(A-Mod)$, we have natural isomorphisms
in $\mathrm{Ho}(SSet)$
$$Map_{Sp(A-Mod)}(M_{*},S_{A}(N))\simeq
Map_{Sp(A-Mod)}(S_{A}(M),S^{n}(S_{A}(N)))\simeq $$
$$\simeq Map_{A-Mod}(M,S^{n}(N)),$$
where $S^{n}$ denotes the suspension functor iterated
$n$-times. Using that $A$ is
$\mathcal{C}_{0}$-good, this shows that
$$\mathbb{R}\underline{h}_{s}^{M_{*}}\simeq
(S^{n})^{*}(\mathbb{R}\underline{h}_{0}^{M}).$$
Moreover, for any $(-n)$-connective objects
$M_{*}$ and $N_{*}$ in $\mathrm{Ho}(Sp(A-Mod))$ one has natural
isomorphisms in $\mathrm{Ho}(SSet)$
$$Map_{Sp(A-Mod)}(M_{*},N_{*})\simeq Map_{A-Mod}(M,N)$$
where $M_{*}\simeq \Omega^{n}(S_{A}(M))$ and
$N_{*}\simeq \Omega^{n}(S_{A}(N))$.
We are therefore reduced to show that for any
$A$-modules $M$ and $N$ in $\mathrm{Ho}(A-Mod)$, the  natural
morphism in $\mathrm{Ho}(SSet)$
$$Map_{A-Mod}(M,N)\longrightarrow
Map_{(A-Mod_{0}^{op})^{\wedge}}((S^{n})^{*}(\mathbb{R}\underline{h}_{0}^{N}),
(S^{n})^{*}(\mathbb{R}\underline{h}_{0}^{M}))$$
is an isomorphism. To see this, we define a morphism in the
opposite direction in the following way. Taking the $n$-th loop
functor on each side gives a morphism
$$Map_{(A-Mod^{op}_{0})^{\wedge}}((S^{n})^{*}(\mathbb{R}\underline{h}_{0}^{N}),
(S^{n})^{*}(\mathbb{R}\underline{h}_{0}^{M}))
\rightarrow $$
$$\longrightarrow Map_{(A-Mod_{0}^{op})^{\wedge}}(\Omega^{n}(S^{n})^{*}(\mathbb{R}\underline{h}_{0}^{N}),
\Omega^{n}(S^{n})^{*}(\mathbb{R}\underline{h}_{0}^{M}).$$
Moreover, there are isomorphisms
$$\Omega^{n}Map_{A-Mod}(M,S^{n}(P))\simeq
Map_{A-Mod}(M,\Omega^{n}(S^{n}(P)))\simeq
Map_{A-Mod}(M,P)\simeq \mathbb{R}\underline{h}_{0}^{M}(P),$$
showing that there exists a natural isomorphism in
$\mathrm{Ho}((A-Mod_{0}^{op})^{\wedge})$ between
$\Omega^{n}(S^{n})^{*}(\mathbb{R}\underline{h}_{0}^{M})$ and
$\mathbb{R}\underline{h}_{0}^{M}$.  One therefore gets a morphism
$$Map_{(A-Mod_{0}^{op})^{\wedge}}((S^{n})^{*}(\mathbb{R}\underline{h}_{0}^{N}),
(S^{n})^{*}(\mathbb{R}\underline{h}_{0}^{M})) \rightarrow $$
$$\longrightarrow Map_{(A-Mod_{0}^{op})^{\wedge}}(\Omega^{n}(S^{n})^{*}(
\mathbb{R}\underline{h}_{0}^{N}),
\Omega^{n}(S^{n})^{*}(\mathbb{R}\underline{h}_{0}^{M})\simeq$$
$$\simeq
Map_{(A-Mod_{0}^{op})^{\wedge}}
(\mathbb{R}\underline{h}_{0}^{N},\mathbb{R}\underline{h}_{0}^{M}).$$
Using that $A$ is $\mathcal{C}_{0}$-good we get the required morphism
$$Map_{(A-Mod_{0}^{op})^{\wedge}}((S^{n})^{*}(\mathbb{R}\underline{h}_{0}^{N}),
(S^{n})^{*}(\mathbb{R}\underline{h}_{0}^{M})) \longrightarrow
Map_{A-Mod}(M,N),$$
and it is easy to check it is an inverse in $\mathrm{Ho}(SSet)$ to the natural morphism
$$Map_{A-Mod}(M,N)\longrightarrow
Map_{(A-Mod_{0}^{op})^{\wedge}}((S^{n})^{*}(\mathbb{R}\underline{h}_{0}^{N}),
(S^{n})^{*}(\mathbb{R}\underline{h}_{0}^{M})).$$
\end{proof}

\section{Descent for modules and stable modules}

In this last section we present some 
definitions concerning descent for modules and stable modules in our general context. 
In a few words, a co-augmented co-simplicial object $A \longrightarrow B_{*}$ in $Comm(\mathcal{C})$,
 is said to \textit{have the descent property for modules}
(resp. \textit{for stable modules}) if the homotopy theory of $A$-modules (resp. of stable  $A$-modules)
is equivalent to the homotopy theory of certain co-simplicial $B_{*}$-modules (resp. stable
$B_{*}$-modules). From the geometric, dual point of view, the object $A \longrightarrow B_{*}$ should be thought
as an augmented simplicial  space $Y_{*} \longrightarrow X$, and having the descent property 
essentially means that the theory of sheaves on $X$ is equivalent to the theory of certain sheaves
on $Y_{*}$ (see \cite[Exp. $V^{bis}$]{sga4II} for more details on this point of view).

The purpose of this section is only to introduce the basic set-up for descent that will be
used later to state that the theory of quasi-coherent modules is local with respect to the
topology, or in other words that hypercovers have the descent property for modules. 
This is a very important property allowing local-to-global arguments. \\

Let $A_{*}$ be a co-simplicial object
in the category $Comm(\mathcal{C})$ of commutative monoids.
Therefore, $A_{*}$ is given by a functor
$$\begin{array}{ccc}
\Delta &  \longrightarrow & Comm(\mathcal{C}) \\
 \left[ n\right]   & \mapsto & A_{n}.
\end{array}$$
A co-simplicial $A_{*}$-module $M_{*}$ is by definition the following datum.
\begin{itemize}
\item A $A_{n}$-module $M_{n} \in A_{n}-Mod$ for any $n\in \Delta$.
\item For any morphism $u : [n] \rightarrow [m]$ in $\Delta$,
a morphism of $A_{n}$-modules $\alpha_{u} : M_{n} \longrightarrow M_{m}$,
such that $\alpha_{v}\circ \alpha_{u}=\alpha_{v\circ u}$ for
any $\xymatrix{[n] \ar[r]^{u} & [m] \ar[r]^{v} & [p]}$ in $\Delta$.
\end{itemize}
In the same way, a morphism of co-simplicial $A_{*}$-modules
$f : M_{*} \longrightarrow N_{*}$ is the data of
morphisms $f_{n} : M_{n} \longrightarrow N_{n}$ for any
$n$, commuting with the $\alpha$'s
$$\alpha^{N}_{u}\circ f_{n}=f_{m}\circ \alpha^{M}_{u}$$
for any $u : [n] \rightarrow [m]$.

The co-simplicial $A_{*}$-modules and morphisms of $A_{*}$-modules form
a category, denoted by $csA_{*}-Mod$. It is furthermore a
$\mathbb{U}$-combinatorial model category
for which the equivalences (resp. fibrations)
are the morphisms $f : M_{*} \longrightarrow N_{*}$ such that
each $f_{n} : M_{n} \longrightarrow N_{n}$
is an equivalence (resp. a fibration) in $A_{n}-Mod$. \\

Let $A$ be a commutative monoid, and
$B_{*}$ be a co-simplicial commutative
$A$-algebra. We can also consider
$B_{*}$ as a co-simplicial commutative monoid
together with a co-augmentation
$$A \longrightarrow B_{*},$$
which is a morphism of co-simplicial objects
when  $A$ is considered as a constant co-simplicial
object. As a co-simplicial commutative monoid
$A$ possesses a category of co-simplicial $A_{*}$-modules
$csA-Mod$ which is nothing else than the
model category of co-simplicial objects in
$A-Mod$ with its projective levelwise model structure.

For a co-simplicial $A$-module $M_{*}$, we define
a co-simplicial $B_{*}$-module $B_{*}\otimes_{A}M_{*}$ by the formula
$$(B_{*}\otimes_{A}M_{*}):=B_{n}\otimes_{A}M_{n},$$
and for which the transitions morphisms are given
by the one of $B_{*}$ and of $M_{*}$.
This construction defines a functor
$$B_{*}\otimes_{A} - : csA-Mod \longrightarrow
csB_{*}-Mod,$$
which has a right adjoint
$$csB_{*}-Mod \longrightarrow csA-Mod.$$
This right adjoint, sends a $B_{*}$-module $M_{*}$
to its underlying co-simplicial $A$-module. Clearly,
this defines a Quillen adjunction.
There exists  another adjunction
$$ct : \mathrm{Ho}(A-Mod) \longrightarrow
\mathrm{Ho}(csA-Mod) \qquad \mathrm{Ho}(A-Mod) \longleftarrow
\mathrm{Ho}(csA-Mod) : Holim,$$
where $Holim$ is defined as the total right derived functor
of the functor lim (see \cite[8.5]{hi}, \cite[10.13]{ds}).

Composing these two adjunctions gives a new adjunction
$$B_{*}\otimes^{\mathbb{L}}_{A}- : \mathrm{Ho}(A-Mod) \longrightarrow
\mathrm{Ho}(csB_{*}-Mod) \qquad \mathrm{Ho}(A-Mod) \longleftarrow
\mathrm{Ho}(csB_{*}-Mod) : \int.$$

\begin{df}\label{ddesc}
\begin{enumerate}
\item
Let $B_{*}$ be a co-simplicial commutative monoid and
$M_{*}$ be a co-simplicial $B_{*}$-module. We say that
$M_{*}$ \emph{is homotopy cartesian}\index{homotopy cartesian! cosimplicial module}
if for any $u : [n] \rightarrow [m]$ in $\Delta$
the morphism, induced by $\alpha_{u}$,
$$M_{n}\otimes_{B_{n}}^{\mathbb{L}}B_{m} \longrightarrow
M_{m}$$
is an isomorphism in $\mathrm{Ho}(B_{m}-Mod)$.
\item Let $A$ be a commutative monoid and
$B_{*}$ be a co-simplicial commutative
$A$-algebra. We say that the co-augmentation morphism
$A \longrightarrow B_{*}$ \emph{satisfies the
descent condition}, if
in the adjunction
$$B_{*}\otimes^{\mathbb{L}}_{A}- : \mathrm{Ho}(A-Mod) \longrightarrow
\mathrm{Ho}(csB_{*}-Mod) \qquad \mathrm{Ho}(A-Mod) \longleftarrow
\mathrm{Ho}(csB_{*}-Mod) : \int$$
the functor $B_{*}\otimes^{\mathbb{L}}_{A}-$ is fully faithful and
induces an equivalence between $\mathrm{Ho}(A-Mod)$ and the
full subcategory of $\mathrm{Ho}(csB_{*}-Mod)$ consisting of
homotopy cartesian objects.
\end{enumerate}
\end{df}

\begin{rmk}
\emph{If a morphism $A \longrightarrow B_{*}$ satisfies the descent
condition, then so does any co-simplicial object
equivalent to it (as a morphism in the model category
of co-simplicial commutative monoids).}
\end{rmk}

\begin{lem}\label{ldesc}
Let $A \longrightarrow B_{*}$ be a co-augmented
co-simplicial commutative monoid in $Comm(\mathcal{C})$
which satisfies the descent condition. Then, the natural
morphism
$$A \longrightarrow Holim_{n}B_{n}$$
is an isomorphism in $\mathrm{Ho}(Comm(\mathcal{C}))$.
\end{lem}

\begin{proof} This is clear as
$\int B_{*}=Holim_{n}B_{n}$. \end{proof}

We now pass to descent for stable modules.
For this, let again $A_{*}$ be a co-simplicial object
in the category $Comm(\mathcal{C})$ of commutative monoids.
A co-simplicial stable $A_{*}$-module $M_{*}$ is by definition the following datum.
\begin{itemize}
\item A stable $A_{n}$-module $M_{n} \in Sp(A_{n}-Mod)$
for any $n\in \Delta$.
\item For any morphism $u : [n] \rightarrow [m]$ in $\Delta$,
a morphism of stable $A_{n}$-modules $\alpha_{u} : M_{n} \longrightarrow M_{m}$,
such that $\alpha_{v}\circ \alpha_{u}=\alpha_{v\circ u}$ for
any $\xymatrix{[n] \ar[r]^{u} & [m] \ar[r]^{v} & [p]}$ in $\Delta$.
\end{itemize}
In the same way, a morphism of co-simplicial stable $A_{*}$-modules
$f : M_{*} \longrightarrow N_{*}$ is the data of
morphisms $f_{n} : M_{n} \longrightarrow N_{n}$ in $Sp(A_{n}-Mod)$ for any
$n$, commuting with the $\alpha$'s
$$\alpha^{N}_{u}\circ f_{n}=f_{m}\circ \alpha^{M}_{u}$$
for any $u : [n] \rightarrow [m]$.

The co-simplicial stable $A_{*}$-modules and morphisms of $A_{*}$-modules form
a category, denoted by $Sp(csA_{*}-Mod)$. It is furthermore a
$\mathbb{U}$-combinatorial model category
for which the equivalences (resp. fibrations)
are the morphisms $f : M_{*} \longrightarrow N_{*}$ such that
each $f_{n} : M_{n} \longrightarrow N_{n}$
is an equivalence (resp. a fibration) in $Sp(A_{n}-Mod)$. We note that
$Sp(csA_{*}-Mod)$ is also naturally equivalent as a category to the
category of $S^{1}_{A}$-spectra in
$csA_{*}Mod$, hence the notation $Sp(csA_{*}-Mod)$
is not ambiguous. \\

Let $A$ be a commutative monoid, and
$B_{*}$ be a co-simplicial commutative
$A$-algebra. We can also consider
$B_{*}$ as a co-simplicial commutative monoid
together with a co-augmentation
$$A \longrightarrow B_{*},$$
which is a morphism of co-simplicial objects
when $A$ is considered as a constant co-simplicial
object. As a co-simplicial commutative monoid
$A$ possesses a category of co-simplicial stable $A_{*}$-modules
$Sp(csA-Mod)$ which is nothing else than the
model category of co-simplicial objects in
$Sp(A-Mod)$ with its projective levelwise model structure.

For a co-simplicial stable $A$-module $M_{*}$, we define
a co-simplicial stable $B_{*}$-module $B_{*}\otimes_{A}M_{*}$ by the formula
$$(B_{*}\otimes_{A}M_{*}):=B_{n}\otimes_{A}M_{n},$$
and for which the transitions morphisms are given
by the one of $B_{*}$ and of $M_{*}$. This construction defines a functor
$$B_{*}\otimes_{A} - : Sp(csA-Mod) \longrightarrow
Sp(csB_{*}-Mod),$$
which has a right adjoint
$$Sp(csB_{*}-Mod) \longrightarrow Sp(csA-Mod).$$
This right adjoint, sends a  co-simplicial stable $B_{*}$-module $M_{*}$
to its underlying co-simplicial stable $A$-module. Clearly,
this is a Quillen adjunction.
One also has an adjunction
$$ct : \mathrm{Ho}(Sp(A-Mod)) \longrightarrow
\mathrm{Ho}(Sp(csA-Mod)) $$
$$ \mathrm{Ho}(Sp(A-Mod)) \longleftarrow
\mathrm{Ho}(Sp(csA-Mod)) : Holim,$$
where $Holim$ is defined as the total right derived functor
of the functor lim (see \cite[8.5]{hi}, \cite[10.13]{ds}).

Composing these two adjunctions gives a new adjunction
$$B_{*}\otimes^{\mathbb{L}}_{A}- : \mathrm{Ho}(Sp(A-Mod)) \longrightarrow
\mathrm{Ho}(Sp(csB_{*}-Mod)) $$ 
$$\mathrm{Ho}(Sp(A-Mod)) \longleftarrow
\mathrm{Ho}(Sp(csB_{*}-Mod)) : \int.$$

\begin{df}\label{ddesc'}
\begin{enumerate}
\item
Let $B_{*}$ be a co-simplicial commutative monoid and
$M_{*}$ be a co-simplicial stable $B_{*}$-module. We say that
$M_{*}$ \emph{is homotopy cartesian}
if for any $u : [n] \rightarrow [m]$ in $\Delta$
the morphism, induced by $\alpha_{u}$,
$$M_{n}\otimes_{B_{n}}^{\mathbb{L}}B_{m} \longrightarrow
M_{m}$$
is an isomorphism in $\mathrm{Ho}(Sp(B_{m}-Mod))$.
\item Let $A$ be a commutative monoid and
$B_{*}$ be a co-simplicial commutative
$A$-algebra. We say that the co-augmentation morphism
$A \longrightarrow B_{*}$ \emph{satisfies the
stable descent condition}, if
in the adjunction
$$B_{*}\otimes^{\mathbb{L}}_{A}- : \mathrm{Ho}(Sp(A-Mod)) \longrightarrow
\mathrm{Ho}(Sp(csB_{*}-Mod)) $$ 
$$\mathrm{Ho}(Sp(A-Mod)) \longleftarrow
\mathrm{Ho}(Sp(csB_{*}-Mod)) : \int$$
the functor $B_{*}\otimes^{\mathbb{L}}_{A}-$ is fully faithful and
induces an equivalence between $\mathrm{Ho}(Sp(A-Mod))$ and the
full subcategory of $\mathrm{Ho}(Sp(csB_{*}-Mod))$ consisting of
homotopy cartesian objects.
\end{enumerate}
\end{df}

The following proposition insures that
the unstable descent condition implies the stable one
under certain conditions.

\begin{prop}\label{pdescstable}
Let $A$ be a commutative monoid and
$B_{*}$ be a co-simplicial commutative
$A$-algebra, such that $A\longrightarrow B_{*}$ satisfies the
descent condition. Assume that the
two following conditions are satisfied.
\begin{enumerate}
\item The suspension functor $S : \mathrm{Ho}(\mathcal{C}) \longrightarrow
\mathrm{Ho}(\mathcal{C})$ is fully faithful.
\item For any $n$, the morphism $A \longrightarrow B_{n}$
is flat in the sense of Def. \ref{d7-}.
\end{enumerate}
Then $A\longrightarrow B_{*}$ satisfies the stable descent property.
\end{prop}

\begin{proof} The second condition insures that
the diagram
$$\xymatrix{
\mathrm{Ho}(Sp(A-Mod)) \ar[r] \ar[d] &
\mathrm{Ho}(Sp(csB_{*}-Mod) \ar[d] \\
\mathrm{Ho}(A-Mod) \ar[r] & \mathrm{Ho}(csB_{*}-Mod)}
$$
commutes up to a natural isomorphism (here the
vertical functors are the right adjoint to the
suspensions inclusion functors, and send a spectrum
to its $0$-th level). So, for any
$M\in Sp(A-Mod)$, the adjunction morphism
$$M\longrightarrow \int(B_{*}\otimes^{\mathbb{L}}M)$$
is easily seen to induces an isomorphism on the
$0$-th level objects in $\mathrm{Ho}(A-Mod)$. As this is true
for any $M$, and in particular for the suspensions of $M$,
we find that this adjunction morphism induces isomorphisms
on each $n$-th level objects in $\mathrm{Ho}(A-Mod)$. Therefore
it is an isomorphism. In the same way, one proves that
for any cartesian object $M_{*}\in \mathrm{Ho}(Sp(csB_{*}-Mod))$
the adjunction morphism
$$B_{*}\otimes_{A}\int(M_{*}) \longrightarrow M_{*}$$
is an isomorphism. \end{proof}

\section{Comparison with the usual notions}

In this last section we present what the general
notions introduced before give, when $\mathcal{C}$ is
the model category of $k$-modules with the trivial model structure.
The other non trivial examples will be given at the beginning
of the various chapters \S \ref{IIder}, \S \ref{IIunb} and \S \ref{IIbnag} where the
case of simplicial modules, complexes and symmetric spectra
will be studied. \\

Let $k$ be a commutative ring in $\mathbb{U}$, and
we let $\mathcal{C}$ be the category of $k$-modules belonging
to the universe $\mathbb{U}$, endowed with its
trivial model structure for which equivalences
are isomorphisms and all morphisms are fibrations
and cofibrations. The category $\mathcal{C}$ is then
a symmetric monoidal model category for the
tensor product $-\otimes_{k}-$. Furthermore, all of our
assumption \ref{ass-1}, \ref{ass0}, \ref{ass1}
and \ref{ass2} are satisfied. The category $Comm(\mathcal{C})$ is of course the category
of commutative $k$-algebras in $\mathbb{U}$, endowed
with its trivial model structure. In the same way,
for any commutative $k$-algebra $A$, the category
$A-Mod$ is the usual category of $A$-modules
in $\mathbb{U}$ together with its trivial model structure.
We set $\mathcal{C}_{0}=\mathcal{C}$ and
$\mathcal{A}=Comm(\mathcal{C})$.
The notions we have presented before restrict essentially
to the usual notions of algebraic geometry with some remarkable caveats.

\begin{itemize}

\item For any morphism of commutative $k$-algebras $A \longrightarrow B$
and any $B$-module $M$, the simplicial set
$\mathbb{D}er_{A}(B,M)$ is discrete and naturally
isomorphic to the set $Der_{A}(B,M)$ of derivations
from $B$ to $M$ over $A$. Equivalently,
the $B$-module $\mathbb{L}_{B/A}$ (defined in Def. \ref{d2}) is simply the usual
$B$-module of K\"ahler differentials $\Omega_{B/A}^{1}$. Note that this $\mathbb{L}_{B/A}$
 is \textit{not} the Quillen-Illusie cotangent complex of $A\rightarrow B$.

\item For any morphism of commutative $k$-algebras $A \longrightarrow B$
the natural morphism $B\longrightarrow THH(B/A)$
is always an isomorphism. Indeed
$$THH(A)\simeq A\otimes_{A\otimes_{k}A}A\simeq A.$$

\item A morphism of commutative $k$-algebras $A \longrightarrow B$
is finitely presented in the sense of Def. \ref{d3} if and only
if it $B$ is a finitely presented $A$-algebra in the usual sense.
In the same way, finitely presented objects
in $A-Mod$ in the sense of Def. \ref{d3} are the finitely
presented $A$-modules. Also, perfect objects
in $A-Mod$ are the projective $A$-modules of finite type.

\item A morphism $A \longrightarrow B$ of commutative
$k$-algebras is a formal covering if and only if
it is a faithful morphism of rings.

\item Let $f : A\longrightarrow B$ be a morphism
of commutative $k$-algebras, and $Spec\, B \longrightarrow
Spec\, A$ the associated morphism of schemes.
We have the following comparison board.\\

$$\begin{array}{c|c}
\mathbf{In\; the \; sense\; of \; Def. \; \ref{d5},\ref{d6},\ref{d7}} &
\mathbf{As\; a\; morphism\; of\; affine\; schemes} \\
 & \\
epimorphism & monomorphism \\
flat & flat \\
formal \; Zariski \; open \; immersion & flat\; monomorphism \\
Zariski\; open\; immersion & open\; immersion \\
(formally)\; unramified & (formally)\; unramified \\
formally\; thh-etale & always \; satisfied \\
formally\; i-smooth & always \; satisfied \\
thh-etale & finitely\; presented \\
\end{array}$$

The reader will immediately notice the absence of (formally) smooth and (formally) \'etale
maps in the previous table. This is essentially due to the fact that there are no easy \textit{general}
characterizations of this maps in terms of the module of K\"ahler differentials alone (which in this
trivial model structure context is our cotangent complex). On the contrary such characterizations do
exist in terms of the ``correct'' cotangent complex which is the Quillen-Illusie one. But this
correct cotangent complex of a morphism of usual commutative $k$-algebras $A\rightarrow B$ will appear as
the cotangent complex according to our definition \ref{d2} only if we consider this morphism in the
category of simplicial $k$-modules, i.e. if we replace the category $\mathcal{C}$ of $k$-modules with
the category of simplicial $k$-modules. In other words our definitions of (formally) smooth and (formally) \'etale
maps reduce to the usual ones (e.g. \cite{egaIV-4}) between commutative $k$-algebras only if we consider them
in the context of \textit{derived algebraic geometry} (see \S \ref{IIder}). This is consistent with the general philosophy
that some aspects of usual algebraic geometry, especially those related to infinitesimal lifting properties and
deformation theory, are conceptually more transparent in (and actually already a part of)
derived algebraic geometry. See also Remark \ref{why} for another instance of this point of view.

\item For commutative monoid $A$,
one has $\pi_{*}(A)=A$, and
for any $a\in A$ the commutative
$A$-algebra $A \longrightarrow A[a^{-1}]$
is the usual localization of $A$ inverting
$a$.

\item One has that $A_{K}\simeq 0$ for any
projective $A$-module $K$ of constant finite rank $n>0$ (for example
if the scheme $Spec\, A$ is connected and $K$ is non zero).

\end{itemize}

\chapter{Geometric stacks: Basic theory}\label{partI.3}

In this chapter, after a brief reminder of \cite{hagI}, we present the key definition of this
work, the notion of \emph{$n$-geometric stack}. The definition we present
here is a generalization of the original notion of \emph{geometric $n$-stack}
introduced by C. Simpson in \cite{s4}. As already remarked in
\cite{s4}, the notion of $n$-geometric stack only depends on a
topology on the opposite category of commutative monoids $Comm(\mathcal{C})^{op}$, and
on a class \textbf{P} of morphisms. Roughly speaking, geomeric stacks are
the stacks obtained by taking quotient of representable stacks
by some equivalence relations in \textit{P}. By choosing different
classes \textbf{P} one gets different notions of geometric stacks.
For example, in the classical situation where $\mathcal{C}=\mathbb{Z}-Mod$, and
the topology is chosen to be the \'etale topology, Deligne-Mumford algebraic stacks
correspond to the case where \textbf{P} is the class of \'etale morphisms, whereas
Artin algebraic stacks correspond to the case where \textbf{P} is the class
of smooth morphisms. We think it is
important to leave the choice of the class 
\textbf{P} open in the general definition, so
that it can be specialized differently depending of the kind of
objects one is willing to consider. \\

From the second section on, we will fix a HA context
$(\mathcal{C},\mathcal{C}_{0},\mathcal{A})$, in the sense of
Def. \ref{dha}.

\section{Reminders on model topoi}\label{remhagI}

To make the paper essentially self-contained, we briefly summarize in this subsection the basic notions and results of the theory of stacks
in homotopical contexts as exposed in \cite{hagI}. We will limit ourselves to recall only the topics
that will be needed in the sequel; the reader is addressed to \cite{hagI} for further details and for proofs.\\

Let $M$ be a $\mathbb{U}$-small model category, and $W(M)$ its class of weak equivalences. We let $SPr(M):=SSets^{M^{op}}$ be the category of simplicial presheaves on $M$ with its projective model structure, i.e. with equivalences and fibrations defined objectwise. 

The model category of \textit{prestacks} $M^{\wedge}$ on $M$ is the model category obtained as the left Bousfield localization of 
$SPr(M)$ at $\left\{h_{u}\;|\; u \in W(M)\right\}$, where $h:M\rightarrow Pr(M)\hookrightarrow SPr(M)$ is the (constant) Yoneda embedding. 
The homotopy category 
 $\mathrm{Ho}(M^{\wedge})$ can be identified with the full subcategory of $\mathrm{Ho}(SPr(M))$ consisting of those simplicial presheaves $F$ on $M$ that preserve weak equivalences (\cite[Def. 4.1.4]{hagI}); any such simplicial presheaf (i.e. any object in $\mathrm{Ho}(M^{\wedge})$) will be called a \textit{prestack on} $M$. $M^{\wedge}$ is a $\mathbb{U}$-cellular and $\mathbb{U}$-combinatorial (\cite[App. A]{hagI}) simplicial model category and its derived simplicial $Hom$'s will be denoted simply by $\mathbb{R}\underline{Hom}$ (denoted as $\mathbb{R}_{w}\underline{Hom}$ in \cite[\S 4.1]{hagI}). \\

If $\Gamma_{*} : M\rightarrow M^{\Delta}$ is a cofibrant resolution functor for $M$ (\cite[16.1]{hi}), and we define $$\underline{h}: M \rightarrow M^{\wedge} \;: \; x\longmapsto (\underline{h}_{x}: y\mapsto Hom_{M}(\Gamma_{*}(y),x)),$$ we have that $\underline{h}$ preserves fibrant objects and weak equivalences between fibrant objects (\cite[Lem. 4.2.1]{hagI}). Therefore we can right-derive $\underline{h}$ to get a functor $\mathbb{R}\underline{h}:=\underline{h}\circ R:\mathrm{Ho}(M)\rightarrow \mathrm{Ho}(M^{\wedge})$, where $R$ is a fibrant replacement functor in $M$; $\mathbb{R}\underline{h}$ is in fact fully faithful (\cite[Thm. 4.2.3]{hagI}) and is therefore called the (\textit{model}) \textit{Yoneda embedding} for the model category $M$ ($\mathbb{R}\underline{h}$, as opposed to $\underline{h}$, does not depend, up to a unique isomorphism, on the choice of the cofibrant resolution functor $\Gamma_{*}$).

We also recall that the canonical morphism $h_{x}\rightarrow \mathbb{R}\underline{h}_{x}$ is always an isomorphism in $\mathrm{Ho}(M^{\wedge})$ (\cite[Lem. 4.2.2]{hagI}), and that with the notations introduced above for the derived simplicial $Hom$'s in $M^{\wedge}$, the model Yoneda lemma (\cite[Cor. 4.2.4]{hagI}) is expressed by the isomorphisms in $\mathrm{Ho}(SSets)$
$$\mathbb{R}\underline{Hom}(\mathbb{R}\underline{h}_{x},F)\simeq \mathbb{R}\underline{Hom}(h_{x},F)\simeq F(x)$$ for any fibrant object $F$ in $M^{\wedge}$.  \\

A convenient homotopical replacement of the notion of a Grothendieck topology in the case of model categories, is the following (\cite[Def. 4.3.1]{hagI})

\begin{df}\label{modtop}
A \emph{model (pre-)topology}\index{model!(pre)-topology} $\tau$ on a $\mathbb{U}$-small model
category $M$, is the datum for any object $x \in M$, of
a set $Cov_{\tau}(x)$ of subsets
of objects in $\mathrm{Ho}(M)/x$, called $\tau$\emph{-covering families} of
$x$, satisfying the following three conditions

\begin{enumerate}

\item \emph{(Stability)} For all $x \in M$ and any isomorphism $y \rightarrow x$ in $\mathrm{Ho}(M)$,
the one-element set $\{y \rightarrow x\}$ is in $Cov_{\tau}(x)$.

\item \emph{(Composition)} If $\{u_{i} \rightarrow x\}_{i \in I} \in Cov_{\tau}(x)$, and for any $i \in I$,
$\{v_{ij} \rightarrow u_{i}\}_{j \in J_{i}} \in Cov_{\tau}(u_{i})$,
the family $\{v_{ij} \rightarrow x\}_{i \in I, j\in J_{i}}$
is in $Cov_{\tau}(x)$.

\item \emph{(Homotopy base change)} Assume the two previous conditions hold. For any
$\{u_{i} \rightarrow x\}_{i \in I} \in Cov_{\tau}(x)$, and any morphism in $\mathrm{Ho}(M)$,
$y \rightarrow x$, the family $\{u_{i}\times^{h}_{x}y \rightarrow y\}_{i \in I}$ is in $Cov_{\tau}(y)$.

\end{enumerate}
A $\mathbb{U}$-small model category $M$
together with a model pre-topology $\tau$ will be called a $\mathbb{U}$-small \emph{model site}\index{model!site}.
\end{df}

By \cite[Prop. 4.3.5]{hagI} a model pre-topology $\tau$ on $M$ induces and is essentially the same thing as a Grothendieck topology, still denoted by $\tau$, on the homotopy category $\mathrm{Ho}(M)$.\\

Given a model site $(M,\tau)$ we have, as in \cite[Thm. 4.6.1]{hagI},
a model category $M^{\sim,\tau}$ ($\mathbb{U}$-combinatorial and left proper) of \textit{stacks} on the model site,
which is defined  as the left Bousfield localization of the model category $M^{\wedge}$ of prestacks on $M$ along a class $H_{\tau}$ of \textit{homotopy} $\tau$\textit{-hypercovers} (\cite[4.4, 4.5]{hagI}).
To any prestack $F$ we can associate
a sheaf $\pi_{0}$ of connected components on the site $(\mathrm{Ho}(M), \tau)$ defined as the associated sheaf to the presheaf $x\longmapsto \pi_{0}(F(x))$.  In a similar way (\cite[Def. 4.5.3]{hagI}),
for any $i>0$, any fibrant object $x\in M$, and any $s\in F(x)_{0}$, we can define a sheaf of homotopy groups $\pi_{i}(F,s)$ on the induced comma site $(\mathrm{Ho}(M/x), \tau)$.
The weak equivalences in $M^{\sim,\tau}$ turn out to be exactly the $\pi_{*}$-sheaves isomorphisms (\cite[Thm. 4.6.1]{hagI}), i.e. those maps $u:F\rightarrow G$ in $M^{\wedge}$ inducing an isomorphism of sheaves $\pi_{0}(F)\simeq \pi_{0}(G)$ on $(\mathrm{Ho}(M), \tau)$, and isomorphisms $\pi_{i}(F,s)\simeq \pi_{i}(G,u(s))$ of sheaves on $(\mathrm{Ho}(M/x), \tau)$ for any $i\geq 0$, for any choice of fibrant $x\in M$ and any base point $s\in F(x)_{0}$.\\

The left Bousfield localization construction defining $M^{\sim,\tau}$ yields a pair of adjoint Quillen functors
$$\mathrm{Id} : M^{\wedge} \longrightarrow M^{\sim,\tau}
\qquad M^{\wedge} \longleftarrow M^{\sim,\tau} : \mathrm{Id}$$
which induces an adjunction pair at the level of homotopy categories
$$a:=\mathbb{L}\mathrm{Id} : \mathrm{Ho}(M^{\wedge}) \longrightarrow \mathrm{Ho}(M^{\sim,\tau})
\qquad \mathrm{Ho}(M^{\wedge}) \longleftarrow \mathrm{Ho}(M^{\sim,\tau}) : j:=\mathbb{R}\mathrm{Id}$$
where $j$ is fully faithful.

\begin{df}\label{dstack0}
\begin{enumerate}
\item A \emph{stack}\index{stack!on a model site} on the model site $(M,\tau)$ is an object
$F\in SPr(M)$ whose image
in $\mathrm{Ho}(M^{\wedge})$ is in the essential image
of the functor $j$.
\item If $F$ and $G$ are stacks on the model site $(M,\tau)$, a \emph{morphism of stacks} is a morphism $F\rightarrow G$ in
in $\mathrm{Ho}(SPr(M))$, or equivalently in
$\mathrm{Ho}(M^{\wedge})$, or equivalently in
$\mathrm{Ho}(M^{\sim,\tau})$.
\item A morphism of stacks $f : F \longrightarrow G$
is a \emph{covering} (or a \emph{cover} or an \emph{epimorphism}) if the induced morphism of sheaves
$$\pi_{0}(f) : \pi_{0}(F) \longrightarrow  \pi_{0}(G)$$
is an epimorphism in the category of sheaves.
\end{enumerate}
\end{df}

Recall that a simplicial presheaf $F:M^{\mathrm{op}}\rightarrow SSet_{\mathbb{V}}$ is a stack if and only if it preserves weak equivalences and satisfy a $\tau$-\textit{hyperdescent} condition (descent, i.e. sheaf-like, condition with respect to the class $H_{\tau}$ of homotopy hypercovers): see \cite[Def. 4.6.5 and Cor. 4.6.3]{hagI}. We will always consider $\mathrm{Ho}(M^{\sim,\tau})$ embedded in $\mathrm{Ho}(M^{\wedge})$ embedded in $\mathrm{Ho}(SPr(M))$, omitting in particular to mention explicitly the functor $j$ above. With this conventions, the functor $a$ above becomes an endofunctor of $\mathrm{Ho}(M^{\wedge})$, called the \textit{associated stack} functor for the model site $(M,\tau)$. The associated stack functor preserves finite homotopy limits and all homotopy colimits (\cite[Prop. 4.6.7]{hagI}).

\begin{df}\label{sottocan}
A model pre-topology $\tau$ on $M$ is \emph{sub-canonical}\index{model!(pre)topology!subcanonical} if for any $x\in M$ the pre-stack $\mathbb{R}\underline{h}_{x}$ is a stack.
\end{df}

$M^{\sim,\tau}$ is a left proper (but not right proper) simplicial model category and its derived simplicial $Hom$'s will be denoted by $\mathbb{R}_{\tau}\underline{Hom}$ (denoted by $\mathbb{R}_{w,\,\tau}\underline{Hom}$ in \cite[Def. 4.6.6]{hagI}); for $F$ and $G$ prestacks on $M$, there is always a morphism in $\mathrm{Ho}(SSet)$ $$\mathbb{R}_{\tau}\underline{Hom}(F,G)\rightarrow \mathbb{R}\underline{Hom}(F,G)$$ which is an isomorphism when $G$ is a stack (\cite[Prop.4.6.7]{hagI}). \\
Moreover $M^{\sim,\tau}$ is a $t$-complete \textit{model topos} (\cite[Def. 3.8.2]{hagI}) therefore possesses important exactness properties. For the readers' convenience we collect below (from \cite{hagI}) the definition of (t-complete) model topoi and the main theorem characterizing them (Giraud-like theorem).\\

For a $\mathbb{U}$-combinatorial model
category $N$, and a $\mathbb{U}$-small set $S$ of
morphisms in $N$, we denote by $L_{S}N$ the left Bousfield
localization of $N$ along $S$. It is a model category, having
$N$ as underlying category, with the same cofibrations
as $N$ and whose equivalences are the $S$-local equivalences
(\cite[Ch. 3]{hi}).

\begin{df}\label{tcompleta}
\begin{enumerate}
    \item An object $x$ in a model category $N$ is \emph{truncated}\index{truncated object} 
 if there exists an integer $n\geq 0$, such that for  
 any $y\in N$ the mapping space $\mathrm{Map}_{N}(y,x)$ (\cite[\S 17.4]{hi}) is a $n$-truncated simplicial set\\
 (i.e. $\pi_{i}(\mathrm{Map}_{N}(x,y),u)=0$ for all $i>n$ and for all base point $u$).
    \item A model category $N$ is \emph{t-complete}\index{t-complete!model category} if a morphism $u:y\rightarrow y'$ in $\mathrm{Ho}(N)$ is an isomorphism if and only if the induced map $u^{*}:[y',x]\rightarrow [y,x]$ is a bijection for any truncated object $x\in N$.
\end{enumerate}
\end{df}

Recall that for any $\mathbb{U}$-small $S$-category $T$ (i.e. a category $T$ enriched over simplicial sets, \cite[Def. 2.1.1]{hagI}), we
can define a $\mathbb{U}$-combinatorial model
category $SPr(T)$ of \textit{simplicial} functors $T^{op} \longrightarrow SSet_{\mathbb{U}}$, in
which equivalences and fibrations
are defined levelwise (\cite[Def. 2.3.2]{hagI}).

\begin{df}\emph{(\cite[\S 3.8]{hagI})}\label{dmodtop}
A \emph{$\mathbb{U}$-model topos}\index{model!topos} is a $\mathbb{U}$-combinatorial model
category $N$ such that there exists
a $\mathbb{U}$-small $S$-category $T$ and a
$\mathbb{U}$-small set of morphisms $\mathcal{S}$ in $SPr(T)$
satisfying the following two conditions.
\begin{enumerate}
\item The model category $N$ is Quillen equivalent to
$L_{\mathcal{S}}SPr(T)$.
\item The identity functor
$$\mathrm{Id} : SPr(T) \longrightarrow L_{\mathcal{S}}SPr(T)$$
preserves homotopy pullbacks.
\end{enumerate}
A $t$\emph{-complete model topos}\index{t-complete model topos} is a $\mathbb{U}$-model topos $N$ which is $t$-complete as a model category.
\end{df}

We need to recall a few special morphisms in the standard simplicial category $\Delta$.
For any $n>0$, and $0\leq i< n$ we let
$$\begin{array}{cccc}
\sigma_{i} : & [1] & \longrightarrow & [n] \\
 & 0 & \mapsto & i \\
 & 1 & \mapsto & i+1.
\end{array}$$

\begin{df}\label{d12}
Let $N$ be a model category.
A \emph{Segal groupoid object in $N$}\index{Segal groupoid} is
a simplicial object
$$X_{*} : \Delta^{op} \longrightarrow N$$
satisfying the following two conditions.
\begin{enumerate}
\item For any $n>0$, the natural
morphism
$$\prod_{0\leq i<n}\sigma_{i} : X_{n} \longrightarrow \underbrace{X_{1}\times^{h}_{X_{0}}X_{1}\times^{h}_{X_{0}}
\dots \times^{h}_{X_{0}}X_{1}}_{n\; times}$$
is an isomorphism in $\mathrm{Ho}(N)$.
\item The morphism
$$d_{0}\times d_{1} :
X_{2} \longrightarrow X_{1}\times^{h}_{d_{0},X_{0},d_{0}}X_{1}$$
is an isomorphism in $\mathrm{Ho}(N)$.
\end{enumerate}
The homotopy category of Segal groupoid objects in $N$
is the full subcategory of $\mathrm{Ho}(N^{\Delta^{op}})$
consisting of Segal groupoid objects. It is denoted
by $\mathrm{Ho}(SeGpd(N))$.
\end{df}

The main theorem characterizing model topoi is the following
analog of Giraud's theorem.

\begin{thm}\emph{(\cite[Thm. 4.9.2]{hagI})}\label{tgiraud}\index{Giraud-like Theorem}
A model category $N$ is a model topos if and only
if it satisfies the following conditions.
\begin{enumerate}
\item The model category $N$ is $\mathbb{U}$-combinatorial.

\item For any $\mathbb{U}$-small family of objects
$\{x_{i}\}_{i\in I}$ in $N$, and any $i\neq j$ in $I$ the following square
$$\xymatrix{
\emptyset \ar[r] \ar[d] & x_{i} \ar[d] \\
x_{j} \ar[r] & \coprod_{k\in I}^{\mathbb{L}} x_{k}}$$
is homotopy cartesian.

\item For any $\mathbb{U}$-small category $I$, 
any morphism $y \rightarrow z$ and any
$I$-diagram $x : I \longrightarrow N/z$, the natural morphism
$$Hocolim_{i\in I}(x_{i}\times^{h}_{z}y) \longrightarrow
(Hocolim_{i\in I}x_{i})\times^{h}_{z}y$$
is an isomorphism in $\mathrm{Ho}(N)$.

\item For any Segal groupoid object (in the sense of Def. \ref{d12})
$$X_{*} : \Delta^{op} \longrightarrow N,$$
the natural morphism
$$X_{1} \longrightarrow X_{0}\times^{h}_{|X_{*}|}X_{0}$$
is an isomorphism in $\mathrm{Ho}(N)$.

\end{enumerate}
\end{thm}

An important consequence is the following

\begin{cor}\label{cgiraudapp}
For any $\mathbb{U}$-model topos $N$ and any fibrant object
$x\in N$, the category $\mathrm{Ho}(N/x)$ is cartesian closed.
\end{cor}

The exactness properties of model topoi will be frequently used all along 
this work. For instance, we will often use that for
any cover of stacks $p : F \longrightarrow G$ (over some model
site $(M,\tau)$), the natural morphism
$$|F_{*}| \longrightarrow G$$
is an isomorphism of stacks, where $F_{*}$ is the
homotopy nerve of $p$ (i.e. the nerve of 
a fibration equivalent to $p$, computed in the category of simplicial
presheaves). This  result is also recalled in 
Lem. \ref{lquot}.  

\section{Homotopical algebraic geometry context}

Let us fix a HA context
$(\mathcal{C},\mathcal{C}_{0},\mathcal{A})$.\\

We denote by $Aff_{\mathcal{C}}$ the opposite of the model
category $Comm(\mathcal{C})$: this will be our base model category $M$ to 
which we will apply the \cite{hagI} constructions recalled in \S \ref{remhagI}.

An object $X\in Aff_{\mathcal{C}}$ corresponding
to a commutative monoid $A \in Comm(\mathcal{C})$ will be symbolically
denoted by $X=Spec\, A$. We will consider
the model category $Aff_{\mathcal{C}}^{\wedge}$ of pre-stacks on
$Aff_{\mathcal{C}}$ as described in \S \ref{remhagI} above. By definition,
it is the left Bousfield localization
of $SPr(Aff_{\mathcal{C}}):=SSet_{\mathbb{V}}^{Aff_{\mathcal{C}}^{op}}$
(the model category of $\mathbb{V}$-simplicial presheaves
on $Aff_{\mathcal{C}}$) along the ($\mathbb{V}$-small) set of equivalences
of $Aff_{\mathcal{C}}$, and the homotopy category $\mathrm{Ho}(Aff_{\mathcal{C}}^{\wedge})$
will be naturally identified with the
full subcategory of $\mathrm{Ho}(SPr(Aff_{\mathcal{C}}))$ consisting of all
functors $F : Aff_{\mathcal{C}}^{op} \longrightarrow SSet_{\mathbb{V}}$
preserving weak equivalences. Objects in
$\mathrm{Ho}(Aff_{\mathcal{C}}^{\wedge})$ will be called \emph{pre-stacks}, and the derived simplicial $Hom$'s of the
simplicial model category $Aff_{\mathcal{C}}^{\wedge}$
will be denoted by $\mathbb{R}\underline{Hom}$.

We will fix once for all
a model pre-topology $\tau$ on $Aff_{\mathcal{C}}$ (Def. \ref{modtop}),
which induces a Grothendieck
topology on $\mathrm{Ho}(Aff_{\mathcal{C}})$, still denoted by the same symbol.
As recalled in \S \ref{modtop},
one can then consider a model category $Aff_{\mathcal{C}}^{\sim,\tau}$, of stacks
on the model site $(Aff_{\mathcal{C}}, \tau)$.
A morphism
$F \longrightarrow G$ of pre-stacks is an equivalence
in $Aff_{\mathcal{C}}^{\sim,\tau}$ if it induces isomorphisms
on all homotopy sheaves (for any choice of $X \in Aff_{\mathcal{C}}$ and any $s\in F(X)$).\\

To ease the notation we will write $\mathrm{St}(\mathcal{C},\tau)$\index{$\mathrm{St}(\mathcal{C},\tau)$} for the homotopy category $\mathrm{Ho}(Aff_{\mathcal{C}}^{\sim,\tau})$ of stacks.\\

The Bousfield localization construction yields an adjunction
$$a : \mathrm{Ho}(Aff_{\mathcal{C}}^{\wedge}) \longrightarrow \mathrm{St}(\mathcal{C},\tau)
\qquad \mathrm{Ho}(Aff_{\mathcal{C}}^{\wedge}) \longleftarrow \mathrm{St}(\mathcal{C},\tau) : j$$
where $j$ is fully faithful.

\begin{df}\label{dstack}
\begin{enumerate}
\item A \emph{stack} is an object
$F\in SPr(Aff_{\mathcal{C}})$ whose image
in $\mathrm{Ho}(Aff_{\mathcal{C}}^{\wedge})$ is in the essential image
of the functor $j$ above.
\item A \emph{morphism of stacks} is a morphism between stacks 
in $\mathrm{Ho}(SPr(Aff_{\mathcal{C}}))$, or equivalently in
$\mathrm{Ho}(Aff_{\mathcal{C}}^{\wedge})$, or equivalently in
$\mathrm{St}(\mathcal{C},\tau)$.
\item A morphism of stacks $f : F \longrightarrow G$
is a \emph{covering} (or a \emph{cover}) if the induced morphism of sheaves
$$\pi_{0}(f) : \pi_{0}(F) \longrightarrow  \pi_{0}(G)$$
is an epimorphism in the category of sheaves.

\end{enumerate}
\end{df}

We will always omit mentioning the functor $j$ and
consider the category $\mathrm{St}(\mathcal{C},\tau)$
as embedded in $\mathrm{Ho}(Aff_{\mathcal{C}}^{\wedge})$, and therefore
 embedded in $\mathrm{Ho}(SPr(Aff_{\mathcal{C}}))$.
With these conventions, the endofunctor $a$ of
$\mathrm{Ho}(Aff_{\mathcal{C}}^{\wedge})$ becomes
the \emph{associated stack functor}, which commutes with finite homotopy limits and arbitrary homotopy colimits.

A functor $F : Aff_{\mathcal{C}}^{op} \longrightarrow SSet_{\mathbb{V}}$
is a stack (\S \ref{remhagI})
if and only if it preserves equivalences
and possesses the descent property with respect to homotopy $\tau$-hypercovers. The derived simplicial
$Hom$'s in the model category
$Aff_{\mathcal{C}}^{\sim,\tau}$ of stacks will be denoted
by
$\mathbb{R}_{\tau}\underline{Hom}$. The natural morphism
$$\mathbb{R}\underline{Hom}(F,G)\longrightarrow
\mathbb{R}_{\tau}\underline{Hom}(F,G)$$
is an isomorphism in $\mathrm{Ho}(SSet)$ when $G$ is
a stack.

The model category $Aff_{\mathcal{C}}^{\sim,\tau}$ is a $t$-complete model
topos (Def. \ref{dmodtop}). We warn the reader that
neither of the model categories $Aff_{\mathcal{C}}^{\wedge}$ and $Aff_{\mathcal{C}}^{\sim,\tau}$
is right proper, though they are both left proper. Because of this certain care
has to be taken when considering homotopy pullbacks and more
generally homotopy limit constructions, as well as
comma model categories of objects over a base object. Therefore,
even when nothing is specified,
adequate fibrant replacement may have been chosen before
considering certain constructions. \\

The model Yoneda embedding (\S \ref{remhagI})
$$\underline{h} : Aff_{\mathcal{C}} \longrightarrow Aff_{\mathcal{C}}^{\wedge}$$
has a total right
derived functor
$$\mathbb{R}\underline{h} :
\mathrm{Ho}(Aff_{\mathcal{C}}) \longrightarrow \mathrm{Ho}(Aff_{\mathcal{C}}^{\wedge}),$$
which is fully faithful.
We also have a
naive Yoneda functor
$$h : Aff_{\mathcal{C}} \longrightarrow Aff_{\mathcal{C}}^{\wedge}$$
sending an object $X\in Aff_{\mathcal{C}}$ to the
presheaf of sets it represents (viewed as a simplicial
presheaf). With these notations, the Yoneda lemma reads
$$\mathbb{R}\underline{Hom}(\mathbb{R}\underline{h}_{X},F)\simeq
\mathbb{R}\underline{Hom}(h_{X},F)\simeq F(X)$$
for any fibrant object $F \in Aff_{\mathcal{C}}^{\wedge}$. The natural morphism $h_{X} \longrightarrow
\mathbb{R}\underline{h}_{X}$ is always
an isomorphism in $\mathrm{Ho}(Aff_{\mathcal{C}}^{\wedge})$ for any $X \in Aff_{\mathcal{C}}$.

We will also use the notation $\underline{Spec}\, A$ for
$\underline{h}_{Spec\, A}$. We warn the reader that
$Spec\, A$ lives in $Aff_{\mathcal{C}}$ whereas
$\underline{Spec}\, A$ is an object of
the model category of stacks $Aff_{\mathcal{C}}^{\sim,\tau}$. \\

We will assume the topology $\tau$ satisfies some
conditions. In order to state them, recall
the category $sAff_{\mathcal{C}}$ of simplicial objects
in $Aff_{\mathcal{C}}$ is a simplicial model category for the
Reedy model structure. Therefore, for any
object $X_{*} \in sAff_{\mathcal{C}}$ and any $\mathbb{U}$-small simplicial set
$K$ we can define an object $X_{*}^{\mathbb{R}\underline{K}}\in \mathrm{Ho}(sAff_{\mathcal{C}})$, by
first taking a Reedy fibrant model for $X_{*}$ and then
the exponential by $K$. The zero-th part of $X_{*}^{\mathbb{R}\underline{K}}$ will be
simply denoted by
$$X_{*}^{\mathbb{R}K}:=(X_{*}^{\mathbb{R}\underline{K}})_{0} \in \mathrm{Ho}(Aff_{\mathcal{C}}).$$
We also refer the reader to \cite[\S 4.4]{hagI} for more details and notations. \\

The following assumption on the pre-topology will be made.

\begin{ass}\label{ass5}
\begin{enumerate}
\item The topology $\tau$ on $\mathrm{Ho}(Aff_{\mathcal{C}})$ is quasi-compact. In other
words, for any covering family
$\{U_{i} \longrightarrow X\}_{i\in I}$ in $Aff_{\mathcal{C}}$ there exists a finite
subset $I_{0} \subset I$ such that the induced family
$\{U_{i} \longrightarrow X\}_{i\in I_{0}}$ is a covering.

\item For any finite family of objects $\{X_{i}\}_{i\in I}$ in $Aff_{\mathcal{C}}$ (including the empty family) 
the family of morphisms
$$\{X_{i} \longrightarrow \coprod^{\mathbb{L}}_{j\in I} X_{j} \}_{i\in I}$$
form a $\tau$-covering family of $\coprod^{\mathbb{L}}_{j\in I} X_{j}$.

\item Let $X_{*} \longrightarrow Y$ be an
augmented simplicial object in $Aff_{\mathcal{C}}$, corresponding
to a co-augmented co-simplicial object
$A \longrightarrow B_{*}$ in $Comm(\mathcal{C})$.
We assume that for any $n$, the one element family of morphisms
$$X_{n} \longrightarrow
X_{*}^{\mathbb{R}\partial \Delta^{n}}\times^{h}_{Y^{\mathbb{R}\partial \Delta^{n}}}Y$$
form a $\tau$-covering family in $Aff_{\mathcal{C}}$. Then the morphism
$$A \longrightarrow B_{*}$$
satisfies the descent condition in the sense
of Def. \ref{ddesc}.
\end{enumerate}
\end{ass}

The previous assumption has several consequences on the
homotopy theory of stacks. They are subsumed in the following
lemma.

\begin{lem}\label{lass5}
\begin{enumerate}
\item For any finite family of objects $X_{i}$ in $Aff_{\mathcal{C}}$ the natural
morphism
$$\coprod_{i}\mathbb{R}\underline{h}_{X_{i}}  \longrightarrow
\mathbb{R}\underline{h}_{\coprod^{\mathbb{L}}_{i}X_{i}}$$
is an equivalence in $Aff_{\mathcal{C}}^{\sim,\tau}$.
\item Let $H$ be the ($\mathbb{V}$-small) set of
augmented simplicial objects $X_{*} \longrightarrow Y$
in $Aff_{\mathcal{C}}$ such that for any $n$ the one element family of morphisms
$$X_{n} \longrightarrow
X_{*}^{\mathbb{R}\partial \Delta^{n}}\times^{h}_{Y^{\mathbb{R}\partial \Delta^{n}}}Y$$
is a $\tau$-covering family in $Aff_{\mathcal{C}}$. Then, the model category
$Aff_{\mathcal{C}}^{\sim,\tau}$ is the left Bousfield localization
of $Aff_{\mathcal{C}}^{\wedge}$ along the set of morphisms
$$|\mathbb{R}\underline{h}_{X_{*}}| \longrightarrow
\mathbb{R}\underline{h}_{Y} \qquad \coprod_{i}\mathbb{R}\underline{h}_{U_{i}}  \longrightarrow
\mathbb{R}\underline{h}_{\coprod^{\mathbb{L}}_{i}U_{i}}$$
where $X_{*} \rightarrow Y$ runs in $H$ and $\{U_{i}\}$ runs through the set
of all finite families of objects in $Aff_{\mathcal{C}}$.
\end{enumerate}
\end{lem}

\begin{proof} 
$(1)$ The case where the set of indices $I$ is empty
follows from our assumption \ref{ass5} $(2)$ with $I$ empty, 
as it states that the empty family covers the initial object on $Aff_{\mathcal{C}}$. 

Let us assume that the set of indices is not empty.
By induction, it is clearly enough to treat the case where
the finite family consists of two objects $X$ and $Y$.
Our assumption \ref{ass5} $(2)$ then implies that
the natural morphism
$$p : \mathbb{R}\underline{h}_{X}\coprod \mathbb{R}\underline{h}_{Y} \longrightarrow
\mathbb{R}\underline{h}_{X\coprod^{\mathbb{L}} Y}$$
is a covering. Therefore, $\mathbb{R}\underline{h}_{X\coprod^{\mathbb{L}} Y}$ is naturally
equivalent to the homotopy colimit of the
homotopy nerve of the morphism $p$. Using this remark we see that
it is enough to prove that
$$\mathbb{R}\underline{h}_{X}\times^{h}_{
\mathbb{R}\underline{h}_{X\coprod^{\mathbb{L}} Y}}
\mathbb{R}\underline{h}_{Y}\simeq \emptyset,$$
as then the homotopy nerve of $p$ will be a constant simplicial object
with values $\mathbb{R}\underline{h}_{X}\coprod \mathbb{R}\underline{h}_{Y}$. As the functor
$\mathbb{R}\underline{h}$ commutes with homotopy pullbacks, it is therefore enough to check that
$$A\otimes^{\mathbb{L}}_{A\times^{h} B}B\simeq 0,$$
for $A$ and $B$ two commutative monoids in $\mathcal{C}$ such that
$X=Spec\, A$ and $Y=Spec\,  B$ (here $0$ is the final object in $Comm(\mathcal{C})$).
For this we can of course suppose that
$A$ and $B$ are fibrant objects in $\mathcal{C}$.

We define a functor
$$F : A\times B-Comm(\mathcal{C}) \longrightarrow A-Comm(\mathcal{C}) \times B-Comm(\mathcal{C})$$
by the formula
$$F(C):=(C\otimes_{A\times B}A,C\otimes_{A\times B}B).$$
The functor $F$ is left Quillen for the product model structures on the right hand
side, and its right adjoint is given by
$G(C,D):=C\times D$
for any $(C,D) \in A-Comm(\mathcal{C}) \times B-Comm(\mathcal{C})$. For any $C \in A\times B-Comm(\mathcal{C})$, one has
$$C\simeq C\otimes^{\mathbb{L}}_{A\times B}\left(A\times B\right) \simeq
\left(C\otimes^{\mathbb{L}}_{A\times B}A\right) \times \left(C\otimes^{\mathbb{L}}_{A\times B}B\right),$$
because of our assumptions \ref{ass-1} and \ref{ass2}, which implies
that the adjunction morphism
$$C \longrightarrow \mathbb{R}G(\mathbb{L}F(C))$$
is an isomorphism in $\mathrm{Ho}(A\times B-Comm(\mathcal{C}))$. As the functor
$G$ reflects equivalences (because of our assumption \ref{ass-1}) this implies
that $F$ and $G$ form a Quillen equivalence. Therefore, the functor $\mathbb{R}G$ commutes with
homotopy push outs, and we have
$$A\otimes^{\mathbb{L}}_{A\times B}B\simeq
\mathbb{R}G(A,0)\otimes^{\mathbb{L}}_{\mathbb{R}G(A,B)}\mathbb{R}G(0,B)\simeq
\mathbb{R}G\left( (A,0)\coprod^{\mathbb{L}}_{(A,B)}(0,B)\right)\simeq
\mathbb{R}G(0)\simeq 0.$$

$(2)$ We know by \cite{hagI} that $Aff_{\mathcal{C}}^{\sim,\tau}$ is the
left Bousfield localization of $Aff_{\mathcal{C}}^{\wedge}$ along the set of
morphisms $|F_{*}| \longrightarrow h_{X}$, where $F_{*}
\longrightarrow h_{X}$ runs in a certain $\mathbb{V}$-small set of
$\tau$-hypercovers. Recall that for each hypercover $F_{*}
\longrightarrow h_{X}$ in this set, each simplicial presheaf
$F_{n}$ is a coproduct of some $h_{U}$. Using the
quasi-compactness assumption \ref{ass5} $(1)$ one sees immediately
that one can furthermore assume that each $F_{n}$ is a finite
coproduct of some $h_{U}$. Finally, using the part $(1)$ of the
present lemma we see that the descent condition
of \cite{hagI} can be stated as two distinct conditions, one
concerning finite coproducts and the other one concerning
representable hypercovers. From this we deduce part $(2)$
of the lemma. \end{proof}

Lemma \ref{lass5} $(2)$ can be reformulated as follows.

\begin{cor}\label{class5}
A simplicial presheaf
$$F : Comm(\mathcal{C}) \longrightarrow SSet_{\mathbb{V}}$$
is a stack if and only if it satisfies the following
three conditions.
\begin{itemize}

\item For any equivalence $A\longrightarrow B$  in $Comm(\mathcal{C})$
the induced morphism $F(A) \longrightarrow F(B)$
is an equivalence of simplicial sets.

\item For any finite family of commutative monoids $\{A\}_{i\in I}$ in $\mathcal{C}$ (including the
empty family), the natural morphism
$$F(\prod_{i\in I} {}^{h}A_{i}) \longrightarrow \prod_{i\in I}F(A_{i})$$
is an isomorphism in $\mathrm{Ho}(SSet)$.

\item For any co-simplicial commutative $A$-algebra
$A\longrightarrow B_{*}$, corresponding to
a $\tau$-hypercover
$$Spec\, B_{*} \longrightarrow Spec\, A$$
in $Aff_{\mathcal{C}}$, the induced morphism
$$F(A) \longrightarrow Holim_{[n]\in \Delta}F(B_{n})$$
is an isomorphism in $\mathrm{Ho}(SSet)$.

\end{itemize}
\end{cor}

Another important consequence of lemma \ref{lass5}
is the following.

\begin{cor}\label{c0'}
The model pre-topology $\tau$ on $Aff_{\mathcal{C}}$ is sub-canonical in the
sense of Def. \ref{sottocan}.
\end{cor}

\begin{proof} We need to show that for any $Z \in Aff_{\mathcal{C}}$ the
object $G:=\mathbb{R}\underline{h}_{Z}$ is a stack, or in other
words is a local object in $Aff_{\mathcal{C}}^{\sim,\tau}$.
For this, we use our lemma \ref{lass5} $(2)$.
The descent property for
finite coproducts is obviously satisfied because of the Yoneda
lemma.
Let $X_{*} \longrightarrow Y$ be a
simplicial object
in $Aff_{\mathcal{C}}$ such that
$$\mathbb{R}\underline{h}_{X_{*}} \longrightarrow
\mathbb{R}\underline{h}_{Y}$$
is a $\tau$-hypercover.
By lemma \ref{ldesc}
the natural morphism
$$Hocolim_{n}X_{n} \longrightarrow Y$$
is an isomorphism in $\mathrm{Ho}(Aff_{\mathcal{C}})$. Therefore,
the Yoneda lemma implies that one has
$$\mathbb{R}\underline{Hom}(\underline{h}_{Y},G)
\simeq Map(Y,Z)\simeq Holim_{n}Map(X_{n} ,Z)\simeq
Holim_{n}\mathbb{R}\underline{Hom}(F_{*},G),$$
showing that $G$ is a stack. \end{proof}

The corollary \ref{c0'} implies that
$\mathbb{R}\underline{h}$ provides a
fully faithful functor
$$\mathbb{R}\underline{h} : \mathrm{Ho}(Aff_{\mathcal{C}}) \longrightarrow
\mathrm{St}(\mathcal{C},\tau).$$
Objects  in the essential image
of $\mathbb{R}\underline{h}$ will be called
\emph{representable objects}. If such an object
corresponds to a commutative monoid $A \in \mathrm{Ho}(Comm(\mathcal{C}))$, it will
also be denoted by $\mathbb{R}\underline{Spec}\, A\in \mathrm{St}(\mathcal{C},\tau)$.
In formula
$$\mathbb{R}\underline{Spec}\, A:=\mathbb{R}\underline{h}_{Spec\, A},$$
for any $A\in Comm(\mathcal{C})$ corresponding to $Spec\, A\in Aff_{\mathcal{C}}$.
As $\mathbb{R}\underline{h}$ commutes with
$\mathbb{U}$-small homotopy limits,
we see that the subcategory of representable
stacks is stable by $\mathbb{U}$-small
homotopy limits. The reader should be careful that
a $\mathbb{V}$-small homotopy limit of
representable stacks is not representable in general.
Lemma \ref{lass5} $(1)$ also implies that
a finite coproduct of representable stacks is a representable stack, and we have
$$\coprod_{i}\mathbb{R}\underline{h}_{U_{i}}\simeq \mathbb{R}\underline{h}_{\coprod^{\mathbb{L}}_{i}U_{i}}.$$

Also, by identifying the category
$\mathrm{Ho}(Comm(\mathcal{C}))^{op}$ with the full subcategory of $\mathrm{St}(\mathcal{C},\tau)$
consisting of representable stacks, one can extend the notions
of morphisms defined in \S \ref{partI.2} (e.g. (formally)
\'etale, Zariski open immersion, flat, smooth \dots) to morphisms between
representable stacks. Indeed, they are all invariant by equivalences and therefore
are properties of morphisms in $\mathrm{Ho}(Comm(\mathcal{C}))$. We will often use
implicitly these
extended notions. In particular, we will use the expression
\emph{$\tau$-covering families of representable stacks} to
denote families of morphisms of representable stacks corresponding
in $\mathrm{Ho}(Comm(\mathcal{C}))^{op}$ to $\tau$-covering families.

We will use the same terminology for the morphisms in $\mathrm{Ho}(Comm(\mathcal{C}))$ and
for the corresponding morphisms of representable stacks, except for
the notion of epimorphism (Def. \ref{d5}), which for obvious reasons 
will be replaced by \emph{monomorphism} in the context of stacks. This is justified since
 a morphism $A \longrightarrow B$ is an epimorphism in the sense
of Def. \ref{d5} if and only if the induced morphism of stacks
$$\mathbb{R}\underline{Spec}\, B \longrightarrow \mathbb{R}\underline{Spec}\, A$$
is a monomorphism in the model category $Aff_{\mathcal{C}}^{\sim,\tau}$ (see Remark \ref{explepi}). 

\begin{rmk} 
\emph{The reader should be warned that we will also
use the expression \emph{epimorphism of stacks}, which will
refer to a morphism of stacks that induces
an epimorphism on the sheaves $\pi_{0}$ (see Def. \ref{dstack0} or \cite{hagI}, where  
they are also called \emph{coverings}).
It is important to notice that
a $\tau$-covering family of representable
stacks $\{X_{i} \longrightarrow X\}$ induces
an epimorphism of stacks
$\coprod X_{i} \longrightarrow X$. On the contrary,
there might very well exist families of morphisms
of representable stacks $\{X_{i} \longrightarrow X\}$
such that $\coprod X_{i} \longrightarrow X$ is an epimorphism
of stacks, but which are not
$\tau$-covering families (e.g.
a morphism between representable stacks that admits a section). }
\end{rmk}

\begin{cor}\label{c0''}
Let $\{u_{i} : X_{i}=Spec\, A_{i} \longrightarrow X=Spec\, A\}_{i\in I}$ be a covering family in
$Aff_{\mathcal{C}}$. Then, the family of base change functors
$$\{\mathbb{L}u_{i}^{*} : \mathrm{Ho}(A-Mod) \longrightarrow \mathrm{Ho}(A_{i}-Mod)\}_{i\in I}$$
is conservative.
In other words, a $\tau$-covering family in $Aff_{\mathcal{C}}$ is
a formal covering in the sense of Def. \ref{dcov}.
\end{cor}

\begin{proof} By the quasi-compactness assumption on $\tau$ we can assume that
the covering family is finite. Also, the morphism
$A \longrightarrow \prod^{h}_{i}A_{i}=B$ is a covering. Therefore,
the descent assumption \ref{ass5} $(3)$ implies that
the base change functor
$$B\otimes^{\mathbb{L}}_{A} - : \mathrm{Ho}(A-Mod) \longrightarrow \mathrm{Ho}(B-Mod)$$
is conservative. Finally, we have seen during the proof of \ref{lass5} $(1)$ that
the product of the base change functors
$$\prod_{i}A_{i}\otimes^{\mathbb{L}}_{B}- : \mathrm{Ho}(B-Mod) \longrightarrow \prod_{i}\mathrm{Ho}(A_{i}-Mod)$$
is an equivalence. Therefore, the composition
$$\prod_{i}A_{i}\otimes^{\mathbb{L}}_{A}- : \mathrm{Ho}(A-Mod) \longrightarrow \prod_{i}\mathrm{Ho}(A_{i}-Mod)$$
is a conservative functor. \end{proof}

We also
recall the Yoneda lemma for stacks, stating that for
any $A\in Comm(\mathcal{C})$, and any fibrant object $F \in Aff_{\mathcal{C}}^{\sim,\tau}$,
there is a natural equivalence of simplicial sets
$$\mathbb{R}\underline{Hom}(\mathbb{R}\underline{Spec}\, A,F)
\simeq \mathbb{R}_{\tau}\underline{Hom}(\mathbb{R}\underline{Spec}\, A,F)\simeq
F(A).$$
For an object $F \in \mathrm{St}(\mathcal{C},\tau)$,
and $A \in Comm(\mathcal{C})$ we will use the
following notation
$$\mathbb{R}F(A):=\mathbb{R}_{\tau}\underline{Hom}(\mathbb{R}\underline{Spec}\, A,F).$$
Note that $\mathbb{R}F(A)\simeq RF(A)$, where $R$ is a fibrant replacement
functor on $Aff_{\mathcal{C}}^{\sim,\tau}$. Note that there is always a natural
morphism $F(A) \longrightarrow \mathbb{R}F(A)$, which is
an equivalence precisely when $F$ is a stack. \\

Finally, another
important consequence of assumption \ref{ass5} is the
local character of representable stacks.

\begin{prop}\label{c0}
Let $G$ be a representable stack and
$F \longrightarrow G$ be any morphism. Assume
there exists a $\tau$-covering family of representable stacks
$$\{G_{i} \longrightarrow G\},$$
such that each
stack $F\times^{h}_{G}G_{i}$
is representable. Then $F$ is a representable stack.
\end{prop}

\begin{proof}
Let $X \in Aff_{\mathcal{C}}$ be a fibrant object such that
$G\simeq \underline{h}_{X}$. We can of course
assume that $G=\underline{h}_{X}$. We can also assume
that $F \longrightarrow G$ is a fibration, and therefore
that $G$ and $F$ are fibrant objects.

By choosing a
refinement of the covering family $\{G_{i} \longrightarrow G\}$,
one can suppose that the covering family is finite and that
each morphism $G_{i} \longrightarrow G$
is the image by $\underline{h}$ of
a fibration $U_{i} \longrightarrow X$ in $Aff_{\mathcal{C}}$.
Finally, considering the
coproduct $U=\coprod^{\mathbb{L}} U_{i} \longrightarrow X$
in $Aff_{\mathcal{C}}$ and using lemma \ref{lass5} $(1)$ one can
suppose that the family $\{G_{i} \longrightarrow G\}$
has only one element and is the image by
$\underline{h}$ of a fibration
$U \longrightarrow X$  in $Aff_{\mathcal{C}}$.

We consider the augmented simplicial object $U_{*} \longrightarrow X$
in $Aff_{\mathcal{C}}$, which is the nerve of the morphism $U \rightarrow X$, and
the corresponding augmented simplicial object
$\underline{h}_{U_{*}} \longrightarrow
\underline{h}_{X}$ in $Aff_{\mathcal{C}}^{\sim,\tau}$.
We form the pullback in $Aff_{\mathcal{C}}^{\sim,\tau}$
$$\xymatrix{
F \ar[r] & \underline{h}_{X} \\
F_{*} \ar[r] \ar[u] & \underline{h}_{U_{*}}\ar[u],}$$
which is a homotopy pullback because of our choices.
In particular, for any $n$, $F_{n}$ is a representable
stack.

Clearly $F_{*}$ is the nerve of the
fibration $F\times_{\underline{h}_{X}}\underline{h}_{U}
\longrightarrow F$. As this last morphism is
an epimorphism in $Aff_{\mathcal{C}}^{\sim,\tau}$ the natural morphism
$$|F_{*}|:=Hocolim_{n} F_{n} \longrightarrow F$$
is an isomorphism in $\mathrm{St}(\mathcal{C},\tau)$. Therefore,
it remains to show that $|F_{*}|$ is a representable
stack.

We will consider the category $sAff_{\mathcal{C}}^{\sim,\tau}$of simplicial
objects in $Aff_{\mathcal{C}}^{\sim,\tau}$, endowed with its
Reedy model structure (see \cite[\S 4.4]{hagI} for details
and notations). In the same way, we will consider
the Reedy model structure on $sAff_{\mathcal{C}}$, the category
of simplicial objects in $Aff_{\mathcal{C}}$.

\begin{lem}\label{lc0}
There exists a simplicial object $V_{*}$ in $Aff_{\mathcal{C}}$ and
an isomorphism $\mathbb{R}\underline{h}_{V_{*}} \simeq F_{*}$
in $\mathrm{Ho}(sAff_{\mathcal{C}}^{\sim,\tau})$.
\end{lem}

\begin{proof}
First of all, our functor $\underline{h}$
extends in the obvious way to a functor on the categories of
simplicial objects
$$\underline{h} : sAff_{\mathcal{C}} \longrightarrow sAff_{\mathcal{C}}^{\sim,\tau},$$
by the formula
$$(\underline{h}_{X_{*}})_{n}:=\underline{h}_{X_{n}}$$
for any $X_{*} \in sAff_{\mathcal{C}}$.
This functor possesses a right derived functor
$$\mathbb{R}\underline{h} : \mathrm{Ho}(sAff_{\mathcal{C}})
\longrightarrow \mathrm{Ho}(sAff_{\mathcal{C}}^{\sim,\tau}),$$
which is easily seen to be fully faithful.

We claim that the essential image of $\mathbb{R}\underline{h}$ consists
of all simplicial objects $F_{*} \in \mathrm{Ho}(sAff_{\mathcal{C}}^{\sim,\tau})$
such that for any $n$, $F_{n}$ is a representable stack. This
will obviously imply the lemma.
Indeed, as $\mathbb{R}\underline{h}$ commutes with
homotopy limits, one sees that this essential
image is stable by ($\mathbb{U}$-small) homotopy limits. Also, any
object $F_{*} \in sAff_{\mathcal{C}}^{\sim,\tau}$ can be written
as a homotopy limit
$$F_{*}\simeq Holim_{n}\mathbb{R}Cosk_{n}(F_{*})$$
of its derived coskeleta (see \cite[\S 4.4]{hagI}). Recall that
for a fibrant object $F_{*}$ in $sAff_{\mathcal{C}}^{\sim,\tau}$
one has
$$Cosk_{n}(F_{*})_{p}\simeq \mathbb{R}Cosk_{n}(F_{*})_{p}\simeq
(F_{*})^{Sk_{n}\Delta^{p}},$$
where, for a simplicial set $K$, $(F_{*})^{K}$ is defined to be
the equalizer of the two natural maps
$$\prod_{[n]}(F_{n})^{K_{n}} \rightrightarrows
\prod_{[n] \rightarrow [m]}(F_{m})^{K_{n}}.$$
Now, for any simplicial set $K$, and any integer $n\geq 0$, there is
 a homotopy push out square
$$\xymatrix{
Sk_{n}K \ar[r] & Sk_{n+1}K \\
\coprod_{K^{\partial \Delta^{n+1}}}\partial \Delta^{n+1} \ar[u] \ar[r] &
\coprod_{K^{\Delta^{n+1}}}\Delta^{n+1}. \ar[u]}$$
Using this, and the fact that $K \mapsto F_{*}^{K}$
sends homotopy push outs to homotopy pullbacks
when $F_{*}$ is fibrant,
we see that for any fibrant object
$F_{*} \in sAff_{\mathcal{C}}^{\sim,\tau}$, and any $n$, we have a homotopy pullback
diagram in $sAff_{\mathcal{C}}^{\sim,\tau}$
$$\xymatrix{
\mathbb{R}Cosk_{n}(F_{*}) \ar[r] & \mathbb{R}Cosk_{n-1}(F_{*})\\
A^{n}_{*} \ar[r] \ar[u] &
B^{n}_{*}. \ar[u]}
$$
Here, $A^{n}_{*}$ and $B^{n}_{*}$ are defined by the following formulas
$$
\begin{array}{cccc}
A^{n}_{*}  : & \Delta^{op} & \longrightarrow & Aff_{\mathcal{C}}^{\sim,\tau} \\
 & [p] & \mapsto & \prod_{(\Delta^{p})^{\Delta^{n+1}}}F_{*}^{\Delta^{n+1}}
\end{array}$$

$$
\begin{array}{cccc}
B^{n}_{*}  : & \Delta^{op} & \longrightarrow & Aff_{\mathcal{C}}^{\sim,\tau} \\
 & [p] & \mapsto & \prod_{(\Delta^{p})^{\partial\Delta^{n+1}}}
F_{*}^{\partial \Delta^{n+1}}.
\end{array}$$

Therefore, by induction on $n$, it is enough to see that
if $F_{*}$ is fibrant then  
$Cosk_{0}(F_{*})$, $A^{n}_{*}$ and
$B^{n}_{*}$ all belongs to the essential image of 
$\mathbb{R}\underline{h}$.

The simplicial object $Cosk_{0}(F_{*})$ is isomorphic
to the nerve of $F_{0} \longrightarrow *$,
and therefore is the image by $\mathbb{R}\underline{h}$ of the
nerve of a fibration $X \longrightarrow *$ in $Aff_{\mathcal{C}}$ representing
$F_{0}$. This shows that $Cosk_{0}(F_{*})$ is in the
image of $\mathbb{R}\underline{h}$.

We have $F_{*}^{\Delta^{n+1}}=F_{n+1}$, which is a representable
stack. Let $X_{n+1} \in Aff_{\mathcal{C}}$ be a fibrant object
such that $F_{n+1}$ is equivalent to $\underline{h}_{X_{n+1}}$. Then,
as $\underline{h}$ commutes with limits, we see that $A^{n}_{*}$
is equivalent the image by $\underline{h}$ of the simplicial object
$$[p] \mapsto \prod_{(\Delta^{p})^{\Delta^{n+1}}}X_{n+1}.$$

In the same way, $F_{*}^{\partial \Delta^{n+1}}$ can be written
as a finite homotopy limit of $F_{i}$'s, and therefore
is a representable stack. Let $Y_{n+1}$ be a fibrant
object in $Aff_{\mathcal{C}}$ such that
$F_{*}^{\partial \Delta^{n+1}}$ is equivalent to $\underline{h}_{Y_{n+1}}$. Then,
$B^{n}_{*}$ is equivalent to the image
by $\underline{h}$ of the
simplicial object
$$[p] \mapsto \prod_{(\Delta^{p})^{\partial \Delta^{n+1}}}Y_{n+1}.$$
This proves the lemma. \end{proof}

We now finish the proof of Proposition \ref{c0}. Let $V_{*}$
be a simplicial object in $sAff_{\mathcal{C}}$ such that
$F_{*}\simeq \mathbb{R}\underline{h}_{V_{*}}$. The augmentation
$F_{*} \longrightarrow \mathbb{R}\underline{h}_{U_{*}}$
gives rise to a well defined morphism in $\mathrm{Ho}(sAff_{\mathcal{C}}^{\sim,\tau})$
$$q : \mathbb{R}\underline{h}_{V_{*}} \longrightarrow \mathbb{R}\underline{h}_{U_{*}}.$$
We can of course suppose that $V_{*}$ is a cofibrant
object in $sAff_{\mathcal{C}}$. As $U_{*}$ is the nerve
of a fibration between fibrant objects it is
fibrant in $sAff_{\mathcal{C}}$. Therefore,
as $\mathbb{R}\underline{h}$ is fully faithful, we can represent
$q$, up to an isomorphism, as
the image by $\underline{h}$ of a morphism in $sAff_{\mathcal{C}}$
$$r : V_{*} \longrightarrow U_{*}.$$
On the level of commutative monoids, the morphism
$r$ is given by a morphism $B_{*} \longrightarrow C_{*}$  of co-simplicial
objects in $Comm(\mathcal{C})$. By construction, for any
morphism $[n] \rightarrow [m]$ in $\Delta$,
the natural morphism
$$F_{m} \longrightarrow F_{n}\times^{h}_{\mathbb{R}\underline{h}_{U_{n}}}\mathbb{R}\underline{h}_{U_{m}}$$
is an isomorphism in $\mathrm{St}(\mathcal{C},\tau)$. This implies that
the underlying co-simplicial $B_{*}$-module of $C_{*}$ is
homotopy cartesian in the sense of Def. \ref{ddesc}.
Our
assumption \ref{ass5} $(3)$ implies that if $Y:=Hocolim_{n}V_{*} \in Aff_{\mathcal{C}}$, then
the natural morphism
$$V_{*} \longrightarrow U_{*}\times^{h}_{X}Y$$
is an isomorphism in $\mathrm{Ho}(sAff_{\mathcal{C}})$. As homotopy colimits
in $Aff_{\mathcal{C}}^{\sim,\tau}$ commute with homotopy pullbacks, this
implies that
$$|F_{*}|\simeq |\mathbb{R}\underline{h}_{V_{*}}|\simeq
|\mathbb{R}\underline{h}_{U_{*}}|\times^{h}_{\mathbb{R}\underline{h}_{X}}
\mathbb{R}\underline{h}_{Y}.$$
But, as $\mathbb{R}\underline{h}_{U_{*}} \longrightarrow \mathbb{R}\underline{h}_{X}$
is the homotopy nerve of an epimorphism we have
$$|\mathbb{R}\underline{h}_{U_{*}}|\simeq \mathbb{R}\underline{h}_{X},$$
showing finally that
$|F_{*}|$ is isomorphic to $\mathbb{R}\underline{h}_{Y}$
and therefore is a representable stack. \end{proof}

Finally, we finish this first section by the following
description of the comma model category $Aff_{\mathcal{C}}^{\sim,\tau}/\underline{h}_{X}$,
for some fibrant object $X\in Aff_{\mathcal{C}}$. This is not a completely trivial
task as the model category $Aff_{\mathcal{C}}^{\sim,\tau}$ is not right proper.

For this, let $A\in Comm(\mathcal{C})$ such that $X=Spec\, A$, so that
$A$ is a cofibrant object in $Comm(\mathcal{C})$. We consider the
comma model category $A-Comm(\mathcal{C})=(Aff_{\mathcal{C}}/X)^{op}$, which is also
the model category $Aff_{A-Mod}^{op}=Comm(A-Mod)$.
The model pre-topology $\tau$ on $Aff_{\mathcal{C}}$ induces in a natural
way a model pre-topology $\tau$ on
$Aff_{\mathcal{C}}/X=Aff_{A-Mod}$. Note that there exists a natural
equivalence of categories, compatible with the
model structures, between
$(Aff_{\mathcal{C}}/X)^{\sim,\tau}$ and $Aff_{\mathcal{C}}^{\sim,\tau}/h_{X}$.
We consider the natural
morphism $h_{X} \longrightarrow \underline{h}_{X}$.
It gives rise to a Quillen adjunction
$$Aff_{\mathcal{C}}^{\sim,\tau}/h_{X} \longrightarrow Aff_{\mathcal{C}}^{\sim,\tau}/\underline{h}_{X}
\qquad
Aff_{\mathcal{C}}^{\sim,\tau}/h_{X} \longleftarrow Aff_{\mathcal{C}}^{\sim,\tau}/\underline{h}_{X}$$
where the right adjoint sends $F \longrightarrow \underline{h}_{X}$
to $F\times_{\underline{h}_{X}}h_{X} \longrightarrow h_{X}$.

\begin{prop}\label{pcomma}
The above Quillen adjunction induces a Quillen equivalence
between $Aff_{\mathcal{C}}^{\sim,\tau}/h_{X}\simeq (Aff_{\mathcal{C}}/X)^{\sim,\tau}$
and $Aff_{\mathcal{C}}^{\sim,\tau}/\underline{h}_{X}$.
\end{prop}

\begin{proof} For $F \longrightarrow \underline{h}_{X}$
a fibrant object, and $Y\in Aff_{\mathcal{C}}/X$, there is
a homotopy cartesian square
$$\xymatrix{
(F\times_{\underline{h}_{X}}h_{X})(Y) \ar[r] \ar[d] &
h_{X}(Y) \ar[d] \\
F(Y) \ar[r] & \underline{h}_{X}(Y).}$$
As the morphism $h_{X}(Y) \longrightarrow \underline{h}_{X}(Y)$
is always surjective up to homotopy when $Y$ is cofibrant, this
implies easily that the derived pullback functor
$$\mathrm{Ho}(Aff_{\mathcal{C}}^{\sim,\tau}/\underline{h}_{X}) \longrightarrow
\mathrm{Ho}(Aff_{\mathcal{C}}^{\sim,\tau}/h_{X})$$
is conservative. Therefore, it is enough to show that the
forgetful functor
$$\mathrm{Ho}(Aff_{\mathcal{C}}/X^{\sim,\tau}) \longrightarrow
\mathrm{Ho}(Aff_{\mathcal{C}}^{\sim,\tau}/\underline{h}_{X})$$
is fully faithful.

Using the Yoneda lemma for $Aff_{\mathcal{C}}/X$, we have
$$Map_{(Aff_{\mathcal{C}}/X)^{\sim,\tau}}(h_{Y},h_{Z})\simeq
Map_{Aff_{\mathcal{C}}/X}(Y,Z)$$
for any two objects $Y$ and $Z$ in $Aff_{\mathcal{C}}/X$. Therefore,
there exists a natural fibration sequence
$$\xymatrix{
Map_{(Aff_{\mathcal{C}}/X)^{\sim,\tau}}(h_{Y},h_{Z}) \ar[r] & Map_{Aff_{\mathcal{C}}}(Y,Z) \ar[r] &
Map_{Aff_{\mathcal{C}}}(Z,X).}$$
In the same way, the Yoneda lemma for $Aff_{\mathcal{C}}$ implies that there exists
a fibration  sequence
$$\xymatrix{
Map_{Aff_{\mathcal{C}}^{\sim,\tau}/\underline{h}_{X}}(h_{Y},h_{Z}) \ar[r] &
Map_{Aff_{\mathcal{C}}}(Y,Z) \ar[r] &
Map_{Aff_{\mathcal{C}}}(Z,X).}$$
This two fibration sequences implies that the forgetful functor
induces equivalences of simplicial sets
$$Map_{Aff_{\mathcal{C}}/X^{\sim,\tau}}(h_{Y},h_{Z}) \simeq
Map_{Aff_{\mathcal{C}}^{\sim,\tau}/\underline{h}_{X}}(Y,Z).$$
In other words, the functor
$$\mathrm{Ho}(Aff_{\mathcal{C}}/X^{\sim,\tau}) \longrightarrow
\mathrm{Ho}(Aff_{\mathcal{C}}^{\sim,\tau}/\underline{h}_{X})$$
is fully faithful when restricted to
the full subcategory of representable stacks.
But, any object in $\mathrm{Ho}(Aff_{\mathcal{C}}/X^{\sim,\tau})$ is a
homotopy colimit of representable stacks. Furthermore,
as the derived pullback
$$\mathrm{Ho}(Aff_{\mathcal{C}}^{\sim,\tau}/\underline{h}_{X}) \longrightarrow
\mathrm{Ho}(Aff_{\mathcal{C}}^{\sim,\tau}/h_{X})$$
commutes with homotopy colimits (as homotopy pullbacks
of simplicial sets do), this implies that the
functor
$$\mathrm{Ho}(Aff_{\mathcal{C}}/X^{\sim,\tau}) \longrightarrow
\mathrm{Ho}(Aff_{\mathcal{C}}^{\sim,\tau}/\underline{h}_{X})$$
is fully faithful on the whole category. \end{proof}

The important consequence of Prop. \ref{pcomma}
comes from the fact that it allows to see
objects in $\mathrm{Ho}(Aff_{\mathcal{C}}^{\sim,\tau}/\underline{h}_{X})$
as functors
$$A-Comm(\mathcal{C}) \longrightarrow SSet.$$
This last fact will be used implicitly in the sequel of this work. \\

As explained in the introduction to this chapter, we will need to fix
a class \textbf{P} of morphisms in $Aff_{\mathcal{C}}$.  Such a class will be then used to glue 
representable stacks to get a 
\emph{geometric stack}. In other words, geometric stacks will 
be the objects obtained by taking some
quotient of representable stacks by equivalence relations 
whose structural morphisms are in  \textbf{P}. Of course, different choices of  \textbf{P} will lead to 
different notions of geometric stacks. To fix his intuition 
the reader may think of  \textbf{P} as being the class of smooth morphisms
(though in some applications \textbf{P} can be something different). \\

From now, and all along this section, we fix a class \textbf{P} of
morphism in $Aff_{\mathcal{C}}$,
that is stable by equivalences. 
As the Yoneda functor
$$\mathbb{R}\underline{h} : \mathrm{Ho}(Aff_{\mathcal{C}})^{op} \longrightarrow
\mathrm{St}(\mathcal{C},\tau)$$
is fully faithful
we can extend the notion of morphisms belonging to \textbf{P}
to its essential image. So, a morphism of representable
objects in $\mathrm{Ho}(SPr(Aff_{\mathcal{C}}^{\sim,\tau}))$ is in \textbf{P} if
by definition it correspond to a morphism in $\mathrm{Ho}(Aff_{\mathcal{C}})$ which is
in \textbf{P}. We will make the following
assumtions on morphisms of \textbf{P} with respect to the topology $\tau$, 
making ``being in \textbf{P}'' into a $\tau$-local property.

\begin{ass}\label{ass4}
\begin{enumerate}

\item Covering families consist of morphisms in \textbf{P} i.e. for any
$\tau$-covering family $\{U_{i} \longrightarrow X\}_{i\in I}$ in 
$Aff_{\mathcal{C}}$, the morphism
$U_{i} \longrightarrow X$ is in \textbf{P} for all $i\in I$. 

\item Morphisms in \textbf{P} are stable by compositions,
equivalences and homotopy pullbacks.

\item Let $f : X \longrightarrow Y$ be a morphism
in $Aff_{\mathcal{C}}$.
If there exists
a $\tau$-covering family
$$\{U_{i}\longrightarrow X\}$$
such that 
each composite morphism $U_{i} \longrightarrow Y$
lies in \textbf{P}, then $f$ belongs to \textbf{P}.

\item For any two objects $X$ and $Y$ in $Aff_{\mathcal{C}}$, the
two natural morphisms
$$X \longrightarrow X\coprod^{\mathbb{L}} Y \qquad
Y \longrightarrow X\coprod^{\mathbb{L}} Y$$
are in \textbf{P}.

\end{enumerate}
\end{ass}

The reader will notice that assumptions \ref{ass5} and \ref{ass4} together imply the
following useful fact.

\begin{lem}\label{lass4}
Let $\{X_{i} \longrightarrow X\}$ be a finite family of morphisms in
\textbf{P}. The total morphism
$$\coprod^{\mathbb{L}}_{i}X_{i} \longrightarrow X$$
is also in \textbf{P}.
\end{lem}

\begin{proof} We consider the family of natural morphisms
$$\{X_{j} \longrightarrow \coprod^{\mathbb{L}}_{i}X_{i}\}_{j}.$$
According to our assumption \ref{ass5} it is a $\tau$-covering family in
$Aff_{\mathcal{C}}$. Moreover, each morphism $X_{j} \longrightarrow
\coprod^{\mathbb{L}}_{i}X_{i}$ and $X_{j} \longrightarrow X$ is in \textbf{P}, so
assumption \ref{ass4} $(3)$ implies that so is
$\coprod^{\mathbb{L}}_{i}X_{i} \longrightarrow X$. \end{proof}

We finish this section by the following definition.

\begin{df}\label{dhag}
A \emph{Homotopical Algebraic Geometry context}
(or simply \emph{HAG context})\index{HAG context} is
a $5$-tuple $(\mathcal{C},\mathcal{C}_{0},\mathcal{A},\tau,\textbf{P})$,
where $(\mathcal{C},\mathcal{C}_{0},\mathcal{A})$ is
a HA context in the sense of Def. \ref{dha}, $\tau$
is a model pre-topology on $Aff_{\mathcal{C}}$, \textbf{P}
is a class of morphism in $Aff_{\mathcal{C}}$, and such that
assumptions \ref{ass5} and \ref{ass4} are satisfied.
\end{df}

\section{Main definitions and standard properties}

From now on, we fix a HAG context $(\mathcal{C},\mathcal{C}_{0},\mathcal{A},\tau,\textbf{P})$
in the sense of Def. \ref{dhag}. We will consider
the model category of stacks $Aff_{\mathcal{C}}^{\sim,\tau}$
as described in the previous section, and introduce the notion
of \textit{geometric} and $n$\textit{-geometric stacks}, which are objects in $Aff_{\mathcal{C}}^{\sim,\tau}$
satisfying certain properties.

The basic geometric idea is that a stack is geometric if
it is obtained by taking the quotient of a representable stack $X$ (or
more generally of a disjoint union of representable stacks) by the action of 
a groupoid object $X_{1}$ acting on $X$, such that $X_{1}$ is itself representable, 
and such that the source morphism $X_{1} \longrightarrow X$ is a morphism in the chosen class \textbf{P}. 
If one thinks of \textbf{P} as being the class of certain smooth morphisms, 
being geometric is thus equivalent of being a quotient by a smooth groupoid action.

It turns out that this notion is not enough for certain applications, as some
naturally arising stacks are obtained as quotients by an action of a groupoid in geometric stacks
rather than in representable stacks (e.g. the quotients by a group-stack action). 
We will therefore also introduce the notion of $n$-geometric stack, 
which is defined inductively as a stack obtained as a quotient by an action of
a groupoid object in $(n-1)$-geometric stacks whose source morphism
is in \textbf{P}. Of course, for this
definition to make sense one must also explain, inductively on $n$, what are the morphisms in \textbf{P} between
$n$-geometric stacks. 

The inductive definition we give below uses a different (though equivalent) point of view, closer to the 
original definition 
of algebraic stacks due to Deligne-Mumford and M. Artin. It says that 
a stack $F$ is $n$-geometric if for any pair of points of $F$, the stack of equivalences between them
is $(n-1)$-geometric, and if moreover it receives a morphism, which is surjective and in \textbf{P},  
from a representable stack (or from a disjoint union of representable stacks).
The equivalence of this definition with the previously mentioned quotient-by-groupoids point of view will be
established in the next section (see Prop. \ref{p13}). \\

\begin{df}\label{d11}
\begin{enumerate}
\item A stack is \emph{$(-1)$-geometric} if it is
representable.
\item
A morphism of stacks $f : F \longrightarrow G$ is
\emph{$(-1)$-representable} if for any representable
stack $X$ and any morphism $X\longrightarrow G$, the
homotopy pullback $F\times^{h}_{G}X$ is a representable
stack.
\item
A morphism of stacks $f : F \longrightarrow G$ is
\emph{in $(-1)$-\textbf{P}} if it is
$(-1)$-representable, and if for any representable
stack $X$ and any morphism $X\longrightarrow G$,
the induced morphism
$$F\times^{h}_{G}X \longrightarrow X$$
is a \textbf{P}-morphism between representable stacks.\\
\end{enumerate}

Now let $n\geq 0$.

\begin{enumerate}

\item Let $F$ be any stack. An \emph{$n$-atlas for $F$}
is a $\mathbb{U}$-small family of morphisms
$\{U_{i} \longrightarrow F\}_{i\in I}$ such that
\begin{enumerate}
\item Each $U_{i}$ is representable.

\item Each morphism $U_{i} \longrightarrow F$
is in $(n-1)$-\textbf{P}.

\item The total morphism
$$\coprod_{i\in I}U_{i} \longrightarrow F$$
is an epimorphism.

\end{enumerate}

\item A stack $F$ is \emph{$n$-geometric}\index{stack!$n$-geometric} if
it satisfies the following two conditions.
\begin{enumerate}
\item The diagonal morphism $F \longrightarrow F\times^{h} F$
is $(n-1)$-representable.
\item The stack $F$ admits an $n$-atlas.
\end{enumerate}
\item A morphism of stacks $F \longrightarrow G$
is \emph{$n$-representable}\index{morphism of stacks!$n$-representable} if for any
representable stack $X$ and any morphism
$X \longrightarrow G$,  the
homotopy pullback
$F\times^{h}_{G}X$ is $n$-geometric.

\item A morphism of stacks $F \longrightarrow G$ is
\emph{in $n$-\textbf{P}}\index{morphism of stacks!in $n$-\textbf{P}}
(or \emph{has the property $n$-\textbf{P}},
or \emph{is a $n$-\textbf{P}-morphism})
if it is $n$-representable and
if for any representable stack $X$, any morphism $X \longrightarrow G$,
there exists an $n$-atlas $\{U_{i}\}$ of $F\times^{h}_{G}X$, such that
each composite morphism $U_{i} \longrightarrow X$ is in \textbf{P}.

\end{enumerate}
\end{df}

\begin{rmk}
\emph{In the above definition, condition (2a) follows from condition 
(2b). This is not immediate now but will be an easy consequence of the description of geometric
stacks as quotients by groupoids given in the next section. We prefer to keep 
the definition of n-geometric stacks with the two conditions
(2a) and (2b) as it is very similar to the usual definition 
of algebraic stacks found in the literature (e.g. in \cite{lm}).}
\end{rmk}

The next Proposition gives the fundamental properties of geometric $n$-stacks.

\begin{prop}\label{p11}
\noindent \begin{enumerate}
\item Any $(n-1)$-representable morphism is $n$-representable.
\item Any $(n-1)$-\textbf{P}-morphism is a $n$-\textbf{P}-morphism.
\item $n$-representable  morphisms
are stable by isomorphisms, homotopy pullbacks and
compositions.
\item $n$-\textbf{P}-morphisms
are stable by isomorphisms, homotopy pullbacks and
compositions.
\end{enumerate}
\end{prop}

\begin{proof}  We use a big induction on $n$. All the assertions
are easily verified for $n=-1$ 
using our assumptions \ref{ass4} on the
morphisms in \textbf{P}.
So, we fix an integer
$n\geq 0$ and suppose that \textit{all} the assertions are true for
any $m<n$; let's prove that they all remain true at the level $n$.\\

$(1)$ By definition \ref{d11} it is
enough to check that a $(n-1)$-geometric stack $F$ is
$n$-geometric. But an $(n-1)$-atlas for $F$ is a $n$-atlas by induction hypothesis
(which tells us in particular that a $(n-2)$-\textbf{P}-morphism is a $(n-1)$-\textbf{P}-morphism),
and moreover the diagonal of $F$ is $(n-2)$-representable thus $(n-1)$-representable,
again by induction hypothesis. Therefore $F$ is $n$-geometric. \\

$(2)$ Let $f : F \longrightarrow G$ be an $(n-1)$-\textbf{P}-morphism.
By definition it is $(n-1)$-representable, hence $n$-representable by $(1)$. Let
$X \longrightarrow G$ be a morphism with $X$ representable.
Since $f$ is in $(n-1)$-\textbf{P}, there exists an $(n-1)$-atlas
$\{U_{i}\}$ of $F\times^{h}_{G}X$ such that
each $U_{i} \longrightarrow X$ is in \textbf{P}. But, as already
observed in $(1)$, our inductive hypothesis, shows that
any $(n-1)$-atlas is an $n$-atlas, and we conclude that
$f$ is also in $n$-\textbf{P}. \\

$(3)$ Stability by isomorphisms and homotopy pullbacks
is clear by definition. To prove the stability
by composition, it is enough to prove that if
$f : F \longrightarrow G$ is an $n$-representable morphism
and $G$ is $n$-geometric then so is $F$.

Let $\{U_{i}\}$ be an $n$-atlas of $G$, and let $F_{i}:=F\times^{h}_{G}U_{i}$.
The stacks $F_{i}$ are $n$-geometric, so we can find
an $n$-atlas $\{V_{i,j}\}_{j}$ for $F_{i}$, for any $i$. By induction hypothesis
(telling us in particular that $(n-1)$-\textbf{P}-morphisms are closed under composition)
we see that the family of morphisms $\{V_{i,j} \longrightarrow F\}$
is an $n$-atlas for $F$. It remains to show that the diagonal
of $F$ is $(n-1)$-representable.

There is a homotopy cartesian square
$$\xymatrix{F\times_{G}^{h}F \ar[r] \ar[d]& F\times^{h} F \ar[d]\\
G \ar[r] & G\times^{h} G.}$$
As $G$ is $n$-geometric, the stability
of $(n-1)$-representable morphisms under homotopy pullbacks (true by induction hypothesis) implies  that
$F\times^{h}_{G}F \longrightarrow F\times^{h} F$
is $(n-1)$-representable. Now, the diagonal
of $F$ factors has
$F \longrightarrow F\times_{G}^{h}F \longrightarrow F\times^{h} F,$
and therefore by stability of $(n-1)$-representable morphisms by
composition (true by induction hypothesis), we see that
it is enough to show that $F \longrightarrow F\times^{h}_{G}F$
is $(n-1)$-representable. Let $X$ be a representable
stack and $X \longrightarrow F\times^{h}_{G}F$ be any
morphism. Then, we have
$$F\times^{h}_{F\times^{h}_{G}F}X\simeq
X\times^{h}_{(F\times^{h}_{G}X)\times^{h}(F\times^{h}_{G}X)}
(F\times^{h}_{G}X).$$
As by hypothesis $F\times^{h}_{G}X$ is $n$-geometric this
shows that $F\times^{h}_{F\times^{h}_{G}F}X$ is $(n-1)$-geometric,
showing that $F\longrightarrow F\times^{h}_{G}F$
is $(n-1)$-representable. \\

$(4)$ Stability by isomorphisms and homotopy pullbacks
is clear by definition. Let $F \longrightarrow G \longrightarrow H$
be two $n$-\textbf{P}-morphisms of stacks. By $(3)$ we already know
the composite morphism to be $n$-representable. By definition
of being in $n$-\textbf{P} one can assume that $H$ is representable.
Then, there exists an $n$-atlas $\{U_{i}\}$ for $G$ such that
each morphism $U_{i} \longrightarrow H$ is in \textbf{P}.
Let $F_{i}:=F\times^{h}_{G}U_{i}$, and let
$\{V_{i,j}\}_{j}$ be an $n$-atlas for
$F_{i}$ such that each $V_{i,j} \longrightarrow U_{i}$ is
in \textbf{P}. Since by induction hypothesis $(n-1)$-\textbf{P}-morphisms are
closed under composition, we
see that $\{V_{i,j}\}$ is indeed an $n$-atlas for $F$ such that
each $V_{i,j} \longrightarrow H$ is in \textbf{P}. 
\end{proof}

An important consequence of our descent assumption
\ref{ass5} and Prop. \ref{c0}
is the following useful proposition.

\begin{prop}\label{p14}
Let $f : F \longrightarrow G$ be any morphism such that
$G$ is an $n$-geometric stack. We suppose that
there exists a $n$-atlas $\{U_{i}\}$ of $G$
such that each stack $F\times^{h}_{G}U_{i}$
is $n$-geometric. Then $F$ is $n$-geometric.

If furthermore each projection $F\times_{X}^{h}U_{i} \longrightarrow U_{i}$
is in $n$-\textbf{P}, then so is $f$.
\end{prop}

\begin{proof} Using the stability of
$n$-representable morphisms by composition
(see Prop. \ref{p11} $(3)$) we see that it is
enough to show that $f$ is $n$-representable.
The proof goes by induction on $n$. For
$n=-1$ this is our corollary
Prop. \ref{c0}. Let us assume
$n\geq 0$ and the proposition proved for rank less than $n$.
Using Prop. \ref{p11} $(3)$ it is enough to suppose that
$G$ is a representable stack $X$.

Let $\{U_{i}\}$ be an $n$-atlas of $X$ as in the statement, and
let $\{V_{i,j}\}$ be an $n$-atlas for $F_{i}:=F\times^{h}_{X}U_{i}$.
Then, the composite family $V_{i,j} \longrightarrow F$
is clearly an $n$-atlas for $F$. It remains to prove that
$F$ has an $(n-1)$-representable diagonal.

The diagonal of $F$ factors as
$$F \longrightarrow F\times_{X}^{h}F \longrightarrow F\times^{h}F.$$
The last morphism being the homotopy pullback
$$\xymatrix{F\times_{X}^{h}F \ar[r] \ar[d] & F\times^{h}F \ar[d] \\
X \ar[r] & X\times^{h}X}$$
is representable and therefore $(n-1)$-representable. Finally,
let $Z$ be any representable stack and $Z \longrightarrow F\times_{X}^{h}F$
be a morphism. Then, the morphism
$$Z\times^{h}_{F\times_{X}^{h}F}F \longrightarrow X$$
satisfies the conditions of the proposition \ref{p14}
for the rank $(n-1)$. Indeed, for any $i$ we have
$$(Z\times^{h}_{F\times_{X}^{h}F}F)\times^{h}_{X}U_{i}\simeq
(Z\times_{X}^{h}U_{i})\times_{F_{i}\times_{U_{i}}^{h}F_{i}}^{h}F_{i}.$$
Therefore, using the induction hypothesis we deduce
that the stack $Z\times^{h}_{F\times_{X}^{h}F}F$ is $(n-1)$-geometric, proving
that $F \longrightarrow F\times_{X}^{h}F$ is $(n-1)$-representable.

The last part of the proposition follows from the fact that
any $n$-atlas $\{V_{i,j}\}$ of $F_{i}$ is such that
each morphism $V_{i,j} \longrightarrow X$ is in $n$-\textbf{P}
by construction. \end{proof}

\begin{cor}\label{c10}
The full subcategory of $n$-geometric stacks
in $\mathrm{St}(\mathcal{C},\tau)$
is stable by homotopy pullbacks, and by
$\mathbb{U}$-small disjoint union 
if $n\geq 0$.
\end{cor}

\begin{proof} Let $\xymatrix{F\ar[r] & H & G\ar[l]}$
be a diagram of stacks.
There are two homotopy  cartesian squares
$$\xymatrix{
F\times^{h}G \ar[r] \ar[d] & G \ar[d] & \qquad & F\times_{H}^{h}G \ar[r] \ar[d] &
H\ar[d] \\
F \ar[r] & \bullet& \qquad & F\times^{h}G \ar[r] & H\times^{h}H,}$$
showing that the stability under homotopy pullbacks follows from
the stability of $n$-representable morphisms under compositions
and homotopy pullbacks.

Let us prove the second part of the corollary, concerning
$\mathbb{U}$-small disjoint union.
Suppose now tha
$n\geq 0$ and
let $F$ be $\coprod_{i}F_{i}$ with each $F_{i}$
an $n$-geometric stack. Then, we have
$$F\times^{h}F\simeq \coprod_{i,j}F_{i}\times^{h}F_{j}.$$
For any representable
stack $X$, and any morphism $X \longrightarrow F\times^{h} F$,
there exists a $0$-atlas $\{U_{k}\}$ of $X$, and
commutative diagrams of stacks
$$\xymatrix{
U_{k} \ar[r] \ar[d] & F_{i(k)}\times^{h}F_{j(k)} \ar[d] \\
X \ar[r] & \coprod_{i,j}F_{i}\times^{h}F_{j}.
}$$
We apply Prop. \ref{p14} to the morphism
$$G:=F\times^{h}_{F\times^{h}F}X \longrightarrow X$$
and for the covering $\{U_{k}\}$. We have
$$G\times_{X}^{h}U_{k}\simeq \emptyset \; if \; i(k)\neq j(k)$$
and
$$G\times_{X}^{h}U_{k}\simeq F_{i(k)}\times^{h}_{F_{i(k)}\times^{h}F_{i(k)}}U_{k}$$
otherwise. Prop. \ref{p14} implies that
$G$ is $(n-1)$-geometric, and therefore that the diagonal
of $F$ is $(n-1)$-representable.
Finally, the same argument and assumption \ref{ass4} show that
the disjoint union of $n$-atlases of the $F_{i}$'s will
form an $n$-atlas for $F$.  \end{proof}

Finally, let us mention  the following important additional property.

\begin{prop}\label{p11'}
Let $f : F \longrightarrow G$ be
an $n$-representable morphism. If $f$
is in $m$-\textbf{P} for $m>n$ then it
is in $n$-\textbf{P}.
\end{prop}

\begin{proof} By induction on $m$ it is enough to treat the
case $m=n+1$. The proof goes then by induction on $n$. For
$n=-1$ 
this is our assumption \ref{ass4} $(3)$. For 
$n\geq 0$
we can by definition assume $G$ is a representable
stack and therefore that $F$ is $n$-geometric.
Then, there exists an $(n+1)$-atlas $\{U_{i}\}$
for $F$ such that each $U_{i} \longrightarrow G$ is in
\textbf{P}. By induction, $\{U_{i}\}$ is also
an $n$-atlas for $F$, which implies that
the morphism $f$ is in fact in \textbf{P}. 
\end{proof}

The last proposition implies in particular
that for an $n$-representable
morphism of stacks the property of being in
$n$-\textbf{P} does not depend on $n$. We will therefore
give the following definition.

\begin{df}\label{d11'}
A morphism in $\mathrm{St}(\mathcal{C},\tau)$ is \emph{in \textbf{P}}\index{morphism of stacks!in \textbf{P}} if it
is in $n$-\textbf{P} for some integer $n$.
\end{df}

\section{Quotient stacks}

We will now present a characterization of
geometric stacks in terms of quotient by groupoid
actions. This point of view is very much similar
to the presentation of manifolds by charts, and
much less intrinsic than definition Def. \ref{d11}.
However, it is sometimes more easy to handle
as several stacks have natural presentations 
as quotients by groupoids.\\

Let $X_{*}$ be a Segal groupoid object in a model category $M$ (Def. \ref{d12}).
Inverting the equivalence
$$X_{2} \longrightarrow X_{1}\times^{h}_{X_{0}}X_{1}$$
and composing with $d_{1} : X_{2} \longrightarrow X_{1}$
gives a well defined morphism in $\mathrm{Ho}(M)$
$$\mu : X_{1}\times^{h}_{X_{0}}X_{1} \longrightarrow X_{1},$$
that is called \emph{composition}. In the same way, inverting
the equivalence
$$X_{2} \longrightarrow X_{1}\times^{h}_{d_{0},X_{0},d_{0}}X_{1},$$
and composing with
$$Id\times^{h}s_{0} : X_{1} \longrightarrow X_{1}\times^{h}_{d_{0},X_{0},d_{0}}X_{1}$$
$$d_{2} : X_{2} \longrightarrow X_{1}$$
gives a well defined isomorphism in $\mathrm{Ho}(M)$
$$i : X_{1} \longrightarrow X_{1},$$
called \emph{inverse}. It is easy to check that
$d_{1}\circ i = d_{0}$ as morphisms in $\mathrm{Ho}(M)$, showing that
the two morphisms $d_{0}$ and $d_{1}$ are
always isomorphic in $\mathrm{Ho}(M)$. Finally, using condition
$(1)$ of Def. \ref{d12} we see that for any $i>0$, all the
face morphisms $$X_{i} \longrightarrow X_{i-1}$$
of a Segal groupoid object
are isomorphic in $\mathrm{Ho}(M)$.

\begin{df}\label{d13}
A Segal groupoid object $X_{*}$ in
$\mathrm{Ho}(SPr(Aff_{\mathcal{C}}^{\sim,\tau}))$ is
\emph{an $n$-\textbf{P} Segal groupoid}\index{Segal groupoid!$n$-\textbf{P} } if
it satisfies the following two conditions.
\begin{enumerate}
\item The stacks $X_{0}$ and
$X_{1}$ are 
disjoint unions of
$n$-geometric stacks.

\item The morphism $d_{0} : X_{1} \longrightarrow X_{0}$
is in
$n$-\textbf{P}.

\end{enumerate}
\end{df}

As $n$-geometric stacks are stable by homotopy pullbacks,
$X_{i}$ is a
disjoint union of $n$-geometric stacks for any
$i$ and any $n$-\textbf{P} Segal groupoid $X_{*}$.
Furthermore, the condition $(2)$ of
Def. \ref{d12} implies that the two morphisms
$$d_{0}, d_{1} : X_{1} \longrightarrow X_{0}$$
are isomorphic as morphisms in $\mathrm{St}(\mathcal{C},\tau)$.
Therefore, for any $n$-\textbf{P} Segal groupoid
$X_{*}$, all the faces $X_{i} \longrightarrow X_{i-1}$
are in \textbf{P}.

\begin{prop}\label{p13}
Let $F\in \mathrm{Ho}(SPr(Aff_{\mathcal{C}}^{\sim,\tau}))$
be a stack and 
$n\geq 0$. The following conditions are
equivalent.
\begin{enumerate}
\item The stack $F$ is $n$-geometric.
\item There exists an $(n-1)$-\textbf{P} Segal
groupoid object $X_{*}$ in $SPr(Aff_{\mathcal{C}}^{\sim,\tau})$,
such that $X_{0}$ is 
a disjoint union of representable stacks, and an isomorphism in
$\mathrm{Ho}(SPr(Aff_{\mathcal{C}}^{\sim,\tau}))$
$$F\simeq |X_{*}|:=Hocolim_{[n]\in \Delta}X_{n}.$$
\item
There exists an $(n-1)$-\textbf{P} Segal
groupoid object $X_{*}$ in $SPr(Aff_{\mathcal{C}}^{\sim,\tau})$,
and an isomorphism in
$\mathrm{Ho}(SPr(Aff_{\mathcal{C}}^{\sim,\tau}))$
$$F\simeq |X_{*}|:=Hocolim_{[n]\in \Delta}X_{n}.$$
\end{enumerate}
If these conditions are satisfied we say that
\emph{$F$ is the quotient stack of the $(n-1)$-\textbf{P} Segal groupoid
$X_{*}$}.
\end{prop}

\begin{proof} We have already seen that
a $0$-geometric stack is $n$-geometric for any $n$. Therefore,
$(2) \Rightarrow (3)$. It remains to show that
$(1)$ implies $(2)$ and $(3)$ implies $(1)$. \\

$(1) \Rightarrow (2)$ Let $F$ be an $n$-geometric stack, and
$\{U_{i}\}$ be an $n$-atlas for $F$. We let
$$p : X_{0}:= \coprod_{i}U_{i} \longrightarrow F$$
be the natural projection. Up to an equivalence, we can  represent
$p$ by a fibration
$X_{0} \longrightarrow F$ between fibrant objects
in $SPr(Aff_{\mathcal{C}}^{\sim,\tau})$. We define
a simplicial object $X_{*}$ to be the nerve of $p$
$$X_{n}:=\underbrace{X_{0}\times_{F}X_{0}\times_{F}\dots
\times_{F}X_{0}}_{n\; times}.$$
Clearly, $X_{*}$ is a groupoid object in $SPr(Aff_{\mathcal{C}}^{\sim,\tau})$
in the usual sense, and as $p$ is a fibration between fibrant objects
it follows
that it is also a Segal groupoid object in the sense
of Def. \ref{d12}. Finally, as $F$ is $n$-geometric,
one has $X_{1}\simeq \coprod_{i,j}U_{i}\times_{F}^{h}U_{j}$
which is therefore an
$(n-1)$-geometric stack by Cor. \ref{c10}. The morphism
$d_{0} : X_{1} \longrightarrow X_{0}\simeq \coprod_{i}U_{i}$
is then given by the projections $U_{i}\times_{F}^{h}U_{j} \longrightarrow U_{i}$
which are in $(n-1)$-\textbf{P} as $\{U_{i}\}$ is an $n$-atlas.
This implies that $X_{*}$ is an $(n-1)$-\textbf{P} Segal groupoid
such that $X_{0}$ is 
a disjoint union of representable.

\begin{lem}\label{lquot}
The natural morphism
$$|X_{*}| \longrightarrow F$$
is an isomorphism in $\mathrm{Ho}(SPr(Aff_{\mathcal{C}}^{\sim,\tau}))$.
\end{lem}

\begin{proof} For any fibration of simplicial sets
$f : X \longrightarrow Y$, we know that
the natural morphism from the geometric realization
of the nerve of $f$ to $Y$ is equivalent to an inclusion
of connected components. This implies
that the morphism $|X_{*}| \longrightarrow F$ induces
an isomorphisms on the sheaves $\pi_{i}$ for $i>0$
and an injection on $\pi_{0}$. As the morphism
$p : X_{0} \longrightarrow F$ is an epimorphism and factors
as $X_{0} \rightarrow |X_{*}| \rightarrow F$,
the morphism $|X_{*}| \longrightarrow F$
is also an epimorphism. This shows that it is surjective
on the sheaves $\pi_{0}$ and therefore is an isomorphism
in $\mathrm{Ho}(SPr(Aff_{\mathcal{C}}^{\sim,\tau})$. \end{proof}

The previous lemma finishes the proof of $(1) \Rightarrow (2)$. \\

$(3) \Rightarrow (1)$ Let $X_{*}$
be an $(n-1)$-\textbf{P} Segal groupoid and
$F=|X_{*}|$. First of all,
we recall the following important fact.

\begin{lem}\label{l13}
Let $M$ be a $\mathbb{U}$-model topos, and
$X_{*}$ be a Segal groupoid object in $M$
with homotopy colimit $|X_{*}|$. Then, for any $n>0$, the
natural morphism
$$X_{n} \longrightarrow \underbrace{X_{0}\times^{h}_{|X_{*}|}X_{0}\times^{h}_{|X_{*}|}
\dots \times^{h}_{|X_{*}|}X_{0}}_{n\; times}$$
is an isomorphism in $\mathrm{Ho}(M)$.
\end{lem}

\begin{proof} This is one of the
standard properties of model topoi. See Thm. \ref{tgiraud}. \end{proof}

Let $\{U_{i}\}$ be an $(n-1)$-atlas
for $X_{0}$, and
let us consider the composed morphisms
$$f_{i} : U_{i} \longrightarrow X_{0} \longrightarrow F.$$
Clearly, $\coprod_{i}U_{i} \longrightarrow F$
is a composition of epimorphism, and is
therefore a epimorphism. In order to prove that $\{U_{i}\}$
form an $n$-atlas for $F$ it is enough to prove that
the morphism $X_{0} \longrightarrow F$ is in $(n-1)$-\textbf{P}.

Let $X$ be any representable stack, $X \longrightarrow F$
be a morphism, and let
$G$ be $X_{0}\times^{h}_{F}X$. As the morphism $X_{0} \longrightarrow F$ is
an epimorphism, we can find a covering family
$\{Z_{j} \longrightarrow X\}$, such that each
$Z_{j}$ is representable, and such that there
exists a commutative diagram in $\mathrm{St}(\mathcal{C},\tau)$
$$\xymatrix{
X_{0} \ar[r] & F \\
 & X \ar[u] \\
 & \coprod_{j} Z_{j}. \ar[u] \ar[uul]}$$
By assumption \ref{ass4} $(1)$ we can also assume that 
each morphism $Z_{j} \longrightarrow X$ is in \textbf{P}, and therefore
that $\{Z_{j} \longrightarrow X\}$ is a $0$-atlas of $X$.

In order to prove that $G \longrightarrow X$ is in $(n-1)$-\textbf{P}
it is enough by Prop. \ref{p14} to prove that each
stack $G_{j}:=G\times^{h}_{X}Z_{j}\simeq X_{0}\times_{F}^{h}Z_{j}$
is $(n-1)$-geometric, and furthermore that each
projection $G_{j} \longrightarrow Z_{j}$ is in
\textbf{P}. We have
$$G_{j}\simeq (X_{0}\times_{F}^{h}X_{0})\times^{h}_{X_{0}}Z_{j}.$$
Therefore, lemma \ref{l13} implies that
$G_{j}\simeq X_{1}\times^{h}_{X_{0}}Z_{j}$, showing that
each $G_{j}$ is $(n-1)$-geometric and finally that
$G$ is $(n-1)$-geometric. Furthermore, this also shows that
each projection $G_{j} \longrightarrow Z_{j}$
is of the form $X_{1}\times^{h}_{X_{0}}Z_{j} \longrightarrow Z_{j}$
which is in \textbf{P} by hypothesis on the
Segal groupoid object $X_{*}$. \end{proof}

\begin{cor}\label{cp13}
Let $f : F \longrightarrow G$ be an epimorphism of stacks and $n\geq 0$.
If $F$ is $n$-geometric and $f$ is $(n-1)$-representable and
in \textbf{P}, then $G$ is $n$-geometric.
\end{cor}

\begin{proof} Indeed, let $U \longrightarrow F$ be an
$n$-atlas, and let $g : U \longrightarrow G$ be the
composition. The morphism $g$ is still an epimorphism and
$(n-1)$-representable, and thus we can assume that
$F$ is $0$-representable.
The morphism $f$ being an epimorphism, $G$ is equivalent to the
quotient stack of the Segal groupoid $X_{*}$ which is the homotopy nerve
of $f$. By assumption this Segal groupoid is an $(n-1)$-Segal groupoid
and thus $G$ is 
$n$-geometric by
Prop. \ref{p13} $(1)$. \end{proof}

\section{Quotient stacks and torsors}

Writing a stack as a quotient stack of a Segal groupoid
is also useful in order to describe associated stacks
to certain pre-stacks. Indeed, it often happens that
a pre-stack is defined as the quotient of a Segal
groupoid, and we are going to show in this section
that the associated stack of such quotients can be described
using an adequate notion of torsor over a Segal groupoid.
This section is in fact completely independent of the
notion of geometricity and concern pure stacky statements
that are valid in any general model topos. We have decided
to include this section as it helps understanding the
quotient stack construction presented previously. However, it
is not needed for the rest of this work and proofs will
be more sketchy than usual. \\

We let $M$ be a general $\mathbb{U}$-model topos
in the sense of Def. \ref{dmodtop}. The main case of application
will be $M=Aff_{\mathcal{C}}^{\sim,\tau}$ but we rather prefer to state
the results in the most general setting (in particular we do
not even assume that $M$ is $t$-complete).

Let $X_{*}$ be a Segal groupoid object in $M$ in the
sense of Def. \ref{d12}, and we assume that
each $X_{n}$ is a fibrant object in $M$.
We will consider
$sM$, the category of simplicial objects in $M$, which
will be endowed with its levelwise projective
model structure, for which fibrations and equivalences
are defined levelwise. We consider
$sM/X_{*}$, the model category of simplicial objects over
$X_{*}$. Finally, we let $Z$ be a fibrant object in $M$, and
$X_{*} \longrightarrow Z$ be a morphism in $sM$ (where
$Z$ is considered as a constant simplicial object), and we assume that
the induced morphism
$$|X_{*}| \longrightarrow Z$$
is an isomorphism in $\mathrm{Ho}(M)$.

We define a Quillen adjunction
$$\phi : sM/X_{*} \longrightarrow M/Z \qquad
sM/X_{*} \longleftarrow M/Z : \psi$$
in the following way. For any $Y \longrightarrow Z$ in $M$, we set
$$\psi(Y):=Y\times_{Z}X_{*}\in sM,$$
or in other words $\psi(Y)_{n}=Y\times_{Z}X_{n}$ and with the obvious
transitions morphisms. The left adjoint to $\psi$ sends
a simplicial object $Y_{*} \longrightarrow X_{*}$ to
its colimit in $M$
$$\phi(Y_{*}):=Colim_{n\in \Delta^{op}}Y_{n} \longrightarrow
Colim_{n}X_{n\in \Delta^{op}} \longrightarrow Z.$$
It is easy to check that $(\phi,\psi)$ is a Quillen
adjunction.

\begin{prop}\label{pdesc}
The functor
$$\mathbb{R}\psi : \mathrm{Ho}(M/Z) \longrightarrow
\mathrm{Ho}(sM/X_{*})$$
is fully faithful. Its essential image consists of
all $Y_{*} \longrightarrow X_{*}$ such that for
any morphism $[n] \rightarrow [m]$ in $\Delta$ the
square
$$\xymatrix{
Y_{m} \ar[r] \ar[d] & Y_{n} \ar[d] \\
X_{m} \ar[r] & X_{n}}$$
is homotopy cartesian.
\end{prop}

\begin{proof} Let $Y \rightarrow Z$ in $M/Z$. Proving that
$\mathbb{R}\psi$ is fully faithful is equivalent to prove that the
natural morphism
$$|Y\times^{h}_{Z}X_{*}| \longrightarrow Y$$
is an isomorphism in $\mathrm{Ho}(M)$. But, using the standard
properties of model topoi (see Thm. \ref{tgiraud}), we have
$$|Y\times^{h}_{Z}X_{*}|\simeq
Y\times^{h}_{Z}|X_{*}|\simeq Y$$
as $|X_{*}| \simeq Z$. This shows that
$\mathbb{R}\psi$ is fully faithful.

By definition of the functors $\phi$ and $\psi$, it is clear that
$\mathbb{R}\psi$ takes its values in the
subcategory described in the statement of the proposition.
Conversely, let
$Y_{*} \longrightarrow X_{*}$ be an object in $\mathrm{Ho}(M)$
satisfying the condition of the proposition. As $X_{*}$ is
a Segal groupoid object, we know that $X_{*}$ is naturally
equivalent to the homotopy nerve
of the augmentation morphism $X_{0} \longrightarrow Z$
(see Thm. \ref{tgiraud}).
Therefore, the object
$\mathbb{R}\psi\mathbb{L}\phi(Y_{*})$ is by definition the
homotopy nerve of the morphism
$$|Y_{*}|\times^{h}_{Z}X_{0} \longrightarrow |Y_{*}|.$$
But, we have
$$|Y_{*}|\times^{h}_{Z}X_{0}\simeq
|Y_{*}\times^{h}_{Z}X_{0}|\simeq Y_{0}$$
by hypothesis on $Y_{*}$. Therefore, the
object $\mathbb{R}\psi\mathbb{L}\phi(Y_{*})$ is naturally
isomorphic in $\mathrm{Ho}(sM/X_{*})$ to the homotopy nerve
of the natural
$$Y_{0} \longrightarrow |Y_{*}|.$$
Finally, as $X_{*}$ is a Segal groupoid object, so is
$Y_{*}$ by assumption. The standard properties
of model topoi (see Thm. \ref{tgiraud})
then tell us that
$Y_{*}$ is naturally equivalent to the
homotopy nerve of $Y_{0} \longrightarrow |Y_{*}|$, and thus
to $\mathbb{R}\psi\mathbb{L}\phi(Y_{*})$ by what we have
just done. \end{proof}

The model category $sM/X_{*}$, or rather its full subcategory
of objects satisfying the conditions of Prop. \ref{pdesc},
can be seen as the category of objects in $M$ together
with an \emph{action of the Segal groupoid $X_{*}$}. Proposition
\ref{pdesc} therefore says that the homotopy theory of
stacks over $|X_{*}|$ is equivalent to the homotopy theory
of stacks together with an action of $X_{*}$. This point of view will
now help us to describe the stack associated to $|X_{*}|$.

For this, let $F$ be a fixed fibrant object in
$M$. We define a new model category
$sM/(X_{*},F)$ in the following way. Its objects are
pairs $(Y_{*},f)$, where $Y_{*} \rightarrow X_{*}$ is
an object in $sM/X_{*}$ and $f : Colim_{n}Y_{n} \longrightarrow F$
is a morphism in $M$. Morphisms $(Y_{*},f) \rightarrow (Y_{*}',f')$
are given by morphisms $Y_{*} \longrightarrow Y_{*}'$ in
$sM/X_{*}$, such that
$$\xymatrix{
Colim_{n}Y_{n} \ar[rr] \ar[rd]_-{f} & & Colim_{n}Y'_{n} \ar[ld]^-{f'}\\
 & F & }$$
commutes in $M$. The model structure on $sM/(X_{*},F)$ is defined
in such a way that fibrations and equivalences are defined
on the underlying objects in $sM$. The model category
$sM/(X_{*},F)$ is also the comma model category
$sM/(X_{*}\times F)$ where $F$ is considered as
a constant simplicial object in $M$.

\begin{df}\label{dtors}
An object $Y_{*} \in sM/(X_{*},F)$ is a \emph{$X_{*}$-torsor\index{torsor}
on $F$} if it satisfies the following two conditions.
\begin{enumerate}
\item
For all morphism $[n] \rightarrow [m]$ in $\Delta$, the
square
$$\xymatrix{
Y_{m} \ar[r] \ar[d] & Y_{n} \ar[d] \\
X_{m} \ar[r] & X_{n}}$$
is homotopy cartesian.
\item The natural morphism
$$|Y_{*}| \longrightarrow F$$
is an isomorphism in $\mathrm{Ho}(M)$.
\end{enumerate}
The \emph{space of $X_{*}$-torsors over $F$}, denoted by
$Tors_{X_{*}}(F)$, is the nerve of
the sub category of fibrant objects $sM/(X_{*},F)^{f}$, consisting of
equivalences between $X_{*}$-torsors
on $F$.
\end{df}

Suppose that $f : F \longrightarrow F'$ is a morphism
between fibrant objects in $M$. We get a pullback
functor
$$sM/(X_{*},F') \longrightarrow sM/(X_{*},F),$$
which is right Quillen, and such that
the induced functor on fibrant objects
$$sM/(X_{*},F')^{f} \longrightarrow sM/(X_{*},F)^{f}$$
sends $X_{*}$-torsors over $F$ to $X_{*}$-torsors
over $F'$ (this uses the commutation of
homotopy colimits with homotopy pullbacks). Therefore,
restricting to the sub categories of equivalences, we get
a well defined morphism between spaces of torsors
$$f^{*} : Tors_{X_{*}}(F') \longrightarrow Tors_{X_{*}}(F).$$
By applying the standard strictification procedure, 
we can always suppose that $(f\circ g)^{*}=g^{*}\circ f^{*}$.
This clearly defines a  functor from $(M^{f})^{op}$, the
opposite full subcategory of fibrant objects in $M^{f}$, to
$SSet$
$$Tors_{X_{*}} : (M^{f})^{op} \longrightarrow SSet.$$
This functor sends equivalences in $M^{f}$ to equivalences
of simplicial sets, and therefore induces a $\mathrm{Ho}(SSet)$-enriched
functor (using for example \cite{dk1})
$$Tors_{X_{*}} : \mathrm{Ho}(M^{f})^{op}\simeq \mathrm{Ho}(M)^{op} \longrightarrow \mathrm{Ho}(SSet).$$
In other words, there are natural morphisms in $\mathrm{Ho}(SSet)$
$$Map_{M}(F,F') \longrightarrow
Map_{SSet}(Tors_{X_{*}}(F'),Tors_{X_{*}}(F)),$$
compatible with compositions.

The main classification result is the following. It gives a way to
describe the stack associated to $|X_{*}|$ for some
Segal groupoid object $X_{*}$ in $M$.

\begin{prop}\label{pdesc2}
Let $X_{*}$ be a Segal groupoid object in $M$ and
$Z$ be a fibrant model for $|X_{*}|$ in $M$.
Then, there exists an element
$\alpha \in \pi_{0}(Tors_{X_{*}}(Z))$, such that for
any fibrant object $F\in M$, the evaluation at $\alpha$
$$\alpha^{*} : Map_{M}(F,Z) \longrightarrow Tors_{X_{*}}(F)$$
is an isomorphism in $\mathrm{Ho}(SSet)$.
\end{prop}

\begin{proof} We can of course assume that $X_{*}$ is
cofibrant in $sM$, as everything is invariant
by changing $X_{*}$ with something equivalent.
Let $Colim_{n}X_{n} \longrightarrow Z$ be a morphism
in $M$ such that $Z$ is fibrant and
such that the induced morphism
$$|X_{*}| \longrightarrow Z$$
is an equivalence. Such a morphism exists as
$Colim_{n}X_{n}$ is cofibrant in $M$ and computes
the homotopy colimit $|X_{*}|$. The element
$\alpha \in \pi_{0}(Tors_{X_{*}}(Z))$ is defined
to be the pair $(X_{*},p) \in \mathrm{Ho}(sM/(X_{*},Z))$, consisting of
the identity of $X_{*}$ and the natural augmentation
$p : Colim_{n}X_{n} \longrightarrow Z$. Clearly,
$\alpha$ is a $X_{*}$-torsor over $Z$.

Applying Cor. \ref{capp2} of Appendix A, we see that there
exists a homotopy fiber sequence
$$\xymatrix{
Map_{M}(F,Z) \ar[r] & N((M/Z)^{f}_{W}) \ar[r] & N(M^{f}_{W})}$$
where the $(M/Z)^{f}_{W}$ (resp. $M^{f}_{W}$)
is the subcategory of equivalences between
fibrant objects in $M/Z$ (resp. in $M$), and the morphism on the right
is the forgetful functor and the fiber is taken
at $F$. In the same way, there exists
a homotopy fiber sequence
$$\xymatrix{
Tors_{X_{*}}(F) \ar[r] & N(((sM/X_{*})^{f}_{W})^{cart}) \ar[r] &
N(M^{f}_{W})}$$
where $(sM/X_{*})^{f}_{W})^{cart}$ is the subcategory of
$sM/X_{*}$ consisting of equivalence between
fibrant objects satisfying condition $(1)$ of Def. \ref{dtors}, and where
the morphism on the right is induced
my the homotopy colimit functor of underlying simplicial objects
and the fiber is taken at $F$. There exists a morphism of
homotopy fiber sequences of simplicial sets
$$\xymatrix{
Map_{M}(F,Z) \ar[r] \ar[d]_-{\alpha^{*}} & N((M/Z)^{f}_{W}) \ar[d] \ar[r] & N(M^{f}_{W})
\ar[d]^-{Id} \\
Tors_{X_{*}}(F) \ar[r] & N(((sM/X_{*})^{f}_{W})^{cart}) \ar[r] &
N(M^{f}_{W})}$$
and the arrow in the middle is an equivalence
because of our proposition \ref{pdesc}.
\end{proof}

\section{Properties of morphisms}

We fix another class \textbf{Q} of morphisms
in $Aff_{\mathcal{C}}$, which is stable by
equivalences, compositions and homotopy pullbacks.
Using the Yoneda embedding the notion of morphisms
in \textbf{Q} (or simply \emph{\textbf{Q}-morphisms})
is extended to a notion of morphisms between
representable stacks.

\begin{df}\label{d20}
We say that \emph{morphisms in \textbf{Q} are
compatible with $\tau$ and \textbf{P}} (or equivalently
that  \emph{morphisms in \textbf{Q} are
local with respect to  $\tau$ and \textbf{P}}) if the following 
two conditions are satisfied:
\begin{enumerate}
	\item If $f : X \longrightarrow Y$ is a morphism
in $Aff_{\mathcal{C}}$ 
such that there exists
a covering family
$$\{U_{i}\longrightarrow X\}$$
with
each composite morphism $U_{i} \longrightarrow Y$ in \textbf{Q}, 
then $f$ belongs to \textbf{Q}.

	\item If $f : X \longrightarrow Y$ is a morphism
in $Aff_{\mathcal{C}}$ and there
exists a covering family
$$\{U_{i}\longrightarrow Y\}$$
such that each homotopy pullback morphism
$$X \times^{h}_{Y}U_{i} \longrightarrow U_{i}$$
is in \textbf{Q}, then $f$ belongs to \textbf{Q}.
  \end{enumerate}
\end{df}

For a class of morphism \textbf{Q}, compatible with
$\tau$ and \textbf{P} in the sense above we can
make the following definition.

\begin{df}\label{d21}
Let \textbf{Q} be a class of morphisms in $Aff_{\mathcal{C}}$, stable
by equivalences, homotopy pullbacks and compositions, and
which is compatible with $\tau$ and \textbf{P}
in the sense above.
A morphism of stacks $f : F \longrightarrow G$
is \emph{in \textbf{Q}} (or equivalently  \emph{is a
\textbf{Q}-morphism})\index{morphism of stacks!\textbf{Q}-morphism} if it is
$n$-representable for some $n$, and if
for any representable stack $X$ and any morphism
$X \longrightarrow G$ there exists an $n$-atlas
$\{U_{i}\}$ of $F\times^{h}_{G}X$ such that
each morphism $U_{i} \longrightarrow X$ between representable stacks
is in \textbf{Q}.
\end{df}

Clearly, because of our definition \ref{d20}, the notion of
morphism in \textbf{Q} of definition \ref{d21} is compatible
with the original notion. Furthermore, it is easy to check, as
it was done for \textbf{P}-morphisms, the following proposition.

\begin{prop}\label{p20}
\begin{enumerate}
\item
Morphisms in \textbf{Q} are stable by equivalences, compositions and
homotopy pullbacks.

\item
Let $f : F \longrightarrow G$ be any morphism between
$n$-geometric stacks. We suppose that
there exists a $n$-atlas $\{U_{i}\}$ of $G$
such that each  projection $F\times_{X}^{h}U_{i} \longrightarrow U_{i}$
is in \textbf{Q}. Then $f$ is in \textbf{Q}.
\end{enumerate}
\end{prop}

\begin{proof} Exercise. \end{proof}

We can also make the following two general definitions of morphisms
of stacks.

\begin{df}\label{d23}
Let $f : F \longrightarrow G$ be a morphism of stacks.
\begin{enumerate}
\item The morphism is \emph{categorically locally finitely presented}\index{morphism of stacks!categorically locally finitely presented} if
for any representable
stack $X=\mathbb{R}\underline{Spec}\, A$, any morphism $X \longrightarrow G$,
and any $\mathbb{U}$-small filtered system of commutative $A$-algebras
$\{B_{i}\}$, the natural morphism
$$\mbox{Hocolim}_{i}
\mbox{Map}_{Aff_{\mathcal{C}}^{\sim,\tau}/X}(\mathbb{R}\underline{Spec}\, B_{i},F\times_{G}^{h}X)
\longrightarrow
\mbox{Map}_{Aff_{\mathcal{C}}^{\sim,\tau}/X}(\mathbb{R}\underline{Spec}\, (\mbox{Hocolim}_{i}B_{i}),F\times_{G}^{h}X)$$
is an isomorphism in $\mathrm{Ho}(SSet)$.

\item The morphism $f$ is \emph{quasi-compact}\index{morphism of stacks!quasi-compact} if for any representable
stack $X$ and any morphism $X \longrightarrow G$ there exists
a finite family of representable stacks $\{X_{i}\}$ and an epimorphism
$$\coprod_{i}X_{i} \longrightarrow F\times_{G}^{h}X.$$

\item The morphism $f$ is \emph{categorically finitely presented}\index{morphism of stacks!categorically finitely presented} if it is categorically  locally finitely presented
and quasi-compact.

\item The morphism $f$ is a \emph{monomorphism}\index{monomorphism of stacks} if the natural morphism
$$F \longrightarrow F\times_{G}^{h}F$$
is an isomorphism in $\mathrm{St}(\mathcal{C},\tau)$.

\item Assume that the class \textbf{Q} of finitely presented morphism
in $Aff_{\mathcal{C}}$ is compatible with the model topology $\tau$ in the
sense of Def. \ref{d20}. The morphism $f$ is \emph{locally finitely presented}
if it is a \textbf{Q}-morphism in the sense of Def. \ref{d21}. It is
a \emph{finitely presented morphism} if it is
quasi-compact and locally finitely presented.

\end{enumerate}
\end{df}

Clearly, the above notions are compatible with the one of definition
\ref{d5}, in the sense that a morphism of
commutative monoids $A \longrightarrow B$ has
a certain property in the sense of Def. \ref{d5} if and only if
the corresponding morphism of stacks
$\mathbb{R}\underline{Spec}\, B \longrightarrow \mathbb{R}\underline{Spec}\, A$ has the
same property in the sense of Def. \ref{d23}.

\begin{prop}
Quasi-compact morphisms, categorically (locally) finitely presented morphisms,
(locally) finitely presented morphisms
and monomorphisms are stable by equivalences,
composition and homotopy pullbacks.
\end{prop}

\begin{proof} Exercise. \end{proof}

\begin{rmk}\label{rlp}
\emph{When the class of finitely presented morphisms in
$Aff_{\mathcal{C}}$ is compatible with the topology $\tau$,
Def. \ref{d23} gives us two different notions of locally finitely presented
morphism which are a priori
rather difficult to compare. Giving precise conditions under which they
coincide is however not so important as much probably 
they are already different in some of our main examples (e.g. for each example
for which the representable stacks are not truncated, as in 
the complicial, or brave new algebraic geometry presented in \S 2.3 and \S 2.4).
In this work we have chosen to
use only the non categorical version of locally finitely presented morphisms, the price
to pay being of course that they do not have easy functorial characterization.}
\end{rmk}

\section{Quasi-coherent modules, perfect modules
and vector bundles}\label{Iqcoh}

For a commutative monoid $A$ in $\mathcal{C}$, we define
a category $A-QCoh$, of quasi-coherent modules
on $A$ (or equivalently on $Spec\, A$) in the following
way. Its objects are the data of
a $B$-module $M_{B}$ for any commutative $A$-algebra
$B\in A-Comm(\mathcal{C})$, together with
an isomorphism
$$\alpha_{u} : M_{B}\otimes_{B}B' \longrightarrow M_{B'}$$
for any morphism $u : B \longrightarrow B'$ in $A-Comm(\mathcal{C})$,
such that one has
$\alpha_{v}\circ (\alpha_{u}\otimes_{B'}B'')= \alpha_{v\circ u}$
for any pair of morphisms
$$\xymatrix{B \ar[r]^-{u} & B' \ar[r]^-{v} & B''}$$
in $A-Comm(\mathcal{C})$. Such data will be denoted
by $(M,\alpha)$.
A morphism in $A-QCoh$, from
$(M,\alpha)$ to $(M',\alpha')$ is given by
a family of morphisms of $B$-modules $f_{B} : M_{B} \longrightarrow M_{B}'$, for
any $B \in A-Comm(\mathcal{C})$, such that for any
$u : B \rightarrow B'$ in $A-Comm(\mathcal{C})$
the diagram
$$\xymatrix{
M_{B}\otimes_{B}B' \ar[r]^-{\alpha_{u}} \ar[d]_-{f_{B}\otimes_{B}B'} &
M_{B'} \ar[d]^-{f_{B'}} \\
(M_{B}')\otimes_{B}B' \ar[r]_-{\alpha_{u}'} & M_{B'}'}$$
commutes. As the categories
$A-Mod$ and $Comm(\mathcal{C})$ are all $\mathbb{V}$-small,
so are the categories $A-QCoh$.

There exists a natural projection $A-QCoh \longrightarrow
A-Mod$, sending $(M,\alpha)$ to $M_{A}$, and it is
straightforward to check that it is an equivalence of categories.
In particular, the model structure on $A-Mod$ will
be transported naturally on $A-QCoh$ through this equivalence.
Fibrations (resp. equivalences) in $A-QCoh$ are simply the morphisms
$f : (M,\alpha) \longrightarrow (M',\alpha')$ such that
$f_{A} : M_{A} \longrightarrow M_{A}'$ is a fibration
(resp. an equivalence).

Let now $f : A \longrightarrow A'$ be a morphism of
commutative monoids in $\mathcal{C}$. There exists a
pullback functor
$$f^{*} : A-QCoh \longrightarrow A'-QCoh$$
defined by $f^{*}(M,\alpha)_{B}:=M_{B}$ for any $B\in A-Comm(\mathcal{C})$, and
for $u : B \longrightarrow B'$ in $A-Comm(\mathcal{C})$ the transition
morphism
$$f^{*}(M,\alpha)_{B}\otimes_{B}B'=M_{B}\otimes_{B}B' \longrightarrow f^{*}(M,\alpha)_{B'}=M_{B'}$$
is given by $\alpha_{u}$. By definition of the model structure on
$A-QCoh$, the functor
$$f^{*} : A-QCoh \longrightarrow A'-QCoh$$
is clearly a left Quillen functor, as the natural
diagram
$$\xymatrix{
A-QCoh \ar[r]^-{f^{*}} \ar[d] & A'-QCoh \ar[d] \\
A-Mod \ar[r]_-{f^{*}}  & A'-Mod}$$
commutes up to a natural isomorphism. Furthermore,
for any pair of morphisms  
$$\xymatrix{A\ar[r]^{f} & A' \ar[r]^{g} & A''} $$
in $Comm(\mathcal{C})$, there is
an equality $(g\circ f)^{*}=g^{*}\circ f^{*}$.
In other words, the rule 
$$A \mapsto A-QCoh \qquad (f : A \rightarrow A') \mapsto
f^{*}$$
defines a $\mathbb{U}$-cofibrantly generated
left Quilllen presheaf on $Aff_{\mathcal{C}}=Comm(\mathcal{C})^{op}$ in the
sense of Appendix B.

We now consider for any $A\in Comm(\mathcal{C})$, the subcategory
$A-QCoh^{c}_{W}$ of $A-QCoh$, consisting of equivalences between
cofibrant objects. As these are preserved by the pullback
functors $f^{*}$, we obtain this way a new presheaf
of $\mathbb{V}$-small categories
$$\begin{array}{cccc}
QCoh^{c}_{W} : & Comm(\mathcal{C})=Aff_{\mathcal{C}}^{op} & \longrightarrow & Cat_{\mathbb{V}} \\
 & A & \mapsto & A-QCoh^{c}_{W}.
\end{array}$$
Composing with the nerve functor
$$N : Cat_{\mathbb{V}} \longrightarrow SSet_{\mathbb{V}}$$
we get  a simplicial presheaf
$$\begin{array}{cccc}
N(QCoh^{c}_{W}) : & Comm(\mathcal{C})=Aff_{\mathcal{C}}^{op} & \longrightarrow & SSet_{\mathbb{V}} \\
 & A & \mapsto & N(A-QCoh^{c}_{W}).
\end{array}$$

\begin{df}\label{d30}
The simplicial presheaf of \emph{quasi-coherent modules}
is\\ $N(QCoh^{c}_{W})$ defined above. It is denoted
by $\mathbf{QCoh}$\index{$\mathbf{QCoh}$}, and is considered
as an object in $Aff_{\mathcal{C}}^{\sim,\tau}$.
\end{df}

It is important to note that for any $A\in Comm(\mathcal{C})$,
the simplicial set $\mathbf{QCoh}(A)$ is naturally
equivalent to the nerve of $A-Mod_{W}^{c}$, the subcategory of
equivalences between cofibrant objects in $A-Mod$, and therefore
also to the nerve of $A-Mod_{W}$, the subcategory of
equivalences in $A-Mod$.
In particular,
$\pi_{0}(\mathbf{QCoh}(A))$ is in bijection with isomorphisms
classes of objects in $\mathrm{Ho}(A-Mod)$ (i.e. equivalence classes
of objects in $A-Mod)$. Furthermore, by \cite{dk3} (see also
Appendix A),
for any $x\in \mathbf{QCoh}(A)$, corresponding to
an equivalence class of $M\in A-Mod$, the connected component
of $\mathbf{QCoh}(A)$ containing $x$ is naturally equivalent to
$BAut(M)$, where $Aut(M)$ is the simplicial monoid of
self equivalences of $M$ in $A-Mod$. In particular,
we have
$$\pi_{1}(\mathbf{QCoh}(A),x)\simeq Aut_{\mathrm{Ho}(A-Mod)}(M)$$
$$\pi_{i+1}(\mathbf{QCoh}(A),x)\simeq [S^{i}M,M]_{A-Mod} \; \forall \; i>1.$$

The main result of this section is the following.

\begin{thm}\label{t2}
The simplicial presheaf $\mathbf{QCoh}$ is a stack.
\end{thm}

\begin{proof} This is a direct application
of the strictification theorem \ref{tstrict}
recalled in Appendix B.

More precisely, we use our lemma \ref{lass5} $(2)$.
Concerning finite direct sums, we have already seen (during the
proof of lemma \ref{lass5}) that
the natural functor
$$\left( -\otimes_{A\times^{h} B}A\right) \times
\left( -\otimes_{A\times^{h} B}B\right) : (A\times^{h} B)-Mod \longrightarrow
A-Mod \times B-Mod$$
is a Quillen equivalence. This implies that for any two
objects $X,Y\in Aff_{\mathcal{C}}$ the natural morphism
$$\mathbf{QCoh}(X\coprod^{\mathbb{L}}Y) \longrightarrow \mathbf{QCoh}(X)\times
\mathbf{QCoh}(Y)$$
is an equivalence. It only remains to show that
$\mathbf{QCoh}$ has the descent property with respect
to hypercovers of the type described in lemma \ref{lass5} $(2)$.
But this is nothing else than
our assumption \ref{ass5} $(3)$ together with Cor. \ref{cstrict}.
\end{proof}

An important consequence of theorem \ref{t2} is the following.

\begin{cor}\label{ct2}
Let $A \in Comm(\mathcal{C})$ be a commutative monoid,
$A-Mod_{W}$ the subcategory of equivalences in
$A-Mod$, and $N(A-Mod_{W})$ be its nerve.
Then, there exists natural isomorphisms
in $\mathrm{Ho}(SSet)$
$$\mathbb{R}\mathbf{QCoh}(A)\simeq \mathbb{R}_{\tau}\underline{Hom}
(\mathbb{R}\underline{Spec}\, A,\mathbf{QCoh})\simeq
N(A-Mod_{W}).$$
\end{cor}

To finish this section, we will describe
two important sub-stacks of $\mathbf{QCoh}$, namely the stack
of \textit{perfect modules} and the stack of \textit{vector bundles}.

For any commutative monoid $A$ in $\mathcal{C}$, we let
$\mathbf{Perf}(A)$ be the sub-simplicial set
of $\mathbf{QCoh}(A)$ consisting of all
connected components corresponding to
perfect objects in $\mathrm{Ho}(A-Mod)$ (in the sense of Def. \ref{d4}).
More precisely,
if $Iso(D)$ denotes the set of isomorphisms
classes of a category $D$, the simplicial set
$\mathbf{Perf}(A)$ is defined as the pullback
$$\xymatrix{
\mathbf{Perf}(A) \ar[r] \ar[d] & Iso(\mathrm{Ho}(A-Mod)^{perf}) \ar[d] \\
\mathbf{QCoh}(A) \ar[r] & \pi_{0}\mathbf{QCoh}(A)\simeq Iso(\mathrm{Ho}(A-Mod))}$$
where $\mathrm{Ho}(A-Mod)^{perf}$ is the full subcategory of
$\mathrm{Ho}(A-Mod)$ consisting of perfect $A$-modules
in the sense of Def. \ref{d4}.

We say that an $A$-module $M \in \mathrm{Ho}(A-Mod)$
is a \emph{rank $n$ vector bundle}, if there exists
a covering family $A \longrightarrow A'$ such that
$M\otimes^{\mathbb{L}}_{A}A'$ is isomorphic in
$\mathrm{Ho}(A'-Mod)$ to $(A')^{n}$. As we have defined
the sub simplicial set $\mathbf{Perf}(A)$ of
$\mathbf{QCoh}(A)$ we define
$\mathbf{Vect}_{n}(A)$ to be the sub simplicial set
of $\mathbf{QCoh}(A)$ consisting of
connected components corresponding to rank $n$ vector bundles.

For any morphism of commutative monoids $u : A \longrightarrow A'$, the
base change functor
$$\mathbb{L}u^{*} : \mathrm{Ho}(A-Mod) \longrightarrow \mathrm{Ho}(A'-Mod)$$
preserves perfect modules as well as rank $n$ vector bundles.
Therefore, the sub simplicial sets
$\mathbf{Vect}_{n}(A)$ and $\mathbf{Perf}$ form in fact
full sub simplicial presheaves
$$\mathbf{Vect}_{n} \subset \mathbf{QCoh} \qquad
\mathbf{Perf} \subset \mathbf{QCoh}.$$
The simplicial presheaves $\mathbf{Vect}_{n}$ and
$\mathbf{Perf}$ then define objects in $Aff_{\mathcal{C}}^{\sim,\tau}$.

\begin{cor}\label{ct22}
The simplicial presheaves $\mathbf{Perf}$ and
$\mathbf{Vect}_{n}$ are stacks.
\end{cor}

\begin{proof} Indeed, as they are full sub-simplicial presheaves
of the stack $\mathbf{QCoh}$, it is clearly enough to show that
being a perfect module and being a vector bundle or rank $n$ is
a local condition for the topology $\tau$. For vector bundles
this is obvious from the definition.

Let $A \in Comm(\mathcal{C})$ be a commutative monoid,
and $P$ be an $A$-module, such that there exists
a $\tau$-covering $A \longrightarrow B$ such that
$P\otimes_{A}^{\mathbb{L}}B$ is a perfect $B$-module.
Assume that $A \rightarrow B$ is
a cofibration, and let $B_{*}$ be its co-nerve,
considered as a
co-simplicial object in $A-Comm(\mathcal{C})$.
Let $Q$ be any $A$-module, and
define two objects in $\mathrm{Ho}(csB_{*}-Mod)$ by
$$\mathbb{R}\underline{Hom}_{A}(P,Q)_{*}:=
\mathbb{R}\underline{Hom}_{A}(P,Q)\otimes_{A}^{\mathbb{L}}B_{*}$$
$$\mathbb{R}\underline{Hom}_{A}(P,Q_{*}):=
\mathbb{R}\underline{Hom}_{A}(P,Q\otimes_{A}^{\mathbb{L}}B_{*}).$$
There is a natural morphism
$$\mathbb{R}\underline{Hom}_{A}(P,Q)_{*} \longrightarrow
\mathbb{R}\underline{Hom}_{A}(P,Q_{*}).$$
These co-simplicial objects are both cartesian, and by assumption \ref{ass5} $(3)$
applied to the $A$-modules $Q$ and $\mathbb{R}\underline{Hom}_{A}(P,Q)$
the induced morphism in
$\mathrm{Ho}(A-Mod)$
$$\int \mathbb{R}\underline{Hom}_{A}(P,Q)_{*} \simeq
\mathbb{R}\underline{Hom}_{A}(P,Q) \simeq $$
$$\simeq Holim_{n\in \Delta}\mathbb{R}\underline{Hom}_{A}(P,Q\otimes_{A}^{\mathbb{L}}B_{n})
\simeq
\int \mathbb{R}\underline{Hom}_{A}(P,Q\otimes_{A}^{\mathbb{L}}B_{*})$$
is an isomorphism.
Therefore, assumption \ref{ass5} $(3)$ implies that
the natural morphism
$$\mathbb{R}\underline{Hom}_{A}(P,Q)_{0} \longrightarrow
\mathbb{R}\underline{Hom}_{A}(P,Q_{0})$$
is an isomorphism. By definition this implies that
$$\mathbb{R}\underline{Hom}_{A}(P,Q) \otimes^{\mathbb{L}}_{A}B
\longrightarrow \mathbb{R}\underline{Hom}_{A}(P,Q\otimes^{\mathbb{L}}_{A}B)$$
is an isomorphism in $\mathrm{Ho}(A-Mod)$. In particular,
when $Q=A$ we find that the natural morphism
$$P^{\vee}\otimes_{A}^{\mathbb{L}}B \longrightarrow
\mathbb{R}\underline{Hom}_{B}(P\otimes_{A}^{\mathbb{L}}B,B)$$
is an isomorphism in $\mathrm{Ho}(B-Mod)$. As $P\otimes_{A}^{\mathbb{L}}B$
is by assumption a perfect $B$-module, we find that the natural
morphism
$$P^{\vee}\otimes_{A}^{\mathbb{L}}Q \longrightarrow
\mathbb{R}\underline{Hom}(P,Q)$$
becomes an isomorphism after base changing to $B$, and this
for any $A$-module $Q$. As $A\longrightarrow B$ is
a $\tau$-covering, we see that this implies that
$$P^{\vee}\otimes_{A}^{\mathbb{L}}Q \longrightarrow
\mathbb{R}\underline{Hom}(P,Q)$$
is always an isomorphism in $\mathrm{Ho}(A-Mod)$, for any $Q$, and thus
that $P$ is a perfect $A$-module.
\end{proof}

\begin{df}\label{d30'}
The \emph{stack of vector bundles of rank $n$} is
$\mathbf{Vect}_{n}$\index{$\mathbf{Vect}_{n}$}. The \emph{stack of perfect modules}
is $\mathbf{Perf}$\index{$\mathbf{Perf}$}.
\end{df}

The same construction can also been done in the stable
context.
For a commutative monoid $A$ in $\mathcal{C}$, we define
a category $A-QCoh^{Sp}$, of stable quasi-coherent modules
on $A$ (or equivalently on $Spec\, A$) in the following
way. Its objects are the data of
a stable $B$-module $M_{B}\in Sp(B-Mod)$ for any commutative $A$-algebra
$B\in A-Comm(\mathcal{C})$, together with
an isomorphism
$$\alpha_{u} : M_{B}\otimes_{B}C \longrightarrow M_{B'}$$
for any morphism $u : B \longrightarrow B'$ in $A-Comm(\mathcal{C})$,
such that one has
$\alpha_{v}\circ (\alpha_{u}\otimes_{B'}B'')= \alpha_{v\circ u}$
for any pair of morphisms
$$\xymatrix{B \ar[r]^-{u} & B' \ar[r]^-{v} & B''}$$
in $A-Comm(\mathcal{C})$. Such data will be denoted
by $(M,\alpha)$.
A morphism in $A-QCoh^{Sp}$, from
$(M,\alpha)$ to $(M',\alpha')$ is given by
a family of morphisms of stable
$B$-modules $f_{B} : M_{B} \longrightarrow M_{B}'$, for
any $B \in A-Comm(\mathcal{C})$, such that for any
$u : B \rightarrow B'$ in $A-Comm(\mathcal{C})$
the diagram
$$\xymatrix{
M_{B}\otimes_{B}B' \ar[r]^-{\alpha_{u}} \ar[d]_-{f_{B}\otimes_{B}B'} &
M_{B'} \ar[d]^-{f_{B'}} \\
(M_{B}')\otimes_{B}B' \ar[r]_-{\alpha_{u}'} & M_{B'}'}$$
commutes. As the categories
$Sp(A-Mod)$ and $Comm(\mathcal{C})$ are all $\mathbb{V}$-small,
so are the categories $A-QCoh^{Sp}$.

There exists a natural projection $A-QCoh^{Sp} \longrightarrow
Sp(A-Mod)$, sending $(M,\alpha)$ to $M_{A}$, and it is
straightforward to check that it is an equivalence of categories.
In particular, the model structure on $Sp(A-Mod)$ will
be transported naturally on $A-QCoh^{Sp}$ through this equivalence.

Let now $f : A \longrightarrow A'$ be a morphism of
commutative monoids in $\mathcal{C}$. There exists a
pullback functor
$$f^{*} : A-QCoh^{Sp} \longrightarrow A'-QCoh^{Sp}$$
defined by $f(M,\alpha)_{B}:=M_{B}$ for any $B\in A-Comm(\mathcal{C})$, and
for $u : B \longrightarrow B'$ in $A-Comm(\mathcal{C})$ the transition
morphism
$$f(M,\alpha)_{B}\otimes_{B}B'=M_{B}\otimes_{B}B' \longrightarrow f(M,\alpha)_{B'}=M_{B'}$$
is given by $\alpha_{u}$. By definition of the model structure on
$A-QCoh^{Sp}$, the functor
$$f^{*} : A-QCoh^{Sp} \longrightarrow A'-QCoh^{Sp}$$
is clearly a left Quillen functor. Furthermore,
for any pair of morphisms $$\xymatrix{
A\ar[r]^{f} & A' \ar[r]^{g} & A''}$$ in $Comm(\mathcal{C})$, there is
an equality $(g\circ f)^{*}=g^{*}\circ f^{*}$.
In other words, the rule
$$A \mapsto A-QCoh^{Sp} \qquad (f : A \rightarrow A') \mapsto
f^{*}$$
defines a $\mathbb{U}$-cofibrantly generated
left Quilllen presheaf on $Aff_{\mathcal{C}}=Comm(\mathcal{C})^{op}$ in the
sense of Appendix B.

We now consider for any $A\in Comm(\mathcal{C})$, the subcategory
$A-QCoh^{Sp,c}_{W}$ of $A-QCoh^{Sp}$, consisting of equivalences between
cofibrant objects. As these are preserved by the pullback
functors $f^{*}$, one gets this way a new presheaf
of $\mathbb{V}$-small categories
$$\begin{array}{cccc}
QCoh^{Sp,c}_{W} : & Comm(\mathcal{C})=Aff_{\mathcal{C}}^{op} & \longrightarrow & Cat_{\mathbb{V}} \\
 & A & \mapsto & A-QCoh^{Sp,c}_{W}.
\end{array}$$
Composing with the nerve functor
$$N : Cat_{\mathbb{V}} \longrightarrow SSet_{\mathbb{V}}$$
one gets a simplicial presheaf
$$\begin{array}{cccc}
N(QCoh^{Sp,c}_{W}) : & Comm(\mathcal{C})=Aff_{\mathcal{C}}^{op} & \longrightarrow & SSet_{\mathbb{V}} \\
 & A & \mapsto & N(A-QCoh^{Sp,c}_{W}).
\end{array}$$

\begin{df}\label{d30stable}
The simplicial presheaf of \emph{stable quasi-coherent modules}
is $N(QCoh^{Sp,c}_{W})$ defined above. It is denoted
by $\mathbf{QCoh}^{Sp}$\index{$\mathbf{QCoh}^{Sp}$}, and is considered
as an object in $Aff_{\mathcal{C}}^{\sim,\tau}$.
\end{df}

It is important to note that for any $A\in Comm(\mathcal{C})$,
the simplicial set $\mathbf{QCoh}^{Sp}(A)$ is naturally
equivalent to the nerve of $A-Mod_{W}^{Sp,c}$, the subcategory of
equivalences between cofibrant objects in $Sp(A-Mod)$, and therefore
also to the nerve of $Sp(A-Mod)_{W}$, the subcategory of
equivalences in $Sp(A-Mod)$.
In particular,
$\pi_{0}(\mathbf{QCoh}^{Sp}(A))$ is in bijection with isomorphisms
classes of objects in $\mathrm{Ho}(Sp(A-Mod))$ (i.e. equivalence classes
of objects in $Sp(A-Mod))$. Furthermore, by \cite{dk3} (see also
Appendix A),
for any $x\in \mathbf{QCoh}^{Sp}(A)$, corresponding to
an equivalence class of $M\in Sp(A-Mod)$, the connected component
of $\mathbf{QCoh}^{Sp}(A)$ containing $x$ is naturally equivalent to
$BAut(M)$, where $Aut(M)$ is the simplicial monoid of
self equivalences of $M$ in $Sp(A-Mod)$. In particular, we have
$$\pi_{1}(\mathbf{QCoh}^{Sp}(A),x)\simeq Aut_{\mathrm{Ho}(Sp(A-Mod))}(M) $$
$$\pi_{i+1}(\mathbf{QCoh}^{Sp}(A),x)\simeq [S^{i}M,M]_{Sp(A-Mod)} \; \forall \; i>1.$$

The same proof as theorem \ref{t2}, but
based on Proposition \ref{pdescstable} gives the
following stable version.

\begin{thm}\label{t2stable}
Assume that the two conditions are satisfied.
\begin{enumerate}
\item  The suspension functor $S : \mathrm{Ho}(\mathcal{C}) \longrightarrow
 \mathrm{Ho}(\mathcal{C})$ is fully faithful.
\item For all $\tau$-covering
family $\{U_{i} \rightarrow X\}$ in $Aff_{\mathcal{C}}$,
each morphism $U_{i} \rightarrow X$ is flat in the
sense of Def. \ref{d7-}.
\end{enumerate}
Then, the simplicial presheaf $\mathbf{QCoh}^{Sp}$ is a stack.
\end{thm}

For any commutative monoid $A$ in $\mathcal{C}$, we let
$\mathbf{Perf}^{Sp}(A)$ be the sub-simplicial set
of $\mathbf{QCoh}^{Sp}(A)$ consisting of all
connected components corresponding to
perfect objects in $\mathrm{Ho}(Sp(A-Mod))$ (in the sense of Def. \ref{d4}).
This defines
a full sub-prestack $\mathbf{Perf}^{Sp}$
of $\mathbf{QCoh}^{Sp}$. Then, the same argument as for Cor. \ref{ct22}
gives the following corollary.

\begin{cor}\label{ct22stable}
Assume that the two conditions are satisfied.
\begin{enumerate}
\item  The suspension functor $S : \mathrm{Ho}(\mathcal{C}) \longrightarrow
 \mathrm{Ho}(\mathcal{C})$ is fully faithful.
\item For all $\tau$-covering
family $\{U_{i} \rightarrow X\}$ in $Aff_{\mathcal{C}}$,
each morphism $U_{i} \rightarrow X$ is flat in the
sense of Def. \ref{d7-}.
\end{enumerate}
The simplicial presheaf $\mathbf{Perf}^{Sp}$ is a stack.
\end{cor}

This justifies the following definition.

\begin{df}\label{d30'stable}
Under the condition of Cor. \ref{ct22stable},
the \emph{stack of stable perfect modules}
is $\mathbf{Perf}^{Sp}$\index{$\mathbf{Perf}^{Sp}$}.
\end{df}

We finish this section by the standard description of
$\mathbf{Vect}_{n}$ as the classifying stack of
the group stack $\mathbf{Gl}_{n}$, of invertible $n$ by $n$ matrices.
For this, we notice that the natural morphism of stacks
$* \longrightarrow \mathbf{Vect}_{n}$, pointing the trivial
rank $n$ vector bundle, induces an isomorphism of sheaves of sets
$*\simeq \pi_{0}(\mathbf{Vect}_{n})$. Therefore,
the stack $\mathbf{Vect}_{n}$ can be written as the geometric
realization of the homotopy nerve of the morphism
$* \longrightarrow \mathbf{Vect}_{n}$. In other words, we can
find a Segal groupoid object $X_{*}$, such that $X_{0}=*$,
and with $|X_{*}|\simeq \mathbf{Vect}_{n}$. Furthermore,
the object $X_{1}$ is naturally equivalent to the
loop stack $\Omega_{*}\mathbf{Vect}_{n}:=*\times^{h}_{\mathbf{Vect}_{n}}*$.
By construction and by Prop. \ref{papp2}, this loop stack can be described
as the simplicial presheaf
$$\begin{array}{cccc}
\Omega_{*}\mathbf{Vect}_{n} : & Comm(\mathcal{C}) & \longrightarrow & SSet \\
 & A & \mapsto & Map_{\mathcal{C}}'(1^{n},A^{n}),
\end{array}$$
where $Map_{\mathcal{C}}'(1^{n},A^{n})$ is the sub simplicial set
of the mapping space $Map_{\mathcal{C}}(1^{n},A^{n})$
consisting of all connected components
corresponding to automorphisms in
$$\pi_{0}Map_{\mathcal{C}}(1^{n},A^{n})\simeq
\pi_{0}Map_{A-Mod}(A^{n},A^{n})\simeq [A^{n},A^{n}]_{A-Mod}.$$

The important fact concerning the stack
$\Omega_{*}\mathbf{Vect}_{n}$ is the following
result.

\begin{prop}\label{p24}
\begin{enumerate}
\item
The stack $\Omega_{*}\mathbf{Vect}_{n}$ is representable
and the morphism
$\Omega_{*}\mathbf{Vect}_{n} \longrightarrow *$ is
formally smooth.

\item
If Moreover $\mathbf{1}$ is finitely presented in $\mathcal{C}$, then the morphism
$\Omega_{*}\mathbf{Vect}_{n} \longrightarrow *$ is finitely presented (and thus smooth by $(1)$).
\end{enumerate}
\end{prop}

\begin{proof} $(1)$ We start by defining a larger stack
$\mathbf{M}_{n}$, of $n$ by $n$ matrices. We set
$$\begin{array}{cccc}
\mathbf{M}_{n} : & Comm(\mathcal{C}) & \longrightarrow & SSet \\
 & A & \mapsto & Map_{\mathcal{C}}(1^{n},A^{n}).
\end{array}$$

This stack is representable, as it is isomorphic in
$\mathrm{St}(\mathcal{C},\tau)$ to
$\mathbb{R}\underline{Spec}\, B$, where
$B=\mathbb{L}F(\mathbf{1}^{n^{2}})$ is the free
commutative monoid generated by the object $\mathbf{1}^{n^{2}} \in \mathcal{C}$.
We claim that the natural inclusion morphism
$$\Omega_{*}\mathbf{Vect}_{n} \longrightarrow \mathbf{M}_{n}$$
is $(-1)$-representable and a formally \'etale morphism.
Indeed, let $A$ be any commutative monoid,
$$x : X:=\mathbb{R}\underline{Spec}\, A \longrightarrow
\mathbf{M}_{n}$$ 
be a morphism of stacks, and
let us consider the stack
$$F:=\Omega_{*}\mathbf{Vect}_{n}\times^{h}_{\mathbf{M}_{n}} X\longrightarrow 
X.$$
The point $x$ corresponds via the Yoneda lemma to
a morphism $u : A^{n} \longrightarrow A^{n}$ in $\mathrm{Ho}(A-Mod)$.
Now, for any commutative monoid $A'$, the natural morphism
$$\mathbb{R}F(A') \longrightarrow (\mathbb{R}\underline{Spec}\, A)(A')$$
identifies $\mathbb{R}F(A')$ with the sub simplicial set
of $(\mathbb{R}\underline{Spec}\, A)(A')\simeq Map_{Comm(\mathcal{C})}(A,A')$
consisting of all connected components corresponding to
morphisms $A \longrightarrow A'$ in $\mathrm{Ho}(Comm(\mathcal{C}))$ such that
$$u\otimes_{A}^{\mathbb{L}}A' : (A')^{n} \longrightarrow
(A')^{n}$$
is an isomorphism in $\mathrm{Ho}(A'-Mod)$. Considering $u$ as an element
of $[A^{n},A^{n}]\simeq M_{n}(\pi_{0}(A))$, we can consider
its determinant $d(u)\in \pi_{0}(A)$. Then, using notations
from Def. \ref{dloc}, we clearly have an isomorphism of stacks
$$F\simeq \mathbb{R}\underline{Spec}\, (A[d(u)^{-1}]).$$
By Prop. \ref{ploc3} this shows that the morphism
$$F \longrightarrow \mathbb{R}\underline{Spec}\, A$$
is a formally \'etale morphism between representable stacks.
As $\mathbf{M}_{n}$ is representable this implies that
$\Omega_{*}\mathbf{Vect}_{n}$ is
a representable stack and that the morphism
$$\Omega_{*}\mathbf{Vect}_{n} \longrightarrow \mathbf{M}_{n}$$
is formally \'etale (it is also a flat monomorphism by
\ref{ploc2}).

Moreover, we have $\mathbf{M}_{n} \simeq
\mathbb{R}\underline{Spec}\, B$, where $B:=\mathbb{L}F(\mathbf{1}^{n^{2}})$
is the derived free commutative monoid over $\mathbf{1}^{n^{2}}$. This
implies that $\mathbb{L}_{B}$ is a free $B$-module of rank
$n^{2}$, and therefore that the morphism $\mathbf{1}
\longrightarrow B$ is formally smooth in the sense of Def.
\ref{d7}. By composition, we find that
$\Omega_{*}\mathbf{Vect}_{n} \longrightarrow *$
is a formally smooth morphism as required. \\

$(2)$ This follows from $(1)$, Prop.  \ref{ploc2} $(2)$ and the fact that
$B=\mathbb{L}F(\mathbf{1}^{n^{2}})$ is a finitely presented object in $Comm(\mathcal{C})$. \end{proof}

\begin{df}\label{d31}
\begin{enumerate}
\item
The stack $\Omega_{*}\mathbf{Vect}_{n}$ is denoted
by $\mathbf{Gl}_{n}$\index{$\mathbf{Gl}_{n}$!linear group stack of rank $n$}, and is called
the \emph{linear group stack of rank $n$}. The stack
$\mathbf{Gl}_{1}$ is denoted by $\mathbb{G}_{m}$\index{$\mathbb{G}_{m}$!multiplicative group stack}, and is called the
\emph{multiplicative group stack}.
\item The stack
$\mathbf{M}_{n}$\index{$\mathbf{M}_{n}$!stack of $n\times n$ matrices} defined during the proof of
\ref{p24} $(1)$ is called the \emph{stack of $n\times n$ matrices}. The stack
$\mathbf{M}_{1}$ is denoted by $\mathbb{G}_{a}$\index{$\mathbb{G}_{a}$!additive group stack}, and is called
the \emph{additive group stack}.
\end{enumerate}
\end{df}

Being a stack of loops, the stack $\mathbf{Gl}_{n}=\Omega_{*}\mathbf{Vect}_{n}$
has a natural group structure, encoded in the fact that
it is the $X_{1}$ of a Segal groupoid object $X_{*}$
with $X_{0}=*$. Symbolically, we will simply write
$$B\mathbf{Gl}_{n}:=|X_{*}|.$$

Our conclusion is that the stack
$\mathbf{Vect}_{n}$ can be written as
$B\mathbf{Gl}_{n}$, where $\mathbf{Gl}_{n}$ is a
\emph{formally smooth representable group stack}. Furthermore this group
stack is smooth when the unit $\mathbf{1}$ is finitely presented.
As a corollary we
get the following geometricity result
on $\mathbf{Vect}_{n}$.

\begin{cor}\label{cp24}
Assume
that the unit $1$ is a finitely presented object in $\mathcal{C}$.
Assume that all smooth morphisms in $Comm(\mathcal{C})$
belong to \textbf{P}. Then, the stack
$\mathbf{Vect}_{n}$ is $1$-geometric, the morphism
$\mathbf{Vect}_{n} \longrightarrow *$ is in \textbf{P} and finitely presented, and furthermore its
diagonal is a $(-1)$-representable morphism.
\end{cor}

\begin{proof} The $1$-geometricity statement is a consequence of Prop. \ref{p13}
and the fact that the natural morphism $* \longrightarrow \mathbf{Vect}_{n}$ is
a $1$-\textbf{P}-atlas. That $*$ is a $1$-\textbf{P}-atlas also implies that
$\mathbf{Vect}_{n} \longrightarrow *$ is in \textbf{P} and finitely presented.
The statement concerning the diagonal follows from the fact that
$\mathbf{Gl}_{n}$ is a representable stack and the locality of
representable objects Prop. \ref{c0}. \end{proof}

We finish with an analogous situation for perfect modules. We let
$K$ be a perfect object in $\mathcal{C}$, and we define a stack
$\mathbb{R}\underline{End}(K)$ in the following way. We chose a
cofibrant replacement $QK$ of $K$, and $\Gamma_{*}$ a simplicial resolution functor
on $\mathcal{C}$. One sets
$$\begin{array}{cccc}
\mathbb{R}\underline{End}(K) : & Comm(\mathcal{C}) & \longrightarrow & SSet_{\mathbb{V}} \\
 & A & \mapsto & Hom(QK,\Gamma_{*}(QK\otimes A)).
\end{array}$$
Note that for any $A$ the simplicial set $\mathbb{R}\underline{End}(K)(A)$ is
naturally equivalent to $$Map_{A-Mod}(K\otimes^{\mathbb{L}}A,K\otimes^{\mathbb{L}}A).$$

\begin{lem}\label{lendperf}
The simplicial presheaf $\mathbb{R}\underline{End}(K)(A) \in Aff_{\mathcal{C}}^{\sim,\tau}$
is a stack. It is furthermore representable.
\end{lem}

\begin{proof} This is clear as $K$ being perfect one sees that there exists
an isomorphism in $\mathrm{Ho}(SPr(Aff_{\mathcal{C}}))$
$$\mathbb{R}\underline{End}(K)\simeq \mathbb{R}\underline{Spec}\, B,$$
where $B:=\mathbb{L}F(K\otimes^{\mathbb{L}}K^{\vee})$ is the derived
free commutative monoid on the object $K\otimes^{\mathbb{L}}K^{\vee}$.
\end{proof}

We now define a sub-stack $\mathbb{R}\underline{Aut}(K)$ of
$\mathbb{R}\underline{End}(K)$. For a commutative monoid $A\in Comm(\mathcal{C})$,
we define $\mathbb{R}\underline{Aut}(K)(A)$ to be the union of connected
components of $\mathbb{R}\underline{End}(K)(A)$ corresponding to isomorphisms
in
$$\pi_{0}(\mathbb{R}\underline{End}(K)(A))\simeq [K\otimes^{\mathbb{L}}A,K\otimes^{\mathbb{L}}A]_{A-Mod}.$$
This clearly defines a full sub-simplicial presheaf $\mathbb{R}\underline{Aut}(K)$
of $\mathbb{R}\underline{End}(K)$, which is a sub-stack as one can see easily
using Cor. \ref{c0''}.

\begin{prop}\label{pautperf}
Assume that $\mathcal{C}$ is a stable model category.
Then, the stack $\mathbb{R}\underline{Aut}(K)$ is representable. Furthermore the morphism
$\mathbb{R}\underline{Aut}(K) \longrightarrow \mathbb{R}\underline{End}(K)$
is a formal Zariski open immersion, and $\mathbb{R}\underline{Aut}(K) \longrightarrow
*$ is fp. If furthermore $\mathbf{1}$ is finitely presented in $\mathcal{C}$ then
$\mathbb{R}\underline{Aut}(K) \longrightarrow*$ is a perfect morphism.
\end{prop}

\begin{proof} It is the same as \ref{p24} but using
the construction $A_{K}$ and Prop. \ref{psupp}, instead of the standard
localization $A[a^{-1}]$. More precisely, for a representable stack
$X:=\mathbb{R}\underline{Spec}\, A$ and a morphism
$x : X \longrightarrow \mathbb{R}\underline{End}(K)$, the homotopy
pullback
$$\mathbb{R}\underline{Aut}(K)\times^{h}_{\mathbb{R}\underline{End}(K)}X \longrightarrow X$$
is isomorphic in $\mathrm{Ho}(Aff_{\mathcal{C}}^{\sim,\tau}/X)$ to
$$\mathbb{R}\underline{Spec}\, A_{E} \longrightarrow \mathbb{R}\underline{Spec}\, A,$$
where $E$ is the homotopy cofiber of the endomorphism
$x : K\otimes^{\mathbb{L}}A \longrightarrow K\otimes^{\mathbb{L}}A$ corresponding to the point $x$.
\end{proof}

\chapter{Geometric stacks: Infinitesimal theory}\label{partI.4}

As in the previous chapter, we fix once for all a
HAG context $(\mathcal{C},\mathcal{C }_{0},\mathcal{A},\tau,\textbf{P})$. \\

In this chapter we will assume furthermore that the suspension
functor
$$S : \mathrm{Ho}(\mathcal{C}) \longrightarrow \mathrm{Ho}(\mathcal{C})$$
is fully faithful. In particular, for any commutative monoid $A\in Comm(\mathcal{C})$,
the stabilization functor
$$S_{A} : \mathrm{Ho}(A-Mod) \longrightarrow \mathrm{Ho}(Sp(A-Mod))$$
is fully faithful (see \ref{lstabmod}). We will therefore forget to mention
the functor $S_{A}$ and simply consider
$A$-modules as objects in $\mathrm{Ho}(Sp(A-Mod))$, corresponding to
$0$-connective objects.

\section{Tangent stacks and cotangent complexes}

We consider the initial commutative monoid $\mathbf{1} \in Comm(\mathcal{C})$.
It can be seen as a module over itself, and gives rise to
a trivial square zero extension $\mathbf{1}\oplus \mathbf{1}$ (see \S \ref{Ider}).

\begin{df}\label{d14}
The \emph{dual numbers over $\mathcal{C}$} is the commutative
monoid $\mathbf{1}\oplus \mathbf{1}$, and is denoted by
$\mathbf{1}[\epsilon]$\index{$\mathbf{1}[\epsilon]$!monoid of dual numbers}. The corresponding representable
stack is denoted by
$$\mathbb{D}_{\epsilon}:=\mathbb{R}\underline{Spec}\, (\mathbf{1}[\epsilon])$$
and is called the \emph{infinitesimal disk}\index{$\mathbb{D}_{\epsilon}$!infinitesimal disk}.
\end{df}

Of course, as every trivial square zero extension the natural morphism
$\mathbf{1} \longrightarrow \mathbf{1}[\epsilon]$ possesses a natural retraction
$\mathbf{1}[\epsilon] \longrightarrow \mathbf{1}$. On the level of representable stacks this
defines a natural global point
$$* \longrightarrow \mathbb{D}_{\epsilon}.$$
We recal that $\mathrm{Ho}(Aff^{\sim,\tau}_{\mathcal{C}})$ being the homotopy category 
of a model topos has internal $Hom$'s objects $\mathbb{R}_{\tau}\underline{\mathcal{H}om}(-,-)$. 
They satisfy the usual adjunction isomorphisms
$$Map_{Aff^{\sim,\tau}_{\mathcal{C}}}(F,\mathbb{R}_{\tau}\underline{\mathcal{H}om}(G,H))\simeq
Map_{Aff^{\sim,\tau}_{\mathcal{C}}}(F\times^{h}G,H).$$

\begin{df}\label{d15}
Let $F \in \mathrm{St}(\mathcal{C},\tau)$ be a stack. The \emph{tangent
stack of $F$}\index{tangent stack} is defined to be
$$TF:=\mathbb{R}_{\tau}\underline{\mathcal{H}om}(\mathbb{D}_{\epsilon},F).$$
The natural morphism $* \longrightarrow  \mathbb{D}_{\epsilon}$ induces
a well defined projection
$$\pi : TF \longrightarrow F,$$
and the projection $\mathbb{D}_{\epsilon} \longrightarrow *$ induces
a natural section of $\pi$
$$e : F \longrightarrow TF.$$
\end{df}

For any commutative monoid $A\in Comm(\mathcal{C})$, it is
clear that there is a natural equivalence
$$A\otimes^{\mathbb{L}}(\mathbf{1}[\epsilon])\simeq A\oplus A,$$
where $A\oplus A$ is the trivial square zero extension
of $A$ by itself. We will simply denote
$A\oplus A$ by $A[\epsilon]$, and
$$\mathbb{D}_{\epsilon}^{A}:=\mathbb{R}\underline{Spec}\, (A[\epsilon])
\simeq \mathbb{R}\underline{Spec}\, (A)\times \mathbb{D}_{\epsilon}$$
the infinitesimal disk over $A$.

With these notations, and for any fibrant object
$F \in Aff_{\mathcal{C}}^{\sim,\tau}$, the stack
$TF$ can be described as the following simplicial presheaf
$$\begin{array}{cccc}
TF : & Comm(\mathcal{C})=Aff_{\mathcal{C}}^{op} & \longrightarrow & SSet_{\mathbb{V}} \\
 & A & \mapsto & F(A[\epsilon]).
\end{array}$$
Note that if $F$ is fibrant, then so is $TF$ as defined above. 
In other words for any $A\in Comm(\mathcal{C})$ there exists 
a natural equivalence of simplicial sets
$$\mathbb{R}TF(A)\simeq \mathbb{R}F(A[\epsilon]).$$

\begin{prop}\label{p15}
The functor $F \mapsto TF$ commutes with $\mathbb{V}$-small homotopy
limits.
\end{prop}

\begin{proof} This is clear as
$\mathbb{R}_{\tau}\underline{\mathcal{H}om}(H,-)$ always commutes with
homotopy limits for any $H$. \end{proof}

Let $F \in \mathrm{St}(\mathcal{C},\tau)$ be a stack,
$A \in Comm(\mathcal{C})$ a commutative monoid and
$$x : \mathbb{R}\underline{Spec}\, A \longrightarrow F$$
be a $A$-point. Let $M$ be an $A$-module, and let
$A\oplus M$ be the trivial square zero extension of $A$ by $M$.
Let us fix the following notations
$$X:=\mathbb{R}\underline{Spec}\, A \qquad
X[M]:=\mathbb{R}\underline{Spec}\, (A\oplus M).$$
The natural augmentation $A\oplus M \longrightarrow A$
gives rise to a
natural morphism of stacks
$X \longrightarrow X[M]$.

\begin{df}\label{d16}
Let $x : X \longrightarrow F$ be as above. The simplicial
set of \emph{derivations from $F$ to $M$ at the point $x$} is
defined by\index{$\mathbb{D}er_{F}(X,M)$}
$$\mathbb{D}er_{F}(X,M):=
Map_{X/Aff_{\mathcal{C}}^{\sim,\tau}}(X[M],F)
\in \mathrm{Ho}(SSet_{\mathbb{V}}).$$
It will also denoted by 
$\mathbb{D}er_{F}(x,M)$.
\end{df}

As the construction $M \mapsto X[M]$ is functorial in $M$, we
get this way a well defined functor
$$\mathbb{D}er_{F}(X,-) : \mathrm{Ho}(A-Mod) \longrightarrow \mathrm{Ho}(SSet_{\mathbb{V}}).$$
This functor is furthermore naturally compatible with the
$\mathrm{Ho}(SSet)$-enrichment, in the sense that there exists
natural morphisms
$$Map_{A-Mod}(M,N) \longrightarrow Map_{SSet}(\mathbb{D}er_{F}(X,M),
\mathbb{D}er_{F}(X,N))$$
which are compatible with compositions.
Note that the morphism $X \longrightarrow X[M]$ has
a natural retraction $X[M] \longrightarrow X$, and therefore
that the simplicial set $\mathbb{D}er_{F}(X,M)$ has a distinguished
base point, the trivial derivation. The functor
$\mathbb{D}er_{F}(X,-)$ takes its values in the homotopy category
of pointed simplicial sets.

We can also describe the functor $\mathbb{D}er_{F}(X,M)$ using
functors of points in the following way. Let $F \in Aff_{\mathcal{C}}^{\sim,\tau}$
be a fibrant object. For any commutative monoid $A$, any $A$-module
$M$ and any point $x\in F(A)$, we consider the standard homotopy fiber\footnote{
The standard homotopy fiber product of a diagram
$\xymatrix{X \ar[r] & Z & \ar[l] Y}$ of fibrant simplicial sets
is defined for example by
$$X\times^{h}_{Z}Y:=(X\times Y)\times_{Z\times Z}Z^{\Delta^{1}}.$$
} of
$$F(A\oplus M) \longrightarrow F(A)$$
at the point $x$. This clearly defines a functor
$$\begin{array}{ccc}
A-Mod & \longrightarrow & SSet_{\mathbb{V}} \\
M & \mapsto & Hofiber\left(F(A\oplus M) \rightarrow F(A)\right)
\end{array}$$
which is a lift of the functor considered above
$$\mathbb{D}er_{F}(X,-) : \mathrm{Ho}(A-Mod) \longrightarrow  \mathrm{Ho}(SSet_{\mathbb{V}}),$$
where $x\in F(A)$ corresponds via the Yoneda lemma to
a morphism $X=\mathbb{R}\underline{Spec}\, A \longrightarrow F$.
The fact that the functor $\mathbb{D}er_{F}(X,-)$ has a natural lift
as above is important, as
it then makes sense to say that it commutes with homotopy limits or homotopy
colimits. The functor $\mathbb{D}er_{F}(X,-)$
can be considered as an object in $\mathrm{Ho}((A-Mod^{op})^{\wedge})$, the homotopy
category of pre-stacks on the model category $A-Mod^{op}$, as defined
in \cite[\S 4.1]{hagI} (see also \S \ref{remhagI}). Restricting to the subcategory
$A-Mod_{0}$ we get an object
$$\mathbb{D}er_{F}(X,-)\in \mathrm{Ho}((A-Mod_{0}^{op})^{\wedge}).$$
In the sequel the functor
$\mathbb{D}er_{F}(X,-)$ will always be considered
as an object in $\mathrm{Ho}((A-Mod_{0}^{op})^{\wedge})$.

\begin{df}\label{d17}
Let $F$ be any stack and let $A \in \mathcal{A}$.
\begin{enumerate}
\item
Let $x : X:=\mathbb{R}\underline{Spec}\, A \longrightarrow F$
be an $A$-point. We say that \emph{$F$ has a
cotangent complex at $x$} if
there exists an integer $n$,
an $(-n)$-connective stable $A$-module $\mathbb{L}_{F,x} \in \mathrm{Ho}(Sp(A-Mod))$, and
an isomorphism in
$\mathrm{Ho}((A-Mod^{op}_{0})^{\wedge})$
$$\mathbb{D}er_{F}(X,-)\simeq
\mathbb{R}\underline{h}_{s}^{\mathbb{L}_{F,x}}.$$

\item If $F$ has a cotangent complex at $x$, the stable $A$-module
$\mathbb{L}_{F,x}$ is then called the
\emph{cotangent complex of $F$ at $x$}\index{cotangent complex! of a stack at a point}.
\item
If $F$ has a cotangent complex at $x$,
the \emph{tangent complex of $F$ at $x$}\index{tangent complex! of a stack at a point} is then the
stable $A$-module
$$\mathbb{T}_{F,x}:=\mathbb{R}\underline{Hom}_{A}^{Sp}(
\mathbb{L}_{F,x},A)\in \mathrm{Ho}(Sp(A-Mod)).$$
\end{enumerate}
\end{df}

In other words, the existence of a cotangent complex
of $F$ at $x : X:=\mathbb{R}\underline{Spec}\, A \longrightarrow F$
is equivalent to the co-representability of the
functor
$$\mathbb{D}er_{F}(X,-) : A-Mod_{0} \longrightarrow  SSet_{\mathbb{V}},$$
by some $(-n)$-connective object $\mathbb{L}_{F,x} \in \mathrm{Ho}(Sp(A-Mod))$.
The fact that $\mathbb{L}_{F,x}$ is well defined
is justified by our Prop. \ref{pstabmod}. \\

The first relation between cotangent complexes and the tangent
stack is given by the following proposition.

\begin{prop}\label{p16}
Let $F$ be a stack and $x : X:=\mathbb{R}\underline{Spec}\, A \longrightarrow F$
be an $A$-point with $A\in \mathcal{A}$.
If $F$ has a cotangent complex $\mathbb{L}_{F,x}$ at the
point $x$ then there exists a natural isomorphism
in $\mathrm{Ho}(SSet_{\mathbb{V}})$
$$\mathbb{R}\underline{Hom}_{Aff_{\mathcal{C}}^{\sim,\tau}/F}(X,TF)\simeq
Map_{Sp(A-Mod)}(A,\mathbb{T}_{F,x})\simeq Map_{Sp(A-Mod)}(\mathbb{L}_{F,x},A).$$
\end{prop}

\begin{proof} By definition of $TF$,
there is a natural isomorphism
$$\mathbb{R}\underline{Hom}_{Aff_{\mathcal{C}}^{\sim,\tau}/F}(X,TF)\simeq
\mathbb{R}\underline{Hom}_{Aff_{\mathcal{C}}^{\sim,\tau}/X}(X[A],TF)
\simeq Map_{Sp(A-Mod)}(\mathbb{L}_{F,x},A).$$
Moreover, by definition of the tangent complex we have
$$Map_{Sp(A-Mod)}(\mathbb{L}_{F,x},A)\simeq Map_{Sp(A-Mod)}(A,\mathbb{T}_{F,x}).$$
\end{proof}

Now, let $F$ be a stack, and $u$ a morphism
in $\mathrm{Ho}(Aff_{\mathcal{C}}^{\sim,\tau}/F)$
$$\xymatrix{
Y\ar[rr]^-{u} \ar[rd]_-{y} & & X \ar[ld]^-{y} \\
 & F, & }$$
with $X=\mathbb{R}\underline{Spec}\, A$ and
$Y=\mathbb{R}\underline{Spec}\, B$
belonging to $\mathcal{A}$.
Let $M\in B-Mod_{0}$, which is also an $A$-module
by the forgetful functor $A-Mod_{0} \longrightarrow B-Mod_{0}$.
There is a commutative diagram of
commutative monoids
$$\xymatrix{A\oplus M \ar[r] \ar[d] &  B\oplus M \ar[d] \\
 A \ar[r] & B,}$$
inducing a commutative square of representable stacks
$$\xymatrix{X[M] &  Y[M] \ar[l] \\
 X\ar[u] & Y. \ar[u] \ar[l]}$$
This implies the existence of a natural morphism
in $\mathrm{Ho}(SSet_{\mathbb{V}})$
$$\mathbb{D}er_{F}(Y,M) \longrightarrow \mathbb{D}er_{F}(X,M).$$
If the stack $F$ has cotangent complexes at both points
$x$ and $y$, Prop. \ref{pstabmod} induces a well defined morphism
in $\mathrm{Ho}(Sp(B-Mod))$
$$u^{*} : \mathbb{L}_{F,x}\otimes^{\mathbb{L}}_{A}B \longrightarrow
\mathbb{L}_{F,y}.$$
Of course, we have $(u\circ v)^{*}=v^{*}\circ u^{*}$
whenever this formula makes sense.

In the same way, the construction of $\mathbb{L}_{F,x}$
is functorial in $F$. Let $f : F \longrightarrow F'$ be
a morphism of stacks, $A\in \mathcal{A}$, and
$x : \mathbb{R}\underline{Spec}\, A \longrightarrow F$
be an $A$-point with image
$f(x) : \mathbb{R}\underline{Spec}\, A \longrightarrow F'$. Then,
for any
$A$-module $M$, there is a natural morphism
$$\mathbb{D}er_{F}(X,M) \longrightarrow \mathbb{D}er_{F'}(X,M).$$
Therefore, if  $F$ has a cotangent complex
at $x$ and $F'$ has a cotangent complex at $x'$, we get
a natural morphism in $\mathrm{Ho}(Sp(A-Mod))$
$$df_{x} : \mathbb{L}_{F',f(x)}\longrightarrow
\mathbb{L}_{F,x},$$
called the \emph{differential of $f$ at $x$}. Once again,
we have $d(f\circ g)_{x}=dg_{x}\circ df_{x}$ each time
this formula makes sense (this is the chain rule).
Dually, we also get
by duality the \emph{derivative of $f$ at $x$}
$$Tf_{x} : \mathbb{T}_{F,x} \longrightarrow
\mathbb{T}_{F',f(x)}.$$

\begin{df}\label{d18}
A stack $F$ \emph{has a global cotangent complex relative
to the HA context $(\mathcal{C},\mathcal{C}_{0},\mathcal{A})$} (or simply
\emph{has a cotangent complex}\index{cotangent complex!(global) of a stack} when there is no
ambiguity on the context) if the following two
conditions are satisfied.
\begin{enumerate}
\item For any $A\in \mathcal{A}$, and any point
$x : \mathbb{R}\underline{Spec}\, A \longrightarrow F$, the stack
$F$ has a cotangent complex $\mathbb{L}_{F,x}$ at $x$.
\item For any morphism $u : A \longrightarrow B$ in $\mathcal{A}$,
and
any morphism in $\mathrm{Ho}(Aff_{\mathcal{C}}^{\sim,\tau}/F)$
$$\xymatrix{
Y:=\mathbb{R}\underline{Spec}\, B\ar[rr]^-{u} \ar[rd]_-{y} & &
X:=\mathbb{R}\underline{Spec}\, A \ar[ld]^-{x} \\
 & F, & }$$
the induced morphism
$$u^{*} : \mathbb{L}_{F,x}\otimes^{\mathbb{L}}_{A}B \longrightarrow
\mathbb{L}_{F,y}$$
is an isomorphism in $\mathrm{Ho}(Sp(B-Mod))$.
\end{enumerate}
\end{df}

As a corollary of Prop. \ref{p1} and the standard
properties of derivations  any
representable stack has a cotangent complex.

\begin{prop}\label{p16-}
Any representable stack
$F=\mathbb{R}\underline{Spec}A$ has a global cotangent complex.
\end{prop}

\begin{proof} This is nothing else than the
existence of a universal derivation as proved in
Prop. \ref{p1}. \end{proof}

If $X=\mathbb{R}\underline{Spec}\, A$
is a representable stack in $\mathcal{A},$
and $x : X \rightarrow X$ is the identity, then
the stable $A$-module $\mathbb{L}_{X,x}$ is naturally
isomorphic in $\mathrm{Ho}(Sp(A-Mod))$ to
$\mathbb{L}_{A}$. More generally, for any morphism
$A \longrightarrow B$ with $B\in \mathcal{A}$, corresponding
to $y : \mathbb{R}\underline{Spec}\, B \longrightarrow X$,
the $B$-module $\mathbb{L}_{X,y}$ is naturally isomorphic
in $\mathrm{Ho}(Sp(B-Mod))$
to $\mathbb{L}_{A}\otimes_{A}^{\mathbb{L}}B$. \\

The next proposition explains the relation between
the tangent stack and the global cotangent complex when it exists.
It is a globalization of Prop.\ref{p16}.

\begin{prop}\label{p17}
Let $F$ be a stack having a cotangent complex. Let
$x : X=\mathbb{R}\underline{Spec}\, A \longrightarrow F$
be any morphism, and
$$TF_{x}:=TF\times_{F}^{h}X \longrightarrow X$$
the natural projection. Let $A\longrightarrow B$ be
a morphism with $B\in \mathcal{A}$, corresponding to  a morphism of representable stacks
$X=\mathbb{R}\underline{Spec}\, A \longrightarrow Y=\mathbb{R}\underline{Spec}\, B$. Then, there
exists a natural isomorphism in $\mathrm{Ho}(SSet)$
$$\mathbb{R}TF_{x}(B)\simeq Map_{Aff_{\mathcal{C}}^{\sim,\tau}/F}(Y,TF)\simeq
Map_{Sp(A-Mod)}(\mathbb{L}_{F,x},B).$$
\end{prop}

\begin{proof} This is a reformulation of Prop. \ref{p16}, and the fact that
$$Map_{Sp(A-Mod)}(\mathbb{L}_{F,x},B)\simeq Map_{Sp(B-Mod)}(\mathbb{L}_{F,x}\otimes_{A}^{\mathbb{L}}B,B).$$
\end{proof}

\begin{rmk}\label{r10}
\emph{Of course, the isomorphism of proposition Prop. \ref{p17}
is functorial in $F$.}
\end{rmk}

\begin{prop}\label{p18}
Let $F$ be an $n$-geometric stack. We assume that
for any  $A\in \mathcal{A}$, and
any point $x : X:=\mathbb{R}\underline{Spec}\, A \longrightarrow F$,
and any $A$-module $M\in A-Mod_{0}$, the natural morphism
$$\mathbb{D}er_{F}(X,M) \simeq \mathbb{D}er_{F}(X,\Omega S(M))\longrightarrow
\Omega \mathbb{D}er_{F}(X,S(M))$$
is an isomorphism in $\mathrm{Ho}(SSet)$.
Then $F$ has a global cotangent complex, which is furthermore
$(-n)$-connective.
\end{prop}

\begin{proof} The proof is by induction on $n$. For
$n=-1$ this is Prop. \ref{p16-} and does not use our
exactness condition on the functor $\mathbb{D}er_{F}(X,-)$.

Let $n\geq 0$ be an integer and $F$ be an $n$-geometric
stack. Let $X=\mathbb{R}\underline{Spec}\, A$ be a representable stack
in $\mathcal{A}$ and
$x : X \longrightarrow F$ be any morphism
in $\mathrm{St}(\mathcal{C},\tau)$. We consider
the natural morphisms
$$d : X \longrightarrow X\times^{h}X \qquad
d_{F} : X \longrightarrow X\times^{h}_{F}X.$$
By induction on $n$, we see that the stacks
$X\times^{h}X$ and $X\times^{h}_{F}X$ both have cotangent complexes
at the point $d$ and $d_{F}$, denoted respectively
by $\mathbb{L}$ and $\mathbb{L}'$. There is moreover a natural
morphism in $\mathrm{Ho}(Sp(A-Mod))$
$$f : \mathbb{L'} \longrightarrow \mathbb{L},$$
induced by
$$\xymatrix{X \ar[r] & X\times_{F}^{h} X \ar[r] & X\times^{h}X.}$$
We set $\mathbb{L}''$ as the homotopy cofiber of $f$
in $Sp(A-Mod)$. By construction, for
any $A$-module $M \in A-Mod_{0}$, the simplicial set
$Map_{Sp(A-Mod)}(\mathbb{L}'',M)$ is naturally equivalent to the
homotopy fiber of
$$\mathbb{D}er_{X}(X,M)\times^{h}_{\mathbb{D}er_{F}(X,M)}
\mathbb{D}er_{X}(X,M)\longrightarrow\mathbb{D}er_{X}(X,M)\times \mathbb{D}er_{X}(X,M),$$
and thus is naturally equivalent to
$\Omega_{*}\mathbb{D}er_{F}(X,M)$.
By assumption we have natural isomorphisms in $\mathrm{Ho}(SSet)$
$$Map_{Sp(A-Mod)}(\Omega(\mathbb{L}''),M)\simeq
Map_{Sp(A-Mod)}(\mathbb{L}'',S(M))$$
$$\simeq \Omega_{*}\mathbb{D}er_{F}(X,S(M))
\simeq \mathbb{D}er_{F}(X,\Omega_{*}S(M))\simeq \mathbb{D}er_{F}(X,M).$$
This implies that
$\Omega_{*}\mathbb{L}''\in \mathrm{Ho}(Sp(A-Mod))$ is a cotangent complex
of $F$ at the point $x$. By induction on $n$ and by
construction we also see that
this cotangent complex is $(-n)$-connective.

Now, let
$$\xymatrix{
Y=\mathbb{R}\underline{Spec}\, B\ar[rr]^-{u} \ar[rd]_-{y} & &
X=\mathbb{R}\underline{Spec}\, A
\ar[ld]^-{x} \\
 & F, & }$$
be a morphism in $\mathrm{Ho}(Aff_{\mathcal{C}}^{\sim,\tau}/F)$, with
$A$ and $B$ in $\mathcal{A}$.
We consider the commutative diagram with
homotopy cartesian squares
$$\xymatrix{
& X \ar[r] & X\times^{h}_{F}X \ar[r] & X\times^{h}X \\
Y \ar[r] \ar[ru] & Y\times^{h}_{X}Y \ar[u] \ar[r] &
Y\times^{h}_{F}Y \ar[u] \ar[r] & Y\times^{h}Y. \ar[u]}$$
By the above explicit construction and an induction on $n$, the fact that the natural morphism
$$u^{*} : \mathbb{L}_{F,x}\otimes^{\mathbb{L}}_{A}B \longrightarrow
\mathbb{L}_{F,y}$$
is an isomorphism simply follows from the
next lemma.

\begin{lem}
Let
$$\xymatrix{
X=\mathbb{R}\underline{Spec}\, A \ar[r]^-{x} & F \ar[r] & G \\
Y=\mathbb{R}\underline{Spec}\, B \ar[u]^-{u} \ar[r]_-{y} & F' \ar[u]\ar[r] &
\ar[u] G' }$$
be a commutative diagram with the right hand square being homotopy cartesian 
in $Aff_{\mathcal{C}}^{\sim,\tau}$.
We assume that $A$ and $B$ are in $\mathcal{A}$ and that
$F$ and $G$ have global cotangent complexes.
Then the natural square
$$\xymatrix{
\mathbb{L}_{G',y} \ar[r] & \mathbb{L}_{F',y} \\
\mathbb{L}_{G,x}\otimes^{\mathbb{L}}_{A}B \ar[u]\ar[r] & \mathbb{L}_{F,x}\otimes^{\mathbb{L}}_{A}B \ar[u]}$$
is homotopy cartesian in $Sp(B-Mod)$.
\end{lem}

\begin{proof} This is immediate from the
definition and the homotopy cartesian square
$$\xymatrix{
\mathbb{D}er_{F}(Y,M)\simeq  \mathbb{D}er_{F}(X,M)\ar[r] &
\mathbb{D}er_{G}(Y,M) \simeq  \mathbb{D}er_{G}(X,M)\\
\mathbb{D}er_{F'}(Y,M) \ar[r]\ar[u] & \mathbb{D}er_{G'}(Y,M) \ar[u]}$$
for any $A$-module $M$. \end{proof}

This finishes the proof of Prop. \ref{p18}. \end{proof}

In fact, the proof of Proposition
\ref{p18} also proves the following

\begin{prop}\label{cp18}
Let $F$ be a stack such that the diagonal
$F \longrightarrow F\times^{h}F$ is $(n-1)$-representable.
We suppose that for an
 $A \in \mathcal{A}$, any
point $x : X:=\mathbb{R}\underline{Spec}\, A \longrightarrow F$,
and any $A$-module $M\in A-Mod_{0}$ the natural morphism
$$\mathbb{D}er_{F}(X,M) \simeq \mathbb{D}er_{F}(X,\Omega S(M))\longrightarrow
\Omega \mathbb{D}er_{F}(X,S(M))$$
is an isomorphism in $\mathrm{Ho}(SSet)$.
Then $F$ has a cotangent complex, which is furthermore
$(-n)$-connective.
\end{prop}

We finish this section by the notion of relative cotangent complex
and its relation with the absolute notion.
Let $f : F \longrightarrow G$ be a morphism of stacks,
$A\in \mathcal{A}$, and
$X:=\mathbb{R}\underline{Spec}\, A \longrightarrow F$ be a morphism.
We define an object $\mathbb{D}er_{F/G}(X,-) \in \mathrm{Ho}((A-Mod_{0}^{op})^{\wedge})$,
to be the standard homotopy fiber of the morphism
of the natural morphism
$$df : \mathbb{D}er_{F}(X,-) \longrightarrow \mathbb{D}er_{G}(X,-).$$
In terms of functors the object $\mathbb{D}er_{F/G}(X,-)$ sends
an $A$-module $M\in A-Mod_{0}$ to the simplicial set
$$\mathbb{D}er_{F/G}(X,M)=
Map_{X/Aff_{\mathcal{C}}^{\sim,\tau}/G}(X[M],F).$$

\begin{df}\label{d17bis}
Let $f : F \longrightarrow G$ be a morphism of stacks.
\begin{enumerate}
\item Let $A \in \mathcal{A}$,
and $x : X:=\mathbb{R}\underline{Spec}\, A \longrightarrow F$
be an $A$-point. We say that \emph{$f$ has a
(relative) cotangent complex at $x$ relative
to $\mathcal{A}$}
(or simply \emph{$f$ has a
(relative) cotangent complex at $x$}\index{cotangent complex! relative of stacks}
when $\mathcal{A}$ is unambiguous)
if there exists an integer $n$, and
an $(-n)$-connective stable $A$-module $\mathbb{L}_{F/G,x} \in \mathrm{Ho}(Sp(A-Mod))$, and
an isomorphism in
$\mathrm{Ho}(A-Mod_{0}^{op})^{\wedge}$
$$\mathbb{D}er_{F/G}(X,-)\simeq
\mathbb{R}\underline{h}_{s}^{\mathbb{L}_{F/G,x}}.$$

\item If $f$ has a cotangent complex at $x$, the stable $A$-module
$\mathbb{L}_{F,x}$ is then called the
\emph{(relative) cotangent complex of $f$ at $x$}\index{cotangent complex! relative of stacks at a point}.
\item
If $f$ has a cotangent complex at $x$,
the \emph{(relative) tangent complex of $f$ at $x$}\index{tangent complex! relative of stacks at a point} is then the
stable $A$-module
$$\mathbb{T}_{F/G,x}:=\mathbb{R}\underline{Hom}_{A}^{Sp}(
\mathbb{L}_{F/G,x},A)\in \mathrm{Ho}(Sp(A-Mod)).$$
\end{enumerate}
\end{df}

Let now be a morphism of stacks $f : F \longrightarrow G$, and a commutative
diagram in $Aff_{\mathcal{C}}^{\sim,\tau}$
$$\xymatrix{
Y\ar[rr]^-{u} \ar[rd]_-{y} & & X \ar[ld]^-{y} \\
 & F, & }$$
with $X=\mathbb{R}\underline{Spec}\, A$ and
$Y=\mathbb{R}\underline{Spec}\, B$ belonging to $\mathcal{A}$.
We have a natural morphism in $\mathrm{Ho}((A-Mod_{0}^{op})^{\wedge})$
$$\mathbb{D}er_{F/G}(Y,) \longrightarrow \mathbb{D}er_{F/G}(X,-).$$
If the morphism $f : F\longrightarrow G$ has cotangent complexes at both points
$x$ and $y$, Prop. \ref{pstabmod} induces a well defined morphism
in $\mathrm{Ho}(Sp(B-Mod))$
$$u^{*} : \mathbb{L}_{F/G,x}\otimes^{\mathbb{L}}_{A}B \longrightarrow
\mathbb{L}_{F/G,y}.$$
Of course, we have $(u\circ v)^{*}=v^{*}\circ u^{*}$
when this formula makes sense.

\begin{df}\label{d18bis}
A morphism of stacks $f : F \longrightarrow G$ \emph{has a relative 
cotangent complex relative to $(\mathcal{C},\mathcal{C}_{0},\mathcal{A})$} (or simply
\emph{has a cotangent complex}\index{cotangent complex! (global) relative of stacks} when
the HA context is clear) if the following two
conditions are satisfied.
\begin{enumerate}
\item For any  $A\in \mathcal{A}$, and any point
$x : \mathbb{R}\underline{Spec}\, A \longrightarrow F$, the morphism
$f$ has a cotangent complex $\mathbb{L}_{F/G,x}$ at $x$.
\item For any morphism $u : A \longrightarrow B$ in $\mathcal{A}$, and
any morphism in $\mathrm{Ho}(Aff_{\mathcal{C}}^{\sim,\tau}/F)$
$$\xymatrix{
Y:=\mathbb{R}\underline{Spec}\, B\ar[rr]^-{u} \ar[rd]_-{y} & &
X:=\mathbb{R}\underline{Spec}\, A \ar[ld]^-{x} \\
 & F, & }$$
the induced morphism
$$u^{*} : \mathbb{L}_{F/G,x}\otimes^{\mathbb{L}}_{A}B \longrightarrow
\mathbb{L}_{F/G,y}$$
is an isomorphism in $\mathrm{Ho}(Sp(B-Mod))$.
\end{enumerate}
\end{df}

The important remark is the following, relating absolute and relative
notions of cotangent complexes.

\begin{lem}\label{lp16-}
Let $f : F \longrightarrow G$ be a morphism
of stacks.
\begin{enumerate}
\item If both stacks $F$ and $G$ have cotangent complexes then
the morphism $f$ has a cotangent complex. Furthermore, for
any  $A\in \mathcal{A}$, and any morphism of stacks
$X=\mathbb{R}\underline{Spec}\, A \longrightarrow F$, there is
a natural homotopy cofiber sequence of stable $A$-modules
$$\xymatrix{
\mathbb{L}_{G,x}\ar[r] & \mathbb{L}_{F,x} \ar[r] & \mathbb{L}_{F/G,x}.}$$

\item If the morphism $f$ has a cotangent complex then for any
stack $H$ and any morphism $H \longrightarrow G$, the
morphism $F\times^{h}_{G}H \longrightarrow H$ has a
relative cotangent
complex and furthermore we have
$$\mathbb{L}_{F/G,x}\simeq \mathbb{L}_{F\times^{h}_{G}H /H,x}$$
for any $A\in \mathcal{A}$, and any morphism of stacks
$X=\mathbb{R}\underline{Spec}\, A \longrightarrow F\times^{h}_{G}H$.

\item If for any $A\in \mathcal{A}$ and any morphism of
stacks $x : X:=\mathbb{R}\underline{Spec}\, A \longrightarrow F$,
the morphism $F\times^{h}_{G}X \longrightarrow X$
has a relative cotangent complex, then the morphism
$f$ has a relative cotangent complex. Furthermore, we have
$$\mathbb{L}_{F/G,x}\simeq \mathbb{L}_{F\times^{h}_{G}X /X,x}.$$

\item If for any $A\in \mathcal{A}$ and any morphism of
stacks $x : X:=\mathbb{R}\underline{Spec}\, A \longrightarrow F$,
the stack $F\times^{h}_{G}X$
has a cotangent complex, then the morphism
$f$ has a relative cotangent complex. Furthermore we have
a natural homotopy cofiber sequence
$$\mathbb{L}_{A} \longrightarrow \mathbb{L}_{F\times^{h}_{G}X,x} \longrightarrow \mathbb{L}_{F/G,x}.$$

\end{enumerate}
\end{lem}

\begin{proof} $(1)$ and $(2)$ follow easily from
the definition. Point $(3)$ follows from $(2)$. Finally,
point $(4)$ follows from $(3)$, $(1)$
and Prop. \ref{p16-}. \end{proof}

\section{Obstruction theory}

Recall from \ref{dext}
that for any commutative monoid $A$, any $A$-module $M$, and
any derivation $d : A \longrightarrow A\oplus M$, we can form
the square zero extension of $A$ by $\Omega M$, denoted by
$A\oplus_{d} \Omega M$, as the homotopy cartesian
square in $Comm(\mathcal{C})$
$$\xymatrix{
A\oplus_{d}\Omega M \ar[r]^-{p} \ar[d] & A \ar[d]^{d} \\
A\ar[r]_-{s} & A\oplus M}$$
where $s : A \longrightarrow A\oplus M$ is the
trivial derivation. In the sequel, the morphism
$p : A\oplus_{d}\Omega M \longrightarrow A$ will be called the
natural projection.

\begin{df}\label{d24}
\begin{enumerate}
\item
A stack $F$ is \emph{infinitesimally cartesian relative to
the HA context $(\mathcal{C},\mathcal{C}_{0},\mathcal{A})$} (or
simply \emph{inf-cartesian}\index{inf-cartesian}\index{infinitesimally cartesian} when
the HA context is unambiguous) if for any
 commutative monoid
$A\in \mathcal{A}$, any $M\in A-Mod_{1}$, and any derivation
$d \in \pi_{0}(\mathbb{D}er(A,M))$, corresponding
to a morphism $d : A \longrightarrow A\oplus M$ in $\mathrm{Ho}(Comm(\mathcal{C})/A)$,
the square
$$\xymatrix{
\mathbb{R}F(A\oplus_{d}\Omega M) \ar[r] \ar[d] & \mathbb{R}F(A)\ar[d]^-{d} \\
\mathbb{R}F(A) \ar[r]_-{s} & \mathbb{R}F(A\oplus M)}$$
is homotopy cartesian.

\item A stack $F$ \emph{has an obstruction theory}\index{obstruction theory} (relative
to $(\mathcal{C},\mathcal{C}_{0},\mathcal{A})$)
if it has a (global) cotangent complex and if it is infinitesimally cartesian
(relative to $(\mathcal{C},\mathcal{C}_{0},\mathcal{A})$).

\end{enumerate}
\end{df}

One also has a relative version.

\begin{df}\label{d24bis}
\begin{enumerate}
\item
A morphism of stacks $F \longrightarrow G$ is \emph{infinitesimally cartesian relative to $(\mathcal{C},\mathcal{C}_{0},\mathcal{A})$} (or
simply \emph{inf-cartesian} if the HA context is clear)\index{inf-cartesian!relative}\index{infinitesimally cartesian!relative}
if for any commutative monoid
$A\in \mathcal{A}$, any $A$-module $M\in A-Mod_{1}$, and any derivation
$d \in \pi_{0}(\mathbb{D}er(A,M))$, corresponding
to a morphism $d : A \longrightarrow A\oplus M$ in $\mathrm{Ho}(Comm(\mathcal{C})/A)$,
the square
$$\xymatrix{
\mathbb{R}F(A\oplus_{d}\Omega M) \ar[r] \ar[d] & \mathbb{R}G(A\oplus_{d}\Omega M) \ar[d] \\
\mathbb{R}F(A)\times^{h}_{\mathbb{R}F(A\oplus M)}\mathbb{R}F(A) \ar[r] & 
\mathbb{R}G(A)\times^{h}_{\mathbb{R}G(A\oplus M)}\mathbb{R}G(A)}$$
is homotopy cartesian.

\item A morphism of stacks $f : F \longrightarrow G$ \emph{has an obstruction theory relative
to $(\mathcal{C},\mathcal{C}_{0},\mathcal{A})$} (or simply
\emph{has an obstruction theory}\index{obstruction theory!relative} if the HA context is clear) if
it has a (global) cotangent complex and if it is infinitesimally cartesian
relative $(\mathcal{C},\mathcal{C}_{0},\mathcal{A})$.

\end{enumerate}
\end{df}
As our HA context $(\mathcal{C},\mathcal{C}_{0},\mathcal{A})$ is fixed
once for all we will from now avoid
to mention the expression \emph{relative to $(\mathcal{C},\mathcal{C}_{0},\mathcal{A})$}
when referring to the property of having an obstruction theory.
The more precise terminology will only be used when two different
HA contexts are involved (this will only happen
in \S \ref{IIunb}). \\

We have the following
generalization of
lemma \ref{lp16-}.

\begin{lem}\label{lp16-bis}
Let $f : F \longrightarrow G$ be a morphism
of stacks.
\begin{enumerate}
\item If both stacks $F$ and $G$ have
an obstruction theory then
the morphism $f$ has an obstruction theory.

\item If the morphism $f$ has an obstruction theory then for any
stack $H$ and any morphism $H \longrightarrow G$, the
morphism $F\times^{h}_{G}H \longrightarrow H$ has a
relative obstruction theory.

\item If for any $B\in Comm(\mathcal{C})$ and any morphism of
stacks $y : Y:=\mathbb{R}\underline{Spec}\, B \longrightarrow G$,
the stack $F\times^{h}_{G}B$
has an obstruction theory , then the morphism
$f$ has a relative obstruction theory.

\end{enumerate}
\end{lem}

\begin{proof} The existence of the cotangent complexes
is done in
Lem. \ref{lp16-}, and it only remains to deal with
the inf-cartesian property. The points $(1)$ and $(2)$
are clear by definition.

$(3)$ Let $A\in \mathcal{A}$, $M\in A-Mod_{1}$ and
$d \in \pi_{0}(\mathbb{D}er(A,M))$. We need to show that
the square
$$\xymatrix{
\mathbb{R}F(A\oplus_{d}\Omega M) \ar[r] \ar[d] & \mathbb{R}G(A\oplus_{d}\Omega M) \ar[d] \\
\mathbb{R}F(A)\times^{h}_{\mathbb{R}F(A\oplus M)}\mathbb{R}F(A) \ar[r] &
\mathbb{R}G(A)\times^{h}_{\mathbb{R}G(A\oplus M)}\mathbb{R}G(A)}$$
is homotopy cartesian. Let $z$ be a point in $
\mathbb{R}G(A\oplus_{d}\Omega M)$.
We need to prove that the morphism induced on the homotopy fibers
of the two horizontal morphisms taken at $z$ is an equivalence.
But this easily follows from the fact that
the pullback of $f$ by the morphism corresponding to $z$,
$$F\times^{h}_{G}X_{d}[\Omega M] \longrightarrow
X_{d}[\Omega M],$$
has an obstruction theory. 
\end{proof}

\begin{prop}\label{p22-}
\begin{enumerate}
\item
Any representable stack has an obstruction theory.
\item
Any representable morphism has an obstruction theory.
\end{enumerate}
\end{prop}

\begin{proof} $(1)$ By Prop. \ref{p16-} we already know that
representable stacks have cotangent complexes.
Using the Yoneda lemma, it is obvious to check that any representable
stack is inf-cartesian. Indeed,
for $F=\mathbb{R}\underline{Spec}\, B$ we 
have
$$\mathbb{R}F(A\oplus_{d}\Omega M)\simeq
Map_{Comm(\mathcal{C})}(B,A\oplus_{d}\Omega M)\simeq$$
$$Map_{Comm(\mathcal{C})}(B,A)\times^{h}_{Map_{Comm(\mathcal{C})}(B,A\oplus M)}
Map_{Comm(\mathcal{C})}(B,A)\simeq 
\mathbb{R}F(A)\times^{h}_{\mathbb{R}F(A\oplus M)}\mathbb{R}F(A).$$
$(2)$ Follows from $(1)$ and Lem. \ref{lp16-bis} $(3)$.
\end{proof}

In general, the expression \emph{has an obstruction theory
relative to $(\mathcal{C},\mathcal{C}_{0},\mathcal{A})$} is justified by the
following proposition.

\begin{prop}\label{p22}
Let $F$ be a stack which has an obstruction theory.
Let $A\in \mathcal{A}$, $M\in A-Mod_{1}$
and let  $d \in \pi_{0}(\mathbb{D}er(A,M))$
be a derivation with $A\oplus_{d}\Omega M$ the corresponding square zero extension.
Let us denote by
$$X:=\mathbb{R}\underline{Spec}\, A
\longrightarrow X_{d}[\Omega M]:=\mathbb{R}\underline{Spec}\, (A\oplus_{d}\Omega M)$$
the morphism of representable stacks corresponding to the natural
projection $A\oplus_{d}\Omega M \longrightarrow A$. Finally, let
$x : X\longrightarrow F$ be an $A$-point of $F$.
\begin{enumerate}
\item
There exists a natural obstruction
$$\alpha(x)\in \pi_{0}(Map_{Sp(A-Mod)}(\mathbb{L}_{F,x},M))=[\mathbb{L}_{F,x},M]_{Sp(A-Mod)}$$
vanishing if and only $x$ extends to a morphism $x'$ in $\mathrm{Ho}(X/Aff_{\mathcal{C}}^{\sim,\tau})$
$$\xymatrix{
X \ar[rr] \ar[dr]_-{x} & & X_{d}[\Omega M] \ar[dl]^-{x'} \\
 & F. & }$$
\item Let us suppose that $\alpha(x)=0$. Then, the simplicial set of lifts of $x$
$$\mathbb{R}\underline{Hom}_{X/Aff_{\mathcal{C}}^{\sim,\tau}}(X_{d}[\Omega M],F)$$
is (non canonically) isomorphic in $\mathrm{Ho}(SSet)$ to the simplicial set
$$Map_{Sp(A-Mod)}(\mathbb{L}_{F,x},\Omega M)\simeq \Omega Map_{Sp(A-Mod)}(\mathbb{L}_{F,x},M).$$
More precisely, it is a simplicial torsor over the simplicial group
$$\Omega Map_{Sp(A-Mod)}(\mathbb{L}_{F,x},M).$$

\end{enumerate}
\end{prop}

\begin{proof} First of all, the space of lifts $x'$ is by definition
$\mathbb{R}\underline{Hom}_{X/Aff_{\mathcal{C}}^{\sim,\tau}}(X_{d}[\Omega M],F)$, which is
naturally equivalent to
the homotopy fiber at $x$ of the natural morphism
$$\mathbb{R}_{\tau}\underline{Hom}(X_{d}[\Omega M],F) \longrightarrow
\mathbb{R}_{\tau}\underline{Hom}(X,F).$$
Using that $F$ is inf-cartesian, we see that there exists a
homotopy cartesian square
$$\xymatrix{
\mathbb{R}_{\tau}\underline{Hom}(X_{d}[\Omega M],F) \ar[r] \ar[d] & \mathbb{R}_{\tau}\underline{Hom}(X,F) \ar[d]^-{d} \\
\mathbb{R}_{\tau}\underline{Hom}(X,F) \ar[r]_-{s} & \mathbb{R}_{\tau}\underline{Hom}(X[M],F).}$$
Therefore, the simplicial set
$\mathbb{R}\underline{Hom}_{X/Aff_{\mathcal{C}}^{\sim,\tau}}(X_{d}[\Omega M],F)$ fits into
a homotopy cartesian square
$$\xymatrix{\mathbb{R}\underline{Hom}_{X/Aff_{\mathcal{C}}^{\sim,\tau}}(X_{d}[\Omega M],F) \ar[r] \ar[d] & \bullet \ar[d]^-{d}  \\
\bullet \ar[r]_-{0} & \mathbb{R}\underline{Hom}_{X/Aff_{\mathcal{C}}^{\sim,\tau}}(X[M],F).}$$
As $F$ has a cotangent complex we have
$$\mathbb{R}\underline{Hom}_{X/Aff_{\mathcal{C}}^{\sim,\tau}}(X[M],F)\simeq Map_{Sp(A-Mod)}(\mathbb{L}_{F,x},M),$$ and
we see that the image of the right vertical arrow
in the last diagram provides the element
$\alpha(x)\in \pi_{0}(Map_{Sp(A-Mod)}(\mathbb{L}_{F,x},M))$, which clearly vanishes
if and only if $\mathbb{R}\underline{Hom}_{X/Aff_{\mathcal{C}}^{\sim,\tau}}(X_{d}[\Omega M],F)$
is non-empty. Furthermore, this last homotopy cartesian diagram also shows that
if $\alpha(x)=0$, then one has an isomorphism in $\mathrm{Ho}(SSet)$
$$\mathbb{R}\underline{Hom}_{X/Aff_{\mathcal{C}}^{\sim,\tau}}(X_{d}[\Omega M],F) \simeq
\Omega Map_{Sp(A-Mod)}(\mathbb{L}_{F,x},M).$$
\end{proof}

One checks immediately that the obstruction of
Prop. \ref{p22} is functorial in $F$ and $X$ in the following sense.
If $f : F \longrightarrow F'$ be a morphism
of stacks, then clearly
$$df(\alpha(x))=\alpha(f(x)) \in \pi_{0}(\mathbb{D}er_{F'}(X,M)).$$
In the same way, if $A \longrightarrow B$ is a morphism
in $\mathcal{A}$,
corresponding to a morphism of representable stacks $Y \longrightarrow X$, and
$y : Y \longrightarrow F$ be the composition, then we have
$$y^{*}(\alpha(x))=\alpha_{y} \in \pi_{0}(\mathbb{D}er_{F}(Y,M\otimes_{A}^{\mathbb{L}}B)).$$

Proposition \ref{p22} also has a relative form, whose proof is
essentially the same. We will also express it in a more
precise way.

\begin{prop}\label{p22bis}
Let $f : F \longrightarrow G$ be a morphism of stacks which has an obstruction theory.
Let $A\in \mathcal{A}$, $M\in A-Mod_{1}$, $d \in \pi_{0}(\mathbb{D}er(A,M))$
a derivation and $A\oplus_{d}\Omega M$ the corresponding square zero extension.
Let $x$ be a point in
$\mathbb{R}F(A)\times^{h}_{\mathbb{R}G(A\oplus_{d} \Omega M)}\mathbb{R}G(A)$ with
projection $y\in \mathbb{R}F(A)$, and let $L(x)$ be the
homotopy fiber, taken at $x$, of the morphism
$$\mathbb{R}F(A\oplus_{d} \Omega M) \longrightarrow
\mathbb{R}F(A)\times^{h}_{\mathbb{R}G(A\oplus_{d} \Omega M)}\mathbb{R}G(A).$$
Then there exists a natural point $\alpha(x)$ in
$Map_{A-Mod}(\mathbb{L}_{F/G,x},M)$, and a natural isomorphism in
$\mathrm{Ho}(SSet)$
$$L(z)\simeq \Omega_{\alpha(x),0}Map_{A-Mod}(\mathbb{L}_{F/G,y},M),$$
where $\Omega_{\alpha(x),0}Map_{A-Mod}(\mathbb{L}_{F/G,y},M)$
is the simplicial set of paths from $\alpha(x)$ to $0$.
\end{prop}

\begin{proof} Essentially the same as for Prop. \ref{p22}. The point
$x$ corresponds to a commutative diagram
in $\mathrm{Ho}(Aff_{\mathcal{C}}^{\sim,\tau}/G)$
$$\xymatrix{
X \ar[r] \ar[d] & F \ar[d] \\
X_{d}[\Omega M] \ar[r] & G},$$
where $X:=\mathbb{R}\underline{Spec}\, A$ and
$X_{d}[\Omega M]:=\mathbb{R}\underline{Spec}\, (A\oplus_{d}\Omega M)$.
Composing with the natural commutative diagram
$$\xymatrix{
X[M] \ar[r]^-{d} \ar[d]_-{s} & X \ar[d] \\
X \ar[r] & X_{d}[\Omega M]}$$
we get a well defined commutative diagram
in $\mathrm{Ho}(Aff_{\mathcal{C}}^{\sim,\tau}/G)$
$$\xymatrix{
X[M] \ar[r]^-{d} \ar[d]_-{s} & X \ar[r] & F \ar[d] \\
X \ar[rr] & & G,}$$
giving rise to a well defined point
$$\alpha(x)\in \mathbb{R}\underline{Hom}_{X/Aff_{\mathcal{C}}^{\sim,\tau}/G}(X[M],F)=
\mathbb{D}er_{F/G}(X,M).$$
Using that the morphism $f$ is inf-cartesian, we easily see  that
the simplicial set $\Omega_{\alpha(x),0}\mathbb{D}er_{F/G}(X,M)$ is naturally
equivalent to the space of lifts
$$L(x)=\mathbb{R}\underline{Hom}_{X/Aff_{\mathcal{C}}^{\sim,\tau}/G}(X_{d}[\Omega M],F)\simeq
\Omega_{\alpha(x),0}\mathbb{D}er_{F/G}(X,M).$$
\end{proof}

\begin{prop}\label{p23}
Let
$F$ be a stack whose diagonal $F \longrightarrow F\times^{h}F$
is $n$-representable for some $n$.
Then $F$ has an obstruction theory if and only if
it is inf-cartesian.
\end{prop}

\begin{proof} It is enough to show that
a stack $F$ that is inf-cartesian satisfies the condition
of proposition \ref{p18}. But this follows easily from
the following homotopy cartesian square
$$\xymatrix{
A\oplus M \ar[r] \ar[d] & A \ar[d]^-{s}\\
A \ar[r]_-{s} & A\oplus S(M)}$$
for any commutative monoid $A\in \mathcal{A}$ and any
$A$-module $M\in A-Mod_{1}$. \end{proof}

\section{Artin conditions}

In this section we will give
conditions on the topology $\tau$ and on the class of morphisms
\textbf{P} ensuring that any stack which is geometric for such a $\tau$ and \textbf{P}  will have an
obstruction theory. We call these conditions \emph{Artin's conditions}, 
though we warn the reader that these are not the rather famous conditions
for a functor to be representable by an algebraic space. We refer instead to 
the fact that an algebraic stack in the sense of Artin (i.e. with a smooth atlas)
has a good infinitesimal and obstruction theory. We think that this has been 
first noticed by M. Artin, since this is precisely one part of the easy direction of his representability criterion
(\cite[I, 1.6]{Ar}  or \cite[Thm. 10.10]{lm}). 

\begin{df}\label{d25}
We will say that \emph{$\tau$ and \textbf{P} satisfy Artin's conditions\index{Artin conditions}
relative to $(\mathcal{C},\mathcal{C}_{0},\mathcal{A})$}
(or simply \emph{satisfy Artin's conditions} if the
HA context is clear)
if there
exists a class \textbf{E} of morphisms in $Aff_{\mathcal{C}}$ such that the
following conditions are satisfied.
\begin{enumerate}
\item Any morphism in \textbf{P} is formally i-smooth
in the sense of Def. \ref{dismooth}.
\item Morphisms in \textbf{E} are formally \'etale, are stable by equivalences, homotopy pullbacks
and composition.
\item For any
morphism $A \longrightarrow B$ in \textbf{E} 
with $A\in \mathcal{A}$, we have $B\in \mathcal{A}$. 
\item For any epimorphism of representable stacks
$Y \longrightarrow X$, which is a \textbf{P}-morphism,
there exists an epimorphism of representable stacks
$X' \longrightarrow X$, which is in \textbf{E}, and
a commutative diagram in $\mathrm{St}(\mathcal{C},\tau)$
$$\xymatrix{
& Y \ar[d] \\
X' \ar[ru] \ar[r] & X.}$$
\item Let $A\in \mathcal{A}$, $M
\in A-Mod_{1}$,
$d\in \pi_{0}(\mathbb{D}er(A,M))$ be a derivation, and
$$X:=\mathbb{R}\underline{Spec}\, A \longrightarrow
X_{d}[\Omega M]:=\mathbb{R}\underline{Spec}, (A\oplus_{d}\Omega M)$$
be the natural morphism in $\mathrm{St}(\mathcal{C},\tau)$.
Then, a formally \'etale morphism of representable stacks
$p : U \longrightarrow X_{d}[\Omega M]$ is in \textbf{E}
if and only if $U\times^{h}_{X_{d}[\Omega M]}X \longrightarrow X$ is so.
Furthermore, if $p$ is in \textbf{E}, then $p$ is an epimorphism of
stacks if and only if $U\times^{h}_{X_{d}[\Omega M]}X \longrightarrow X$ is so.
\end{enumerate}
\end{df}

The above definition might seem a bit technical and somehow hard to follow. In order to fix his intuition, we suggest the reader
to think in terms of standard algebraic geometry with $\tau$ being the etale 
topology, \textbf{P} the class of smooth morphisms, and 
\textbf{E} the class of \'etale morphisms. This is only meant to convey some classical geometric 
intuition because this classical situation in  algebraic geometry 
does not really fit into the above definition; in fact in this case the base category is the category of 
$k$-modules and thus the suspension functor is not fully faithful. \\

In order to simplify notations we will say that a morphism
$A\longrightarrow B$ in $Comm(\mathcal{C})$ is an \emph{\textbf{E}-covering}\index{\textbf{E}-covering}
if it is in \textbf{E} and if
the corresponding morphism of stacks
$$\mathbb{R}\underline{Spec}\, B \longrightarrow
\mathbb{R}\underline{Spec}\, A$$
is an epimorphism of stacks.

\begin{thm}\label{t1}
Assume $\tau$ and
\textbf{P} satisfy Artin's conditions.
\begin{enumerate}
\item
Any $n$-representable morphism of stacks has an obstruction theory.
\item Let $f : F \longrightarrow G$ be an $n$-representable
morphism of stacks. If $f$ is in \textbf{P} then
for any $A\in \mathcal{A}$, and any morphism
$x : X:=\mathbb{R}\underline{Spec}\, A \longrightarrow F$
there exists an  \textit{E}-covering
$$x' : X':=\mathbb{R}\underline{Spec}\, A'  \longrightarrow X$$
such that for any $M\in A'-Mod_{1}$ the natural morphism
$$[\mathbb{L}_{X'/G,x'},M] \longrightarrow [\mathbb{L}_{F/G,x},M]_{A-Mod}$$
is zero.
\end{enumerate}
\end{thm}

\begin{proof} Before starting the proof we will need
the following general fact, that will be used
all along the proof of the theorem.

\begin{lem}\label{lcart}
Let $D$ be a pointed model category for which the suspension
functor
$$S : \mathrm{Ho}(D) \longrightarrow \mathrm{Ho}(D)$$
is fully faithful. Then, a homotopy co-cartesian square
in $D$ is also homotopy cartesian.
\end{lem}

\begin{proof} When $D$ is a stable model category this
is well known since homotopy fiber sequences are also
homotopy cofiber sequences (see \cite{ho}). The general case
is proved in the same way. When $D$ is furthermore
$\mathbb{U}$-cellular (which will be our case), one can even deduce
the result from the stable case by using the
left Quillen functor $D \longrightarrow Sp(D)$, from
$D$ to the model category of spectra in $D$ as
defined in \cite{ho2}, and using the fact that it
is homotopically fully faithful. \end{proof}

The previous lemma can be applied to $\mathcal{C}$, but also
to the model categories $B-Mod$ of modules over
some commutative monoid $B$. In particular, homotopy
cartesian square of $B$-modules which are also
homotopy co-cartesian will remain homotopy cartesian after a
derived tensor product by any $B$-module. \\

Let us now start the proof
of theorem \ref{t1}. \\

\bigskip

\noindent\textit{Some topological invariance statements} \\

We start with several results concerning topological invariance of
formally \'etale morphisms and morphisms in \textbf{E}.

\begin{lem}\label{l15}
Let $A$ be a commutative monoid, $M$ an $A$-module
and $A\oplus M \longrightarrow A$ the natural augmentation.
Then, the homotopy push out functor
$$A\otimes^{\mathbb{L}}_{A\oplus M} - : \mathrm{Ho}((A\oplus M)-Comm(\mathcal{C})) \longrightarrow
\mathrm{Ho}(A-Comm(\mathcal{C}))$$
induces an equivalence between the full sub-categories
consisting of
formally \'etale commutative $A\oplus M$-algebras and formally \'etale
commutative $A$-algebras.
\end{lem}

\begin{proof} We see that the functor is essentially surjective
as the morphism $A\oplus M \longrightarrow A$ possesses
a section $A \longrightarrow A\oplus M$.
Next, we show that the  functor $A\otimes^{\mathbb{L}}_{A\oplus M} - $ is
conservative. For this, let us consider the
commutative square
$$\xymatrix{
A\oplus M  \ar[r] \ar[d] & A \ar[d] \\
A\ar[r] & A\oplus S(M),}$$
which, as a commutative square in $\mathcal{C}$, is homotopy cocartesian and homotopy cartesian.
Therefore, this square is homotopy cocartesian in
$(A\oplus M)-Mod$, and thus lemma \ref{lcart} implies that
for any commutative $A\oplus M$-algebra $B$, the natural morphism
$$B \longrightarrow \left(A\otimes_{A\oplus M}^{\mathbb{L}}B\right)\times^{h}_{
(A\oplus S(M))\otimes_{A\oplus M}^{\mathbb{L}}B}
\left(A\otimes_{A\oplus M}^{\mathbb{L}}B\right)
\simeq \left(A\times^{h}_{A\oplus S(M)}A\right)\otimes_{A\oplus M}^{\mathbb{L}}B$$
is then an isomorphism in $\mathrm{Ho}(Comm(\mathcal{C}))$. This clearly implies that the
functor $A\otimes^{\mathbb{L}}_{A\oplus M} - $ is conservative.

Now, let $A\oplus M \longrightarrow B$ be a formally \'etale
morphism of commutative monoids. The diagonal morphism $M
\longrightarrow M\times^{h} M\simeq M\oplus M$ in $\mathrm{Ho}(A-Mod)$, induces
a well defined morphism in $\mathrm{Ho}(Comm(\mathcal{C})/A\oplus M)$
$$\xymatrix{
 & (A\oplus M)\oplus M \ar[d] \\
A\oplus M \ar[ru] \ar[r]_-{Id} & A\oplus M,}$$
and therefore a natural element in $\pi_{0}\mathbb{D}er(A\oplus M,M)$.
Composing with $M\longrightarrow M\otimes_{A}^{\mathbb{L}}B$, we get
a well defined element in $\pi_{0}\mathbb{D}er(A\oplus M,M\otimes_{A}^{\mathbb{L}}B)$.
Using that $A\oplus M \longrightarrow B$ is formally \'etale, this element extends uniquely to an
element in $\pi_{0}\mathbb{D}er(B,M\otimes_{A}^{\mathbb{L}}B)$. This last derivation
gives rise to a well defined morphism in $\mathrm{Ho}((A\oplus M)-Comm(\mathcal{C}))$
$$u : B \longrightarrow \left(A\otimes_{A\oplus M}^{\mathbb{L}}B\right)\oplus M\otimes_{A}^{\mathbb{L}}B.$$
Furthermore, by construction, this morphism is sent to the identity of
$A\otimes_{A\oplus M}^{\mathbb{L}}B$ by the functor
$A\otimes^{\mathbb{L}}_{A\oplus M} -$, and
as we have seen this implies that $u$ is an isomorphism in $\mathrm{Ho}((A\oplus M)-Comm(\mathcal{C}))$.

We now finish the proof of the lemma by showing that the functor
$A\otimes^{\mathbb{L}}_{A\oplus M} -$ is fully faithful. For this, let
$A\oplus M \longrightarrow B$ and $A\oplus M \longrightarrow B'$ be two formally \'etale morphisms
of commutative monoids. As we have seen, $B'$ can be written as
$$B'\simeq (A\oplus M)\otimes_{A}^{\mathbb{L}}A'\simeq A'\oplus M'$$
where $A \longrightarrow A'$ is formally \'etale, and
$M':=M\otimes_{A}^{\mathbb{L}}A'$. We consider the natural morphism
$$Map_{(A\oplus M)-Comm(\mathcal{C})}(B,B') \simeq Map_{(A\oplus M)-Comm(\mathcal{C})}(B,A'\oplus M') \longrightarrow$$
$$
Map_{A-Comm(\mathcal{C})}(A\otimes_{A\oplus M}^{\mathbb{L}}B,A\otimes_{A\oplus M}^{\mathbb{L}}B')\simeq
Map_{A-Comm(\mathcal{C})}(A\otimes_{A\oplus M}^{\mathbb{L}}B,A').$$
The homotopy fiber of this morphism at a point $B \rightarrow A'$ is identified with
$\mathbb{D}er_{A\oplus M}(B,M')$, which is contractible as
$B$ is formally \'etale over $A\oplus M$. This shows that
the morphism
$$Map_{(A\oplus M)-Comm(\mathcal{C})}(B,B') \longrightarrow
Map_{A-Comm(\mathcal{C})}(A\otimes_{A\oplus M}^{\mathbb{L}}B,A\otimes_{A\oplus M}^{\mathbb{L}}B')$$
has contractible homotopy fibers and therefore is an
isomorphism in $\mathrm{Ho}(SSet)$, and finishes the
proof of the lemma. \end{proof}

\begin{lem}\label{l18'}
Let $A$ be a commutative monoid and $M$ an $A$-module.
Then
the homotopy push out functor
$$A\otimes^{\mathbb{L}}_{A\oplus M} - : \mathrm{Ho}((A\oplus M)-Comm(\mathcal{C})) \longrightarrow
\mathrm{Ho}(A-Comm(\mathcal{C}))$$
induces an equivalence between the full sub-categories
consisting of
\textbf{E}-coverings of $A\oplus M$ and of
\textbf{E}-coverings of $A$.
\end{lem}

\begin{proof} Using lemma \ref{l15} it is enough to show that
a formally \'etale morphism $f : A\oplus M \longrightarrow B$
is in \textbf{E} (respectively an \textbf{E}-covering) if
and only if $A \longrightarrow A':=A\otimes_{A\oplus M}^{\mathbb{L}}B$
is
in \textbf{E} (respectively an \textbf{E}-covering). But, as $f$ is formally \'etale,
lemma \ref{l15} implies that it can
written as
$$A\oplus M \longrightarrow A'\oplus M'\simeq (A\oplus M)\otimes_{A}^{\mathbb{L}}A',$$
with $M':=M\otimes_{A}^{\mathbb{L}}A'$.
Therefore the lemma simply follows from the stability of
epimorphisms and morphisms in \textbf{E} by homotopy pullbacks.
\end{proof}

\begin{lem}\label{l16}
Let $A\in Comm(\mathcal{C})$, $M\in A-Mod_{1}$,
and $d \in \pi_{0}(\mathbb{D}er(A,M))$ be a derivation.
Let $B:=A\oplus_{d}\Omega M \longrightarrow A$ be the natural
augmentation, and let us consider the base change functor
$$A\otimes_{B}^{\mathbb{L}}- :
\mathrm{Ho}(B-Comm(\mathcal{C})) \longrightarrow \mathrm{Ho}(A-Comm(\mathcal{C})).$$
Then, $A\otimes_{B}^{\mathbb{L}} -$ restricted
to the full subcategory consisting of
formally \'etale commutative $B$-algebras
is fully faithful.
\end{lem}

\begin{proof} Let
$$\xymatrix{
B \ar[r] \ar[d] & A \ar[d]^-{d} \\
A \ar[r]_-{s} & A\oplus M}$$
be the standard homotopy cartesian square
of commutative monoids, which
is also homotopy co-cartesian in $\mathcal{C}$
as $M\in A-Mod_{1}$. We represent it
as a fibered square in $Comm(\mathcal{C})$
$$\xymatrix{
B \ar[r] \ar[d] & A \ar[d]^-{d} \\
A' \ar[r]_-{s'} & A\oplus M,}$$
where $s' : A' \longrightarrow A\oplus M$
is fibrant replacement of the trivial
section $s : A \longrightarrow A\oplus M$.

We define
a model category $D$ whose objects are
$5$-plets $(B_{1},B_{2},B_{3},a,b)$, where
$B_{1}\in A-Comm(\mathcal{C})$, $B_{2}\in A'-Comm(\mathcal{C})$,
$B_{3}\in (A\oplus M)-Comm(\mathcal{C})$, and
$a$ and $b$ are morphisms of commutative $(A\oplus M)$-algebras
$$\xymatrix{
(A\oplus M)\otimes_{A}B_{1} \ar[r]^-{a} & B_{3} & \ar[l]_-{b}
(A\oplus M)\otimes_{A'}B_{2}}$$
(where the co-base change on the left is taken
with respect of the morphism $s' : A' \rightarrow A\oplus M$ and
the one on the right with respect to
$d : A \rightarrow A\oplus M$). For an object $(B_{1},B_{2},B_{3},a,b)$
in $D$, the morphisms $a$ and $b$ can also be understood
as  $B_{1} \longrightarrow B_{3}$ in $A-Comm(\mathcal{C})$ and
$B_{2} \longrightarrow B_{3}$ in $A'-Comm(\mathcal{C})$.
The morphisms
$$(B_{1},B_{2},B_{3},a,b) \longrightarrow (B_{1}',B_{2}',B_{3}',a',b')$$
in $D$ are defined
in the obvious way, as families of morphisms
$\{B_{i} \rightarrow B_{i}'\}$
commuting with the $a$'s and $b$'s. A morphism in $D$
is defined to be an equivalence or a cofibration if each
morphism $B_{i} \rightarrow B_{i}'$ is so. A morphism
in $D$ is defined to be a fibration if each
morphism $B_{i} \rightarrow B_{i}'$ is a fibration in $\mathcal{C}$, and
if the natural morphisms
$$B_{1} \longrightarrow B_{1}'\times_{B_{3}'}B_{3}
\qquad B_{2} \longrightarrow B_{2}'\times_{B_{3}'}B_{3}$$
are fibrations in $\mathcal{C}$.
This defines a model category structure on $D$
which is a Reedy type model structure.
An important fact concerning $D$ is the
description of its mapping spaces as
the following homotopy cartesian square
$$\xymatrix{
Map_{D}((\underline{B},a,b),(\underline{B'},a',b'))
\ar[r] \ar[d] & Map_{A-Comm(\mathcal{C})}(B_{1},B_{1}')\times
Map_{A'-Comm(\mathcal{C})}(B_{2},B_{2}')
\ar[d] \\
Map_{(A\oplus M)-Comm(\mathcal{C})}(B_{3},B_{3}') \ar[r] &
Map_{A-Comm(\mathcal{C})}
(B_{1},B_{3}')\times Map_{A'-Comm(\mathcal{C})}
(B_{2},B_{3}')}$$ where we have denoted $\underline{B}:=(B_{1},B_{2},B_{3})$ and $\underline{B'}:=(B_{1}',B_{2}',B_{3}')$.
There exists a natural functor
$$F : B-Comm(\mathcal{C}) \longrightarrow  D$$
sending a commutative $B$-algebra $B'$ to
the object
$$F(B'):=\left( A\otimes_{B}B',A'\otimes_{B}B',(A\oplus M)\otimes_{B}B',
a,b \right)$$
where
$$a : (A\oplus M)\otimes_{A}(A\otimes_{B}B') \simeq (A\oplus M)\otimes_{B}B' \qquad
b : (A\oplus M)\otimes_{A'}(A'\otimes_{B}B') \simeq (A\oplus M)\otimes_{B}B'$$
are the two natural isomorphisms in $(A\oplus M)-Comm(\mathcal{C})$. The
functor $F$ has a right adjoint $G$, sending
an object $(B_{1},B_{2},B_{3},a,b)$ to the pullback in
$B-Comm(\mathcal{C})$
$$\xymatrix{
G(B_{1},B_{2},B_{3},a,b) \ar[r] \ar[d] & B_{1} \ar[d]^-{a} \\
B_{2} \ar[r]_-{b} & B_{3}.}$$
Clearly the adjunction $(F,G)$ is a Quillen adjunction. Furthermore,
lemma \ref{lcart} implies that
for any commutative $B$-algebra $B \rightarrow B'$ the adjunction morphism
$$B' \longrightarrow \mathbb{R}G\mathbb{L}F(B')=
A\otimes_{B}^{\mathbb{L}}B'\times^{h}_{
(A\oplus M)\otimes_{B}^{\mathbb{L}}B'}
A'\otimes_{B}^{\mathbb{L}}B'\simeq (A\times^{h}_{A \oplus M}A')\otimes^{\mathbb{L}}_{B}B'$$
is an isomorphism in $\mathrm{Ho}(B-Comm(\mathcal{C}))$. This implies
in particular that
$$\mathbb{L}F : \mathrm{Ho}(B-Comm(\mathcal{C})) \longrightarrow \mathrm{Ho}(D)$$
is fully faithful.

We now consider the functor
$$D \longrightarrow A-Comm(\mathcal{C})$$
sending $(B_{1},B_{2},B_{3},a,b)$ to $B_{1}$.
Using our lemma
\ref{l15}, and the description of the mapping spaces
in $D$, it is not hard to see that the
induced functor
$$\mathrm{Ho}(D) \longrightarrow \mathrm{Ho}(A-Comm(\mathcal{C}))$$
becomes fully faithful when restricted to the full subcategory
of $\mathrm{Ho}(D)$ consisting of objects $(B_{1},B_{2},B_{3},a,b)$
such that $A \rightarrow B_{1}$
and $A' \rightarrow B_{2}$ are formally \'etale and
the induced morphism
$$a : (A\oplus M)\otimes^{\mathbb{L}}_{A}B_{1} \longrightarrow B_{3} \qquad
b : (A\oplus M)\otimes^{\mathbb{L}}_{A'}B_{2} \longrightarrow B_{3}$$
are isomorphisms in $\mathrm{Ho}((A\oplus M)-Comm(\mathcal{C}))$. Putting all of this
together we deduce that the functor
$$A\otimes_{B}^{\mathbb{L}}- :
\mathrm{Ho}(B-Comm(\mathcal{C})) \longrightarrow \mathrm{Ho}(A-Comm(\mathcal{C}))$$
is fully faithful when restricted to the full subcategory
of formally \'etale morphisms. \end{proof}

\begin{lem}\label{l16'}
Let $A$ be a commutative monoid, $M$ an $A$-module
and $d : \mathbb{D}er(A,M)$ be a derivation.
Let $B:=A\oplus_{d}\Omega M \longrightarrow A$ be the natural
augmentation. Then, there exists a natural
homotopy cofiber sequence of $A$-modules
$$\xymatrix{
\mathbb{L}_{B}\otimes_{B}^{\mathbb{L}}A \ar[r] &
\mathbb{L}_{A} \ar[r]^-{d} & \mathbb{L}QZ(M),}$$
where
$$Q : A-Comm_{nu}(\mathcal{C}) \longrightarrow A-Mod \qquad
A-Comm_{nu}(\mathcal{C}) \longleftarrow A-Mod : Z$$
is the Quillen adjunction described during the proof
of \ref{p1}.
\end{lem}

\begin{proof} When $d$ is the trivial derivation, we know
the lemma is correct as by Prop. \ref{p2} $(4)$ we have
$$\mathbb{L}_{A\oplus M}\otimes_{A\oplus M}^{\mathbb{L}}A\simeq
\mathbb{L}_{A}\coprod \mathbb{L}QZ(M).$$

For the general situation, we use our
left Quillen functor
$$F : B-Comm(\mathcal{C}) \longrightarrow D$$
defined during the proof of lemma \ref{l16}.
The commutative monoid $A$ is naturally an
$(A\oplus M)$-algebra, and can be considered
as a natural object $(A,A,A)$ in $D$ with the obvious
transition morphisms
$$A\otimes_{A}(A\oplus M)\simeq A\oplus M \rightarrow A,$$
$$A\otimes_{A'}(A\oplus M)\rightarrow
A\otimes_{A}(A\oplus M)\simeq A\oplus M \rightarrow A.$$
For any $A$-module $N$, one can consider $A\oplus N$
as a commutative $A$-algebra, and therefore as as an object
$(A\oplus N,  A\oplus N, A\oplus N)$
in $D$ (with the obvious transition morphisms). 
We will simply denote by $A$ the object $(A,A,A)\in D$, and
by $A\oplus N$ the object $(A\oplus N,  A\oplus N, A\oplus N)\in D$.
The
left Quillen property of $F$ implies that
$$\mathbb{D}er(B,N)\simeq Map_{D/A}(F(B),A\oplus N).$$
This shows that the morphism
$$\mathbb{D}er(A,N) \longrightarrow \mathbb{D}er(B,N)$$
is equivalent to the morphism
$$\mathbb{D}er(A,N)  \longrightarrow
\mathbb{D}er(A,N)\times^{h}_{\mathbb{D}er(A\oplus M,N)}\mathbb{D}er(A,N).$$
This implies the the morphism
$$\mathbb{L}_{B}\otimes^{\mathbb{L}}_{B}A \longrightarrow
\mathbb{L}_{A}$$
is naturally equivalent to
the morphism of
$A$-modules
$$\mathbb{L}_{A}\coprod^{\mathbb{L}}_{\mathbb{L}_{A\oplus M}
\otimes^{\mathbb{L}}_{A\oplus M}A} \mathbb{L}_{A} \longrightarrow \mathbb{L}_{A}.$$
Using the already known result for the trivial
extension $A\oplus M$ we get the required natural cofiber sequence
$$\mathbb{L}_{B}\otimes^{\mathbb{L}}_{B}A \longrightarrow \mathbb{L}_{A}
\longrightarrow \mathbb{L}QZ(M).$$
\end{proof}

\begin{lem}\label{l17}
Let $A$ be a commutative monoid, $M\in A-Mod_{1}$
and $d \in \pi_{0}(\mathbb{D}er(A,M))$ be a derivation.
Let $B:=A\oplus_{d}\Omega M \longrightarrow A$ be the natural
augmentation, and let us consider the base change functor
$$A\otimes_{B}^{\mathbb{L}}- :
\mathrm{Ho}(B-Comm(\mathcal{C})) \longrightarrow \mathrm{Ho}(A-Comm(\mathcal{C})).$$
Then, $A\otimes_{B}^{\mathbb{L}} -$ induces an equivalence between
the full sub-categories consisting of
formally \'etale commutative $B$-algebras
and of formally \'etale commutative $A$-algebras.
\end{lem}

\begin{proof} By lemma \ref{l16} we already know that
the functor is fully faithful, and it only remains to show that
any formally \'etale $A$-algebra $A\rightarrow A'$
is of the form $A\otimes^{\mathbb{L}}_{B}B'$ for some
formally \'etale morphism $B \rightarrow B'$.

Let $A \longrightarrow A'$ be a formally \'etale morphism.
The derivation $d \in \pi_{0}\mathbb{D}er(A,M)$ lifts uniquely to
a derivation $d' \in \pi_{0}\mathbb{D}er(A',M')$
where $M':=M\otimes^{\mathbb{L}}_{A}A'$. We form
the corresponding square zero extension
$$\xymatrix{
B':=A'\oplus_{d}\Omega M' \ar[r] \ar[d] & A'\ar[d]^-{d'} \\
A' \ar[r] & A'\oplus M'}$$
which comes equipped with a natural morphism
$B \longrightarrow B'$ fitting in a homotopy commutative square
$$\xymatrix{
B' \ar[r] & A' \\
B \ar[r] \ar[u] & A. \ar[u]}$$
We claim that $B \longrightarrow B'$ is formally \'etale
and that the natural morphism
$A\otimes^{\mathbb{L}}_{B}B' \longrightarrow A'$ is
an isomorphism in $\mathrm{Ho}(A-Comm(\mathcal{C}))$.

There are natural cofiber sequences in $\mathcal{C}$
$$B' \longrightarrow A'\longrightarrow M' $$
$$B \longrightarrow A \longrightarrow M,$$
as well as a natural morphism of cofiber sequences
in $\mathcal{C}$
$$\xymatrix{
A\otimes_{B}^{\mathbb{L}}B'\ar[r] \ar[d] & \ar[r] A\otimes_{B}^{\mathbb{L}}A' \ar[d] &
 A\otimes_{B}^{\mathbb{L}}  M'\ar[d] \\
B\otimes_{B}^{\mathbb{L}}A'\simeq A' \ar[r] &
A\otimes_{B}^{\mathbb{L}}A' \ar[r] & M \otimes_{B}^{\mathbb{L}}A',}$$
which by our lemma \ref{lcart} is also a morphism of
fiber sequences.
The two right vertical morphisms are equivalences and
thus so is the arrow on the left. This shows that
$A\otimes^{\mathbb{L}}_{B}B' \simeq A'$.

Finally, lemma \ref{l16'} implies the existence of
a natural morphism of cofiber sequences of $A'$-modules
$$\xymatrix{
\mathbb{L}_{B}\otimes_{B}^{\mathbb{L}}A \ar[r] \ar[d] &
\mathbb{L}_{A} \ar[d] \ar[r] & \mathbb{L}QZ(M) \ar[d] \\
\mathbb{L}_{B'}\otimes_{B'}^{\mathbb{L}}A' \ar[r] &
\mathbb{L}_{A'} \ar[r] & \mathbb{L}QZ(M')\simeq  \mathbb{L}QZ(M)\otimes^{\mathbb{L}}_{A}A'.}$$
As $A\rightarrow A'$ is \'etale, this implies that
the natural morphism
$$(\mathbb{L}_{B}\otimes_{B}^{\mathbb{L}}B')\otimes^{\mathbb{L}}_{B'}A'
\simeq (\mathbb{L}_{B}\otimes_{B}^{\mathbb{L}}A)\otimes^{\mathbb{L}}_{A}A'
\longrightarrow \mathbb{L}_{B'}\otimes_{B'}^{\mathbb{L}}A'$$
is an isomorphism in $\mathrm{Ho}(A'-Mod)$. This would show that
$B \longrightarrow B'$ is formally \'etale if one knew that
the base change functor
$$A'\otimes_{B'}^{\mathbb{L}} - : \mathrm{Ho}(B'-Mod) \longrightarrow
\mathrm{Ho}(A'-Mod)$$
were conservative. However, this is the case
as lemma \ref{lcart} implies that for any
$B'$-module $N$ we have a natural isomorphism
in $\mathrm{Ho}(B'-Mod)$
$$N\simeq (A'\otimes^{\mathbb{L}}_{B'}N)\times^{h}_{
(A'\oplus M')\otimes^{\mathbb{L}}_{B'}N}(A'\otimes^{\mathbb{L}}_{B'}N).$$
\end{proof}

\begin{lem}\label{l17'}
Let $A\in \mathcal{A}$, $M\in A-Mod_{1}$
and $d \in \pi_{0}(\mathbb{D}er(A,M))$ be a derivation.
Let
$$A\oplus_{d}\Omega M \longrightarrow A$$
be the natural morphism and let us consider the homotopy push-out functor
$$(A\oplus_{d}\Omega M)\otimes_{A}^{\mathbb{L}} - :
\mathrm{Ho}((A\oplus_{d}\Omega M)-Comm(\mathcal{C})) \longrightarrow
\mathrm{Ho}(A-Comm(\mathcal{C})).$$
Then, $(A\oplus_{d}\Omega M)\otimes_{A}^{\mathbb{L}} -$ induces an equivalence between
the full sub-categories consisting of
\textbf{E}-covers of $A\oplus_{d}\Omega M$ and
of \textbf{E}-covers of $A$.
\end{lem}

\begin{proof} This is immediate from
Lem. \ref{l17} and condition $(4)$ of
Artin's conditions \ref{d25}. \end{proof}

\textit{Proof of Theorem \ref{t1}} \\

We are now ready to prove that $F$
has an obstruction theory. To simplify notations
we assume that $F$ is fibrant in $Aff^{\sim,\tau}_{\mathcal{C}}$, 
so $F(A)\simeq \mathbb{R}F(A)$ for any $A \in Comm(\mathcal{C})$.  
We then argue by induction
on the integer $n$. For $n=-1$ Theorem \ref{t1}
follows from Prop. \ref{p22-}, hypothesis
Def. \ref{d25} $(5)$ and Prop. \ref{pismooth}.
We now assume that $n\geq 0$ and that both statement of theorem
\ref{t1} are true for all $m<n$. \\

We start by proving Thm. \ref{t1} $(1)$ for rank $n$. For this, we
use Lem. \ref{lp16-bis} and Prop. \ref{p23}, which show that
we only need to prove that any $n$-geometric
stack is inf-cartesian.

\begin{lem}\label{l18}
Let $F$ be an $n$-geometric stack. Then
$F$ is inf-cartesian.
\end{lem}

\begin{proof}
Let $A\in \mathcal{A}$, $M\in A-Mod_{1}$, and
$d\in \pi_{0}(\mathbb{D}er(A,M))$.
Let $x$ be a point in
$\pi_{0}(F(A)\times^{h}_{F(A\oplus M)}F(A))$, whith
projection $x_{1}\in \pi_{0}(F(A))$ on the first factor.
We need to show that
the homotopy fiber, taken at $x$, of the morphism
$$F(A\oplus_{d}\Omega M) \longrightarrow F(A)\times^{h}_{F(A\oplus M)}F(A)$$
is contractible. For this, we  replace the
homotopy cartesian diagram
$$\xymatrix{A\oplus_{d}\Omega M \ar[r] \ar[d] & A \ar[d] \\
A \ar[r] & A\oplus M}$$
be an equivalent commutative diagram in $Comm(\mathcal{C})$
$$\xymatrix{B \ar[r] \ar[d] & B_{1} \ar[d] \\
B_{2} \ar[r] & B_{3},}$$
in such a way that each morphism is a cofibration
in $Comm(\mathcal{C})$.  The point $x$ can be represented as a
point in the standard homotopy pullback
$F(B_{1})\times^{h}_{F(B_{3})}F(B_{2})$.
We then
define a functor
$$S : B-Comm(\mathcal{C})=(Aff_{\mathcal{C}}/Spec\, B)^{op} \longrightarrow SSet_{\mathbb{V}},$$
in the following way. For any morphism of commutative monoids
$B \longrightarrow B'$, the simplicial set $S(B')$ is defined to be
the standard homotopy fiber, taken at $x$,  of the natural morphism
$$F(B') \longrightarrow
F(B_{1}\otimes_{B}B')\times^{h}_{
F(B_{3}\otimes_{B}B')}F(B_{2}\otimes_{B}B').$$
Because of our choices on the $B_{i}$'s,
it is clear that the simplicial presheaf $S$ is a stack
on the comma model site $Aff_{\mathcal{C}}/Spec\, B$. Therefore,
in order to show that $S(B)$ is contractible it is
enough to show that $S(B')$ is contractible for
some morphism $B \longrightarrow B'$ such that
$\mathbb{R}\underline{Spec}\, B' \longrightarrow
\mathbb{R}\underline{Spec}\, B$ is an epimorphism of stacks.
In particular, we are allowed to homotopy base change by some
\textbf{E}-covering of $B$. Also, using our lemma
\ref{l17} (or rather its proof), we see that for any
\textbf{E}-covering $B\rightarrow B'$, the homotopy cartesian square
$$\xymatrix{B' \ar[r] \ar[d] & B_{1}\otimes_{B}B' \ar[d] \\
B_{2}\otimes_{B}B' \ar[r] & B_{3}\otimes_{B}B',}$$
is in fact equivalent to some
$$\xymatrix{A'\oplus_{d'}\Omega M' \ar[r] \ar[d] & A' \ar[d]^-{d'} \\
A' \ar[r] & A'\oplus M',}$$
for some \textbf{E}-covering $A \longrightarrow A'$
(and with
$M'\simeq M\otimes^{\mathbb{L}}_{A}A'$, and where $d'$ is
the unique derivation $d' \in \pi_{0}\mathbb{D}er(A',M')$
extending $d$). This shows that we can
always replace $A$ by $A'$, $d$ by $d'$ and $M$ by $M'$.
In particular, for an $(n-1)$-atlas $\{U_{i} \longrightarrow F\}$ 
$F$, we can assume that
the point $x_{1}\in \pi_{0}(F(A))$, image of
the point $x$, lifts to a
point in $y_{1}\in \pi_{0}(U_{j}(A))$ for some
$j$. We will denote $U:=U_{j}$. 

\begin{sublem}
The point $x \in \pi_{0}(F(A)\times^{h}_{F(A\oplus M)}F(A))$
lifts to  point $y\in \pi_{0}(U(A)\times^{h}_{U(A\oplus M)}U(A))$
\end{sublem}

\begin{proof} We consider the commutative
diagram of simplicial sets
$$\xymatrix{
U(A)\times^{h}_{U(A\oplus M)}U(A) \ar[r]^{f} \ar[d]_-{p}  & F(A)\times^{h}_{F(A\oplus M)}F(A) \ar[d]^-{q} \\
U(A) \ar[r] & F(A)}$$
induced by the natural projection $A\oplus M \rightarrow A$.
Let $F(p)$ and $F(q)$ be the homotopy fibers
of the morphisms $p$ and $q$ taken at $y_{1}$ and $x_{1}$.
We have a natural morphism $g : F(p) \longrightarrow F(q)$.
Moreover, the homotopy fiber of the morphism $f$, taken at the point $x$,
receives a natural morphism from the
homotopy fiber of the morphism $g$. It is therefore
enough to show that the homotopy fiber of
$g$ is not empty. But, by definition of derivations, the
morphism $g$ is equivalent to the morphism
$$\Omega_{d,0}\mathbb{D}er_{U}(X,M) \longrightarrow \Omega_{d,0}\mathbb{D}er_{F}(X,M),$$
where the derivation $d$ is  given
by the image of the point $y_{1}$ by $d : A \rightarrow A\oplus M$.
Therefore, the homotopy fiber of the morphism $g$ is
equivalent to
$$\Omega_{d,0}\mathbb{D}er_{U/F}(X,M)\simeq
\Omega_{d,0}Map(\mathbb{L}_{U/F,y_{1}},M).$$
But, using Thm. \ref{t1} $(2)$ at rank $(n-1)$ and for
the morphism $U\longrightarrow F$ we obtain that
$\Omega_{d,0}\mathbb{D}er_{U/F}(X,M)$
is non empty. This finishes the proof
of the sub-lemma. \end{proof}

We now consider the commutative diagram
$$\xymatrix{
U(A\oplus_{d}\Omega M) \ar[r]^-{a} \ar[d]_-{b'} & U(A)\times^{h}_{U(A\oplus M)}U(A) \ar[d]^-{b} \\
F(A\oplus_{d}\Omega M) \ar[r]_-{a'}  & F(A)\times^{h}_{F(A\oplus M)}F(A).}$$
The morphism $a$ is an equivalence because $U$ is
representable. Furthermore, by our inductive assumption
the above square is homotopy cartesian. This implies
that the homotopy fiber
of $a'$ at $x$ is either contractible or empty. But, by the
above sub-lemma the
point $x$ lifts, up to homotopy, to a point
in $U(A)\times^{h}_{U(A\oplus M)}U(A)$, showing that this homotopy fiber is non empty.
\end{proof}

We have finished the proof of lemma \ref{l18} which implies that
any $n$-geometric stack is inf-cartesian, and thus that any
$n$-representable morphism has an obstruction theory. 
It only remain to
show part $(2)$ of Thm. \ref{t1}
at rank $n$. For this, we use Lem. \ref{lp16-} $(3)$ which
implies that we can assume that $G=*$. Let $U \longrightarrow
F$ be an $n$-atlas, $A\in \mathcal{A}$ and
$x : X:=\mathbb{R}\underline{Spec}\, A \longrightarrow F$ be a point.
By passing to an epimorphism of representable stacks
$X' \longrightarrow X$ which is in \textbf{E}, we can
suppose that the point $x$ factors through
a point $u : X \longrightarrow U$, where $U$ is representable
and $U\longrightarrow F$ is in \textbf{P}. By composition and the
hypothesis that morphisms in \textbf{P} are formally i-smooth, we see
that $U\longrightarrow *$ is a formally i-smooth morphism. We then have
a diagram
$$\mathbb{L}_{F,x} \longrightarrow
\mathbb{L}_{U,u} \longrightarrow \mathbb{L}_{X,x},$$
which obviously implies that
for any $M\in A-Mod_{1}$ the natural morphism
$$[\mathbb{L}_{X,x},M]\longrightarrow [\mathbb{L}_{F,x},M]$$
factors through
 the morphism
$$[\mathbb{L}_{X,x},M]\longrightarrow [\mathbb{L}_{U,u},M]$$
which is itself equal to zero by Prop. \ref{pismooth}. \end{proof}

We also extract from
Lem. \ref{l17} and its proof the following important corollary.

\begin{cor}\label{ct12}
Let $A\in Comm(\mathcal{C})$, $M\in A-Mod$ and
$d \in \pi_{0}(\mathbb{D}er(A,M))$ be a derivation.
Assume that the square
$$\xymatrix{
B=A\oplus_{d}\Omega M \ar[r] \ar[d] & A \ar[d] \\
A \ar[r] & A\oplus M}$$
is homotopy co-cartesian in $\mathcal{C}$, then
the base change functor
$$A\otimes_{B}^{\mathbb{L}} - : \mathrm{Ho}(B-Comm(\mathcal{C}))
\longrightarrow \mathrm{Ho}(A-Comm(\mathcal{C}))$$
induces an equivalence between the full sub-categories of formally
\'etale commutative $B$-algebras and formally \'etale
commutative $A$-algebras. The same statement holds
with \emph{formally \'etale} replaced by \emph{\'etale}.
\end{cor}

\begin{proof} Only the assertion with \emph{formally \'etale}
replaced by \emph{\'etale} requires an argument. For this, we only need
to prove that if a formally \'etale morphism
$f : B \longrightarrow B'$ is such that
$A \longrightarrow A\otimes_{B}^{\mathbb{L}}B'=A'$
is finitely presented, then so is $f$. For this we use the fully faithful functor
$$\mathbb{L}F : \mathrm{Ho}(B-Comm(\mathcal{C})) \longrightarrow \mathrm{Ho}(D)$$
defined during the proof of Lem. \ref{l16}. Using the description of
mappping spaces in $D$ in terms of a certain homotopy pullbacks, and
using the fact that filtered homotopy colimits in $SSet$ commutes with
homotopy pullbacks, we deduce the statement. \end{proof}
 
\part{Applications}

\chapter*{Introduction to Part 2}

In this second part we apply the theory developed in the first part to study the geometry
of stacks in various HAG contexts (Def. \ref{dhag}). In particular we will specialize our
base symmetric monoidal model category $\mathcal{C}$ to the following cases:
\begin{itemize}
  \item $\mathcal{C}=\mathbb{Z}-Mod$, the category of $\mathbb{Z}$-modules to get a theory
  of \textit{geometric stacks in} (classical) \textit{Algebraic Geometry} (\S \ref{IIsim});
	\item $\mathcal{C}=sMod_{k}$, the category of simplicial modules over an arbitrary base commutative ring $k$ to get a theory of \textit{derived} or $D^{-}$\textit{-geometric stacks} (\S \ref{IIder});  
	\item $\mathcal{C}=C(k)$, the category of unbounded cochain complexes of modules over a characteristic zero base commutative ring $k$ to get a theory of \textit{geometric stacks in} \textit{complicial algebraic geometry}, also called geometric $D$-stacks (\S \ref{IIunb});
	\item $\mathcal{C}=Sp^{\Sigma}$, the category of symmetric spectra (\cite{hss, shi}) to get a theory of \textit{geometric	stacks in brave new algebraic geometry} (\S \ref{IIbnag}).
\end{itemize}

In \S \ref{IIsim} we are concerned with \textbf{classical algebraic geometry}, the base category $\mathcal{C}=\mathbb{Z}-Mod$ being endowed with the trivial model structure.
We verify that if $k-Aff:=Comm(\mathcal{C})^{\mathrm{op}}=(k-Alg)^{\mathrm{op}}$ is endowed with its 
\'etale  Grothendieck topology and $\mathbf{P}$ is the class of smooth morphisms between (usual) commutative rings then $$(\mathcal{C},\mathcal{C}_{0},\mathcal{A},\tau,\textbf{P}):=(k-Mod,k-Mod,k-Alg,\textrm{\'et},\mathbf{P})$$
is a HAG-context according to Def. \ref{dhag}. $n$-geometric stacks in this context will be called
\textit{Artin} $n$\textit{-stacks}, and essentially coincide with geometric $n$-stacks as defined in \cite{s4}; their model category
will be denoted by $k-Aff^{\sim,\textrm{\'et}}$ and its homotopy category by $\mathrm{St}(k)$\index{$\mathrm{St}(k)$}.

After having established a coherent dictionary (Def. \ref{dII-2}), we show in 
\S  \ref{comparison} how the theory of schemes, of algebraic spaces, and of Artin's algebraic stacks
in groupoids (\cite{lm}) embeds in our theory of geometric stacks (Prop. \ref{pII-2}).
We also remark
(Rmk. \ref{why}) that the general infinitesimal theory for geometric stacks developed in Part I does not apply
to this context. The reason for this is that the category $\mathcal{C}=k-Mod$ (with its trivial model structure)
is as unstable as it could possibly be: the suspension functor $S:\mathrm{Ho}(\mathcal{C})=k-Mod 
\rightarrow k-Mod=\mathrm{Ho}(\mathcal{C})$ is trivial (i.e. sends each $k$-module to the zero $k$-module). The 
explanation for this is the following. Usual infinitesimal theory that applies to schemes, algebraic spaces
or to (some classes of) Artin's algebraic stacks in groupoids is in fact (as made clear e.g. by the definition of
cotangent complex of a scheme (\cite[II.2]{ill}) which uses simplicial resolutions of objects $Comm(\mathcal{C})$),
already conceptually part of \textit{derived algebraic geometry} in the sense that its classical
definition already requires to embed $Comm(\mathcal{C})=k-Alg$ into the category of simplicial $k$-algebras.
And in fact (as we will show in \S \ref{IIder}) when schemes, algebraic spaces and Artin algebraic stacks
in groupoids are viewed as \textit{derived stacks}, then their classical infinitesimal theory can be recovered 
(and generalized, see Cor. \ref{cpII-3}) and interpreted
 geometrically within our general formalism of Chapter \ref{partI.4} (in particular in Prop. \ref{p16}).\\
 
Therefore we are naturally brought to \S \ref{IIder} where we treat the case of \textbf{derived algebraic geometry}, i.e.
the case where $\mathcal{C}:=sk-Mod$, the category of simplicial modules over an arbitrary commutative ring $k$.

In \S \ref{IIder.1} we describe the model categories 
$\mathcal{C}=sk-Mod$ and $Comm(\mathcal{C})=sk-Alg$ whose opposite is denoted by $k-D^{-}Aff$, in particular finite
cell and finitely presented objects, suspension and loop functors, Postnikov towers and stable modules. 
We also show that $(\mathcal{C}, \mathcal{C}_{0}, \mathcal{A}):=(sk-Mod, sk-Mod, sk-Alg)$ is a HA context
in the sense of Def. \ref{dha}.

In \S \ref{IIder.2} we show how the general definitions of properties of modules (e.g. projective, flat,
perfect) and of morphisms between commutative rings in $\mathcal{C}$ (e.g. finitely presented, flat, (formally) smooth,
(formally) \'etale, Zariski open immersion) given in Chapter \ref{partI.2} translates concretely in the present context. The basic idea here is that of \textit{strongness} which says that a module $M$ over a simplicial $k$-algebra
$A$ (respectively, a morphism $A\rightarrow B$ in $sk-Alg$) has the property 
$\mathcal{P}$, defined in the abstract setting of Chapter \ref{partI.2}, if and only if $\pi_{0}(M)$ has the corresponding classical property as a 
$\pi_{0}(A)$-module and $\pi_{0}(M)\otimes_{\pi_{0}(A)} \pi_{*}(A) \simeq \pi_{*}(M)$ 
(respectively, the induced morphism $\pi_{0}(A)\rightarrow \pi_{0}(B)$ has the corresponding classical property, and
$\pi_{0}(B)\otimes_{\pi_{0}(A)} \pi_{*}(A) \simeq \pi_{*}(B)$). A straightforward extension
of the \'etale topology to simplicial $k$-algebras, then provides us with an \'etale model site $(k-D^{-}Aff,\textrm{\'et})$ satisfying assumption \ref{ass5} 
(Def. \ref{dII-5} and Lemma \ref{lII-5}) and with the corresponding model 
category $k-D^{-}Aff^{\sim, \textrm{\'et}}$ of \textrm{\'et} $D^{-}$-stacks (Def. \ref{dII-6}). 
The homotopy category $\mathrm{Ho}(k-D^{-}Aff^{\sim, \textrm{\'et}})$ of \textrm{\'et}
$D^{-}$-stacks will be simply denoted by $D^{-}\mathrm{St}(k)$. We conclude the section with two 
useful corollaries about topological invariance of \'etale and Zariski open immersions (Cor. \ref{ctII-1} and \ref{ctII-1zar}),
stating that ths small \'etale and Zariski site of a simplicial ring $A$ is
equivalent to the corresponding site of $\pi_{0}(A)$.  

In \S \ref{IIder.3} we describe our HAG context (Def. \ref{dhag}) for derived algebraic geometry by choosing 
the class \textbf{P} to be the class of smooth morphisms in $sk-Alg$ and the model topology to be the \'etale topology. This
HAG context $(\mathcal{C},\mathcal{C}_{0}, \mathcal{A}, \tau, \textbf{P}):=(sk-Mod, sk-Mod, sk-Alg, \textrm{\'et}, \mathrm{smooth})$ 
is shown to satisfy Artin's conditions (Def. \ref{pII-4}) relative to the HA context 
$(\mathcal{C}, \mathcal{C}_{0}, \mathcal{A}):=(sk-Mod, sk-Mod, sk-Alg)$ in Prop. \ref{d25}; as a corollary of the general theory
of Part I, this gives (Cor. \ref{cII-4}) an obstruction theory (respectively, a relative obstruction theory) for any $n$-geometric $D^{-}$-stack (resp., for any $n$-representable morphism between $D^{-}$-stacks), and in particular a (relative) cotangent complex for any $n$-geometric $D^{-}$-stack (resp., for any $n$-representable morphism). We
finish the section showing that the properties of being flat, smooth, \`etale and finitely presented can be extended to
$n$-representable morphisms between $D^{-}$-stacks (Lemma \ref{lII-7}), and by the definition of 
open and closed immersion of $D^{-}$-stacks (Def. \ref{dII-7}).

In \S \ref{IIder.4} we study  \textit{truncations} of derived stacks. The inclusion functor
$j:k-Aff \hookrightarrow sk-Alg$, that sends a commutative $k$-algebra $R$ to the constant simplicial $k$-algebra $R$,
is Quillen right adjoint to the functor $\pi_{0}:sk-Alg \rightarrow k-Alg$, and this Quillen adjunction induces an adjunction 
$$i:=\mathbb{L}j_{!} : \mathrm{St}(k) \longrightarrow D^{-}\mathrm{St}(k)$$
$$\mathrm{St}(k) \longleftarrow D^{-}\mathrm{St}(k) : \mathbb{R}i^{*}=:t_{0}$$ between the (homotopy) categories of derived and un-derived stacks.  The functor $i$ is fully faithful and commutes with homotopy colimits (Lemma \ref{lII-8}) and embeds the theory of stacks into the theory of derived stacks, while the functor
$t_{0}$, called the truncation functor, sends the affine stack corresponding to a simplicial $k$-algebra $A$
to the affine scheme $\mathrm{Spec}\;\pi_{0}(A)$, and commutes with homotopy limits and colimits (Lemma \ref{lII-8'}).
Both the inclusion and the truncation functor preserve $n$-geometric stacks and flat, smooth, \'etale 
morphisms between them (Prop. \ref{pII-5}). This gives a nice compatibility between the theories in \S \ref{IIsim}
and \S \ref{IIder}, and therefore between moduli spaces and their derived analogs.

In particular, we get the that for \textit{any} Artin algebraic stack in groupoids
(actually, any Artin $n$-stack) $\mathcal{X}$, $i(\mathcal{X})$ has an obstruction theory, and therefore a cotangent complex; in other words
viewing Artin stacks as derived stacks simplifies and clarifies a lot their infinitesimal theory, as already 
remarked in this introduction. We finish the section by showing (Prop. \ref{pII-6}) that for any geometric 
$D^{-}$-stack $F$, its truncation, viewed again as a derived stack (i.e. $it_{0}(F)$) sits inside $F$
as a closed sub-stack, and one can reasonably think of $F$ as behaving like a formal thickening of its
truncation.

In \S \ref{IIder.5} we give useful criteria for a $n$-representable morphism between $D^{-}$-stacks
being smooth (respectively, \'etale) in terms of locally finite presentation of the induced morphism on truncations, and infinitesimal lifting properties (Prop. \ref{pII-11}, resp., Prop. \ref{pII-12}) 
or properties of the cotangent complex (Cor. \ref{cpII-11}, resp., Cor. \ref{cpII-12}).

The final \S \ref{IIder.6} of Chapter \ref{IIder} contains applications of derived algebraic geometry to the construction of various derived versions of moduli spaces as $D^{-}$-stacks. In \ref{IIder.6.1} we first show (Lemma \ref{lII-10}) that
the stack of rank $n$ vector bundles when viewed as a $D^{-}$-stack using the inclusion $i : \mathrm{St}(k) \rightarrow D^{-}\mathrm{St}(k)$ is indeed isomorphic to the derived
stack $\mathbf{Vect}_{n}$ of rank $n$ vector bundles defined as in \S \ref{Iqcoh}. Then, for any simplicial set $K$, we define the derived stack $\mathbb{R}\mathbf{Loc}_{n}(K)$  as the derived exponentiation of 
$\mathbf{Vect}_{n}$ with respect to $K$ (Def. \ref{dII-11}), show that when $K$ is finite dimensional then 
$\mathbb{R}\mathbf{Loc}_{n}(K)$ is a finitely presented $1$-geometric $D^{-}$-stack (Lemma \ref{lII-11}), and
identify its truncation with the usual Artin stack of rank $n$ local systems on $K$ (Lemma \ref{lII-12}).
Finally we give a more concrete geometric interpretation of $\mathbb{R}\mathbf{Loc}_{n}(K)$ as a moduli space of derived geometric objects (\textit{derived rank $n$ local systems}) on the topological realization $\left|K\right|$
(Prop. \ref{pII-7}) and show that the tangent space of $\mathbb{R}\mathbf{Loc}_{n}(K)$ at a global point correponding
to a rank $n$ local system $E$ on $K$ is the cohomology complex $C^{*}(K,E\otimes_{k}E^{\vee})[1]$ (Prop. \ref{pII-8}).
The latter result shows in particular that the $D^{-}$- stack $\mathbb{R}\mathbf{Loc}_{n}(K)$ depends on strictly
more than the fundamental groupoid of $K$ (because its tangent spaces can be nontrivial even if $K$ is simply connected) and therefore carries higher homotopical informations as opposed to 
the usual (i.e underived) Artin stack of rank $n$ local systems. In subsection \ref{IIder.6.2} we treat the
case of algebras over an operad. If $\mathcal{O}$ is an operad in the category of $k$-modules, we consider
a simplicial presheaf $\mathbf{Alg}^{\mathcal{O}}_{n}$ on $k-D^{-}Aff$ which associates to any simplicial
$k$-algebra $A$ the nerve of the subcategory $\mathcal{O}-Alg(A)$ of weak equivalences in the category of cofibrant algebras over the operad $\mathcal{O}\otimes _{k}A$ (which is an operad in simplicial $A$-modules) whose underlying $A$-module is a vector bundle of rank $n$. In Prop. \ref{pII-15} we show that $\mathbf{Alg}^{\mathcal{O}}_{n}$
is a $1$-geometric quasi-compact $D^{-}$-stack, and in Prop. \ref{pII-16} we identify its tangent space
in terms of derived derivations. In subsection \ref{IIder.6.3}, using a special case of
J.Lurie's representability criterion (see Appendix $C$), we give sufficient conditions for the 
$D^{-}$-stack $\mathbf{Map}(\mathcal{X},F)$ of morphisms of derived stacks $i(\mathcal{X})\rightarrow F$ 
(for $\mathcal{X}\in \mathrm{St}(k)$) to be $n$-geometric (Thm. \ref{pII-17}),
and compute its tangent space in two particular cases (Cor. \ref{cpII-17} and \ref{cpII-17'}). The latter 
of these, i.e. the case $\mathbb{R}\mathcal{M}_{DR}(X):=\mathbf{Map}(i(X_{DR}),\mathbf{Vect}_{n})$, where
$X$ is a complex smooth projective variety, is particularly important because it is the first step
in the construction of a \textit{derived} version of \textit{non-abelian Hodge theory} which will be investigated in a future work.\\

In \S \ref{IIunb} we treat the case of \textbf{complicial} (or \textbf{unbounded}) \textbf{algebraic geometry},
i.e. homotopical algebraic geometry over the base category $\mathcal{C}:= C(k)$ of (unbounded) complexes
of modules over a commutative $\mathbb{Q}$-algebra $k$ \footnote{The case of $k$ of positive characteristic can be
treated as a special case of \textit{brave new algebraic geometry} (\S \ref{IIbnag}) over the base $Hk$, $H$ being the
Eilenberg-Mac Lane functor.}.

Here $Comm(\mathcal{C})$ is the model category $k-cdga$ of commutative differential graded $k$-algebras, 
(called shortly cdga's) whose opposite model category will be denoted by $k-DAff$. A new feature of complicial algebraic
geometry with respect to derived algebraic geometry is the existence of \textit{two} interesting HA contexts
(Lemma \ref{lcont}):
the \textit{weak} one $(\mathcal{C},\mathcal{C}_{0},\mathcal{A}):=(C(k),C(k),k-cdga)$,
 and the \textit{connective} one $(\mathcal{C},\mathcal{C}_{0},\mathcal{A}):=(C(k),C(k)_{\leq 0},k-cdga_{0})$
where $C(k)_{\leq 0}$ is the full subcategory of $C(k)$ consisting of complexes which are cohomologically
trivial in positive degrees and $k-cdga_{0}$ is the full subcategory of $k-cdga$ of connected algebras (i.e.
cohomologically trivial in non-zero degrees). A related phenomenon is the fact that the notion of 
\textit{strongness} that describes completely the properties of modules and of morphisms between
between simplicial algebras in derived algebraic geometry, is less strictly connected with standard properties 
of modules and of morphisms between
between cdga's. For example, morphisms between cdga's which have strongly the property $\mathcal{P}$ 
(e.g. $\mathcal{P}$= flat, \'etale etc.; see Def. \ref{unbstrong}) have the property  $\mathcal{P}$ 
but the converse is true only if additional hypotheses are met (Prop. \ref{pII-28}). 
Using the strong version of \'etale morphisms, we endow the category $k-DAff$ with the \textit{strongly \'etale
model topology} s-\'et (Def. \ref{dII-25} and Lemma \ref{lII-26} ), and define the model category 
of $D$\textit{-stacks} as $k-DAff^{\sim, \textrm{s-\'et}}$ (Def. \ref{dII-26}). The homotopy category of $k-DAff^{\sim, \textrm{s-\'et}}$
will be simply denoted by $D\mathrm{St}(k)$. 

While the notion
of stack does not depend on the HA or HAG contexts chosen but only the base model site,
the notion of geometric stack depends on the HA and HAG contexts. In \S \ref{IIunb.2}
we complete the weak HA context above to the HAG context 
$(\mathcal{C},\mathcal{C}_{0},\mathcal{A},\tau, \textbf{P}):=(C(k),C(k),k-cdga,\textrm{s-\'et},\mathbf{P}_{w})$,
where $\mathbf{P}_{w}$ is the class of formally perfect morphisms (Def. \ref{d7'}),  and call the corresponding
geometric stacks \textit{weakly geometric stacks} (Def. \ref{dII-27}).

Section \ref{IIunb.3} contains
some interesting examples of weakly geometric stacks. We first show (Prop. \ref{pperf}) 
that the stack $\mathbf{Perf}$ of perfect modules 
(defined in an abstract context in Def. \ref{d30'}) is weakly $1$-geometric, categorically locally of finite 
presentation and its diagonal is $(-1)$-representable. In subsection \ref{IIunb.3.2} we first define 
the simplicial presheaf $\mathbf{Ass}$ sending a cdga $A$ to the nerve of the subcategory of weak
equivalences in the category of associative and unital (not necessarily commutative) cofibrant $A$-dg algebras
whose underlying $A$-dg module is perfect, and show that this is a weakly $1$-geometric $D$-stack 
(Prop. \ref{pass} and Cor. \ref{cpass}). Then we define (using the model structure on dg-categories
of \cite{tab}) a dg-categorical variation of $\mathbf{Ass}$, denoted by $\mathbf{Cat}_{*}$, by 
sending a cdga $A$ to the nerve of the subcategory of weak equivalences
of the model category of $A$-dg categories $\mathcal{D}$ which are connected and have perfect and
cofibrant $A$-dg-modules of morphisms, and prove (Prop. \ref{pcat}) that the canonical classifying functor 
$\mathbf{Ass}\rightarrow \mathbf{Cat}_{*}$ is a weakly $1$-representable, $fp$ epimorphism of
$D$-stacks; it follows that $\mathbf{Cat}_{*}$ is a weakly $2$-geometric $D$-stack (Cor. \ref{pcat}). 
 
In \S \ref{IIunb.4} we switch to the connective HA context $(\mathcal{C},\mathcal{C}_{0},\mathcal{A}):=(C(k),C(k)_{\leq 0},k-cdga_{0})$ and complete it 
to the connective HAG context $$(C(k),C(k)_{\leq 0},k-cdga_{0},\textrm{s-\'et},\textrm{fip-smooth}),$$ where fip-smooth is the class of formally 
perfect (Def. \ref{d7}) and formally i-smooth morphisms (Def. \ref{dismooth}). Geometric
stacks in this HAG context will be simply called \textit{geometric} $D$\textit{-stacks} (Def. \ref{dII-27'});
since fip-smooth morphisms are in $\mathbf{P}_{w}$, any geometric $D$-stack is weakly geometric. As opposed to the
weak HAG context, this connective context indeed satisfies Artin's condition of 
Def. \ref{d25} (Prop. \ref{lII-27''}),
and therefore any geometric $D$-stack has an obstruction theory (and a cotangent complex).

In \S \ref{IIunb.5} we  give some examples of geometric $D$-stacks. We first observe (subsection \ref{IIunb.5.1})
 that the normalization functor $N: sk-Alg \rightarrow k-cdga$ induces a fully faithful functor 
$j:=\mathbb{L}N_{!}:D^{-}\mathrm{St}(k)\hookrightarrow D\mathrm{St}(k)$. 
This provides us with lots of examples of (geometric) $D$-stacks. 
 In subsection \ref{IIunb.5.2} we study the $D$-stack
 of CW-perfect modules. After having defined, for any cdga $A$, the notion of CW-A-dg-module of amplitude
 in $[a,b]$ (Def. \ref{dcwmod}) and proved some stability properties of this notion (Lemma \ref{cwstab}),
 we define the sub-$D$-stack $\mathbf{Perf}^{CW}_{[a,b]}\subset
\mathbf{Perf}$, consisting of all
perfect modules locally equivalent to some CW-dg-modules of amplitude
contained in $[a,b]$. We prove that $\mathbf{Perf}^{CW}_{[a,b]}$ is $1$-geometric (Prop. \ref{pcwmod}),
and that its tangent space at a point corresponding to a perfect CW-$A$-dg-module $E$ is given
by the complex $(E^{\vee}\otimes^{\mathbb{L}}_{A}E)[1]$ (Cor. \ref{cpcwmod}).
 In subsection 
\ref{IIunb.5.3} we define the $D$-stack of CW-dg-algebras as the homotopy pullback of $\mathbf{Ass}\rightarrow \mathbf{Perf}$ along the inclusion $\mathbf{Perf}^{CW}_{[a,b]}\rightarrow \mathbf{Perf}$, 
prove that it is $1$-geometric and compute its tangent space at a 
global point in terms of derived derivations (Cor. \ref{ccwass}). Finally subsection \ref{IIunb.5.4}
is devoted to the analysis of the $D$-stack $\mathbf{Cat}^{CW}_{*,\, [n,0]}$ 
of CW-dg-categories of perfect amplitude in $[n,0]$ (Def. \ref{dcwdgcat}). As opposite to the
weakly geometric $D$-stack $\mathbf{Cat}_{*}$, that cannot have a reasonable infinitesimal theory,
its full sub-$D$-stack $\mathbf{Cat}^{CW}_{*,\, [n,0]}$ is not only a $2$-geometric $D$-stack
but has a tangent space that can be computed in terms of Hochschild homology (Thm. \ref{tdg-cat}). 
As corollaries of this important result we can prove a folklore statement (see e.g. \cite[p. 266]{kos} in the
case of $A_{\infty}$-categories with one object) regarding the deformation theory
of certain negative dg-categories being controlled by the Hochschild complex  
of dg-categories (Cor. \ref{ctdg-cat}), and a result showing that if one wishes to keep
the existence of the cotangent complex, the restriction to non-positively graded dg-categories 
is unavoidable (Cor. \ref{ctdg-cat2}).\\

In the last, short \S \ref{IIbnag}, we establish the basics of \textbf{brave new algebraic 
geometry}, i.e. of homotopical algebraic geometry over the base category $\mathcal{C}=Sp^{\Sigma}$ 
of symmetric spectra (\cite{hss, shi}). We consider $\mathcal{C}$ endowed with 
the positive model structure of \cite{shi}, which is better behaved than the usual 
one when dealing with commutative monoid objects an modules over them. 
We denote $Comm(\mathcal{C})$ by $S-Alg$ (and call its objects commutative $S$-algebras,
$S$ being the sphere spectrum, or sometimes bn-rings), and its opposite model category by $S-Aff$.
Like in the 
case of complicial algebraic geometry, we consider two HA contexts here (Lemma \ref{connbnag}): 
$(\mathcal{C},\mathcal{C}_{0},\mathcal{A}):=(Sp^{\Sigma},Sp^{\Sigma}, S-Alg)$ and $(\mathcal{C},\mathcal{C}_{0},\mathcal{A}):=(Sp^{\Sigma},Sp^{\Sigma}_{c}, S-Alg_{0})$ where 
$Sp^{\Sigma}_{c}$ is the subcategory of connective symmetric spectra, and $S-Alg_{0}$ 
the subcategory of $S$-algebras with homotopy groups concentrated in degree zero.
After giving some examples of formally \'etale and formally $thh$-\'etale maps
between bn-rings, we define (Def. \ref{bnagstrong}) \textit{strong versions} of flat, 
(formally) \'etale, (formally) smooth, 
and Zariski open immersions, exactly like in chapters \ref{IIder} and \ref{IIunb}, and
give some results relating them to the corresponding non-strong notions (Prop. \ref{pII-28'}).  
An interesting exception to this relationship occurs in the case of smooth morphism:
the Eilenberg-MacLane functor $H$ from commutative rings to bn rings does not
preserve (formal) smoothness in general, though it preserves (formal) strong smoothness, due to
the presence of non-trivial Steenrod operations in characteristic $p>0$ (Prop. \ref{affineline}).
We endow the category $S-Aff$ with the \textit{strong \'etale} model topology s-\'et (Def. \ref{dII-25bis}
 and Lemma \ref{lII-26bis}),
 and define the model category $S-Aff^{\sim,\,\textrm{s-\'et}}$ of $S$\textit{-stacks} (Def. \ref{dII-26'}). 
 The homotopy category of  $S-Aff^{\sim,\,\textrm{s-\'et}}$ will be simply denoted by $\mathrm{St}(S)$. 
The two HA contexts defined above are completed (Cor. \ref{clII-27bis}) to two 
different HAG contexts by choosing for both HA contexts the s-\'et model topology, and the class  $\mathbf{P}_{\textrm{s-\'et}}$ of strongly 
\'etale morphsims for the first HA context (respectively, the class $\mathbf{P}$ of fip-morphisms  
for the second HA context). We call \textit{geometric Deligne-Mumford} $S$\textit{-stacks} (respectively, 
 \textit{geometric} $S$\textit{-stacks}) the geometric stacks in the first (resp, second) HAG
  context, and observe that both contexts satisfy Artin's condition of 
Def. \ref{d25} (Prop. \ref{artinbnag}).

Finally in \S \ref{IIbnag.2}, we use the definition of topological modular forms
of Hopkins-Miller to build a $1$-geometric Deligne-Mumford stack $\overline{\mathcal{E}}_{\mathbf{S}}$
that is a ``bn-derivation'' of the usual stack $\overline{\mathcal{E}}$ of generalized elliptic curves 
\footnote{See e.g. \cite[IV]{dera}, where it is denoted by $\mathcal{M}_{(1)}$.} (i.e., the truncation
of $\overline{\mathcal{E}}_{\mathbf{S}}$ is $\overline{\mathcal{E}}$), and such 
that the spectrum $\mathsf{tmf}$ coincide with the spectrum of \textit{functions on} 
$\overline{\mathcal{E}}_{\mathbf{S}}$. We conclude by the remark that
a moduli theoretic interpretation of $\overline{\mathcal{E}}_{\mathbf{S}}$ (or, most
probably some variant of it), i.e. finding out which are the brave new objects that it
classifies, could give not only interesting new geometry over bn rings but also
new insights on classical objects of algebraic topology.

\chapter{Geometric $n$-stacks in algebraic geometry (after C. Simpson)}\label{IIsim}

All along this chapter we fix an associative commutative ring $k\in \mathbb{U}$ with unit.

\section{The general theory}

We consider $\mathcal{C}=k-Mod$, the category of
$k$-modules in the universe $\mathbb{U}$. We endow
the category $\mathcal{C}$ with the trivial model structure for which
equivalences are isomorphisms and all morphisms
are cofibrations and fibrations. The category $\mathcal{C}$ is
furthermore a symmetric monoidal model category
for the monoidal structure given by the tensor
product of $k$-modules.
The
assumptions \ref{ass-1}, \ref{ass1}, \ref{ass0} and
\ref{ass2} are all trivially satisfied. The category
$Comm(\mathcal{C})$ is identified with the category
$k-Alg$, of commutative (associative and unital)
$k$-algebras in $\mathbb{U}$, endowed with the trivial
model structure. Objects in $k-Alg$ will simply be
called commutative $k$-algebras, without any reference
to the universe $\mathbb{U}$.
For any $A\in k-Alg$, the category
$A-Mod$ is the usual symmetric monoidal category of
$A$-modules in $\mathbb{U}$, also endowed with its
trivial model structure. Furthermore,
we have $\mathbb{R}\underline{Hom}_{A}(M,N)\simeq
\underline{Hom}_{A}(M,N)$ for any two objects $A \in k-Alg$, and
is the usual $A$-module of morphisms from $M$ to $N$.
We set $\mathcal{C}_{0}:=k-Mod$, and
$\mathcal{A}:=k-Alg$. The triplet
$(k-Mod,k-Mod,k-Alg)$
is then a HA context in the sense of Def. \ref{dha}.

The category $Aff_{\mathcal{C}}$ is identified with
$k-Alg^{op}$, and therefore to the category
of affine $k$-schemes in $\mathbb{U}$. It will simply be denoted
by $k-Aff$ and its objects will simply be called
affine $k$-schemes, without any reference to the
universe $\mathbb{U}$. The model category of pre-stacks $k-Aff^{\wedge}=Aff_{\mathcal{C}}^{\wedge}$
is simply the model category of $\mathbb{V}$-simplicial
presheaves on $k-Aff$, for which equivalences and fibrations
are defined levelwise. The Yoneda functor
$$h : k-Aff \longrightarrow k-Aff^{\wedge}$$
is the usual one, and sends an affine $k$-scheme $X$ to the
presheaf of sets it represents $h_{X}$ (considered
as a presheaf of constant simplicial sets). Furthermore we have natural
isomorphisms of functors
$$h\simeq \underline{h} \simeq \mathbb{R}\underline{h} :
k-Aff \longrightarrow \mathrm{Ho}(k-Aff^{\wedge}),$$
which is nothing else than the natural composition
$$\xymatrix{k-Aff \ar[r]^-{h} & k-Aff^{\wedge} \ar[r] & \mathrm{Ho}(k-Aff^{\wedge}).}$$

We let $\tau=\textrm{\'et}$, the usual \'etale pre-topology
on $k-Aff$ (see e.g. \cite{mil}). Recall that a family of morphisms
$$\{X_{i}=Spec\, A_{i} \longrightarrow X=Spec\, A\}_{i\in I}$$
is an $\textrm{\'et}$-covering family\index{$\textrm{\'et}$-covering family!in $k-Aff$} if and only if it contains a finite
sub-family $\{X_{i} \longrightarrow X\}_{i\in J}$, $J\subset I$,
such that the corresponding morphism of commutative $k$-algebras
$$A \longrightarrow \prod_{i\in J}A_{i}$$
is a faithfully flat and \'etale morphism of commutative rings.

\begin{lem}\label{lII-1}
The \'etale topology on $k-Aff$ satisfies assumption
\ref{ass5}.
\end{lem}

\begin{proof} Points $(1)$ and $(2)$ of \ref{ass5}
are clear. Point $(3)$ is induced by the
faithfully flat descent of quasi-coherent modules
for affine $ffqc$-hypercovers (see e.g. \cite[Exp. $V^{\mathrm{bis}}$]{sga4II})). \end{proof}

The model category of stacks $k-Aff^{\sim,\textrm{\'et}}$ is the
projective model structure for simplicial presheaves
on the Grothendieck site $(k-Aff,\textrm{\'et})$, as defined
for example in \cite{bl} (see also \cite[\S 1]{to1} ). Its homotopy category, denoted simply as
$\mathrm{St}(k)$\index{$\mathrm{St}(k)$}, is then identified with the
full subcategory of $\mathrm{Ho}(k-Aff^{\wedge})$ consisting
of simplicial presheaves
$$F : k-Alg=k-Aff^{op} \longrightarrow SSet_{\mathbb{V}}$$
satisfying the following two
conditions
\begin{itemize}
\item For any two commutative $k$-algebra
$A$ and $B$, the natural morphism
$$F(A\times B) \longrightarrow F(A)\times F(B)$$
is an isomorphism in $\mathrm{Ho}(SSet)$.
\item For any co-augmented
co-simplicial commutative $k$-algebra, $A \longrightarrow B_{*}$, such that
the augmented simplicial object
$$Spec\, B_{*} \longrightarrow Spec\, A$$
is an \'etale hypercover, the natural morphism
$$F(A) \longrightarrow Holim_{n \in \Delta}F(B_{n})$$
is an isomorphism in $\mathrm{Ho}(SSet)$.
\end{itemize}

It is well known (and also a consequence of Cor. \ref{c0}) that the \'etale topology is sub-canonical;
therefore there exists a fully faithful functor
$$h : k-Aff \longrightarrow \mathrm{St}(k)\subset \mathrm{Ho}(k-Aff^{\sim,\textrm{\'et}}).$$
Furthermore, we have
$$Spec\, A\simeq \underline{Spec}\, A\simeq \mathbb{R}\underline{Spec}\, A$$
for any $A\in k-Alg$.

We set \textbf{P} to be the class of \textit{smooth morphisms}
in $k-Aff$ in the sense of \cite[17.3.1]{egaIV-4}. 
It is well known that
our assumption \ref{ass4} 
is satisfied (e.g. that smooth morphisms are \'etale-local in the source and target, 
see for example \cite[17.3.3, 17.3.4, 17.7.3]{egaIV-4}).
In particular, we get that
$(k-Mod,k-Mod,k-Alg,\textrm{\'et},\mathbf{P})$ is
a HAG context in the sense of Def. \ref{dhag}.
The general definition \ref{d11} can then be applied, and provides
a notion of $n$-geometric stack in $\mathrm{St}(k)$.
A first important observation (the lemma below) is that
$n$-geometric stacks are $n$-stacks in the sense of
\cite{s4}. This is a special feature of standard algebraic geometry,
and the same would be true for any theory for which the
model structure on $\mathcal{C}$ is trivial: the
\emph{geometric complexity} is a bound for the \emph{stacky
complexity}.

Recall that a stack $F \in \mathrm{St}(k)$
is \emph{$n$-truncated} if for any
$X\in k-Aff$, any $s\in \pi_{0}(F(X))$ and any
$i>n$, the sheaf $\pi_{i}(F,s)$ is trivial. By \cite[3.7]{hagI}, this
is equivalent to say that
for any stack $G$ the simplicial set
$\mathbb{R}_{\tau}\underline{Hom}(G,F)$
is $n$-truncated.

\begin{lem}\label{lII-2}
Let $F$ be an $n$-geometric stack in $\mathrm{St}(k)$.
Then $F$ is 
$(n+1)$-truncated.
\end{lem}

\begin{proof} The proof is by induction on $n$.
Representable stacks are nothing else than
affine schemes, and therefore are $0$-truncated.
Suppose that the lemma is known for
$m<n$. Let $F$ be an $n$-geometric stack,
$X\in k-Aff$ and $s\in \pi_{0}(F(X))$. We have a natural isomorphism of sheaves
on $k-Aff/X$
$$\pi_{i}(F,s)\simeq \pi_{i-1}(X\times_{F}^{h}X,s),$$
where $X \longrightarrow F$ is the morphisms of stacks
corresponding to $s$. As the diagonal of $F$ is $(n-1)$-representable,
$X\times_{F}^{h}X\simeq F\times_{F\times^{h}F}X\times^{h}X$ is
$(n-1)$-geometric. By induction we find
that $\pi_{i}(F,s)\simeq *$ for any $i>n$. \end{proof}

Lemma \ref{lII-2} justifies the following terminology, closer
to the usual terminology one can find in the literature.

\begin{df}\label{dII-1}
An \emph{Artin $n$-stack}\index{Artin $n$-stack} is an $n$-truncated stack which is
$m$-geometric for some integer $m$.
\end{df}

The general theory of Artin $n$-stacks could then be
pursued in a similar fashion as for
Artin stacks in \cite{lm}. A part of this is done
in \cite{s4} and will not be reproduced here, as many 
of these statements will be settled down in the more general 
context of geometric $D^{-}$-stacks (see \S \ref{IIder}).
Let us mention however, that as explained
in Def. \ref{d21}, we can define the notions of
flat, smooth, \'etale, unramified, regular, Zariski open
immersion \dots morphisms
between Artin $n$-stacks. These kinds of morphisms are
as usual stable by homotopy pullbacks, compositions
and equivalences. In particular this allows
the following definition.

\begin{df}\label{dII-2}
\begin{itemize}
\item
An Artin $n$-stack is a \emph{Deligne-Mumford $n$-stack}\index{Deligne-Mumford $n$-stack} if
there exists an $n$-atlas $\{U_{i}\}$ for $F$ such that
each morphism $U_{i} \longrightarrow F$
is an \'etale morphism.
\item An Artin $n$-stack is an \emph{algebraic space}\index{algebraic space} if
it is a Deligne-Mumford $n$-stack, and if
furthermore
the diagonal $F \longrightarrow F\times^{h}F$
is a monomorphism in the sense of Def. \ref{d23}.
\item An Artin $n$-stack $F$ is an \emph{scheme}\index{scheme} if
there exists an $n$-atlas $\{U_{i}\}$ for $F$ such that
each morphism $U_{i} \longrightarrow F$
is a monomorphism.
\end{itemize}
\end{df}

\begin{rmk}\label{remnstack}
\begin{enumerate}
\item \emph{An algebraic space in the sense of the definition above which is automatically a $1$-geometric stack, 
and is nothing else than
an algebraic space in the usual sense. Indeed, this can be shown by induction on $n$:
an algebraic space which is also $n$-geometric is by definition the quotient of a union of affine schemes $X$ by some
\'etale equivalence relation $R \subset X\times X$ where $R$ is an algebraic space which is $(n-1)$-geometric. In particular, $R$ being a subobject in 
$X\times X$ we see that $R$ is a separated algebraic space, and thus is a 
$0$-geometric stack. This implies that 
$X/R$ is a $1$-geometric stack. In the same way, any scheme is automatically a $1$-geometric stack. Moreover, 
algebraic spaces (resp. schemes) which are $0$-geometric stacks are precisely algebraic spaces (resp. schemes) with 
an affine diagonal. }
\item \emph{Thought there is a small discrepancy between the notion
of Artin $n$-stack and the notion of $n$-geometric stack in $\mathrm{St}(k)$, 
our notion of Artin $n$-stack is equivalent to the notion of \emph{slightly geometric
n-stacks} of \cite{s4}. 
 }
\end{enumerate}
\end{rmk}

\section{Comparison with Artin's algebraic stacks}\label{comparison}

Artin $n$-stacks as defined in the last section are simplicial presheaves, whereas
Artin stacks are usually presented in the literature using the theory of
fibered categories (se e.g. \cite{lm}). In this section we 
briefly explain how the theory of fibered categories in groupoids 
can be embedded in the theory of simplicial presheaves, and
how this can be used in order to compare
the original definition of Artin stacks to our
definition of Artin $n$-stacks. \\

In \cite{hol}, it is shown that there exists a model
category $Grpd/\mathcal{S}$, of cofibered categories in
groupoids over a Grothendieck site $\mathcal{S}$. The fibrant
objects for this model structure are precisely the
stacks in groupoids in the sense of \cite{lm}, and
the equivalences in $Grpd/\mathcal{S}$ are the morphisms
of cofibered categories becoming equivalences on the
associated stacks (i.e. local equivalences). There exists furthermore
a Quillen equivalence
$$p : P(\mathcal{S},Grpd) \longrightarrow Grpd/\mathcal{S} \qquad
P(\mathcal{S},Grpd) \longleftarrow Grpd/\mathcal{S} : \Gamma,$$
where $P(\mathcal{S},Grpd)$ is the local projective model category of
presheaves of groupoids on $\mathcal{S}$. Finally, there exists a
Quillen adjunction
$$\Pi_{1} : SPr_{\tau}(\mathcal{S}) \longrightarrow P(\mathcal{S},Grpd) \qquad
SPr_{\tau}(\mathcal{S}) \longleftarrow P(\mathcal{S},Grpd) : B,$$
where $\Pi_{1}$ is the natural extension of the functor
sending a simplicial set to its fundamental groupoid, $B$ is the
natural extension of the nerve functor, and the
$SPr_{\tau}(\mathcal{S})$ is the local projective model
structure of simplicial presheaves on the site $\mathcal{S}$
(also denoted by $\mathcal{S}^{\sim,\tau}$ in our context, at 
least when $\mathcal{S}$ has limits and
colimits and thus can be considered as a model
category with the trivial model structure).
The functor $B$ preserves equivalences and the
induced functor
$$B : \mathrm{Ho}(P(\mathcal{S},Grpd)) \longrightarrow \mathrm{Ho}(SPr_{\tau}(\mathcal{S}))$$
is fully faithful and its image consists of all
$1$-truncated objects (in the sense of \cite[\S 3.7]{hagI}). Put in another way, the model
category $P(\mathcal{S},Grpd)$ is Quillen equivalent to the
$S^{2}$-nullification of $SPr_{\tau}(\mathcal{S})$ (denoted
by $SPr_{\tau}^{\leq 1}(\mathcal{S})$ in \cite[\S 3.7]{hagI}). In conclusion, there
exists a chain of Quillen equivalences
$$Grpd/\mathcal{S}\rightleftarrows P(\mathcal{S},Grpd) \rightleftarrows SPr_{\tau}^{\leq 1}(\mathcal{S}),$$
and therefore a well defined adjunction
$$t : \mathrm{Ho}(SPr_{\tau}(\mathcal{S})) \longrightarrow \mathrm{Ho}(Grpd/\mathcal{S}) \qquad
\mathrm{Ho}(SPr_{\tau}(\mathcal{S})) \longleftarrow \mathrm{Ho}(Grpd/\mathcal{S}) : i,$$
such that the right adjoint
$$i : \mathrm{Ho}(Grpd/\mathcal{S}) \longrightarrow \mathrm{Ho}(SPr_{\tau}(\mathcal{S}))$$
is fully faithful and its image consists of all $1$-truncated objects.

The category $\mathrm{Ho}(Grpd/\mathcal{S})$ can also be described as the
category whose objects are stacks in groupoids in the sense
of \cite{lm}, and whose morphisms are given by $1$-morphisms of stacks
up to $2$-isomorphisms. In other words, for two given stacks $F$ and $G$ in
$Grpd/\mathcal{S}$, the set of morphisms from $F$ to $G$ in
$\mathrm{Ho}(Grpd/\mathcal{S})$ is the set of isomorphism classes
of the groupoid $\mathcal{H}om(F,G)$, of morphisms of stacks. This implies that the
usual category of stacks in groupoids, up to $2$-isomorphisms, can be
identified through the functor $i$ with the full subcategory of
$1$-truncated objects in
$\mathrm{Ho}(SPr_{\tau}(\mathcal{S}))$. Furthermore, the functor $i$ being
defined as the composite of right derived functors and derived
Quillen equivalences will commutes with homotopy limits. As homotopy limits
in $Grpd/\mathcal{S}$ can also be identified with the $2$-limit of stacks
as defined in \cite{lm}, the functor $i$ will send
$2$-limits to homotopy limits. As a particular case we obtain that
$i$ sends the $2$-fiber product of stacks in groupoids to the homotopy fiber product.

The $2$-categorical structure of stacks in groupoids can also be recovered
from the model category $SPr_{\tau}(\mathcal{S})$. Indeed, applying the
simplicial localization techniques of \cite{dk1} to the Quillen adjunctions
described above, we get a well defined diagram of $S$-categories
$$L(Grpd/\mathcal{S}) \longrightarrow LP(\mathcal{S},Grpd) \longrightarrow
LSPr_{\tau}(\mathcal{S}),$$
which is fully faithful in the sense of \cite[Def. 2.1.3]{hagI}. In particular, the $S$-category
$L(Grpd/\mathcal{S})$ is naturally equivalent to the full sub-$S$-category
of $LSPr_{\tau}(\mathcal{S})$ consisting of $1$-truncated objects.
Using \cite{dk3}, the $S$-category $L(Grpd/\mathcal{S})$ is also equivalent to the
$S$-category whose objects are stacks in groupoids, cofibrant as objects in
$Grpd/\mathcal{S}$, and whose morphisms simplicial sets are given by
the simplicial $Hom$'s sets of $Grpd/\mathcal{S}$. These simplicial $Hom$'s sets are
simply the nerves of the groupoid of functors between cofibered categories in
groupoids. In other words, replacing the simplicial sets of morphisms
in $L(Grpd/\mathcal{S})$ by their fundamental groupoids, we find
a $2$-category naturally $2$-equivalent to the usual $2$-category of
stacks in groupoids on $\mathcal{S}$. Therefore, we see that
the $2$-category of stacks in groupoids can be identified, up to
a natural $2$-equivalence, as the $2$-category obtained from
the full sub-$S$-category of $LSPr_{\tau}(\mathcal{S})$ consisting
of $1$-truncated objects, by replacing its simplicial sets of morphisms by
their fundamental groupoids. \\

We now come back to the case where $\mathcal{S}=(k-Aff,\textrm{\'e}t)$, the Grothendieck
site of affine $k$-schemes with the \'etale pre-topology. We have seen that there exists a fully faithful functor
$$i : \mathrm{Ho}(Grpd/k-Aff^{\sim,\textrm{\'e}t}) \longrightarrow \mathrm{St}(k),$$
from the category of stacks in groupoids up to $2$-isomorphisms, to the
homotopy category of stacks. The image of this functor consists of all
$1$-truncated objects and it is compatible with the simplicial
structure (i.e. possesses a natural lifts as a morphism of $S$-categories).
We also have seen that $i$ sends $2$-fiber products of stacks to
homotopy fiber products.

Using the functor $i$, every stack in groupoids can be seen as
an object in our category of stacks $\mathrm{St}(k)$. For example,
all examples of stacks presented in \cite{lm} give rise to
stacks in our sense. The proposition below subsumes the main
properties of the functor $i$, relating the usual notion of scheme, algebraic space and
stack to the one of our definition Def. \ref{dII-2}.
Recall that a stack in groupoids $X$ is separated (resp. quasi-separated) if its
diagonal is a proper (resp. separated) morphism.

\begin{prop}\label{pII-2}
\begin{enumerate}
\item For any commutative $k$-algebra $A$, there exists a natural isomorphism
$$i(Spec\, A)\simeq Spec\, A.$$

\item If $X$ is a scheme (resp. algebraic space, resp.
Deligne-Mumford stack, resp. Artin stack) with an
affine diagonal in the sense of \cite{lm}, then
$i(X)$ is an Artin $0$-stack which is $0$-geometric (resp. an Artin $0$-stack which is $0$-geometric, resp.
a Deligne-Mumford $1$-stack which is $0$-geometric, resp. an Artin $1$-stack which is $0$-geometric)
in the sense of Def. \ref{dII-2}. 

\item If $X$ is a 
scheme (resp. algebraic space, resp.
Deligne-Mumford stack, resp. Artin stack) in the sense of \cite{lm}, then
$i(X)$ is an Artin $0$-stack which is $1$-geometric (resp. an Artin $0$-stack which is $1$-geometric, resp.
a Deligne-Mumford $1$-stack which is $1$-geometric, resp. an Artin $1$-stack which is $1$-geometric) in the sense of Def. \ref{dII-2}.

\item Let $f : F \longrightarrow G$ be a morphism between
Artin stacks in the sense of \cite{lm}. Then
the morphism $f$ is flat (resp. smooth, resp. \'etale,
resp. unramified, resp. Zariski open immersion)
if and only if $i(f) : i(F) \longrightarrow i(G)$ is so.

\end{enumerate}
\end{prop}

\begin{proof} This readily follows from the definition using the fact that
$i$ preserves affine schemes, epimorphisms of stacks, and sends $2$-fiber  products to
homotopy pullbacks. \end{proof}

\begin{rmk}\label{why}\emph{If we try to apply the general infinitesimal and obstruction theory developed in
\S \ref{partI.4} to the present HAG-context, we immediately see that this is impossible because the
suspension functor $S:\mathrm{Ho}(\mathcal{C})= k-Mod \rightarrow k-Mod=\mathrm{Ho}(\mathcal{C})$ is trivial.
On the other hand, the reader might object that there is already a well established infinitesimal theory,
at least in the case of schemes, algebraic spaces and for a certain class of algebraic stacks in groupoids, and
that our theory does not seem to be able to reproduce it. The answer to this question turns out to be both conceptually
and technically relevant. First of all, if we look at e.g. the definition of the cotangent complex of a scheme (\cite[2.1.2]{ill}) we realize 
that a basic and necessary step is to enlarge the category of rings to the category
of simplicial rings  in order to be able to consider free (or more generally cofibrant) resolutions of maps between rings. In our setup, this can be reformulated by saying that in order to get the correct infinitesimal 
theory, even for ordinary schemes, it is necessary to view them as geometric objects in \textit{derived algebraic geometry}, i.e. on homotopical algebraic geometry over the base category $\mathcal{C}=sk-Mod$
of simplicial $k$-modules (so that $Comm(\mathcal{C})$ is exactly the category of commutative simplicial 
$k$-algebras). In other words, the usual infinitesimal theory of schemes is already ``secretly'' a part of
derived algebraic geometry, that will be studied in detail in the next chapter \ref{IIder}. 
Moreover, as it will be
shown, this approach has, even for classical objects like schemes or Artin stacks in groupoids, both conceptual advantages (like e.g. the fact that the cotangent complex of
a scheme can be interpreted geometrically as a \textit{genuine} cotangent space to the scheme when viewed as a derived stack,
satisfying a natural universal property, while the cotangent complex of a scheme do not have any
universal property inside the theory of schemes),
and technical advantages (like the fact, proved in Cor. \ref{cpII-3}, that \textit{any} Artin stack in groupoids has an
obstruction theory). }
\end{rmk}

\begin{conv}\emph{From now on we will omit mentioning the functor $i$, and
will simply view stacks in groupoids as objects
in $\mathrm{St}(k)$. In particular, we will allow
ourselves to use the standard notions and vocabulary of
the general theory of schemes.}
\end{conv}

\chapter{Derived algebraic geometry}\label{IIder}

All along this chapter $k$ will be a fixed commutative (associative and unital) ring.

\section{The HA context}\label{IIder.1}

In this section we specialize our general theory of Part I 
to the case where $\mathcal{C}=sk-Mod$, is the category
of simplicial $k$-modules in the universe $\mathbb{U}$. The category
$sk-Mod$ is endowed with its standard model
category structure, for which the fibrations and equivalences are
defined on the underlying simplicial sets (see
for example \cite{gj}). The tensor product of
$k$-modules extends naturally to a levelwise tensor product on
$sk-Mod$, making it into a symmetric monoidal model
category. Finally, $sk-Mod$ is known to be a
$\mathbb{U}$-combinatorial proper and simplicial model category.

The model category $sk-Mod$ is known to satisfy
assumptions \ref{ass-1}, \ref{ass1} and \ref{ass0} (see \cite{schw-shi}).
Finally, it follows easily from \cite[II.4, II.6]{q0} that
$sk-Mod$ also satisfies assumption
\ref{ass2}. \\

The category $Comm(sk-Mod)$ will be denoted by $sk-Alg$, and its
objects will be called simplicial commutative  $k$-algebras. More
generally, for $A\in sk-Alg$, the category $A-Comm(\mathcal{C})$ will be
denoted by $A-Alg_{s}$. For any $A\in sk-Alg$ we will denoted by
$A-Mod_{s}$ the category of $A$-modules in $sk-Mod$, which is
nothing else than the category of simplicial modules over the
simplicial ring $A$. The model structure on $sk-Alg$, $A-Alg_{s}$ and
$A-Mod_{s}$ is the usual one, for which the equivalences and
fibrations are defined on the underlying simplicial sets. For an
object $A\in sk-Alg$, we will denote by $\pi_{i}(A)$ its homotopy
group (pointed at $0$). The graded abelian group $\pi_{*}(A)$
inherits a structure of a commutative graded algebra from $A$,
which defines a functor $A \mapsto \pi_{*}(A)$ from $sk-Alg$ to
the category of commutative graded $k$-algebras. More generally,
if $A$ is a simplicial commutative $k$-algebra, and $M$ is an
$A$-module, the graded abelian group $\pi_{*}(M)$ has a natural
structure of a graded $\pi_{*}(A)$-module.

There exists a Quillen adjunction
$$\pi_{0} : sk-Alg \longrightarrow k-Alg \qquad
sk-Alg \longleftarrow k-Alg : i,$$
where $i$ sends a commutative $k$-algebra to the corresponding
constant
simplicial commutative $k$-algebra. This Quillen adjunction induces a fully faithful functor
$$i : k-Alg \longrightarrow \mathrm{Ho}(sk-Alg).$$
From now on we will omit to mention the functor $i$, and always
consider $k-Alg$ as embedded in $sk-Alg$, except if the contrary
is specified. Note that when $A\in k-Alg$, also considered
as an object in $sk-Alg$, we have two different notions of $A$-modules,
$A-Mod$, and $A-Mod_{s}$. The first one is the usual category of
$A$-modules, whereas the second one is the category of simplicial objects
in $A-Mod$. 

For any morphism of simplicial commutative $k$-algebras
$A \longrightarrow B$, the $B$-module
$\mathbb{L}_{B/A}$ constructed in \ref{p1} is
naturally isomorphic in $\mathrm{Ho}(B-Mod)$ to
D. Quillen's cotangent complex introduced in
\cite{q}. In particular, if $A \longrightarrow B$ is a morphism
between (non-simplicial) commutative $k$-algebras, then we have
$\pi_{0}(\mathbb{L}_{B/A})\simeq \Omega^{1}_{B/A}$. More generally,
we find by adjunction
$$\pi_{0}(\mathbb{L}_{B/A})\simeq \Omega^{1}_{\pi_{0}(B)/\pi_{0}(A)}.$$
Recall also that $A \longrightarrow B$ in $k-Alg$ is \'etale in the
sense of \cite[17.1.1]{egaIV-4} if and only if
$\mathbb{L}_{B/A}\simeq 0$ and $B$ is finitely presented as a commutative $A$-algebra. 
In the same way, a morphism
$A \longrightarrow B$ in $k-Alg$ is smooth
in the sense of \cite[17.1.1]{egaIV-4} if and only if
$\Omega^{1}_{B/A}$ is a projective $B$-module,
$\pi_{i}(\mathbb{L}_{B/A})\simeq 0$ for $i>0$ and $B$ is finitely presented as a commutative $A$-algebra..
Finally, recall the existence of a natural first quadrant spectral sequence
(see \cite[II.6 Thm. 6(b)]{q0})
$$Tor^{p}_{\pi_{*}(A)}(\pi_{*}(M),\pi_{*}(N))_{q} \Rightarrow
\pi_{p+q}(M\otimes^{\mathbb{L}}_{A}N),$$
for $A\in sk-Alg$ and any objects $M$ and $N$ in $A-Mod_{s}$.

We set $\mathcal{C}_{0}:=\mathcal{C}=sk-Mod$, and
$\mathcal{A}:=sk-Alg$. Then, clearly,
assumption \ref{ass7} is also satisfied. The triplet
$(sk-Mod,sk-Mod,sk-Alg)$ is then a HA context
in the sense of Def. \ref{dha}.
Note that
for any $A\in sk-Alg$ we have $(A-Mod_{s})_{0}=A-Mod$, whereas
$(A-Mod_{s})_{1}$ consists of all $A$-modules $M$ such that
$\pi_{0}(M)=0$, also called \emph{connected modules}. \\

For an integer $n\geq 0$ we define
the $n$-th sphere $k$-modules by
$S^{n}_{k}:=S^{n}\otimes k \in sk-Mod$. The free commutative monoid on
$S^{n}_{k}$ is an object $k[S^{n}] \in sk-Alg$, such that for any
$A\in sk-Alg$ there are functorial isomorphisms
$$[k[S^{n}],A]_{sk-Alg}\simeq \pi_{n}(A).$$
In the same way we define $\Delta^{n}\otimes k\in sk-Mod$ and
its associated free commutative monoid $k[\Delta^{n}]$.
There are natural morphisms $k[S^{n}] \longrightarrow k[\Delta^{n+1}]$, coming from the
natural inclusions
$\partial\Delta^{n+1}=S^{n} \hookrightarrow \Delta^{n+1}$.
The set of morphisms
$$\{k[S^{n}] \longrightarrow k[\Delta^{n+1}]\}_{n\geq 0}$$
form a generating set of cofibrations in $sk-Alg$. 
The model category $sk-Alg$ is then easily checked to be
compactly generated in the sense of Def. \ref{dI-cell}.
A finite cell object in $sk-Alg$ is then any object $A\in sk-Alg$
for which there exists
a finite sequence in $sk-Alg$
$$\xymatrix{
A_{0}=k \ar[r] & A_{1} \ar[r] & A_{2} \ar[r] & \dots \ar[r] & A_{m}=A,}$$
such that for any $i$ there exists a push-out square in $sk-Alg$
$$\xymatrix{
A_{i} \ar[r] & A_{i+1} \\
k[S^{n_{i}}] \ar[u] \ar[r] & k[\Delta^{n_{i}+1}] \ar[u].}$$
Our Prop. \ref{pI-cell} implies that an
object $A\in sk-Alg$ is finitely presented in the sense of
Def. \ref{d3} if and only if it is equivalent to a retract of
a finite cell object (see also \cite[III.2]{ekmm} or \cite[Thm. III.5.7]{km} for other
proofs). More generally,
for $A\in sk-Alg$, there exists a notion of finite cell
object in $A-Alg_{s}$ using the elementary morphisms
$$A[S^{n}]:=A\otimes_{k}k[S^{n}] \longrightarrow A[\Delta^{n+1}]:=A\otimes_{k}k[\Delta^{n+1}].$$
In the same way, a morphism $A \longrightarrow B$ in $sk-Alg$ is finitely presented
in the sense of Def. \ref{d3} if and only if $B$ is equivalent to a retract of
a finite cell objects in $A-Alg_{s}$. Prop. \ref{pI-cell} also implies that
any morphism $A\longrightarrow B$, considered as an object in
$A-Alg_{s}$, is equivalent to a filtered colimit of
finite cell objects, so in particular to a filtered
homotopy colimit of finitely presented objects.

The Quillen adjunction between $k-Alg$ and $sk-Alg$ shows that the functor
$$\pi_{0} : sk-Alg \longrightarrow k-Alg$$
does preserve finitely presented morphisms. On the contrary, the inclusion functor
$i : k-Alg \longrightarrow sk-Alg$ does not preserve finitely presented objects, and
the finite presentation condition in $sk-Alg$ is in general stronger than
in $k-Alg$.

For $A\in sk-Alg$, we also have a notion of finite cell objects in $A-Mod_{s}$, based 
the generating set for cofibrations consisting of morphisms of the form
$$S^{n}_{A}:=A\otimes_{k}S^{n}_{k} \longrightarrow \Delta^{n+1}_{A}:=A\otimes_{k}\Delta^{n+1}_{A}.$$
Using Prop. \ref{pI-cell} we see that the finitely presented objects in 
$A-Mod_{s}$ are the objects equivalent to a retract
of a finite cell objects (see also \cite[III.2]{ekmm} or \cite[Thm. III.5.7]{km}). Moreover, the  functor
$$\pi_{0} : A-Mod_{s} \longrightarrow \pi_{0}(A)-Mod$$
is left Quillen, so preserves finitely presented objects. On the contrary, for
$A\in k-Alg$, the natural inclusion functor $A-Mod \longrightarrow A-Mod_{s}$
from $A$-modules to simplicial $A$-modules does not preserve
finitely presented objects in general.

The category $sk-Mod$ is also Quillen equivalent (actually equivalent) to the model category
$C^{-}(k)$ of non-positively graded cochain complexes of $k$-modules, through the Dold-Kan correspondence (\cite[8.4.1]{wei}). In particular, the
suspension functor
$$S : \mathrm{Ho}(sk-Mod) \longrightarrow \mathrm{Ho}(sk-Mod)$$
corresponds to the shift functor $E \mapsto E[1]$ on the level of complexes, and is
a fully faithful functor. This implies that for any
$A\in sk-Alg$, the suspension functor
$$S : \mathrm{Ho}(A-Mod_{s}) \longrightarrow \mathrm{Ho}(A-Mod_{s})$$
is also fully faithful. We have furthermore
$\pi_{i}(S(M))\simeq \pi_{i+1}(M)$ for all $M\in A-Mod_{s}$.
The suspension and loop functors will be denoted respectively by
$$M[1]:=S(M) \qquad M[-1]:=\Omega(M).$$

For any $A \in sk-Alg$, we can construct a functorial tower in $sk-Alg$,
called the \emph{Postnikov tower},
$$\xymatrix{
A \ar[r] & \dots \ar[r] & A_{\leq n} \ar[r] & A_{\leq n-1} \ar[r] & \dots \ar[r]
& A_{\leq 0}=\pi_{0}(A)}$$
in such a way that $\pi_{i}(A_{\leq n})=0$ for all $i>n$, and the morphism
$A \longrightarrow A_{\leq n}$ induces isomorphisms on the $\pi_{i}$'s for all 
$i\leq n$. The morphism $A\longrightarrow A_{\leq n}$ is characterized by the
fact that for any $B\in sk-Alg$ which is $n$-truncated (i.e. $\pi_{i}(B)=0$
for all $i>n$), the induced morphism
$$Map_{sk-Alg}(A_{\leq n},B) \longrightarrow Map_{sk-Alg}(A,B)$$
is an isomorphism in $\mathrm{Ho}(SSet)$. This implies in particular that
the Postnikov tower is furthermore unique up to equivalence (i.e. unique as an object
in the homotopy category of diagrams). There exists a natural isomorphism in $\mathrm{Ho}(sk-Alg)$
$$A\simeq Holim_{n}A_{\leq n}.$$
For any integer $n$, the homotopy fiber of the morphism
$$A_{\leq n} \longrightarrow A_{\leq n-1}$$
is isomorphic in $\mathrm{Ho}(sk-Mod)$ to $S^{n}\otimes_{k}\pi_{n}(A)$, and is also denoted by
$\pi_{n}(A)[n]$. The $k$-module $\pi_{n}(A)$ has a natural
structure of a $\pi_{0}(A)$-module, and this induces a natural structure of
a simplicial $\pi_{0}(A)$-module on each $\pi_{n}(A)[i]$ for all $i$. Using the natural projection
$A_{\leq n-1} \longrightarrow \pi_{0}(A)$, we thus see the object
$\pi_{n}(A)[i]$ as an object in $A_{\leq n-1}-Mod_{s}$. Note that
there is a natural isomorphism in $\mathrm{Ho}(A_{\leq n-1}-Mod_{s})$
$$S(\pi_{n}(A)[i])\simeq \pi_{n}(A)[i+1] \qquad
\Omega (\pi_{n}(A)[i])\simeq \pi_{n}(A)[i-1],$$
where $\pi_{n}(A)[i]$ is understood to be $0$ for $i<0$. We recall the
following important and well known fact.

\begin{lem}\label{lpostn}
With the above notations, there exists a unique
derivation
$$d_{n}\in \pi_{0}(\mathbb{D}er_{k}(A_{\leq n-1},\pi_{n}(A)[n+1]))$$
such that the natural projection
$$A_{\leq n-1}\oplus_{d_{n}}\pi_{n}(A)[n] \longrightarrow A_{\leq n-1}$$
is isomorphic in $\mathrm{Ho}(sk-Alg/A_{\leq n-1})$ to the natural morphism
$$A_{\leq n} \longrightarrow A_{\leq n-1}.$$
\end{lem}

\begin{proof}[Sketch of Proof] (See also \cite{ba}
for more details).

The uniqueness of $d_{n}$ follows easily from our lemma \ref{l16'}, and the fact that the natural
morphism
$$\mathbb{L}QZ(\pi_{n}(A)[n+1]) \longrightarrow \pi_{n}(A)[n+1]$$
induces an isomorphism on $\pi_{i}$ for all $i\leq n+1$
(this follows from our lemma \ref{lII-4-} below).
To prove the existence
of $d_{n}$, we consider the homotopy push-out diagram in
$sk-Alg$
$$\xymatrix{
A_{\leq n-1} \ar[r] & B \\
A_{\leq n} \ar[u] \ar[r] & A_{\leq n-1}. \ar[u]}$$
The identity of $A_{\leq n-1}$ induces a morphism
$B \longrightarrow A_{\leq n-1}$, which is a retraction of
$A_{\leq n-1} \longrightarrow B$. Taking the $(n+1)$-truncation gives a commutative diagram
$$\xymatrix{
A_{\leq n-1} \ar[r] & B_{\leq n+1} \\
A_{\leq n} \ar[u] \ar[r] & A_{\leq n-1}. \ar[u]_-{s}}$$
in such a way that $s$ has a retraction. This easily implies that
the morphism $s$ is isomorphic, in a non-canonical way, to the zero derivation
$A_{\leq n-1} \longrightarrow A_{\leq n-1}\oplus \pi_{n}(A)[n+1]$. The top horizontal
morphism of the previous diagram then gives rise to a derivation
$$d_{n} : A_{\leq n-1} \longrightarrow A_{\leq n-1}\oplus \pi_{n}(A)[n+1].$$
The diagram
$$\xymatrix{
A_{\leq n-1} \ar[r]^{d_{n}} & B_{\leq n+1} \\
A_{\leq n} \ar[u] \ar[r] & A_{\leq n-1}, \ar[u]_-{s}}$$
is then easily checked to be homotopy cartesian, showing that
$A_{\leq n} \longrightarrow A_{\leq n-1}$ is isomorphic to
$A_{\leq n-1}\oplus_{d_{n}}\pi_{n}(A)[n] \longrightarrow A_{\leq n-1}$.
\end{proof}

Finally, the truncation construction also exists for modules. For
any $A\in sk-Alg$, and $M\in A-Mod_{s}$, there exists a natural
tower of morphisms in $A-Mod_{s}$
$$\xymatrix{
M \ar[r] & \dots \ar[r] & M_{\leq n} \ar[r] & M_{\leq n-1} \ar[r] & \dots \ar[r]
& M_{\leq 0}=\pi_{0}(M),}$$
such a way that $\pi_{i}(M_{\leq n})=0$ for all $i>n$, and the morphism
$M \longrightarrow M_{\leq n}$ induces an isomorphisms on $\pi_{i}$ for
$i\leq n$. The natural morphism
$$M \longrightarrow Holim_{n}M_{\leq n}$$
is an isomorphism in $\mathrm{Ho}(A-Mod_{s})$. Furthermore, the $A$-module
$M_{\leq n}$ is induced by a natural $A_{\leq n}$-module, still
denoted by $M_{\leq n}$, through the natural morphism
$A \longrightarrow A_{\leq n}$. The natural projection
$M\longrightarrow M_{\leq n}$ is again characterized by the fact that for
any $A$-module $N$ which is $n$-truncated, the induced morphism
$$Map_{A-Mod_{s}}(M_{\leq n},N) \longrightarrow Map_{A-Mod_{s}}(M,N)$$
is an isomorphism in $\mathrm{Ho}(SSet)$. 

For an object $A\in sk-Alg$, the homotopy category
$\mathrm{Ho}(Sp(A-Mod_{s}))$, of stable $A$-modules can be described in
the following way. By normalization, the commutative
simplicial $k$-algebra $A$ can be transformed into
a commutative $dg$-algebra over $k$, $N(A)$ (because $N$ is lax symmetric monoidal). We can therefore
consider its model category of unbounded $N(A)$-dg-modules, and its
homotopy category $\mathrm{Ho}(N(A)-Mod)$. The two categories
$\mathrm{Ho}(Sp(A-Mod_{s}))$ and $\mathrm{Ho}(N(A)-Mod)$ are then
naturally equivalent. In particular, when $A$ is a commutative $k$-algebra, then
$N(A)=A$, and one finds that $\mathrm{Ho}(Sp(A-Mod_{s}))$ is simply the
unbounded derived category of $A$, or equivalently the homotopy category
of the model category $C(A)$ of unbounded complexes of $A$-modules 
$$\mathrm{Ho}(Sp(A-Mod_{s}))\simeq D(A)\simeq \mathrm{Ho}(C(A)).$$
Finally, using our Cor. \ref{cI-cell} (see also \cite[III.7]{ekmm}), we see that the
perfect objects in the symmetric monoidal model category $\mathrm{Ho}(Sp(A-Mod_{s}))$
are exactly the finitely presented objects. \\

We now let $k-D^{-}Aff$ be the opposite model category
of $sk-Alg$. We use our general notations, $Spec\, A\in k-D^{-}Aff$ being the
object corresponding to $A\in sk-Alg$. We will also sometimes use the notation
$$t_{0}(Spec\, A):=Spec\, \pi_{0}(A).$$

\section{Flat, smooth, \'etale and Zariski open morphisms}\label{IIder.2}

According to our general definitions presented in \S \ref{partI.2} 
we have various notions of projective and perfect modules,
flat, smooth, \'etale, unramified \dots morphisms in $sk-Alg$. Our first task, before
visiting our general notions of stacks, will be
to give concrete descriptions of these notions. \\

Any object $A \in sk-Alg$ gives rise to a commutative graded
$k$-algebra of homotopy $\pi_{*}(A)$, which is functorial in $A$. In particular,
$\pi_{i}(A)$ is always endowed with a natural structure of a
$\pi_{0}(A)$-module, functorially in $A$. For a morphism $A\longrightarrow B$ in
$sk-Alg$ we obtain a natural morphism
$$\pi_{*}(A)\otimes_{\pi_{0}(A)}\pi_{0}(B) \longrightarrow \pi_{*}(B).$$
More generally, for $A \in sk-Alg$ and $M$ an $A$-module, one has a natural morphism
of $\pi_{0}(A)$-modules $\pi_{0}(M) \longrightarrow \pi_{*}(M)$, giving rise to a natural morphism
$$\pi_{*}(A)\otimes_{\pi_{0}(A)}\pi_{0}(M) \longrightarrow \pi_{*}(M).$$
These two morphisms are the same when the commutative $A$-algebra $B$ is considered as an
$A$-module in the usual way.

\begin{df}\label{dII-3}
Let $A \in sk-Alg$ and $M$ be an $A$-module. The $A$-module $M$ is
\emph{strong}\index{strong! simplicial module} if the natural morphism
$$\pi_{*}(A)\otimes_{\pi_{0}(A)}\pi_{0}(M) \longrightarrow \pi_{*}(M)$$
is an isomorphism.
\end{df}

\begin{lem}\label{lII-3}
Let $A \in sk-Alg$ and $M$ be an $A$-module.
\begin{enumerate}
\item The $A$-module $M$ is projective if and only if
it is strong and $\pi_{0}(M)$ is a projective
$\pi_{0}(A)$-module.

\item The $A$-module $M$ is flat if and only if it is
strong and $\pi_{0}(M)$ is a flat
$\pi_{0}(A)$-module.

\item The $A$-module $M$ is perfect if and only if
it is strong and $\pi_{0}(M)$ is a projective
$\pi_{0}(A)$-module of finite type.

\item The $A$-module $M$ is projective and finitely presented if and only
if it is perfect.

\end{enumerate}
\end{lem}

\begin{proof} $(1)$ Let us suppose that $M$ is projective.We first notice that a
retract of a strong module $A$-module is again a strong $A$-module. This allows us to
suppose that $M$ is free, which clearly implies that $M$ is strong and that
$\pi_{0}(M)$ is a free $\pi_{0}(A)$-module (so in particular projective). Conversely,
let $M$ be a strong $A$-module with $\pi_{0}(M)$ projective over
$\pi_{0}(A)$. We write $\pi_{0}(M)$ as a retract of a free $\pi_{0}$-module
$$\xymatrix{\pi_{0}(M) \ar[r]^-{i}& \pi_{0}(A)^{(I)}=\oplus_{I}\pi_{0}(A) \ar[r]^-{r} & \pi_{0}(M).}$$
The morphism $r$ is given by a family of elements $r_{i} \in \pi_{0}(M)$ for
$i\in I$, and therefore can be seen as a morphism $r' : A^{(I)} \longrightarrow
M$, well defined in $\mathrm{Ho}(A-Mod_{s})$. In the same way, the projector
$p=i\circ r$ of $\pi_{0}(A)^{(I)}$, can be seen as a projector
$p'$ of $A^{(I)}$ in the homotopy category $\mathrm{Ho}(A-Mod_{s})$. By construction, this
projector gives rise to a split fibration sequence
$$K \longrightarrow A^{(I)} \longrightarrow C,$$
and the morphism $r'$ induces a well defined morphism
in $\mathrm{Ho}(A-Mod_{s})$
$$r': C \longrightarrow M.$$
By construction, this morphism induces an isomorphisms on $\pi_{0}$, and as
$C$ and $M$ are strong modules, $r'$ is an isomorphism in $\mathrm{Ho}(A-Mod_{s})$. \\

$(2)$ Let $M$ be a strong $A$-module with $\pi_{0}(M)$ flat over
$\pi_{0}(A)$, and $N$ be any $A$-module. Clearly,
$\pi_{*}(M)$ is flat as a $\pi_{*}(A)$-module. Therefore, the Tor spectral
sequence of \cite{q}
$$Tor^{*}_{\pi_{*}(A)}(\pi_{*}(M),\pi_{*}(N)) \Rightarrow
\pi_{*}(M\otimes^{\mathbb{L}}_{A}N)$$
degenerates and gives a natural isomorphism
$$\pi_{*}(M\otimes^{\mathbb{L}}_{A}N) \simeq
\pi_{*}(M)\otimes_{\pi_{*}(A)}\pi_{*}(N)\simeq $$
$$\simeq (\pi_{0}(M)
\otimes_{\pi_{0}(A)}\pi_{*}(A))\otimes_{\pi_{*}(A)}\pi_{*}(N)\simeq
\pi_{0}(M)\otimes_{\pi_{0}(A)}\pi_{*}(N).$$
As $\pi_{0}(M)\otimes_{\pi_{0}(A)}-$ is an exact functor its transform
long exact sequences into long exact sequences. This easily implies that
$M\otimes^{\mathbb{L}}_{A}-$ preserves homotopy fiber
sequences, and therefore that
$M$ is a flat $A$-module.

Conversely, suppose that $M$ is a flat
$A$-module. Any short exact sequence $0\rightarrow N \rightarrow P$
of $\pi_{0}(A)$-modules can also be seen as a homotopy fiber sequence
of $A$-modules, as any morphism $N\rightarrow P$ is always a fibration.
Therefore, we obtain a homotopy fiber sequence
$$0 \longrightarrow M\otimes^{\mathbb{L}}_{A}N \longrightarrow M\otimes^{\mathbb{L}}_{A}P$$
which on $\pi_{0}$ gives a short exact sequence
$$0 \longrightarrow \pi_{0}(M)\otimes_{\pi_{0}(A)}N \longrightarrow
\pi_{0}(M)\otimes_{\pi_{0}(A)}P.$$
This shows that $\pi_{0}(M)\otimes_{\pi_{0}(A)}-$ is an exact functor, and therefore
that $\pi_{0}(M)$ is a flat $\pi_{0}(A)$-module.
Furthermore, taking $N=0$ we get that for any
$\pi_{0}(A)$-module $P$ one has
$\pi_{i}(M\otimes^{\mathbb{L}}_{A}P)=0$ for any $i>0$. In other words, we have
an isomorphism in $\mathrm{Ho}(A-Mod_{s})$
$$M\otimes^{\mathbb{L}}_{A}P\simeq \pi_{0}(M)\otimes_{\pi_{0}(A)}P.$$
By shifting $P$ we obtain that for any $i\geq 0$ and any $\pi_{0}(A)$-module
$P$ we have
$$M\otimes^{\mathbb{L}}_{A}(P[i])\simeq (\pi_{0}(M)\otimes_{\pi_{0}(A)}P)[i].$$
Passing to Postnikov towers we see that this implies that for any
$A$-module $P$ we have
$$\pi_{i}(M\otimes^{\mathbb{L}}_{A}P)\simeq
\pi_{0}(M)\otimes_{\pi_{0}(A)}\pi_{i}(P).$$
Applying this to $P=A$ we find that
$M$ is a strong $A$-module. \\

$(3)$ Let us first suppose that $M$ is strong with $\pi_{0}(M)$ projective of finite
type over $\pi_{0}(A)$. By point $(1)$ we know that $M$ is
a projective $A$-module. Moreover, the proof of $(1)$ also shows that
$M$ is a retract in $\mathrm{Ho}(A-Mod_{s})$ of some $A^{I}$ with $I$ finite.
By our general result Prop. \ref{pd7-} this implies that $A$ is perfect.
Conversely, let $M$ be a perfect $A$-module. By $(2)$ and
Prop. \ref{pd7-} we know that $M$ is strong and that
$\pi_{0}(M)$ is flat over $\pi_{0}(A)$. Furthermore, the unit
$k$ of $sk-Mod$ is finitely presented, so
by Prop. \ref{pd7-} $M$ is  finitely presented object
in $A-Mod_{s}$. Using the left Quillen functor $\pi_{0}$ from $A-Mod_{s}$ to
$\pi_{0}(A)$-modules, we see that this implies that $\pi_{0}(M)$ is a finitely presented
$\pi_{0}(A)$-module, and therefore is projective of finite type.\\

$(4)$ Follows from $(1)$, $(3)$ and Prop. \ref{pd7-}. \end{proof}

\begin{df}\label{dII-4}
\begin{enumerate}
\item A morphism $A \longrightarrow B$ in $sk-Alg$ is
\emph{strong}\index{strong!morphism between simplicial algebras} if the natural morphism
$$\pi_{*}(A)\otimes_{\pi_{0}(A)}\pi_{0}(B) \longrightarrow \pi_{*}(B)$$
is an isomorphism (i.e. $B$ is strong as an $A$-module).
\item
A morphism $A \longrightarrow B$ in $sk-Alg$ is
\emph{strongly flat}\index{strongly flat!morphism between simplicial algebras} (resp. \emph{strongly smooth}, resp.
\emph{strongly \'etale}\index{strongly \'etale!morphism between simplicial algebras}, resp.
\emph{a strong Zariski open immersion})\index{strong Zariski open immersion!between simplicial algebras} if it is strong and
if the morphism of affine schemes
$$Spec\, \pi_{0}(B) \longrightarrow Spec\, \pi_{0}(A)$$
is flat (resp. smooth, resp. \'etale, resp. a Zariski open immersion).

\end{enumerate}
\end{df}

We start by a very useful criterion in order to
recognize finitely presented morphisms. We have learned this
proposition from J. Lurie (see \cite{lu}).

\begin{prop}\label{plocfp}
Let $f : A \longrightarrow B$ be a morphism in $sk-Alg$. Then,
$f$ is finitely presented if and only if it satisfies the following
two conditions.
\begin{enumerate}
\item The morphism $\pi_{0}(A) \longrightarrow \pi_{0}(B)$
is a finitely presented morphism of commutative rings.
\item The cotangent complex $\mathbb{L}_{B/A} \in \mathrm{Ho}(B-Mod_{s})$ is
finitely presented.
\end{enumerate}
\end{prop}

\begin{proof} Let us assume first that $f$ is finitely presented.
Then $(1)$ and $(2)$ are easily seen to be true by fact that
$\pi_{0}$ is left adjoint and by definition
of derivations.

Let us now assume that $f : A \longrightarrow B$ is a morphism in $sk-Alg$ such that
$(1)$ and $(2)$ are satisfied.
Let $k\geq 0$ be an integer, and let $P(k)$ be the following property:
for any filtered diagram $C_{i}$ in $A-Alg_{s}$, such that
$\pi_{j}(C_{i})=0$ for all $j>k$, the natural morphism
$$Hocolim_{i}Map_{A-Alg_{s}}(B,C_{i}) \longrightarrow
Map_{A-Alg_{s}}(B,Hocolim_{i}C_{i})$$
is an isomorphism in $\mathrm{Ho}(SSet)$.

We start to prove by induction on $k$ that $P(k)$ holds for all $k$. For $k=0$ this
is hypothesis $(1)$. Suppose $P(k-1)$ holds, and let
$C_{i}$ be a any filtered diagram in $A-Alg_{s}$, such that
$\pi_{j}(C_{i})=0$ for all $j>k$. We consider $C=Hocolim_{i}C_{i}$, as well as the
$k$-th Postnikov towers
$$C \longrightarrow C_{\leq k-1} \qquad
(C_{i}) \longrightarrow (C_{i})_{\leq k-1}.$$
There is a commutative square of simplicial sets
$$\xymatrix{
Hocolim_{i}Map_{A-Alg_{s}}(B,C_{i}) \ar[d] \ar[r] & Hocolim_{i}Map_{A-Alg_{s}}(B,(C_{i})_{k-1}) \ar[d]\\
Map_{A-Alg_{s}}(B,C) \ar[r] & Map_{A-Alg_{s}}(B,C_{k-1}).}$$
By induction, the morphism on the right is an equivalence. Furthermore, using
Prop. \ref{p22bis}, Lem. \ref{lpostn} and the fact that the cotangent complex $\mathbb{L}_{B/A}$
is finitely presented, we see that the morphism induced on the
homotopy fibers of the horizontal morphisms is also an equivalence. By the five
lemma this implies that the morphism
$$Hocolim_{i}Map_{A-Alg_{s}}(B,C_{i}) \longrightarrow Map_{A-Alg_{s}}(B,C)$$
is an equivalence. This shows that $P(k)$ is satisfied.

As $\mathbb{L}_{B/A}$ is finitely presented, it is a retract of
a finite cell $B$-module (see Prop. \ref{pI-cell}). In particular, there is
an integer $k_{0}>0$, such that $[\mathbb{L}_{B/A},M]_{B-Mod_{s}}=0$ for any
$B$-module $M$ such that $\pi_{i}(M)=0$ for all $i<k_{0}$ (one can chose
$k_{0}$ strictly bigger than the dimension of the cells of a module
of which $\mathbb{L}_{B/A}$ is a retract).
Once again, Prop. \ref{p22bis} and Lem. \ref{lpostn} implies that
for any commutative $A$-algebra $C$, the natural projection
$C\longrightarrow C_{\leq k_{0}}$ induces a bijection
$$[B,C]_{A-Alg_{s}}\simeq [B,C_{\leq k_{0}}]_{A-Alg_{s}}.$$
Therefore, as the property $P(k_{0})$ holds, we find
that for any filtered system $C_{i}$ in $A-Alg_{s}$, the natural morphism
$$Hocolim_{i}Map_{A-Alg_{s}}(B,C_{i}) \longrightarrow
Map_{A-Alg_{s}}(B,Hocolim_{i}C_{i})$$
induces an isomorphism in $\pi_{0}$. As this is valid for any
filtered system, this shows that the morphism
$$Hocolim_{i}Map_{A-Alg_{s}}(B,C_{i}) \longrightarrow
Map_{A-Alg_{s}}(B,Hocolim_{i}C_{i})$$
induces an isomorphism on all the $\pi_{i}$'s, and therefore
is an equivalence. This shows that $f$ is finitely presented.
\end{proof}

\begin{prop}\label{pII-3}
Let $f : A \longrightarrow B$ be a morphism in $sk-Alg$.
\begin{enumerate}
\item The morphism $f$ is smooth if and only if
it is perfect. The morphism $f$ is formally smooth if it is
formally perfect. The morphism $f$ is formally
i-smooth if and only if it is formally smooth.

\item The morphism $f$ is (formally) unramified if and only if
it is (formally) \'etale.

\item The morphism $f$ is (formally) \'etale if and only if
it is (formally) thh-\'etale.

\item
A morphism $A \longrightarrow B$ in $sk-Alg$ is
\emph{flat} (resp.
\emph{a Zariski open immersion}) if
and only if it is \emph{strongly flat} (resp.
\emph{a strong Zariski open immersion}).

\end{enumerate}
\end{prop}

\begin{proof} $(1)$ The first two assumptions follow from Lem. \ref{lII-3} $(4)$.
For the comparison between formally i-smooth and formally smooth morphism
we notice that a $B$-module $P\in B-Mod$ is projective if and only if
$[P,M[1]]=0$ for any $M\in B-Mod$. This and Prop. \ref{pismooth}
imply the statement (note that in our present context
$\mathcal{A}=sk-Alg$ and $\mathcal{C}_{0}=\mathcal{C}$). \\

$(2)$ Follows from the fact that the suspension functor of $sk-Mod$ is
fully faithful and from Prop. \ref{p6} $(1)$. \\

$(3)$ By Prop. \ref{p6} $(2)$ (formally) thh-\'etale morphisms are (formally)
\'etale. Conversely, we need to prove that a
formally \'etale morphism $A \longrightarrow B$ is thh-\'etale. For this, we use the
well known spectral sequence
$$\pi_{*}(Sym^{*}(\mathbb{L}_{B/A}[1]))\Rightarrow \pi_{*}(THH(B/A))$$
of \cite[8]{q}. Using this spectral sequence we see that
$B\simeq THH(B/A)$ if and only if $\mathbb{L}_{B/A}\simeq 0$. In particular
a formally \'etale morphism is always formally thh-\'etale. \\

$(4)$ For flat morphism this is Lem. \ref{lII-3} $(2)$. Let $f : A \longrightarrow B$ be
a Zariski open immersion. By definition, $f$ is a flat morphism, and therefore
is strongly flat by what we have seen. Moreover,
for any commutative $k$-algebra $C$, considered as an object
$C\in sk-Alg$ concentrated in degree $0$, we have
natural isomorphisms in $\mathrm{Ho}(SSet)$
$$Map_{sk-Alg}(A,C)\simeq
Hom_{k-Alg}(\pi_{0}(A),C)$$ 
$$Map_{sk-Alg}(B,C)\simeq
Hom_{k-Alg}(\pi_{0}(B),C).$$
In particular, as $f$ is a epimorphism in $sk-Alg$ the induced morphism of affine
schemes $\varphi: Spec\, \pi_{0}(B) \longrightarrow Spec\, \pi_{0}(A)$ is a monomorphism of schemes. This last morphism is therefore
a flat monomorphism of affine schemes. Moreover, as $f$ is finitely presented
so is $\varphi$, which is
therefore a finitely presented flat monomorphism. By \cite[2.4.6]{egaIV-4}, a finitely
presented flat morphism is open, so $\varphi$ is an open flat monomorphism and thus a Zariski open immersion.
This implies that $f$ is a strong Zariski open immersion.

Conversely, let $f : A \longrightarrow B$ be a strong Zariski open immersion.
By $(1)$ it is a flat morphism. It only remains to show that it is also
an epimorphism in $sk-Alg$ and that it is finitely presented. For the first of these properties,
we use the $Tor$ spectral sequence to see that the natural
morphism $B \longrightarrow B\otimes_{A}^{\mathbb{L}}B$ induces an isomorphism on
homotopy groups
$$\pi_{*}(B\otimes_{A}^{\mathbb{L}}B)\simeq
\pi_{*}(B)\otimes_{\pi_{*}(A)}\pi_{*}(B)\simeq (\pi_{0}(B)\otimes_{\pi_{0}(A)}\pi_{0}(B))\otimes_{\pi_{0}(A)}\pi_{*}(A)
\simeq \pi_{*}(B).$$
In other words, the natural morphism
$B\otimes^{\mathbb{L}}_{A}B \longrightarrow B$ is an isomorphism in
$\mathrm{Ho}(sk-Alg)$, which implies that $f$ is an epinomorphism.
It remain to be shown that $f$ is furthermore finitely presented, but this
follows from Prop. \ref{plocfp} and Prop. \ref{p6}.
\end{proof}

\begin{thm}\label{tII-1}
\begin{enumerate}
\item A morphism in $sk-Alg$ is \'etale if and
only if it is strongly  \'etale.

\item A morphism in $sk-Alg$ is  smooth  if and
only if it is strongly smooth.

\end{enumerate}
\end{thm}

\begin{proof} We start by two fundamental lemmas. \\

Recall from \S \ref{Ider} the Quillen adjunction
$$Q : A-Comm^{nu}(\mathcal{C}) \longrightarrow A-Mod \qquad
A-Comm^{nu}(\mathcal{C}) \longleftarrow A-Mod : Z.$$

\begin{lem}\label{lII-4-}
Let $A\in sk-Alg$, and $M\in A-Mod_{s}$ be a $A$-module
such that $\pi_{i}(M)=0$ for all $i<k$, for some fixed integer $k$.
Then, the adjunction morphism
$$\mathbb{L}QZ(M) \longrightarrow M$$
induces isomorphisms
$$\pi_{i}(\mathbb{L}QZ(M)) \simeq \pi_{i}(M)$$
for all $i\leq k$. In particular
$$\pi_{i}(\mathbb{L}QZ(M)) \simeq 0 \; for \; i<k \qquad
\pi_{k}(\mathbb{L}QZ(M))\simeq \pi_{k}(M).$$
\end{lem}

\begin{proof} The non-unital $A$-algebra $Z(M)$ being $(k-1)$-connected, we can write
it, up to an equivalence, as a CW object in $A-Comm^{nu}(sk-Mod)$
with cells of dimension at least $k$ (in the sense
of \cite{ekmm}). In other words $Z(M)$ is equivalent to
a filtered colimit
$$\xymatrix{\dots  \ar[r] & A_{i} \ar[r] & A_{i+1} \ar[r] & \dots }$$
where at each step there is a push-out diagram in $A-Comm^{nu}(sk-Mod)$
$$\xymatrix{
A_{i} \ar[r] & A_{i+1} \\
\coprod k[S^{n_{i}}]^{nu} \ar[u] \ar[r] & \coprod k[\Delta^{n_{i}+1}]^{nu}, \ar[u]}$$
where $n_{i+1}>n_{i}\geq k-1$, and $k[K]^{nu}$ is the free non-unital
commutative simplicial $k$-algebra generated by a simplicial set $K$.
Therefore, $Q$ being left Quillen, the $A$-module
$\mathbb{L}QZ(M)$ is the homotopy colimit of
$$\xymatrix{\dots  \ar[r] & Q(A_{i}) \ar[r] & Q(A_{i+1}) \ar[r] & \dots }$$
where at each step there exists a homotopy push-out diagram in $A-Mod_{s}$
$$\xymatrix{
Q(A_{i}) \ar[r] & Q(A_{i+1}) \\
\oplus Q(k[S^{n_{i}}]^{nu}) \ar[u] \ar[r] & 0. \ar[u]}$$
Computing homotopy groups using long exact sequences, we see that
the statement of the lemma can be reduced to prove that
for a free non-unital commutative simplicial $k$-algebra $k[S^{n}]^{nu}$, we have
$$\pi_{i}(Q(k[S^{n}]^{nu}))\simeq 0 \; for \; i<n \qquad
\pi_{n}(Q(k[S^{n}]^{nu}))\simeq k.$$
But, using that $Q$ is left Quillen we find
$$Q(k[S^{n}]^{nu})\simeq S^{n}\otimes k,$$
which implies the result. \end{proof}

\begin{lem}\label{lII-4}
Let $A$ be any object in $sk-Alg$ and let us consider the $k$-th stage of its
Postnikov tower
$$A_{\leq k}\longrightarrow A_{\leq k-1}.$$
There exist natural isomorphisms
$$\pi_{k+1}(\mathbb{L}_{A_{\leq k-1}/A_{\leq k}})\simeq
\pi_{k}(A)$$
$$\pi_{i}(\mathbb{L}_{A_{\leq k-1}/A_{\leq k}})\simeq 0 \; for \;
i\leq k.$$
\end{lem}

\begin{proof} This follows from lemma \ref{lpostn},
lemma \ref{l16'}, and  lemma \ref{lII-4-}. \end{proof}

Let us now start the proof of theorem \ref{tII-1}. \\

$(1)$ Let $f : A\longrightarrow B$ be a strongly \'etale
morphism. By definition of strongly \'etale the morphism $f$ is
strongly flat. In particular, the square
$$\xymatrix{
A \ar[r] \ar[d] & B \ar[d] \\
\pi_{0}(A) \ar[r] & \pi_{0}(B)}$$
is homotopy cocartesian in $sk-Alg$. Therefore, Prop. \ref{p2} $(2)$ implies
$$\mathbb{L}_{B/A}\otimes_{B}\pi_{0}(B)\simeq
\mathbb{L}_{\pi_{0}(B)/\pi_{0}(A)}\simeq 0.$$
Using the Tor spectral sequence
$$Tor^{*}_{\pi_{*}(B)}(\pi_{0}(B),\pi_{*}(\mathbb{L}_{B/A})) \Rightarrow
\pi_{*}(\mathbb{L}_{B/A}\otimes_{B}\pi_{0}(B))=0$$
we find by induction on $k$ that
$\pi_{k}(\mathbb{L}_{B/A})=0$. This shows that the morphism
$f$ is formally \'etale. The fact that it is also finitely presented
follows then from
Prop. \ref{plocfp}.

Conversely, let $f : A\longrightarrow B$ be an \'etale morphism.
As $\pi_{0} : sk-Alg \longrightarrow k-Alg$ is left Quillen, we
deduce immediately that $\pi_{0}(A) \longrightarrow \pi_{0}(B)$
is an \'etale morphism of commutative rings. We will prove by
induction that all the truncations
$$f_{k} : A_{\leq k} \longrightarrow B_{\leq k}$$
are strongly \'etale and thus \'etale as well by what we have seen before. 
We 
thus assume that $f_{k-1}$ is strongly \'etale (and thus also \'etale
by what we have seen in the first part of the proof).
By adjunction, 
the truncations
$$f_{k} : A_{\leq k} \longrightarrow B_{\leq k}$$
are such that
$$(\mathbb{L}_{B_{\leq k}/A_{\leq k}})_{\leq k} \simeq (\mathbb{L}_{B/A})_{\leq k}\simeq 0.$$
Furthermore, let $M$ be any $\pi_{0}(B)$-module, and 
$d$ be any morphism in $\mathrm{Ho}(B_{\leq k}-Mod_{s})$
$$d : \mathbb{L}_{B_{\leq k}/A_{\leq k}} \longrightarrow M[k+1].$$
We consider the commutative diagram
in $\mathrm{Ho}(A_{\leq k}/sk-Alg/B_{\leq k})$
$$\xymatrix{
A_{\leq k} \ar[r] \ar[d] & B_{\leq k}\oplus_{d}M[k] \ar[d] \\
B_{\leq k} \ar[r] & B_{\leq k}.}$$
The obstruction of the existence of a morphism
$$u : B_{\leq k} \longrightarrow B_{\leq k}\oplus_{d}M[k]$$
in $\mathrm{Ho}(A_{\leq k}/sk-Alg/B_{\leq k})$ is precisely 
$d\in [\mathbb{L}_{B_{\leq k}/A_{\leq k}},M[k+1]]$. On the other hand, by 
adjunction, the existence of such a morphism $f$ is equivalent to 
the existence of a morphism
$$u' : B \longrightarrow B_{\leq k}\oplus_{d}M[k]$$
in $\mathrm{Ho}(A/sk-Alg/B_{\leq k})$. The obstruction of the existence of $u'$ 
itself is
$d'\in [\mathbb{L}_{B/A},M[k+1]]$, the image of $d$ by the natural
morphism, which vanishes as
$A \rightarrow B$ is formally \'etale. This implies that $u'$ and thus $u$ exists, 
and therefore that for any $\pi_{0}(B)$-module $M$ we have
$$[\mathbb{L}_{B_{\leq k}/A_{\leq k}},M[k+1]]=0.$$
As the object $\mathbb{L}_{B_{\leq k}/A_{\leq k}}$  is already known to be
$k$-connected, we conclude that it is furthermore
is $(k+1)$-connected. 

Now, there exists a morphism of homotopy cofiber sequences in
$sk-Mod$
$$\xymatrix{\mathbb{L}_{A_{\leq k}}\otimes_{A_{\leq k}}^{\mathbb{L}}A_{\leq k-1} \ar[r] \ar[d] &
\mathbb{L}_{A_{\leq k-1}} \ar[r] \ar[d] &    \mathbb{L}_{A_{\leq k-1}/A_{\leq k}} \ar[d] \\
\mathbb{L}_{B_{\leq k}}\otimes_{B_{\leq k}}^{\mathbb{L}}B_{\leq k-1} \ar[r] &
\mathbb{L}_{B_{\leq k-1}} \ar[r] &  \mathbb{L}_{B_{\leq k-1}/B_{\leq k}}. }$$
Base changing the first row by $A_{\leq k-1} \longrightarrow B_{\leq k-1}$ gives another
morphism of homotopy cofiber sequences in $sk-Mod$
$$\xymatrix{\mathbb{L}_{A_{\leq k}}\otimes_{A_{\leq k}}^{\mathbb{L}}B_{\leq k-1} \ar[r] \ar[d] &
\mathbb{L}_{A_{\leq k-1}}\otimes_{A_{\leq k-1}}^{\mathbb{L}}B_{\leq k-1} \ar[r] \ar[d] &    
\mathbb{L}_{A_{\leq k-1}/A_{\leq k}}\otimes_{A_{\leq k-1}}^{\mathbb{L}}B_{\leq k-1} \ar[d] \\
\mathbb{L}_{B_{\leq k}}\otimes_{B_{\leq k}}^{\mathbb{L}}B_{\leq k-1} \ar[r] &
\mathbb{L}_{B_{\leq k-1}} \ar[r] &  \mathbb{L}_{B_{\leq k-1}/B_{\leq k}}. }$$
Passing to the long exact sequences in homotopy, and using that
$f_{k-1}$ is \'etale, as well as 
the fact that $\mathbb{L}_{B_{\leq k}/A_{\leq k}}\otimes_{A_{\leq k-1}}^{\mathbb{L}}B_{\leq k-1}$ 
is $(k+1)$-connected, it is easy to see that
the vertical morphism on the right induces an isomorphism
$$\pi_{k+1}(\mathbb{L}_{A_{\leq k-1}/A_{\leq k}}\otimes^{\mathbb{L}}_{A_{\leq k-1}}B_{\leq k-1})
\simeq \pi_{k+1}(\mathbb{L}_{B_{\leq k-1}/B_{\leq k}}).$$
Therefore, Lemma \ref{lII-4} implies that the natural morphism
$$\pi_{k}(A)\otimes_{\pi_{0}(A)}\pi_{0}(B) \longrightarrow
\pi_{k}(B)$$
is an isomorphism. By induction on $k$ this shows that
$f$ is strongly \'etale. \\

$(2)$ Let $f : A \longrightarrow B$ be a strongly smooth morphism. As the morphism $f$ is
strongly flat, there is  a homotopy push-out square in $sk-Alg$
$$\xymatrix{
A \ar[r] \ar[d] & B\ar[d] \\
\pi_{0}(A) \ar[r] & \pi_{0}(B).}$$
Together with Prop. \ref{p2} $(2)$, we thus have a natural isomorphism
$$\mathbb{L}_{B/A}\otimes_{B}^{\mathbb{L}}\pi_{0}(B)\simeq \mathbb{L}_{\pi_{0}(B)/\pi_{0}(A)}.$$
As the morphism $\pi_{0}(A) \longrightarrow \pi_{0}(B)$ is smooth,
we have $\mathbb{L}_{\pi_{0}(B)/\pi_{0}(A)}\simeq \Omega^{1}_{\pi_{0}(B)/\pi_{0}(A)}[0]$.
Using that $\Omega^{1}_{\pi_{0}(B)/\pi_{0}(A)}$ is a projective $\pi_{0}(B)$-module, we see that 
the isomorphism
$$\Omega^{1}_{\pi_{0}(B)/\pi_{0}(A)}[0]\simeq \mathbb{L}_{B/A}\otimes_{B}^{\mathbb{L}}\pi_{0}(B)$$
can be lifted to a morphism
$$P \longrightarrow \mathbb{L}_{B/A}$$
in $\mathrm{Ho}(B-Mod_{s})$, where $P$ is a projective $B$-module
such that $P\otimes_{B}^{\mathbb{L}}\pi_{0}(B)\simeq \Omega^{1}_{\pi_{0}(B)/\pi_{0}(A)}[0]$. 
We let 
$K$ be the homotopy cofiber of this last morphism. By construction, 
we have 
$$K\otimes_{B}^{\mathbb{L}}\pi_{0}(B)\simeq 0,$$
which easily implies by induction on $k$ that $\pi_{k}(K)=0$. Therefore, 
the morphism 
$$P \longrightarrow \mathbb{L}_{B/A}$$
is in fact an isomorphism in $\mathrm{Ho}(B-Mod_{s})$, showing that 
$\mathbb{L}_{B/A}$ is a projective $B$-module. 
Moreover, the homotopy cofiber sequence
in $B-Mod_{s}$
$$\xymatrix{
\mathbb{L}_{A}\otimes_{A}^{\mathbb{L}}B \ar[r] & \mathbb{L}_{B} \ar[r]&
\mathbb{L}_{B/A}}$$
gives rise to a homotopy cofiber sequence
$$\xymatrix{
\mathbb{L}_{B} \ar[r] &
\mathbb{L}_{B/A} \ar[r] & S(\mathbb{L}_{A}\otimes_{A}^{\mathbb{L}}B).}$$
But, $\mathbb{L}_{B/A}$ being a retract of a free $B$-module, we see that
$[\mathbb{L}_{B/A},S(\mathbb{L}_{A}\otimes_{A}^{\mathbb{L}}B)]$ is a retract
of a product of $\pi_{0}(S(\mathbb{L}_{A}\otimes_{A}^{\mathbb{L}}B))=0$ and thus
is trivial.
This implies that the morphism
$\mathbb{L}_{B/A} \longrightarrow S(\mathbb{L}_{A}\otimes_{A}^{\mathbb{L}}B)$
is trivial in $\mathrm{Ho}(B-Mod_{s})$, and therefore that the homotopy cofiber sequence
$$\xymatrix{
\mathbb{L}_{A}\otimes_{A}^{\mathbb{L}}B \ar[r] & \mathbb{L}_{B}
\ar[r]& \mathbb{L}_{B/A},}$$ which is also a homotopy fiber
sequence, splits. In particular, the morphism
$\mathbb{L}_{A}\otimes_{A}^{\mathbb{L}}B \longrightarrow
\mathbb{L}_{B}$ has a retraction. We have thus shown that $f$ is a
formally smooth morphism. The fact that $f$ is furthermore
finitely presented follows from Prop. \ref{plocfp} and the fact
that $\mathbb{L}_{B/A}$ if finitely presented because
$\Omega_{B/A}^{1}$ is so (see Lem. \ref{lII-3}).

Conversely, let $f : A \longrightarrow B$ be a smooth morphism in
$sk-Alg$ and let us prove it is strongly smooth. First of all, using that
$\pi_{0} : sk-Alg \longrightarrow k-Alg$ is left Quillen, we see that
$\pi_{0}(A) \longrightarrow \pi_{0}(B)$ has the required
lifting property for being a formally smooth morphism. Furthermore, it is a finitely presented
morphism, so is a smooth morphism of commutative rings.
We then form the
homotopy push-out square in $sk-Alg$
$$\xymatrix{
A\ar[r] \ar[d] & B \ar[d] \\
\pi_{0}(A) \ar[r] & C.}$$
By base change, $\pi_{0}(A) \longrightarrow C$ is a smooth morphism. We will start to prove that
the natural morphism $C \longrightarrow \pi_{0}(C) \simeq  \pi_{0}(B)$
is an isomorphism. Suppose it is not, and let $i$ be the smallest
integer $i>0$ such that $\pi_{i}(C)\neq 0$. Considering the
homotopy cofiber sequence
$$\xymatrix{
\mathbb{L}_{C_{\leq i}/\pi_{0}(A)}\otimes^{\mathbb{L}}_{C_{\leq i}}\pi_{0}(C) \ar[r] &
\mathbb{L}_{\pi_{0}(C)/\pi_{0}(A)}\simeq \Omega^{1}_{\pi_{0}(C)/\pi_{0}(A)}[0]
\ar[r] & \mathbb{L}_{C_{\leq i}/\pi_{0}(C)},}$$
and using lemma \ref{lII-4}, we see that
$$\pi_{i}(\mathbb{L}_{C/\pi_{0}(A)}\otimes^{\mathbb{L}}_{C}\pi_{0}(C))\simeq
\pi_{i}(\mathbb{L}_{C_{\leq i}/\pi_{0}(A)}\otimes^{\mathbb{L}}_{C_{\leq i}}\pi_{0}(C))\simeq
\pi_{i+1}(\mathbb{L}_{C_{\leq i}/\pi_{0}(C)})\simeq \pi_{i}(C)
\neq 0.$$
But this contradicts the fact that
$\mathbb{L}_{C/\pi_{0}(A)}\otimes^{\mathbb{L}}_{C}\pi_{0}(C)$ is a projective (and thus strong by
lemma \ref{lII-3} $(1)$)
$\pi_{0}(C)$-module, and thus the fact that $\pi_{0}(A) \longrightarrow C$
is a smooth morphism. We therefore have a homotopy push-out diagram
in $sk-Alg$
$$\xymatrix{
A\ar[r] \ar[d] & B \ar[d] \\
\pi_{0}(A) \ar[r] & \pi_{0}(B).}$$
Using that the bottom horizontal morphism in flat and the Tor spectral sequence, we get
by induction that the natural morphism
$$\pi_{i}(A)\otimes_{\pi_{0}(A)}\pi_{0}(B) \longrightarrow \pi_{i}(B)$$
is an isomorphism. We thus have seen that $f$ is a strongly smooth morphism.
\end{proof}

An important corollary of theorem \ref{tII-1} is the following
topological invariance of \'etale morphisms.

\begin{cor}\label{ctII-1}
Let $A\in sk-Alg$ and $t_{0}(X)=Spec\, (\pi_{0}A)\longrightarrow X=Spec\, A$ be the natural morphism.
Then, the base change functor
$$\mathrm{Ho}(k-D^{-}Aff/X) \longrightarrow \mathrm{Ho}(k-D^{-}Aff/t_{0}(X))$$
induces an equivalence from the full subcategory of \'etale
morphism $Y \rightarrow X$ to the full subcategory of
\'etale morphisms $Y'\rightarrow t_{0}(X)$. Furthermore, this equivalence
preserves epimorphisms of stacks.
\end{cor}

\begin{proof} We consider the Postnikov tower
$$\xymatrix{A \ar[r] & \dots \ar[r] & A_{\leq k} \ar[r] & A_{\leq k-1} \ar[r] & \dots \ar[r] &
A_{\leq 0}=\pi_{0}(A),}$$
and the associated diagram in $k-D^{-}Aff$
$$\xymatrix{X_{\leq 0}=t_{0}(X) \ar[r] & X_{\leq 1} \ar[r] &  \dots \ar[r] &
X_{\leq k-1} \ar[r] & X_{\leq k}\ar[r] & \dots \ar[r] & X.}$$
We define a model category $k-D^{-}Aff/X_{\leq *}$, whose objects are families of
objects $Y_{k} \longrightarrow X_{\leq k}$ in $k-D^{-}Aff/X_{\leq k}$, together with
transitions morphisms
$Y_{k-1} \longrightarrow Y_{k}\times_{X_{\leq k}}X_{\leq k-1}$
in $k-D^{-}Aff/X_{\leq k-1}$. The morphisms in $k-D^{-}Aff/X_{\leq *}$
are simply the families of morphisms $Y_{k} \longrightarrow Y_{k}'$ in
$k-D^{-}Aff/X_{\leq k}$ that commute with the transition morphisms.
The model structure on $k-D^{-}Aff/X_{\leq *}$ is such that
a morphism $f : Y_{*} \longrightarrow Y_{*}'$ in $k-D^{-}Aff/X_{\leq *}$ is
an equivalence (resp. a fibration) if all morphisms $Y_{k} \longrightarrow Y_{k}'$
are equivalences (resp. fibrations) in $k-D^{-}Aff$.

There exists a Quillen adjunction
$$G : k-D^{-}Aff/X_{\leq *} \longrightarrow
k-D^{-}Aff/X \qquad
k-D^{-}Aff/X_{\leq *} \longleftarrow
k-D^{-}Aff/X : F,$$
defined by $G(Y_{*})$ to be
the colimit in $k-D^{-}Aff/X$ of the diagram
$$\xymatrix{Y_{0} \ar[r] & Y_{1} \ar[r] & \dots \ar[r] Y_{k} \ar[r] & \dots \ar[r] & X.}$$
The right adjoint $F$ sends a morphism $Y\longrightarrow X$ to the various
pullbacks
$$F(Y)_{k}:=Y\times_{X}X_{\leq k} \longrightarrow X_{k}$$
together with the obvious transition isomorphisms. Using
Thm. \ref{tII-1} $(1)$ it is not hard to check that the derived adjunction
$$\mathbb{L}G : \mathrm{Ho}(k-D^{-}Aff/X_{\leq *}) \longrightarrow
\mathrm{Ho}(k-D^{-}Aff/X)$$
$$\mathrm{Ho}(k-D^{-}Aff/X_{\leq *}) \longleftarrow
\mathrm{Ho}(k-D^{-}Aff/X) : \mathbb{R}F,$$
induces an equivalence from the full subcategory of $\mathrm{Ho}(k-D^{-}Aff/X)$
consisting of \'etale morphisms $Y \longrightarrow X$, and the full subcategory of
$\mathrm{Ho}(k-D^{-}Aff/X_{\leq *})$ consisting of objects $Y_{*}$ such that each
$Y_{k} \longrightarrow X_{\leq k}$ is \'etale and each transition morphism
$$Y_{k-1} \longrightarrow Y_{k}\times^{h}_{X_{\leq k}}X_{\leq k-1}$$
is an isomorphism in $\mathrm{Ho}(k-D^{-}Aff/X_{\leq k-1})$. Let us denote these
two categories respectively by
$\mathrm{Ho}(k-D^{-}Aff/X)^{et}$ and $\mathrm{Ho}(k-D^{-}Aff/X_{\leq *})^{cart,et}$.
Using Cor. \ref{ct12}
and Lem. \ref{lpostn}, we know that
each base change functor
$$\mathrm{Ho}(k-D^{-}Aff/X_{\leq k}) \longrightarrow \mathrm{Ho}(k-D^{-}Aff/X_{\leq k-1})$$
induces an equivalence on the full sub-categories of \'etale morphisms. This easily implies that
the natural projection functor
$$\mathrm{Ho}(k-D^{-}Aff/X_{\leq *})^{cart,et} \longrightarrow \mathrm{Ho}(k-D^{-}Aff/X_{\leq 0})^{et}$$
is an equivalence of categories. Therefore, by composition, we find that the base change
functor
$$\mathrm{Ho}(k-D^{-}Aff/X)^{et} \longrightarrow \mathrm{Ho}(k-D^{-}Aff/X_{\leq *})^{cart,et} \longrightarrow
\mathrm{Ho}(k-D^{-}Aff/X_{\leq 0})^{et}$$
is an equivalence of categories.

Finally, the statement concerning epimorphism of stacks is obvious, as
a flat morphism $Y \longrightarrow X$ in $k-D^{-}Aff$ induces
an epimorphism of stacks if and only if $t_{0}(Y) \longrightarrow t_{0}(X)$
is a surjective morphism of affine schemes. \end{proof}

A direct specialization of \ref{ctII-1} is the following.

\begin{cor}\label{ctII-1zar}
Let $A\in sk-Alg$ and $t_{0}(X)=Spec\, (\pi_{0}A)\longrightarrow X=Spec\, A$ be the natural morphism.
Then, the base change functor
$$\mathrm{Ho}(k-D^{-}Aff/X) \longrightarrow \mathrm{Ho}(k-D^{-}Aff/t_{0}(X))$$
induces an equivalence from full sub-categories of Zariski open
immersions $Y \rightarrow X$ to the full subcategory of
Zariski open immersions $Y'\rightarrow t_{0}(X)$. Furthermore, this equivalence
preserves epimorphisms of stacks.
\end{cor}

\begin{proof} Indeed, using \ref{ctII-1} it is enough to see
that an \'etale morphism $A \longrightarrow B$ is a Zariski open
immersion if and only if $\pi_{0}(A) \longrightarrow \pi_{0}(B)$
is. But this true by Thm. \ref{tII-1} and Prop. \ref{pII-3}. \end{proof}

From the proof of Thm. \ref{tII-1} we also extract the following
more precise result.

\begin{cor}\label{ctII-1'}
Let $f : A \longrightarrow B$ be a morphism in $sk-Alg$.
The following are equivalent.
\begin{enumerate}
\item The morphism $f$ is smooth (resp. \'etale).
\item The morphism $f$ is flat and
$\pi_{0}(A) \longrightarrow \pi_{0}(B)$ is smooth (resp. \'etale).
\item The morphism $f$ is formally smooth (resp.
formally \'etale) and $\pi_{0}(B)$ is a
finitely presented $\pi_{0}(A)$-algebra.
\end{enumerate}
\end{cor}

We are now ready to define the \'etale model topology (Def. \ref{modtop}) on $k-D^{-}Aff$.

\begin{df}\label{dII-5}
A family of morphisms $\{Spec\, A_{i} \longrightarrow Spec\, A\}_{i\in I}$
in $k-D^{-}Aff$ is an \emph{\'etale covering family}
(or simply \emph{\'et-covering family})\index{\'et-covering family} if it satisfies the following two
conditions.
\begin{enumerate}
\item Each morphism $A \longrightarrow A_{i}$
is \'etale.
\item There exists a finite sub-set $J\subset I$ such that the
family $\{A \longrightarrow A_{i}\}_{i\in J}$ is a
formal covering family in the sense of \ref{dcov}.
\end{enumerate}
\end{df}

Using that \'etale morphisms are precisely the strongly \'etale morphisms (see Corollary \ref{ctII-1'}) we immediately deduce
that a family of morphisms $\{Spec\, A_{i} \longrightarrow Spec\, A\}_{i\in I}$
in $k-D^{-}Aff$ is an \'et-covering family if and only if there exists a finite
sub-set $J\subset I$ satisfying the following two conditions.

\begin{itemize}

\item For all $i\in I$, the natural morphism
$$\pi_{*}(A)\otimes_{\pi_{0}(A)}\pi_{0}(A_{i}) \longrightarrow
\pi_{*}(A_{i})$$
is an isomorphism.

\item The morphism of affine schemes
$$\coprod_{i\in J} Spec\, \pi_{0}(A_{i}) \longrightarrow Spec\, \pi_{0}(A)$$
is \'etale and surjective.

\end{itemize}

\begin{lem}\label{lII-5}
The \'et-covering families define a model topology (Def. \ref{modtop}) on $k-D^{-}Aff$,
which satisfies assumption \ref{ass5}.
\end{lem}

\begin{proof} That \'et-covering families defines a model topology simply follows
from the general properties of \'etale morphisms and formal coverings described in
propositions \ref{pcov} and \ref{p4'}. It only remain to show that
the \'etale topology satisfies assumption \ref{ass5}.

The \'etale topology is quasi-compact by definition, so $(1)$ of \ref{ass5}
is satisfied. In the same way, property $(2)$ of \ref{ass5} is obviously
satisfied according to the explicit definition of \'etale coverings
given above. Finally, let us check property $(3)$ of \ref{ass5}.
For this, let $A \longrightarrow B_{*}$ be a co-simplicial object
in $sk-Alg$, corresponding to a
representable \'etale-hypercover in $k-D^{-}Aff$ in the sense of
\ref{ass5} $(3)$. We consider the adjunction
$$B_{*}\otimes_{A}^{\mathbb{L}} - : \mathrm{Ho}(A-Mod_{s}) \longrightarrow \mathrm{Ho}(csB_{*}-Mod)$$
$$ \mathrm{Ho}(A-Mod_{s}) \longleftarrow \mathrm{Ho}(csB_{*}-Mod) : \int$$
defined in \S \ref{partI.2}. We restrict it to the full subcategory
$\mathrm{Ho}(csB_{*}-Mod)^{cart}$ of $\mathrm{Ho}(csB_{*}-Mod)$ consisting of cartesian objects
$$B_{*}\otimes_{A}^{\mathbb{L}} - : \mathrm{Ho}(A-Mod_{s}) \longrightarrow \mathrm{Ho}(csB_{*}-Mod)^{cart}$$
$$ \mathrm{Ho}(A-Mod_{s}) \longleftarrow \mathrm{Ho}(csB_{*}-Mod)^{cart} : \int$$
and we need to prove that this is an equivalence. By definition
of formal coverings, the base change functor
$$B_{*}\otimes_{A}^{\mathbb{L}} - : \mathrm{Ho}(A-Mod_{s}) \longrightarrow \mathrm{Ho}(csB_{*}-Mod)^{cart}$$
is clearly conservative, so it only remains to show that
the adjunction map $$Id \longrightarrow \int \circ (B_{*}\otimes_{A}^{\mathbb{L}} -)$$
is an isomorphism.

For this, let $E_{*} \in \mathrm{Ho}(csB_{*}-Mod)^{cart}$, and let us consider the
adjunction morphism
$$E_{1} \longrightarrow (Holim_{n}E_{n})\otimes_{A}^{\mathbb{L}}B_{1}.$$
We need to show that this morphism is an isomorphism in $\mathrm{Ho}(B_{1}-Mod_{s})$.
For this, we first use that $A \longrightarrow B_{1}$ is
a strongly flat morphism, and thus
$$\pi_{*}((Holim_{n}E_{n})\otimes_{A}^{\mathbb{L}}B_{1})\simeq
\pi_{*}((Holim_{n}E_{n}))\otimes_{\pi_{0}(A)}\pi_{0}(B_{1}).$$
We then apply the spectral sequence
$$H^{p}(Tot(\pi_{q}(E_{*})))\Rightarrow \pi_{q-p}((Holim_{n}E_{n}).$$
The object $\pi_{q}(E_{*})$ is now a co-simplicial module over
the co-simplicial commutative ring $\pi_{0}(B_{*})$, which is
furthermore cartesian
as all the coface morphisms $B_{n} \longrightarrow B_{n+1}$
are flat. The morphism $\pi_{0}(A) \longrightarrow
\pi_{0}(B_{*})$ being an \'etale, and thus faithfully flat, hypercover in the usual sense,
we find by the usual flat descent for quasi-coherent sheaves that
$H^{p}(Tot(\pi_{q}(E_{*})))\simeq 0$ for $p\neq 0$. Therefore, the above spectral sequence
degenerates and gives an isomorphism
$$\pi_{p}((Holim_{n}E_{n})\simeq Ker\left(
\pi_{p}(E_{0}) \rightrightarrows \pi_{p}(E_{1}) \right).$$
In other words, $\pi_{p}((Holim_{n}E_{n}))$ is the
$\pi_{0}(A)$-module obtained by
descent from $\pi_{p}(E_{*})$ on $\pi_{p}(B_{*})$. In particular, the natural morphism
$$\pi_{p}((Holim_{n}E_{n})\otimes_{\pi_{0}(A)}\pi_{0}(B_{1}) \longrightarrow
\pi_{p}(E_{1})$$
is an isomorphism. Putting all of this together we find that
$$E_{1} \longrightarrow (Holim_{n}E_{n})\otimes_{A}^{\mathbb{L}}B_{1}$$
is an isomorphism in $\mathrm{Ho}(B_{1}-Mod_{s})$. \end{proof}

We have now the model site (Def. \ref{modtop}) $(k-D^{-}Aff, \textrm{\'e}t)$, with the
\'etale model topology, and we make the following definition.

\begin{df}\label{dII-6}
\begin{enumerate}
\item A \emph{$D^{-}$-stack}\index{$D^{-}$-stack} is an object
$F\in k-D^{-}Aff^{\sim,\textrm{\'e}t}$ which is a stack
in the sense of Def. \ref{dstack}.
\item
The \emph{model category of $D^{-}$-stacks} is
$k-D^{-}Aff^{\sim,\textrm{\'e}t}$\index{$k-D^{-}Aff^{\sim,\textrm{\'e}t}$}, and its homotopy category will be simply denoted by $D^{-}\mathrm{St}(k)$.\index{$D^{-}\mathrm{St}(k)$} 
\end{enumerate}
\end{df}

The following result is a corollary of Proposition \ref{plocfp}. It states that 
the property of being finitely presented is local for the \'etale topology 
defined above.

\begin{cor}\label{clocfp}
Let $f : A \longrightarrow B$ be a morphism in $sk-Alg$.
\begin{enumerate}
\item If there exists an \'etale covering $B \longrightarrow B'$ such that
$A \longrightarrow B'$ is finitely presented, then $f$ is finitely presented.
\item If there exists an \'etale covering $A \longrightarrow A'$, such that
$$A' \longrightarrow A'\otimes_{A}^{\mathbb{L}}B$$
is finitely presented, then $f$ is finitely presented.
\end{enumerate}
\end{cor}

\begin{proof} This follows from proposition \ref{plocfp}. Indeed,
it suffices to prove that both conditions of proposition \ref{plocfp}
have the required local property. The first one is well known, and the
second one is a consequence of corollary \ref{ct22stable},
as finitely presented objects in $\mathrm{Ho}(Sp(B-Mod_{s}))$
are precisely the perfect objects. \end{proof}

\section{The HAG context: Geometric $D^{-}$-stacks}\label{IIder.3}

We now let \textbf{P} be the class of smooth morphisms in $sk-Alg$.

\begin{lem}\label{lII-6}
The class \textbf{P} of smooth morphisms and the \'etale model topology satisfy
assumption \ref{ass4}.
\end{lem}

\begin{proof} As \'etale morphisms are also smooth we see that 
assumption \ref{ass4} $(1)$ is satisfied.
Assumption \ref{ass4} $(2)$ is satisfied as smooth morphisms are stable
by homotopy pullbacks, compositions and equivalences. Let us prove that
smooth morphisms satisfy \ref{ass4} $(3)$ for the \'etale model topology.

Let $X \longrightarrow Y$ be a morphism in $k-D^{-}Aff$, and
$\{U_{i} \longrightarrow X\}$ a finite \'etale covering family such that
each $U_{i} \longrightarrow X$ and $U_{i} \longrightarrow Y$
is smooth. First of all, using that $\{U_{i} \longrightarrow X\}$ is a flat formal
covering, and that each morphism $U_{i} \longrightarrow Y$ is flat, we see that the
morphism $X \longrightarrow Y$ is flat. Using our results
Prop. \ref{pII-3} and Thm. \ref{tII-1} it only remain to show that the
morphism $t_{0}(X) \longrightarrow t_{0}(Y)$
is a smooth morphism between affine schemes. But, passing to $t_{0}$ we find a
smooth covering family
$\{t_{0}(U_{i}) \longrightarrow t_{0}(X)\}$, such that each morphism
$t_{0}(U_{i}) \longrightarrow t_{0}(Y)$ is smooth, and we know (see e.g.
\cite{egaIV-4}) that this implies that $t_{0}(X) \longrightarrow t_{0}(Y)$
is a smooth morphism of affine schemes.

Finally, property $(4)$ of \ref{ass4} is obvious. \end{proof}

We have verified our assumptions \ref{ass5} and \ref{ass4} for the \'etale  model topology and
\textbf{P} the class of smooth morphisms. We thus have that
$$(sk-Mod,sk-Mod,sk-Alg,\textrm{\'e}t,\mathbf{P})$$
is a HAG context in the sense of Def. \ref{dhag}.
We can therefore apply our general definitions to obtain
a notion of $n$-geometric $D^{-}$-stacks in $D^{-}\mathrm{St}(k)$, as well as the notion of
$n$-smooth morphisms. We then check that Artin's conditions of
Def. \ref{d25} are satisfied.

\begin{prop}\label{pII-4}
The \'etale model topology and the smooth morphisms satisfy Artin's conditions
relative to the HA context $(sk-Mod,sk-Mod,sk-Alg)$ in the sense of Def. \ref{d25}.
\end{prop}

\begin{proof} We will show that the class \textbf{E} of \'etale morphisms
satisfies conditions $(1)$ to $(5)$ of Def. \ref{d25}. \\

$(1)$ is clear as \textbf{P} is exactly the class of all
i-smooth morphisms.

$(2)$ and $(3)$ are clear by the choice of \textbf{E} and
$\mathcal{A}$. \\

To prove $(4)$, let $p : Y \longrightarrow X$ be a smooth
morphism in $k-D^{-}Aff$, and let us consider the smooth morphism of affine
schemes $t_{0}(p) : t_{0}(Y) \longrightarrow t_{0}(X)$. As $p$ is
a smooth epimorphism of stacks, $t_{0}(p)$ is a smooth and surjective morphism of affine schemes.
It is known (see e.g.
\cite{egaIV-4}) that there exists a \'etale covering of affine schemes
$X'_{0} \longrightarrow t_{0}(X)$
and a commutative diagram
$$\xymatrix{ & t_{0}(Y) \ar[d] \\
X'_{0} \ar[r] \ar[ru] & t_{0}(X).}$$
By Cor. \ref{ctII-1}, we know that there exists an \'etale covering
$X' \longrightarrow X$ in $k-D^{-}Aff$, inducing $X'_{0} \longrightarrow t_{0}(X)$.
Taking the homotopy pullback
$$Y\times_{X}^{h}X' \longrightarrow X'$$
we can replace $X$ by $X'$, and therefore assume that the morphism
$t_{0}(p) : t_{0}(Y) \longrightarrow t_{0}(X)$ has a section. We are going to show that this
section can be extended to a section of $p$, which will be enough to prove what we want.

We use the same trick as in the proof of Cor. \ref{ctII-1}, and
consider once again the model category
$k-D^{-}Aff/X_{\leq *}$. Using Thm. \ref{tII-1} $(2)$ we see that the
functor
$$\mathrm{Ho}(k-D^{-}Aff/X) \longrightarrow \mathrm{Ho}(k-D^{-}Aff/X_{\leq *})$$
induces an equivalences from the full subcategory
$\mathrm{Ho}(k-D^{-}Aff/X)^{sm}$ of smooth morphisms $Z \rightarrow X$ to
the full subcategory $\mathrm{Ho}(k-D^{-}Aff/X_{\leq *})^{cart,sm}$
consisting of objects $Z_{*}$ such that each $Z_{k}
\longrightarrow X_{\leq k}$ is smooth, and each morphism
$$Z_{k-1}
\longrightarrow Z_{k}\times_{Z_{\leq k}}^{h}Z_{\leq k-1}$$
is an isomorphism in $\mathrm{Ho}(k-D^{-}Aff)$. From this we deduce that the
space of sections of $p$ can be described in the following way
$$Map_{k-D^{-}Aff/X}(X,Y)\simeq Holim_{k}Map_{k-D^{-}Aff/X_{\leq k}}(X_{\leq k},Y_{k}),$$
where $Y_{k}\simeq Y_{\leq k}\simeq Y\times_{X}^{h}X_{\leq k}$.
But, as there always exists a surjection
$$\pi_{0}(Holim_{k}Map_{k-D^{-}Aff/X_{\leq k}}(X_{\leq k},Y_{k})) \longrightarrow
Lim_{k}\pi_{0}(Map_{k-D^{-}Aff/X_{\leq k}}(X_{\leq k},Y_{k})),$$
we only need to check that $Lim_{k}\pi_{0}(Map_{k-D^{-}Aff/X_{\leq k}}(X_{\leq k},Y_{k}))\neq \emptyset$.
For this, it is enough to prove that for any $k\geq 1$, the restriction map
$$\pi_{0}(Map_{k-D^{-}Aff/X_{\leq k}}(X_{\leq k},Y_{k})) \longrightarrow
\pi_{0}(Map_{k-D^{-}Aff/X_{\leq k-1}}(X_{\leq k},Y_{k-1}))$$
is surjective. In other words, we need to prove that a morphism
in $\mathrm{Ho}(k-D^{-}Aff/X_{\leq k})$
$$\xymatrix{
 X_{\leq k-1} \ar[d] \ar[r]& Y_{k}  \\
X_{\leq k}  &}$$
can be filled up to a commutative diagram in $\mathrm{Ho}(k-D^{-}Aff/X_{\leq k})$
$$\xymatrix{
 X_{\leq k-1} \ar[d] \ar[r] & Y_{k}  \\
X_{\leq k}. \ar[ru] &}$$
Using Prop. \ref{p22bis} and Lem. \ref{lpostn}, we see that the obstruction for the
existence of this lift lives in the group
$$[\mathbb{L}_{B_{k}/A_{\leq k}},\pi_{k}(A)[k+1]]_{A_{\leq k}-Mod_{s}},$$
where $Y_{k}=Spec\, B_{k}$ and $X=Spec\, A$. But, as
$Y_{k} \longrightarrow X_{\leq k}$ is a smooth morphism, the $A_{\leq k}$-module
$\mathbb{L}_{B_{k}/A_{\leq k}}$ is projective, and therefore
$$[\mathbb{L}_{B_{k}/A_{\leq k}},\pi_{k}(A)[k+1]]_{A_{\leq k}-Mod_{s}}\simeq 0,$$
by Lem. \ref{lII-3} $(5)$. \\

It remains to prove that property $(5)$ of Def. \ref{d25} is satisfied.
Let $X\longrightarrow X_{d}[\Omega M]$ be as in the statement of Def. \ref{d25} $(4)$,
with $X=Spec\, A$ and some connected $A$-module $M\in A-Mod_{s}$,
and let $U \longrightarrow X_{d}[\Omega M]$ be an \'etale morphism. We know
by Lem. \ref{l17} (or rather its proof) that
this morphism is of the form
$$X'_{d'}[\Omega M'] \longrightarrow X_{d}[\Omega M]$$
for some \'etale morphism $X'=Spec\, A'\longrightarrow X=Spec\, A$, and furthermore
$X'\longrightarrow X$ is equivalent to the homotopy pullback
$$U\times_{X_{d}[\Omega M]}X \longrightarrow X.$$
Therefore, using Cor. \ref{ctII-1'}, it is enough to prove that
$$t_{0}(X'_{d'}[\Omega M']) \longrightarrow
t_{0}(X_{d}[\Omega M])$$
is a
surjective morphism of affine schemes if and only if
$$t_{0}(X') \longrightarrow
t_{0}(X)$$
is so. But, 
$M$ being connected, $t_{0}(X_{d}[\Omega M])_{red}\simeq t_{0}(X)$, and
in the same way $t_{0}(X'_{d'}[\Omega M'])_{red}\simeq t_{0}(X')$, this gives the result, since
being a surjective morphism is topologically invariant. \end{proof}

\begin{cor}\label{cII-4}
\begin{enumerate}
\item
Any $n$-geometric $D^{-}$-stack has an obstruction theory.
In particular, any $n$-geometric $D^{-}$-stack
has a cotangent complex.
\item Any $n$-representable morphism of $D^{-}$-stacks
$f : F \longrightarrow G$ has a relative
obstruction theory. In particular,
any $n$-representable morphism of $D^{-}$-stacks has a relative
cotangent complex.
\end{enumerate}
\end{cor}

\begin{proof} Follows from Prop. \ref{pII-4} and theorem
\ref{t1}. \end{proof}

We finish by some properties of morphisms, as in Def. \ref{d21}.

\begin{lem}\label{lII-7}
Let \textbf{Q} be one the following class of morphisms in $k-D^{-}Aff$.
\begin{enumerate}
\item Flat.
\item Smooth.
\item Etale.
\item Finitely presented.
\end{enumerate}
Then, morphisms in \textbf{Q} are compatible with the \'etale topology and the class \textbf{P}
of smooth morphisms in the sense of Def. \ref{d20}.
\end{lem}

\begin{proof} Using the explicit description of flat, smooth and \'etale
morphisms given in Prop. \ref{pII-3} $(1)$, $(2)$, and $(3)$ reduce to the analog well known facts
for morphism between affine schemes. Finally, Cor. \ref{clocfp} implies that
the class of finitely presented morphisms in $k-D^{-}Aff$ is compatible
with the \'etale topology and the class \textbf{P}.  \end{proof}

Lemma \ref{lII-7} and definition \ref{d21} allows us to define the notions of
flat, smooth, \'etale and locally finitely presented morphisms of $D^{-}$-stacks, which are all
stable by equivalences, compositions and homotopy pullbacks. Recall that by
definition, a flat (resp. smooth, resp. \'etale, resp. locally finitely presented)
is always $n$-representable for some $n$. Using our general definition
Def. \ref{d23} we also have notions of
quasi-compact morphisms, finitely presented morphisms,
and monomorphisms between $D^{-}$-stacks. We also make the following definition.

\begin{df}\label{dII-7}
\begin{enumerate}
\item
A morphism of $D^{-}$-stacks is a \emph{Zariski open immersion}\index{Zariski open immersion!of $D^{-}$-stacks} if it is
a locally finitely presented flat monomorphism.
\item A morphism of $D^{-}$-stacks $F \longrightarrow G$ is a \emph{closed immersion}\index{closed immersion!of $D^{-}$-stacks}
if it is representable, and if for any $A\in sk-Alg$ and any
morphism $X=\mathbb{R}\underline{Spec}\, A \longrightarrow G$  the induced morphism of
representable $D^{-}$-stacks
$$F\times^{h}_{G}X\simeq \mathbb{R}\underline{Spec}\, B  \longrightarrow \mathbb{R}\underline{Spec}\, A$$
induces an epimorphism of rings $\pi_{0}(A) \longrightarrow \pi_{0}(B)$.
\end{enumerate}
\end{df}

\section{Truncations}\label{IIder.4}

We consider the natural inclusion functor
$$i : k-Aff \longrightarrow k-D^{-}Aff$$
right adjoint to the functor
$$\pi_{0} : k-D^{-}Aff \longrightarrow k-Aff.$$
The pair $(\pi_{0},i)$ is a Quillen adjunction (for the trivial
model structure on $k-Aff$), and as usual we will omit to mention
the inclusion functor $i$, and simply consider
commutative $k$-algebras as constant simplicial objects.
Furthermore, both functors preserve equivalences and
thus induce a Quillen adjunction on the model category of
pre-stacks (using notations of \cite[\S 4.8]{hagI})
$$i_{!} : k-Aff^{\wedge}=SPr(k-Aff) \longrightarrow k-D^{-}Aff^{\wedge} \qquad
k-Aff^{\wedge} \longleftarrow k-D^{-}Aff^{\wedge} : i^{*}.$$
The functor $i$ is furthermore continuous in the sense of
\cite[\S 4.8]{hagI}, meaning that the right derived functor
$$\mathbb{R}i^{*} : \mathrm{Ho}(k-D^{-}Aff^{\wedge}) \longrightarrow \mathrm{Ho}(k-Aff^{\wedge})$$
preserves the sub-categories of stacks. Indeed, by
Lem. \ref{lass5} $(2)$ and adjunction this follows from the fact that
$i : k-Aff \longrightarrow k-D^{-}Aff$ preserves co-products, equivalences and
\'etale hypercovers. By the general properties of left Bousfield localizations
we therefore get a Quillen adjunction on the model categories of stacks
$$i_{!} : k-Aff^{\sim,\textrm{\'e}t} \longrightarrow k-D^{-}Aff^{\sim,\textrm{\'e}t} \qquad
k-Aff^{\sim,\textrm{\'e}t} \longleftarrow k-D^{-}Aff^{\sim,\textrm{\'e}t} : i^{*}.$$

From this we get a derived adjunction on the homotopy categories of stacks
$$\mathbb{L}i_{!} : \mathrm{St}(k) \longrightarrow D^{-}\mathrm{St}(k)$$
$$\mathrm{St}(k) \longleftarrow D^{-}\mathrm{St}(k) : \mathbb{R}i^{*}.$$

\begin{lem}\label{lII-8}
The functor $\mathbb{L}i_{!}$ is fully faithful.
\end{lem}

\begin{proof} We need to show that for any $F\in \mathrm{St}(k)$ the
adjunction morphism
$$F \longrightarrow \mathbb{R}i^{*}\circ \mathbb{L}i_{!}(F)$$
is an isomorphism.
The functor $\mathbb{R}i^{*}$ commutes with
homotopy colimits, as these are computed in the model category of simplicial
presheaves and thus levelwise. Moreover, as $i_{!}$ is left Quillen
the functor $\mathbb{L}i_{!}$ also commutes with homotopy colimits. Now,
any stack $F\in \mathrm{St}(k)$ is a homotopy colimit of
representable stacks (i.e. affine schemes), and therefore we can
suppose that $F=Spec\, A$, for $A\in k-Alg$. But then
$\mathbb{L}i_{!}(Spec\, A)\simeq \mathbb{R}\underline{Spec}\, A$. Furthermore,
for any $B\in k-Alg$ there are natural isomorphisms in $\mathrm{Ho}(SSet)$
$$\mathbb{R}\underline{Spec}\, A(B)\simeq Map_{sk-Alg}(A,B)\simeq Hom_{k-Alg}(A,B)\simeq
(Spec\, A)(B).$$
This shows that the adjunction morphism
$$Spec\, A \longrightarrow \mathbb{R}i^{*}\circ \mathbb{L}i_{!}(Spec\, A)$$
is an isomorphism. \end{proof}

Another useful remark is the following 

\begin{lem}\label{lII-8'}
The functor $i^{*} : k-D^{-}Aff^{\sim,\textrm{\'e}t} \longrightarrow k-Aff^{\sim,\textrm{\'e}t}$
is right and left Quillen. In particular it preserves equivalences.
\end{lem}

\begin{proof} The functor
$i^{*}$ has a right adjoint
$$\pi_{0}^{*} : k-Aff^{\sim,\textrm{\'e}t} \longrightarrow k-D^{-}Aff^{\sim,\textrm{\'e}t}.$$
Using lemma \ref{lass5} $(2)$ we see that $\pi_{0}^{*}$
is a right Quillen functor. Therefore $i^{*}$ is left Quillen. \end{proof}

\begin{df}\label{dII-8}
\begin{enumerate}
\item
The \emph{truncation functor}\index{truncation functor} is
$$t_{0}:=i^{*} : D^{-}\mathrm{St}(k) \longrightarrow \mathrm{St}(k).$$

\item The \emph{extension functor}\index{extension functor} is the left adjoint to $t_{0}$
$$i:=\mathbb{L}i_{!} : \mathrm{St}(k) \longrightarrow D^{-}\mathrm{St}(k).$$
\item A $D^{-}$-stack $F$ is \emph{truncated} if the adjunction morphism
$$it_{0}(F) \longrightarrow F$$
is an isomorphism in $D^{-}\mathrm{St}(k)$.
\end{enumerate}
\end{df}

By lemmas \ref{lII-8}, \ref{lII-8'} we know that the truncation functor
$t_{0}$ commutes with homotopy limits and homotopy colimits. The extension functor
$i$ itself commutes with homotopy colimits and is fully faithful. An important
remark is that the extension $i$ does not commute with homotopy limits, as the inclusion functor
$k-Alg \longrightarrow sk-Alg$ does not preserve homotopy push-outs.

Concretely, the truncation functor $t_{0}$ sends a functor
$$F : sk-Alg \longrightarrow SSet_{\mathbb{V}}$$
to
$$\begin{array}{cccc}
t_{0}(F) : & k-Alg & \longrightarrow & SSet_{\mathbb{V}} \\
 & A & \mapsto & F(A).
\end{array}$$
By adjunction we clearly have
$$t_{0}\mathbb{R}\underline{Spec}\, A\simeq Spec\, \pi_{0}(A)$$
for any $A\in sk-Alg$, showing that the notation is compatible with the one
we did use before for $t_{0}(Spec\, A)=Spec\, \pi_{0}(A)$ as objects in $k-D^{-}Aff$.  
The extension functor $i$ is characterized by
$$i(Spec\, A)\simeq \mathbb{R}\underline{Spec}\, A,$$
and the fact that it commutes with homotopy colimits. 

\begin{prop}\label{pII-5}
\begin{enumerate}
\item The functor $t_{0}$ preserves epimorphisms of stacks.
\item The functor $t_{0}$ sends $n$-geometric
$D^{-}$-stacks to $n$-geometric stacks, and flat (resp. smooth, resp. \'etale)
morphisms between $D^{-}$-stacks to flat (resp. smooth, resp. \'etale) morphisms
between stacks.
\item The functor $i$ preserves homotopy pullbacks of $n$-geometric stacks
along a flat morphism, sends $n$-geometric stacks to $n$-geometric
$D^{-}$-stacks, and flat (resp. smooth, resp. \'etale) morphisms between $n$-geometric stacks to
flat (resp. smooth, resp. \'etale) morphisms between $n-$-geometric $D^{-}$-stacks.
\item Let $F \in \mathrm{St}(k)$ be an $n$-geometric stack, and
$F' \longrightarrow i(F)$ be a flat morphism of $n$-geometric $D^{-}$-stacks. Then
$F'$ is truncated (and therefore is the image by $i$ of an $n$-geometric stack by $(2)$).
\end{enumerate}
\end{prop}

\begin{proof} $(1)$ By adjunction, this follows from the fact that
$i : k-Alg \longrightarrow sk-Alg$ reflects
\'etale covering families (by Prop. \ref{lII-3}). \\

$(2)$ The proof is by induction on $n$. For $n=-1$, 
this simply follows from the formula
$$t_{0}\mathbb{R}\underline{Spec}\, A\simeq Spec\, \pi_{0}(A),$$
and Prop. \ref{pII-3} and Thm. \ref{tII-1}. Assume that the
property is known for any $m<n$ and let us prove it for $n$. Let
$F$ be an $n$-geometric $D^{-}$-stack, which by Prop. \ref{p13}
can be written as $|X_{*}|$ for some $(n-1)$-smooth Segal groupoid
$X_{*}$ in $k-D^{-}Aff^{\sim,\textrm{\'e}t}$. Using our property at rank
$n-1$ and that $t_{0}$ commutes with homotopy limits shows that
$t_{0}(X_{*})$ is also a $(n-1)$-smooth Segal groupoid object in
$k-Aff^{\sim,\textrm{\'e}t}$. Moreover, as $t_{0}$ commutes with homotopy
colimits we have $t_{0}(F)\simeq |t_{0}(X_{*})|$, which by Prop.
\ref{p13} is an 
$n$-geometric stack. This shows that $t_{0}$ sends
$n$-geometric $D^{-}$-stacks to 
$n$-geometric stacks and therefore
preserves $n$-representable morphisms. Let $f : F \longrightarrow
G$ be a flat (resp. smooth, resp. \'etale) morphism of
$D^{-}$-stacks which is $n$-representable, and let us prove by
induction on $n$ that $t_{0}(f)$ is flat (resp. smooth, resp.
\'etale). We let $X\in \mathrm{St}(k)$ be an affine scheme, and
$X \longrightarrow t_{0}(G)$ be a morphism of stacks. By adjunction
between $i$ and $t_{0}$ we have
$$X\times_{t_{0}(G)}^{h}t_{0}(F)\simeq
t_{0}(X\times_{G}^{h}F),$$
showing that we can assume that $G=X$ is an affine scheme and thus $F$ to be
an $n$-geometric $D^{-}$.
By definition of being flat (resp. smooth, resp. \'etale) there exists a smooth
$n$-atlas $\{U_{i}\}$ of $F$ such that each composite morphism
$U_{i} \longrightarrow X$ is $(n-1)$-representable and flat (resp.
smooth, resp. \'etale). By induction and $(1)$ the family $\{t_{0}(U_{i})\}$
is a smooth $n$-atlas for $t_{0}(F)$, and by induction each composition
$t_{0}(U_{i}) \longrightarrow t_{0}(X)=X$ is flat (resp. smooth, resp. \'etale). This
implies that $t_{0}(F) \longrightarrow t_{0}(X)$ is flat (resp. smooth, resp. \'etale). \\

$(3)$ The proof is by induction on $n$. For $n=-1$ 
this follows from
the formula
$i(Spec\, A)\simeq \mathbb{R}\underline{Spec}\, A$, the description
of flat, smooth and \'etale morphisms (Prop. \ref{pII-3} and Thm. \ref{tII-1}), and the fact that
for any flat morphism of commutative $k$-algebras $A \longrightarrow B$
and any commutative $A$-algebra $C$ there is a natural isomorphism in $\mathrm{Ho}(A-Alg_{s})$
$$B\otimes_{A}C\simeq B\otimes_{A}^{\mathbb{L}}C.$$
Let us now assume the property is proved for $m<n$ and let us prove it for
$n$. Let $F$ be an  $n$-geometric stack, and by Prop. \ref{p13} let us
write it as $F\simeq |X_{*}|$ for
some 
$(n-1)$-smooth Segal groupoid object $X_{*}$ in $k-Aff^{\sim,\textrm{\'e}t}$. By induction,
$i(X_{*})$ is again a 
$(n-1)$-smooth Segal groupoid objects in $k-D^{-}Aff^{\sim,\textrm{\'e}t}$, and
as $i$ commutes with homotopy colimits we have
$i(F)\simeq |i(X_{*})|$. Another application of Prop. \ref{p13} shows that
$F$ is an 
$n$-geometric $D^{-}$-stacks. We thus have seen that
$i$ sends $n$-geometric stacks to 
$n$-geometric $D^{-}$-stacks.

Now, let $F \longrightarrow G$ be a flat morphism between  $n$-geometric stacks, and
$H \longrightarrow G$ any morphism between $n$-geometric stacks. We want to show that the
natural morphism
$$i(F\times_{G}^{h}H) \longrightarrow i(F)\times^{h}_{i(G)}i(H)$$
is an isomorphism in $\mathrm{St}(k)$. For this, we write
$G\simeq |X_{*}|$ for some $(n-1)$-smooth Segal groupoid object in
$k-Aff^{\sim,\textrm{\'e}t}$, and we consider the Segal groupoid objects
$$F_{*}:=F\times_{G}^{h}X_{*}\qquad H_{*}:=H\times_{G}^{h}X_{*},$$
where $X_{*} \longrightarrow |X_{*}|=G$ is the natural augmentation in
$\mathrm{St}(k)$. The Segal groupoid objects
$F_{*}$ and $H_{*}$ are again $(n-1)$-smooth Segal groupoid objects
in $k-Aff^{\sim,\textrm{\'e}t}$ as $G$ is an  $n$-geometric stack. The natural morphisms of
Segal groupoid objects
$$F_{*} \longrightarrow X_{*} \longleftarrow H_{*}$$
gives rise to another 
$(n-1)$-smooth Segal groupoid object
$F_{*}\times^{h}_{X_{*}}H_{*}$. Clearly, we have
$$F\times^{h}_{G}H\simeq |F_{*}\times^{h}_{X_{*}}H_{*}|.$$
As $i$ commutes with homotopy colimits, and by induction on $n$ we have
$$i(F\times^{h}_{G}H)\simeq
|i(F_{*})\times^{h}_{i(X_{*})}i(H_{*})|\simeq
|i(F_{*})|\times^{h}_{|i(X_{*})|}|i(H_{*})|\simeq
i(F)\times^{h}_{i(G)}i(H).$$
It remains to prove that if $f : F \longrightarrow G$ is a flat (resp. smooth, resp.
\'etale) morphism between  $n$-geometric stacks then $i(f) : i(F) \longrightarrow i(G)$
is a flat (resp. smooth, resp. \'etale) morphism between $n$-geometric $D^{-}$-stacks.
For this, let $\{U_{i}\}$ be a smooth $n$-atlas for $G$. We have seen before that
$\{i(U_{i})\}$ is a smooth $n$-atlas for $i(G)$. As we have seen that $i$ commutes with
homotopy pullbacks along flat morphisms, and because of the local properties
of flat (resp. smooth, resp. \'etale) morphisms (see Prop. \ref{p20} and Lem. \ref{lII-7}), we can
suppose that $G$ is one of the $U_{i}$'s and thus is an affine scheme. Now, let
$\{V_{i}\}$ be a smooth $n$-atlas for $F$. The family
$\{i(V_{i})\}$ is a smooth $n$-atlas for $F$, and furthermore
each morphism $i(V_{i}) \longrightarrow i(G)$ is the image by $i$ of a flat (resp. smooth, resp.
\'etale) morphism between affine schemes and therefore is a flat (resp. smooth, resp.
\'etale morphism) of $D^{-}$-stacks. By definition this implies that
$i(F) \longrightarrow i(G)$ is flat (resp. smooth, resp. \'etale). \\

$(4)$ The proof is by induction on $n$. For 
$n=-1$ 
this is simply
the description of flat morphisms of Prop. \ref{pII-3}. Let us assume the property
is proved for $m<n$ and let us prove it for $n$. Let $F' \longrightarrow i(F)$ be
a flat morphism, with $F$ an  $n$-geometric stack and $F'$ an 
$n$-geometric $D^{-}$-stack.
Let $\{U_{i}\}$ be a smooth $n$-atlas for $F$. Then, $\{i(U_{i})\}$ is a
smooth 
$n$-atlas for $i(F)$. We consider the commutative diagram of $D^{-}$-stacks
$$\xymatrix{
it_{0}(F') \ar[r]^-{f} \ar[rd] & F' \ar[d] \\
 & i(F).}$$
We need to prove that $f$ is an isomorphism in $D^{-}\mathrm{St}(k)$. As this is
a local property on $i(F)$, we can take the homotopy pullback over the atlas
$\{i(U_{i})\}$, and thus suppose that $F$ is an affine scheme. Let now $\{V_{i}\}$ be
a smooth $n$-atlas for $F'$, and we consider the
homotopy nerve of the morphism
$$X_{0}:=\coprod_{i}V_{i} \longrightarrow F'.$$
This is a 
$(n-1)$-smooth Segal groupoid objects in $D^{-}Aff^{\sim,\textrm{\'e}t}$, which is such that
each morphism $X_{i} \longrightarrow i(G)$ is flat. Therefore, by induction on $n$, the natural
morphism of Segal groupoid objects
$$it_{0}(X_{*}) \longrightarrow X_{*}$$
is an equivalence. Therefore, as $i$ and $t_{0}$ commutes with homotopy colimits we find
that the adjunction morphism
$$it_{0}(F')\simeq |it_{0}(X_{*})| \longrightarrow |X_{*}| \simeq F'$$
is an isomorphism in $D^{-}\mathrm{St}(k)$.
\end{proof}

An important corollary of Prop. \ref{pII-3} is the following fact.

\begin{cor}\label{cpII-3}
For any Artin $n$-stack, the $D^{-}$-stack $i(F)$ has an
obstruction theory.
\end{cor}

\begin{proof} Follows from \ref{pII-3} $(3)$ and
Cor. \ref{cII-4}. \end{proof}

One also deduces from Prop. \ref{pII-3} and Lem. \ref{lII-2}
the following corollary.

\begin{cor}\label{cII-4bis}
Let $F$ be an $n$-geometric $D^{-}$-stack. Then, for any
$A\in sk-Alg$, such that $\pi_{i}(A)=0$ for $i>k$, the simplicial
set $\mathbb{R}F(A)$ is $(n+k+1)$-truncated.
\end{cor}

\begin{proof} This is by induction on $k$. For $k=0$
this is Lem. \ref{lII-2} and the fact that $t_{0}$ preserves
$n$-geometric stacks. To pass from $k$ to $k+1$, we consider
for any $A\in sk-Alg$ with $\pi_{i}(A)=0$ for $i>k+1$, the
natural morphisms
$$
\mathbb{R}F(A) \longrightarrow \mathbb{R}F(A_{\leq k}),$$
whose homotopy fibers can be described using
Prop. \ref{p22} and Lem. \ref{lpostn}. We find that this homotopy fiber
is either empty, or equivalent to
$$Map_{A_{\leq k}-Mod}(\mathbb{L}_{F,x},\pi_{k+1}(A)[k+1]),$$
which is $(k+1)$-truncated. By induction, $\mathbb{R}F(A_{\leq k})$ is
$(k+1+n)$-truncated and the homotopy fibers of
$$\mathbb{R}F(A) \longrightarrow \mathbb{R}F(A_{\leq k}),$$
are $(k+1)$-truncated, and therefore $\mathbb{R}F(A)$
is $(k+n+2)$-truncated. \end{proof}

Another important property of the truncation functor is the following
local description of the truncation $t_{0}(F)$ sitting inside
the $D^{-}$-stack $F$ itself.

\begin{prop}\label{pII-6}
Let $F$ be an $n$-geometric $D^{-}$-stack. The adjunction morphism
$it_{0}(F) \longrightarrow F$ is a representable morphism. Moreover, for any
$A\in sk-Alg$, and any flat morphism $\mathbb{R}\underline{Spec}\, A \longrightarrow F$, the square
$$\xymatrix{
it_{0}(F) \ar[r] & F \\
\mathbb{R}\underline{Spec}\, \pi_{0}(A) \ar[u] \ar[r] & \mathbb{R}\underline{Spec}\, A \ar[u]}$$
is homotopy cartesian. In particular, the morphism
$it_{0}(F) \longrightarrow F$ is a closed immersion in the sense of
Def. \ref{dII-7}.
\end{prop}

\begin{proof} By the local character of representable morphisms it is
enough to prove that for any flat morphism $X=\mathbb{R}\underline{Spec}\, A \longrightarrow F$, the square
$$\xymatrix{
it_{0}(F) \ar[r] & F \\
it_{0}(X)=\mathbb{R}\underline{Spec}\, \pi_{0}(A) \ar[u] \ar[r] & X \ar[u]}$$
is homotopy cartesian. The morphism
$$it_{0}(F)\times_{F}^{h}X \longrightarrow it_{0}(F)$$
is flat, and by Prop. \ref{pII-3} $(4)$ this implies
that the $D^{-}$-stack $it_{0}(F)\times_{F}^{h}X$ is truncated. In other words the natural
morphism
$$it_{0}(it_{0}(F)\times_{F}^{h}X)\simeq it_{0}(F)\times_{it_{0}(F)}^{h}it_{0}(X)\simeq
it_{0}(X) \longrightarrow it_{0}(F)\times_{F}^{h}X $$
is an isomorphism. \end{proof}

Using our embedding
$$i : \mathrm{St}(k) \longrightarrow
D^{-}\mathrm{St}(k)$$
we will see stacks in $\mathrm{St}(k)$, and in
particular Artin $n$-stacks, as
$D^{-}$-stacks. However, as the functor $i$ does not
commute with homotopy pullbacks we will
still mention it in order to avoid confusions. \\

\section{Infinitesimal criteria for smooth and \'etale morphisms}\label{IIder.5}

Recall that $sk-Mod_{1}$ denotes the full subcategory of
$sk-Mod$ of connected simplicial $k$-modules. It consists of all
$M\in sk-Mod$ for which the adjunction $S(\Omega M) \longrightarrow M$
is an isomorphism in $\mathrm{Ho}(sk-Mod)$, or equivalently for 
which $\pi_{0}(M)=0$. 

\begin{prop}\label{pII-11}
Let $f : F \longrightarrow G$ be an $n$-representable morphism between
$D^{-}$-stacks.
The morphism $f$ is smooth if and only if it satisfies the following two conditions
\begin{enumerate}
\item The morphism $t_{0}(f) : t_{0}(F) \longrightarrow t_{0}(G)$
is a locally finitely presented morphism in $k-Aff^{\sim,\textrm{\'e}t}$.

\item For any $A\in sk-Alg$, any connected $M\in A-Mod_{s}$, and any
derivation $d\in \pi_{0}(\mathbb{D}er(A,M))$, the natural
projection
$A\oplus_{d} \Omega M \longrightarrow A$
induces a surjective morphism
$$\pi_{0}(\mathbb{R}F(A\oplus_{d}\Omega M))\longrightarrow
\pi_{0}\left(\mathbb{R}G(A\oplus_{d}\Omega M)\times^{h}_{\mathbb{R}F(A)}\mathbb{R}G(A)\right).$$

\end{enumerate}
\end{prop}

\begin{proof} First of all we can suppose that $F$ and $G$ are fibrant
objects in $k-D^{-}Aff^{\sim,\textrm{\'e}t}$.

Suppose first that the morphism
$f$ is smooth and let us prove that is satisfies the two conditions
of the proposition. We know by
\ref{pII-5} $(3)$ that $t_{0}(f)$ is then a smooth morphism in
$k-Aff^{\sim,\textrm{\'e}t}$, so condition $(1)$ is satisfied.
The proof that $(2)$ is also satisfied goes by induction on $n$. Let us start with the case
$n=-1$, and in other words when $f$ is a smooth and representable morphism.

We fix a point $x$ in
$\pi_{0}\left(G(A\oplus_{d}\Omega M)\times^{h}_{F(A)}G(A)\right)$, and
we need to show that the homotopy fiber taken at $x$ of the morphism
$$F(A\oplus_{d}\Omega M)\longrightarrow
G(A\oplus_{d}\Omega M)\times^{h}_{F(A)}G(A)$$
is non empty. The point $x$ corresponds via Yoneda to a commutative diagram
in $\mathrm{Ho}(k-D^{-}Aff^{\sim,\textrm{\'e}t}/G)$
$$\xymatrix{
X \ar[r] \ar[d] & F \ar[d] \\
X_{d}[\Omega M] \ar[r] & G,}$$
where $X:=\mathbb{R}\underline{Spec}\, A$, and
$X_{d}[\Omega M]:=\mathbb{R}\underline{Spec}\, (A\oplus_{d}\Omega M)$.
Making a homotopy base change
$$\xymatrix{
X \ar[d] \ar[r] & F\times^{h}_{G}X_{d}[\Omega M] \ar[d] \\
X_{d}[\Omega M] \ar[r] & X_{d}[\Omega M],}$$
we see that we can replace $G$ be $X_{d}[\Omega M]$ and $f$ by
the projection $F\times^{h}_{G}X_{d}[\Omega M] \longrightarrow X_{d}[\Omega M]$.
In particular, we can assume that $G$ is a representable $D^{-}$-stack.
The morphism $f$ can then be written as
$$f : F\simeq \mathbb{R}\underline{Spec}\, C \longrightarrow \mathbb{R}\underline{Spec}\, B\simeq G,$$
and corresponds to a morphism of commutative simplicial
$k$-algebras $B \longrightarrow C$. Then, using Prop. \ref{p22bis}
and Cor. \ref{cII-4} we see that the obstruction for the point $x$ to lifts
to a point in $\pi_{0}(F(A\oplus_{d}\Omega M))$ lives in the
abelian group $[\mathbb{L}_{C/B}\otimes_{C}^{\mathbb{L}}A,M]$.
But, as $B \longrightarrow C$ is assumed to be smooth, the
$A$-module $\mathbb{L}_{C/B}\otimes_{C}^{\mathbb{L}}A$ is a
retract of a free $A$-module. This implies that
$[\mathbb{L}_{C/B}\otimes_{C}^{\mathbb{L}}A,M]$ is a retract of a
product of $\pi_{0}(M)$, and therefore is $0$ by hypothesis on
$M$. This implies that condition $(2)$ of the proposition is
satisfied when $n=-1$.

Let us now assume that condition $(2)$ is satisfied for all
smooth $m$-representable morphisms for $m<n$, and let us prove it for
a smooth $n$-representable morphism $f : F \longrightarrow G$.
Using the same trick as above, se see that we can always assume that
$G$ is a representable $D^{-}$-stack, and therefore that
$F$ is an $n$-geometric $D^{-}$-stack. Then, let us chose
a point $x$ in
$\pi_{0}\left(G(A\oplus_{d}\Omega M)\times^{h}_{F(A)}G(A)\right)$, and
we need to show that $x$ lifts to a point in
$\pi_{0}(F(A\oplus_{d}\Omega M))$.
For this,
we use
Cor. \ref{cII-4} for $F$, and consider its cotangent complex
$\mathbb{L}_{F,y}\in \mathrm{Ho}(Sp(A-Mod_{s}))$, where
$y\in F(A)$ is the image of $x$. There exists a natural functoriality morphism
$$\mathbb{L}_{G,f(y)} \longrightarrow \mathbb{L}_{F,y}$$
whose homotopy cofiber is $\mathbb{L}_{F/G,y}\in \mathrm{Ho}(Sp(A-Mod_{s}))$.
Then, Prop. \ref{p22bis} tell us that the obstruction
for the existence of this lift lives in
$[\mathbb{L}_{F/G,y},M]$. It is therefore enough to show that
$[\mathbb{L}_{F/G,y},M]\simeq 0$ for any $M\in A-Mod_{s}$ such that $\pi_{0}(M)=0$.

\begin{lem}\label{lII-13}
Let $F \longrightarrow G$ be a smooth morphism between $n$-geometric
$D^{-}$-stacks with $G$ a representable stack. Let $A\in sk-Alg$ and
$y : Y=\mathbb{R}\underline{Spec}\, A \longrightarrow F$ be a point.
Then the object
$$\mathbb{L}_{F/G,y} \in \mathrm{Ho}(Sp(A-Mod_{s}))$$
is perfect, and its dual
$\mathbb{L}_{F/G,y}^{\vee} \in \mathrm{Ho}(Sp(A-Mod_{s}))$
is $0$-connective (i.e. belongs to the image of
$\mathrm{Ho}(A-Mod_{s}) \hookrightarrow \mathrm{Ho}(Sp(A-Mod_{s}))$.
\end{lem}

\begin{proof} Recall first that
an object in $\mathrm{Ho}(Sp(A-Mod_{s}))$ is perfect if and only if
it is finitely presented, and if and only if it is
a retract of a finite cell stable $A$-module (see Cor. \ref{cI-cell}, and also
 \cite[III.2]{ekmm} or  \cite[Thm. III.5.7]{km}).

The proof is then by induction on $n$. When $F$ is representable, this
is by definition of smooth morphisms, as then
$\mathbb{L}_{F/G,y}$ is a projective $A$-module of finite presentation, and
so is its dual. Let us suppose the lemma proved for all $m<n$, and lets
prove it for $n$.
First of all, the conditions on $\mathbb{L}_{F/G,y}$ we need to prove are
local for the \'etale topology on $A$, because of Cor. \ref{ct22stable}.
Therefore, one can assume that the point $y$ lifts to a point of
an $n$-atlas for $F$. One can thus suppose that there exists
a representable $D^{-}$-stack $U$, a smooth morphism $U \longrightarrow F$,
such that $y\in \pi_{0}(\mathbb{R}F(A))$ is the image of a point $z\in \mathbb{R}U(A)$.
There exists an fibration sequence of stable $A$-modules
$$\xymatrix{
\mathbb{L}_{F/G,y} \longrightarrow \mathbb{L}_{U/G,z} \longrightarrow
\mathbb{L}_{U/F,z}.}$$
As $U \longrightarrow G$
is smooth, $\mathbb{L}_{U/G,z}$ is a projective $A$-module. Furthermore,
the stable $A$-module
$\mathbb{L}_{U/F,y}$ can be identified with
$\mathbb{L}_{U\times_{F}^{h}Y/Y,s}$, where $s$ is the natural
section $Y \longrightarrow U\times_{F}^{h}Y$ induced by the point $z : Y \longrightarrow U$.
The morphism $U\times_{F}^{h}Y \longrightarrow Y$ being a smooth and $(n-1)$-representable
morphism, induction tells us that the stable $A$-module $\mathbb{L}_{U/F,y}$
satisfies the conditions of the lemma. Therefore, $\mathbb{L}_{F/G,y}$ is the homotopy fiber
of a morphism between stable $A$-module satisfying the conditions of the lemma, and is
easily seen to satisfies itself these conditions. \end{proof}

By the above lemma, we have
$$[\mathbb{L}_{F/G,y},M]\simeq \pi_{0}(\mathbb{L}_{F/G,y}^{\vee}\otimes_{A}^{\mathbb{L}}M)
\simeq \pi_{0}(\mathbb{L}_{F/G,y}^{\vee})\otimes_{\pi_{0}(A)}\pi_{0}(M)\simeq 0$$
for any $A$-module $M$ such that $\pi_{0}(M)=0$. This finishes the proof of the
fact that $f$ satisfies the conditions of Prop. \ref{pII-11} when it is smooth.

Conversely, let us assume that $f : F \longrightarrow G$ is
a morphism satisfying the lifting property of \ref{pII-11}, and let us show that
$f$ is smooth. Clearly, one can suppose that $G$ is a representable stack, and thus
that $F$ is an $n$-geometric $D^{-}$-stack. We need to show that for
any representable $D^{-}$-stack $U$ and any smooth morphism $U \longrightarrow F$ the
composite morphism $U \longrightarrow G$ is smooth. By what we have
seen in the first part of the proof, we known that
$U \longrightarrow F$ also satisfies the lifting properties, and thus so does
the composition $U \longrightarrow G$. We are therefore reduced to the case
where $f$ is a morphism between  representable $D^{-}$-stacks, and thus
corresponds to a morphism of commutative simplicial $k$-algebras
$A \longrightarrow B$. By hypothesis on $f$, $\pi_{0}(A) \longrightarrow \pi_{0}(B)$
is a finitely presented morphism of commutative rings. Furthermore, Prop. \ref{p22bis} and Cor. \ref{cII-4}
show that
for any $B$-module $M$ with $\pi_{0}(M)=0$, we have
$[\mathbb{L}_{B/A},M]=0$. Let $B^{(I)} \longrightarrow \mathbb{L}_{B/A}$
be a morphism of $B$-modules, with $B^{(I)}$ free over some set $I$, and
such that the induced morphism $\pi_{0}(B)^{(I)} \longrightarrow \pi_{0}(\mathbb{L}_{B/A})$
is surjective. Let $K$ be the homotopy fiber of the morphism
$B^{(I)} \longrightarrow \mathbb{L}_{B/A}$, that, according to our choice,
induces a homotopy fiber sequence of $A$-modules
$$\xymatrix{B^{(I)} \ar[r] & \mathbb{L}_{B/A} \ar[r] & K[1].}$$
The short exact sequence
$$[\mathbb{L}_{B/A},B^{(I)}] \longrightarrow
[\mathbb{L}_{B/A},\mathbb{L}_{B/A}] \longrightarrow [\mathbb{L}_{B/A},K[1]]=0,$$
shows that $\mathbb{L}_{B/A}$ is a retract of $B^{(I)}$, and thus is
a projective $B$-module. Furthermore, the homotopy cofiber sequence
$$\xymatrix{\mathbb{L}_{A}\otimes_{A}^{\mathbb{L}}B \ar[r] & \mathbb{L}_{B} \ar[r] & \mathbb{L}_{B/A},}$$
induces a short exact sequence
$$\xymatrix{
[\mathbb{L}_{B},\mathbb{L}_{A}\otimes_{A}^{\mathbb{L}}B]\ar[r] &
[\mathbb{L}_{A}\otimes_{A}^{\mathbb{L}}B,\mathbb{L}_{A}\otimes_{A}^{\mathbb{L}}B] \ar[r] &
[\mathbb{L}_{B/A},\mathbb{L}_{A}\otimes_{A}^{\mathbb{L}}B[1]]=0,}$$
shows that the morphism $\mathbb{L}_{A}\otimes_{A}^{\mathbb{L}}B \longrightarrow
\mathbb{L}_{B}$ has a retraction. We conclude that $f$ is a formally smooth
morphism such that $\pi_{0}(A) \longrightarrow \pi_{0}(B)$ is finitely presented, and
by Cor. \ref{ctII-1'} that $f$ is a smooth morphism. \end{proof}

From the proof of Prop. \ref{pII-11} we extract the following corollary.

\begin{cor}\label{cpII-11}
Let $F \longrightarrow G$ be an $n$-representable morphism of
$D^{-}$-stacks such that the morphism $t_{0}(F) \longrightarrow t_{0}(G)$ is
a locally finitely presented morphism of stacks.
The following three conditions are equivalent.
\begin{enumerate}
\item The morphism $f$ is smooth.
\item
For any $A\in sk-Alg$ and any morphism of stacks
$x : X=\mathbb{R}\underline{Spec}\, A \longrightarrow F$,
the object
$$\mathbb{L}_{F/G,x} \in \mathrm{Ho}(Sp(A-Mod_{s}))$$
is perfect, and its dual
$\mathbb{L}_{F/G,y}^{\vee} \in \mathrm{Ho}(Sp(A-Mod_{s}))$
is $0$-connective (i.e. belongs to the image of
$\mathrm{Ho}(A-Mod_{s}) \hookrightarrow \mathrm{Ho}(Sp(A-Mod_{s}))$.
\item For any $A\in sk-Alg$, any morphism of stacks
$x : X=\mathbb{R}\underline{Spec}\, A \longrightarrow F$,
and any $A$-module $M$ in $sk-Mod_{1}$, we have
$$[\mathbb{L}_{F/G,x},M]=0.$$
\end{enumerate}
\end{cor}

\begin{proof} That $(1)$ implies $(2)$ follows from lemma \ref{lII-13}. Conversely,
if $\mathbb{L}_{F/G,x}$ satisfies the conditions of the corollary, then for
any $A$-module $M$ such that $\pi_{0}(M)$ one has
$[\mathbb{L}_{F/G,x},M]=0$. Therefore, Prop. \ref{p22bis} shows that
the lifting property of Prop. \ref{pII-11} holds, and thus
that $(2)$ implies $(1)$. Furthermore, clearly $(2)$ implies
$(3)$, and conversely $(3)$ together with
Prop. \ref{p22bis} implies the lifting property of Prop. \ref{pII-11}.
\end{proof}

\begin{prop}\label{pII-12}
Let $f : F \longrightarrow G$ be an $n$-representable morphism between
$D^{-}$-stacks.
The morphism $f$ is \'etale if and only if it satisfies the following two conditions
\begin{enumerate}
\item The morphism $t_{0}(f) : t_{0}(F) \longrightarrow t_{0}(G)$
is locally finitely presented as a morphism in $k-Aff^{\sim,\tau}$.

\item For any $A\in sk-Alg$, any $M\in A-Mod_{s}$
whose underlying $k$-module is in $sk-Mod_{1}$, and any
derivation $d\in \pi_{0}(\mathbb{D}er(A,M))$, the natural
projection
$A\oplus_{d} \Omega M \longrightarrow A$
induces an isomorphism in $\mathrm{Ho}(SSet)$
$$\mathbb{R}F(A\oplus_{d}\Omega M) \longrightarrow
\mathbb{R}G(A\oplus_{d}\Omega M)\times^{h}_{\mathbb{R}F(A)}\mathbb{R}G(A).$$

\end{enumerate}
\end{prop}

\begin{proof} First of all one can suppose that $F$ and $G$ are fibrant
objects in $k-D^{-}Aff^{\sim,\textrm{\'e}t}$.

Suppose first that the morphism
$f$ is \'etale and let us prove that is satisfies the two conditions
of the proposition. We know by
\ref{pII-5} $(3)$ that $t_{0}(f)$ is then a \'etale morphism in
$k-Aff^{\sim,\textrm{\'e}t}$, so condition $(1)$ is satisfied.
The proof that $(2)$ is also satisfied goes by induction on $n$. Let us start with the case
$n=-1$, and in other words when $f$ is an \'etale and representable morphism.
In this case the result follows from Cor. \ref{cpismooth} $(2)$. We now assume the result
for all $m<n$ and prove it for $n$. Let $A\in sk-Alg$, $M$ an $A$-module
with $\pi_{0}(M)=0$, and
$d\in \pi_{0}(\mathbb{D}er(A,M))$ be a derivation.
Let $x$ be a point in $\pi_{0}(\mathbb{R}G(A\oplus_{d}\Omega M)\times^{h}_{\mathbb{R}F(A)}\mathbb{R}G(A))$,
with image $y$ in $\pi_{0}(\mathbb{R}F(A))$.
By
Prop. \ref{p22bis} the homotopy fiber of the morphism
$$\mathbb{R}F(A\oplus_{d}\Omega M) \longrightarrow
\mathbb{R}G(A\oplus_{d}\Omega M)\times^{h}_{\mathbb{R}F(A)}\mathbb{R}G(A)$$
at the point $x$ is non empty if and only if a certain
obstruction in $[\mathbb{L}_{F/G,y},M]$ vanishes. Furthermore, if nonempty this
homotopy fiber is then equivalent to $$Map_{A-Mod}(\mathbb{L}_{F/G,y},\Omega (M)).$$
It is then enough to prove that $\mathbb{L}_{F/G,y}=0$, and this is contained in the following lemma.\\

\begin{lem}\label{lII-14}
Let $F \longrightarrow G$ be a \'etale morphism between $n$-geometric
$D^{-}$-stacks. Let $A\in sk-Alg$ and
$y : Y=\mathbb{R}\underline{Spec}\, A \longrightarrow F$ be a point.
Then $\mathbb{L}_{F/G,y}\simeq 0$.
\end{lem}

\begin{proof} We easily reduce to the case when $G$ is a representable $D^{-}$-stack.\\
The proof is then by induction on $n$. When $F$ is representable, this
is by definition of \'etale morphisms.
Let us suppose the lemma proved for all $m<n$, and lets
prove it for $n$.
First of all, the vanishing of $\mathbb{L}_{F/G,y}$ is clearly
a local condition  for the \'etale topology on $A$.
Therefore, we can assume that the point $y$ lifts to a point of
an $n$-atlas for $F$. We can thus suppose that there exists
a representable $D^{-}$-stack $U$, a smooth morphism $U \longrightarrow F$,
such that $y\in \pi_{0}(\mathbb{R}F(A))$ is the image of a point $z\in \mathbb{R}U(A)$.
There exists an fibration sequence of stable $A$-modules
$$\xymatrix{
\mathbb{L}_{F/G,y} \longrightarrow \mathbb{L}_{U/G,z} \longrightarrow
\mathbb{L}_{U/F,z}.}$$
 As $U \longrightarrow G$
is \'etale, $\mathbb{L}_{U/G,z}\simeq 0$. Furthermore,
the stable $A$-module
$\mathbb{L}_{U/F,y}$ can be identified with
$\mathbb{L}_{U\times_{F}^{h}Y/Y,s}$, where $s$ is the natural
section $Y \longrightarrow U\times_{F}^{h}Y$ induced by the point $z : Y \longrightarrow U$.
The morphism $U\times_{F}^{h}Y \longrightarrow Y$ being an \'etale and $(n-1)$-representable
morphism, induction tells us that the stable $A$-module $\mathbb{L}_{U/F,y}$
vanishes. We conclude that $\mathbb{L}_{F/G,y}\simeq 0$. \end{proof}

The lemma finishes the proof that if $f$ is \'etale then
it satisfies the two conditions of the proposition. 

Conversely, let
us assume that $f$ satisfies the properties $(1)$ and $(2)$ of the proposition.
To prove that $f$ is \'etale we can suppose that $G$ is a representable $D^{-}$-stack.
We then need to show that for any representable $D^{-}$-stack $U$ and any
smooth morphism $u : U \longrightarrow F$, the induced morphism
$v : U \longrightarrow G$ is \'etale. But Prop. \ref{p22bis} and our
assumption $(2)$ easily implies
$\mathbb{L}_{F/G,v}\simeq 0$. Furthermore, the obstruction for the homotopy cofiber sequence
$$\xymatrix{\mathbb{L}_{U/G,v} \ar[r] & \mathbb{L}_{F/G,u} \ar[r] & \mathbb{L}_{U/F,u},}$$
to splits lives in $[\mathbb{L}_{U/F,u},S(\mathbb{L}_{U/G,v})]$, which is zero
by Cor. \ref{cpII-11}. Therefore, $\mathbb{L}_{U/G,v}$ is a retract
of $\mathbb{L}_{F/G,u}$ and thus vanishes.
This implies that $U\longrightarrow G$ is a formally \'etale
morphism of representable $D^{-}$-stacks. Finally, our assumption $(1)$ and
Cor. \ref{ctII-1'} implies that $U \longrightarrow G$
is an \'etale morphism as required. \end{proof}

\begin{cor}\label{cpII-12}
Let $F \longrightarrow G$ be an $n$-representable morphism of
$D^{-}$-stacks such that the morphism $t_{0}(F) \longrightarrow t_{0}(G)$ is
a locally finitely presented morphism of stacks.
The following two conditions are equivalent.
\begin{enumerate}
\item The morphism $f$ is \'etale.
\item
For any $A\in sk-Alg$ and any morphism of stacks
$x : X=\mathbb{R}\underline{Spec}\, A \longrightarrow F$,
one has
$\mathbb{L}_{F/G,x} \simeq 0$.
\end{enumerate}
\end{cor}

\begin{proof} This follows from Prop. \ref{pII-12} and
Prop. \ref{p22bis}. \end{proof}

\section{Some examples of geometric $D^{-}$-stacks}\label{IIder.6}

We present here some basic examples of geometric $D^{-}$-stacks. Of course we do not claim
to be exhaustive, and many other
interesting examples will not be discussed here and will appear in future works
(see e.g. \cite{go,tv}).\\

\subsection{Local systems}\label{IIder.6.1}

Recall from Def. \ref{d30'} the existence of the $D^{-}$-stack $\mathbf{Vect}_{n}$, of
rank $n$ vector bundles. Recall by \ref{lII-3} that for $A\in sk-Alg$, an $A$-module
$M\in A-Alg_{s}$ is a rank $n$ vector bundle, if and only if
$M$ is a strong $A$-module and $\pi_{0}(M)$ is a projective $\pi_{0}(A)$-module of rank $n$. Recall
also from Lem. \ref{lII-3} that vector bundles are precisely the locally perfect modules.
The conditions of \ref{cp24} are all satisfied in the present context and therefore
we know that $\mathbf{Vect}_{n}$ is a smooth $1$-geometric $D^{-}$-stack.
As a consequence of Prop. \ref{pII-5} $(4)$ we deduce
that the $D^{-}$-stack $\mathbf{Vect}_{n}$ is truncated in the sense of
\ref{dII-8}. We also have a stack of rank $n$ vector bundles
$\mathbf{Vect}_{n}$ in $\mathrm{St}(k)$. Using the same notations for these
two different objects is justified by the following lemma.

\begin{lem}\label{lII-10}
There exists a natural isomorphism in $D^{-}\mathrm{St}(k)$
$$i(\mathbf{Vect}_{n})\simeq \mathbf{Vect}_{n}.$$
\end{lem}

\begin{proof} As we know that the $D^{-}$-stack
$\mathbf{Vect}_{n}$ is truncated, it is equivalent to show that there exists
a natural isomorphism
$$\mathbf{Vect}_{n}\simeq t_{0}(\mathbf{Vect}_{n})$$
in $\mathrm{Ho}(k-Aff^{\sim,\tau})$.

We start by defining a morphism of stacks
$$\mathbf{Vect}_{n}\longrightarrow t_{0}(\mathbf{Vect}_{n}).$$
For this, we construct
for a commutative $k$-algebra $A$, a natural functor
$$\phi_{A} : A-QCoh^{c}_{W} \longrightarrow i(A)-QCoh^{c}_{W},$$
where $i(A) \in sk-Alg$ is the constant simplicial commutative $k$-algebra
associated to $A$, and where $A-QCoh^{c}_{W}$ and $i(A)-QCoh^{c}_{W}$
are defined in \S \ref{Iqcoh}. Recall that
$A-QCoh^{c}_{W}$ is the category whose objects are
the data of a $B$-module $M_{B}$ for any morphism $A \longrightarrow B$ in $k-Alg$, together with
isomorphisms $M_{B}\otimes_{B}B' \simeq M_{B'}$ for any
$A\longrightarrow B\longrightarrow B'$ in $k-Alg$, satisfying the usual
cocycle conditions. The morphisms $M\rightarrow M'$ in $A-QCoh^{c}_{W}$ are simply the families of
isomorphisms $M_{B}\simeq M'_{B}$ of $B$-modules which commute with the transitions isomorphisms.
In the same way, the objects in $i(A)-QCoh^{c}_{W}$ are the data
a cofibrant $B$-module $M_{B}\in B-Mod_{s}$ for any morphism $i(A) \longrightarrow B$ in $sk-Alg$, together with
equivalences $M_{B}\otimes_{B}B' \longrightarrow M_{B'}$ for any
$i(A)\longrightarrow B\longrightarrow B'$ in $sk-Alg$, satisfying the usual
cocycle conditions. The morphisms $M\rightarrow M'$ in $i(A)-QCoh^{c}_{W}$ are simply the
equivalences $M_{B} \longrightarrow M'_{B}$ which commutes with the transition
equivalences. The functor
$$\phi_{A} : A-QCoh^{c}_{W} \longrightarrow i(A)-QCoh^{c}_{W},$$
sends an object $M$ to the object $\phi_{A}(M)$, where $\phi_{A}(M)_{B}$
is the simplicial $B$-module defined by
$$(\phi_{A}(M)_{B})_{n}:=M_{B_{n}}$$
(note that $B$ can be seen as a simplicial commutative $A$-algebra). The functor
$\phi_{A}$ is clearly functorial in $A$ and thus defines a morphism of simplicial presheaves
$$\mathbf{QCoh} \longrightarrow t_{0}(\mathbf{QCoh}).$$
We check easily that the sub-stack
$\mathbf{Vect}_{n}$ of $\mathbf{QCoh}$ is sent to the sub-stack
$t_{0}(\mathbf{Vect}_{n})$ of $t_{0}(\mathbf{QCoh})$, and therefore we get a morphism of stacks
$$\mathbf{Vect}_{n} \longrightarrow t_{0}(\mathbf{Vect}_{n}).$$
To see that this morphism is an isomorphism of stacks we construct a morphism in the other
direction by sending a simplicial $i(A)$-module
$M$ to the $\pi_{0}(A)$-module $\pi_{0}(M)$. By \ref{lII-3} we easily see that this
defines an inverse of the above morphism.
\end{proof}

For a simplicial set $K\in SSet_{\mathbb{U}}$, and an object
$F\in k-D^{-}Aff^{\sim,\textrm{\'e}t}$, we can use the simplicial structure of the category
$k-D^{-}Aff^{\sim,\textrm{\'e}t}$ in order to define the exponential
$F^{K}\in k-D^{-}Aff^{\sim,\textrm{\'e}t}$. The model category $k-D^{-}Aff^{\sim,\textrm{\'e}t}$ being a
simplicial model category the functor
$$(-)^{K} : k-D^{-}Aff^{\sim,\textrm{\'e}t} \longrightarrow k-D^{-}Aff^{\sim,\textrm{\'e}t}$$
is right Quillen, and therefore can be derived on the right. Its right derived
functor will be denoted by
$$\begin{array}{ccc}
D^{-}\mathrm{St}(k) & \longrightarrow & D^{-}\mathrm{St}(k) \\
F & \mapsto & F^{\mathbb{R}K}.
\end{array}$$

Explicitly, we have
$$F^{\mathbb{R}K}\simeq (RF)^{K}$$
where $RF$ is a fibrant replacement of $F$ in $k-D^{-}Aff^{\sim,\textrm{\'e}t}$.

\begin{df}\label{dII-11}
Let $K$ be a $\mathbb{U}$-small simplicial set. The \emph{derived moduli
stack of rank $n$ local systems on $K$}\index{$\mathbb{R}\mathbf{Loc}_{n}(K)$!derived moduli
stack of rank $n$ local systems on $K$} is defined to be
$$\mathbb{R}\mathbf{Loc}_{n}(K):=\mathbf{Vect}_{n}^{\mathbb{R}K}.$$
\end{df}

We start by the following easy observation.

\begin{lem}\label{lII-11}
Assume that $K$ is a finite dimensional simplicial set.
Then, the $D^{-}$-stack $\mathbb{R}\mathbf{Loc}_{n}(K)$ is a
finitely presented $1$-geometric
$D^{-}$-stack.
\end{lem}

\begin{proof}  We consider the following
homotopy co-cartesian square of simplicial sets
$$\xymatrix{
Sk_{i}K \ar[r] & Sk_{i+1}K \\
\coprod_{Hom(\partial\Delta^{i+1},K)}\partial\Delta^{i+1} \ar[r] \ar[u] &
\coprod_{Hom(\Delta^{i+1},K)}\Delta^{i+1},\ar[u]}$$
where $Sk_{i}K$ is the skeleton of dimension of $i$ of $K$. This gives a homotopy pullback
square of $D^{-}$-stacks
$$\xymatrix{
\mathbb{R}\mathbf{Loc}_{n}(Sk_{i+1}K)\ar[d] \ar[r] & \mathbb{R}\mathbf{Loc}_{n}(Sk_{i}K) \ar[d] \\
\prod^{h}_{Hom(\Delta^{i+1},K)}\mathbb{R}\mathbf{Loc}_{n}(\Delta^{i+1}) \ar[r] &
\prod^{h}_{Hom(\partial\Delta^{i+1},K)}\mathbb{R}\mathbf{Loc}_{n}(\partial\Delta^{i+1}).}$$
As finitely presented $1$-geometric $D^{-}$-stacks are stable
by homotopy pullbacks, we see
by induction on the skeleton that it only remains to show that
$\mathbb{R}\mathbf{Loc}_{n}(\partial\Delta^{i+1})$ is a finitely presented $1$-geometric
$D^{-}$-stack. But there is
an isomorphism in $\mathrm{Ho}(SSet)$,
$$\partial\Delta^{i+1}\simeq *\coprod_{\partial\Delta^{i}}^{\mathbb{L}}*,$$
giving rise to an isomorphism of $D^{-}$-stacks
$$\mathbb{R}\mathbf{Loc}_{n}(\partial\Delta^{i+1})\simeq
\mathbb{R}\mathbf{Loc}_{n}\times^{h}_{\mathbb{R}\mathbf{Loc}_{n}(\partial\Delta^{i})}\mathbb{R}\mathbf{Loc}_{n}.$$
By induction on $i$, we see that it is enough to show that
$\mathbf{Vect}_{n}$ is a finitely presented $1$-geometric $D^{-}$-stack
which is known from Cor. \ref{cp24}. \end{proof}

Another easy observation is the description of the truncation
$t_{0}\mathbb{R}\mathbf{Loc}_{n}(K)$. For this, recall that
the simplicial set $K$ has a fundamental groupoid $\Pi_{1}(K)$. The usual
Artin stack of rank $n$ local systems on $K$ is the stack in groupoids defined by
$$\begin{array}{cccc}
\mathbf{Loc}_{n}(K) : &  k-Alg & \longrightarrow & \{Groupoids\} \\
 & A & \mapsto & \underline{Hom}(\Pi_{1}(K),\mathbf{Vect}_{n}(A)).
\end{array}$$
In other words, $\mathbf{Loc}_{n}(K)(A)$ is the groupoid
of functors from $\Pi_{1}(K)$ to the groupoid of rank $n$ projective
$A$-modules. As usual, this Artin stack is considered as
an object in $\mathrm{St}(k)$.

\begin{lem}\label{lII-12}
There exists a natural isomorphism in $\mathrm{St}(k)$
$$\mathbf{Loc}_{n}(K)\simeq t_{0}\mathbb{R}\mathbf{Loc}_{n}(K).$$
\end{lem}

\begin{proof} The truncation $t_{0}$ being the right derived functor of a right Quillen
functor commutes with derived exponentials. Therefore, we have
$$t_{0}\mathbb{R}\mathbf{Loc}_{n}(K)\simeq
(t_{0}\mathbb{R}\mathbf{Loc}_{n})^{\mathbb{R}K}\simeq
(\mathbf{Vect}_{n})^{\mathbb{R}K}\in \mathrm{St}(k).$$
The stack $\mathbf{Vect}_{n}$ is $1$-truncated, and therefore
we also have natural isomorphisms
$$(\mathbf{Vect}_{n})^{\mathbb{R}K}(A)\simeq
Map_{SSet}(K,\mathbf{Vect}_{n}(A))\simeq
\underline{Hom}(\Pi_{1}(K),\Pi_{1}(\mathbf{Vect}_{n}(A))).$$
This equivalence is natural in $A$ and provides the isomorphism of the lemma. \end{proof}

Our next step is to give a more geometrical interpretation
of the $D^{-}$-stack $\mathbb{R}\mathbf{Loc}_{n}(K)$, in terms of certain
local systems of objects on the topological realization of $K$.

Let now $X$ be a $\mathbb{U}$-small topological space.
Let $A\in sk-Mod$, and let us define a model category
$A-Mod_{s}(X)$, of $A$-modules over $X$. The category
$A-Mod_{s}(X)$ is simply the category of
presheaves on $X$ with values in $A-Mod_{s}$. The model structure
on $A-Mod_{s}(X)$ is of the same type as the local projective model
structure on simplicial presheaves. We first define an intermediate
model structure on $A-Mod_{s}(X)$ for which equivalences (resp. fibrations)
are morphism $\mathcal{E} \longrightarrow \mathcal{F}$ in $A-Mod_{s}(X)$
such that for any open subset $U\subset X$ the induced morphism
$\mathcal{E}(U) \longrightarrow \mathcal{F}(U)$ is an equivalence (resp.
a fibration). This model structure exists as $A-Mod_{s}$ is a
$\mathbb{U}$-cofibrantly generated model category, and let us call
it the \emph{strong} model structure. The final model structure on
$A-Mod_{s}(X)$ is the one for which
cofibrations are the same cofibrations as for the strong model structure, and equivalences
are the morphisms $\mathcal{E} \longrightarrow \mathcal{F}$ such that for any
point $x\in X$ the induced morphism on the stalks
$\mathcal{E}_{x} \longrightarrow \mathcal{F}_{x}$ is an equivalence in $A-Mod_{s}$.
The existence of this model structure is proved the same way as for the case
of simplicial presheaves (we can also use the forgetful functor
$A-Mod_{s}(X) \longrightarrow SPr(X)$ to lift the local projective model structure
on $SPr(X)$ in a standard way).

For a commutative simplicial $k$-algebra $A$, we consider
$A-Mod_{s}(X)^{c}_{W}$, the subcategory of
$A-Mod_{s}(X)$ consisting of cofibrant objects and equivalences between them.
For a morphism of commutative simplicial $k$-algebras $A \longrightarrow B$, we have
a base change functor
$$\begin{array}{ccc}
A-Mod_{s}(X)&  \longrightarrow & B-Mod_{s}(X)\\
\mathcal{E} & \mapsto & \mathcal{E}\otimes_{A}B
\end{array}$$
which is a left Quillen functor, and therefore induces a well defined functor
$$-\otimes_{A}B : A-Mod_{s}(X)^{c}_{W} \longrightarrow B-Mod_{s}(X)^{c}_{W}.$$
This defines a lax functor $A \mapsto A-Mod_{s}(X)^{c}_{W}$, from
$sk-Alg$ to $Cat$, which can be strictified in the usual way. We will omit to mention
explicitly this strictification here and will do as if
$A \mapsto A-Mod_{s}(X)^{c}_{W}$ does define a genuine functor
$sk-Alg \longrightarrow Cat$.

We then define a sub-functor of $A \mapsto
A-Mod_{s}(X)^{c}_{W}$ in the following way. For $A\in sk-Alg$,
let $A-Loc_{n}(X)$ be the full subcategory of
$A-Mod_{s}(X)^{c}_{W}$ consisting of objects $\mathcal{E}$, such
that there exists an open covering $\{U_{i}\}$ on $X$, such that
each restriction $\mathcal{E}_{|U_{i}}$ is isomorphic
in $\mathrm{Ho}(A-Mod_{s}(U_{i}))$ to a constant
presheaf with fibers a projective $A$-module of rank $n$
(i.e. projective $A$-module $E$ such that $\pi_{0}(E)$ is a projective
$\pi_{0}(A)$-module of rank $n$). This defines a sub-functor of
$A\mapsto A-Mod_{s}(X)^{c}_{W}$, and thus a functor
from $sk-Alg$ to $Cat$. Applying the nerve functor we obtain a simplicial presheaf
$\mathbb{R}\mathbf{Loc}_{n}(X) \in k-D^{-}Aff^{\sim,\textrm{\'e}t}$, defined by
$$\mathbb{R}\mathbf{Loc}_{n}(X)(A):=N(A-Loc_{n}(X)).$$

\begin{prop}\label{pII-7}
Let $K$ be a simplicial set in $\mathbb{U}$ and $|K|$ be its topological
realization.
The simplicial presheaf $\mathbb{R}\mathbf{Loc}_{n}(|K|)$ is a $D^{-}$-stack, and there
exists an isomorphism in $D^{-}\mathrm{St}(k)$
$$\mathbb{R}\mathbf{Loc}_{n}(|K|)\simeq \mathbb{R}\mathbf{Loc}_{n}(K).$$
\end{prop}

\begin{proof} We first remark that if $f : X \longrightarrow X'$ is a
homotopy equivalence of topological spaces, then the induced morphism
$$f^{*} : \mathbb{R}\mathbf{Loc}_{n}(X') \longrightarrow \mathbb{R}\mathbf{Loc}_{n}(X')$$
is an equivalence of simplicial presheaves. Indeed, a standard argument reduces
to the case where $f$ is the projection $X\times [0,1] \longrightarrow X$, and then
to $[0,1] \longrightarrow *$, for which one can use the same argument as
in \cite[Lem. 2.16]{to4}.

Let us first prove that
$\mathbb{R}\mathbf{Loc}_{n}(|K|)$ is a $D^{-}$-stack, and for this we will prove
that the simplicial presheaf $\mathbb{R}\mathbf{Loc}_{n}(|K|)$ can be written,
in $SPr(k-D^{-}Aff)$, as a certain homotopy limit of $D^{-}$-stacks. As $D^{-}$-stacks
in $SPr(k-D^{-}Aff)$ are stable by homotopy limits this will prove what we want.

Let $U_{*}$ an open hypercovering of $|K|$ such that each
$U_{i}$ is a coproduct of contractible open subsets in $|K|$ (such
a hypercover exists as $|K|$ is a locally contractible space).
It is not hard to show using Cor. \ref{cstrict} and standard
cohomological descent, that for any
$A\in sk-Alg$, the
natural morphism
$$N(A-Mod_{s}(|K|)^{c}_{W}) \longrightarrow
Holim_{m\in \Delta}N(A-Mod_{s}(U_{m})^{c}_{W})$$
is an isomorphism in $\mathrm{Ho}(SSet)$. From this, we easily deduce that the
natural morphism
$$\mathbb{R}\mathbf{Loc}_{n}(|K|) \longrightarrow
Holim_{m\in \Delta}\mathbb{R}\mathbf{Loc}_{n}(U_{m})$$
is an isomorphism in $\mathrm{Ho}(SPr(k-D^{-}Aff))$. Therefore, we are reduced to show that
$\mathbb{R}\mathbf{Loc}_{n}(U_{m})$ is a $D^{-}$-stack. But, $U_{m}$ being
a coproduct of contractible topological spaces, $\mathbb{R}\mathbf{Loc}_{n}(U_{m})$
is a product of some $\mathbb{R}\mathbf{Loc}_{n}(U)$ for some contractible space
$U$. Moreover $\mathbb{R}\mathbf{Loc}_{n}(U)$ is naturally isomorphic
in $\mathrm{Ho}(SPr(k-D^{-}Aff))$ to $\mathbb{R}\mathbf{Loc}_{n}(*)=\mathbf{Vect}_{n}$.
As we know that $\mathbf{Vect}_{n}$ is a $D^{-}$-stack, this shows that
the simplicial presheaf $\mathbb{R}\mathbf{Loc}_{n}(|K|)$ is a homotopy limit of $D^{-}$-stacks
and thus is itself a $D^{-}$-stack.

We are left to prove that $\mathbb{R}\mathbf{Loc}_{n}(|K|)$ and
$\mathbb{R}\mathbf{Loc}_{n}(K)$ are isomorphic. But, we have seen that
$$\mathbb{R}\mathbf{Loc}_{n}(|K|) \simeq Holim_{m\in \Delta}\mathbb{R}\mathbf{Loc}_{n}(U_{m}),$$
for an open hypercover $U_{*}$ of $|K|$, such that each
$U_{m}$ is a coproduct of contractible open subsets. We let
$K':=\pi_{0}(U_{*})$ be the simplicial set of connected components of $U_{*}$, and thus
$$Holim_{m\in \Delta}\mathbb{R}\mathbf{Loc}_{n}(U_{m})\simeq
Holim_{m\in \Delta}\mathbb{R}\mathbf{Loc}_{n}(K'_{m}),$$
where $K'_{m}$ is considered as a discrete topological space. We have thus proved
that
$$\mathbb{R}\mathbf{Loc}_{n}(|K|)\simeq Holim_{m\in \Delta}\prod^{h}_{K'_{m}}
\mathbf{Vect}_{n}\simeq (\mathbf{Vect}_{n})^{\mathbb{R}K'}.$$
But, by \cite[Lem. 2.10]{to4}, we know that $|K|$ is homotopically equivalent to
$|K'|$, and thus that $K$ is equivalent to $K'$. This implies that
$$\mathbb{R}\mathbf{Loc}_{n}(|K|)\simeq (\mathbf{Vect}_{n})^{\mathbb{R}K'}\simeq
(\mathbf{Vect}_{n})^{\mathbb{R}K}\simeq \mathbb{R}\mathbf{Loc}_{n}(|K|).$$
\end{proof}

We will now describe the cotangent complex of $\mathbb{R}\mathbf{Loc}_{n}(K)$. For this, we
fix a global point
$$E : * \longrightarrow \mathbb{R}\mathbf{Loc}_{n}(K),$$ which by
Lem. \ref{lII-12} corresponds to a functor
$$E : \Pi_{1}(K) \longrightarrow \Pi_{1}(\mathbf{Vect}_{n}(k)),$$
where $\Pi_{1}(\mathbf{Vect}_{n}(k))$ can be identified with the
groupoid of rank $n$ projective $k$-modules. The object $E$ is
thus a local system of rank $n$ projective $k$-modules on $K$ in
the usual sense. We will compute the cotangent
complex $\mathbb{L}_{\mathbb{R}\mathbf{Loc}_{n}(K),E} \in
\mathrm{Ho}(Sp(sk-Mod))$ (recall that
$\mathrm{Ho}(Sp(sk-Mod))$ can be naturally identified with the
unbounded derived category of $k$). For this, we let
$E\otimes_{k}E^{\vee}$ be the local system on $K$ of endomorphisms
of $E$, and $C_{*}(K,E\otimes_{k}E^{\vee})$ will be the complex of
homology of $K$ with coefficients in the local system
$E\otimes_{k}E^{\vee}$. We consider
$C_{*}(K,E\otimes_{k}E^{\vee})$ as an unbounded complex of
$k$-modules, and therefore as an object in $\mathrm{Ho}(Sp(sk-Mod))$.

\begin{prop}\label{pII-8}
There exists an isomorphism in
$\mathrm{Ho}(Sp(sk-Mod))\simeq \mathrm{Ho}(C(k))$
$$\mathbb{L}_{\mathbb{R}\mathbf{Loc}_{n}(K),E}\simeq
C_{*}(K,E\otimes_{k}E^{\vee})[-1].$$
\end{prop}

\begin{proof} Let $M\in sk-Mod$, and let us consider the
simplicial set $$\mathbb{D}er_{E}(\mathbb{R}\mathbf{Loc}_{n}(K),M),$$ of derivations
of $\mathbb{R}\mathbf{Loc}_{n}(K)$ at the point $E$ and with coefficients in $M$.
By definition of $\mathbb{R}\mathbf{Loc}_{n}(K)$ we have
$$\mathbb{D}er_{E}(\mathbb{R}\mathbf{Loc}_{n}(K),M)\simeq
Map_{SSet/\mathbf{Vect}_{n}(k)}(K,\mathbf{Vect}_{n}(k\oplus M)),$$
where $K \longrightarrow \mathbf{Vect}_{n}(k)$ is given by the
object $E$, and $\mathbf{Vect}_{n}(k\oplus M) \longrightarrow \mathbf{Vect}_{n}(k)$
is the natural projection. At this point we use Prop. \ref{papp2}
in order to describe, functorially in $M$, the
morphism $\mathbf{Vect}_{n}(k\oplus M) \longrightarrow \mathbf{Vect}_{n}(k)$.
For this, we let $\mathcal{G}(k)$ to be the groupoid
of projective $k$-modules of rank $n$. We also define
an $S$-category $\mathcal{G}(k\oplus M)$ in the following way.
Its objects are projective $k$-modules of rank $n$. The
simplicial set of morphisms in $\mathcal{G}(k\oplus M)$ between two
such $k$-modules $E$ and $E'$ is defined to be
$$\mathcal{G}(k\oplus M)_{(E,E')}:=\underline{Hom}^{Eq}_{(k\oplus M)-Mod_{s}}(E\oplus (E\otimes_{k}M),E'\oplus
(E'\otimes_{k}M)),$$
the simplicial set of equivalences from $E\oplus (E\otimes_{k}M)$ to $E'\oplus (E'\otimes_{k}M)$, in the model category
$(k\oplus M)-Mod_{s}$. It is important to note that
$E\oplus (E\otimes_{k}M)$ is isomorphic to $E\otimes_{k}(k\oplus M)$, and therefore is a cofibrant
object in $(k\oplus M)-Mod_{s}$ (as the base change of a cofibrant
object $E$ in $sk-Mod$). There exists a natural
morphism of $S$-categories
$$\mathcal{G}(k\oplus M) \longrightarrow \mathcal{G}(k)$$
being the identity on the set of objects, and the composition of natural morphisms
$$\underline{Hom}^{Eq}_{(k\oplus M)-Mod_{s}}(E\oplus (E\otimes_{k}M),E'\oplus
(E'\otimes_{k}M)) \longrightarrow \underline{Hom}^{Eq}_{k-Mod_{s}}(E,E')\longrightarrow $$
$$\pi_{0}(\underline{Hom}^{Eq}_{k-Mod_{s}}(E,E'))\simeq \mathcal{G}(k)_{(E,E')}.$$
on the simplicial sets of morphisms. Clearly, Prop. \ref{papp2} and its functorial
properties, show that the morphism
$$\mathbf{Vect}_{n}(k\oplus M) \longrightarrow \mathbf{Vect}_{n}(k)$$
is equivalent to
$$N(\mathcal{G}(k\oplus M)) \longrightarrow N(\mathcal{G}(k)),$$
and in a functorial way in $M$.

For any $E \in \mathcal{G}(k)$, we can consider the classifying simplicial set
$K(E\otimes_{k}E^{\vee}\otimes_{k}M,1)$ of the simplicial abelian
group $E\otimes_{k}E^{\vee}\otimes_{k}M$, and for any isomorphism of projective $k$-modules
of rank $n$, $E \simeq E'$, the corresponding
isomorphism of simplicial set
$$K(E\otimes_{k}E^{\vee}\otimes_{k}M,1) \simeq K(E'\otimes_{k}(E')^{\vee}\otimes_{k}M,1).$$
This defines a local system $L$ of simplicial sets on the groupoid
$\mathcal{G}(k)$, for which the total space
$$Hocolim_{\mathcal{G}(k)}L \longrightarrow N(\mathcal{G}(k))$$
is easily seen to be isomorphic to the projection
$$N(\mathcal{G}(k\oplus M)) \longrightarrow N(\mathcal{G}(k)).$$
The conclusion is that the natural projection
$$\mathbf{Vect}_{n}(k\oplus M) \longrightarrow \mathbf{Vect}_{n}(k)$$
is equivalent, functorially in $M$, to the total space of the
local system $E \mapsto K(E\otimes_{k}E^{\vee}\otimes_{k}M,1)$ on
$\mathbf{Vect}_{n}(k)$.
Therefore, one finds a natural equivalence of simplicial sets
$$\mathbb{D}er_{E}(\mathbb{R}\mathbf{Loc}_{n}(K),M)\simeq
Map_{SSet/\mathbf{Vect}_{n}(k)}(K,\mathbf{Vect}_{n}(k\oplus M))\simeq $$
$$C^{*}(K,E\otimes_{k}E^{\vee}\otimes_{k}M[1])\simeq Map_{C(k)}(C_{*}(K,E\otimes_{k}E^{\vee})[-1],M).$$
As this equivalence is functorial in $M$, this shows that
$$\mathbb{L}_{\mathbb{R}\mathbf{Loc}_{n}(K),E}\simeq
C_{*}(K,E\otimes_{k}E^{\vee})[-1],$$
as required. \end{proof}

\begin{rmk}
\emph{An important consequence of Prop. \ref{pII-8} is that the
$D^{-}$-stack depends on strictly more than
the fundamental groupoid of $K$. Indeed, the tangent space
of $\mathbb{R}\mathbf{Loc}_{n}(K)$ at a global point corresponding to a local
system $E$ on $K$ is $C^{*}(K,E\otimes_{k}E^{\vee})[1]$, which can be non-trivial
even when $K$ is simply connected. In general, there is a closed
immersion of $D^{-}$-stacks (Prop. \ref{pII-6})
$$it_{0}(\mathbb{R}\mathbf{Loc}_{n}(K)) \longrightarrow \mathbb{R}\mathbf{Loc}_{n}(K),$$
which on the level of tangent spaces induces the natural morphism
$$\tau_{\leq 0}(C^{*}(K,E\otimes_{k}E^{\vee})[1])\longrightarrow
C^{*}(K,E\otimes_{k}E^{\vee})[1]$$
giving an isomorphism on $H^{0}$ and $H^{1}$. This shows that
the $D^{-}$-stack $\mathbb{R}\mathbf{Loc}_{n}(K)$ contains
strictly more information than the usual Artin stack of
local systems on $K$, and does encode some
higher homotopical invariants of $K$.}
\end{rmk}

\subsection{Algebras over an operad}\label{IIder.6.2}

Recall (e.g. from \cite{re,sp}) the notions of operads and algebras over them, as
well as their model structures.
Let $\mathcal{O}$ be an operad in $k-Mod$, the category of
$k$-modules. We assume that for any $n$, the $k$-module
$\mathcal{O}(n)$ is projective. Then, for any $A\in sk-Alg$, we
can consider $\mathcal{O}\otimes_{k}A$, which is an operad in the
symmetric monoidal category $A-Mod_{s}$ of $A$-modules. We can
therefore consider $\mathcal{O}\otimes_{k}A-Alg_{s}$, the category
of $\mathcal{O}$-algebras in $A-Mod_{s}$. According to \cite{hin}
the category $\mathcal{O}\otimes_{k}A-Alg_{s}$ can be endowed with
a natural structure of a $\mathbb{U}$-cofibrantly generated model
category for which equivalences and fibrations are defined on the
underlying objects in $sk-Mod$. For a morphism $A\longrightarrow B$
in $sk-Alg$, there is a Quillen adjunction
$$B\otimes_{A} - : \mathcal{O}\otimes_{k}A-Alg_{s} \longrightarrow
\mathcal{O}\otimes_{k}B-Alg_{s} \qquad
\mathcal{O}\otimes_{k}A-Alg_{s} \longleftarrow
\mathcal{O}\otimes_{k}B-Alg_{s} : F,$$
where $F$ is the forgetful functor. The rule $A \mapsto \mathcal{O}\otimes_{k}A-Alg_{s}$
together with base change functors $B\otimes_{A} - $ is almost a left Quillen
presheaf in the sense of Appendix B, except that the associativity of composition
of base change is only valid up to a natural isomorphism. However, the standard strictification
techniques can be applied in order to replace, up to a natural equivalence, this by a genuine
left Quillen presheaf. We will omit to mention this replacement and will proceed as if
$A \mapsto \mathcal{O}\otimes_{k}A-Alg_{s}$ actually defines a left Quillen presheaf
on $k-D^{-}Aff$.

For any $A\in sk-Alg$, we consider $\mathcal{O}-Alg(A)$, the category of
cofibrant objects in $\mathcal{O}\otimes_{k}A-Alg_{s}$ and equivalences between them.
Finally, we let $\mathcal{O}-Alg_{n}(A)$ be the full subcategory consisting of
objects $B\in \mathcal{O}-Alg(A)$, such that the underlying $A$-module of $B$ is
a vector bundle of rank $n$ (i.e. that $B$ is a strong $A$-module, and
$\pi_{0}(B)$ is a projective $\pi_{0}(A)$-module of rank $n$). The base change functors
clearly preserves the sub-categories $\mathcal{O}-Alg_{n}(A)$, and we get this way
a well defined presheaf of $\mathbb{V}$-small categories
$$\begin{array}{ccc}
k-D^{-}Aff & \longrightarrow & Cat_{\mathbb{V}} \\
A & \mapsto & \mathcal{O}-Alg_{n}(A).
\end{array}$$
Applying the nerve functor we obtain a simplicial presheaf
$$\begin{array}{ccc}
k-D^{-}Aff & \longrightarrow & SSet_{\mathbb{V}} \\
A & \mapsto & N(\mathcal{O}-Alg_{n}(A)).
\end{array}$$
This simplicial presheaf will be denoted by $\mathbf{Alg}_{n}^{\mathcal{O}}$\index{$\mathbf{Alg}_{n}^{\mathcal{O}}$}, and is considered
as an object in $k-D^{-}Aff^{\sim,\textrm{\'e}t}$.

\begin{prop}\label{pII-15}
\begin{enumerate}
\item The object $\mathbf{Alg}_{n}^{\mathcal{O}} \in k-D^{-}Aff^{\sim,\textrm{\'e}t}$
is a $D^{-}$-stack.

\item The $D^{-}$-stack $\mathbf{Alg}_{n}^{\mathcal{O}}$
is $1$-geometric and quasi-compact.

\end{enumerate}
\end{prop}

\begin{proof} $(1)$ The proof relies on the standard argument
based on Cor. \ref{cstrict}, and is left to the reader. \\

$(2)$ We consider the natural morphism of $D^{-}$-stacks
$$p : \mathbf{Alg}_{n}^{\mathcal{O}} \longrightarrow \mathbf{Vect}_{n},$$
defined by forgetting the $\mathcal{O}$-algebra structure. Precisely, it
sends an object $B\in \mathcal{O}-Alg_{n}(A)$ to its underlying
$A$-modules, which is a rank $n$-vector bundle by definition. We are going to prove that
the morphism $p$ is a representable morphism, and this will imply the result
as $\mathbf{Vect}_{n}$ is already known to be $1$-geometric and quasi-compact.
For this, we consider the natural morphism $* \longrightarrow \mathbf{Vect}_{n}$, which
is a smooth $1$-atlas for $\mathbf{Vect}_{n}$, and consider the
homotopy pullback
$$\widetilde{\mathbf{Alg}_{n}^{\mathcal{O}}}:=\mathbf{Alg}_{n}^{\mathcal{O}}\times_{\mathbf{Vect}_{n}}^{h}*.$$
It is enough by Prop. \ref{p14} to show that $\widetilde{\mathbf{Alg}_{n}^{\mathcal{O}}}$
is a representable stack. For this, we use \cite[Thm. 1.1.5]{re} in order to show that
the $D^{-}$-stack $\widetilde{\mathbf{Alg}_{n}^{\mathcal{O}}}$ is isomorphic in
$D^{-}\mathrm{St}(k)$ to $Map(\mathcal{O},\underline{End}(k^{n}))$, defined
as follows. For any $A\in sk-Alg$, we consider
$\underline{End}(A^{n})$, the operad in $A-Mod_{s}$ of endomorphisms
of the object $A^{n}$ (recall that for $M\in A-Mod_{s}$, the operad
$\underline{End}(M)$ is defined by $\underline{End}(M)(n):=\underline{Hom}_{A}(M^{\otimes_{k} n},M)$).
We let $Q\mathcal{O}$ be a cofibrant replacement of the operad
$\mathcal{O}$ in the model category $sk-Mod$. Finally, $Map(\mathcal{O},\underline{End}(k^{n}))$ is defined
as
$$\begin{array}{ccc}
sk-Alg & \longrightarrow & SSet \\
A & \mapsto &
Map(\mathcal{O},\underline{End}(k^{n}))(A):=\underline{Hom}_{Op}(Q\mathcal{O},\underline{End}(A^{n})),
\end{array}$$
where $\underline{Hom}_{Op}$ denotes the simplicial set of morphism in the simplicial category
of operads in $sk-Mod$.
As we said, \cite[Thm. 1.1.5]{re} implies that $\widetilde{\mathbf{Alg}_{n}^{\mathcal{O}}}$ is isomorphic
to $Map(\mathcal{O},\underline{End}(k^{n}))$. Therefore, it remain to show that
$Map(\mathcal{O},\underline{End}(k^{n}))$ is a representable $D^{-}$-stack.

For this, we can write $\mathcal{O}$, up to an equivalence, as the homotopy colimit of free operads
$$\mathcal{O}\simeq Hocolim_{i}\mathcal{O}_{i}.$$
Then, we have
$$Map(\mathcal{O},\underline{End}(k^{n}))\simeq
Holim_{i}Map(\mathcal{O}_{i},\underline{End}(k^{n}))$$
and as representable $D^{-}$-stacks are stable by homotopy limits, we are reduced to the case
where $\mathcal{O}$ is a free operad. This means that there exists an integer $m\geq 0$, such that
for any other operad $\mathcal{O}'$ in $sk-Mod$, we have a natural isomorphism of sets
$$Hom_{Op}(\mathcal{O},\mathcal{O}')\simeq Hom(k,\mathcal{O}'(m)).$$
In particular, we find that the $D^{-}$-stack
$Map(\mathcal{O},\underline{End}(k^{n}))$ is isomorphic to the $D^{-}$-stack sending
$A\in sk-Alg$ to the simplicial set
$\underline{Hom}_{A-Mod_{s}}((A^{n})^{\otimes_{k} m},A)$, of morphisms
from the $A$-module $(A^{n})^{\otimes_{k} m}$ to $A$. But this last $D^{-}$-stack
is clearly isomorphic to
$\mathbb{R}\underline{Spec}\, B$, where $B$ is a the free commutative simplicial
$k$-algebra on $k^{n^{m}}$, or in other words a polynomial algebra over $k$ with
$n^{m}$ variables. \end{proof}

We still fix an operad $\mathcal{O}$ in $k-Mod$ (again requiring $\mathcal{O}(m)$ to be a projective
$k$-module for any $m$), and we let
$B$ be an $\mathcal{O}$-algebra in $k-Mod$, such that
$B$ is projective of rank $n$ as a $k$-module. This defines a well defined
morphism of stacks
$$B : * \longrightarrow \mathbf{Alg}_{n}^{\mathcal{O}}.$$
We are going to describe the cotangent complex
of $\mathbf{Alg}_{n}^{\mathcal{O}}$ at $B$ using the notion of
(derived) derivations for $\mathcal{O}$-algebras.

For this, recall from
\cite{gh} the notion of $\mathcal{O}$-derivations from $B$ and  with coefficients in
a $B$-module. For any $B$-module $M$, one can define
the square zero extension $B\oplus M$ of $B$ by $M$, which
is another $\mathcal{O}$-algebra together with a natural
projection $B\oplus M \longrightarrow M$. The $k$-module of derivations from
$B$ to $M$ is defined by
$$Der_{k}^{\mathcal{O}}(B,M):=Hom_{\mathcal{O}-Alg/B}(B,B\oplus M).$$
In the same way, for any $\mathcal{O}$-algebra $B'$ with a
morphism $B' \longrightarrow B$ we set
$$Der_{k}^{\mathcal{O}}(B',M):=Hom_{\mathcal{O}-Alg/B}(B',B\oplus M).$$
The functor $B' \mapsto Der_{k}^{\mathcal{O}}(B',M)$ can be derived on the left,
to give a functor
$$\mathbb{R}Der_{k}^{\mathcal{O}}(-,M) : \mathrm{Ho}(\mathcal{O}-Alg_{s}/B)^{op} \longrightarrow
\mathrm{Ho}(SSet_{\mathbb{V}}).$$
As shown in \cite{gh}, the functor $M \mapsto
\mathbb{R}Der_{k}^{\mathcal{O}}(B,M)$ is co-represented
by an object $\mathbb{L}^{\mathcal{O}}_{B}\in \mathrm{Ho}(sB-Mod)$, thus there is
a natural isomorphism in $\mathrm{Ho}(SSet_{\mathbb{V}})$
$$\mathbb{R}Der_{k}^{\mathcal{O}}(B,M)\simeq
Map_{sB-Mod}(\mathbb{L}_{B}^{\mathcal{O}},M).$$
The category $sB-Mod$ of simplicial $B$-modules is a
closed model category for which equivalences and
fibrations are detected in $sk-Mod$. Furthermore,  the
category $sB-Mod$ is tensored and co-tensored
over $sk-Mod$ making it into a
$sk-Mod$-model category in the sense of \cite{ho}.
Passing to model categories of spectra, we obtain
a model category $Sp(sB-Mod)$ which is
a $Sp(sk-Mod)$-model category.
We will then set
$$\mathbb{R}\underline{Der}_{k}^{\mathcal{O}}(B,M):=
\mathbb{R}\underline{Hom}_{Sp(sB-Mod)}(\mathbb{L}_{B}^{\mathcal{O}},M)
\in \mathrm{Ho}(Sp(sk-Mod)),$$
where $\mathbb{R}\underline{Hom}_{Sp(sB-Mod)}$ denotes the
$\mathrm{Ho}(Sp(sk-Mod))$-enriched derived $Hom$ of $Sp(sB-Mod)$.
Note that $\mathrm{Ho}(sB-Mod)$ is equivalent to
the unbounded derived category of $B$-modules, and as
well that $\mathrm{Ho}(Sp(sk-Mod))$ is equivalent to the
unbounded derived category of $k$-modules.

\begin{prop}\label{pII-16}
With the notations as above, the tangent complex of the $D^{-}$-stack
$\mathbf{Alg}_{n}^{\mathcal{O}}$ at the point $B$ is given by
$$\mathbb{T}_{\mathbf{Alg}_{n}^{\mathcal{O}},B}\simeq
\mathbb{R}\underline{Der}_{k}^{\mathcal{O}}(B,B)[1] \in \mathrm{Ho}(Sp(sk-Mod)).$$
\end{prop}

\begin{proof} We have an isomorphism in $\mathrm{Ho}(Sp(sk-Mod))$
$$\mathbb{T}_{\mathbf{Alg}_{n}^{\mathcal{O}},B}\simeq
\mathbb{T}_{\Omega_{B}\mathbf{Alg}_{n}^{\mathcal{O}},B}[1].$$
Using Prop. \ref{papp2}, we see that the $D^{-}$-stack
$\Omega_{B}\mathbf{Alg}_{n}^{\mathcal{O}}$ can be described as
$$\begin{array}{ccc}
sk-Alg & \longrightarrow & SSet \\
A & \mapsto & Map^{eq}_{\mathcal{O}-Alg_{s}}(B,B\otimes_{k}A)
\end{array}$$
where $Map^{eq}_{\mathcal{O}-Alg_{s}}$ denotes the
mapping space of equivalences in the category of $\mathcal{O}$-algebras
in $sA-Mod$. In particular, we see that for $M\in sk-Mod$, 
the simplicial set $\mathbb{D}er_{\Omega_{B}\mathbf{Alg}_{n}^{\mathcal{O}}}(B,M)$
is naturally equivalent to the homotopy fiber, taken at the identity, of the morphism
$$Map_{\mathcal{O}-Alg_{s}}(B,B\otimes_{k}(k\oplus M))
\longrightarrow Map_{\mathcal{O}-Alg_{s}}(B,B).$$
The $\mathcal{O}$-algebra $B\otimes_{k}(k\oplus M)$ can be identified with
$B\oplus (M\otimes_{k}B)$, the square zero extension of $B$ by $M\otimes_{k}B$,
as defined in \cite{gh}.
Therefore,  by definition of derived derivations this
homotopy fiber is naturally equivalent to
$\mathbb{R}Der_{k}^{\mathcal{O}}(B,M)$. This shows that
$$\mathbb{D}er_{\Omega_{B}\mathbf{Alg}_{n}^{\mathcal{O}}}(B,M)
\simeq \mathbb{R}Der_{k}^{\mathcal{O}}(B,M\otimes_{k}B)\simeq
Map_{Sp(sB-Mod)}(\mathbb{L}_{B}^{\mathcal{O}},M\otimes_{k}B).$$
This implies that for any $N\in sk-Mod$, we have
$$Map_{Sp(sk-Mod)}(N,\mathbb{T}_{\Omega_{B}\mathbf{Alg}_{n}^{\mathcal{O}},B})
\simeq \mathbb{D}er_{\Omega_{B}\mathbf{Alg}_{n}^{\mathcal{O}}}(B,N^{\vee})\simeq
Map_{Sp(sB-Mod)}(\mathbb{L}_{B}^{\mathcal{O}},
N^{\vee}\otimes_{k}B)\simeq $$
$$Map_{Sp(sk-Mod)}(N,
\mathbb{R}\underline{Hom}_{Sp(sB-Mod)}(\mathbb{L}_{B}^{\mathcal{O}},B))
\simeq Map_{Sp(sk-Mod)}(N,\mathbb{R}\underline{Der}_{k}(B,B)).$$
The Yoneda lemma implies the existence of a natural
isomorphism in $\mathrm{Ho}(Sp(sk-Mod))$
$$\mathbb{T}_{\Omega_{B}\mathbf{Alg}_{n}^{\mathcal{O}},B}\simeq
\mathbb{R}\underline{Der}_{k}(B,B).$$
We thus we have
$$\mathbb{T}_{\mathbf{Alg}_{n}^{\mathcal{O}},B}
\simeq \mathbb{R}\underline{Der}_{k}(B,B)[1].$$
\end{proof}

\begin{rmk}
\emph{The proof of Prop. \ref{pII-16} actually shows that
for any $M\in sk-Mod$ we have
$$\mathbb{R}\underline{Hom}_{k}^{Sp}
(\mathbb{L}_{\mathbf{Alg}_{n}^{\mathcal{O}},B},M)\simeq
\mathbb{R}Der_{k}(B,B\otimes_{k}M)[1].$$}
\end{rmk}

\subsection{Mapping $D^{-}$-stacks}\label{IIder.6.3}

We let $X$ be a stack in $\mathrm{St}(k)$, and
$F$ be an $n$-geometric $D^{-}$-stack. We are going to investigate the
geometricity of the $D^{-}$-stack of morphisms from $i(X)$ to $F$\index{$\mathbf{Map}(X,F)$}.
$$\mathbf{Map}(X,F):=\mathbb{R}_{\textrm{\'e}t}\underline{\mathcal{H}om}(i(X),F)\in D^{-}\mathrm{St}(k).$$
Recall that $\mathbb{R}_{\textrm{\'e}t}\underline{\mathcal{H}om}$ denotes the internal
Hom of the category $D^{-}\mathrm{St}(k)$. 

\begin{thm}\label{pII-17}
With the notations above, we assume the following three conditions are satisfied.
\begin{enumerate}
\item The stack
$$t_{0}(\mathbf{Map}(X,F))\simeq \mathbb{R}_{\textrm{\'e}t}\underline{\mathcal{H}om}(X,t_{0}(F))\in \mathrm{St}(k)$$
is $n$-geometric.
\item The $D^{-}$-stack $\mathbf{Map}(X,F)$ has a cotangent complex.
\item The stack $X$ can be written in $\mathrm{St}(k)$ has a homotopy colimit
$Hocolim_{i}U_{i}$, where $U_{i}$ is a affine scheme, flat over $Spec\, k$.
\end{enumerate}
Then the $D^{-}$-stack $\mathbf{Map}(X,F)$ is $n$-geometric.
\end{thm}

\begin{proof} The only if part is clear. Let us suppose that $\mathbf{Map}(X,F)$
satisfies the three conditions. To prove that it is an $n$-geometric $D^{-}$-stack we
are going to lift an $n$-atlas of $t_{0}(\mathbf{Map}(X,F))$
to an $n$-atlas of $\mathbf{Map}(X,F)$. For this we use the following special case of
J.Lurie's representability criterion, which can be proved using the material of this
paper (see Appendix C).

\begin{thm}{(J. Lurie, see \cite{lu} and 
Appendix $C$)}\label{tII-4}
Let $F$ be a $D^{-}$-stack. The following conditions are equivalent.
\begin{enumerate}
\item $F$ is an $n$-geometric $D^{-}$-stack.
\item $F$ satisfies the following three conditions.
\begin{enumerate}
\item The truncation $t_{0}(F)$ is an Artin 
$(n+1)$-stack.
\item $F$ has an obstruction theory.
\item For any $A\in sk-Alg$, the natural morphism
$$\mathbb{R}F(A) \longrightarrow Holim_{k}\mathbb{R}F(A_{\leq k})$$
is an isomorphism in $\mathrm{Ho}(SSet)$.
\end{enumerate}
\end{enumerate}
\end{thm}

We need to prove that $\mathbf{Map}(X,F)$ satisfies the conditions
$(a)-(c)$ of theorem \ref{tII-4}. Condition $(a)$ is clear by assumption.
The existence of
a cotangent complex for $\mathbf{Map}(X,F)$ is guaranteed by assumption. The fact
that $\mathbf{Map}(X,F)$ is moreover inf-cartesian follows
from the general fact.

\begin{lem}\label{lpII-17}
Let $F$ be a $D^{-}$-stack which is inf-cartesian. Then, for any $D^{-}$-stack $F'$, the $D^{-}$-stack
$\mathbb{R}_{\textrm{\'e}t}\underline{\mathcal{H}om}(F',F)$ is inf-cartesian.
\end{lem}

\begin{proof} Writing $F'$ has a homotopy colimit of
representable $D^{-}$-stacks
$$F'\simeq Hocolim_{i}U_{i},$$
we find
$$\mathbb{R}_{\textrm{\'e}t}\underline{\mathcal{H}om}(F',F)\simeq
Holim_{i}\mathbb{R}_{\textrm{\'e}t}\underline{\mathcal{H}om}(U_{i},F).$$
As being inf-cartesian is clearly stable by homotopy limits, we reduce to the
case where $F'=\mathbb{R}\underline{Spec}\, B$ is a representable $D^{-}$-stack.
Let $A\in sk-Alg$, $M$ be an $A$-module with $\pi_{0}(M)=0$,
and $d\in \pi_{0}(\mathbb{D}er(A,M))$ be a derivation. Then, the commutative square
$$\xymatrix{
\mathbb{R}_{\textrm{\'e}t}\underline{\mathcal{H}om}(F',F)(A\oplus_{d}\Omega M) \ar[r] \ar[d] &
\mathbb{R}_{\textrm{\'e}t}\underline{\mathcal{H}om}(F',F)(A) \ar[d] \\
\mathbb{R}_{\textrm{\'e}t}\underline{\mathcal{H}om}(F',F)(A) \ar[r] &
\mathbb{R}_{\textrm{\'e}t}\underline{\mathcal{H}om}(F',F)(A\oplus M)}$$
is equivalent to the commutative square
$$\xymatrix{
\mathbb{R}F((A\oplus_{d}\Omega M)\otimes_{k}^{\mathbb{L}}B) \ar[r] \ar[d] &
 \mathbb{R}F(A\otimes_{k}^{\mathbb{L}}B)\ar[d] \\
\mathbb{R}F(A\otimes_{k}^{\mathbb{L}}B) \ar[r] &
\mathbb{R}F((A\oplus M)\otimes_{k}^{\mathbb{L}}B).}$$
Using the $F$ is inf-cartesian with respect to
$A\otimes_{k}^{\mathbb{L}}B\in sk-Alg$, and the derivation
$d\otimes_{k}B \in \pi_{0}(\mathbb{D}er(A\otimes_{k}^{\mathbb{L}}B,M\otimes_{k}^{\mathbb{L}}B))$
implies that this last square is homotopy cartesian. This shows that
$\mathbb{R}_{\textrm{\'e}t}\underline{\mathcal{H}om}(F',F)$ is inf-cartesian.
\end{proof}

In order to finish the proof of Thm. \ref{pII-17} it only remain to show that
$\mathbf{Map}(X,F)$ also satisfies the condition $(c)$ of Thm. \ref{tII-4}. For this,
we write $X$ as $Hocolim_{i}U_{i}$ with $U_{i}$ affine and flat over $k$,
and therefore $\mathbf{Map}(X,F)$ can be written
as the homotopy limit $Holim_{i}\mathbf{Map}(U_{i},F)$. In order to check condition $(c)$ we can therefore
assume that $X$ is an affine scheme, flat over $k$. Let us write $X=Spec\, B$, with
$B$ a commutative flat $k$-algebra. Then, for any $A\in sk-Mod$, the morphism
$$\mathbb{R}\mathbf{Map}(X,F)(A) \longrightarrow Holim_{k}\mathbb{R}\mathbf{Map}(X,F)(A_{\leq k})$$
is equivalent to
$$\mathbb{R}F(A\otimes_{k}B) \longrightarrow Holim_{k}\mathbb{R}F((A_{\leq k})\otimes_{k}B).$$
But, as $B$ is flat over $k$, we have $(A_{\leq k})\otimes_{k}B\simeq (A\otimes_{k}B)_{\leq k}$, and therefore
the above morphism is equivalent to
$$\mathbb{R}F(A\otimes_{k}B) \longrightarrow Holim_{k}\mathbb{R}F((A\otimes_{k}B)_{\leq k})$$
and is therefore an equivalence by $(1) \Rightarrow (2)$ of Thm. \ref{tII-4} applied to $F$.
\end{proof}

The following corollaries are direct consequences of Thm. \ref{pII-17}. The only
non trivial part consists of proving the existence of a cotangent complex, which
we will assume in the present version of this work.

\begin{cor}\label{cpII-17}
Let $X$ be a projective and flat scheme over $Spec\, k$, and $Y$ a projective and smooth scheme over
$Spec\, k$. Then, the $D^{-}$-stack
$\mathbf{Map}(i(X),i(Y))$ is a $1$-geometric $D^{-}$-scheme. Furthermore, for any
morphism of schemes $f : X \longrightarrow Y$, the cotangent complex
of $\mathbf{Map}(i(X),i(Y))$ at the point $f$ is given by
$$\mathbb{L}_{\mathbf{Map}(i(X),i(Y)),f}\simeq
(C^{*}(X,f^{*}(T_{Y})))^{\vee},$$
where $T_{Y}$ is the tangent sheaf of $Y \longrightarrow Spec\, k$.
\end{cor}

Let us now suppose that $k=\mathbb{C}$ is the field of complex numbers, and
let $X$ be a smooth and projective variety. We will be interested in the
sheaf $X_{DR}$ of \cite{s2}, defined by
$X_{DR}(A):=A_{red}$, for a commutative $\mathbb{C}$-algebra $A$.
Recall that the stack $\mathcal{M}_{DR}(X)$ is defined as the stack of morphisms
from $X_{DR}$ to $\mathbf{Vect}_{n}$, and is identified with the
stack of flat bundles on $X$ (see \cite{s2}). It is known that
$\mathcal{M}_{DR}(X)$ is an Artin $1$-stack.

\begin{cor}\label{cpII-17'}
The $D^{-}$-stack\index{$\mathbb{R}\mathcal{M}_{DR}(X)$} $\mathbb{R}\mathcal{M}_{DR}(X):=\mathbf{Map}(i(X_{DR}),\mathbf{Vect}_{n})$
is $1$-geometric. For a point $E : * \longrightarrow \mathbb{R}\mathcal{M}_{DR}(X)$, corresponding
to a flat vector bundle $E$ on $X$, the cotangent complex of $\mathbb{R}\mathcal{M}_{DR}(X)$
at $E$ is given by
$$\mathbb{L}_{\mathbb{R}\mathcal{M}_{DR}(X),E}\simeq (C^{*}_{DR}(X,E\otimes E^{*}))^{\vee}[-1],$$
where $C^{*}_{DR}(X,E\otimes E^{*})$ is the algebraic de Rham cohomology of $X$ with coefficients
in the flat bundle $E\otimes E^{*}$.
\end{cor}

Corollary \ref{cpII-17'} is only the beginning of the story; in fact we can also produce, in a similar way, 
$\mathbb{R}\mathcal{M}_{Dol}(X)$ and $\mathbb{R}\mathcal{M}_{Hod}(X)$, which are 
derived versions of the moduli stacks of Higgs bundles and $\lambda$-connections
of \cite{s2}, and this would lead us to a \textit{derived version of non-abelian Hodge theory}.
We think this is very interesting research direction because these derived
moduli also encode higher homotopical data in their tangent complexes. We hope
to come back to this topic in a future work.

\chapter{Complicial algebraic geometry}\label{IIunb}

In this chapter we present a second context of application of the general formalism of Part I, 
in which the base model category is $C(k)$, the category of unbounded complexes
over some ring $k$ of characteristic zero. Contrary to the previously considered applications, 
the general notions presented in \S 1.2 does not produce here notions which are very close to 
the usual ones for commutative rings. As a consequence the geometric intuition is here only a very loose 
guide. 

We will present two different HAG contexts over $C(k)$. The first one is very weak in the sense
that it is very easy for a stack to be geometric in this context (these geometric stacks will be called \emph{weakly geometric}). 
The price to pay for this abundance of geometric stacks is that this context 
does not satisfy Artin's conditions and thus there is no good infinitesimal theory. 

The second HAG context we consider 
is a bit closer to the geometric intuition, and satisfies the Artin's conditions so it behaves well 
infinitesimally. 

Both of these ``unbounded'' contexts seem interesting as we are able to produce examples
of geometric stacks which cannot be represented by geometric $D^{-}$-stacks, i.e. by the kind 
of geometric stacks studied in the previous chapter. \\

In this chapter $k$ will be a commutative $\mathbb{Q}$-algebra.

\section{Two HA contexts}\label{IIunb.1}

We let $\mathcal{C}:=C(k)$, the model category of unbounded
complexes of $k$-modules in $\mathbb{U}$. The model structure on
$C(k)$ is the projective one, for which fibrations are epimorphisms and
equivalences are quasi-isomorphisms. The model category $C(k)$
is a symmetric monoidal model category for the tensor product of complexes.
Furthermore, it is well known that our assumptions
\ref{ass-1}, \ref{ass0}, \ref{ass1} and \ref{ass2} are satisfied.

The category $Comm(C(k))$ is the usual model category of unbounded commutative
differential graded algebras over $k$, for which fibrations are epimorphisms and
equivalences are quasi-isomorphisms. The category $Comm(C(k))$ will be denoted
by $k-cdga$, and its objects will simply  be called cdga's. For $A\in k-cdga$, the category
$A-Mod$ is the category of unbounded $A$-dg-modules, again with its natural model
structure. In order to avoid confusion of notations
we will denote the category $A-Mod$ by
$A-Mod_{dg}$. Note that for a usual commutative $k$-algebra
$k'$ we have
$k'-Mod_{dg}=C(k')$, whereas $k'-Mod$ will denote the usual
category of $k'$-modules. Objects in $A-Mod_{dg}$ will be
called \emph{$A$-dg-modules}.

As for the case of simplicial algebras, we will set for any $E\in C(k)$
$$\pi_{i}(E):=H^{-i}(E).$$
When $A$ is a $k$-cdga, the $\mathbb{Z}$-graded $k$-module $\pi_{*}(A)$
has a natural structure of a commutative graded $k$-algebra. In the same
way for $M$ an $A$-dg-module, $\pi_{*}(M)$ becomes a
graded $\pi_{*}(A)$-module. Objects
$A\in k-cdga$ such that $\pi_{i}(A)=0$ for any $i<0$ will be
called $(-1)$-connected. Any $A\in k-cdga$ possesses
a $(-1)$-connected cover
$A' \longrightarrow A$, which is such that
$\pi_{i}(A')\simeq \pi_{i}(A)$ for any $i\geq 0$.

As in the context of derived algebraic geometry we will
denote by $M\mapsto M[1]$ the suspension functor.
In the same way, $M\mapsto M[n]$ is the $n$-times
iterated suspension functor (when $n$ is negative
this means the $n$-times iterated loop functor). \\

The first new feature of complicial algebraic geometry is the existence of two interesting
choices for the subcategory $\mathcal{C}_{0}$, both of them
of particular interest depending on the context.
We let $C(k)_{\leq 0}$ be the full subcategory of
$C(k)$ consisting of complexes $E$ such that
$\pi_{i}(E)=0$ for any $i<0$ (or equivalently, such that 
$H^{i}(E)=0$ for all $i>0$, which explains better the notation). We also set
$k-cdga_{0}$ to be the full subcategory
of $k-cdga$ consisting of $A\in k-cdga$ such that
$\pi_{i}(A)=0$ for any $i\neq 0$. In the same way we
denote by $k-cdga_{\leq 0}$ for the full subcategory of
$k-cdga$ consisting of $A$ such that $\pi_{i}(A)\equiv H^{-i}(A)=0$ for any $i<0$.

\begin{lem}\label{lcont}
\begin{enumerate}
\item The triplet $(C(k),C(k),k-cdga)$ is a HA context.
\item The triplet $(C(k),C(k)_{\leq 0},k-cdga_{0})$
is a HA context.
\end{enumerate}
\end{lem}

\begin{proof} The only non-trivial point is to show that
any $A\in k-cdga_{0}$ is $C(k)_{\leq 0}$-good
in the sense of Def. \ref{dgood}. By definition
$A$ is equivalent to some usual commutative $k$-algebra, so we
can replace $A$ by $k$ itself. We are then left to
prove that the natural functor
$$\mathrm{Ho}(C(k))^{op} \longrightarrow \mathrm{Ho}((C(k)^{op}_{\leq 0})^{\wedge})$$
is fully faithful.

To prove this, we consider the restriction functor
$$i^{*} : \mathrm{Ho}((C(k)^{op})^{\wedge})\longrightarrow \mathrm{Ho}((C(k)^{op}_{\leq 0})^{\wedge})$$
induced by the inclusion $i : C(k)_{\leq 0} \subset C(k)$. We restrict the functor $i^{*}$ to the
full sub-categories of corepresentable objects
$$i^{*} : \mathrm{Ho}((C(k)^{op})^{\wedge})^{\mathrm{corep}}\longrightarrow \mathrm{Ho}((C(k)^{op}_{\leq 0})^{\wedge})^{\mathrm{corep}}.$$
A precision here: we say that
a functor $F : C(k)_{\leq 0} \longrightarrow SSet$ is \emph{corepresentable} 
if it is of the form $$D \mapsto Map_{C(k)}(E,i(D))$$ for some object $E \in C(k)$ (i.e. 
it belongs the essential image of $\mathrm{Ho}(C(k))^{op} \longrightarrow \mathrm{Ho}((C(k)^{op}_{\leq 0})^{\wedge})$.
We claim that this restricted functor $i^{*}$ is an equivalence of categories. Indeed, an inverse $f$ can be 
constructed by sending a functor $F : C(k)_{\leq 0} \longrightarrow SSet$
to the functor
$$\begin{array}{cccc}
f(F) : & C(k) & \longrightarrow & SSet \\
 & D & \mapsto & Holim_{n \geq 0}\Omega^{n}F(D(\leq n)[n]),
\end{array}$$ 
where $D(\leq n)$ is the naive truncation of $D$
defined by $D(\leq n)_{m}=D_{m}$ if $m\leq n$ and
$D(\leq n)_{m}=0$ if $m>n$, and $\Omega^{n}F(D(\leq n)[n])$
is the $n$-fold loop space of the simplicial set $F(D(\leq n)[n])$, 
based at the natural point $*\simeq F(0) \rightarrow F(D(\leq n)[n])$. 
As any $D\in C(k)$ is functorially equivalent to the homotopy limit $Holim_{n}D(\leq n)$, 
it is easy to check that the functor $f$ and $i^{*}$ are inverse to each other. 

To finish the proof, it is enough to notice that there exists a commutative diagram (up to a natural 
isomorphism)
$$\xymatrix{
\mathrm{Ho}(C(k))^{op} \ar[r]^-{\mathbb{R}\underline{h}} \ar[rd] & \mathrm{Ho}((C(k))^{\wedge})^{corep} \ar[d]^-{i^{*}} \\
 & \mathrm{Ho}((C(k)^{op}_{\leq 0})^{\wedge})^{corep}. }$$
The Yoneda embedding $\mathbb{R}\underline{h}$ being fully faithful, this implies that 
the functor 
$$\mathrm{Ho}(C(k))^{op} \longrightarrow \mathrm{Ho}((C(k)^{op}_{\leq 0})^{\wedge})$$ 
is also fully faithful. \end{proof}

\begin{rmk}
\emph{There are objects $A\in k-cdga$ which are \textit{not}
$C(k)_{\leq 0}$-good. For example, we can
take $A$ with $\pi_{*}(A)\simeq k[T,T^{-1}]$
where $T$ is in degree $2$. Then, there is no 
nontrivial $A$-dg-module $M$ such that $M\in C(k)_{\leq 0}$: this 
clearly implies that $A$ can not be $C(k)_{\leq 0}$-good.}

\emph{Another example is given by $A=k[T]$ where $T$ is in degree $-2$, and
$M$ be the $A$-module $A[T^{-1}]$. As there are no morphisms from
$M$ to $A$-modules in $C(k)_{\leq 0}$, then $M$ is sent to zero 
by the functor}
$$\mathrm{Ho}(A-Mod)^{op} \longrightarrow \mathrm{Ho}((A-Mod^{op}_{\leq 0})^{\wedge}).$$
\end{rmk}

\begin{df}\label{unbstrong}
\begin{enumerate}
\item
Let $A\in k-cdga$, and $M$ be an $A$-dg-module. The $A$-dg-module $M$ is \emph{strong}\index{strong!dg-module}
if the natural morphism
$$\pi_{*}(A)\otimes_{\pi_{0}(A)}\pi_{0}(M)\longrightarrow \pi_{*}(M)$$
is an isomorphism.
\item A morphism $A\longrightarrow B$ in $k-cdga$ is \emph{strongly flat}\index{morphism between cdga's!strongly flat}
(resp. \emph{strongly smooth}\index{morphism between cdga's!strongly smooth}, resp. \emph{strongly \'etale}\index{morphism between cdga's!strongly \'etale},
resp. \emph{a strong Zariski open immersion}\index{strong Zariski open immersion!between cdga's}) if $B$ is strong as an $A$-dg-module, and
if the morphism of affine schemes
$$Spec\, \pi_{0}(B) \longrightarrow Spec\, \pi_{0}(A)$$
is flat (resp. smooth, resp. \'etale, resp. a Zariski open immersion).
\end{enumerate}
\end{df}

One of the main difference between derived algebraic geometry and
complicial  algebraic geometry lies in the fact that
the strong notions of flat, smooth, \'etale and Zariski open immersion are not
as easily related to the general notions presented in
\S \ref{partI.2}. We have the following partial comparison result.

\begin{prop}\label{pII-28}
Let $f : A \longrightarrow B$ be a morphism in $k-cdga$.
\begin{enumerate}
\item If $A$ and $B$ are $(-1)$-connected, the morphism
$f$ is smooth (resp. i-smooth, resp. resp. \'etale, resp. a Zariski open immersion) in the sense of Def. \ref{d5}, \ref{d7},
if and only if
$f$ is strongly smooth (resp. strongly \'etale, resp. a strong Zariski open immersion).
\item If the morphism $f$ is strongly flat (resp. strongly smooth,
resp. strongly \'etale, resp. a strong Zariski open immersion), then
it is flat (resp. smooth, resp. \'etale, resp. a Zariski open immersion) in the sense
of Def. \ref{d5}, \ref{d7}.
\end{enumerate}
\end{prop}

\begin{proof} $(1)$ The proof is precisely the same as for
Thm. \ref{tII-1}, as the homotopy theory of $(-1)$-connected
$k-cdga$ is equivalent to
the one of commutative simplicial $k$-algebras, and as this
equivalence preserves cotangent complexes. \\

$(2)$ For flat morphisms there is nothing to prove, as 
all morphisms are flat in the sense of Def. \ref{d5} since the model category $C(k)$ is stable. Let us suppose that
$f : A \longrightarrow B$ is strongly smooth (resp. strongly \'etale, resp. a strong Zariski open immersion).
we can consider the morphism induced on the $(-1)$-connected covers
$$f' : A '\longrightarrow B',$$
where we recall that the $(-1)$-connected cover $A' \longrightarrow A$
induces isomorphims on $\pi_{i}$ for all $i\geq 0$, and is such that 
$\pi_{i}(A')=0$ for all $i<0$. As the morphism $f$ is
strongly flat that the square
$$\xymatrix{
A' \ar[r]^-{f'} \ar[d] & B' \ar[d] \\
A \ar[r]_-{f} & B}$$
is homotopy co-cartesian. Therefore, as our notions of smooth, \'etale and Zariski open
immersion are stable by homotopy push-out, it is enough to show that
$f'$ is smooth (resp. \'etale, resp. a Zariski open immersion). But this follows from
$(1)$. \end{proof}

\begin{ex}\label{exudag}
\emph{Before going further into complicial algebraic geometry we would like to 
present an example illustrating the difference between the strong notions of flat, smooth, \'etale and Zariski open immersion 
and the general notions presented in
\S \ref{partI.2}, showing in particular that proposition \ref{pII-28} $(2)$ does not have a converse.}

\emph{Let $A$ be any commutative $k$-algebra, $X=Spec\, A$ the corresponding affine scheme,
and $U \subset X$ be a quasi-compact open subscheme. It is easy to see that 
there exists a perfect complex of $A$-modules $K$, such that 
$U$ is the open subscheme of $X$ on which $K$ is acyclic. By Prop. \ref{psupp} there
exists then a morphism $A \longrightarrow A_{K}$ in $k-cdga$ which is a Zariski open immersion.
Moreover, the universal property of $A \longrightarrow A_{K}$ shows that  
if $A_{K}$ is cohomologically concentrated in degree $0$ then 
$U$ is affine and we have $U\simeq Spec\, \pi_{0}(A_{K})$. Therefore, as soon 
as $U$ is not affine, $A_{K}$ cannot be concentrated in degree $0$. The morphism
$A \longrightarrow A_{K}$ is thus a Zariski open immersion, and therefore 
\'etale, but is not a strong morphism.  }

\emph{This example also shows that if the scheme $U$ is considered as 
a scheme over $C(k)$ (see \S \ref{IIunb.5.1}), then $U$ is 
equivalent to $\mathbb{R}\underline{Spec}\, A_{K}$, and thus
is affine as a stack over $C(k)$, even though $U$ is not necessarily an affine
subscheme of $X$ over $k$.} 
\end{ex}

\bigskip

The opposite model category $k-cdga^{op}$ will be denoted
by $k-DAff$. We will endow it with the strong \'etale model
topology.

\begin{df}\label{dII-25}
A family of morphisms $\{Spec\, A_{i} \longrightarrow Spec\, A\}_{i\in I}$
in $k-DAff$ is a \emph{strong \'etale covering family}\index{strong \'etale covering family!in $k-DAff$}
(or simply \emph{s-\'et covering family}) if it satisfies the following two
conditions.
\begin{enumerate}
\item Each morphism $A \longrightarrow A_{i}$
is strongly \'etale.
\item There exists a finite subset $J\subset I$ such that the
family $\{A \longrightarrow A_{i}\}_{i\in J}$ is a
formal covering family in the sense of \ref{dcov}.
\end{enumerate}
\end{df}

Using the definition of strong \'etale morphisms, we immediately check
that a family of morphisms $\{Spec\, A_{i} \longrightarrow Spec\, A\}_{i\in I}$
in $k-DAff$ is a s-\'et covering family if and only if there exists a finite
sub-set $J\subset I$ satisfying the following two conditions.

\begin{itemize}

\item For all $i\in I$, the natural morphism
$$\pi_{*}(A)\otimes_{\pi_{0}(A)}\pi_{0}(A_{i}) \longrightarrow
\pi_{*}(A_{i})$$
is an isomorphism.

\item The morphism of affine schemes
$$\coprod_{i\in J} Spec\, \pi_{0}(A_{i}) \longrightarrow Spec\, \pi_{0}(A)$$
is \'etale and surjective.

\end{itemize}

\begin{lem}\label{lII-26}
The s-\'et covering families define a model topology on $k-DAff$, which satisfies
assumption \ref{ass5}.
\end{lem}

\begin{proof} The same as for Lem. \ref{lII-5}. \end{proof}

The model topology s-\'et gives rise to a model category
of stacks $k-DAff^{\sim,\textrm{s-\'et}}$.

\begin{df}\label{dII-26}
\begin{enumerate}
\item A \emph{$D$-stack}\index{$D$-stack} is an object
$F\in k-DAff^{\sim,\textrm{s-\'et}}$ which is a stack
in the sense of Def. \ref{dstack}.
\item
The \emph{model category of $D$-stacks} is
$k-DAff^{\sim,\textrm{s-\'et}}$\index{$k-DAff^{\sim,\textrm{s-\'et}}$}, and its homotopy category is simply denoted by $D\mathrm{St}(k)$\index{$D\mathrm{St}(k)$}.
\end{enumerate}
\end{df}

From Prop. \ref{pII-28} we get the following generalization of Cor. \ref{ctII-1}.

\begin{cor}\label{cII-set}
Let $A\in k-cdga$ and $A'\longrightarrow A$ be its
$(-1)$-connected cover.
Let us consider the natural morphisms
$$t_{0}(X):=Spec\, (\pi_{0}A)\longrightarrow X'=Spec\, A'
\qquad  X=Spec\, A \longrightarrow X':=Spec A'.$$
Then, the homotopy base change functors
$$\mathrm{Ho}(k-DAff/X') \longrightarrow \mathrm{Ho}(k-DAff/t_{0}(X))$$
$$\mathrm{Ho}(k-DAff/X') \longrightarrow \mathrm{Ho}(k-DAff/X)$$
induces equivalences between the full sub-categories
of strongly \'etale morphisms.
Furthermore, these equivalences
preserve epimorphisms of stacks.
\end{cor}

\begin{proof} Using Cor. \ref{ctII-1} we see that it is enough to show that
the base change along the connective cover $A' \longrightarrow A$
$$\mathrm{Ho}(k-DAff/X') \longrightarrow \mathrm{Ho}(k-DAff/X)$$
induces an equivalences from the full sub-categories of strongly \'etale
morphism $Y' \rightarrow X'$ to the full subcategory of
\'etale morphisms $Y\rightarrow X$. But an inverse to this functor is given by
sending a strongly \'etale morphism $A \longrightarrow B$ to its
connective cover $A' \longrightarrow B'$. \end{proof}

\section{Weakly geometric $D$-stacks}\label{IIunb.2}

We now let \textbf{P}$_{w}$
be the class of formally perfect morphisms in $k-DAff$,
also called in the present context \emph{weakly smooth morphisms}.

\begin{lem}\label{lII-27}
The class \textbf{P}$_{w}$ of fp morphisms and the s-\'et model topology satisfy assumptions \ref{ass4}.
\end{lem}

\begin{proof} The only non-trivial thing to show is that the notion of
fp morphism is local for the s-\'et topology on the source. As 
strongly \'etale morphisms are also formally \'etale, 
this easily reduces to showing that being a perfect module is
local for the s-\'et model topology. But this follows from
corollary \ref{ct22}. \end{proof}

We can now state that
$(C(k),C(k),k-cdga,\textrm{s-\'et},\mathbf{P}_{w})$
is a HAG context in the sense of Def. \ref{dhag}.
From our general definitions we obtain a
first notion of geometric $D$-stacks.

\begin{df}\label{dII-27}
\begin{enumerate}
\item
A \emph{weakly $n$-geometric $D$-stack}\index{$D$-stack!weakly $n$-geometric}\index{weakly $n$-representable!$D$-stack} is a
a $D$-stack $F\in D\mathrm{St}(k)$ which is
$n$-geometric for the HAG context
$(C(k),C(k),k-cdga,\textrm{s-\'et},\mathbf{P}_{w})$.
\item
A \emph{weakly $n$-representable morphism of $D$-stacks}\index{weakly $n$-representable!morphism of $D$-stacks}
is an $n$-representable morphism for the
HAG context $(C(k),C(k),k-cdga,\textrm{s-\'et},\mathbf{P}_{w})$.
\end{enumerate}
\end{df}

In the context $(C(k),C(k),k-cdga,\textrm{s-\'et},\mathbf{P}_{w})$ the
i-smooth morphisms are the formally etale morphisms, and therefore
Artin's conditions of Def. \ref{d25} can not be satisfied.
There are actually many interesting weakly geometric $D$-stacks which
do not have cotangent complexes, as we will see.

\section{Examples of weakly geometric $D$-stacks}\label{IIunb.3}

\subsection{Perfect modules}\label{IIunb.3.1}

We consider the $D$-stack $\mathbf{Perf}$, as defined
in Def. \ref{d30'}. As we have seen during the proof of
Lem. \ref{lII-27} the notion of being perfect is
local for the s-\'et topology. In particular,
for any $A\in k-cdga$, the simplicial set
$\mathbf{Perf}(A)$ is naturally equivalent to the
nerve of the category of equivalences between
perfect $A$-dg-modules.

\begin{prop}\label{pperf}
\begin{enumerate}
\item The $D$-stack $\mathbf{Perf}$ is categorically
locally of finite presentation in the sense
of Def. \ref{d23}.
\item
The $D$-stack $\mathbf{Perf}$ is weakly $1$-geometric. Furthermore
its diagonal is $(-1)$-representable.

\end{enumerate}
\end{prop}

\begin{proof} $(1)$ Let $A=Colim_{i\in I}A_{i}$ be a filtered
colimit of objects in $k-cdga$. We need to prove that the
morphism
$$Colim_{i}\mathbf{Perf}(A_{i}) \longrightarrow \mathbf{Perf}(A)$$
is an equivalence. Let $j\in I$, and two perfect $A_{j}$-dg-modules $P$ and $Q$.
Then, we have
$$Map_{A-Mod}(A\otimes_{A_{j}}^{\mathbb{L}}P,A\otimes_{A_{j}}^{\mathbb{L}}Q)
\simeq
Map_{A_{j}-Mod}(P,A\otimes_{A_{j}}^{\mathbb{L}}Q)\simeq $$
$$ \simeq Map_{A_{j}-Mod}(P,Colim_{i\in j/J}A_{i}\otimes_{A_{j}}^{\mathbb{L}}Q).$$
As $P$ is perfect it is finitely presented in the sense of Def. \ref{d3}, and thus
we get
$$Map_{A-Mod}(A\otimes_{A_{j}}^{\mathbb{L}}P,A\otimes_{A_{j}}^{\mathbb{L}}Q)
\simeq Colim_{i\in j/I}Map_{A_{i}-Mod}(A_{i}\otimes_{A_{j}}^{\mathbb{L}}P,
A_{i}\otimes_{A_{j}}^{\mathbb{L}}Q).$$
Furthermore, as during the proof of Prop. \ref{pautperf} we can show that
the same formula holds when $Map$ is replaced by the
sub-simplicial set $Map^{eq}$ of equivalences. Invoking the relations
between mapping spaces and loop spaces of nerves of model categories
(see Appendix B), this clearly implies that
the morphism
$$Colim_{i}\mathbf{Perf}(A_{i}) \longrightarrow \mathbf{Perf}(A)$$
induces isomorphisms on $\pi_{i}$ for $i>0$ and an injective morphism
on $\pi_{0}$.

It only remains to show that for any perfect
$A$-dg-module $P$, there exists $j\in I$ and
a perfect $A_{j}$-dg-module $Q$ such that $P\simeq A\otimes_{A_{j}}^{\mathbb{L}}Q$.
For this, we use that perfect modules are precisely the retract of
finite cell modules. By construction, it is clear that a finite
cell $A$-dg-module is defined over some $A_{j}$. We can therefore
write $P$ as a direct factor of $A\otimes_{A_{j}}^{\mathbb{L}}Q$
for some $j\in I$ and $Q$ a perfect $A_{j}$-dg-module. The direct
factor $P$ is then determined by a projector
$$p\in [A\otimes_{A_{j}}^{\mathbb{L}}Q,A\otimes_{A_{j}}^{\mathbb{L}}Q]\simeq
Colim_{i\in j/I}[A_{i}\otimes_{A_{j}}^{\mathbb{L}}Q,A_{i}\otimes_{A_{j}}^{\mathbb{L}}Q].$$
Thus, $p$ defines a projector $p_{i}$ in
$[A_{i}\otimes_{A_{j}}^{\mathbb{L}}Q,A_{i}\otimes_{A_{j}}^{\mathbb{L}}Q]$
for some $i$, corresponding to a direct
factor $P_{i}$ of $A_{i}\otimes_{A_{j}}^{\mathbb{L}}Q$. Clearly, we have
$$P\simeq A\otimes_{A_{i}}^{\mathbb{L}}P_{i}.$$

$(2)$ We start by showing that
the diagonal of $\mathbf{Perf}$ is $(-1)$-representable. In other words, we
need to prove that for any $A\in k-cdga$, and any
two perfect $A$-dg-modules $M$ and $N$, the $D$-stack
$$\begin{array}{cccc}
Eq(M,N) :  & A-cdga & \longrightarrow & SSet_{\mathbb{V}} \\
 & (A\rightarrow B) & \mapsto & Map_{B-Mod}^{eq}(M\otimes_{A}^{\mathbb{L}}B,
N\otimes_{A}^{\mathbb{L}}B),
\end{array}$$
(where $Map^{eq}$ denotes the sub-simplicial set of $Map$ consisting
of equivalences) is a representable $D$-stack. This is true, and
the proof is exactly the same as the proof of Prop. \ref{pautperf}.

To construct a $1$-atlas, let us chose a $\mathbb{U}$-small
set $\mathcal{F}$ of representative for the isomorphism classes
of finitely presented objects in $\mathrm{Ho}(k-cdga)$. Then, for any
$A\in \mathcal{F}$, let $\mathcal{M}_{A}$ be a $\mathbb{U}$-small
set of representative for the isomorphism
classes of perfect objects in $\mathrm{Ho}(A-Mod_{dg})$. The fact that
these $\mathbb{U}$-small sets $\mathcal{F}$ and $\mathcal{M}$ exists
follows from the fact finitely presented objects are
retracts of finite cell objects (see Cor. \ref{cI-cell}).

We consider the natural morphism
$$p : U:=\coprod_{A\in \mathcal{F}}\coprod_{P\in \mathcal{M}_{A}}
U_{A,P}:=\mathbb{R}\underline{Spec}\, A \longrightarrow
\mathbf{Perf}.$$
By construction the morphism $p$ is an epimorphism of stacks. Furthermore,
as the diagonal of $\mathbf{Perf}$ is $(-1)$-representable
each morphism $U_{A,P}=\mathbb{R}\underline{Spec}\, A \longrightarrow
\mathbf{Perf}$ is $(-1)$-representable. Finally,
as $U_{A,P}=\mathbb{R}\underline{Spec}\, A$ and $\mathbf{Perf}$ are both
categorically locally of finite presentation, we see that for any
$B\in k-cdga$ and any $Y:=\mathbb{R}\underline{Spec}\, B \longrightarrow
\mathbf{Perf}$, the morphism
$$U_{A,P}\times_{\mathbf{Perf}}^{h}Y \longrightarrow Y$$
is a finitely presented morphism between representable
$D$-stacks. In particular, it is a perfect morphism. We therefore
conclude that
$$p : U:=\coprod_{A\in \mathcal{F}}\coprod_{P\in \mathcal{M}_{A}}
U_{A,P}:=\mathbb{R}\underline{Spec}\, A \longrightarrow
\mathbf{Perf}$$
is a $1$-atlas for $\mathbf{Perf}$. This finishes the proof that
$\mathbf{Perf}$ is weakly $1$-geometric.
\end{proof}

\subsection{The $D$-stacks of dg-algebras and dg-categories}\label{IIunb.3.2}

We
consider for any $A\in k-cdga$ the model category $A-dga$, of associative and unital
$A$-algebras (in $\mathbb{U}$). The model structure on $A-dga$ is the usual one for which
equivalences are quasi-isomorphisms and fibrations are epimorphisms.
We consider $Ass(A)$ to be the subcategory
of $A-dga$ consisting of equivalences between
objects $B\in A-dga$ satisfying the following two
conditions.
\begin{itemize}
\item The object $B$ is cofibrant in $A-dga$.
\item The underlying $A$-dg-module of $B$ is perfect.
\end{itemize}

For a morphism $A \longrightarrow A'$ in $k-cdga$, we have a
base change functor
$$A'\otimes_{A} - : Ass(A) \longrightarrow Ass(A'),$$
making $A\mapsto Ass(A)$ into a pseudo-functor
on $k-cdga$. Using the usual strictification procedure, and
passing to the nerve we obtain a simplicial presheaf
$$\begin{array}{cccc}
\mathbf{Ass} : & k-DAff & \longrightarrow & SSet_{\mathbb{V}} \\
 & A & \mapsto & N(Ass(A)).
\end{array}$$

\begin{prop}\label{pass}
\begin{enumerate}
\item
The simplicial presheaf $\mathbf{Ass}$\index{$\mathbf{Ass}$} is a $D$-stack.
\item
The natural projection $\mathbf{Ass} \longrightarrow
\mathbf{Perf}$, which forget the algebra structure,
is $(-1)$-representable.
\end{enumerate}
\end{prop}

\begin{proof} $(1)$ This is exactly the same proof as
Thm. \ref{t2}. \\

$(2)$ Let $A\in k-cdga$ and $\mathbb{R}\underline{Spec}\, A \longrightarrow \mathbf{Perf}$ be a point corresponding to
a perfect $A$-dg-module $E$. We denote by
$\widetilde{\mathbf{Ass}_{E}}$ the homotopy fiber
taken at $E$ of the morphism $\mathbf{Ass} \longrightarrow
\mathbf{Perf}$.

\begin{lem}
The $D$-stack $\widetilde{\mathbf{Ass}_{E}}$ is representable.
\end{lem}

\begin{proof} This is the same argument as for Prop. \ref{pII-15}. We see
using \cite[Thm. 1.1.5]{re} that $\widetilde{\mathbf{Ass}_{E}}$ is the $D$-stack (over
$Spec\, A$)
$Map(\mathcal{A}ss,\mathbb{R}\underline{End}(E))$, of morphisms from
the associative operad to the the (derived) endomorphism operad of
$E$, defined the same way as in the proof of \ref{pII-15}. Again the same argument
as for \ref{pII-15}, consisting
of writing $\mathcal{A}ss$ as a homotopy colimit of free operads, reduces the statement of the
lemma to prove that for a perfect $A$-dg-module $K$ of, the $D$-stack
$$\begin{array}{ccc}
A-cdga & \longrightarrow & SSet \\
A' & \mapsto & Map_{(C(k))}(K,A')
\end{array}$$
is representable. But this is true as it is equivalent to
$\mathbb{R}\underline{Spec}\, B$, where $B$ is the (derived) free $A$-cdga
over $K$. \end{proof}

The previous lemma shows that $\mathbf{Ass} \longrightarrow
\mathbf{Perf}$ is a $(-1)$-representable morphism, and finishes
the proof of Prop. \ref{pass}. \end{proof}

\begin{cor}\label{cpass}
The $D$-stack $\mathbf{Ass}$ is
weakly $1$-geometric.
\end{cor}

\begin{proof} Follows from Prop. \ref{pperf} and
Prop. \ref{pass}. \end{proof}

We now consider a slight modification of
$\mathbf{Ass}$, by considering
dg-algebras as
dg-categories with only one objects. For this, let
$A\in k-cdga$. Recall that a $A$-dg-category is by definition a
category enriched over the symmetric monoidal category $A-Mod_{dg}$, of $A$-dg-modules.
More precisely, a $A$-dg-category $\mathcal{D}$ consists of the following data
\begin{itemize}
\item A set of objects $Ob(\mathcal{D})$.
\item For any pair of objects $(x,y)$ in $Ob(\mathcal{D})$
an $A$-dg-module $\mathcal{D}(x,y)$.
\item For any triple of objects $(x,y,z)$ in $Ob(\mathcal{D})$ a composition morphism
$$\mathcal{D}(x,y)\otimes_{A}\mathcal{D}(y,z) \longrightarrow \mathcal{D}(x,z),$$
which satisfies the usual unital and associativity conditions.
\end{itemize}

The $A$-dg-categories (in the universe $\mathbb{U}$) form a category
$A-dgCat$, with the obvious notion of morphisms. For an $A-dg$-category
$\mathcal{D}$, we can form a category $\pi_{0}(\mathcal{D})$, sometimes
called the \emph{homotopy category of} $\mathcal{D}$, whose objects are
the same as $\mathcal{D}$ and for which morphisms from $x$ to $y$ is the set
$\pi_{0}(\mathcal{D}(x,y))$ (with the obvious induced compositions).
The construction $\mathcal{D} \mapsto \pi_{0}(\mathcal{D})$ defines
a functor from $A-dgCat$ to the category of $\mathbb{U}$-small categories.
Recall that a
morphism $f : \mathcal{D} \longrightarrow \mathcal{E}$ is then called
a \emph{quasi-equivalence} (or simply an \emph{equivalence})
if it satisfies the following two conditions.
\begin{itemize}
\item For any pair of objects $(x,y)$ in $\mathcal{D}$ the induced morphism
$$f_{x,y} : \mathcal{D}(x,y) \longrightarrow
\mathcal{E}(f(x),f(y))$$
is an equivalence in $A-Mod_{dg}$.

\item The induced functor
$$\pi_{0}(f) : \pi_{0}(\mathcal{D}) \longrightarrow \pi_{0}(\mathcal{E})$$
is an equivalence of categories.
\end{itemize}

We also define a notion of fibration, as the morphisms
$f : \mathcal{D} \longrightarrow \mathcal{E}$ satisfying the following two
conditions.
\begin{itemize}
\item For any pair of objects $(x,y)$ in $\mathcal{D}$ the induced morphism
$$f_{x,y} : \mathcal{D}(x,y) \longrightarrow
\mathcal{E}(f(x),f(y))$$
is a fibration in $A-Mod_{dg}$.

\item For any object $x$ in $\mathcal{D}$, and any
isomorphism $v : f(x) \rightarrow z$ in $\pi_{0}(\mathcal{E})$,
there exists an isomorphism $u : x \rightarrow y$ in
$\pi_{0}(\mathcal{D})$ such that
$\pi_{0}(f)(u)=v$.
\end{itemize}

With these notions of fibrations and equivalences, the
category $A-dgCat$ is a model category. This is
proved in \cite{tab} when $A$ is a commutative ring. The general
case of categories enriched in a well behaved
monoidal model category has been worked out recently by J. Tapia (private communication).

For $A\in k-cdga$, we denote by $Cat_{*}(A)$ the category of
equivalences between objects $\mathcal{D}\in A-dgCat$ satisfying the
following two conditions.
\begin{itemize}
\item For any two objects
$x$ and $y$ in $\mathcal{D}$ the $A$-dg-module $\mathcal{D}(x,y)$ is
perfect and cofibrant in $A-Mod_{dg}$.
\item The category $\pi_{0}(\mathcal{D})$ possesses a unique
object up to isomorphism.
\end{itemize}

For a morphism $A \longrightarrow A'$ in $k-cdga$, we have
a base change functor
$$-\otimes_{A}A' : Cat_{*}(A) \longrightarrow Cat_{*}(A')$$
obtained by the formula
$$Ob(\mathcal{D}\otimes_{A}A'):=Ob(\mathcal{D}) \qquad
(\mathcal{D}\otimes_{A}A')(x,y):=\mathcal{D}(x,y)\otimes_{A}A'.$$
This makes $A \mapsto Cat_{*}(A)$ into a pseudo-functor
on $k-cdga$. Strictifying and applying the nerve construction we get
a well defined simplicial presheaf\index{$\mathbf{Cat}_{*}$}
$$\begin{array}{cccc}
\mathbf{Cat}_{*} : & k-cdga & \longrightarrow & SSet_{\mathbb{V}} \\
 & A & \mapsto & N(Cat_{*}(A)).
\end{array}$$

It is worth mentioning that $\mathbf{Cat}_{*}$ is not a stack, since
there are non trivial twisted forms of objects in $Cat_{*}(k)$ with respect to  
the \'etale topology on $k$. These twisted forms can be
interpreted as certain \emph{stacks in dg-categories} having \textit{locally}
a unique object up to equivalences, but they might have
either no global objects or several nonequivalent global objects. 
We will not explicitly describe the stack associated to
$\mathbf{Cat}_{*}$, as this will be irrelevant for the sequel, and will simply
consider $\mathbf{Cat}_{*}$ as an object in $D\mathrm{St}(k)$. 

There exists a morphism of simplicial presheaves
$$B : \mathbf{Ass} \longrightarrow \mathbf{Cat}_{*},$$
which sends an associative $A$-algebra $C$ to the
$A$-dg-category $BC$, having a unique object $*$ and $C$ as the
endomorphism $A$-algebra of $*$. The morphism $B$ will be considered
as a morphism in
$D\mathrm{St}(k)$.

\begin{prop}\label{pcat}
The morphism $B : \mathbf{Ass} \longrightarrow \mathbf{Cat}_{*}$
is weakly $1$-representable, fp and an epimorphism of $D$-stacks.
\end{prop}

\begin{proof} We will prove a more precise result, giving explicit description
of the homotopy fibers of $B$. For this, we start by some
model category considerations relating associative
dga to dg-categories. Recall that for any $A\in k-cdga$, we have
two model categories, $A-dga$ and $A-dgCat$, of associative
$A$-algebras and $A$-dg-categories. We consider
$\mathbf{1}/A-dgCat$, the model category of dg-categories together 
with a distinguised object. More precisely, 
$\mathbf{1}$ is the $A$-dg-category with a unique object and $A$ as its
endomorphism dg-algebra, and
$\mathbf{1}/A-dgCat$ is the comma model category.
The functor $B : A-dga \longrightarrow \mathbf{1}/A-dgCat$
is a left Quillen functor. Indeed, its right adjoint $\Omega_{*}$, sends
a pointed $A$-dg-category $\mathcal{D}$ to the
$A$-algebra $\mathcal{D}(*,*)$ of endomorphisms of the distinguished
point $*$. Clearly, the adjunction morphism
$\Omega_{*}B\Rightarrow Id$ is an equivalence, and thus the
functor
$$B : \mathrm{Ho}(A-dga) \longrightarrow \mathrm{Ho}(\mathbf{1}/A-dgCat)$$
is fully faithful. As a consequence, we get that for any
$C,C'\in A-dga$ there exists a natural homotopy fiber sequence of simplicial sets
$$Map_{A-dgCat}(\mathbf{1},BC') \longrightarrow Map_{A-dga}(B,B') \longrightarrow
Map_{A-dgCat}(BC,BC').$$

Now, let $A\in k-cgda$ and $B$ be an
associative $A$-algebra, corresponding to a morphism
$$x : X:=\mathbb{R}\underline{Spec}\, A \longrightarrow \mathbf{Ass}.$$
We consider the $D$-stack
$$F:=\mathbf{Ass}\times^{h}_{\mathbf{Cat}_{*}}X.$$
Using the relations between mapping spaces in $A-dga$ and $A-dgCat$ above we
easily see that the $D$-stack $F$ is connected (i.e. $\pi_{0}(F)\simeq *$).
In particular, in order to show that $F$ is weakly $1$-geometric, it is
enough to prove that $\Omega_{B}F$ is a representable $D$-stack. Using again
the homotopy fiber sequence of mapping spaces above we see that
the $D$-stack $\Omega_{B}F$ can be described as
$$\begin{array}{cccc}
\Omega_{B}F : & A-cdga & \longrightarrow & SSet_{\mathbb{V}} \\
 & (A\rightarrow A') & \mapsto & Map_{\mathbf{1}/A'-dgCat}(S^{1}\otimes^{\mathbb{L}}\mathbf{1},B(B\otimes_{A}^{\mathbb{L}}A'))
\end{array}$$
where $S^{1}\otimes^{\mathbb{L}}\mathbf{1}$ is computed in the model category
$A-dgCat$. One can easily check that one has an isomorphism in $\mathbf{1}/A-dgCat$
$$S^{1}\otimes^{\mathbb{L}}\mathbf{1}\simeq B(A[T,T^{-1}]),$$
where $A[T,T^{-1}]:=A\otimes_{k}k[T,T^{-1}]$.
Therefore, there is a natural equivalence
$$Map_{\mathbf{1}/A'-dgCat}(S^{1}\otimes^{\mathbb{L}}\mathbf{1},B(B\otimes_{A}^{\mathbb{L}}A'))\simeq
Map_{A-dga}(A[T,T^{-1}],B\otimes_{A}^{\mathbb{L}}A'),$$
and thus the $D$-stack $\Omega_{B}F$ can also be described by
$$\begin{array}{cccc}
\Omega_{B}F : & A-cdga & \longrightarrow & SSet_{\mathbb{V}} \\
 & (A\rightarrow A') & \mapsto & Map_{A-dga}(A[T,T^{-1}],B\otimes_{A}^{\mathbb{L}}A').
\end{array}$$
The morphism of $k$-algebras $k[T] \longrightarrow k[T,T^{-1}]$ induces
natural morphisms (here $B^{\vee}$ is the dual of the perfect $A$-dg-module $B$)
$$Map_{A-dga}(A[T,T^{-1}],B\otimes_{A}^{\mathbb{L}}A') \longrightarrow
Map_{A-dga}(A[T],B\otimes_{A}^{\mathbb{L}}A')\simeq
Map_{A-Mod}(B^{\vee},A').$$
It is not difficult to check that this gives a morphism of $D$-stacks
$$\Omega_{B}F \longrightarrow \mathbb{R}\underline{Spec}\, \mathbb{L}F(B^{\vee}),$$
where $\mathbb{L}F(B^{\vee})$ is the derived free $A$-cdga over the $A$-dg-module
$B^{\vee}$. Furthermore, applying Prop.
\ref{psupp} this morphism is easily seen to be representable by an open Zariski immersion, showing that
$\Omega_{B}F$ is thus a representable $D$-stack. \end{proof}

\begin{cor}\label{ccat}
The $D$-stack $\mathbf{Cat}_{*}$ is weakly $2$-geometric.
\end{cor}

\begin{proof} This follows immediately from
Cor. \ref{cpass}, Prop. \ref{pcat} and the general criterion of Cor. \ref{cp13}.
\end{proof}

Let $A\in k-cdga$ and $B$ be an associative $A$-algebra corresponding
to a morphism of $D$-stacks
$$B : X:=\mathbb{R}\underline{Spec}\, A \longrightarrow \mathbf{Ass}.$$
We define a $D$-stack $B^{*}$ on $A-cdga$ in the following way
$$\begin{array}{cccc}
B^{*} : & A-cdga & \longrightarrow & SSet_{\mathbb{V}} \\
 & (A\rightarrow A') & \mapsto & Map_{A-dga}(A[T,T^{-1}],B\otimes_{A}^{\mathbb{L}}A').
\end{array}$$
The $D$-stack $B^{*}$ is called the \emph{$D$-stack of invertible
elements in $B$}.
The $D$-stack $B^{*}$ possesses in fact a natural loop stack structure (i.e.
the above functor factors naturally, up to equivalence, through a functor
from $A-cdga$ to the category of loop spaces), induced by the
Hopf algebra structure on $A[T,T^{-1}]$. This loop stack structure is also
the one induced by the natural equivalences
$$B^{*}(A')\simeq \Omega_{*}Map_{A-dgCat}(*,B\otimes_{A}^{\mathbb{L}}A').$$
Delooping gives
another $D$-stack
$$\begin{array}{cccc}
K(B^{*},1) :  & A-cdga & \longrightarrow & SSet_{\mathbb{V}} \\
 & (A\rightarrow A') & \mapsto & K(B^{*}(A'),1).
\end{array}$$

The following corollary is a reformulation of Prop. \ref{pcat} and of its proof.

\begin{cor}\label{ccat2}
Let $A\in k-cdga$ and $B\in \mathbf{Ass}(A)$ corresponding to
an associative $A$-algebra $B$. Then, there is a natural homotopy cartesian square
 of $D$-stacks
$$\xymatrix{
\mathbf{Ass} \ar[r] & \mathbf{Cat}_{*} \\
K(B^{*},1) \ar[r] \ar[u] & \mathbb{R}\underline{Spec}\, A \ar[u]^{B}.}$$
\end{cor}

\section{Geometric $D$-stacks}\label{IIunb.4}

We now switch to the HA context
$(C(k),C(k)_{\leq 0},k-cdga_{0})$. Recall that
$C(k)_{\leq 0}$ is the subcategory of $C(k)$ consisting
of $(-1)$-connected object, and $k-cdga_{0}$ is the
subcategory of $k-cdga$ of objects cohomologically concentrated
in degree $0$.
Within this HA context we let \textbf{P} to
be the class of formally perfect and
formally i-smooth morphisms in $k-DAff$. Morphisms in
\textbf{P} will simply be called \emph{fip-smooth morphisms}.
There are no easy description of them, but
Prop. \ref{pismooth} implies that
a morphism $f : A \longrightarrow B$ be a morphism
is fip-smooth if it
satisfies the following two conditions.
\begin{itemize}
\item The cotangent complex $\mathbb{L}_{B/A}\in \mathrm{Ho}(A-Mod_{dg})$
is perfect.
\item For any $R\in k-cdga_{0}$, any connected module
$M\in R-Mod_{dg}$ and any morphism $B \longrightarrow R$, one has
$$\pi_{0}(\mathbb{L}_{B/A}^{\vee}\otimes_{B}^{\mathbb{L}}M)=0.$$
\end{itemize}
The converse is easily shown to be true if $A$ and $B$ are both
$(-1)$-connected.

\begin{lem}\label{lII-27'}
The class \textbf{P} of fip-smooth morphisms
and the s-\'et model topology satisfy assumptions \ref{ass4}.
\end{lem}

\begin{proof} We see that the only non
trivial part is to show that the notion of fip-smooth
morphism is local for the s-\'et topology. For fp-morphisms this
is Cor. \ref{ct22}. For fi-smooth morphisms this is an
easy consequence of the definitions. \end{proof}

From Lem. \ref{lII-27'} we get a HAG
context $(C(k),C(k)_{\leq 0},k-cdga_{0},s-et,\mathbf{P})$.

\begin{df}\label{dII-27'}
\begin{enumerate}
\item
An \emph{$n$-geometric $D$-stack}\index{$D$-stack!$n$-geometric} is
a $D$-stack $F\in D\mathrm{St}(k)$ which is
$n$-geometric for the HAG context
$(C(k),C(k)_{\leq 0},k-cdga_{0},s-et,\mathbf{P})$.
\item A \emph{strongly $n$-representable morphism
of $D$-stacks}\index{strongly $n$-representable morphism
of $D$-stacks} is an $n$-representable morphism
of $D$-stacks for the HAG context
$(C(k),C(k)_{\leq 0},k-cdga_{0},s-et,\mathbf{P})$.
\end{enumerate}
\end{df}

Note that \textbf{P} $\subset$ \textbf{P}$_{w}$
and therefore that any $n$-geometric
$D$-stack is weakly $n$-geometric.

\begin{lem}\label{lII-27''}
The s-\'et topology and the class \textbf{P} of fip-smooth morphisms
satisfy Artin's conditions of Def. \ref{d25}.
\end{lem}

\begin{proof} This is essentially the same proof as Prop. \ref{pII-4}, and
is even more simple as one uses here the
strongly \'etale topology. \end{proof}

\begin{cor}\label{clII-27''}
Any $n$-geometric $D$-stack possesses
an obstruction theory (relative to the
HA context $(C(k),C(k)_{\leq 0},k-cdga_{0})$).
\end{cor}

\begin{proof} This follows from Lem. \ref{lII-27''}
and Thm. \ref{t1}. \end{proof}

\section{Examples of geometric $D$-stacks}\label{IIunb.5}

\subsection{$D^{-}$-stacks and $D$-stacks}\label{IIunb.5.1}

We consider the normalization functor
$N : sk-Alg \longrightarrow k-cdga$,
sending a simplicial commutative $k$-algebra $A$ to its
normalization $N(A)$, with its induced structure of
commutative differential graded algebra.
The pullback functor gives a Quillen adjunction
$$N_{!} : k-D^{-}Aff^{\sim,et} \longrightarrow k-DAff^{\sim,s-et} \qquad
k-D^{-}Aff^{\sim,et} \longleftarrow k-DAff^{\sim,s-et} : N^{*}.$$

As the functor $N$ is known to be homotopically fully faithful, the
left derived functor
$$j:=\mathbb{L}N_{!} : D^{-}\mathrm{St}(k) \longrightarrow
D\mathrm{St}(k)$$
is fully faithful. We can actually characterize the essential image
of the functor $j$ as consisting of all $D$-stacks $F$
for which for any $A\in k-cdga$ with $(-1)$-connected cover $A'\longrightarrow A$, the
morphism $\mathbb{R}F(A') \longrightarrow \mathbb{R}F(A)$ is an equivalence.
In other words, for any $F\in k-D^{-}Aff^{\sim,et}$, and any $A\in k-cdga$, we have
$$\mathbb{R}j(F)(A)\simeq \mathbb{R}F(D(A')),$$
where $D(A')\in sk-Alg$ is a denormalization of $A'$, the
$(-1)$-connected cover of $A$.

\begin{prop}\label{DvsD-}
\begin{enumerate}
\item For any $A\in sk-Alg$, we have
$$j(\mathbb{R}\underline{Spec}\, A)\simeq
\mathbb{R}\underline{Spec}\, N(A).$$
\item The functor $j$ commutes with homotopy limits and homotopy colimits.
\item
The functor $j$ sends $n$-geometric $D^{-}$-stacks to
$n$-geometric $D$-stacks.
\end{enumerate}
\end{prop}

\begin{proof} This is clear. \end{proof}

The previous proposition shows that any $n$-geometric
$D^{-}$-stack gives rise to an $n$-geometric $D$-stack, and thus
provides us with a lot of examples of those.\\

\subsection{CW-perfect modules}\label{IIunb.5.2}

Let $A\in k-cdga$. We define by induction on $n=b-a$ the notion
of a perfect CW-$A$-dg-module of amplitude contained in $[a,b]$.

\begin{df}\label{dcwmod}
\begin{enumerate}
\item A \emph{perfect CW-$A$-dg-module of amplitude contained
in $[a,a]$} is an $A$-dg-module $M$ isomorphic in
$\mathrm{Ho}(A-Mod_{dg})$ to $P[-a]$, with $P$ a projective and
finitely presented $A$-dg-module (as usual in the sense of definitions \ref{d3} and \ref{d7-}).

\item Assume that the notion of perfect CW-$A$-dg-module of amplitude contained
in $[a,b]$ has been defined for any $a\leq b$ such that  $b-a=n-1$.
A \emph{perfect CW-$A$-dg-module\index{perfect CW-dg-module} of amplitude contained
in $[a,b]$}, with $b-a=n$, is an $A$-dg-module $M$ isomorphic in $\mathrm{Ho}(A-Mod_{dg})$ to the
homotopy cofiber of a morphism
$$P[-a-1] \longrightarrow N,$$
where $P$ is projective and finitely presented, and
$N$ is a perfect CW-$A$-dg-module of amplitude contained
in $[a+1,b]$.

\end{enumerate}
\end{df}

The perfect CW-$A$-dg-modules satisfy the following stability
conditions.

\begin{lem}\label{cwstab}
\begin{enumerate}
\item If $M$ is a perfect CW-$A$-dg-module of amplitude
contained in $[a,b]$, and
$A\longrightarrow A'$ is a morphism in $k-cdga$, then
$A'\otimes_{A}^{\mathbb{L}}M$ is a perfect
CW-$A$-dg-module of amplitude
contained in $[a,b]$.

\item Let $A\in k-cdga_{\leq 0}$ be a $(-1)$-connected
$k-cdga$. Then, any perfect $A$-dg-module is
a perfect CW-$A$-dg-module of amplitude $[a,b]$ for some
integer $a\leq b$.

\item Let $A\in k-cdga_{\leq 0}$ be a $(-1)$-connected
$k-cdga$, and $M\in A-Mod_{dg}$. If there exists
a s-et covering $A\longrightarrow A'$ such that
$A'\otimes_{A}^{\mathbb{L}}M$ is a perfect CW-$A$-dg-module
of amplitude contained in $[a,b]$, then so is
$M$.

\end{enumerate}
\end{lem}

\begin{proof} Only $(2)$ and $(3)$ requires a proof. Furthermore, $(3)$ clearly follows
from $(2)$ and the local nature of perfect modules (see Cor. \ref{ct22}), and
thus it only remains to prove $(2)$. But this is proved in
\cite[III.7]{ekmm}. \end{proof}

We define a sub-$D$-stack $\mathbf{Perf}^{CW}_{[a,b]}\subset
\mathbf{Perf}$, consisting of all
perfect modules locally equivalent to some CW-dg-modules of amplitude
contained in $[a,b]$. Precisely,
for $A\in k-cdga$, $\mathbf{Perf}^{CW}_{[a,b]}(A)$\index{$\mathbf{Perf}^{CW}_{[a,b]}$} is the
sub-simplicial set
of $\mathbf{Perf}(A)$ which is the union of all
connected components corresponding to $A$-dg-modules
$M$ such that there is an s-\'et covering
$A\longrightarrow A'$ with  $A'\otimes_{A}^{\mathbb{L}}M$
 a perfect CW-$A'$-dg-module of amplitude contained
in $[a,b]$.

\begin{prop}\label{pcwmod}
The $D$-stack
$\mathbf{Perf}^{CW}_{[a,b]}$ is $1$-geometric.
\end{prop}

\begin{proof} We proceed by induction
on $n=b-a$. For $n=0$, the $D$-stack
$\mathbf{Perf}^{CW}_{[a,a]}$ is simply
the stack of vector bundles $\mathbf{Vect}$, which is
$1$-geometric by our general result Cor. \ref{cp24}.
Assume that $\mathbf{Perf}^{CW}_{[a,b]}$
is known to be $1$-geometric for $b-a< n$.
Shifting if necessary, we can clearly assume that $a=0$, and thus
that $b=n-1$.

That the diagonal of $\mathbf{Perf}^{CW}_{[a,b]}$ is
$(-1)$-representable comes from the fact that
$\mathbf{Perf}^{CW}_{[a,b]} \longrightarrow
\mathbf{Perf}$ is a monomorphisms and from Prop. \ref{pperf}.
It remains to construct a $1$-atlas for
$\mathbf{Perf}^{CW}_{[a,b]}$.

Let $A\in k-cdga$. We consider the model category
$Mor(A-Mod_{dg})$, whose objects are morphisms in $A-Mod_{dg}$
(and with the usual model category structure, see
e.g. \cite{ho}). We consider the subcategory
$Mor(A-Mod_{dg})'$ consisting cofibrant objects
$u : M \rightarrow N$ in $Mor(A-Mod_{dg})$ (i.e.
$u$ is a cofibration between cofibrant $A$-dg-modules)
such that $M$ is projective of finite presentation
and $N$ is locally a perfect CW-$A$-dg-module of amplidtude
contained in $[0,n-1]$. Morphisms
in $Mor(A-Mod_{dg})'$ are taken to be equivalences in
$Mor(A-Mod_{dg})$. The correspondence $A \mapsto Mor(A-Mod_{dg})'$ defines
a pseudo-functor on $k-cdga$, and after strictification and
passing to the nerve we get this way a simplicial presheaf
$$\begin{array}{cccc}
F : & k-cdga & \longrightarrow & SSet_{\mathbb{V}} \\
 & A & \mapsto & N(Mor(A-Mod_{dg})').
\end{array}$$
There exists a morphism of $D$-stacks
$$F \longrightarrow \mathbf{Vect}\times^{h}\mathbf{Perf}^{CW}_{[0,n-1]}$$
that sends an object $M\rightarrow N$ to $(M,N)$.
This morphism is easily seen to be $(-1)$-representable,
as its fiber at an $A$-point $(M,N)$ is the $D$-stack
of morphisms from $M$ to $N$, or equivalently
is $\mathbb{R}\underline{Spec}\, B$ where
$B$ is the derived free $A$-cdga over $M^{\vee}\otimes_{A}^{\mathbb{L}}N$.
By induction hypothesis we deduce that
the $D$-stack $F$ itself if $1$-geometric.

Finally, there exists a natural morphism of $D$-stacks
$$p : F \longrightarrow \mathbf{Perf}^{CW}_{[-1,n-1]}$$
sending a morphism  $u : M \rightarrow N$ in $Mor(A-Mod_{dg})'$ to
its homotopy cofiber. The morphism $p$ is an epimorphism of $D$-stacks
by definition of being locally a perfect CW-dg-module,
thus it only remains to show that
$p$ is fip-smooth. Indeed, this would imply the existence
of a $1$-atlas for $\mathbf{Perf}^{CW}_{[-1,n-1]}$, and thus
that it is $1$-geometric. Translating we get
that $\mathbf{Perf}^{CW}_{[a,b]}$ is $1$-geometric
for $b-a=n$.

In order to prove that $p$ is fip-smooth, let
$A\in k-cdga$, and $K : \mathbb{R}\underline{Spec}\, A \longrightarrow
\mathbf{Perf}^{CW}_{[-1,n-1]}$ be an $A$-point, corresponding
to a perfect CW-$A$-dg-module of amplitude contained in $[-1,n-1]$.
We consider the homotopy cartesian square
$$\xymatrix{
F \ar[r]^-{p} & \mathbf{Perf}^{CW}_{[-1,n-1]} \\
F' \ar[u] \ar[r]_-{p'} & \mathbb{R}\underline{Spec}\, A. \ar[u]}$$
We need to prove that $p$ is a fip-smooth morphism.
The $D$-stack $F'$ has natural projection
$F' \longrightarrow \mathbf{Vect}$ given
by $F '\longrightarrow F \longrightarrow \mathbf{Vect}$, where
the second morphism sends a morphism $M\rightarrow N$ to
the vector bundle $M$. We therefore get a morphism
of $D$-stacks
$$F' \longrightarrow \mathbf{Vect}\times^{h}
\mathbb{R}\underline{Spec}\, A \longrightarrow
\mathbb{R}\underline{Spec}\, A.$$
As the morphism $\mathbf{Vect} \longrightarrow *$ is smooth, it is
enough to show that the morphism
$$F' \longrightarrow \mathbf{Vect}\times^{h}
\mathbb{R}\underline{Spec}\, A$$
is fip-smooth. For this we consider the homotopy cartesian square
$$\xymatrix{
F' \ar[r] & \mathbf{Vect}\times^{h}
\mathbb{R}\underline{Spec}\, A \\
F_{0}' \ar[u] \ar[r] & \mathbb{R}\underline{Spec}\, A, \ar[u]}$$
where the section $\mathbb{R}\underline{Spec}\, A \longrightarrow
\mathbf{Vect}$ correspond to a trivial rank $r$ vector bundle
$A^{r}$. It only remains to show that the morphism
$F_{0}' \longrightarrow\mathbb{R}\underline{Spec}\, A$
is fip-smooth.
The $D$-stack $F_{0}'$ over $A$ can then be easily described
(using for example our general Cor. \ref{cstrict})
as
$$\begin{array}{cccc}
F_{0}' : & A-cdga & \longrightarrow & SSet_{\mathbb{V}} \\
 & (A\rightarrow A') & \mapsto & Map_{A-Mod}(K[-1],A^{r}).
\end{array}$$
In other words, we can write
$F_{0}'\simeq \mathbb{R}\underline{Spec}\, B$, where
$B$ is the derived free $A-cdga$ over $(K[-1])^{r}$.
But, as $K$ is a perfect CW-$A$-dg-module of amplitude contained
in $[-1,n-1]$, $(K[-1])^{r}$ is a perfect CW $A$-dg-module
of amplitude contained in $[0,n]$. The proposition will then
follow from the general lemma.

Let $A \in k-cdga$, and recall that the forgetful functor $A-cdga \longrightarrow A-Mod_{dg}$ is right Quillen, 
and that its derived left adjoint 
$$\mathrm{Ho}(A-Mod_{dg}) \longrightarrow \mathrm{Ho}(A-cdga)$$
sends, by definition, an $A$-module $E$ to the \emph{derived free $A$-cdga over $E$}.

\begin{lem}\label{lcwsmooth}
Let $A \in k-cdga$, and $K$ be a perfect CW-$A$-dg-module of
amplitude contained in $[0,n]$. Let
$B$ be the derived free $A-cdga$ over $K$, and
$$p : Y:=\mathbb{R}\underline{Spec}\, B \longrightarrow X:=\mathbb{R}\underline{Spec}\, A$$
be the natural projection.
Then $p$ is fip-smooth.
\end{lem}

\begin{proof} Let $k'$ be any commutative $k$-algebra, and
$B \longrightarrow k$ be any morphism in $\mathrm{Ho}(k-cdga)$, corresponding to a point
$x : Spec\, k \longrightarrow X$. Then, we have
an isomorphism in $D(k')$
$$\mathbb{L}_{Y/X,x}\simeq
K\otimes_{A}^{\mathbb{L}}k'.$$
We thus see that $\mathbb{L}_{Y/X,x}$ is
a perfect complex of $k'$-modules of Tor amplitude
concentrated in degrees $[0,n]$. In particular, it is
clear that for any complex $M$ of $k'$-modules
such that $\pi_{i}(M)=0$ for any $i>0$, then we have
$$[\mathbb{L}_{Y/X,x},M]_{D(k')}\simeq
\pi_{0}(\mathbb{L}_{Y/X,x}^{\vee}\otimes_{k'}^{\mathbb{L}}M)\simeq
0.$$
This shows that the $A$-algebra $B$ is fip-smooth, and thus
shows the lemma. \end{proof}
The previous lemma finishes the proof of the proposition.
\end{proof}

A direct consequence of Prop. \ref{pcwmod} and Thm. \ref{t1}
is the following.

\begin{cor}\label{cpcwmod}
The $D$-stack $\mathbf{Perf}^{CW}_{[a,b]}$ has
an obstruction theory
(relative to the HA context $(C(k),C(k)_{\leq 0},k-cdga_{0})$.
For any $A\in k-cdga_{0}$ and any
point $E : X:=\mathbb{R}\underline{Spec}\, A \longrightarrow
\mathbf{Perf}^{CW}_{[a,b]}$ corresponding to
a perfect CW-$A$-dg-module $E$, there are natural isomorphisms
in $\mathrm{Ho}(A-Mod_{dg})$
$$\mathbb{L}_{\mathbf{Perf}^{CW}_{[a,b]},E}\simeq
E^{\vee}\otimes_{A}^{\mathbb{L}}E[-1]$$
$$\mathbb{T}_{\mathbf{Perf}^{CW}_{[a,b]},E}\simeq
E^{\vee}\otimes_{A}^{\mathbb{L}}E[1].$$
\end{cor}

\begin{proof} The first part of the corollary follows
from our Thm. \ref{t1} and Prop. \ref{pcwmod}. Let $A$ and
$E : X:=\mathbb{R}\underline{Spec}\, A \longrightarrow
\mathbf{Perf}^{CW}_{[a,b]}$
as in the statement. We have
$$\mathbb{L}_{\mathbf{Perf}^{CW}_{[a,b]},E}\simeq
\mathbb{L}_{\Omega_{E}\mathbf{Perf}^{CW}_{[a,b]},E}[-1].$$
Moreover, $\Omega_{E}\mathbf{Perf}^{CW}_{[a,b]}\simeq
\mathbb{R}\underline{Aut}(E)$, where
$\mathbb{R}\underline{Aut}(E)$ is the $D$-stack of self-equivalences
of the perfect module $E$ as defined in \S \ref{Iqcoh}. By Prop. \ref{pautperf}
we know that $\mathbb{R}\underline{Aut}(E)$ is representable, and
furthermore that the natural inclusion morphism
$$\mathbb{R}\underline{Aut}(E) \longrightarrow \mathbb{R}\underline{End}(E)$$
is a formally \'etale morphism of representable $D$-stacks. Therefore,
we have
$$\mathbb{L}_{\Omega_{E}\mathbf{Perf}^{CW}_{[a,b]},E}\simeq
\mathbb{L}_{\mathbb{R}\underline{End}(E),Id}.$$
Finally, we have
$$\mathbb{R}\underline{End}(E)\simeq \mathbb{R}\underline{Spec}\, B,$$
where $B$ is the derived free $A$-cdga over $E^{\vee}\otimes_{A}^{\mathbb{L}}E$.
This implies that
$$\mathbb{L}_{\mathbb{R}\underline{End}(E),Id}\simeq
E^{\vee}\otimes_{A}^{\mathbb{L}}E,$$
and by what we have seen that
$$\mathbb{L}_{\mathbf{Perf}^{CW}_{[a,b]},E}\simeq
E^{\vee}\otimes_{A}^{\mathbb{L}}E[-1].$$
\end{proof}

\begin{rmk}\label{rcwmod}
\begin{enumerate}
\item
\emph{It is important to note that, if
$$N^{*} : D\mathrm{St}(k) \longrightarrow D^{-}\mathrm{St}(k)$$
denotes the restriction functor, we have
$$N^{*}(\mathbf{Perf}^{CW}_{[a,b]})\simeq \mathbf{Perf}_{[a,b]},$$
where $\mathbf{Perf}_{[a,b]}$ is the sub-stack of
$\mathbf{Perf} \in D^{-}\mathrm{St}(k)$ consisting of
perfect modules of Tor-amplitude contained in
$[a,b]$. However, the two $D$-stacks $\mathbf{Perf}^{CW}_{[a,b]}$ and
$j(\mathbf{Perf}_{[a,b]})$ are not the same. Indeed, 
for any $A \in k-cdga$ with $A'$ as $(-1)$-connective cover, 
we have
$$j(\mathbf{Perf}_{[a,b]})(A) \simeq \mathbf{Perf}_{[a,b]}(A').$$
In general the natural morphism
$$-\otimes_{A}^{\mathbb{L}}A' : \mathbf{Perf}_{[a,b]}(A') \longrightarrow \mathbf{Perf}^{CW}_{[a,b]}(A)$$
is not an equivalence. For example, let us suppose that there exist
a non zero element $x \in \pi_{-1}(A)\simeq [A,A[1]]$, then the matrix
$$\left(\begin{array}{cc}Id_{A} & x \\ 0 & Id_{A[1]}\end{array}\right)$$
defines an equivalence $A \bigoplus A[1] \simeq A\bigoplus A[1]$
of $A$-dg-modules which is not induced by an equivalence
$A'\bigoplus A'[1] \simeq A'\bigoplus A'[1]$ of $A'$-dg-modules. }
\item \emph{We can also show that the $D^{-}$-stack
$\mathbf{Perf}_{[a,b]} \in D^{-}\mathrm{St}(k)$ is
$(n+1)$-geometric for $n=b-a$. The proof is essentially the
same as for Prop. \ref{pcwmod} and will not be reproduced
here. Of course, the formula for the cotangent
complex remains the same. See \cite{tv} for more details.}
\end{enumerate}
\end{rmk}

\subsection{CW-dg-algebras}\label{IIunb.5.3}

Recall the existence of the following diagram of
$D$-stacks
$$\xymatrix{
 & \mathbf{Ass} \ar[d] \\
\mathbf{Perf}^{CW}_{[a,b]} \ar[r] & \mathbf{Perf}.}$$

\begin{df}\label{dcwass}
The \emph{$D$-stack of CW-dg-algebras of amplitude contained
in $[a,b]$}\index{$\mathbf{Ass}_{[a,b]}^{CW}$} is defined by the following homotopy cartesian square
$$\xymatrix{
\mathbf{Ass}_{[a,b]}^{CW} \ar[r] \ar[d] & \mathbf{Ass} \ar[d] \\
\mathbf{Perf}^{CW}_{[a,b]} \ar[r] & \mathbf{Perf}.
}$$
\end{df}

Let $B$ be an associative $k$-dga, and let
$M$ be a $B$-bi-dg-module (i.e.
a $B\otimes_{k}^{\mathbb{L}}B^{op}$-dg-module). We can form
the square zero extension
$B\oplus M$, which is another associative $k$-dga together with
a natural projection $B\oplus M \longrightarrow B$.
The simplicial set of (derived) derivations
from $B$ to $M$ is then defined as
$$\mathbb{R}Der_{k}(B,M):=Map_{k-dga/B}(B,B\oplus M).$$
The same kind of proof as for Prop. \ref{p1} shows that
there exists an object $\mathbb{L}_{B}^{Ass}\in \mathrm{Ho}(B\otimes_{k}^{\mathbb{L}}B^{op}-Mod_{dg})$
an natural isomorphisms in $\mathrm{Ho}(SSet_{\mathbb{V}})$
$$\mathbb{R}Der_{k}(B,M)\simeq
Map_{B\otimes_{k}^{\mathbb{L}}B^{op}-Mod}(\mathbb{L}_{B}^{Ass},M).$$
We then set
$$\mathbb{R}\underline{Der}_{k}(B,M):=
\mathbb{R}\underline{Hom}_{B\otimes_{k}^{\mathbb{L}}B^{op}-Mod}(\mathbb{L}_{B}^{Ass},M)
\in \mathrm{Ho}(C(k)),$$
where $\mathbb{R}\underline{Hom}_{B\otimes_{k}^{\mathbb{L}}B^{op}-Mod}$ denotes
the $\mathrm{Ho}(C(k))$-enriched derived Hom's of
the $C(k)$-model category $B\otimes_{k}^{\mathbb{L}}B^{op}-Mod_{dg}$.

\begin{cor}\label{ccwass}
The $D$-stack $\mathbf{Ass}^{CW}_{[a,b]}$ is
$1$-geometric. For any global point
$B : Spec\, k \longrightarrow \mathbf{Ass}^{CW}_{[a,b]}$, corresponding to
an associative $k$-dga $B$, one has a natural isomorphism
in $\mathrm{Ho}(C(k))$
$$\mathbb{T}_{\mathbf{Ass}^{CW}_{[a,b]},B}\simeq
\mathbb{R}\underline{Der}_{k}(B,B)[1].$$
\end{cor}

\begin{proof} Using Prop. \ref{pass} we see that the natural
projection
$$\mathbf{Ass}^{CW}_{[a,b]} \longrightarrow
\mathbf{Perf}^{CW}_{[a,b]}$$
is $(-1)$-representable. Therefore, by Prop. \ref{pcwmod} we know
that $\mathbf{Ass}^{CW}_{[a,b]}$ is $1$-geometric. In patricular,
Thm. \ref{t1} implies that it has an obstruction theory
relative to the HA context $(C(k),C(k)_{\leq 0},k-cdga_{0})$.
We thus have
$$\mathbb{T}_{\mathbf{Ass}^{CW}_{[a,b]},B}\simeq
\mathbb{T}_{\Omega_{B}\mathbf{Ass}^{CW}_{[a,b]},B}[1].$$
The identification
$$\mathbb{T}_{\Omega_{B}\mathbf{Ass}^{CW}_{[a,b]},B}\simeq
\mathbb{R}\underline{Der}_{k}(B,B)$$
follows from the exact same argument
as Prop. \ref{pII-16}, using $\mathbb{L}_{B}^{Ass}$ instead
of $\mathbb{L}_{B}^{\mathcal{O}}$. 
\end{proof}

\subsection{The $D$-stack of negative CW-dg-categories}\label{IIunb.5.4}

Recall from \S \ref{IIunb.3.2} the existence of the morphism of $D$-stacks
$$B : \mathbf{Ass} \longrightarrow \mathbf{Cat}_{*},$$
sending an associative $k$-algebra $C$ to the dg-category
$BC$ having a unique object and $C$ as endomorphisms of
this object.

\begin{df}\label{dcwdgcat}
Let $n\leq 0$ be an integer.
The \emph{$D$-stack of CW-dg-categories of amplitude
contained in $[n,0]$} is defined as the full sub-$D$-stack
$\mathbf{Cat}_{*,[n,0]}^{CW}$\index{$\mathbf{Cat}_{*,[n,0]}^{CW}$} of $\mathbf{Cat}_{*}$ consisting
of the essential image of the morphism
$$\mathbf{Ass}^{CW}_{[n,0]} \longrightarrow
\mathbf{Ass} \longrightarrow \mathbf{Cat}_{*}.$$
\end{df}

More precisely, for $A\in k-cdga$, one sets
$\mathbf{Cat}_{*,[n,0]}^{CW}(A)$ to be the sub-simplicial
set of $\mathbf{Cat}_{*}(A)$ consisting of $A$-dg-categories
$\mathcal{D}$ such that for any $*\in Ob(\mathcal{D})$
the $A$-dg-module $\mathcal{D}(*,*)$ is (locally)
a perfect CW-$A$-dg-module of amplitude contained in
$[n,0]$.

Recall that for any associative $k$-dga $B$,
one has a model
category $B\otimes_{k}^{\mathbb{L}}B^{op}-Mod_{dg}$ of
$B$-bi-dg-modules. This model category is naturally
tensored and co-tensored over the symmetric
monoidal model category $C(k)$, making it into
a $C(k)$-model category in the sense of \cite{ho}.
The derived $\mathrm{Ho}(C(k))$-enriched Hom's
of $B\otimes_{k}^{\mathbb{L}}B^{op}-Mod_{dg}$ will then
be denoted by
$\mathbb{R}\underline{Hom}_{B\otimes_{k}^{\mathbb{L}}B^{op}}$.

Finally, for any associative $k$-dga $B$, one sets
$$\mathbb{HH}_{k}(B,B):=\mathbb{R}\underline{Hom}_{B\otimes_{k}^{\mathbb{L}}B^{op}}
(B,B)\in \mathrm{Ho}(C(k)),$$
where $B$ is considered as $B$-bi-dg-module in the obvious way.

\begin{thm}\label{tdg-cat}
\begin{enumerate}
\item The morphism
$$B : \mathbf{Ass}^{CW}_{[n,0]}
\longrightarrow \mathbf{Cat}_{*,[n,0]}^{CW}$$
is a $1$-representable fip-smooth covering of $D$-stacks.
\item The associated $D$-stack to $\mathbf{Cat}^{CW}_{*,[n,0]}$
is a $2$-geometric $D$-stack.
\item If $C$ is an associative $k$-dg-algebra, corresponding to a point
$C : Spec\, k\longrightarrow \mathbf{Ass}_{[n,0]}^{CW}$, then one has
a natural isomorphism in $\mathrm{Ho}(C(k))$
$$\mathbb{T}_{\mathbf{Cat}_{*,[n,0]}^{CW},BC}\simeq \mathbb{HH}_{k}(C,C)[2].$$
\end{enumerate}
\end{thm}

\begin{proof} $(1)$ Using Cor. \ref{ccat2},
it is enough to show that for any $A\in k-cdga$,
and any associative $A$-algebra $C$, which is
a perfect CW-$A$-dg-module of amplitude contained in $[n,0]$, the
morphism of $D$-stack $K(C^{*},1) \longrightarrow \mathbb{R}\underline{Spec}\, A$
is $1$-representable and fip-smooth. For this it is clearly enough to
show that the $D$-stack $C^{*}$ is representable and that the
morphism $C^{*} \longrightarrow \mathbb{R}\underline{Spec}\, A$
is fip-smooth. We already have seen during the proof
of Prop. \ref{pcat} that $C^{*}$ is representable and
a Zariski open sub-$D$-stack
of $\mathbb{R}\underline{Spec}\, F$, where
$F$ is the free $A$-cdga over $C^{\vee}$.
It is therefore enough to see that
$A \longrightarrow F$ is fip-smooth, which
follows from Lem. \ref{lcwsmooth} as by assumption
$C^{\vee}$ is a perfect CW-$A$-dg-module of amplitude
contained in $[0,n]$. \\

$(2)$ Follows from Cor. \ref{ccwass} and Cor. \ref{cp13}. \\

$(3)$ We consider a point $C : Spec\, k \longrightarrow
\mathbf{Ass}^{CW}{[n,0]}$, and the homotopy cartesian square
$$\xymatrix{
\mathbf{Ass}^{CW}_{[n,0]} \ar[r] & \mathbf{Cat}_{*,[n,0]}^{CW} \\
K(C^{*},1) \ar[r] \ar[u] & Spec\, k \ar[u]^{C}.}$$
By $(2)$,  Cor. \ref{ccwass} and Thm. \ref{t1} we know that
all the stacks in the previous square
have an obstruction theory, and thus a cotangent complex
(relative to the HA context $(C(k),C(k)_{\leq 0},k-cdga_{0})$).
Therefore, one finds a homotopy fibration sequence of complexes
of $k$-modules
$$\xymatrix{
\mathbb{T}_{K(C^{*},1),*} \ar[r] & \mathbb{T}_{\mathbf{Ass}^{CW}_{[n,0]},C}
\ar[r] & \mathbb{T}_{\mathbf{Cat}^{CW}_{*,[n,0]},BC}}$$
that can also be rewritten as
$$\xymatrix{
C[1] \ar[r] & \mathbb{R}\underline{Der}_{k}(C,C)[1]
\ar[r] & \mathbb{T}_{\mathbf{Cat}^{CW}_{*,[n,0]},BC}.}$$
The morphism $C \longrightarrow \mathbb{R}\underline{Der}_{k}(C,C)$
can be described in the following way.
The $C$-bi-dg-module $\mathbb{L}_{C}^{Ass}$ can be easily identified
with the homotopy fiber (in the model category of
$C\otimes_{k}^{\mathbb{L}}C^{op}$-dg-modules) of the multiplication morphism
$$C\otimes_{k}^{\mathbb{L}}C^{op} \longrightarrow C.$$
The natural morphism
$\mathbb{L}_{C}^{Ass} \longrightarrow C\otimes_{k}^{\mathbb{L}}C^{op}$
then induces our morphism on the level of derivations
$$C\simeq \mathbb{R}\underline{Hom}_{C\otimes_{k}^{\mathbb{L}}C^{op}-Mod}(
C\otimes_{k}^{\mathbb{L}}C^{op},C) \longrightarrow
\mathbb{R}\underline{Hom}_{C\otimes_{k}^{\mathbb{L}}C^{op}-Mod}(
\mathbb{L}_{C}^{Ass},B) \simeq \mathbb{R}\underline{Der}_{k}(C,C).$$
In particular, we see that there exists a natural
homotopy fiber sequence in $C(k)$
$$\xymatrix{
\mathbb{HH}_{k}(C,C)=
\mathbb{R}\underline{Hom}_{C\otimes_{k}^{\mathbb{L}}C^{op}-Mod}(C,C) \ar[r]
& C \ar[r] & \mathbb{R}\underline{Der}_{k}(C,C).}$$
We deduce that there exists a natural
isomorphism in $\mathrm{Ho}(C(k))$
$$\mathbb{T}_{\mathbf{Cat}^{CW}_{*,[n,0]},BC}[-2]\simeq
\mathbb{HH}_{k}(C,C).$$
\end{proof}

An important corollary of Thm. \ref{tdg-cat} is given by the following
fact. It appears in many places in the literature but we 
know of no references including a proof of it. For this, we recall
that for any commutative $k$-algebra $k'$ we denote by
$Cat_{*}(k')$ the category of
equivalences between $k'$-dg-categories satisfying the
following two conditions.
\begin{itemize}
\item For any two objects
$x$ and $y$ in $\mathcal{D}$ the complex of $k'$-module $\mathcal{D}(x,y)$ is
perfect (and cofibrant in $C(k')$).
\item The category $\pi_{0}(\mathcal{D})$ possesses a unique
object up to isomorphism.
\end{itemize}
We finally let $Cat_{*}^{[n,0]}(k')$ be the full subcategory
of $Cat_{*}(k')$ consisting of objects
$\mathcal{D}$ such that for any two objects
$x$ and $y$ the perfect complex $\mathcal{D}(x,y)$ has
Tor amplitude
contained in $[n,0]$ for some $n\leq 0$.

\begin{cor}\label{ctdg-cat}
Let $\mathcal{D} \in Cat_{*}^{[n,0]}(k')$. Then, the homotopy fiber
$\mathbb{D}ef_{\mathcal{D}}$, taken
at the point $\mathcal{D}$, of the
morphism of simplicial sets
$$N(Cat_{*}^{[n,0]}(k[\epsilon])) \longrightarrow N(Cat_{*}^{[n,0]}(k))$$
is given by
$$\mathbb{D}ef_{\mathcal{D}}\simeq Map_{C(k)}(k,\mathbb{HH}_{k}(\mathcal{D},\mathcal{D})[2]).$$
In particular, we have
$$\pi_{i}(\mathbb{D}ef_{\mathcal{D}})\simeq
\mathbb{HH}_{k}^{2-i}(\mathcal{D},\mathcal{D}).$$
\end{cor}

In the above corollary we have used $\mathbb{HH}_{k}(\mathcal{D},\mathcal{D})$,
the Hochschild complex of a dg-category $\mathcal{D}$. It is defined
the same way as for associative dg-algebras, and when
$\mathcal{D}$ is equivalent to $BC$ for an associative dg-alegbra $C$ we have
$$\mathbb{HH}_{k}(\mathcal{D},\mathcal{D})\simeq
\mathbb{HH}_{k}(C,C).$$

Finally, we would like to mention that restricting to
negatively graded dg-categories seems difficult to avoid if we
want to keep the existence of a cotangent complex.

\begin{cor}\label{ctdg-cat2}
Assume that $k$ is a field.
Let $C\in \mathbf{Ass}^{CW}_{[a,b]}(k)$ be a $k$-point corresponding to
an associative dg-algebra $C$. If we have $H^{i}(C)\neq 0$ for
some $i> 0$ then the $D$-stack
$\mathbf{Cat}_{*}$ does not have a cotangent complex
at the point $BC$.
\end{cor}

\begin{proof} Suppose that $\mathbf{Cat}_{*}$ does  have a cotangent complex
at the point $BC$. Then, as so does the $D$-stack
$\mathbf{Ass}^{CW}_{[a,b]}$ (by Cor. \ref{ccwass}), we see that the
homotopy fiber, taken at the point $BC$, of the morphism $B : \mathbf{Ass} \longrightarrow
\mathbf{Cat}_{*}$ has a cotangent complex at $C$.
This homotopy fiber is $K(C^{*},1)$ (see Cor. \ref{ccat2}),
and thus we would have
$$\mathbb{L}_{K(C^{*},1),C}\simeq C^{\vee}[-1].$$
Let $k[\epsilon_{i-1}]$ be the square zero extension
of $k$ by $k[i-1]$ (i.e. $\epsilon_{i-1}$ is in degree
$-i+1$), for some $i>0$ as in the statement.
We now consider the homotopy fiber sequence
$$\xymatrix{
Map_{C(k)}(C^{\vee}[-1],k[i-1]) \ar[r] & K(C^{*},1)(k[\epsilon_{i-1}]) \ar[r] &
K(C^{*},1)(k).}$$
Considering the long exact sequence in homotopy we find
$$\xymatrix{
\pi_{1}(K(C^{*},1)(k[\epsilon_{i-1}]))=H^{0}(C)^{*}\oplus H^{i-1}(C) \ar[r] &
\pi_{1}(K(C^{*},1)(k))=H^{0}(C)^{*}
\ar[r] & }$$
$$\xymatrix{\pi_{0}(Map_{C(k)}(C^{\vee}[-1],k[i-1]))=H^{i}(C) \ar[r] &
\pi_{0}(K(C^{*},1)(k[\epsilon_{i-1}]))=*.}$$
This shows that the last morphism must be injective, which can not
be the case as soon as $H^{i}(C)\neq 0$. \end{proof}

\begin{rmk}\label{rmkcat}
\begin{enumerate}
\item \emph{Some of the results in Thm. \ref{tdg-cat} were announced
in \cite[Thm. 5.6]{hagdag}. We need to warn the reader that
\cite[Thm. 5.6]{hagdag} is not correct for the description of
$\widetilde{\mathbb{R}Cat_{\mathcal{O}}}$ briefly given
before that theorem
(the same mistake appears in
\cite[Thm. 4.4]{to3}). Indeed, $\widetilde{\mathbb{R}Cat_{\mathcal{O}}}$
would correspond to} isotrivial \emph{deformations
of dg-categories, for which the underlying
complexes of morphisms stays locally constant. Therefore, the tangent
complex of  $\widetilde{\mathbb{R}Cat_{\mathcal{O}}}$ can not be
the full Hochschild complex as stated in \cite[Thm. 5.6]{hagdag}.
Our theorem \ref{tdg-cat} corrects this mistake.}

\item \emph{
We like to consider our Thm. \ref{tdg-cat} $(3)$ and
Cor. \ref{ctdg-cat} as
a possible explanation of the following sentence in
\cite[p. 266]{kos}:} \\``In some sense, the full Hochschild complex controls deformations
of the $A_{\infty}$-category with one object, such that its endomorphism
space is equal to $A$.'' \\

\emph{
Furthermore, our homotopy fibration sequence
$$\xymatrix{
K(C^{*},1) \ar[r] & \mathbf{Ass}^{CW}_{[n,0]} \ar[r] & \mathbf{Cat}_{*,[n,0]}^{CW}}$$
is the geometric global counter-part of the
well known exact triangles of complexes
(see e.g. \cite[p. 59]{ko2})
$$\xymatrix{
C[1] \ar[r] & \mathbb{R}\underline{Der}_{k}(C,C)[2] \ar[r] & \mathbb{HH}_{k}(C,C)[2]
\ar[r]^-{+1} & }$$
as we pass from the former to the latter by taking
tangent complexes at the point $C$. }

\item \emph{We saw in Cor. \ref{ctdg-cat2} that
the full $D$-stack $\mathbf{Cat}_{*}$ can not have
a reasonable infinitesimal theory.
We think it is important to mention that
even Cor. \ref{ctdg-cat} cannot reasonably be true if we
remove the assumption that $\mathcal{D}(x,y)$ is
of Tor-amplitude contained in $[n,0]$ for some
$n\leq 0$. Indeed, for any commutative $k$-algebra
$k'$, the morphism
$$B : \mathbf{Ass}(k') \longrightarrow \mathbf{Cat}_{*}(k')$$
is easily seen to induce an isomorphism on $\pi_{0}$
and a surjection on $\pi_{1}$. From this and the
fact that $\mathbb{T}_{\mathbf{Ass},C}=\mathbb{R}\underline{Der}_{k}(C,C)[1]$,
we easily deduce that the natural morphism
$$\pi_{0}(\mathbb{R}\underline{Der}_{k}(C,C)[1]) \longrightarrow
\pi_{0}(\mathbb{D}ef_{BC})$$
is surjective. Therefore, if Cor. \ref{ctdg-cat} were true for
$\mathcal{D}=BC$, we would have that the morphism
$$H^{1}(\mathbb{R}\underline{Der}_{k}(C,C)) \longrightarrow
\mathbb{HH}^{2}(C,C)$$
is surjective, or equivalently that the natural morphism
$$H^{2}(C) \longrightarrow H^{2}(\mathbb{R}\underline{Der}_{k}(C,C))$$
is injective. But this is not the case in general, as the morphism
$C \longrightarrow \mathbb{R}\underline{Der}_{k}(C,C)$ can be
zero (for example when $C$ is commutative). It is therefore not
strictly correct to state that the Hochschild cohomology of
an associative dg-algebra controls its deformation
as a dg-category, contrary to what appears in several
references (including some of the authors !).  }

\item \emph{For $n=0$, we can easily show that the
restriction $N^{*}(\mathbf{Cat}_{*,[0,0]}^{CW})\in D^{-}\mathrm{St}(k)$ is
a $2$-geometric $D^{-}$-stack in the sense of \S \ref{IIder}. Furthermore, its
truncation $t_{0}N^{*}(\mathbf{Cat}_{*,[0,0]}^{CW}) \in \mathrm{Ho}(k-Aff^{\sim,et})$
is naturally equivalent to the Artin $2$-stack of
$k$-linear categories with one object. The tangent complex
of $N^{*}(\mathbf{Cat}_{*,[0,0]}^{CW})$ at a point $C$ corresponding
to a $k$-alegbra projective of finite type over $k$ is then
the usual Hochschild complex $HH(C,C)$ computing the Hochschild
cohomology of $C$.}

\emph{
Of course $\mathbf{Cat}_{*,[n,0]}^{CW}$ is only a rough
approximation to what should be the stack 
of dg-categories, and in particular we think that $\mathbf{Cat}_{*,[n,0]}^{CW}$ 
is not suited for dealing with dg-categories coming from algebraic geometry.
Let for example $X$ be a smooth and projective variety over a field
$k$; it is known that the the derived category $D_{qcoh}(X)$ is of the form
$\mathrm{Ho}(B-Mod_{dg})$ for some associative $dg$-algebra $B$ which is
perfect as a complex of $k$-modules. It is however very unlikely that
$B$ can be chosen to be concentrated in degrees $[n,0]$ for
some $n\leq 0$ (by construction
$H^{i}(B)=Ext^{i}(E,E)$ for some compact generator
$E\in D_{qcoh}(X)$). So the dg-algebra $B$ when considered as 
a dg-category will \textit{not} define a point in $\mathbf{Cat}_{*,[n,0]}^{CW}$. Another, 
more serious problem comes from the fact that the dg-algebra $B$ 
is not uniquely determined, but is only unique up to Morita equivalence. 
As a consequence, the variety $X$ can deform and $B$ might not follow this deformation
(though another, Morita equivalent, dg-algebra will follow the deformation), and 
thus $X \mapsto B$ will not be a morphism of stacks (even locally around $X$). Therefore, 
it seems very important to consider the stack $\mathbf{Cat}_{*}$ modulo  
Morita equivalences. We also think that passing to Morita equivalences will solve the
problem of the non geometricity of $\mathbf{Cat}_{*}$ mentioned above. 
This direction is currently being investigated by M. Anel.}
\end{enumerate}
\end{rmk}

\chapter{Brave new algebraic geometry}\label{IIbnag}

In this final chapter we briefly present brave new algebraic geometry\footnote{The term ``brave new rings'' was invented by F. Waldhausen to describe structured ring spectra; we have only adapted it to our situation.}, 
i.e. algebraic geometry over ring spectra. We will emphasize
the main differences with derived algebraic geometry, and
the subject will be studied in more details in future works.

As in the case of complicial algebraic geometry, we will present two distinct  
HAG contexts (see Cor. \ref{clII-27bis}) which essentially differs in the choice of the class \textbf{P}.
The first one, where \textbf{P} is chosen to be the class of strongly \'etale morphisms (see Def. \ref{bnagstrong}), is suited for defining \emph{brave new} Deligne-Mumford stacks, and, as all the contexts based on ``strong'' morphisms, 
it is geometrically very close to usual algebraic geometry. The second HAG context, where \textbf{P} is chosen to be the class of fip-smooth morphisms 
(i.e. formally perfect and
formally i-smooth morphisms), is weaker and is similar to the corresponding weak context already presented in complicial algebraic geometry:
it allows to define \emph{brave new} geometric stacks which are not 
Deligne-Mumford. Since the notion of fip-smooth morphism 
for brave new rings behaves differently from 
commutative rings (see Prop. \ref{affineline}), the geometric intuition in this context is once again a bit far from 
standard algebraic geometry (e.g. smooth morphisms are not necessary flat). 
Nonetheless we think that this context is very interesting, as it is not only geometrically reasonable 
(e.g. it satisfies Artin's conditions), but it is also able to ``see'' some of the interesting 
new phenomena arising in the theory of structured ring spectra.

\section{Two HAG contexts}\label{IIbnag.1}

We let $\mathcal{C}:=Sp^{\Sigma}$, the category of symmetric spectra
in $\mathbb{U}$. The model structure we are going to use on
$Sp{\Sigma}$ is the so called positive stable model structure
described in \cite{shi}. This model structure is
Quillen equivalent to the usual model structure, but is much better behaved with
respect to homotopy theory of monoids and modules objects.
The model category $Sp^{\Sigma}$
is a symmetric monoidal model category for the smash product of
symmetric spectra.
Furthermore, all our assumptions
\ref{ass-1}, \ref{ass0}, \ref{ass1} and \ref{ass2} are satisfied
thanks to \cite[Theorem 14.5]{mmss},
\cite[Thm. 3.1, Thm. 3.2]{shi}, and \cite[Cor. 4.3]{shi} in conjunction with Lemma \cite[5.4.4]{hss}.

The category $Comm(Sp^{\Sigma})$ is the category of commutative symmetric
ring spectra, together with the positive stable model structure.
 The category $Comm(Sp^{\Sigma})$ will be denoted
by $S-Alg$, and its objects will simply  be called
\emph{commutative $S$-algebras} or also \emph{bn rings} (where \emph{bn} stands for
\emph{brave new}).
For any $E\in Sp^{\Sigma}$, we will set
$$\pi_{i}(E):=\pi^{stab}_{i}(RE),$$
where $RE$ is a fibrant replacement of $E$ in $Sp^{\Sigma}$, and
$\pi^{stab}_{i}(RE)$ are the naive stable homotopy groups
of the $\Omega$-spectrum $RE$. Note that if $E\in Sp^{\Sigma}$ is fibrant, then there is a natural isomorphism
$\pi^{stab}_{*}E\simeq \pi_{*}E$, and that a map $f:E'\rightarrow E''$ is a weak
equivalence in $Sp^{\Sigma}$ if and only if $\pi_{*}f$ is an isomorphism.

When $A$ is a commutative $S$-algebra, the $\mathbb{Z}$-graded abelian group $\pi_{*}(A)$
has a natural structure of a commutative graded algebra. In the same
way, when $M$ an $A$-module, $\pi_{*}(M)$ becomes a
graded $\pi_{*}(A)$-module. An object $E$ will be called connective, or
$(-1)$-connected, if
$\pi_{i}(E)=0$ for all $i<0$. We let $\mathcal{C}_{0}$ be
$Sp^{\Sigma}_{c}$, be full subcategory of connective objects in $Sp^{\Sigma}$. We let
$\mathcal{A}$ be $S-Alg_{0}$, the full subcategory of $Comm(Sp^{\Sigma})$ consisting of
commutative $S$-algebras $A$ with $\pi_{i}(A)=0$ for any $i\neq 0$.
If we denote by $H: CommRings\longrightarrow Comm(\mathcal{C})$ the Eilenberg-MacLane functor,
then $S-Alg_{0}$ is the subcategory of $S-Alg$ formed by all the
commutative $S$-algebras equivalent to some $Hk$ for some
commutative ring $k$.

\begin{lem}\label{connbnag}
The triplet $(Sp^{\Sigma},Sp^{\Sigma}_{c},S-Alg_{0})$ is a HA context.
\end{lem}

\begin{proof} The only thing to check is that
any object $A\in S-Alg_{0}$ is $Sp^{\Sigma}_{c}$-good in the sense
of Def. \ref{dgood}. Thanks to the equivalence between the homotopy theory
of $Hk$-modules and of complexes of $k$-modules (see \cite[Thm. IV.2.4]{ekmm}), this has
already been proved during the proof of Lem. \ref{lcont} $(2)$ . \end{proof}

The following example lists some classes of formally \'etale maps in brave new algebraic geometry,
according to our general definitions in Chapter \ref{partI.1}.

\begin{ex}\label{ex1}
\emph{
\begin{enumerate}
  \item If $A$ and $B$ are connective $S$-algebras,
a morphism $A\rightarrow B$ is formally $thh$-\'etale if and only if
  it is formally \'etale (\cite[Cor. 2.8]{min}).
\item A morphism of (discrete) commutative rings $R\rightarrow R'$ is formally \'etale if and only if
    the associated morphism $HR\rightarrow HR'$ of bn rings is formally \'etale if and only if
    the associated morphism $HR\rightarrow HR'$ of bn rings is formally $thh$-\'etale  (\cite[\S 5.2]{hagI}).
\item the complexification
map $KO\rightarrow KU$ is $thh$-formally \'etale (by \cite[p. 3]{ro}) hence formally
    \'etale. More generally, the same argument shows that any \textit{Galois extension}
of bn rings, according to J. Rognes \cite{ro}, is formally $thh$-\'etale, hence formally \'etale.
    \item There exist examples of formally \'etale morphisms of bn-rings
which are not $thh$-\'etale (see \cite{min} or \cite[\S 5.2]{hagI}).
\end{enumerate}
}
\end{ex}

As in the case of complicial algebraic geometry
(Ch. \ref{IIunb}) we find it useful to introduce also \textit{strong} versions
for properties of morphisms between commutative $S$-algebras.

\begin{df}\label{bnagstrong}
\begin{enumerate}
\item
Let $A\in S-Alg$, and $M$ be an $A$-module. The $A$-module $M$ is \emph{strong}\index{strong!$A$-module for $A\in S-Alg$}
if the natural morphism
$$\pi_{*}(A)\otimes_{\pi_{0}(A)}\pi_{0}(M)\longrightarrow \pi_{*}(M)$$
is an isomorphism.
\item A morphism $A\longrightarrow B$ in $S-Alg$ is \emph{strongly flat}\index{morphism in $S-Alg$!strongly flat}
(resp. \emph{strongly (formally) smooth}\index{morphism in $S-Alg$!strongly (formally) smooth}, resp. \emph{strongly (formally) \'etale}\index{morphism in $S-Alg$!strongly (formally) \'etale},
resp. \emph{a strong Zariski open immersion}\index{morphism in $S-Alg$!strong Zariski open immersion}) if $B$ is strong as an $A$-module, and
if the morphism of affine schemes
$$Spec\, \pi_{0}(B) \longrightarrow Spec\, \pi_{0}(A)$$
is flat (resp. (formally) smooth, resp. (formally) \'etale, resp. a Zariski open immersion).
\end{enumerate}
\end{df}

One of the main difference between derived algebraic geometry and
unbounded derived algebraic geometry was that
the strong notions of flat, smooth, \'etale and Zariski open immersion are not
as easily related to the corresponding general notions presented in
\S \ref{partI.2}. In the present situation, the comparison is even more loose as
typical phenomena arising from the existence of Steenrod operations
in characteristic $p$, make the notion of smooth morphisms of $S$-algebras rather subtle, and
definitely different from the above notion of strongly smooth morphisms.
We do not think this is a problem of the theory, but rather we think of this
as an interesting new feature of brave new algebraic geometry,
as compared to derived algebraic geometry, and well worthy of investigation.

\begin{prop}\label{pII-28'}
Let $f : A \longrightarrow B$ be a morphism in $S-Alg$.
\begin{enumerate}
\item If $A$ and $B$ are connective, the morphism
$f$ is \'etale (resp. a Zariski open immersion) in the sense of Def. \ref{d5},
if and only if
$f$ is strongly \'etale (resp. a strong Zariski open immersion).
\item If the morphism $f$ is strongly flat (
resp. strongly \'etale, resp. a strong Zariski open immersion), then
it is flat (resp. \'etale, resp. a Zariski open immersion) in the sense
of Def. \ref{d5}.
\end{enumerate}
\end{prop}

\begin{proof} $(1)$ The proof is the same as for
Thm. \ref{tII-1}. \\

$(2)$ The proof is the same as for Prop. \ref{pII-28}. \end{proof}

The reader will notice that the proof of
Thm. \ref{tII-1} $(2)$ does not apply to the present context
as a smooth morphism of commutative rings is in general
not smooth when considered
as a morphism of commutative $S$-algebras. The typical example
of this phenomenon is the following.

\begin{prop}\label{affineline}
\begin{itemize}
    \item The canonical map $H\mathbb{Q}\rightarrow H(\mathbb{Q}[T])$ is smooth, i-smooth and perfect.
    \item The canonical map $H\mathbb{F}_{p}\rightarrow H(\mathbb{F}_{p}[T])$ is strongly smooth but
\emph{not} formally smooth, nor formally i-smooth (Def. \ref{d7}).
\end{itemize}
\end{prop}

\begin{proof} For any discrete commutative ring $k$, we have a canonical
map $a_{k}:Hk[T]:=F_{Hk}(Hk)\rightarrow H(k[T])$
of commutative $Hk$-algebras, corresponding to the map $H(k\rightarrow k[T])$
pointing the element $T$. Now
$$\pi_{*}(Hk[T])\simeq \bigoplus_{r\geq 0}\mathrm{H}_{*}(\Sigma_{r},k)$$
where in the group homology $\mathrm{H}_{*}(\Sigma_{r},k)$, $k$ is a trivial $\Sigma_{r}$-module.
Since $\mathrm{H}_{n}(\Sigma_{r},\mathbb{Q})=0$, for $n\neq 0$,
and $\mathrm{H}_{0}(\Sigma_{r},\mathbb{Q})\simeq \mathbb{Q}$, for any $r\geq 0$
we see that $a_{\mathbb{Q}}$ is a stable homotopy equivalence, and therefore a weak equivalence (\cite[Thm. 3.1.11]{hss}).
In other words $H(\mathbb{Q}[T])$ ``is'' the free commutative $H\mathbb{Q}$-algebra on one generator;
therefore it is finitely presented over $H\mathbb{Q}$ and, since for any (discrete) commutative ring $k$ we have
$$\mathbb{L}_{Hk[T]/Hk}\simeq Hk[T],$$ the cotangent complex $\mathbb{L}_{H(\mathbb{Q}[T])/H\mathbb{Q}}$ is free
of rank one over $H(\mathbb{Q}[T])$, hence projective and perfect. So $H\mathbb{Q}\rightarrow H(\mathbb{Q}[T])$ is smooth and perfect.

Let's move to the char $p>0$ case. It is clear that $H\mathbb{F}_{p}\rightarrow H(\mathbb{F}_{p}[T])$ is strongly
smooth; let's suppose that it is formally smooth.
In particular $\mathbb{L}_{H(\mathbb{F}_{p}[T])/H\mathbb{F}_{p}}$ is a projective $H(\mathbb{F}_{p}[T])$-module.
Therefore $\pi_{*}\mathbb{L}_{H(\mathbb{F}_{p}[T])/H\mathbb{F}_{p}}$ injects into
$$\pi_{*}\coprod_{E}H(\mathbb{F}_{p}[T])\simeq \prod_{E}\mathbb{F}_{p}[T]$$ (concentrated in degree $0$).
But, by \cite[Thm. 4.2]{bamcc} and \cite[Thm. 4.1]{rr},
$$\pi_{*}\mathbb{L}_{H(\mathbb{F}_{p}[T])/H\mathbb{F}_{p}}\simeq (H(\mathbb{F}_{p}[T]))_{*}(H\mathbb{Z})$$
and the last ring has $(H\mathbb{F}_{p})_{*}(H\mathbb{Z})$ as a direct summand
(using the augmentation $H(\mathbb{F}_{p}[T])\rightarrow H\mathbb{F}_{p}$).
Now, it is known that $(H\mathbb{F}_{p})_{*}(H\mathbb{Z})$ is not concentrated in degree $0$
(it is a polynomial $\mathbb{F}_{p}$-algebra in positive degrees generators, for $p=2$,
and the tensor product of such an algebra with an exterior $\mathbb{F}_{p}$-algebra for odd $p$).
Therefore $H\mathbb{F}_{p}\rightarrow H(\mathbb{F}_{p}[T])$ cannot be formally smooth.
\end{proof}

We remark again that, as now made clear by the proof above, the conceptual reason for the non-smoothness of
$H\mathbb{F}_{p}\rightarrow H(\mathbb{F}_{p}[T])$ is essentially the existence of (non-trivial)
Steenrod operations in characteristic $p>0$. Since formal smoothness is stable under base-change,
we also conclude that $H\mathbb{Z}\rightarrow H(\mathbb{Z}[T])$ is not formally smooth. The same
argument also shows that this morphism is not formally i-smooth. \\

The following example shows that the converse of Prop. \ref{pII-28'} (2)
is false in general.

\begin{ex} \emph{The complexification map $m:KO\rightarrow KU$ is
formally \'etale but not strongly formally
\'etale. In fact, we have $$\pi_{*}m: \pi_{*}(KO)=\mathbb{Z}[\eta, \beta,
\lambda^{\pm 1}]/(\eta^{3},2\eta, \eta \beta, \beta^{2}-4\lambda)\longrightarrow
\pi_{*}(KO)=\mathbb{Z}[\nu^{\pm 1}],
$$ with $deg (\eta)=1, deg(\beta)= 4, deg(\lambda)=8, deg(\nu)=2$, $\pi_{*}m(\eta)= 0$,
$\pi_{*}m(\beta)= 2\nu^{2}$ and
$\pi_{*}m(\lambda) = \nu^{4}$. In particular $\pi_{0}m$ is an isomorphism
(hence \'etale) but $m$ is not strong.\\
We address the reader to \cite[Rmk. 5.2.9]{hagI} for an example, due to
M. Mandell, of a non connective formally \'etale extension of
$H\mathbb{F}_{p}$, which is therefore not strongly formally \'etale.
There also exist examples of Zariski open immersion
$HR \longrightarrow A$, here $R$ is a smooth commutative $k$-algebra, 
such that $A$ possesses non trivial negative homotopy groups (see
\cite[\S 5.2]{hagI} for more details).
}
\end{ex}

The opposite model category $S-Alg$ will be denoted
by $SAff$. We will endow it with the following \textit{strong \'etale model
topology}.

\begin{df}\label{dII-25bis}
A family of morphisms $\{Spec\, A_{i} \longrightarrow Spec\, A\}_{i\in I}$
in $SAff$ is a \emph{strong \'etale covering family}\index{strong \'etale covering family!in $SAff$}
(or simply \emph{s-\'et covering family}) if it satisfies the following two
conditions.
\begin{enumerate}
\item Each morphism $A \longrightarrow A_{i}$
is strongly \'etale.
\item There exists a finite sub-set $J\subset I$ such that the
family $\{A \longrightarrow A_{i}\}_{i\in J}$ is a
formal covering family in the sense of \ref{dcov}.
\end{enumerate}
\end{df}

Using the definition of strong \'etale morphisms, we immediately check
that a family of morphisms $\{Spec\, A_{i} \longrightarrow Spec\, A\}_{i\in I}$
in $SAff$ is a s-\'et covering family if and only if there exists a finite
sub-set $J\subset I$ satisfying the following two conditions.

\begin{itemize}

\item For all $i\in I$, the natural morphism
$$\pi_{*}(A)\otimes_{\pi_{0}(A)}\pi_{0}(A_{i}) \longrightarrow
\pi_{*}(A_{i})$$
is an isomorphism.

\item The morphism of affine schemes
$$\coprod_{i\in J} Spec\, \pi_{0}(A_{i}) \longrightarrow Spec\, \pi_{0}(A)$$
is \'etale and surjective.

\end{itemize}

\begin{lem}\label{lII-26bis}
The s-\'et\index{s-\'et!model topology on $SAff$} covering families define a model topology on $SAff$, that satisfies
assumption \ref{ass5}.
\end{lem}

\begin{proof} The same as for Lem. \ref{lII-5}. \end{proof}

The model topology s-\'et gives rise to a model category
of stacks $SAff^{\sim,\textrm{s-\'et}}$\index{$SAff^{\sim,\textrm{s-\'et}}$}.

\begin{df}\label{dII-26'}
\begin{enumerate}
\item An \emph{$S$-stack}\index{$S$-stack} is an object
$F\in SAff^{\sim,\textrm{s-\'et}}$ which is a stack
in the sense of Def. \ref{dstack}.
\item
The \emph{model category of $S$-stacks} is
$SAff^{\sim,\textrm{s-\'et}}$\index{$SAff^{\sim,\textrm{s-\'et}}$}, and its homotopy category will be simply
denoted by $\mathrm{St}(S)$.\index{$\mathrm{St}(S)$}
\end{enumerate}
\end{df}

We now set \textbf{P} to be the class of fip-smooth morphisms (i.e. formally perfect and
formally i-smooth morphisms) in
$SAlg$, and \textbf{P}$_{\textrm{s-\'et}}$ be the class of
strongly \'etale morphisms.

\begin{lem}\label{lII-27bis}
\begin{enumerate}
\item
The class \textbf{P}$_{\textrm{s-\'et}}$ of strongly \'etale morphisms and the s-\'et model topology satisfy
assumptions \ref{ass4}.
\item
The class \textbf{P} of fip-smooth morphisms and the s-\'et model topology satisfy
assumptions \ref{ass4}.
\end{enumerate}
\end{lem}

\begin{proof} It is essentially the same as for Lem. \ref{lII-27}. \end{proof}

\begin{cor}\label{clII-27bis}
\begin{enumerate}
\item The 5-tuple $(Sp^{\Sigma}, Sp^{\Sigma}, S-Alg, \textrm{s-\'et}, \mathbf{P}_{\textrm{s-\'et}})$
is a HAG context.
\item The 5-tuple $(Sp^{\Sigma}, Sp^{\Sigma}_{c}, S-Alg_{0}, \textrm{s-\'et}, \mathbf{P})$
is a HAG context.
\end{enumerate}
\end{cor}

According to our general theory, the notions of morphisms in
\textbf{P} and \textbf{P}$_{\textrm{s-\'et}}$ gives two notions of geometric stacks in
$SAff^{\sim,\textrm{s-\'et}}$.

\begin{df}\label{dII-27bis}
\begin{enumerate}
\item \emph{A $n$-geometric Deligne-Mumford $S$-stack}\index{$S$-stack!$n$-geometric Deligne-Mumford} is an $n$-geometric $S$-stack
with respect to the class \textbf{P}$_{\textrm{s-\'et}}$ of strongly \'etale morphisms.

\item \emph{An $n$-geometric $S$-stack}\index{$S$-stack!$n$-geometric} is an $n$-geometric $S$-stack
with respect to the class \textbf{P} of fip-smooth morphisms.

\end{enumerate}
\end{df}

Of course, as \textbf{P}$_{\textrm{s-\'et}}$ is included in \textbf{P} any strong $n$-geometric
Deligne-Mumford $S$-stack is an $n$-geometric $S$-stack. \\

Finally, the reader can easily check the following proposition.

\begin{prop}\label{artinbnag}

\begin{enumerate}
\item The topology $\textrm{s-\'et}$ and the class \textbf{P}$_{\textrm{s-\'et}}$ satisfy Artin condition
(for the HA context $(Sp^{\Sigma},Sp^{\Sigma},S-Alg)$.

\item The topology $\textrm{s-\'et}$ and the class \textbf{P} satisfy Artin's condition
(for the HA context $(Sp^{\Sigma},Sp^{\Sigma}_{c},S-Alg_{0})$.

\end{enumerate}
\end{prop}

In particular, we obtain as a corollary of Thm. \ref{t1} that
any $n$-geometric $S$-stack has an obstruction theory relative to
the HA context $(Sp^{\Sigma},Sp^{\Sigma}_{c},S-Alg_{0})$. In the same way, any
strong $n$-geometric Deligne-Mumford $S$-stack has an obstruction theory
relative to the context $(Sp^{\Sigma},Sp^{\Sigma},S-Alg)$.

Without going into details, we mention that
all the examples of geometric $D$-stacks given in the previous chapter
can be generalized to examples of geometric $S$-stacks.
One fundamental example is $\mathbf{Perf}^{CW}_{[a,b]}$
of perfect CW-modules of amplitude contained in $[a,b]$, which
by a similar argument as for Prop. \ref{pcwmod} is a $1$-geometric
$S$-stack.

\section{Elliptic cohomology as a Deligne-Mumford $S$-stack}\label{IIbnag.2}

In this final section we present the construction
of a $1$-geometric Deligne-Mumford $S$-stack using the sheaf of spectra of topological modular
forms. \\

The Eilenberg-MacLane spectrum construction (see \cite[Ex. 1.2.5]{hss})
gives rise to a fully faithful functor
$$\mathbb{L}H_{!} : \mathrm{St}(\mathbb{Z}) \longrightarrow \mathrm{St}(S),$$
which starts from the homotopy category of stacks on the
usual \'etale site of affine schemes,\index{$\mathrm{St}(\mathbb{Z})$} i.e. $$\mathrm{St}(\mathbb{Z}):=\mathrm{Ho}(\mathbb{Z}-Aff^{\sim,\textrm{\'et}}).$$ This functor
has a right adjoint, called the \emph{truncation functor}
$$h^{0}:=H^{*} : \mathrm{St}(S) \longrightarrow \mathrm{St}(\mathbb{Z}),$$
simply given by composing a simplicial presheaf $F : SAff^{op} \longrightarrow SSet$
with the functor $H : Aff \longrightarrow S-Aff$. \\

Let us denote by $\overline{\mathcal{E}}$ the moduli stack of
generalized elliptic curves with integral geometric fibers,
which is the standard
compactification of the moduli stack of elliptic curves by adding
the nodal curves at infinity (see e.g. \cite[IV]{dera}, where it is denoted by $\mathcal{M}_{(1)}$); recall that
$\overline{\mathcal{E}}$ is a Deligne-Mumford stack, proper and smooth over $Spec\; \mathbb{Z}$ (\cite[Prop. 2.2]{dera}).

As shown by
recent works of M. Hopkins, H. Miller, P. Goerss,
N. Strickland, C. Rezk and M. Ando, there exists a natural
presheaf of commutative $S$-algebras on the small
\'etale site $\overline{\mathcal{E}}_{\textrm{\'et}}$ of
$\overline{\mathcal{E}}$. We will denote this presheaf by
$\mbox{\em tmf}$. Recall that by construction, if $U=Spec\, A \longrightarrow \overline{\mathcal{E}}$
is an \'etale morphism, corresponding to an elliptic curve $E$ over
the ring $A$, then $\mbox{\em tmf}(U)$ is the (connective) elliptic cohomology theory
associated to the formal group of $E$ (in particular, one
has $\pi_{0}(\mbox{\em tmf}(U))=A$). Recall also that the (derived) global sections
$\mathbb{R}\Gamma(\overline{\mathcal{E}},\mbox{\em tmf})$, form a
commutative $S$-algebra, well defined in $\mathrm{Ho}(S-Alg)$,
called the \emph{spectrum of topological modular forms}, and denoted by $\mathsf{tmf}$ \footnote{The notation here is a bit nonstandard: 
what is usually called the spectrum of topological modular forms is actually the connective cover of the spectrum we have denoted by $\mathsf{tmf}$}.

Let $U \longrightarrow \overline{\mathcal{E}}$ be a surjective \'etale morphism with
$U$ an affine scheme, and let us consider its nerve
$$\begin{array}{cccc}
U_{*} : & \Delta^{op} & \longrightarrow & Aff \\
 & [n] & \mapsto & U_{n}:=\underbrace{U\times_{\overline{\mathcal{E}}}U\times_{\overline{\mathcal{E}}}\dots
\times_{\overline{\mathcal{E}}}U}_{n\; times}.
\end{array}$$

This is a simplicial object in $\overline{\mathcal{E}}_{\textrm{\'et}}$, and by
applying $\mbox{\em tmf}$ we obtain a
co-simplicial object in $S-Alg$
$$\begin{array}{cccc}
\mbox{\em tmf}(U_{*}) : & \Delta & \longrightarrow & S-Alg \\
 & [n] & \mapsto & \mbox{\em tmf}(U_{n}). 
\end{array}$$
Taking $\mathbb{R}\underline{Spec}$ of this diagram we obtain
a simplicial object in the model category $S-Aff^{\sim,\textrm{s-\'et}}$
$$\begin{array}{cccc}
\mathbb{R}\underline{Spec}\; (\mbox{\em tmf}(U_{*})) : & \Delta^{op} & \longrightarrow & S-Aff^{\sim,\textrm{s-\'et}} \\
 & [n] & \mapsto & \mathbb{R}\underline{Spec}\; (\mbox{\em tmf}(U_{n})).
\end{array}$$

The homotopy colimit of this diagram will be denoted by
$$\overline{\mathcal{E}}_{\mathbf{S}}:=\mathrm{hocolim}_{n\in \Delta^{op}}
\mathbb{R}\underline{Spec}\; (\mbox{\em tmf}(U_{*})) \in \mathrm{St}(S).$$

The following result is technically just a remark as there is
essentially nothing to prove; however, we prefer to state it
as a theorem to emphasize its importance.

\begin{thm}\label{ttmf}
The stack $\overline{\mathcal{E}}_{\mathbf{S}}$ defined above is
a strong Deligne-Mumford $1$-geometric $S$-stack. Furthermore $\overline{\mathcal{E}}_{\mathbf{S}}$ is a
 ``brave new derivation'' of the moduli stack $\overline{\mathcal{E}}$ of elliptic curves, i.e. there exists a natural
isomorphism in $\mathrm{St}(\mathbb{Z})$
$$h^{0}(\overline{\mathcal{E}}_{\mathbf{S}})\simeq \overline{\mathcal{E}}.$$
\end{thm}

\begin{proof} To prove that $\overline{\mathcal{E}}_{\mathbf{S}}$ is geometric, it is enough to check that
the simplicial object $\mathbb{R}\underline{Spec}\; (\mbox{\em tmf}(U_{*}))$ is a strongly \'etale
Segal groupoid. For this, recall that for
any morphism $U=Spec\, B \rightarrow V=Spec\, A$ in $\overline{\mathcal{E}}_{\textrm{\'et}}$, the natural morphism
$$\pi_{*}(\mbox{\em tmf}(V))\otimes_{\pi_{0}(\mbox{\em tmf}(V))}
\pi_{0}(\mbox{\em tmf}(U))\simeq \pi_{*}(\mbox{\em tmf}(V))\otimes_{A}B \longrightarrow
\pi_{*}(\mbox{\em tmf}(U))$$
is an isomorphism. This shows that the functor
$$\mathbb{R}\underline{Spec}\; (\mbox{\em tmf}(-)) : \overline{\mathcal{E}}_{\textrm{\'et}} \longrightarrow
S-Aff^{\sim,\textrm{s-\'et}}$$
preserves homotopy fiber products and therefore sends Segal groupoid objects
to Segal groupoid objects. In particular,
$\mathbb{R}\underline{Spec}\; (\mbox{\em tmf}(U_{*}))$ is a Segal groupoid object. The same fact also shows that
for any morphism $U=Spec\, B \rightarrow V=Spec\, A$ in $\overline{\mathcal{E}}_{\textrm{\'et}}$,
the induced map $\mbox{\em tmf}(V) \longrightarrow \mbox{\em tmf}(U)$ is
a strong \'etale morphism. This implies that
$\mathbb{R}\underline{Spec}\; (\mbox{\em tmf}(U_{*}))$ is a strongly \'etale
Segal groupoid object in representable $S$-stacks, and thus
shows that $\overline{\mathcal{E}}_{\mathbf{S}}$ is indeed a strong Deligne-Mumford $1$-geometric
$S$-stack.

The truncation functor $h^{0}$ clearly commutes with homotopy colimits, and therefore
$$h^{0}(\overline{\mathcal{E}}_{\mathbf{S}})\simeq
\mathrm{hocolim}_{n\in \Delta^{op}}h^{0}(\mathbb{R}\underline{Spec}\;
(\mbox{\em tmf}(U_{n}))) \in \mathrm{St}(\mathbb{Z}).$$
Furthermore, for any connective representable $S$-stack, $\mathbb{R}\underline{Spec}\, A$, one
has a natural isomorphism $h^{0}(\mathbb{R}\underline{Spec}\, A)\simeq Spec\, \pi_{0}(A)$.
Therefore, one sees immediately that there is a natural isomorphism of
simplicial objects in $\mathbb{Z}-Aff^{\sim,\textrm{\'et}}$
$$h^{0}(\mathbb{R}\underline{Spec}\; (\mbox{\em tmf}(U_{*})))\simeq U_{*}.$$
Therefore, we get
$$h^{0}(\overline{\mathcal{E}}_{\mathbf{S}})\simeq
\mathrm{hocolim}_{n\in \Delta^{op}}h^{0}(\mathbb{R}\underline{Spec}\; (\mbox{\em tmf}(U_{n})))\simeq
\mathrm{hocolim}_{n\in \Delta^{op}}U_{n}\simeq \overline{\mathcal{E}},$$
as $U_{*}$ is the nerve of an \'etale covering of $\overline{\mathcal{E}}$. \end{proof}

Theorem \ref{ttmf} tells us that the presheaf of topological modular forms
$\mbox{\em tmf}$ provides a natural geometric $S$-stack $\overline{\mathcal{E}}_{\mathbf{S}}$
whose truncation is the usual stack of elliptic curves
$\overline{\mathcal{E}}$. Furthermore, as one can show that the small strong \'etale topoi
of $\overline{\mathcal{E}}_{\mathbf{S}}$ and $\overline{\mathcal{E}}$
coincide (this is a general fact about strong \'etale model topologies), we see that
$$\mathsf{tmf}:=\mathbb{R}\Gamma(\overline{\mathcal{E}},\mbox{\em tmf}) \simeq
\mathbb{R}\Gamma(\overline{\mathcal{E}}_{\mathbf{S}},\mathcal{O}),$$
and therefore that topological modular forms can be simply interpreted
as \emph{functions on the geometric $S$-stack $\overline{\mathcal{E}}_{\mathbf{S}}$}.
Of course, our construction of $\overline{\mathcal{E}}_{\mathbf{S}}$ has essentially been rigged to
make this true, so this is not a surprise. However, we have gained a bit from the conceptual
point of view: since after all $\overline{\mathcal{E}}$ is a moduli stack, now
that we know the existence of the geometric $S$-stack
$\overline{\mathcal{E}}_{\mathbf{S}}$ we can ask for a \textit{modular interpretation}
of it, or in other words for a direct geometric description of the corresponding
simplicial presheaf on $S-Aff$. An answer to this question not only would
provide a direct construction
of $\mathsf{tmf}$, but would also give a conceptual interpretation of it in a geometric
language closer the usual notion of modular forms.

\begin{Q}\label{q1}
Find a modular interpretation of the
$S$-stack
$\overline{\mathcal{E}}_{\mathbf{S}}$.
\end{Q}

Essentially, we are asking for the brave new ``objects'' that the
$S$-stack $\overline{\mathcal{E}}_{\mathbf{S}}$ classifies.
We could also consider the non-connective
version of the $S$-stack $\overline{\mathcal{E}}_{\mathbf{S}}$ (defined through the non-connective version of $\mbox{\em tmf}$) for
which a modular interpretation seems much more accessible.\\
Very recent work by J. Lurie (see \cite{luelcoom} for a detailed announcement of his results) answers in fact to Question \ref{q1}; he shows that such a variant of 
$\overline{\mathcal{E}}_{\mathbf{S}}$ classifies brave new versions of \cite{ahs}'s elliptic spectra plus additional data (called orientations).
This moduli-theoretic point of view 
 makes use of some very interesting notions of \emph{brave new abelian varieties},
\emph{brave new formal groups} and
their geometry. The complete picture (possibly extended to higher chromatic levels) does not only give
an alternative construction of the spectrum $tmf$ (and a better functoriality) but
it could be the starting point of a rather new\footnote{J. Lurie's approach has as a byproduct also a natural construction of $G$-\textit{equivariant} versions of elliptic cohomology, for any compact $G$.} and deep interaction between
stable homotopy theory and homotopical algebraic geometry, involving 
many new questions and objects, and probably also 
new insights on classical objects of algebraic topology.

\begin{appendix}

\chapter{Classifying spaces of model categories}

The classifying space of a model category $M$ is
defined to be $N(M_{W})$, the nerve of its subcategory of equivalences.
More generally, if $C\subset M_{W}$ is a full subcategory
of the category of equivalences in $M$, which is
closed by equivalences in $M$, the classifying space
of $C$ is $N(C)$, the nerve of $C$.

We fix a $\mathbb{V}$-small model category $M$ and a full subcategory
$C \subset M_{W}$ closed by equivalences. We consider
the model category $(C,C)^{\wedge}$, defined in \cite[\S 2.3.2]{hagI}.
Recall that the underlying category of $(C,C)^{\wedge}$
is the category $SSet_{\mathbb{V}}^{C^{op}}$, of
$\mathbb{V}$-small simplicial presheaves on $C$.
The model structure of $(C,C)^{\wedge}$ is defined as  the
left Bousfield localization of the levelwise projective model
structure on $SSet_{\mathbb{V}}^{C^{op}}$, by inverting
all the morphisms in $C$. The important fact is that
local objects in $(C,C)^{\wedge}$ are functors
$F : C^{op} \longrightarrow SSet_{\mathbb{V}}$
sending all morphisms in $C$ to equivalences.

We define an adjunction
$$N : (C,C)^{\wedge} \longrightarrow SSet/N(C) \qquad
(C,C)^{\wedge} \longleftarrow SSet/N(C) : S$$
in the following way. A functor $F : C^{op} \longrightarrow
SSet$ is sent to the simplicial set $N(F)$, for which
the set of $n$-simplices is the set of parirs
$$N(F)_{n}:=\{(c_{0} \rightarrow c_{1} \rightarrow \dots \rightarrow c_{n},\alpha)\}$$
where $(c_{0} \rightarrow c_{1} \rightarrow \dots \rightarrow c_{n})$ is an $n$-simplex
in $N(C)$ and $\alpha \in F(c_{n})_{n}$ is an $n$-simplex in
$F(c_{n})$. Put in an other way,
$N(F)$ is the diagonal of the bi-simplicial set
$$(n,m) \mapsto N(C/F_{n})_{m}$$
where $C/F_{n}$ is the category of objects of the presheaf
of $n$-simplices in $F$. The functor $N$ has a right adjoint
$$S : SSet/N(C) \longrightarrow (C,C)^{\wedge}$$
sending $X \longrightarrow N(C)$ to the simplicial presheaf
$$\begin{array}{cccc}
S(X) : & C^{op} & \longrightarrow & SSet \\
 & x & \mapsto & \underline{Hom}_{SSet/N(C)}(N(h_{x}),X),
\end{array}$$
where $h_{x}$ is the presheaf of sets represented
by $x \in C$ (note that
$N(h_{x})\rightarrow N(C)$ is isomorphic
to $N(C/x)\rightarrow N(C)$).

\begin{prop}\label{papp1}
The adjunction
$(N,S)$ is a Quillen equivalence.
\end{prop}

\begin{proof} First of all we need to check that
$(N,S)$ is a Quillen adjunction. For this we use the
standard properties of left Bousfield localizations, and we see
that it is enough to check that
$$N : SPr(C) \longrightarrow SSet/N(C) \qquad
SPr(C) \longleftarrow SSet/N(C) : S$$
is a Quillen adjunction (where $SPr(C)$ is the projective
model structure of simplicial presheaves on $C$), and that
$S$ preserves fibrant objects. These two facts are clear
by definition of $S$ and the description of fibrant objects
in $N(C,C)^{\wedge}$.

For any $x\in C$, the
morphism $N(h_{x})=N(C/x) \longrightarrow N(C)$ is
isomorphic in $\mathrm{Ho}(SSet/N(C))$ to
$x \longrightarrow N(C)$.  Therefore,
for $X\in \mathrm{Ho}(SSet/N(C))$ and $x\in C$, the simplicial set
$$\mathbb{R}S(X)(x) \simeq \mathbb{R}\underline{Hom}_{SSet/N(C)}(
x,X)$$
is naturally isomorphic in $\mathrm{Ho}(SSet)$ to the homotopy fiber
of $X \longrightarrow N(C)$ taken at $x$.
This clearly implies that the right derived functor
$$\mathbb{R}S : \mathrm{Ho}(SSet/N(C)) \longrightarrow \mathrm{Ho}((C,C)^{\wedge})$$
is conservative. In particular, it only remains to show that the
adjunction morphism $Id \longrightarrow \mathbb{R}S\mathbb{L}N$
is an isomorphism. But this last assumption follows
from the definition the functor $N$
and from a standard lemma (se for example \cite{q2}), which
shows that the homotopy fiber at $x\in C$ of
$N(F) \longrightarrow N(C)$ is naturally equivalent to
$F(x)$ when $F$ is fibrant in $(C,C)^{\wedge}$. \end{proof}

Recall from \cite[Lem. 4.2.2]{hagI} that for any
object $x \in C$, one can construct a local
model for $h_{x}$ as $\underline{h}_{R(x)}$, sending
$y\in C$ to $Hom_{C}(\Gamma^{*}(y),R(x))$, where
$\Gamma^{*}$ is a co-simplicial replacement functor in $M$.
The natural morphism $h_{x} \longrightarrow \underline{h}_{R(x)}$
being an equivalence in $(C,C)^{\wedge}$, one finds
using Prop. \ref{papp1} natural equivalences
of simplicial sets
$$\underline{Hom}(h_{y},\underline{h}_{R(x)})\simeq
Map^{eq}_{M}(y,x) \simeq \mathbb{R}\underline{Hom}_{SSet/N(C)}(
N(h_{y}),N(h_{x}))\simeq y\times^{h}_{N(C)}x,$$
where $Map^{eq}_{M}(y,x)$ is the sub-simplicial set
of the mapping space $Map_{M}(y,x)$ consisting of
equivalences. As a corollary of this we find the
important result due to Dwyer and Kan. For this, we recall that
the simplicial monoid of self equivalences of an
object $x\in M$ can be defined as
$$Aut(x):=\underline{Hom}_{(C,C)^{\wedge}}(\underline{h}_{R(x)},\underline{h}_{R(x)}).$$

\begin{cor}\label{capp1}
Let $C \subset M_{W}$ be a full subcategory of
equivalences in a model category $M$, which
is stable by equivalences. Then, one has a natural
isomorphism in $\mathrm{Ho}(SSet)$
$$
N(C)\simeq \coprod_{x\in \pi_{0}(N(C)))}BAut(x)$$
where $Aut(x)$ is the simplicial monoid of
self equivalences of $x$ in $M$.
\end{cor}

Another important consequence of Prop. \ref{papp1} is the
following interpretation of mapping spaces in term
of homotopy fibers between classifying spaces of certain model categories.

\begin{cor}\label{capp2}
Let $M$ be a model category and $x,y\in M$ be two fibrant and
cofibrant objects in $M$. Then, there exists a natural
homotopy fiber sequence of simplicial sets
$$Map_{M}(x,y) \longrightarrow N((x/M)_{W}) \longrightarrow
N(M_{W}),$$
where the homotopy fiber is taken at $y \in M$.
\end{cor}

\begin{proof} This follows easily from Cor. \ref{capp1}.
\end{proof}

We will need a slightly more functorial interpretation of Cor.
\ref{capp1} in the particular case where the model category $M$ is
simplicial. We assume now that $M$ is a $\mathbb{V}$-small
simplicial $\mathbb{U}$-model category (i.e. the model category
$M$ is a $SSet_{\mathbb{U}}$-model category in the sense of \cite{ho}).
We still let $C\subset M_{W}$ be a full subcategory of
the category of equivalences in $M$, and still assume that
$C\subset M_{W}$ is stable by equivalences.

We define an $S$-category $\mathcal{G}(C)$ is the following way (recall that
an $S$-category is a simplicially enriched category, see for example
\cite[\S 2.1]{hagI} for more details and notations). The objects of $\mathcal{G}(C)$
are the objects of $C$ which are furthermore fibrant and cofibrant in $M$.
For two objects $x$ and $y$, the simplicial set of morphisms is defined to be
$$\mathcal{G}(C)_{(x,y)}:=\underline{Hom}^{eq}_{M}(x,y),$$
where by definition $\underline{Hom}^{eq}_{M}(x,y)_{n}$ is the set
of equivalences in $M$ from $\Delta^{n}\otimes x$ to $y$
(i.e. $\underline{Hom}^{eq}_{M}(x,y)$
is the sub-simplicial set of $\underline{Hom}_{M}(x,y)$
consisting of equivalences). Clearly, the $S$-category
$\mathcal{G}(C)$ is \emph{groupoid like}, in the sense that its
category of connected components $\pi_{0}(\mathcal{G}(C))$
(also denoted by $\mathrm{Ho}(\mathcal{G}(C)$) is a groupoid (or equivalently
every morphism in $\mathcal{G}(C)$ has an inverse up to homotopy).
Let $C^{c,f}$ be the full subcategory of $C$ consisting of fibrant and
cofibrant objects in $C$. There exist two natural morphisms of $S$-categories
$$C \longleftarrow C^{c,f} \longrightarrow \mathcal{G}(C),$$
where a category is considered as an $S$-category with discrete
simplicial sets of morphisms. Passing to the nerves, one gets
a diagram of simplicial sets
$$N(C) \longleftarrow N(C^{c,f}) \longrightarrow N(\mathcal{G}(C)),$$
where the nerve functor is extended diagonally to $S$-categories (see e.g.
\cite{dk1}). Another interpretation of corollary \ref{capp1} is the following
result.

\begin{prop}\label{papp2}
With the above notations, the two morphisms
$$N(C) \longleftarrow N(C^{c,f}) \longrightarrow N(\mathcal{G}(C)),$$
are equivalences of simplicial sets.
\end{prop}

\begin{proof} It is well known that the left arrow is an equivalence as
a fibrant-cofibrant replacement functor gives an inverse up to homotopy.
For the right arrow, we let
$N(\mathcal{G}(C))_{n}$ be the category of $n$-simplicies in $N(\mathcal{G}(C))$,
defined by having the same objects and with
$$(N(\mathcal{G}(C))_{n})_{x,y}:=(N(\mathcal{G}(C))_{x,y})_{n}.$$
By definition of the nerve, one has a natural equivalence
$$N(\mathcal{G}(C))\simeq Hocolim_{n\Delta^{op}}(N(\mathcal{G}(C))_{n}).$$
Furthermore, it is clear that each functor
$$C^{c,f}=\mathcal{G}(C)_{0} \longrightarrow \mathcal{G}(C)_{n}$$
induces an equivalence of the nerves, as the $0$-simplex $[0] \rightarrow [n]$ clearly induces
a functor
$$N(\mathcal{G}(C)_{n}) \longrightarrow N(C^{c,f})$$
which is a homotopy inverse. \end{proof}

Proposition \ref{papp2} is another interpretation of Prop. \ref{papp1}, as
the $S$-category $\mathcal{G}(C)$ is groupoid-like,
the delooping theorem of G. Segal implies that there exists a natural
equivalence of simplicial sets
$$x\times^{h}_{N(C)}y\longleftarrow x\times^{h}_{N(C^{c,f})}y \longrightarrow
x\times^{h}_{N(\mathcal{G}(C))}y\longleftarrow \mathcal{G}(C)_{(x,y)}=\underline{Hom}^{eq}_{M}(x,y).$$

The advantage of Prop. \ref{papp2} over the more general proposition \ref{papp1}
is that it is more easy to state a functorial property of the
equivalences in the following particular context (the equivalence
in Prop. \ref{papp1} can also be made functorial, but it requires
some additional work, using for example simplicial localization techniques).

We assume that $G : M \longrightarrow N$ is a simplicial, left Quillen functor between
$\mathbb{V}$-small simplicial $\mathbb{U}$-model categories. We let $C\subset M_{W}$ and
$D\subset N_{W}$ be two full sub-categories stable by equivalences, and we suppose that
all objects in $M$ and $N$ are fibrant. Finally, we assume that the functor $G$ restricted
to cofibrant objects sends $C^{c}:=M^{c}\cap C$ to $D^{c}:=N^{c}\cap D$. In this situation,
we define an $S$-functor
$$G : \mathcal{G}(C) \longrightarrow \mathcal{G}(D)$$
simply by using the simplicial enrichment of $G$. Then, one has a commutative diagram
of $S$-categories
$$\xymatrix{
C  & C^{c} \ar[r] \ar[d]^-{G} \ar[l] & \mathcal{G}(C) \ar[d]^-{G} \\
D & D^{c} \ar[r] \ar[l] & \mathcal{G}(D),}$$
and thus a commutative diagram of simplicial sets
$$\xymatrix{
N(C)  & N(C^{c}) \ar[r] \ar[d]^-{G} \ar[l] & N(\mathcal{G}(C)) \ar[d]^-{G} \\
N(D) & N(D^{c}) \ar[r] \ar[l] & N(\mathcal{G}(D)).}$$
The important fact here is that this construction is associative
with respect to composition of the simplicial left Quillen functor $G$.
In other words, if one has a diagram of simplicial model categories
$M_{i}$ and simplicial left Quillen functors, together with sub-categories
$C_{i}\subset (M_{i})_{W}$ satisfying the required properties, then one
obtain a commutative diagram of diagrams of simplicial sets
$$\xymatrix{
N(C_{i})  & N(C_{i}^{c}) \ar[r] \ar[d]^-{G_{i}} \ar[l] & N(\mathcal{G}(C_{i})) \ar[d]^-{G_{i}} \\
N(D_{i}) & N(D_{i}^{c}) \ar[r] \ar[l] & N(\mathcal{G}(D_{i})).}$$

\chapter{Strictification}

Let $I$ be a $\mathbb{U}$-small category. For any $i\in I$
we let $M_{i}$ be a $\mathbb{U}$-cofibrantly generated
model category, and for any $u : i\rightarrow j$ morphism in $I$
we let
$$u^{*} : M_{j} \longrightarrow M_{i} \qquad
M_{j} \longleftarrow M_{i} : u_{*}$$
be a Quillen adjunction. We suppose furthermore that for
any composition $$\xymatrix{i \ar[r]^-{u} & j \ar[r]^-{v} & k}$$
one has an equality of functors
$$u^{*}\circ v^{*} = (v\circ u)^{*}$$
Such a data $(\{M_{i}\}_{i},\{u^{*}\}_{u})$
will be called, according to \cite{sh}, a \emph{($\mathbb{U}$-)cofibrantly
generated left Quillen presheaf over $I$}, and will be
denoted simply by the letter $M$.

For any cofibrantly generated left Quillen presheaf $M$ on $I$, we
consider the category $M^{I}$, of $I$-diagrams in
$M$ in the following way. Objects in $M^{I}$ are
given by the data of objects $x_{i} \in M_{i}$ for any
$i\in I$, together with morphisms
$\phi_{u} : u^{*}(x_{j}) \longrightarrow x_{i}$ for any
$u : i \rightarrow j$ in $I$, making the following diagram commutative
$$\xymatrix{
u^{*}v^{*}(x_{k}) \ar[r]^-{u^{*}(\phi_{v})} \ar[d]_-{Id} &
u^{*}(x_{j}) \ar[d]^-{\phi_{u}} \\
(v\circ u)^{*}(x_{k}) \ar[r]_-{\phi_{v\circ u}} & x_{i}}$$
for any $\xymatrix{i \ar[r]^-{u} & j \ar[r]^-{v} & k}$
in $I$. Morphisms in $M^{I}$ are simply given by families
morphisms $f_{i} : x_{i} \longrightarrow y_{i}$, such that
$f_{i}\circ \phi_{u}=\phi_{u}\circ f_{j}$ for
any $i\rightarrow j$ in $I$.

The category $M^{I}$ is endowed with a model structure for which the
fibrations or equivalences are the morphisms $f$ such that
for any $i\in I$ the induced morphism $f_{i}$ is a fibration
or an equivalence in $M_{i}$. As all model categories
$M_{i}$ are $\mathbb{U}$-cofibrantly generated, it is not hard
to adapt the general argument of \cite[11.6]{hi} in order to prove that
$M^{I}$ is also a $\mathbb{U}$-cofibrantly generated model category.

We define an object
$x \in M^{I}$ to be \emph{homotopy cartesian} if for any
$u : i\rightarrow j$ in $I$ the induced morphism
$$\mathbb{L}u^{*}(x_{i}) \longrightarrow x_{j}$$
is an isomorphism in $\mathrm{Ho}(M_{j})$. The full subcategory
of cartesian objects in $M^{I}$ will be denoted by
$M^{I}_{cart}$. \\

There exists a presheaf of categories
$(-/I)^{op}$ over $I$, having the opposite comma category
$(i/I)^{op}$ as value over the object $i$, and the
natural functor $(i/I)^{op} \longrightarrow (j/I)^{op}$ for
any morphism $j \longrightarrow i$ in $I$. For any
$\mathbb{U}$-cofibrantly left Quillen presheaf $M$ over $I$, we
define a morphism of presheaves of categories over $I$
$$M^{I} \times (-/I)^{op} \longrightarrow M,$$
where $M^{I}$ is seen as a constant presheaf of categories,
in the following way. For an object $i\in I$, the functor
$$M^{I} \times (i/I)^{op} \longrightarrow M_{i}$$
sends an object $(x,u : i\rightarrow j)$
to $u^{*}(x_{j}) \in M_{i}$, and a
morphism $(x,u : i\rightarrow j) \longrightarrow
(y,v : i\rightarrow k)$, given by a morphism $x \rightarrow y$
in $M^{I}$ and a
commutative diagram in $I$
$$\xymatrix{
i\ar[d]_-{v} \ar[rd]^-{u} \\
k \ar[r]^-{w} & j,}$$
is sent to the morphism in $M_{i}$
$$\xymatrix{
u^{*}(x_{j})\simeq v^{*}(w^{*}(x_{j})) \ar[r] & v^{*}(x_{k}).}$$
For a diagram $u : i\rightarrow j$ in $I$, the following
diagram
$$\xymatrix{
M^{I} \times (j/I)^{op} \ar[r] \ar[d] & M_{j} \ar[d]^-{u^{*}} \\
M^{I} \times (i/I)^{op} \ar[r] & M_{i}}$$
clearly commutes, showing that the above definition actually defines
a morphism
$$M^{I}\times(-/I)^{op} \longrightarrow M.$$

We let $(M^{I}_{cart})^{cof}_{W}$ be the subcategory of
$M^{I}$ consisting of homotopy cartesian and cofibrant
objects in $M^{I}$ and equivalences between them.
In the same way we consider the sub-presheaf
$M^{c}_{W}$ whose value at $i\in I$ is the subcategory
of $M_{i}$ consisting of cofibrant objects in $M_{i}$ and
equivalences between them. Note that the
functors $u^{*} : M_{j} \longrightarrow M_{i}$ being
left Quillen for any $u : i\longrightarrow j$, preserves
the sub-categories of equivalences between cofibrant objects.

We have thus defined a morphism of presheaves of categories
$$(M^{I}_{cart})^{c}_{W} \times (-/I)^{op} \longrightarrow
M^{c}_{W},$$
and we now consider the corresponding morphism of
simplicial presheaves obtained by applying the nerve functor
$$N((M^{I}_{cart})^{c}_{W}) \times N((-/I)^{op}) \longrightarrow
N(M^{c}_{W}),$$
that is considered as a morphism in the homotopy category
$\mathrm{Ho}(SPr(I))$, of simplicial presheaves over $I$.
As for any $i$ the category $(i/I)^{op}$ has a final object,
its nerve $N((i/I)^{op})$ is contractible, and therefore
the natural projection
$$N((M^{I}_{cart})^{c}_{W}) \times N((-/I)^{op}) \longrightarrow
N((M^{I}_{cart})^{c}_{W})$$
is an isomorphism in $\mathrm{Ho}(SPr(I))$. We therefore have constructed
a well defined morphism i $\mathrm{Ho}(SPr(I))$, from the
constant simplicial presheaf $N((M^{I}_{cart})^{c}_{W})$ to
the simplicial presheaf $N(M_{W}^{c})$. By adjunction this
gives a well defined morphism in $\mathrm{Ho}(SSet)$
$$N((M^{I}_{cart})^{c}_{W}) \longrightarrow
Holim_{i\in I}N(M^{c}_{W}).$$
The strictification theorem asserts that this last morphism
is an isomorphism in $\mathrm{Ho}(SSet)$. As this seems to be
a folklore result we will not include a proof.

\begin{thm}\label{tstrict}
For any $\mathbb{U}$-small category $I$ and any
$\mathbb{U}$-cofibrantly generated left Quillen
presheaf $M$ on $I$, the natural morphism
$$N((M^{I}_{cart})^{c}_{W}) \longrightarrow
Holim_{i\in I^{op}}N(M^{c}_{W})$$
is an isomorphism in $\mathrm{Ho}(SSet)$.
\end{thm}

\begin{proof} When $M$ is the constant Quillen presheaf of simplicial sets
this is proved in \cite{dk3}. The general case can treated in a similar
way. See also  \cite[Thm. 18.6]{sh} for a stronger result. \end{proof}

Let $I$ and $M$ be as in the statement of Thm. \ref{tstrict}, and
let $M_{0}$ be a $\mathbb{U}$-cofibrantly model
category. We consider $M_{0}$ as a constant left Quillen presheaf on $I$, for which
all values are equal to $M_{0}$, and all transition functors are
identities. We assume that there exists
a left Quillen natural transformation $\phi : M_{0} \longrightarrow M$. By this, we mean
the data of left Quillen functors $\phi_{i} : M_{0} \longrightarrow M_{i}$ for any
$i\in I$, such that for any $u : i\rightarrow j$ one has
$u^{*}\circ \phi_{j}=\phi_{i}$. In this case, we define a functor
$$\phi : M_{0} \longrightarrow M^{I},$$
by the obvious formula $\phi(x)_{i}:=\phi_{i}(x)$ for $x\in M_{0}$, and the
transition morphisms of $\phi(x)$ all being identities. The functor
$\phi$ is not a left Quillen functor, but preserves
equivalences between cofibrant objects, and thus possesses a left derived functor
$$\mathbb{L}\phi : \mathrm{Ho}(M_{0}) \longrightarrow \mathrm{Ho}(M^{I}).$$
One can even show that this functor possesses a right adjoint, sending
an object $x\in M^{I}$ to the homotopy limit of
the diagram in $M_{0}$, $i\mapsto \mathbb{R}\psi_{i}(x_{i})$, where
$\psi_{i}$ is the right adjoint to $\phi_{i}$.

One also has a natural transformation of presheaves of categories
$$(M_{0})_{W}^{c} \longrightarrow M^{c}_{W}$$
inducing a natural morphism of simplicial presheaves on $I$
$$N((M_{0})_{W}^{c}) \longrightarrow N(M^{c}_{W}),$$
and thus a natural morphism in $\mathrm{Ho}(SSet)$
$$N((M_{0})_{W}^{c}) \longrightarrow Holim_{i\in I^{op}}N(M^{c}_{W}).$$

\begin{cor}\label{cstrict}
Let $I$, $M$ and $M_{0}$ be as above, and assume that the functor
$$\mathbb{L}\phi : \mathrm{Ho}(M_{0}) \longrightarrow \mathrm{Ho}(M^{I})$$
is fully faithful and that its image consists of all
homotopy cartesian objects in $\mathrm{Ho}(M^{I})$. Then the induced morphism
$$N((M_{0})_{W}^{c}) \longrightarrow Holim_{i\in I^{op}}N(M^{c}_{W})$$
is an isomorphism in $\mathrm{Ho}(SSet)$.
\end{cor}

\begin{proof} Indeed, we consider the functor
$$\begin{array}{cccc}
G : & (M_{0})^{c}_{W} & \longrightarrow & (M^{I}_{cart})^{c}_{W} \\
 & x & \mapsto & \phi(Qx),
\end{array}$$
where $Qx$ is a functorial cofibrant replacement of $x$. By hypothesis,
the induced morphism on the nerves
$$N((M_{0})_{W}^{c}) \longrightarrow N(M^{I}_{cart})^{c}_{W}$$
is an isomorphism in $\mathrm{Ho}(SSet)$.
We now consider the commutative diagram in $\mathrm{Ho}(SSet)$
$$\xymatrix{
N((M_{0})_{W}^{c}) \ar[r] \ar[rd] & N(M^{I}_{cart})^{c}_{W} \ar[d] \\
 & Holim_{i\in I^{op}}N(M^{c}_{W}).}$$
The right vertical arrow being an isomorphism by Thm. \ref{tstrict}
we deduce the corollary. \end{proof}

\chapter{Representability criterion (after J. Lurie)}

The purpose of this appendix is to give a sketch of a proof of the
following special case of J. Lurie's representability theorem.
Lurie's theorem is much deeper and
out of the range of this work. We will not need it in its full generality and will
content ourselves with this special case, largely enough for
our applications.

\begin{thm}{(J. Lurie, see \cite{lu})}\label{tII-4app}
Let $F$ be a $D^{-}$-stack. The following conditions are equivalent.
\begin{enumerate}
\item $F$ is an $n$-geometric $D^{-}$-stack.
\item $F$ satisfies the following three conditions.
\begin{enumerate}
\item The truncation $t_{0}(F)$ is an Artin 
$(n+1)$-stack.
\item $F$ has an obstruction theory relative to $sk-Mod_{1}$.
\item For any $A\in sk-Alg$, the natural morphism
$$\mathbb{R}F(A) \longrightarrow Holim_{k}\mathbb{R}F(A_{\leq k})$$
is an isomorphism in $\mathrm{Ho}(SSet)$.
\end{enumerate}
\end{enumerate}
\end{thm}

\begin{proof}[Sketch of proof] The only if part is the easy part. $(a)$
is true by Prop. \ref{pII-5} and $(b)$ by Cor. \ref{cII-4}. For $(c)$ one proves the
following more general lemma.

\begin{lem}
Let $f : F \longrightarrow G$ be an $n$-representable morphism. Then for
any $A\in sk-Alg$, the natural square
$$\xymatrix{
\mathbb{R}F(A) \ar[r] \ar[d] & Holim_{k}\mathbb{R}F(A_{\leq k}) \ar[d] \\
\mathbb{R}G(A) \ar[r]  & Holim_{k}\mathbb{R}G(A_{\leq k})}$$
is homotopy cartesian.
\end{lem}

\begin{proof} We prove this by induction on $n$. For $n=-1$, one reduces
easily to the case of a morphism between representable $D^{-}$-stacks, for which
the result simply follows from the fact that
$A\simeq Holim_{k}A_{\leq k}$. Let us now assume that $n\geq 0$ and the
result prove for $m<n$. Let $A\in sk-Alg$ and
$$x\in \pi_{0}(Holim_{k}\mathbb{R}G(A)\times_{\mathbb{R}G(A_{\leq k})}^{h}\mathbb{R}F(A_{\leq k}))$$
with
projections
$$x_{k}\in \pi_{0}(\mathbb{R}G(A)\times_{\mathbb{R}G(A_{\leq k})}^{h}\mathbb{R}F(A_{\leq k})).$$
We need to prove that the homotopy fiber $H$ of
$$\mathbb{R}F(A) \longrightarrow Holim_{k}\mathbb{R}G(A)\times_{\mathbb{R}G(A_{\leq k})}^{h}
\mathbb{R}F(A_{\leq k})$$
at $x$ is contractible. Replacing $F$ by $F\times_{G}^{h}X$ where
$X:=\mathbb{R}\underline{Spec}\, A$, and $G$ by $X$, one can assume that
$G$ is a representable stack and $F$ is an $n$-geometric stack.
As $G$ is representable, the morphism
$$\mathbb{R}G(A) \longrightarrow Holim_{k}\mathbb{R}G(A_{\leq k})$$
is an equivalence. Therefore, we are reduced to the case
where $G=*$. The point $x$ is then a point
in $\pi_{0}(Holim_{k}\mathbb{R}F(A_{\leq k}))$, and we need to prove that the homotopy
fiber $H$, taken at $x$, of the morphism
$$\mathbb{R}F(A) \longrightarrow Holim_{k}\mathbb{R}F(A_{\leq k})$$
is contractible. Using Cor. \ref{ctII-1} one sees easily that
this last statement is local on the small \'etale site of
$A$. By a localization argument we can therefore assume that
each projection $x_{k}\in \pi_{0}(\mathbb{R}F(A_{\leq k}))$ of $x$ is the image of
a point $y_{k}\in \pi_{0}(\mathbb{R}U(A_{\leq k}))$, for some
representable $D^{-}$-stack $U$ and a smooth morphism $U \longrightarrow F$.

Using Prop. \ref{p22bis} and Lem. \ref{lpostn} we see that the homotopy fiber of the morphism
$$\mathbb{R}U(A_{\leq k+1}) \longrightarrow
\mathbb{R}U(A_{\leq k})\times_{\mathbb{R}F(A_{\leq k})}^{h}\mathbb{R}F(A_{\leq k+1})$$
taken at $y_{k}$ is equivalent to $Map_{A_{\leq k}-Mod}(\mathbb{L}_{U/F,y_{k}},\pi_{k+1}(A)[k+1])$.
Cor. \ref{cpII-11} then implies that when $k$ is big enough, the homotopy fibers
of the morphisms
$$\mathbb{R}U(A_{\leq k+1}) \longrightarrow
\mathbb{R}U(A_{\leq k})\times_{\mathbb{R}F(A_{\leq k})}^{h}\mathbb{R}F(A_{\leq k+1})$$
are simply connected, and thus this morphism is surjective on connected
components.
This easily implies that the points $y_{k}$ can be thought as a point
$y\in \pi_{0}(Holim_{k}\mathbb{R}U(A_{\leq k}))$ whose image in
$\pi_{0}(Holim_{k}\mathbb{R}D(A_{\leq k}))$ is equal to $x$.

We then consider the diagram
$$\xymatrix{
\mathbb{R}U(A) \ar[r] \ar[d] & Holim_{k}\mathbb{R}U(A_{\leq k}) \ar[d] \\
\mathbb{R}F(A) \ar[r] & Holim_{k}\mathbb{R}F(A_{\leq k}).}$$
By induction on $n$ we see that this diagram is homotopy cartesian, and that the
top horizontal morphism is an equivalence. There, the morphism induced on the
homotopy fibers of the horizontal morphisms is an equivalence, showing that
$H$ is contractible as required. 
\end{proof}

Conversely, let $F$ be a $D^{-}$-stack satisfying conditions $(a)-(c)$ of
\ref{tII-4app}. The proof goes by induction on $n$. Let us first
$n=-1$. We start by a lifting lemma.

\begin{lem}\label{liftatlas}
Let $F$ be a $D^{-}$-stack satisfying the conditions
$(a)-(c)$ of Thm. \ref{tII-4app}. Then, for any
affine scheme $U_{0}$, and any
\'etale morphism $U_{0} \longrightarrow t_{0}(F)$, there exists a
representable $D^{-}$-stack $U$, a morphism $u : U\longrightarrow F$,
with $\mathbb{L}_{F,u}\simeq 0$, and a
homotopy cartesian square in $k-D^{-}Aff^{\sim,\textrm{\'et}}$
$$\xymatrix{
U_{0} \ar[r] \ar[d] & t_{0}(F) \ar[d] \\
U \ar[r] & F.}$$
\end{lem}

\begin{proof} We are going to construct by induction a sequence of representable $D^{-}$-stacks
$$\xymatrix{U_{0} \ar[r] & U_{1} \ar[r] \dots & U_{k} \ar[r] & U_{k+1} \ar[r] \dots & F}$$
in $k-D^{-}Aff^{\sim,\textrm{\'et}}/F$ satisfying the following properties.
\begin{itemize}
\item One has $U_{k}=\mathbb{R}\underline{Spec}\, A_{k}$
with $\pi_{i}(A_{k})=0$ for all $i>k$.
\item The
corresponding morphism $A_{k+1} \longrightarrow A_{k}$ induces isomorphisms
on $\pi_{i}$ for all $i\leq k$.
\item The morphism $u_{k} : U_{k} \longrightarrow F$ are such that
$$\pi_{i}(\mathbb{L}_{U_{k}/F,u_{k}})=0 \; \forall \; i\leq k+1.$$

\end{itemize}

Assume for the moment that this sequence is constructed, and let
$A:=Holim_{A_{k}}$, and $U:=\mathbb{R}\underline{Spec}\, A$.
The points $u_{k}$, defines a well defined point in
$\pi_{0}(Holim_{k}\mathbb{R}F(A_{k}))$, which by condition
$(d)$ induces a well defined morphism of stacks
$u : U\longrightarrow F$. Clearly, one has a homotopy cartesian square
$$\xymatrix{
U_{0} \ar[r] \ar[d] & t_{0}(F) \ar[d] \\
U \ar[r] & F.}$$
Let $M$ be any $A$-module. Again, using
condition $(d)$, one sees that
$$\mathbb{D}er_{F}(U,M)\simeq Holim_{k}\mathbb{D}er_{F}(U_{k},M_{\leq k})\simeq Holim_{k}
Map_{A_{k}-Mod_{s}}(\mathbb{L}_{U_{k}/F,u_{k}},M_{k})\simeq 0.$$
This implies that $\mathbb{L}_{U/F,u}=0$.

It remains to explain how to construct the sequence of $U_{k}$. This is done
by induction. For $k=0$, the only thing to check is that
$$\pi_{i}(\mathbb{L}_{U_{0}/F,u_{0}})=0 \; \forall \; i\leq 1.$$
This follows easily from Prop. \ref{p22bis} and the fact that
$u_{0} : U_{0} \longrightarrow t_{0}(F)$ is \'etale.

Assume now that all the $U_{i}$ for $i\leq k$ have been constructed.
We consider $u_{k} : U_{k} \longrightarrow F$, and the natural morphism
$$\mathbb{L}_{U_{k}}\longrightarrow \mathbb{L}_{U_{k}/F,u_{k}}
\longrightarrow (\mathbb{L}_{U_{k}/F,u_{k}})_{\leq k+2}=N_{k+1}[k+2],$$
where $N_{k+1}:=\pi_{k+2}(\mathbb{L}_{U_{k}/F,u_{k}})$. The morphism
$$\mathbb{L}_{U_{k}}\longrightarrow N_{k+1}[k+2]$$
defines a square zero extension of $A_{k}$ by $N_{k+1}[k+1]$,
$A_{k+1} \longrightarrow A_{k}$. By construction,
and using Prop. \ref{p22} there exists
a well defined point in
$$\mathbb{R}\underline{Hom}_{U_{k}/k-D^{-}Aff^{\sim,\textrm{\'et}}/F}(U_{k+1},F),$$
corresponding to a morphism in $\mathrm{Ho}(k-D^{-}Aff^{\sim,\textrm{\'et}}/F)$
$$\xymatrix{
U_{k} \ar[dr] \ar[r] & U_{k+1} \ar[d] \\
 & F.}$$
It is not hard to check that
the corresponding morphism
$\mathbb{R}\underline{Spec}\, A_{k+1}=U_{k+1} \longrightarrow F$
has the required properties. 
\end{proof}

Let us come back to the case $n=-1$. Lemma
 \ref{liftatlas} implies that there exists
a representable $D^{-}$-stack $U$, a morphism $U \longrightarrow F$
such that $\mathbb{L}_{F,u}\simeq 0$, and inducing an
isomorphism $t_{0}(U) \simeq t_{0}(F)$. Using Prop. \ref{p22bis} and
Lem. \ref{lpostn}, it is not hard to show by induction on $k$, that
for any $A\in sk-Alg$ with $\pi_{i}(A)=0$ for $i>k$, the induced morphism
$$\mathbb{R}U(A) \longrightarrow \mathbb{R}F(A)$$
is an isomorphism. Using condition $(c)$ for $F$ and $U$ one sees that
this is also true for any $A\in sk-Alg$. Therefore, $U\simeq F$, and $F$
is a representable $D^{-}$-stack. \\

We now finish the proof of Thm. \ref{tII-4app} by induction on $n$.
Let us suppose we know Thm. \ref{tII-4app} for $m<n$ and let
$F$ be a $D^{-}$-stack satisfying conditions
$(a)-(c)$ for rank $n$. By induction we see that the
diagonal of $F$ is $(n-1)$-representable and the hard point is to
prove that $F$ has an $n$-atlas. For this
we lift an $n$-atlas of $it_{0}(F)$ in the following way.
Starting from a smooth morphism $U_{0} \longrightarrow it_{0}(F)$,
we to construct by induction a sequence of representable $D^{-}$-stacks
$$\xymatrix{U_{0} \ar[r] & U_{1} \ar[r] \dots & U_{k} \ar[r] & U_{k+1} \ar[r] \dots & F}$$
in $k-D^{-}Aff^{\sim,\textrm{\'et}}/F$ satisfying the following properties.
\begin{itemize}
\item One has $U_{k}=\mathbb{R}\underline{Spec}\, A_{k}$
with $\pi_{i}(A_{k})=0$ for all $i>k$.
\item The
corresponding morphism $A_{k+1} \longrightarrow A_{k}$ induces isomorphisms
on $\pi_{i}$ for all $i\leq k$.
\item The morphism $u_{k} : U_{k} \longrightarrow F$ is such that
for any $M\in A_{k}-Mod_{s}$ with
$\pi_{i}(M)=0$ for all $i>k+1$ and $\pi_{0}(M)=0$, one has
$$[\mathbb{L}_{U_{k}/F,u_{k}},M]_{Sp(A_{k}-Mod_{s})}=0.$$
\end{itemize}
Such a sequence can be constructed by induction on $k$ using
obstruction theory as in Lem. \ref{liftatlas}. We then let $A=Holim_{k}A_{k}$, and
$U=\mathbb{R}\underline{Spec}\, A$. Then, condition
$(c)$ implies the existence of a well defined
morphism $U \longrightarrow F$, which is seen to be smooth
using Cor. \ref{cpII-11}. Using this lifting of smooth morphisms from $it_{0}$
to $F$, one produces an $n$-atlas for $F$ by lifting an $n$-atlas of
$it_{0}(F)$. 
\end{proof}

\end{appendix}

\backmatter



\begin{thebibliography}{90}

\bibitem[EGAI]{ega1} A. Grothendieck, J. Dieudonn\'e, \textit{El\'ements de G\'eom\'etrie Alg\'ebrique},
I, Springer-Verlag, New York $1971$.


\bibitem[AHS]{ahs} M. Ando, M. J. Hopkins, N. P. Strickland, \textit{Elliptic spectra, the Witten genus, and the theorem of the cube}, Inv. Math. \textbf{146}, (2001), 595-687.


\bibitem[Ar]{Ar} M. Artin, \textit{Algebraization of formal moduli I}, in ``Global Analysis'' (papers in honor of K. Kodaira),
University of Tokyo Press, 1969, p. 21-71.

\bibitem[SGA4-I]{sga4} M. Artin, A. Grothendieck, J. L. Verdier, \textit{Th\'eorie des topos et cohomologie
\'etale des sch\'emas- Tome 1}, Lecture Notes in Math \textbf{269}, Springer Verlag, Berlin, 1972.

\bibitem[SGA4-II]{sga4II} M. Artin, A. Grothendieck, J. L. Verdier, \textit{Th\'eorie des topos et cohomologie
\'etale des sch\'emas- Tome 2}, Lecture Notes in Math \textbf{270}, Springer Verlag, Berlin, 1972.

\bibitem[Ba]{ba} M. Basterra, \textit{Andr\'e-Quillen cohomology of commutative $S$-algebras},
J. Pure Appl. Algebra \textbf{144}, 1999, 111-143.

\bibitem[Ba-McC]{bamcc} M. Basterra, R. McCarthy, \textit{Gamma Homology, Topological Andr\'e-Quillen Homology and Stabilization}, Topology and its Applications \textbf{121}/3, 2002, 551-566.

\bibitem[Be1]{be1} K. Behrend, \textit{Differential Graded Schemes I: Perfect Resolving Algebras}, Preprint math.AG/0212225.

\bibitem[Be2]{be2} K. Behrend, \textit{Differential Graded Schemes II: The 2-category of Differential Graded Schemes}, Preprint math.AG/0212226.



\bibitem[Bl]{bl} B. Blander, \textit{Local projective model structure on simpicial presheaves},
$K$-theory \textbf{24} (2001) No. $3$, 283-301.



\bibitem[Ci-Ka$1$]{ck1} I. Ciocan-Fontanine, M. Kapranov,
\textit{Derived Quot schemes}, Ann. Sci. Ecole Norm. Sup. (4) \textbf{34} (2001), 403-440.

\bibitem[Ci-Ka$2$]{ck2} I. Ciocan-Fontanine, M. Kapranov,
\textit{Derived Hilbert Schemes}, preprint available at \textsf{math.AG/$0005155$}.



\bibitem[Del$1$]{del1} P. Deligne, \textit{Cat\'egories Tannakiennes},
in \textit{Grothendieck Festschrift Vol. $II$}, Progress in Math. $\mathbf{87}$, Birkhauser, Boston $1990$.

\bibitem[Del$2$]{del2} P. Deligne, \textit{Le groupe fondamental de la droite
projective moins trois points}, in \textit{Galois groups over
$\mathbb{Q}$}, Math. Sci. Res. Inst. Publ., $16$, Springer Verlag,
New York, $1989$.

\bibitem[Del-Rap]{dera} P. Deligne, M. Rapoport, \textit{Les sch\'emas de modules de courbes elliptiques}, 143-317
in \textit{Modular functions of one variable II}, LNM 349, Springer, Berlin 1973.

\bibitem[Dem-Gab]{dg} M. Demazure, P. Gabriel, \textit{Groupes alg\'ebriques, Tome $I$},
Masson \& Cie. Paris, North-Holland
publishing company, $1970$.

\bibitem[Du]{du} D. Dugger \textit{Universal homotopy theories}, Adv. Math. \textbf{164} (2001), 144-176.

\bibitem[Du2]{du2} D. Dugger \textit{Combinatorial model categories have presentations},
Adv. in Math. $\mathbf{164}$ (2001), 177-201.


\bibitem[D-K$1$]{dk1} W. Dwyer,  D. Kan, \textit{Simplicial localization of categories},
J. Pure and Appl. Algebra $\mathbf{17}$ (1980), 267-284.

\bibitem[D-K$2$]{dk2} W. Dwyer,  D. Kan, \textit{Equivalences between homotopy theories of diagrams},
in \textit{Algebraic topology and algebraic $K$-theory}, \textit{Annals of Math. Studies} $\mathbf{113}$,
Princeton University Press, Princeton, 1987, 180-205.

\bibitem[D-K$3$]{dk3} W. Dwyer, D. Kan, \textit{Homotopy commutative diagrams and their realizations}, J. Pure Appl. Algebra
\textbf{57} (1989) No. 1, 5-24.

\bibitem[DHK]{dkh} W. Dwyer, P. Hirschhorn, D. Kan, \textit{Model categories and more general abstract homotopy theory},
Book in preparation, available at \textsf{http://www-math.mit.edu/$^{\sim}$psh}.

\bibitem[DS]{ds} W. Dwyer, J. Spalinski, \textit{Homotopy theories and model categories}, Handbook of Algebraic Topology, edited by I. M. James, Elsevier, 1995, 73-126. 

\bibitem[EKMM]{ekmm} A.D. Elmendorf, I. Kriz, M.A. Mandell, J.P. May, \textit{Rings, modules, and
algebras in stable homotopy theory}, Mathematical Surveys and Monographs, vol. $47$,
American Mathematical Society, Providence, $RI$, $1997$.

\bibitem[Goe-Ja]{gj} P. Goerss, J.F. Jardine, \textit{Simplicial homotopy theory},
Progress in Mathematics, Vol. $\mathbf{174}$,
Birkhauser Verlag $1999$.

\bibitem[Goe-Hop]{gh} P. Goerss, M. Hopkins, \textit{Andr\'e-Quillen (co)homology for
simplicial algebras over simplicial operads}, Une D\'egustation Topologique [Topological Morsels]: Homotopy Theory in the Swiss Alps (D. Arlettaz and K. Hess, eds.), Contemp. Math. 265, Amer. Math. Soc., Providence, RI, 2000, pp. 41-85.

\bibitem[Go]{go} J. Gorski, \textit{Representability of
the derived Quot functor}, in preparation.

\bibitem[Gr]{gr} A.Grothendieck, \textit{Cat\'egories cofibr\'ees additives et complexe cotangent relatif},
Lecture Note in Mathematics $79$, Springer-Verlag, Berlin, 1968.

\bibitem[EGAI]{egaI} A. Grothendieck, J. Dieudonn\'e,  \textit{El\'ements de G\'eom\'etrie Alg\'ebrique $I$}, Springer Verlag, Berlin, 1971.

\bibitem[EGAIV]{egaIV-4} A. Grothendieck,  \textit{El\'ements de G\'eom\'etrie Alg\'ebrique $IV$. Etude locale des sch\'emas
et des morphismes de sch\'emas},
Publ. Math. I.H.E.S.,  \textbf{20}, \textbf{24}, \textbf{28}, \textbf{32} (1967).

\bibitem[Ha]{hak} M. Hakim, \textit{Topos annel\'es et sch\'emas relatifs},
Ergebnisse der Mathematik und ihrer Grenzgebiete, Band $64$. Springer-Verlag
Berlin-New York, $1972$.

\bibitem[Hin1]{hin} V. Hinich, \textit{Homological algebra of homotopical algebras}, Comm. in Algebra
$\mathbf{25}$ (1997), 3291-3323.

\bibitem[Hin2]{hin2} V. Hinich, \textit{Formal stacks as dg-coalgebras}, J. Pure Appl. Algebra $\mathbf{162}$
(2001), No. 2-3, 209-250.

\bibitem[Hi]{hi} P. S. Hirschhorn, \textit{Model Categories and Their Localizations}, Math. Surveys and Monographs Series 99, AMS, Providence, 2003.

\bibitem[H-S]{sh} A. Hirschowitz, C. Simpson, \textit{Descente pour les $n$-champs},
preprint available at \textsf{math.AG/$9807049$}.

\bibitem[Hol]{hol} S. Hollander, \textit{A homotopy theory for stacks},
preprint available at \textsf{math.AT/0110247}.

\bibitem[Ho1]{ho} M. Hovey, \textit{Model categories}, Mathematical surveys and monographs, Vol. $\mathbf{63}$,
Amer. Math. Soc., Providence 1998.

\bibitem[Ho2]{ho2} M. Hovey, \textit{Spectra and symmetric spectra in general model categories}, J. Pure Appl. Alg. \textbf{165} (2001), 63-127.

\bibitem[HSS]{hss} M. Hovey, B.E. Shipley, J. Smith, \textit{Symmetric spectra},
J. Amer. Math. Soc. $\mathbf{13}$ (2000), no. 1, 149-208.

\bibitem[Ill]{ill} L. Illusie, \textit{Complexe cotangent et d\'eformations I}, Lectur Notes in Mathematics \textbf{239}, Springer Verlag, Berlin, 1971. 

\bibitem[Ja1]{ja} J. F. Jardine, \textit{Simplicial presheaves}, J. Pure and Appl. Algebra $\mathbf{47}$ (1987),
35-87.

\bibitem[Ja2]{ja2} J. F. Jardine, \textit{Stacks and the homotopy theory of simplicial sheaves},
in \textit{Equivariant stable homotopy theory and related areas} (Stanford, CA, 2000).
Homology Homotopy Appl. $\mathbf{3}$ (2001), No. 2, 361-384.

\bibitem[Jo1]{jo} A. Joyal, Letter to Grothendieck.

\bibitem[Jo2]{jo2} A. Joyal, unpublished manuscript.

\bibitem[Jo-Ti]{joti} A. Joyal, M. Tierney, \textit{Strong stacks and classifying spaces},
in \textit{Category theory (Como, 1990)}, Lecture Notes in Mathematics $\mathbf{1488}$, Springer-Verlag New York, 1991, 213-236.

\bibitem[Ka2]{ka2} M. Kapranov, \textit{Injective resolutions of $BG$ and derived moduli spaces of
local systems}, J. Pure Appl. Algebra $\mathbf{155}$ (2001), No. 2-3, 167-179.



\bibitem[K-P-S]{kps} L. Katzarkov, T. Pantev, C. Simspon, \textit{Non-abelian mixed Hodge structures}, preprint
\textsf{math.AG/0006213}.



\bibitem[Ko]{ko2} M. Kontsevich, \textit{Operads and motives in deformation quantization},
Mosh\'e Flato (1937-1998), Lett. Math. Phys. $\mathbf{48}$, (1999), No. $1$, 35-72.

\bibitem[Ko-So]{kos} M. Kontsevich, Y. Soibelman, \textit{Deformations of algebras over operads and
the Deligne conjecture}, Conf\'erence Mosh\'e Flato 1999, Vol. 1 (Dijon),
255-307, Math. Phys. Stud. $\mathbf{21}$, Kluwer Acad. Publ, Dordrecht, 2000.

\bibitem[Kr-Ma]{km} I. Kriz, J. P. May, \textit{Operads, algebras, modules and motives}, Ast\'erisque $\mathbf{233}$, 1995.

\bibitem[La-Mo]{lm} G. Laumon and L. Moret-Bailly,
\textit{Champs alg\'ebriques}, A Series of Modern Surveys in Mathematics vol. $\mathbf{39}$, Springer-Verlag 2000.

\bibitem[Lu1]{lu} J.Lurie, PhD Thesis, MIT, Boston, 2004.

\bibitem[Lu2]{luelcoom} J. Lurie, \textit{A survey of elliptic cohomology}, Preprint December 2005 (available at http://www.math.harvard.edu/$\sim$lurie/papers/survey.pdf).
 
\bibitem[MMSS]{mmss} M. Mandell, J. P. May, S. Schwede, B. Shipley \textit{Model categories of diagram spectra},
Proc. London Math. Soc. $\mathbf{82}$ (2001), 441--512.

\bibitem[May]{may} J.P. May, \textit{Pairings of categories and spectra}, JPAA \textbf{19} (1980), 299-346.

\bibitem[May2]{m} J.P. May, \textit{Picard groups, Grothendieck rings, and Burnside rings of categories},
Adv. in Math. \textbf{163}, 2001, 1-16.

\bibitem[Mil]{mil} J. S. Milne, \textit{\'Etale cohomology}, Princeton University Press, 1980.

\bibitem[Min]{min} V. Minasian, \textit{Andr\'e-Quillen spectral sequence for THH},
Topology and Its Applications, \textbf{129}, (2003) 273-280.

\bibitem[MCM]{min2} R. Mc Carthy, V. Minasian, \textit{HKR theorem for smooth S-algebras}, Journal of Pure and Applied Algebra, Vol \textbf{185}, 2003, 239-258, 2003.


\bibitem[Q1]{q0} D. Quillen, \textit{Homotopical algebra},
Lecture Notes in Mathematics \textbf{43}, Springer Verlag, Berlin, 1967.

\bibitem[Q2]{q} D. Quillen, \textit{On the (co-)homology of commutative rings}, Applications
of Categorical Algebra (Proc. Sympos. Pure Math., Vol XVII, New York, 1964),
65-87. Amer. Math. Soc., Providence, P.I.

\bibitem[Q3]{q2} D. Quillen, \textit{Higher algebraic
K-theory I}, in \textit{Algebraic K-theory I-Higher K-theories},
Lecture Notes in Mathematics \textbf{341}, Springer Verlag, Berlin.

\bibitem[Re]{re} C. Rezk, \textit{Spaces of algebra structures and cohomology of operads}, Thesis 1996, available
at \textsf{http://www.math.uiuc.edu/~rezk}.



\bibitem[Ri-Rob]{rr} B. Richter, A. Robinson, \textit{Gamma-homology of group algebras and of polynomial algebras},
To appear in the ``Proceedings of the Northwestern conference'' 2002.

\bibitem[Ro]{ro} J. Rognes, \textit{Galois extensions of structured ring spectra},
Preprint \textsf{math.AT/0502183}.



\bibitem[Schw-Shi]{schw-shi} S. Schwede, B. Shipley, \textit{Stable model categories are categories of modules},
Topology $\mathbf{42}$ ($2003$), $103--153$.


\bibitem[Shi]{shi} B. Shipley, \textit{A convenient model category for commutative ring spectra}, Preprint 2002.

\bibitem[S$1$]{s2} C. Simpson, \textit{Homotopy over the complex
numbers and generalized cohomology theory}, in \textit{Moduli
of vector bundles (Taniguchi Symposium, December 1994)}, M.
Maruyama ed., Dekker Publ. (1996), 229-263.

\bibitem[S$2$]{s3} C. Simpson, \textit{A Giraud-type characterization of the simplicial categories
associated to closed model categories as $\infty$-pretopoi}, Preprint \textsf{math.AT/$9903167$}.

\bibitem[S$3$]{s4} C. Simpson, \textit{Algebraic (geometric) $n$-stacks}, Preprint  \textsf{math.AG/$9609014$}.

\bibitem[S$4$]{s5} C. Simpson, \textit{The Hodge filtration on non-abelian cohomology},
Preprint \textsf{math.AG/9604005}.

\bibitem[Sm]{sm} J. Smith, \textit{Combinatorial model categories}, unpublished.



\bibitem[Sp]{sp} M. Spitzweck, \textit{Operads, algebras and modules in model categories and motives},
Ph.D. Thesis, Mathematisches Instit\"ut, Friedrich-Wilhelms-Universit\"at Bonn (2001), available at
\textsf{http://www.uni-math.gwdg.de/spitz/}.



\bibitem[Tab]{tab} G. Tabuada, \textit{Une structure de cat\'egorie
de mod\`eles de Quillen sur la cat\'egorie des dg-cat\'egories},
Preprint math.KT/0407338.


\bibitem[To1]{to1} B. To\"en, \textit{Champs affines},
Preprint math.AG/0012219.

\bibitem[To2]{to3} B. To\"en, \textit{Homotopical and higher categorical
structures in algebraic geometry}, Hablitation Thesis available
at math.AG/0312262

\bibitem[To3]{to4} B. To\"en, \textit{Vers une interpr\'etation Galoisienne de la th\'eorie
de l'homotopie}, Cahiers de topologie et geometrie differentielle categoriques, Volume XLIII (2002),
257-312.

\bibitem[To-Va1]{tv} B. To\"en, M. Vaqui\'e, \textit{Moduli of objects
in dg-categories}, Preprint math.AG/0503269.

\bibitem[To-Va2]{sousz} B. To\"en, M. Vaqui\'e, \textit{Au-dessous de $Spec\, \mathbb{Z}$}, 
Preprint math.AG/0509684.


\bibitem[HAGI]{hagI} B. To\"en, G. Vezzosi, \textit{Homotopical algebraic geometry I: Topos theory}, Advances in Mathematics, 
\textbf{193}, Issue 2 (2005), p. 257-372.

\bibitem[To-Ve1]{msri} B. To\"en, G. Vezzosi, \textit{Segal topoi and stacks over Segal categories},
December 25, 2002, to appear in Proceedings of the Program \textit{``Stacks, Intersection theory and Non-abelian Hodge Theory''},
MSRI, Berkeley, January-May 2002 (also available as Preprint math.AG/0212330).

\bibitem[To-Ve2]{hagdag} B. To\"en, G. Vezzosi, \textit{From HAG to DAG: derived moduli spaces}, p. 175-218, in \textit{``Axiomatic, Enriched and Motivic Homotopy Theory''}, Proceedings of the NATO Advanced Study Institute, Cambridge, UK, (9-20 September 2002), Ed. J.P.C. Greenlees, NATO Science Series II, Volume 131
Kluwer, 2004.

\bibitem[To-Ve3]{newton} B. To\"en, G. Vezzosi, \emph{``Brave New'' algebraic geometry and global derived moduli spaces of ring spectra}, to appear in Proceedings of the
Euroworkshop \textit{``Elliptic Cohomology and Higher Chromatic Phenomena''}
(9 - 20 December 2002), Isaac Newton Institute for Mathematical Sciences (Cambridge, UK), H. Miller, D. Ravenel eds.
(also available as Preprint \textsf{math.AT/0309145}).

\bibitem[To-Ve4]{web} B. To\"en, G. Vezzosi, \textit{Algebraic geometry over model categories.
A general approach to Derived Algebraic Geometry},
Preprint \textsf{math.AG/$0110109$}.

\bibitem[We]{wei} C. Weibel, \textit{An introduction to homological algebra}, Cambridge Univ. Press, Cambridge, 1995.

\end{thebibliography}

\printindex

\end{document}